  \documentclass[reqno, 11 pt]{amsart}
\usepackage{amssymb}       
\usepackage{amsmath}        
\usepackage{euscript}        
\usepackage{amsthm}         
\usepackage{dsfont}     
\usepackage{amsfonts}       
\usepackage{marginnote}    
\usepackage{latexsym} 
\usepackage{amsopn} %per definire nuovi operatori "operator name"
\usepackage{geometry} 
\usepackage{amscd} %per fare diagrammi rettangolari (no frecce diagonali) e mettere le label sulle frecce 
\usepackage{mathrsfs}  
\usepackage{caption}
\usepackage{aurical}   
\usepackage[english]{babel}
\usepackage[utf8]{inputenc}
\usepackage{graphicx}
\usepackage[dvipsnames]{xcolor}
\usepackage{multicol}
\usepackage{nomencl}
\makeglossary
\makenomenclature  
\usepackage{refcount}
\usepackage{upgreek}
\usepackage{float}
\usepackage{hyperref}
\usepackage{times}
\usepackage[titletoc]{appendix}
\newcommand{\comment}[1]{}

\newcommand{\pa}{\partial}
\newcommand\norm[1]{\left\lVert#1\right\rVert}
\newcommand{\normk}[2]{\| #1 \|_{#2}^{k_0,\gamma}}
\newcommand{\di}{{\rm d}}

\newcommand{\teta}{\theta}

\newcommand{\al}{\alpha}

\newcommand{\g}{\gamma}
\newcommand{\eps}{\varepsilon}
\newcommand{\e}{{\varepsilon}}
\newcommand{\vphi}{{\varphi}}
\newcommand{\s}{{\sigma}}
\newcommand{\oo}{{\omega}}
\newcommand{\vs}{\varsigma}
\newcommand{\id}{\operatorname{Id}}

\newcommand{\op}{{\hbox{Op}}}
\newcommand{\opw}{{\hbox{Op}^{\rm W}}}
\newcommand{\opb}{{\hbox{Op}^\mathrm{B}}}
\newcommand{\opbw}{{\hbox{Op}^\mathrm{BW}}}
\providecommand{\vect}[2]{{\bigl[\begin{smallmatrix}#1\\#2\end{smallmatrix}\bigr]}}   
\providecommand{\sm}[4]{{\bigl(\begin{smallmatrix}#1&#2\\#3&#4\end{smallmatrix}\bigr)}}

\newtheorem{Teo}{Theorem}

\newtheorem{thm}{Theorem}[section]
\newtheorem*{thm*}{Theorem}
\newtheorem{prop}[thm]{Proposition}
\newtheorem{lemma}[thm]{Lemma}

\newtheorem{cor}[thm]{Corollary}
\newtheorem*{cor*}{Corollary}
\newtheorem{rmk}[thm]{Remark}
\newtheorem{ex}{Example}
\newtheorem{hypo}[thm]{Hypothesis}
\newtheorem{defn}[thm]{Definition}

\numberwithin{equation}{section}
\newcommand{\wtz}{{\widetilde{z}}}

\newcommand{\wtG}{{\widetilde{G}}}

\newcommand{\ii}{{\rm i}}

\newcommand{\om}{{\omega}}

\newcommand{\whi}{{\widehat i}}
\newcommand{\whz}{{\widehat z}}
\newcommand{\why}{{\widehat y}}

\newcommand{\wh}[1]{{\widehat{#1}}}

\newcommand{\whI}{{\widehat I}}

\newcommand{\x}{\xi}

\newcommand{\ov}{\overline}

\newcommand{\C}{{\mathbb C}}
\newcommand{\D}{{\mathbb D}}

\newcommand{\N}{{\mathbb N}}

\newcommand{\R}{{\mathbb R}}

\newcommand{\T}{{\mathbb T}}
\newcommand{\Z}{{\mathbb Z}}

\newcommand{\cA}{{\mathcal A}}
\newcommand{\cB}{{\mathcal B}}
\newcommand{\cC}{{\mathcal C}}
\newcommand{\cD}{{\mathcal D}}
\newcommand{\cE}{{\mathcal E}}
\newcommand{\cF}{{\mathcal F}}
\newcommand{\cG}{{\mathcal G}}
\newcommand{\cH}{{\mathcal H}}
\newcommand{\cI}{{\mathcal I}}

\newcommand{\cK}{{\mathcal K}}
\newcommand{\cL}{{\mathcal L}}
\newcommand{\cM}{{\mathcal M}}
\newcommand{\cN}{{\mathcal N}}
\newcommand{\cO}{{\mathcal O}}
\newcommand{\calO}{\mathcal{O}}
\newcommand{\cP}{{\mathcal P}}
\newcommand{\cQ}{{\mathcal Q}}
\newcommand{\cR}{{\mathcal R}}
\newcommand{\cS}{{\mathcal S}}
\newcommand{\cT}{{\mathcal T}}
\newcommand{\cU}{{\mathcal U}}

\newcommand{\cW}{{\mathcal W}}

\newcommand{\fm}{{\mathfrak{m}}}

\newcommand{\fM}{{\mathfrak{M}}}

\newcommand{\fr}{{\mathfrak{r}}}

\newcommand{\ta}{{\mathtt{a}}}
\newcommand{\tb}{{\mathtt{b}}}
\newcommand{\tc}{{\mathtt{c}}}
\newcommand{\td}{{\mathtt{d}}}

\renewcommand{\th}{{\mathtt{h}}}

\newcommand{\tk}{{\mathtt{k}}}

\newcommand{\tm}{{\mathtt{m}}}

\newcommand{\tq}{{\mathtt{q}}}
\newcommand{\tr}{{\mathtt{r}}}

\newcommand{\tv}{{\mathtt{v}}}
\newcommand{\tw}{{\mathtt{w}}}

\newcommand{\tC}{{\mathtt{C}}}
\newcommand{\tD}{{\mathtt{D}}}

\newcommand{\tT}{{\mathtt{T}}}

\newcommand{\tV}{{\mathtt{V}}}

\newcommand{\bL}{{\bf L}}

\newcommand{\bT}{{\bf T}}

\usepackage{bm}

\renewcommand{\d}{\td}

\newcommand{\grad}{\nabla}

\newcommand{\abs}[1]{\left| #1 \right|}
\newcommand{\absk}[2]{| #1 |_{#2}^{k_0,\gamma}}
\newcommand{\0}{{(0)}}

\newcommand{\la}{\left\langle}
\newcommand{\ra}{\right\rangle}

\newcommand{\im}{{\rm i}}
\newcommand{\jap}[1]{\langle #1 \rangle}

\newcommand{\uno}{{\mathbb I}}

\newcommand{\nnorm}[1]{{\left\vert\kern-0.25ex\left\vert\kern-0.25ex\left\vert #1 
    \right\vert\kern-0.25ex\right\vert\kern-0.25ex\right\vert}}
\usepackage{framed,enumitem} 

\renewcommand{\div}{\nabla\cdot}
\renewcommand{\t}{\tau}

\providecommand{\vect}[2]{{\bigl[\begin{smallmatrix}#1\\#2\end{smallmatrix}\bigr]}} 
  
\providecommand{\sm}[4]{{\bigl[\begin{smallmatrix}#1&#2\\#3&#4\end{smallmatrix}\bigr]}}
\newcommand{\dps}{\displaystyle}
\newcommand{\uphi}{\underline{\phi}}

\newcommand{\Lip}{\mathrm{Lip}} 
\newcommand{\kug}{{k+1,\gamma}} % MACRO per norme Whitney
\newcommand{\ompaph}{\omega\cdot \partial_\vphi}
\DeclareMathOperator{\dist}{dist}

\makeatletter

\renewcommand{\tocsection}[3]{%
\indentlabel{\@ifnotempty{#2}{\bfseries\ignorespaces#1 #2\quad}}\bfseries#3}
\def\l@subsection{\@tocline{2}{0pt}{2.5pc}{5pc}{}}
\def\l@subsubsection{\@tocline{3}{0pt}{4.5pc}{5pc}{}}
\renewcommand\tocchapter[3]{%
  \indentlabel{\@ifnotempty{#2}{\ignorespaces#2.\quad}}#3%
}
\begin{document} 
 
\title{Time Quasi-Periodic Three-dimensional Traveling Gravity Water Waves}
%Quasi-periodic traveling three-dimensional pure gravity Water Waves}

% Quasi-periodic traveling wave solutions for three-dimensional pure gravity water waves (o pure gravity water waves in 3D)
%Time quasi-periodic traveling three dimensional gravity water waves 
% 

\date{}

\author{Roberto Feola}
\address{\scriptsize{Dipartimento di Matematica e Fisica, Universit\`a degli Studi RomaTre, 
Largo San Leonardo Murialdo 1, 00146, Roma Italy}}
\email{roberto.feola@uniroma3.it}

\author{Riccardo Montalto}
\address{\scriptsize{Dipartimento di Matematica, Universit\`a degli Studi di Milano, Via Saldini 50, I-20133, Milano, Italy}}
\email{riccardo.montalto@unimi.it}

\author{Shulamit Terracina} 
\address{\scriptsize{SISSA, Via Bonomea 265, 34136, Trieste, Italy}}
\email{shterra@sissa.it}

%\thanks{ {\em Acknowledgements.} 
%The authors have been  supported by the  research project 
%PRIN 2020XBFL ``Hamiltonian and dispersive PDEs" of the 
%Italian Ministry of Education and Research (MIUR). S.T. acknowledge the support of the project ERC STARTING GRANT 2021, "Hamiltonian Dynamics, Normal Forms and Water Waves" (HamDyWWa), Project Number: 101039762.

\makeatletter
\@namedef{subjclassname@2020}{\textup{2020} Mathematics Subject Classification}
\makeatother
 


%\maketitle  
   
\begin{abstract}   
Starting with the pioneering computations of Stokes in 1847, the search of traveling waves in fluid mechanics has always been a fundamental topic, since they can be seen as building blocks to determine the long time dynamics (which is a widely open problem).  
In this paper we prove the existence of time quasi-periodic traveling wave solutions for three-dimensional  pure gravity water waves in finite depth, on flat tori, with an arbitrary number of speeds of propagation. 
These solutions are global in time, they do not reduce to stationary solutions in any moving reference frame and they are approximately given by finite sums of Stokes waves traveling with  rationally independent 
speeds of propagation.
% (they satisfies suitable quantitative non-resonance conditions)
This is a very hard  small divisors problem for Partial Differential Equations due to the fact that one deals with a dispersive quasi-linear PDE in higher dimension with a very  complicated geometry of the resonances. 
% shall describe a novel strategy that allows to perform the spectral analysis of the linearized equations (at any approximate traveling wave solutions) which is suited for higher dimensional dispersive PDEs with sublinear dispersion and the usage of conservation laws such as the momentum in order to get some cancellation of resonant terms. 
Our result is the first KAM (Kolmogorov-Arnold-Moser) result for an autonomous, dispersive, quasi-linear PDE in dimension greater than one and it is the first example of global solutions, which do not reduce to steady ones in any moving reference frame, for 3D water waves equations on compact domains. %solves a long standing open problem in Fluid Mechanics and in the community of KAM for PDEs. 
\end{abstract}

\maketitle

\setcounter{tocdepth}{3}
\tableofcontents

\section{Introduction} 

A problem of fundamental importance in mathematical fluid mechanics
regards 
the search for traveling 
 surface waves.
Since the pioneering work of Stokes \cite{stokes} in 1847, 
a huge literature has established the existence of  steady traveling  waves.
Such solutions, either periodic or localized in space,  look stationary in a moving frame, we refer to them as "Stokes waves".  
In the recent years, a huge effort has been made in order to investigate more general time periodic and quasi-periodic surface waves, which are not stationary in any moving frame. They are approximately given by the interaction of traveling waves with different speeds of propagation.
It is well known that this is a difficult small divisors problem (the spectrum of the linearized equation at the flat ocean surface accumulates to zero). The picture for two-dimensional fluids is nowadays quite clear, at least for small amplitude waves. Several methods and tools (like Nash-Moser theory, spectral theory, micro-local analysis, the theory of pseudo-differential operators, normal forms) have been implemented successfully for proving the existence of 2D periodic and quasi-periodic standing surface waves and 2D quasi-periodic traveling water waves (see \cite{IPT, AB15, BM20, BBHM, BFM, BFM2, FGtrave}). On the other hand, the problem for 3D surface waves is still widely open. 
%The 3D Water Waves equations are PDEs in higher space dimension 
%and hence much more complicated resonance phenomena occur.  
Indeed, the only available results for 3D Water Waves concern the existence of {\it bi-periodic traveling waves}, which look {\it steady in a suitable reference frame}. This is the natural extension of Stokes waves in dimension three. The first result is due to Craig-Nicholls \cite{CN} in the case of 3D gravity-capillary water waves (it is not a small divisors problem) and the second one is due to Iooss-Plotnikov \cite{IP-Mem-2009}, \cite{IP2} in the case of 3D pure-gravity water waves (it is a small divisors problem). In particular, in \cite{IP-Mem-2009} the authors consider symmetric diamond patterns by taking wave vectors into a lattice $\Gamma^*$ spanned by two vectors of $\R^2$ with the same length and velocity of propagation which is the bisectrix of these two vectors. The choice of the periodic lattice is part of the problem, indeed in \cite{IP2} they considered non symmetric patterns w.r.t. the direction of propagation of the wave (for instance considering wave vectors with different length). 
%Moreover, as we mentioned above, they look for { bi-periodic traveling waves} 
%which looks { steady in the reference frame, following the wave motion}. 

\noindent
 The goal of our paper is to construct small amplitude, 3D quasi-periodic traveling wave solutions for the three-dimensional, pure gravity water waves equations in finite depth with an arbitrary number of frequencies (bigger than the space dimension),  ``for a large measure set" of the depth parameter 
 \begin{equation}\label{mathtt h1 h2 fissati}
 \mathtt h \in [\mathtt h_1, \mathtt h_2] \quad \text{for some fixed} \quad 0 < \mathtt h_1 < \mathtt h_2\,. 
 \end{equation}
 These traveling wave solutions {\it are not steady in any reference frame},
 generalizing \cite{CN, IP-Mem-2009, IP2}. We refer to such solutions as {\bf quasi-periodic three-dimensional Stokes Waves.}
 As domain we consider a generic torus $\T^2_\Gamma = \R^2 / \Gamma$ 
 where $\Gamma$ is a lattice spanned by a basis $\{ {\bf e}_1, {\bf e}_2 \}$ of $\R^2$ and $\Gamma^*$ is the corresponding dual lattice spanned by $\{ {\bf K}_1, {\bf K}_2 \}$. 
 We make on the basis ${\bf K}_1, {\bf K}_2$ a ``genericity assumption" (see Hypothesis \ref{hyp:lattice} below). 
 Due to the flexibility of our approach, we actually work in  space dimension $d \geq 2$ with no
substantial  differences, 
 even if the physical dimension is $d = 2$.
 
 \noindent
 This result is a very hard small divisors problem for Partial Differential Equations. 
 The existence of quasi-periodic solutions for quasi-linear
 PDEs in dimension greater or equal than two is still widely open. The major difficulty is that the 3D water waves system is a dispersive PDE in higher space dimension and containing derivatives of maximal order in the nonlinearity (quasi-linear PDE), hence very strong resonance phenomena occur. %We mention that previous results on construction of traveling waves for 3D Water Waves equations are due to Craig-Nicholls and Iooss-Plotnikov

The theorem that we prove is given below (informal statement).
\begin{Teo}{\bf (Quasi-periodic 3D traveling water waves in finite depth).}
\label{thm:mainmain}
For a generic choice of the lattice $\Gamma^*$ (in the sense of the hypothesis \ref{hyp:lattice}), for any choice of wave vectors   $\bar{\jmath}_{i} \in \Gamma^*$,  $i=1\,,\ldots,\nu$ with $0 < |\bar{\jmath}_1| < |\bar{\jmath}_2| < \ldots < |\bar{\jmath}_\nu|$ and for a large measure set of depths 
$\mathtt{h}$, there exist a diophantine frequency vector $\omega=\omega(\mathtt{h})=
(\omega_{\bar{\jmath}_{1}}(\mathtt{h}), \ldots,\omega_{\bar{\jmath}_{\nu}}(\mathtt{h}))$
and a non-trivial small amplitude, linearly stable quasi-periodic traveling waves solution, with frequency $\omega$,
of the three dimensional pure gravity water waves problem.
% in finite depth.
\end{Teo}

It is well know that the problem of studying the long time dynamics for the 3D water waves system
for periodic initial data 
 is essentially open. The theorem above provides the first example of global in time
solutions whose time dependence is not trivial, i.e. for instance they do not 
reduce to stationary ones in any moving reference frame.
These solutions have small amplitude and its Sobolev norm its controlled uniformly in time.
Moreover, they generalize 
the \emph{asymmetrical three-dimensional traveling waves } constructed in 
\cite{IP2} where the authors deal with solutions  with frequencies $(\omega_1, \omega_2)$, whereas in our case 
we allow an arbitrary number of wave vectors and frequencies.
As already said, a small divisor problem arises (both in presence  or not of surface tension).
In the pure gravity case considered here, small divisors are particularly dangerous.
This is due to the higher-dimensional nature of the problem, where time and space are independent, 
and the dispersion relation, in absence of surface tension, is very weak.
In addition to this, the resonant phenomena are particularly strong also due 
to the presence of the derivatives in the nonlinearity which produces 
a strong perturbation (at the highest order) of the eigenvalues of the linearized operator.
The crucial difficulty in order to control resonances is that, in our context, we do not have 
a sharp asymptotic expansion of the eigenvalues in decreasing powers of the spectral parameter, which instead is typical in one space dimensional problems. This is also due to the complicated structure of the nonlinearity
which depends on  the quite involved pseudo-differential structure of the Dirichlet-Neumann operator, which is much more simple in one dimension.

%which are
%particularly difficult to control in the higher dimensional framework.
%In particular, in the pure gravity case considered here
%such small divisors are particularly dangerous since
%it is an higher dimension problem in which time and space are 
%independent and the dispersion relation is very weak in absence of surface tension.
%
%, due to the \emph{weak} dispersion relation of the equation,
%the differences of the eigenvalues of the linearized equation at the flat surface
%accumulates to zero very fast.
In order to control the space-time resonances 
we use the depth parameter to impose non-resonance conditions and we take
wave vectors in a non-resonant lattice (see Hyp. \ref{hyp:lattice}).
The main novelty is a combination of such non-resonance conditions with a refined use of the momentum conservation 
in order to prove a \emph{normal form reduction}  of the linearized operator at any approximate solution and get suitable  ``tame'' estimates on its inverse. 
The invertibility is obtained by a reducibility argument that consists into two main steps:
$(i)$ a pseudo-differential reduction in decreasing orders that transforms
the linearized operator to a \emph{Fourier multiplier} up to (arbitrarily) regularizing in space remainders;
$(ii)$ a KAM-like reducibility scheme (for \emph{bounded} operators) based on ``second order Melnikov'' conditions on the eigenvalues, which completely diagonalize the operator (both in time and space).
This approach has been successfully implemented for quasi-linear one-dimensional  PDEs. The sharp spectral analysis of the linearized equations is exactly the main obstacle in dealing with quasi-linear PDEs in higher space dimension and it requires substantially new ideas. 
Step $(i)$ of this reduction is highly non-trivial. Indeed, up to now, the few results in this direction in high space dimension
are confined to small perturbations of  non-resonant ``transport type'' equations, whereas 
we deal with \emph{singular} perturbations (of order one) of the dispersive 
operator $\omega(\x)\sim |\x|^{1/2}$ which has much more complicated structure of the resonances. 

This very strong issue can be interpreted 
even from the point of view of classical mechanics. Indeed, in the context of perturbation theory for Hamiltonian systems, one usually wants to analyze the global dynamics of quasi-integrable Hamiltonian systems of the form $H(x, \xi) = H_0(\xi) + \e P(x, \xi)$, $(x, \xi) \in \T^d_\Gamma \times \R^d$ ($x \in \T^d_\Gamma$ being the angle variable and $\xi \in \R^d$ being the action variable). This is usually achieved by implementing global perturbation theory, namely trying to eliminate the $x$-dependence from the perturbation (globally in $\xi \in \R^d$ and not only locally) up to an arbitrary small correction. It is known that this is much easier in the case where $H_0(\xi) = \zeta \cdot \xi$ where $\zeta \in \R^d$ is a given irrational vector whereas if $\nabla_\xi H_0(\xi)$ is non-constant in $\xi$, the non-resonance conditions required depend on the point of the phase space, since one has to provide suitable lower bounds on $\nabla_\xi H_0(\xi) \cdot k$, $k \in \Z^d \setminus \{ 0 \}$. 
%it is developing perturbation theory of an Hamiltonian 
%of the form $\zeta\cdot\x+\e P(x,\x)$ with $\zeta$  a diophantine vector  
%is much more easier  than considering an Hamiltonian $H(x,\x)=H_0(\x)+\e P(x,\x)$
%where $\nabla_\x H_0(\x)$ is not a constant function. In such a case 
%there would be resonances depending on the point of the phase space.
In order to overcome this issue, we develop a  normal form strategy in higher space dimension
based on the \emph{sub-linear} nature of the dispersion relation combined with a 
sharp use of pseudo-differential calculus and the use of conservation of momentum 
which allows to remove exact dangerous resonances. 

Regarding the \emph{perturbative} diagonalization scheme of step $(ii)$, we remark that
very few results are known in literature in this direction. Mostly, KAM results on $\T^{d}$, 
$d\geq 2$, are restricted to \emph{semilinear}
PDEs, whereas in the case of 
\emph{unbounded} perturbations they are available only for 
very \emph{well-cooked} models. The main difficulties are related with possible multiplicity of the eigenvalues 
which lead to very ``weak'' lower bounds in the second order Melnikov conditions required for the diagonalization 
scheme. This happens in particular, in our case where
 differences of eigenvalues  \emph{accumulate} to zero very fast. 
 We control the loss of derivatives (in time and space) arising from such weak
 Melnikov conditions by taking advantage of the  reduction of step $(i)$ 
 and of the conservation of momentum. As a byproduct of this reducibility approach we
 obtain that the linearized operator at the (approximate) quasi-periodic solutions is diagonal with purely imaginary spectrum
 in the directions normal to the embedded torus. This provides the linear stability of solutions.

In the next parts of this introduction, we shall give the precise statement of our result and we shall explain more in detail the difficulties and the main novelties and ideas of the proof.

%To the best of our knowledge the above theorem is the first existence result concerning blbla
%
%
% \textcolor{red}{Ripensare quasta ultima parte}
% This is the first KAM result for a dispersive, autonomous quasi-linear PDE in dimension greater or equal than two. More precisely, the major difficulties are: $(i)$ The dispersion relation is sublinear $\omega(j) \simeq |j|^{\frac12}$ where $j$ is in a discrete set of $\R^d$, hence the differences $\omega(j) - \omega(j')$ can be zero for many $j, j'$ or it can accumulate to zero very fast. 
% 
% difficulties
%\red{Commenti su coronamenteo di tanti sforzi}

\bigskip

{\bf Acknowledgements}

%\vspace{0.9em}
%{\bf Acknowledgements}.
R. Montalto  is  supported by the ERC STARTING GRANT 2021 “Hamiltonian Dynamics, Normal Forms and Water Waves” (HamDyWWa), Project Number: 101039762. 

S. Terracina  is supported by the European Union ERC CONSOLIDATOR GRANT 2023 GUnDHam, Project Number: 101124921. 

 The Views and opinions expressed are however those of the authors only and do not necessarily reflect those of the European Union or the European Research Council. Neither the European Union nor the granting authority can be held responsible for them.

R. Feola and R. Montalto acknowledge the support of PRIN 2022 “Turbulent effects vs Stability in Equations from Oceanography” (TESEO), project number: 2022HSSYPN. 

R. Feola is supported by GNAMPA-INdAM project "Stable and unstable phenomena in propagation of waves in dispersive media" CUP-E5324001950001

S. Terracina is supported
GNAMPA-INdAM project "Deterministic and probabilistic evolution out-of-equilibrium Hamiltonian systems" 
CUP-E5324001950001

\noindent
During the writing of the present paper we benefited a lot of suggestions and discussions. 

We warmly thank Professors Pietro Baldi, Massimiliano Berti, Zaher Hani, Alberto Maspero and Michela Procesi  for their continuous encouragements in pursuing this difficult project and their valuable comments and suggestions that allowed to improve the presentation of our result. We also warmly thank Federico Murgante for many discussions concerning the sections of the Dirichlet-Neumann operator.

\subsection{Formulation of the problem and main result}
Let $\Gamma$ be a lattice of dimension 
$d\geq1$ in $\mathbb{R}^{d}$, 
with basis ${\bf e}_1, {\bf e}_2, \ldots, {\bf e}_d$, namely
\begin{equation}\label{def:lattice}
\Gamma:=\Big\{\sum_{i=1}^{d}k_i {\bf e}_i\; :\; k_1,\ldots,k_{d}\in \mathbb{Z}\Big\}\,,
\qquad{\rm and \; define}
\qquad
\mathbb{T}^{d}_{\Gamma}:=\mathbb{R}^{d}/\Gamma\,.
\end{equation}
%and define 
%\begin{equation}\label{torogamma}
%\mathbb{T}^{d}_{\Gamma}:=\mathbb{R}^{d}/\Gamma\,.
%\end{equation}
The dual lattice $\Gamma^{*}$ is defined by
\begin{equation}\label{def:lattice2}
\Gamma^{*}:=\Big\{b\in \mathbb{R}^{d}\;:\; b\cdot k\in 2\pi \mathbb{Z}\,,\;\; \forall k\in\Gamma\Big\}\,.
\end{equation}
We denote by $\{ {\bf K}_1, \ldots, {\bf K}_d \}$ the dual basis of $\Gamma^*$. 
We consider an
incompressible and irrotational perfect fluid, under the action of gravity
occupying, at time
$t$, a three dimensional domain with finite depth $\mathtt{h}>0$, given by 
\begin{equation}\label{Deta}
{\mathcal D}_{\eta} := 
\big\{ (x,z)\in \T_{\Gamma}^{d} \times\R \, ;  \ - \th <z<\eta(t,x) \big\}\,,
\end{equation} 
where $\eta$ %:\R\times\T\to \R$ 
is a smooth function. 
Clearly, the physical problem is with $d=2$, but our methods (as we will see) apply to any $d\geq 1$.
So we decide to keep $d$ as a free parameter to underline that the results hold in any dimension. 
The velocity field in the time dependent domain 
$ {\mathcal D}_{\eta} $ is the gradient of a harmonic function $\underline{\phi}$, called the velocity potential. 
The time-evolution of the fluid is determined by a system of equations for 
the two functions $(t,x)\to \eta(t,x) $, $ (t,x,z)\to \underline{\phi}(t,x,z)$.
Following Zakharov \cite{Zak1} and Craig-Sulem \cite{CrSu}, we denote by
$\psi(t,x) = \underline{\phi}(t,x,\eta(t,x))$ the restriction 
of the velocity potential to the free interface.
Given the shape $\eta(t,x)$ of the domain $ {\mathcal{D}}_\eta $ 
and the Dirichlet value $\psi(t,x)$ of the velocity potential at the top boundary, one can recover 
$\underline{\phi}(t,x,z)$ as the unique solution of the elliptic problem 
\begin{equation} \label{BoundaryPr}
(\pa_{zz}+\Delta) \underline{\phi} = 0  \ \text{in } 
{\mathcal{D}}_\eta  \, , \quad 
\partial_z \underline{\phi}= 0  \ \text{at } y = - \th \, , \quad 
\underline{\phi} = \psi \ \;\; \text{on }\;\; \{z= \eta(t,x)\}. 
\end{equation}
The $(\eta,\psi)$ variables then satisfy the gravity water waves system
\begin{equation}\label{eq:113}
\left\{\begin{aligned}
    \partial_t \eta &= G(\eta)\psi \\
\partial_t\psi &=  -g\eta  -\frac{1}{2} |\nabla\psi|^2 
+  \frac{1}{2}\frac{(\nabla\eta\cdot  \nabla\psi + G(\eta)\psi)^2}{1+|\nabla\eta|^2}
\end{aligned}\right.
\end{equation}
where $g>0$ is the acceleration of gravity, $ G(\eta)$ is the 
Dirichlet-Neumann operator, defined by 
\begin{equation}\label{eq:112a}
 G(\eta)\psi:=G(\eta,\th) \psi:= \sqrt{1+|\nabla\eta|^2}(\partial_n\underline{\phi})_{| z=\eta(t,x)}
 = (\partial_z\uphi -\nabla\eta\cdot \nabla\uphi)_{| z=\eta(t,x)}
\end{equation}
and $n$ is the outward unit normal at the free interface $y= \eta(t,x)$. 
Without loss of generality, we set the gravity constant to $g = 1$.

It was first observed by Zakharov \cite{Zak1} that \eqref{eq:113} is a Hamiltonian system with respect to the symplectic form $d\psi\wedge d\eta$ and it can be written as
\begin{equation}\label{HS}
\begin{aligned}
& \qquad \pa_t \eta = \nabla_\psi H (\eta, \psi) \, , \quad  \pa_t \psi = - \nabla_\eta H (\eta, \psi)  \, , 
\end{aligned}
\end{equation}
where $ \nabla $ denotes the $ L^2 $-gradient, with Hamiltonian
\begin{equation}\label{HamiltonianWW}
H(\eta, \psi) := \frac12 \int_{\T^{d}_{\Gamma}} \psi \, G(\eta ) \psi \, dx 
+ \frac{1}{2} \int_{\T^{d}_{\Gamma}} \eta^2  \, dx
\end{equation}
given by the sum of the kinetic 
and potential energy of the fluid.
The invariance of the system \eqref{eq:113} in the $y$ and $x$ 
variable implies the existence of $d + 1$ prime integrals, respectively the ``mass'' $\int_{\T_{\Gamma}^d} \eta \, dx$ 
and the momenta
\begin{equation}\label{Hammomento}
\mathtt{M}_{i}:= \int_{\T^{d}_{\Gamma}} \eta_{x_i} (x) \psi (x)  \, dx\,,\qquad i=1,\ldots,d\,.
\end{equation}

The Hamiltonian \eqref{HamiltonianWW}  is defined on the spaces
\begin{equation}\label{phase-space}
(\eta, \psi) \in H^s_0 (\T^{d}_{\Gamma}; \mathbb{R}) 
\times {\dot H}^s (\T^{d}_{\Gamma}; \mathbb{R})  
\end{equation}
$ s\in \R $, 
where 
${\dot H}^s (\T^{d}_{\Gamma};\mathbb{R})
:= 
H^s (\T^{d}_{\Gamma};\mathbb{R}) \slash{ \sim}$ 
is the homogeneous Sobolev space
obtained by the equivalence relation
$\psi_1 (x) \sim \psi_2 (x)$ if and only if $ \psi_1 (x) - \psi_2 (x) $ is a constant\footnote{The 
fact that $ \psi \in \dot H^s(\T^{d}_{\Gamma};\mathbb{R}) $ is coherent with the fact 
that only the velocity field $ \nabla_{x,y} \Phi $ has physical meaning, 
and the velocity potential $ \Phi $ is defined up to a constant.
For simplicity of notation 
we denote the equivalence class $ [\psi] $ by $ \psi $ and,  
since the quotient map induces an isometry of $ {\dot H}^s (\T^{d}_{\Gamma};\mathbb{R}) $ onto 
$ H^s_0 (\T^{d}_{\Gamma};\mathbb{R}) $, 
we will conveniently identify $ \psi $ with a function with zero average.},
and $H^s_0(\T^{d}_{\Gamma};\mathbb{R})$ is the subspace of 
$H^s(\T^{d}_{\Gamma};\mathbb{R})$ of zero average functions.

\vspace{0.4em}
%\subsection{Main result} 
\noindent
{\bf Linearized water waves at the flat surface.}
Small amplitude solutions are close to the solutions of the linearized system of \eqref{eq:113}
at the equilibrium $(\eta,\psi)=(0,0)$, namely
\begin{equation} \label{eq:113linear}
\begin{cases}
    \partial_t \eta =G(0)\psi \cr
\partial_t\psi = \dps - \eta  
\end{cases}\,,
\qquad G(0)= |D|\tanh(\mathtt{h}|D|)\,.
\end{equation}
Here $G(0)$ is 
the Dirichlet-Neumann operator at the flat surface $\eta=0$.
The Hamiltonian of the system \eqref{eq:113linear} is
\begin{equation}\label{HamLinearReal}
H_{L}:=\frac{1}{2}({\bf \Omega} u,u)_{L^{2}}=
 \frac12 \int_{\T^{d}_{\Gamma}} \psi \, G(0) \psi \, dx 
+ \frac{1}{2} \int_{\T^{d}_{\Gamma}} \eta^2  \, dx\,,
\qquad u=\vect{\eta}{\psi}\,,\quad {\bf \Omega}=\sm{1}{0}{0}{G(0)} \,.
\end{equation}
Let us introduce the Fourier multiplier
\begin{equation}\label{wild}
\begin{aligned}
\Omega(D)&:=\Omega_{\mathtt{h}}(D):=\sqrt{|D|\tanh(\mathtt{h}|D|)}\,,
\\{\rm i.e.}\quad  
\Omega(D)e^{\ii j\cdot x}&=\Omega(j)e^{\ii j\cdot x}:=
\Omega(j;\mathtt{h})e^{\ii j\cdot x}:=
\sqrt{|j|\tanh(\th |j|)}
e^{\ii j\cdot x}\,,\qquad \forall\, j\in \Gamma^{*}\,.
\end{aligned}
\end{equation}
In the complex symplectic variables 
the system \eqref{eq:113linear}  reads
\begin{equation}\label{blackwhite}
\pa_{t}u=-\ii \Omega(D) u\,,\qquad 
 u=\frac{1}{\sqrt{2}}\big(\Omega^{-1/2}(D)\eta+\ii \Omega^{1/2}(D)\psi\big)\,,
\end{equation}
whose solutions %, recalling  \eqref{solLineari}, 
are given by 
\begin{equation}\label{solLinearicomplex}
\begin{aligned}
u(t,x)&=\sum_{j\in \Gamma^*\setminus\{0\}}
u_{j}(0)e^{-\ii\Omega(j)t+\ii jx}\,,
\qquad u_j(0):=\frac{1}{\sqrt{2}}\big( \Omega^{-\frac{1}{2}}(j)\eta_j(0)
+\ii\Omega^{\frac{1}{2}}(j)\psi_{j}(0) \big)\,.
\end{aligned}
\end{equation}
We note that that they can be either periodic, quasi-periodic or almost-periodic 
depending on the Fourier support and on the arithmetic properties of
 $\Omega(j)=\Omega(j; \th)$, which is
the frequency of oscillation of the $j$-th mode.
We refer to the map 
\begin{equation}\label{dispersionLaw}
j\to \Omega(j)=\sqrt{|j|\tanh(\th |j|)}\,,\quad j\in \Gamma^*\setminus\{0\}\,,
\end{equation}
as the \emph{dispersion law} of \eqref{eq:113linear}.
%
%
%The solution of the linear system \eqref{eq:113linear} 
%belonging to
%$H^0_0 (\T^{2}_{\Gamma}; \mathbb{R}) 
%\times {\dot H}^0 (\T^{2}_{\Gamma}; \mathbb{R})  $
%are 
%\begin{equation}\label{solLineari}
%\begin{aligned}
%\eta(t,x)&=\sum_{j\in \Gamma^{*}\setminus\{0\}}
%\Big(
%\eta_{j}(0)\cos(\Omega(j) t)+
%\Omega(j)\psi_{j}(0)\sin(\Omega(j) t)\Big)e^{\ii jx}\,,\\
%\psi(t,x)&=\sum_{j\in \Gamma^*\setminus\{0\}} \Big( 
%\psi_{j}(0)\cos(\Omega(j)t)- 
%\Omega^{-1}(j)\eta_{j}(0)\sin(\Omega(j)t)\Big)e^{\ii jx}
%\end{aligned}
%\end{equation}
%with
%coefficients satisfying $\ov{\eta_{-j}(0)}=\eta_{j}(0)$
%and $\ov{\psi_{-j}(0)}=\psi_{j}(0)$.

\vspace{0.5em}
\noindent
{\bf Traveling quasi-periodic solutions.}
We shall construct solutions for \eqref{eq:113} 
which are quasi-periodic and traveling according to the following definition.

\begin{defn}{\bf (Quasi-periodic traveling waves).}\label{def:quasitravelling}
Let $\nu\in\N$ and let $\omega\in \R^{\nu}$ be an irrational vector, i.e.
$\omega\cdot\ell\neq0$ for any $\ell\in\Z^{\nu}\setminus\{0\}$.
Consider also 
$\tv_1, \dots, \tv_\nu \in\Gamma^{*}$  and 
a function
\[
\begin{aligned}
(\eta,\psi) : \R\times\T^{d}_{\Gamma} & \to \R^{2},\quad 
( t,x)  \mapsto ({\eta}(t, x), {\psi}(t, x))\,.
 \end{aligned}
\]

\noindent
$(i)$ We say that $(\eta,\psi)$ 
 is a \emph{quasi-periodic} wave with frequency vector $\omega$
 if there exists a function $U : \T^{\nu}\times\T_{\Gamma}^{d}\to \R^{2}$
 such that 
 \begin{equation}\label{quasiperiodicWave}
 (\eta(t,x),\psi(t,x))=U(\omega t,x)\,.
 \end{equation}

\noindent
$(ii)$ We say that a quasi-periodic wave $(\eta,\psi)$ is \emph{traveling}
with ``wave vectors'' $\tv_i$, $i=1,\ldots,\nu$,
if there exists a function 
$\breve U : \T^{\nu}\to \R^{2}$ such that\footnote{Here $\mathtt V$ is represented by a matrix with $i$-th column $\mathtt{V}^{i}=(\tv_{j,i})_{j=1,\ldots,\nu}$ and we set
$\tv_{j}=(\tv_{j,i})_{i=1,\ldots,d}$.} 
\begin{equation}\label{travellingWaves}
 %(\eta(t,x),\psi(t,x))=
 \begin{aligned}
&  U(\omega t,x)=  \breve{U}( \omega_1  t- \tv_{1}  \cdot x ,\ldots, \omega_\nu t- \tv_\nu \cdot x ) = \breve{U}(\omega t- \mathtt{V}x)
 \,, \\
 & \mathtt V x := (\mathtt v_1 \cdot x, \mathtt v_2 \cdot x, \ldots, \mathtt v_\nu \cdot x) \qquad  
 \end{aligned}
\end{equation}

\noindent
$(iii)$  A quasi-periodic traveling wave
%, with frequency and velocity vectors $\omega,\tV^{(i)}$, $i=1,\ldots,2$ respectively, 
solving equations \eqref{eq:113}
is called \emph{quasi-periodic traveling solution.}
\end{defn}

Thanks to the 
irrationality of $\omega\in \R^{\nu}$, 
one can check that $(\eta,\psi)$ is a 
quasi-periodic traveling solution if and only if 
the embedding $\T^{\nu}\ni\vphi\mapsto U(\vphi,x)$ in \eqref{quasiperiodicWave}
satisfies
the following $d+1$ PDEs:
\begin{align}
\omega\cdot\pa_{\vphi}U&-X_{H}(U)=0\,,
\label{equazioni embedding1}
\\
\mathtt{V}^{T}\cdot\pa_{\vphi}U&+\nabla U=0\,,
\label{equazioni embedding2}
\end{align}
where\footnote{equivalently $\tV^{i}\cdot\pa_{\vphi}U+\pa_{x_i}U=0$ for $i=1,\ldots,d$.}
$X_H$ corresponds to the r.h.s of the system \eqref{eq:113} and where
\[
\mathtt{V}^{T}\cdot\pa_{\vphi}U=\sum_{j=1}^{\nu}\mathtt{v}_{j} \pa_{\vphi_j}U\,.
\]
The $d$-equations in \eqref{equazioni embedding2}
imply that $U$ can be written as in \eqref{travellingWaves},
and are equivalent to require that $U$ is invariant for the flow
of the momentum Hamiltonians in \eqref{Hammomento}.

We shall construct 
such solutions localized  
 in Fourier space at 
$\nu$ distinct  \textit{tangential sites}. We first need to make a 
crucial assumption on the dual 
lattice $\Gamma^{*}$ in \eqref{def:lattice2} generated by the vectors
${\bf K}_1, \ldots, {\bf K}_d$.

\begin{hypo}{\bf (Non-resonant lattice).}\label{hyp:lattice}
Define the vector 
\[
\vec{P}:=\big({\bf K}_j\cdot{\bf K}_k\big)_{1\leq j\leq k\leq d}\in \R^{d(d+1)/2}\,,
\]
where $\{ {\bf K}_1, \ldots, {\bf K}_d \}$ is the dual basis of the dual lattice $\Gamma^*$ and $\cdot$ denotes the standard scalar product in $\R^{d}$.
We assume that $\vec{P}$ has rationally independent components, i.e. $\vec{P}\cdot\ell\neq0$ 
for all $\ell\in \Z^{d(d+1)/2}\setminus\{0\}$.
\end{hypo}
Under the assumption \ref{hyp:lattice} we consider the set
of tangential sites 
\begin{equation}\label{TangentialSitesWW}
S:=\{\overline{\jmath}_1, \dots, \overline{\jmath}_\nu \}\subset \Gamma^{*}\setminus\{0\}\,,
\qquad
0< |\overline{\jmath}_1| <  |\overline{\jmath}_2| < \dots < | \overline{\jmath}_\nu|\,,\;\;\;\nu\geq1\,.
\end{equation}
Hypothesis \ref{hyp:lattice} %and condition \eqref{TangentialSitesWW} 
provides the precise meaning of 
the expression 
``generic choice'' in the informal statement of the main Theorem \ref{thm:mainmain}.
The solutions of \eqref{eq:113linear} that originate by exciting the tangential modes 
have the form
(recall \eqref{solLinearicomplex})
\begin{equation}\label{solLineari}
\begin{aligned}
\eta(t,x)&=\sum_{j\in S}
\sqrt{2\zeta_{j}}\Omega^{\frac{1}{2}}(j)\Big(\cos(\Omega(j)t)\cos(j \cdot x)+
\sin(\Omega(j)t)\sin(j \cdot x)\Big)\,,\\
\psi(t,x)&=\sum_{j\in S}
\sqrt{2\zeta_{j}}
\Omega^{-\frac{1}{2}}(j)\Big(
\cos(\Omega(j)t)\sin(j \cdot x)-
\sin(\Omega(j)t)\cos(j \cdot x)\Big)
\end{aligned}
\end{equation}
where $\zeta=(\zeta_j)_{j\in S}$ with $\zeta_j>0$. They
are superposition of periodic traveling linear waves with wave vector  $\overline{\jmath}_i\in \Gamma^{*}$ 
and frequency $\Omega(\overline{\jmath}_{i})=\sqrt{|\overline{\jmath}_i|\tanh(\th|\overline{\jmath}_i|)}$,
for $i=1,\ldots,\nu$. 
Such motions are quasi-periodic (or periodic) traveling waves of the form \eqref{travellingWaves} 
with (recall \eqref{wild})
\begin{align}
&\bullet \; {\rm {\bf frequency\, vector: }}\quad
\overline{\omega}:=\overline{\omega}(\th):=
\left(
\omega_{\overline{\jmath}_{i}}(\th)
\right)_{i=1,\ldots,\nu}\in\mathbb{R}^{\nu}\,,
\quad \omega_{\overline{\jmath}_{i}}(\th)
:=\Omega(\overline{\jmath}_i)\,,\label{LinearFreqWW}
\\
&\bullet \; {\rm {\bf wave\, vectors: }}\qquad
\tv_{i}:=\overline{\jmath}_{i}\in \Gamma^{*}\,,\quad i=1,\ldots,\nu\,,\label{velocityvec1}
\\
&\bullet \; {\rm {\bf velocity\, vectors: }}\qquad 
\mathtt{V}^{j}=(\tv_{i,j})_{i=1,\ldots,\nu}=(\overline{\jmath}_{i,j})_{i=1,\ldots,\nu}\,,
\quad j=1,\ldots,d\,.\label{velocityvec2} 
%\text{\color{red} FORSE QUESTA SI PUO TOGLIERE MI SA CHE NON SERVE A NULLA}
\end{align}

We  shall construct traveling quasi-periodic  
solutions of the nonlinear equations  \eqref{eq:113} 
with frequency vector $ \omega $ 
belonging to an open bounded subset 
$ {\mathtt \Omega} $ in $ \R^\nu $
which is diophantine according to the following definition.

\begin{defn}{\bf (Diophantine vectors).}
Fix 
$ \gamma \in (0,1) $, $ \tau >  \nu - 1 $ and fix an open bounded set 
$\mathtt{\Omega}\subset\R^{\nu}$.
We define the set of $(\gamma,\tau)-$\emph{diophantine vectors} as
\begin{equation}\label{def:DCgt}
\tD\tC (\gamma, \tau) := 
\Big\{ \omega\in {\mathtt \Omega} \subset\R^\nu \ : \ 
\abs{\omega\cdot \ell}\geq \gamma \jap{\ell}^{-\tau} \ , 
\ \forall\,\ell\in\Z^\nu \setminus \{0\} \Big\}\,, 
\end{equation}
where  $\langle \ell \rangle :=  \max\{1, |\ell|\} $.
\end{defn}
Moreover our solutions will be ``close'' to linear ones 
(see \eqref{solLineari}) supported on $S$ in the Sobolev spaces
of functions
$H^{s}(\T^{\nu}\times \T_{\Gamma}^{d};\R^{2})$.
The main result of the paper is the following.

\begin{thm}{\bf (Quasi-periodic 3D Stokes Waves in finite depth).} \label{thm:main0}
Consider finitely many tangential sites $ S \subset \Gamma^{*}\setminus\{0\} $
as in \eqref{TangentialSitesWW} and recall \eqref{mathtt h1 h2 fissati}. 
Then there exists $ \bar \sigma >  0 $ such that for any 
$\mathfrak{s} \geq \bar \sigma$,  there exists 
$ \varepsilon_0 \equiv \varepsilon_0(\mathfrak{s}) \in (0,1) $ such that,  
for any $ |\zeta |   \leq \varepsilon_0^2  $, 
$  \zeta := (\zeta_{ \overline{\jmath}_{i}} )_{i = 1, \ldots, \nu} \in \R_+^\nu $,
there exists a Cantor-like set  $ {\cG}_{\zeta} \subset [\th_1, \th_2] $ 
with asymptotically full measure as $ \zeta \to 0 $, i.e. 
$ \lim_{\zeta \to 0} | {\cG}_\zeta|  = {\th}_2- {\th}_1 $ such that
the following holds.
For any  $ \th \in {\cG}_\zeta $, 
the pure gravity water waves equations 
\eqref{eq:113} have a 
quasi-periodic traveling wave solution 
(according to Definition \ref{def:quasitravelling})
of the form
\begin{equation}\label{SoluzioneEsplicitaWW}
\begin{aligned}
\eta(t,x)&=
\sum_{j\in S}
\sqrt{2\zeta_{j}}\Omega^{\frac{1}{2}}(j)\Big(\cos(\omega_{j}t)\cos(j\cdot x)+
\sin(\omega_{j}t)\sin(j\cdot x)\Big)
+r_{1}(\omega t,x)\,,
\\
\psi(t,x)&=
\sum_{j\in S}
\sqrt{2\zeta_{j}}
\Omega^{-\frac{1}{2}}(j)\Big(
\cos(\omega_{j}t)\sin(j\cdot x)-
\sin(\omega_{j}t)\cos(j\cdot x)\Big)+
r_2(\omega t,x)\,,
\end{aligned}
\end{equation}
where the frequency $\omega$ and the functions $r_{i}$, $i=1,2$, satisfy 
the following properties:

\noindent
$\bullet$
the frequency vector $\omega$ belongs to $\tD\tC (\gamma, \tau)$
for some $\tau\gg\nu$ and some $\gamma=\gamma(\zeta) \to 0$ as $|\zeta|\to0$,
and are ``close'' to the linear frequency $\overline{\omega}(\th)$, namely 
\[
 \omega=\omega(\th,\zeta) \to \overline{\omega}(\th) \;\;{\rm as}\;\;
|\zeta|\to 0\,;
\]

\noindent
$\bullet$
the functions  $r_{i}$, $i=1,2$, are quasi-periodic traveling 
waves with frequency vector  $\omega$
and wave vectors $\overline{\jmath}_{i}$, $i=1,\ldots,\nu$,
and are $o(\sqrt{|\zeta|})$ in the space 
$H^{\mathfrak{s}}(\mathbb{T}^{\nu}\times\mathbb{T}_{\Gamma}^{d};\R)$.

\noindent
In addition, these quasi-periodic solutions are linearly stable.
\end{thm}

Some comments are in order.
 
 \smallskip
 $\bullet$
To the best of our knowledge, Theorem \ref{thm:main0} is the first existence result of global solutions 
for \eqref{eq:113} under \emph{periodic}
 boundary conditions (that do not reduce to steady solutions) which provides the existence of families of initial data exhibiting 
 a time-\emph{quasi-periodic} behavior.
 Some global well-posedness results are available when
\eqref{eq:113} is posed on the \emph{Euclidean space} $\R^d$.
 We mention for instance, Wu \cite{Wu97}  Germain-Masmoudi-Shatah  \cite{GMS} and Wu \cite{Wu2} 
 for the three dimensional case
 and Ionescu-Pusateri \cite{IP}, Alazard-Delort \cite{ADel2}, 
 Ifrim-Tataru \cite{IFRT} for the (more difficult since there is less dispersive effect) two dimensional case.
 In the papers mentioned above the dispersive nature of the equation play a fundamental role.
 On the other hand, when the horizontal variable is considered on tori, there are no dispersive effects 
 which allow to control the solutions for long times.
In this case, we construct our solutions by taking into account possible resonant effects among the linear frequencies of oscillations. 
This paradigm has been already implemented in the papers mentioned before
\cite{CN, IP-Mem-2009, IP2},
and, in parallel, in 
studying long time existence for water waves,
by means of a Birkhoff normal form approach.
%mostly in dimension one. 
We refer, for instance to
Hunter-Ifrim-Tataru
\cite{HIT1}, Berti-Feola-Franzoi \cite{BFF}, for cubic life span ($T\sim \e^{-2}$ for data of size $\e\ll1$), 
Berti-Feola-Pusateri \cite{BFP} for time of stability $\sim \e^{-3}$ in the completely resonant case of pure
gravity waves in infinite depth. For almost global existence results (time of stability $\sim \e^{- N}$ for any $N \in \N$), we mention Berti-Delort \cite{BD}, Berti-Maspero-Murgante \cite{BeMaMur} 
which deeply exploit the superlinear dispersion relation of the equation in presence of surface tension.
Concerning the much more difficult case of the 3D model we mention the remarkable paper
by Ionescu-Pusateri \cite{IP6} proving the lifespan of existence time of order $\e^{-5/3+}$
for the gravity-capillary waves in three dimensions.

The result above regards the evolution of \emph{open} sets of initial data close to the origin for 
\emph{finite} time. Of course to obtain a control over an infinite time interval is more difficult.
Besides this, we must deal with the combined problem of a \emph{weak dispersion} of the pure gravity case
and the presence of \emph{many non-trivial} resonances
due to the high space dimensional setting.

\smallskip
$\bullet$
It is well known that 
the water waves vector field $X_{H}$ in  \eqref{eq:113} possesses,  in addition to
 the Hamiltonian structure and the $x$-translation invariance, some further symmetries. Indeed
 one has that 
\begin{itemize}
%\vspace{0.5em}
%\noindent
\item[(i)] $X_{H}$ is  \emph{reversible} with respect to  the involution
 \[\label{involutionReal}
 S: (\eta,\psi)(x)\mapsto (\eta,-\psi)(-x)\,,
\]
i.e. it satisfies $X_{H}\circ S=-S\circ X_{H}$;
%\begin{equation}\label{involution2Real}
%X_{H}\circ S=-S\circ X_{H}\,;
%\end{equation}
\item[(ii)]
$X_{H }$ is  \emph{even-to-even}, 
i.e. maps $x$-even functions into $x$-even functions. 
%In other words it preserves the subspace 
%\begin{equation}\label{evenFuncReal}
%\Big\{\eta(x)=\eta(-x) \,, \;\; \psi(x)=\psi(-x)\Big\}\,.
%\end{equation}
\end{itemize}
Such properties have been extensively exploited in several papers such as
\cite{BBHM}, \cite{BM20}, \cite{IPT}, \cite{BFM}, \cite{BFM2} 
in order to remove some degeneracies by working on subspaces of functions 
satisfying some parity conditions in time and space.
In this paper we are interested in the existence of \emph{traveling} waves which hence
 do not satisfy the eveness property, indeed even solutions are \emph{standing} waves. 
Another  feature of our result is that we do not exploit the reversible structure but we construct traveling waves
by only using  the 
 the Hamiltonian structure and the $x$-translation invariance of \eqref{eq:113},
 as done, for 2D Water Waves,  in \cite{FGtrave}.
We then construct solutions which are not necessarily restricted to the subspace of 
 \emph{reversible} solutions. Of course, if we exploit the reversible structure in our proof, we get reversible solutions. 

\smallskip
$\bullet$ The dispersion relation \eqref{dispersionLaw} is \emph{sub-linear}. As a consequence,
the differences between eigenvalues of the linearized operator at zero fulfill very ``bad''
separation properties which is a very strong obstruction in performing the spectral analysis. 
More precisely, in our procedure,  we must deal with  ``losses of derivatives'' (both w.r.t. time Fourier indices and space Fourier indices) coming from the small 
divisors (see for instance the ``weak'' lower bounds in \eqref{Cantor set}).
This problem of course is much more difficult in the high dimensional case.
% and it represents one of the main obstruction in performing the spectral analysis of linearized equations at approximate quasi-periodic solutions.
 %\cite{IP6} and in the one dimensional cases \cite{BBHM, BFM2} among the others.

However, in one space dimension there is a quite well understood strategy to deal with this problem, which is based on the fact that one has a good asymptotic expansion of the eigenvalues $\mu(j)$, $j \in \Z$ in decreasing powers of the Fourier index $j \in \Z$.
Unfortunately, this does not hold in higher space dimensions and therefore, the control of the resonances is much more complicated. 
%The absence of such an asymptotic in high space dimension is well known in the literature. In particular for Schr\"odinger type operators with unbounded perturbations, this problem has been studied in \cite{BLM2, BLM3}, where some "directional expansion" of the eigenvalues has been proved and it has been applied in \cite{BLM4} (see also the extension \cite{BL5}) to study the growth of Sobolev norms for time dependent linear Schr\"odinger equations}.
%\red{Aggiungere qui le espansioni di langella?}
%On the contrary, in our case high dimensional case, we have also to face an additional problem,
In our case, we basically have to deal with operators 
of the form (recall \eqref{wild})
\begin{equation}\label{verdeoliva}
\omega\cdot\pa_{\vphi}+\e V(\vphi,x)\cdot\nabla+\ii \Omega(D)+{\rm l.o.t}
\end{equation}
where $\vphi\in \T^{\nu}, x\in \T^{d}_{\Gamma}$, $V$ is a smooth vector valued function, $\e\ll1$ a small parameter
and where by $\emph{l.o.t}$ we denoted lower order terms which are 
pseudo-differential operators of order $1/2$ and of size $\e$.
Hence we developed a flexible approach
to obtain 
normal form reduction and good properties of the eigenvalues  of perturbations 
of dispersive operators of order less than one in high dimension.
We remark that operators of such type arise in many equations of fluid mechanics.
Note that, for the moment, our approach does not cover the case of gravity capillary water waves
(super-linear dispersion relation)
for which completely new ideas are required.

We will discuss more in details this step later in this introduction. For the moment, we underline 
the key ingredients of the approach:
$(i)$
the normal form reduction  we perform is based on the fact that the operator $\omega\cdot\pa_{\vphi}$
(of order one in \emph{time}) is dominant with respect to the 
dispersive operator (in space) $\Omega$ which has order $1/2$;
$(ii)$ some cancellation of exact resonances coming from the conservation of momentum,
that  provide a constraint between that space Fourier indexes and the time ones.
%reduces the degree of freedom of the problem

%\medskip
%This is a very strong issue 
%even in classical mechanics. Indeed developing perturbation theory of an Hamiltonian 
%of the form $\zeta\cdot\x+\e P(x,\x)$ with $\zeta$  a diophantine vector  
%is much more easier  than considering an Hamiltonian $H(x,\x)=H_0(\x)+\e P(x,\x)$
%where $\nabla_\x H_0(\x)$ is not a constant function. In such a case 
%there would be resonances depending on the point of the phase space.

\smallskip
$\bullet$ 
Recalling \eqref{dispersionLaw},
we have that the linear frequencies of oscillation admit the expansion
\[
\sqrt{|j|\tanh(\th |j|)}=\sqrt{|j|}+r(j,\mathtt{h}) 
\qquad {\rm where}
\quad 
|\pa_{\mathtt{h}}^{k}r(j,\mathtt{h})|\leq C_{k}e^{-\mathtt{h}|j|}
\]
for any $k\in \N$, $j\in \Gamma^{*}$ with $|j|\geq1$. In the bound above 
%is  uniform in $\mathtt{h}\in [\mathtt{h}_1,\mathtt{h_2}]$ and 
the constant $C_{k}$ 
depends only on $k$ and $\mathtt{h}_1$ (recall that $\mathtt h \in [\mathtt h_1, \mathtt h_2]$, $0 < \mathtt h_1 < \mathtt h_2$). We use the depth  parameter $\mathtt{h}$
to modulate the frequencies of oscillation in order to extend the \emph{degerate KAM theory}
approach
used in the one dimensional case \cite{BBHM}.
 While
in the gravity-capillary case the capillarity parameter $\kappa$
moves the linear frequencies of polynomial quantities in the index $|j|$ (see for instance \cite{BM20}),
in the pure gravity case we have that  frequencies vary in $\mathtt{h}$
only by an exponentially small quantity. 
We prove the \emph{transversality properties} (see Proposition \ref{lem:transversality})
required by the degenerate KAM by deeply exploiting 
the momentum conservation together with the genericity of the lattice, see Hypothesis \ref{hyp:lattice}.
These arguments are developed in Section \ref{sec:degener}.
While the general strategy is the one in \cite{BBHM} there are serious differences in the proof coming from the high dimensional context.
%
%In addition to this, due the our high dimensional context, we cannot
%adapt the arguments in \cite{BBHM} which are strongly base on the dimension one, i.e. $j\in \Z$.

\smallskip
$\bullet$ Our solutions are arbitrarily regular, namely the traveling wave profile is in 
$H^{s}$ for any $s$ large enough, both in time and space. 
We believe that the presence of so many regular and \emph{global} 
solutions could help in the next future 
 to 
tackle the difficult problem of studying 
%One of the possible next frontiers 
%is to use such quasi-periodic solutions to study 
long time dynamics for open sets of initial data, at least close to the origin.

\smallskip
$\bullet$
As a byproduct of our approach, based on a reducibility argument 
of the linearized equation at a quasi-periodic function, we obtain 
the linear stability of the solutions $U(\omega t,x)$ (as done in \cite{BBHM}, \cite{FGtrave} for instance).
More precisely,  we prove that there exists a set of coordinates in which the linearized problem 
at the the solution $U(\omega t,x)$ is diagonal, with purely imaginary spectrum, in the directions normal to the embedded torus.
This implies that 
for the linearized problem, initial data starting close to the torus remains close to it for all times.
This point will be discussed in details in section \ref{sec:linearstability}.
We remark that to obtain such reducibility results is not trivial at all in high space dimension.
There is a quite extended literature a many efforts have been done in extending 
KAM theory for PDEs 
in high space dimension. The crucial difficulty, that arise also in performing just a reducibility scheme for linearized operators,  comes from the presence of many possible resonances. Therefore second Melnikov conditions
(lower bounds
on differences of eigenvalues required by the scheme, see \eqref{secondeINTRO}),
are typically violated  or produces a loss of regularity both in time and space.
To overcome this problem much more difficult schemes are required, for instance taking into account 
only ``block-diagonal'' normal forms which are related to multiplicity of eigenvalues.
We refer for instance to \cite{gengyou1, EGK2016, EK2010, PP2015} and reference therein for KAM theory 
on $\T^{d}$.
Due to the difficult structure of the resonances, the existence of quasi-periodic solutions in high dimension has been addressed also with different  techniques without requiring second Melnikov conditions.
In this case the invertibility of the linearized operators is obtained by means of ``multiscale analysis'' 
and do not rely on reducibility schemes. As a drawback, no information about the linear stability of the solution is 
obtained. For a detailed discussion on this point, which is quite far from the aim of the present introduction,
 we refer to the monograph \cite{BertiBollewaveBook}.
All the mentioned results deal with \emph{semilinear} perturbations. Reducibility results 
in the case of unbounded perturbations of diagonal operators are available only 
for specific models/domains usually with a special and simple geometry of the resonances, 
like perturbations of transport-like operators or operators on Zoll manifolds.
Without trying to be exhaustive we mention   \cite{FGMP19, Mon1, BLM2013, BGMR1, FGsphere}. 
For other higher dimensional examples in fluids, that still regards  reducibility for transport-like operators, 
we refer the reader to the next paragraph in this introduction. It is known that for the reducibility scheme, it is fundamental to have a sharp asymptotic expansion of the eigenvalues, which is available for perturbations of transport like operators whereas in general is known to be false. In particular, for Schr\"odinger type operators with unbounded perturbations, this problem has been studied in \cite{BLM2, BLM3}, where some "directional expansion" of the eigenvalues has been proved and it has been applied in \cite{BLM4} (see also the extension \cite{BL5}) to study the growth of Sobolev norms for time dependent linear Schr\"odinger equations. However, it seems that  such "directional expansions"  are not enough to obtain reducibility in the case of unbounded perturbations.
In our case, we must deal with linear operators like in \eqref{verdeoliva}, which are deeply different from Schr\"odinger type operators (since the dispersion is sublinear). 
In particular, we treat  perturbations of order $1$ of a dispersive operator $\Omega(D)$ of order $1/2$ whose symbol is then note a pure transport ($\nabla_\xi \Omega(\xi)$ is non-constant w.r.t. $\xi$).
We obtain the reducibility result, and hence the linear stability of the quasi-periodic solution,
by implementing a normal form
strategy deeply based on the sub-linear nature of the dispersion and the conservation of momentum,
which allows us to require very weak second Melnikov conditions (see \eqref{secondeINTRO}).
This strategy will be discussed more in detail in the scheme of the proof below.

\smallskip
$\bullet$  	A crucial ingredient of our proof of Theorem \ref{thm:main0} is to have precise and quantitative
information on the  pseudo-differential structure on the linearized equation of \eqref{eq:113} at any 
approximate quasi-periodic solutions.
This requires a refined "algebraic" and "quantitative" analysis of the Dirichlet-Neumann operator appearing in \eqref{eq:112a}.
Of course  it is well-known that the Dirichlet-Neumann is a pseudo-differential operator of order one,
but however the information present in literature are not enough for our aims.
In particular, 
we need a very  precise control of the dependence on time of the symbols and of smoothing remainders
appearing in the pseudo-differential expansion of the Dirichlet-Neumann operator.
An important technical 
part of the proof it has been  to 
provide sharp \emph{space-time} tame estimates for the elliptic problem defining $G(\eta)\psi$
and for the symbols appearing in the expansion.
For more details on this point we refer to sections \ref{sec:generalStrategy} and 
\ref{sec:DNsezione}-\ref{sec:pseudoDN}. We remark that in one dimension, the structure of the Dirichlet-Neumann operator is very simple, since it is a smoothing perturbation of the Fourier multiplier $|D| \tanh(\mathtt h |D|)$ (see \cite{BBHM}  for details)

\vspace{0.5em}

\noindent\textbf{Time quasi-periodic waves in fluids.}
We now briefly discuss some related results about the 
the existence of \emph{special} periodic or quasi-periodic waves in Fluid Mechanics, which is  the research line in which this paper enters.

\vspace{0.3em}
\noindent
{\emph{2D Water Waves.}} As already remarked in this introduction, the first early results about
\emph{traveling} waves date back to $19$-th century. We mention the pioneering work by Stokes 
\cite{stokes} on nonlinear approximation of the flow and 
 Levi-Civita \cite{levi}, and 
Struik \cite{struik} where some rigorous
bifurcation of small amplitude bi-dimensional traveling waves are proved.
These result are essentially an application of infinite dimensional standard bifurcation theory 
(Crandall-Rabinowitz) since 
it is not a small divisor problem.
In the recent years there have been lots of generalization of such results, 
for instance to water waves with vorticity  or global bifurcation results, see  for instance 
 \cite{CS,CIP, Wh0, Wh}.

The generalization of such results to the case of standing waves, which presents small divisors, 
are quite more recent. 
See Plotnikov-Toland \cite{PlTo}, Iooss-Plotnikov-Toland \cite{IPT},
Alazard-Baldi \cite{AB15}.
The quoted papers regard \emph{time periodic} waves. 
The quasi-periodic case for the water waves system has been treated
by Berti-Montalto \cite{BM20} for capillary standing waves, 
Baldi-Berti-Haus-Montalto \cite{BBHM} for gravity standing waves in finite depth.
For traveling quasi-periodic waves with constant vorticity see 
Berti-Franzoi-Maspero
 \cite{BFM, BFM2} and 
 Feola-Giuliani \cite{FGtrave} on pure gravity water waves in infinite depth.
% Regarding the model with vorticity we also mention
% \cite{CS,CIP, Wh0, Wh}
% \red{mettere le cose di whalen}

% We also mention the remarkable 
%results on \emph{completely resonant} models by Baldi-Berti-Montalto \cite{KdVAut} on Kortweg-deVries,
%Feola-Giuliani-Procesi \cite{FGP20} on Degasperis-Procesi equation, 
%Feola-Giuliani \cite{FGtrave} on pure gravity water waves in infinite depth, 
%and G\'omez-Serrano--Ionescu--Park
%\cite{GomIonPark} on the generalized SQG equation.

\vspace{0.3em}
\noindent
{\emph{3D Water Waves.}}
For the three dimensional water waves system there are only few results, restricted to the case of 
waves which are stationary in a moving frame.
%\emph{periodic} time behaviour.
%The three-dimensional case has been successfully approached more recently.
We quote the paper  \cite{CN} (and reference therein)
by  Craig-Nicholls where the
 existence of
traveling wave solutions
 is proved
in the gravity-capillary case with space periodic boundary conditions.
We also refer to the recent paper \cite{GNPW1} for the case with non-constant vorticity.
In absence of capillarity the existence of traveling waves is 
a small divisor problem. 
This case has been treated by 
Iooss and Plotnikov in \cite{IP-Mem-2009, IP2}.
We also mention \cite{Wh3,wh4}.

\vspace{0.3em}
\noindent
{\emph{Other fluid models.}}
In the very recent years, a particular attention 
has been dedicated to the 
construction of periodic or quasi-periodic nonlinear waves in Euler and Navier-Stokes equations.
An important topic is to study the evolution of vortex patches
which are solutions of the 2d Euler equation (or $\alpha$SQG) with initial data whose vorticity is 
the characteristic function of a domain with smooth boundary.
In this direction we mention
Berti-Hassainia-Masmoudi \cite{BhasMas1}, 
Hassainia-Hmidi-Masmoudi \cite{HHM2},
Roulley \cite{roulley1}, 
Hmidi-Roulley 
\cite{hmidiroulley}, 
Hassainia-Hmidi-Roulley \cite{HHR1}, 
G\'omez-Serrano-Ionescu-Park
\cite{GomIonPark}.
We also mention the first rigorous justification of the classical leapfrogging motion in 
\cite{HHM} (see also  \cite{HHR}).
 
 All the results mentioned above reduce to the study of a one dimensional PDE,
 whereas Euler equation with periodic boundary conditions  ($x \in \T^d$, $d =2,3$) or in the channel $\R \times [-1,1]$, is a PDE in dimension 2 or 3.
% For the Euler equation with periodic boundary condition 
 In the case in which $x \in \T^d$, the construction of nonlinear quasi-periodic waves has been 
 obtained in presence of a forcing term which is quasi-periodic in time.
 We mention
Baldi-Montalto \cite{BM21},
Montalto \cite{monNavier}, Franzoi-Montalto \cite{framonNavier} on Navier-Stokes.
These results have been extended in case without forcing by Enciso-Peralta-Torres de Lizaur
\cite{EncPeraltaTorres}. See also the extension \cite{FrMontalmost}.
For examples of large amplitude quasi-periodic waves in fluids, we refer 
to \cite{ BiFMT, CMT24}. Concerning quasi-periodic traveling waves bifurcating from shear flows, for the 2D Euler equation in the channel, we also mention Franzoi-Masmoudi-Montalto \cite{FMaMo}. 

\vspace{0.9em}
\noindent
{\bf Notations.}
Given some parameters $m_1,\ldots,m_{n}$, 
the notation $a\lesssim_{m_1,\ldots,m_{n}} b$
means that there exists a constant $C(m_1,\ldots,m_{n})>0$, depending on such parameters, such that
$a\leq C(m_1,\ldots,m_{n}) b$. Along the paper these parameters $m_1,\ldots,m_{n}$
will be, for instance,  the Sobolev regularity index $s$, the order of pseudo-differential operators, 
the number of derivatives in $\vphi\in \T^{\nu}, x\in \T^d_{\Gamma},\x\in \R^{d}$.
Sometimes, we shall omit to write the dependence $\lesssim_{s_0,k_0}$ on the
low Sobolev regularity $s_0$ (see \eqref{Sone})  and the index of non-degeneracy $k_0$ in \eqref{k0_def}.
In general we shall omit also the dependence on the dimensions $\nu,d$ and on the Diophantine 
constant $\tau>0$ appearing in \eqref{def:DCgt}, which are considered fixed constants.
Along the paper we denote by $\mathbb{N}$ the set of \emph{natural} numbers included $n=0$.
Given $z\in \R^{p}$, $\beta\in \N^{p}$, $p\in\N$ we adopt the standard multi-index notation
to denote
$\pa_{z}^{p}=\pa_{z_1}^{\beta_1}\cdots\pa_{z_{p}}^{\beta_p}$.
%Sometimes, along the paper, we omit to write the dependence s0,k0 with respect to s0, k0, because s0 (defined in (1.21)) and k0 (determined in Sect. 3) are considered as fixed constants. Similarly, the set S+ of tangential sites is considered as fixed along the paper.

\vspace{0.5em}
\subsection{Scheme of the proof}\label{sec:generalStrategy}
  Here we 
 resume the steps of the proof highlighting the difficulties and the key ingredients.

%\smallskip
%\noindent
%$(i)$ {\bf Nash-Moser theorem of hypothetical conjugation.} 
The overall approach for dealing with small divisors problems in quasi-linear PDEs is the Nash-Moser scheme in Berti-Bolle \cite{BB}, 
we refer to sections \ref{sec:nash}-\ref{sec:approx_inv}.
All the work and all the crucial novelties are  concentrated in the analysis of the linearized problem
(see \eqref{Lomegatrue})
 and of the resonances in order to prove the invertibility property \eqref{almi4}.
 Such an estimate is used to prove the convergence of the scheme 
 in section \ref{sec:NaM}.
 In the analysis of the linearized operator we deal with three main issues:
% one one has to deal with the main features and difficulties of the problem, namely
\begin{itemize}
\item{ the very complicated pseudo-differential structure of the Dirichlet-Neumann operator};
\item{ the presence of derivatives in the nonlinearity (the PDE is quasi-linear) which is particularly difficult to manage in higher space dimension};
\item{ the presence of very strong resonance phenomena (due to the fact that the PDE is in more than one space dimension) and the "bad asymptotic expansion" of the eigenvalues of the linearized equations}.
\end{itemize}

%In this section we will concentrate on the main ideas we used to deal with these three issues.

\smallskip
\vspace{0.9em}
\noindent
$(i)$ {\bf Dirichlet-Neumann operator.}
In order to understand precisely the pseudo-differential structure of the linearized operator 
(at a quasi-periodic function)
we need preliminary results on the Dirichlet-Neumann operator in \eqref{eq:112a}.
More precisely, we shall consider a function $\T^{\nu}\times \T^{d}_{\Gamma}\ni(\vphi,x)\mapsto \eta(\vphi,x)$
smooth both in the time and space variables and we prove the following expansion:
\begin{equation}\label{claimoGeta}
G(\eta) =  |D| \tanh(\mathtt h |D|) +   \opw(a_{G}(\vphi,x,\x)) + {\mathcal R}_G\,,
\end{equation}
where $a_{G}$ is a symbol of order $1$ 
(see section \ref{sec:pseudodiffsimboli} for details on pseudo-differential operators), while
$\T^{\nu}\ni\vphi\mapsto
\mathcal{R}_{G}:=\mathcal{R}_{G}(\vphi)$ is a regularizing operator in space, i.e.
for any fixed $\vphi$ it maps $H^{s}(\T^{d}_{\Gamma})$ into $H^{s+\rho}(\T^{d}_{\Gamma})$
for any $\rho$, for $s>0$ large enough.
When ignoring the $\vphi$ dependence in the operator $G(\eta)$, such an expansion is 
essentially known in literature (see for instance \cite{ADel1,AlM}). 
Nevertheless, in our case, in order to implement
normal form scheme on the linearized operator we need additional sharp estimates
of symbols and remainders  in the space-time Sobolev space
$H^{s}(\T^{\nu}\times\T^{d}_{\Gamma})$. This requires to enter much more precisely in the analysis of the structure of the symbol $a_G$ and the remainder ${\mathcal R}_G$. 

More specifically, we prove in Theorem
\ref{lemma totale dirichlet neumann}, the following bounds:
for some $\s>0$, if $\|\eta\|_{H^{s_0+\s}}$ is small enough (for some $s_0\gtrsim \nu+d$) 
then, the symbol $a$ satisfies  for any $s\geq s_0$, $\alpha \in \N$
\begin{equation}\label{alloro1}
\|a_{G}\|_{1,s,\alpha}:=
\max_{0\le |\beta|\le \alpha} \sup_{\xi\in\R^d} 
\|\partial_\xi^\beta a_{G}(\cdot, \cdot, \xi) \|_{H^{s}(\T^{\nu}\times\T^{d}_{\Gamma})} 
\jap{\xi}^{-1+|\beta|}\lesssim_{s,\alpha}
\|\eta\|_{H^{s+\s}(\T^{\nu}\times\T^{d}_{\Gamma})}\,.
\end{equation}
The remainder ${\mathcal R}_G$ satisfies the following estimate.
For any $\beta_0\in \N$, and $M\geq \mathtt{c}\gg s_0+\beta_0$, one has
for any $|\beta|\leq \beta_0$, the estimate
\begin{equation}\label{stimaRRGG}
\begin{aligned}
\|\langle D \rangle^{M} \partial_\vphi^\beta {\mathcal R}_G \langle D \rangle^{-  \mathtt{c}}h
\|_{H^{s}(\T^{\nu}\times\T^{d}_{\Gamma})}&\lesssim_{s,M}
\|\eta\|_{H^{s_0+\s}(\T^{\nu}\times\T^{d}_{\Gamma})}\|h\|_{_{H^{s}(\T^{\nu}\times\T^{d}_{\Gamma})}}
\\&\qquad +
\|\eta\|_{H^{s+\s}(\T^{\nu}\times\T^{d}_{\Gamma})}\|h\|_{_{H^{s_0+\s}(\T^{\nu}\times\T^{d}_{\Gamma})}}\,,
\end{aligned}
\end{equation}
where $\langle D\rangle^{m}$ for  $m\in \R$ is the Fourier multiplier defined as (recall \eqref{def:lattice2})
\[
\langle D\rangle^{m}e^{\ii j\cdot x}=\langle j\rangle^{m}e^{\ii j\cdot x}\,,\qquad \langle j\rangle:=\sqrt{1+|j^2|}\,,
\qquad j\in \Gamma^{*}
\]
 A \emph{tame} estimates as above on the remainder $\mathcal{R}_{G}$ is necessary for implementing 
a normal form scheme on the linearized operator for two main reasons:
$(i)$ tame estimates are necessary to prove the \emph{convergence} of the KAM reducibility scheme
in finite regularity (on Sobolev spaces); $(ii)$ by taking $M\geq \mathtt{c}$ arbitrarily large,
estimate
\eqref{stimaRRGG} implies that $\mathcal{R}_{G}$ ``gains'' in space $M-\mathtt{c}$ derivatives.
This is crucial in order to control the loss of derivatives coming from the small divisors in the normal form scheme.
To prove such an estimate is one of the most difficult and new part of the proof and it is developed in sections
\ref{sec:DNsezione}-\ref{sec:pseudoDN}.
The major difficulties come from the fact we work in high dimension (in $d=1$ the pseudo-differential
expansion simplifies a lot, as we already mentioned), and the fact that we work in finite depth (which is much harder technically than the infinite depth case).

\medskip
We now briefly explain our approach and the specific estimates we obtained on symbols and remainders.

\medskip
\emph{(1) Green's function estimates.} First of all, in subsection \ref{sec:tameLaplace}, 
we provide a priori estimates on the solution of the 
free boundary elliptic problem \eqref{BoundaryPr}. The key simple idea is to straighten the
domain $\mathcal{D}_{\eta}$ in \eqref{Deta} through the diffeomorphism \eqref{change1}
and study the transformed elliptic problem \eqref{elliptic2} where $(x,y)\in \T^{d}_{\Gamma}\times[-1,0]$.
In this setting the simplification comes from the fact that the problem is set on a straight strip, 
while the difficulty is that one has to treat a non constant 
coefficients (depending on $x$ and $y$)
second order elliptic operator $\widetilde{\mathcal{L}}=\widetilde{\mathcal{L}}(\eta;y)$, 
see \eqref{op:tildeL}. 
Such coefficients are small as the embedding $\eta(\vphi,x)$ is small.
By exploiting such  smallness, we look for the solution $\phi(x,y)$ 
of the  problem \eqref{elliptic2} in the following way.
We set
\[
\phi(x,y)=\mathcal{L}_0[\psi]+u(x,y)\,,\qquad 
\mathcal{L}_0:=\frac{\cosh((1 +y)\th|D|)}{\cosh(\mathtt h |D|)}[\cdot]\,,\quad y\in[-1,0]
\]
where $\vphi_0:=\mathcal{L}_0[\psi]$  is the solution of the problem at the flat surface $\eta\equiv0$ 
(see \eqref{problema vphi 0 riscalato laplace}),
\[
\begin{aligned}
& \partial_{yy} \vphi_0 + \mathtt h^{2} \Delta \vphi_0 = 0\,,\quad {\rm in}\quad \T^{d}_{\Gamma}\times(-1,0)\,,
 \\
 &\vphi_0(x, 0) = \psi(x)\,,\qquad \partial_y \vphi_0(x, - 1) = 0\,,
\end{aligned}
\]
while $u(x,y)$ solves the forced  problem (see \eqref{prob omogeneop forzato})
\begin{equation}\label{squalo1}
\begin{aligned}
&\partial_{yy} u + \th^2 \Delta u + \widetilde{F}[u]= -\widetilde{F}\circ\mathcal{L}_0[\psi]\,,
\quad {\rm in}\quad \T^{d}_{\Gamma}\times(-1,0)\,,
\\&
u( 0, x) = 0 \,,\qquad \partial_y u(- 1, 0) = 0\,,
\end{aligned}
\end{equation}
where $\widetilde{F}:=\widetilde{\mathcal{L}}-\pa_{yy}-\mathtt{h}^2\Delta$, see \eqref{def op F vphi}.
The solution of \eqref{squalo1} is written as $u(x,y)=\mathcal{K}(\vphi,y)[\psi]$.
We provide space-time tame estimates on the Green function $\mathcal{K}(\vphi,y)$ 
as a linear operator from $H^{s}(\T^{d}_{\Gamma})$ into $H^{s}([-1,0]\times \T^{d}_{\Gamma})$
for any fixed $\vphi\in \T^{\nu}$ and 
taking into account regularity in $\vphi$.
This is done
in Proposition \ref{stima cal Kn equazione di laplace} of section \ref{sec:stimegreen}.

\medskip
\emph{(2) The non-homogeneous problem.}
By the Green's function estimates in Proposition \ref{stima cal Kn equazione di laplace}
we cannot deduce any pseudo-differential structure on the solution $\phi(x,y)$ of \eqref{elliptic2}.
In order to do this, we study the elliptic problem in a different setting introduced as follows.
Let us consider the diffeomorphism of the strip
$\mathcal{A}$ defined as (see \eqref{cambio di variabile striscia})
\begin{equation}\label{cambio di variabile strisciaINTRO}
\cA (\vphi) : u(x, y) \mapsto u(x, \alpha(\vphi, x) y), \quad \alpha(\vphi, x) 
:= \frac{\mathtt h}{\mathtt h + \eta(\vphi, x)}\,.
\end{equation}
We provides estimates on such  diffeomorphism  in section  \ref{sec:stimeA}
(see Lemmata 
\ref{lemma cambio di variabile striscia strano}-\ref{stime cal A vphi - Id ca,bio variabile strano striscia}).
Secondly, we introduce the function (see \eqref{func:tildevphi} and \eqref{strutturaPhiPhi})
\begin{equation}\label{elliptic55}
\widetilde{\Phi}(x,y):=\vphi_0(x,y)+\chi(y)(\mathcal{A}(\vphi)-\uno)\cL_0[\psi]+
\chi(y)\mathcal{A}(\vphi)\circ \mathcal{K}(\vphi,y)[\psi]\,,
y\in[-1,0]\,,
\end{equation}
where $\chi(y)$ is a cut off function (see \eqref{cut-off}) such that $\chi\equiv1$ for $y\sim0$.
Thanks to this definition, which is adapted to the finite depth case we are studying\footnote{we introduce the cut-off function $\chi(y)$ following the idea in \cite{BD} in such a way the problem for \eqref{elliptic55} 
is posed on the fixed segment 
$[-1,0]$. The simplest diffeomorphism in \eqref{change1bis} which flatten the free surface $z=\eta(\vphi,x)$ to $y=0$ would lead to the problem \eqref{elliptic2bis}
in which the second order elliptic operator $\mathcal{L}$ has coefficients independent of $y$, but the 
the bottom of the strip $z=-h$ is transformed to $y=-1-\th^{-1}\eta(\vphi,x)$.
This is why we study the problem for \eqref{elliptic55}. This issue does not arise in the infinite depth case.
}, 
we have the following properties:
\begin{enumerate}
\item[a)] the function $\widetilde{\Phi}$ solves the elliptic problem (see \eqref{elliptic5})
\begin{equation}\label{elliptictildeintro}
\left\{\begin{aligned}
{\mathcal{L}} \widetilde{\Phi}&=g\,,\qquad x\in\T^{d}_\Gamma\,,\;\; -1<y<0\,,\\
\widetilde{\Phi}(x,0)&=\psi(x)\,,\\
(\pa_{y}\widetilde{\Phi})(x,-1)&=0\,,
\end{aligned}\right.
\end{equation}
where (see \eqref{op:L})
\begin{equation}\label{op:Lintro}
\mathcal{L}:={\mathcal{L}}(\eta):= \partial_{yy} + \mathtt h^{2} \Delta + \beta_1(\vphi,x) \cdot \nabla \partial_y 
+ \beta_2(\vphi,x) \partial_y + \beta_3(\vphi,x) \Delta
\end{equation}
for some smooth functions $\beta_{i}$ depending only on $\eta(\vphi,x)$ 
(see their definition in \eqref{op laplace raddrizzata buona}),
and $g$ is a known forcing term depending implicitly on the whole solution $\phi(x,y)$ of \eqref{elliptic2}.
The advantage (w.r.t. to the problem \eqref{elliptic2})
is that we still have the to work on the straight strip $\T^{d}_{\Gamma}\times[-1,0]$
but we deal with a second order elliptic operator $\mathcal{L}$ whose coefficients are independent of the variable 
$y\in[-1,0]$.
Moreover, in order to control the resonances present in the problem, 
we need to control a large number  of derivatives in the parameter $\th$ 
which modulates the frequencies.
It turns out that to do this it is more convenient to work on the strip $[-1,0]$
and not on $[-\th,0]$.

\item[b)] The Dirichlet-Neumann operator $G(\eta)[\psi]$ in \eqref{eq:112a} can be expressed in terms of 
$\widetilde{\Phi}$ as
\begin{equation}\label{liberto1}
 G(\eta)\psi 
 = \th^{-1}(1+|\nabla\eta|^2)(\pa_{y}\widetilde{\Phi})(x,0)-\nabla\eta\cdot\nabla \psi(x) \,,
\end{equation}
so we only need to show that the derivative $\pa_{y}\widetilde{\Phi}$ admits a pseudo-differential expansion
when $y$ is close to $0$.
\end{enumerate}
The Green's function estimates for the problem \eqref{elliptictildeintro} are proved in Proposition 
\ref{lemma tame laplace partial y u langle D rangle u}. 
Now, by exploiting the special structure of the second order elliptic operator $\mathcal{L}$
we shall prove a  pseudo-differential expansion for $G(\eta)[\psi]$.
This is done in several steps.

\medskip
\emph{(3) Splitting of the second order operator.}
In the following we implement the overall strategy of \cite{AlM} of splitting the second order elliptic problem into two 
first order problems. Of course, in our context we need to take into account the dependence on the angles variables
$\vphi\in \T^{\nu}$. 

In section \ref{sec:decoupling} (see Lemma \ref{divo1}) we show that the operator 
$\mathcal{L}$ in \eqref{op:Lintro} can be written, up to regularizing remainders, as the composition of two 
first order operators. More precisely, we show that there exist
two symbols $a=a(\vphi,x,\x)$, $A=A(\vphi,x,\x)$ of order one 
such that

$\bullet$ $a(\vphi,x,\x)-|\x|$, $A(\vphi,x,\x)-|\x|$ satisfy estimates like \eqref{alloro1};

\smallskip
$\bullet$ one has (see Def. \ref{standard} for the quantization of symbols)
\[
\mathcal{L}=(\pa_{y}+\op(a))\circ (\pa_{y}-\op(A))
\]
up to regularizing remainders satisfying estimates like \eqref{stimaRRGG}.

As a first consequence of this splitting we have that the equation in \eqref{elliptictildeintro}
reads
\[
(\pa_{y}+\op(a))\circ (\pa_{y}-\op(A))\widetilde{\Phi}=\widetilde{g}\,,\qquad x\in \T^{d}_{\Gamma}\,,
\quad y\in[-1,0]\,,
\]
for some new non-homogeneous term $\widetilde{g}$. See \eqref{elliptic44} for more details.

As a second step one introduces the function 
\begin{equation*}%\label{funz:ww}
w:=\chi(y)\big(\pa_{y}-{\rm Op}(A)\big)\widetilde{\Phi}\,,
\end{equation*}
where $\chi$ is the cut-off function  given in \eqref{cut-off}, so that one shall write
\[
(\pa_{y}\widetilde{\Phi})(x,0)=\op(A)[\psi]+w(x,0)\,.
\]
This means that \eqref{liberto1}
 becomes (see Lemma \ref{lem:strutturaDNconw} for details)
 \[
 G(\eta)[\psi]=\th^{-1}(1+|\nabla\eta|^2))\op(A) [\psi] -\nabla\eta\cdot\nabla \psi(x)+
 \th^{-1}(1+|\nabla\eta|^2)) w_{| y = 0} \,.
 \]
 Formula above, together with explicit definition of the symbol 
$A$ in Lemma \ref{divo1}, would implies  the  pseudo-differential expansion 
\eqref{claimoGeta}
provided that the term $w(x,0)$ can be written as 
$w(x,0)=\mathcal{R}(\vphi)[\psi]$
for some   smoothing operator  $\mathcal{R}(\vphi)$ satisfying tame estimates like \eqref{stimaRRGG}.
This is the last and the most technical challenging part of the proof of Theorem 
\ref{lemma totale dirichlet neumann}, since one has to go deeply inside the structure of $w(x, 0)$ which was not done in previous papers. Moreover, in order to have some algebraic properties on the symbol $a_G$ and on the remainder ${\mathcal R}_G$, we exploit the self-adjointness of $G(\eta)$ and we prove that ${\rm Op}^W(a_G)$ and ${\mathcal R}_G$ are self-adjoint separately. This is the content of section \ref{self-adjoint-pezzi-DN}.

\medskip
\emph{(4) The parabolic problem and analysis of the remainders.}
In order to understand the properties of the function $w(x,y)$, we exploit the fact that it solves the 
\emph{parabolic problem}
\begin{equation}\label{Cauchy problem backward caloreINTRO}
\left\{\begin{aligned}
&(\pa_{y} + \op(a))w= F\,,
\\ 
&w(x,-1)=0\,,
\end{aligned}\right.
\end{equation}
where $F$
is an explicit  forcing term (see \eqref{forcingFF}) which essentially depends on $g$ in \eqref{elliptictildeintro} 
and the solution 
$\phi$ of \eqref{elliptic2}. Its precise structure is crucial for our aims. It is analyzed in 
Lemma \ref{lemma tame F Gamma vphi y che palle}. Moreover,
the solution of the equation \eqref{Cauchy problem backward caloreINTRO} has the form 
\begin{equation}\label{soluzione parametrica dirichlet NeumannIINTRO}
w(\cdot, y) = \int_{- 1}^y\cU(\vphi, y - z) F(\cdot, z)\, d z\,,\qquad \forall\,y \in [- 1, 0]\,,
\end{equation}
where the operator $\cU(\vphi, \tau)$, $\vphi \in \T^\nu$ 
is the flow of the (homogeneous) parabolic equation 
\begin{equation}\label{calore strisciaINTRO}
\begin{cases}
\partial_{\tau} \cU(\vphi, \tau) +  {\rm Op}\big(a(\vphi, x, \xi) \big) \cU(\vphi, \tau) = 0\,, 
\quad  0 \leq \tau \leq 1 \\
\cU(\vphi, 0) = {\rm Id} \,. 
\end{cases}
\end{equation}
The remainder $\mathcal{R}_{G}$ contains many terms 
which are discussed in details in subsections
\ref{sec:analisiW}-\ref{sec:analisiS}
which essentially arise in formula \eqref{soluzione parametrica dirichlet NeumannIINTRO}.
To understand them one need information on $F$, which are provided in Lemma 
\ref{lemma tame F Gamma vphi y che palle}
and on the parabolic flow $ \cU(\vphi, \tau)$ in \eqref{calore strisciaINTRO}.
In Lemma \ref{stima derivate lambda vphi flusso parabolico} we provide space-time a priori 
tame estimates on the flow $\cU(\vphi,\tau)$.
In subsections \ref{subsec:Poissonsymbols}-\ref{sec:parametrixEVO}
(see Proposition \ref{proposizione totale parametrica calore})
we prove a  pseudo-differential expansion of $\mathcal{U}(\vphi,\tau)$ 
with arbitrarily smoothing remainders by providing 
a quantitative version of the celebrated ``parametrix methods'' for elliptic equation of M. Taylor
(see \cite{Taylor}).

In the following we illustrate how to treat 
the most difficult remainder appearing in \eqref{soluzione parametrica dirichlet NeumannIINTRO} 
which has the following form:
by recalling \eqref{opS1}, \eqref{restiGammaSelle}, we need to study an operator of the form
\[
\cS_1(\vphi)  :=  \int_{- 1}^{\bar y} {\rm Op}\big(e^{y a(\vphi, x, \xi)} \big) \circ \Gamma(\vphi, y)\, d y\,, 
\qquad \bar{y}<0\,.
\]
Here $a(\vphi,x,\x)$ is a symbol of order one which is close (in the norm in \eqref{alloro1}) to the symbol $|\x|$
for $|\x|\geq 1/2$. 
It comes from a parametrix expansion of the parabolic flow of the problem \eqref{calore strisciaINTRO} 
studied in sections \ref{sec:decoupling}-\ref{sec:parametrixEVO}.

Moreover $\Gamma(\vphi,y)$ is an operator vanishing identically close to the boundary $y=0$, and which depends 
on the Green Function of the elliptic problem 
\eqref{elliptic4tris}. The problem \eqref{elliptic4tris} is obtained 
by \eqref{BoundaryPr} under some changes of variables constructed in order to 
reduce the elliptic problem  in the strip $\{(y,x) : x\in \T^{d}_{\Gamma}\,,\; -\mathtt{h}<y<\eta(x)\}$
into a problem posed on the fixed strip $\{(y,x) : x\in \T^{d}_{\Gamma}\,,\; -1<y<0\}$.
A part from the explicit form of $\Gamma$ the most dangerous terms come from considering
\begin{equation}\label{operatorSintro}
\cS(\vphi)=\int_{- 1}^{\bar y} {\rm Op}\big(e^{y a(\vphi, x, \xi)} \big) \circ\zeta(y)\mathcal{A}(\vphi)\circ\mathcal{K}(\vphi,y)dy
\end{equation}
where $\zeta(y)$ is a cut off function vanishing close to $y=0$,
$\mathcal{A}$ is the diffeomorphism of the strip defined 
in \eqref{cambio di variabile strisciaINTRO}
and $\mathcal{K}$ is the Green function of the problem \eqref{squalo1} in step $(1)$.

In view of the space-time tame estimates on $\mathcal{K}$
 and $\mathcal{A}$ given respectively in subsections  \ref{sec:stimegreen} and \ref{sec:stimeA},
% 
% 
%  
%First of all on $\mathcal{K}$ we have to perform the Green function space-time estimates
%as a linear operator from $H^{s}(\T^{d}_{\Gamma})$ into $H^{s}([-1,0]\times \T^{d}_{\Gamma})$
%for any fixed $\vphi\in \T^{\nu}$ and 
%taking into account regularity in $\vphi$.
%This is done
%in Proposition \ref{stima cal Kn equazione di laplace} of section \ref{sec:stimegreen}.
%We also provides estimates on the diffeomorphism $\mathcal{A}(\vphi)$ in section  \ref{sec:stimeA}
%(see Lemmata 
%\ref{lemma cambio di variabile striscia strano}-\ref{stime cal A vphi - Id ca,bio variabile strano striscia}).
%We conclude providing estimates on the composition
we are able to provide estimates on the composition
$\mathcal{A}(\vphi)\circ\mathcal{K}(\vphi,y)$ in Lemma \ref{corollario derivata cal A vphi cal K vphi}
(see \eqref{operatorRRmm}).
More precisely, we obtain that 
the operator
\begin{equation}\label{operatoreRintro}
  \partial_\vphi^\beta  \zeta(y) \circ \cA(\vphi) \circ \cK(\vphi, y)\circ
  \langle D \rangle^{- |\beta| }
\end{equation}
is bounded as linear operators from $H^{s}([-1,0]\times\T^{d}_{\Gamma})$ into itself.
Actually it is not only bounded but satisfies sharp space-time estimates, see in particular estimates \eqref{stimasuRRm}.
The important feature of such an estimates is that the dependence on the variables
$\vphi,x,y$ are mixed (due to the diffeomorphism $\mathcal{A}$). As a consequence,
one gets (see estimates \eqref{lagrima3})
\[
\|\pa_{\vphi}^{\beta} \zeta(y) \cA(\vphi) u\|_{H^{s}([-1,0]\times \T^{d}_{\Gamma})}
\lesssim \|u\|_{H^{s+|\beta|}([-1,0]\times \T^{d}_{\Gamma})}\,,
\]
and then
\[
\|\pa_{\vphi}^{\beta} \zeta(y) \cA(\vphi) \circ \mathcal{K}(\vphi,y)\psi\|_{H^{s}([-1,0]\times \T^{d}_{\Gamma})}
\lesssim 
\|\mathcal{K}(\vphi,y)\psi\|_{H^{s+|\beta|}([-1,0]\times \T^{d}_{\Gamma})}
\lesssim\|\psi\|_{H^{s+|\beta|}(\T^{d}_\Gamma)}\,.
\]
This is why, in order to obtain the boundedness of the operator in \eqref{operatoreRintro}
on the corrects spaces we need to add the weight $  \langle D \rangle^{- |\beta| }$
on the right hand side.
This is the main reasons why the estimates \eqref{stimaRRGG} on the reminders in the pseudo-differential expansion of $G(\eta)$ are \emph{unbalanced}. More precisely, 
we cannot have that the operator $\langle D\rangle^{n_1}\pa_{\vphi}^{\beta}\mathcal{R}_{G}
\langle D\rangle^{n_2} $
with $n_1+n_2=M-\mathtt{c}$ is bounded. We can only 
admit a loss on derivatives on the left with a gain, which balance the derivatives on $\vphi$, on the right.
We are allow to multiply by $\langle D\rangle^{M}$ on the left  in \eqref{operatorSintro}
due to the following facts.

$\bullet$ Since $|{\rm Re}(a)|\gtrsim|\x|$ and $y<0$ one can prove that $y^Me^{ya(\vphi,x,\x)}$
for any $M>0$
is a symbol of order $-M$. This is proved in Lemma \ref{prop exp simbolo}.

$\bullet$ Due to the presence of the cut-off $\zeta(y)$ (vanishing when $y\sim0$) the integral in 
\eqref{operatorSintro} is taken on the domain $[-1,\bar{y}]$ with $\bar{y}<0$ far from $0$.
Therefore we can apply the following simple estimate: for any $\alpha\geq0$ and 
$f \in L^2([- 1, 0], H^s(\T^{\nu}\times \T^{d}_{\Gamma}))$ one has 
\[
\Big\| \int_{- 1}^{\bar y} \frac{1}{y^\alpha} f(y)\, d y\Big\|_{H^s(\T^{\nu}\times \T^{d}_{\Gamma})}
\lesssim \| f \|_{L^2_y H^s(\T^{\nu}\times \T^{d}_{\Gamma})}\,. 
\]
Combining these arguments one can show that $\langle D\rangle^{M}\mathcal{S}(\vphi)$
is bounded on  $H^{s}(\T^\nu \times \T^{d}_{\Gamma})$.
This is the content of Lemma \ref{lemma stima resto cal S DN}.

\smallskip
\vspace{0.9em}
\noindent
$(iii)$ {\bf Invertibility and spectral analysis
of the linearized operator.} 
As explained in sections \ref{sec:approx_inv} and
\ref{sec:NaM}, the convergence of the Nash-Moser scheme is obtained by the approximate inverse Theorem \ref{alm.approx.inv}
which is based on the invertibility assumption 
\eqref{almi4} on the linearized operator in the normal directions (see \eqref{Lomegatrue}).
In the following we discuss the main ingredients to prove estimates like  
\eqref{almi4} on the operator (see also \eqref{representation Lom}-\eqref{linWW})\footnote{in this qualitative explanation we do not discuss the presence of the $L^{2}$ projector on the normal variables. We discuss the normal form strategy for the whole operator $\mathcal{L}$ without restriction to symplectic subspaces in \eqref{decomposition}.}
\begin{equation}\label{linWWINTRO}
\mathcal{L}:=\mathcal{L}(\vphi):=\omega\cdot\pa_{\vphi}+
\left(
\begin{matrix}
V\cdot\nabla+{\rm div}(V)+G(\eta)B & -G(\eta) \\
1 +B{\rm div}(V)+BG(\eta)B & V\cdot\nabla-BG(\eta)
\end{matrix}
\right)\,,
\end{equation}
where $G(\eta)$ is the Dirirchlet-Neumann operator, and $V,B$ denote
 the horizontal and vertical components of the velocity field at the free interface
 (see \eqref{def:V}-\eqref{form-of-B}).
 
 The overall approach to invert the linearized operator ${\mathcal L}$ is to conjugate 
 it to a diagonal operator with a sharp control of its eigenvalues. 
 This is absolutely not trivial in dimension greater or equal than two, 
 whereas in dimension one is more natural to expect a sharp asymptotic of the eigenvalues. 
 Indeed a quite general method has been developed  in dimension one 
 (see for instance \cite{BBHM,FGtrave}), 
 but it was totally unclear in dimension greater or equal than two. 
 Roughly speaking, the most difficult thing is the interaction between the 
 {\bf strong resonance phenomena (which are due to the dispersion relation and the high space dimension)} 
 and the {\bf quasi-linear nature of the equation}. 
 We develop a normal form method (which looks quite general for quasi-linear, dispersive, 
 translation invariant PDEs with sublinear dispersion in dimension greater or equal than two), 
 which is based on a normal form reduction whose main ideas are the following: 
 since the dispersion relation is  sublinear, one has $\omega \cdot \partial_\vphi$ 
 is of order one whereas $\Omega(D) = \sqrt{|D| \tanh(\mathtt h |D|)}$ 
 is of order $1/2$, therefore the effect of the time derivative is stronger. 
 This allows to develop a reduction scheme, based on pseudo-differential calculus 
 that removes the dependence on time from the linearized operator. 
 In principle the reduced operator is not diagonal in space. 
 Then one exploits the conservation of momentum, that at the linear 
 level can be seen as the invariance of the linearized 
 operator on the space of traveling wave functions, to have some exact cancellations. 
 In particular, one shows that if a momentum preserving operator is independent 
 of time then it is necessarily diagonal. 
 
 \noindent
 Of course, we observe additional problems. 
 From the analytical point of view, the pseudo-differential operators 
 appearing in the perturbation have symbols which are essentially general, 
 namely they do not admit expansion in decreasing symbols which are homogeneous 
 w.r.t. the variable $\xi$. This is basically due to the complicated structure 
 of the Dirichlet-Neumann operator (which is much easier in dimension one, as it is well known). 
 In addition to this, in the reducibility scheme which diagonalizes $\mathcal{L}$ in section 
\ref{sec:KAMreducibility} we have to deal with the possible multiplicity of the eigenvalues 
and presence of many non trivial resonances. It turns out that there are some algebraic 
exact cancellations 
(by exploiting again the conservation of momentum and the irrationality of the torus). 
On the other hand, there is still a bad control of the small divisors that 
have loss of derivatives both in time and space and the eigenvalues 
$\mu(j)$, $j \in \Gamma^*$ have no asymptotic expansion in decreasing powers of $|j|$.

In the following we shall describe
 in details this strategy by underlining the main novelties.  
% \emph{(i)}
% The first main difference is due to the structure of $G(\eta)$ in \eqref{claimoGeta}.
% In particular, differently from the one dimensional case, the symbols 
%in the pseudo-differential expansion of $G(\eta)$ do {\bf NOT} admit an expansion in decreasing symbols 
%homogeneous in the variable $\x$. This makes much more difficult the pseudo-differential reduction 
%of $\mathcal{L}$ to a diagonal operator up to smoothing remainders.
%
%\emph{(ii)} 
To simplify the notation, 
in this introduction we will say that an operator is a \emph{smoothing remainder}
if it satisfies tame estimates like \eqref{stimaRRGG}.

The normal form reduction is divided  into several steps.

\smallskip
\emph{(1) Good unknown and complex variables.} 
This step is based on a quite standard symmetrization argument, 
which do not present many big differences with the one dimensional case.
We just remark that  we make sharp quantitative estimates on the norms of the pseudo-differential operators
which are fundamental for the next steps.

First, we symmetrize the operator 
$\mathcal{L}$ in \eqref{linWWINTRO} at the highest order by introducing (following \cite{AB15}, \cite{BM20}) the \emph{linear} 
good unknown of Alinhac (see \eqref{flussoG})
%\[
%\mathcal{L}_1=\mathcal{G}^{-1}\mathcal{L}_{\omega}
%\mathcal{G}=\left( 
%\begin{matrix} 1 & 0\\ B & 1
%\end{matrix}
%\right)\,,
%\]
so that one has
\[
\mathcal{L}_1=\mathcal{G}^{-1}\mathcal{L}
\mathcal{G}=
\omega \cdot \partial_\vphi +\left(\begin{matrix}
V\cdot\nabla+ b(\vphi, x) & -G(\eta) \\
1 + a(\vphi, x) & V\cdot\nabla
\end{matrix}\right)+ {\mathcal R}_1\,, 
\qquad \mathcal{G}=\left( 
\begin{matrix} 1 & 0\\ B & 1
\end{matrix}
\right)\,,
\]
where $\mathcal{R}_1$ is smoothing is space and
\[
a(\vphi,x):=(\omega\cdot\pa_{\vphi}B)+V\cdot\nabla B\,, \quad b(\vphi, x) := {\rm div}(V)\,.
\]
In Proposition \ref{prop operatore cal L2},
passing to the complex coordinates through the map $\mathcal{C}$ in \eqref{CVWW}
and using the pseudo-differential expansion of $G(\eta)$ in \eqref{claimoGeta}, we obtain
the operator (see \eqref{Weyl} for the definition of the \emph{Weyl} quantization of symbols)
\begin{equation}\label{elle2INTRO}
\begin{aligned}
\mathcal{L}_2&= {\mathcal C} \circ {\mathcal L}_1 \circ {\mathcal C}^{- 1}
\\&=
\omega\cdot\pa_{\vphi}\uno+\opw(\ii V\cdot\x)\uno+\ii E\Omega(D)+\opw(A_2+B_2)+\mathcal{R}_2\,,
\end{aligned}
\end{equation}
where\footnote{In \eqref{elle2INTRO} we used 
the matrices defined as
\[
 \uno:=\sm{1}{0}{0}{1}\,,\;\;\;E:=\sm{1}{0}{0}{-1}\,.
 \]} $\mathcal{R}_2$ is  smoothing remainder, $\Omega(D)$ is the dispersive operator in \eqref{wild},
$A_2,B_2$ are $2\times2$ matrices of symbols respectively of order $1/2$ and $0$ 
(see \eqref{cal A2 Hamiltoniano 0} for more details and additional properties of symbols).
 We make some comments on the operator $\mathcal{L}_2$ above:
 
 \smallskip
 \noindent
 $(i)$ we first remark
that in these coordinates  $\mathcal{L}_2$ is \emph{Hamiltonian} 
and \emph{momentum preserving}
according to Definition \ref{operatoreHam}. 

\smallskip
\noindent
$(ii)$ Thanks to the estimates on symbols given in Theorem \ref{lemma totale dirichlet neumann} and 
quantitative bounds on compositions of pseudo-differential operators (see Lemma \ref{lemmacomposizione}),
we have a sharp control on the norms of symbols 
$A_2,B_2$, see \eqref{stima cal A 2 new}-\eqref{stima cal B 2 new}. 
In particular, the regularity in $\x$ is independent of the regularity in $x,\vphi$.
This property
is not trivial at all since we have to deal with non-homogeneous symbols in $\x$. 
Moreover
it will  be  fundamental in the following  to apply the quantitative version of an Egorov theorem 
given in Section \ref{sec:Egorov}.

\smallskip
\noindent
$(iii)$ At the highest order $\mathcal{L}_2$ is already diagonal and has the special form
\begin{equation}\label{figaro10}
\omega\cdot\pa_{\vphi}\uno+\opw(\ii V\cdot\x)\uno=
\big(\omega\cdot\pa_{\vphi}+V\cdot\nabla+\tfrac{1}{2}{\rm div}(V)\big)\sm{1}{0}{0}{1}\,.
\end{equation}

\vspace{0.9em}
\medskip
\emph{(2) Straightening of first order vector fields.} 
In this step (see Section \ref{sec:almoststraightening}) our aim is to conjugate the operator in \eqref{figaro10}
to the constant coefficient one
\begin{equation}\label{figaro11}
\omega\cdot\pa_{\vphi}\uno+\mathtt{m}_1\cdot\nabla\uno\,,
\end{equation}
for some constant vector $\mathtt{m}_1\in \R^{d}$.
To do this one has to deal with a small divisors problem (see lower bounds in \eqref{tDtCn}).
By exploiting the momentum conservation we follow the idea  of straightening 
of weakly perturbed constant vector field on tori given in \cite{FGMP19}
through a diffeomorphism of the torus $\mathbb{T}_{\Gamma}^{d}$, $x\mapsto x+\beta(\vphi,x)$. 
This approach has been 
fruitfully used also in \cite{BFM2}, \cite{FGtrave}, but we refer to \cite{BiFMT} 
for details in the high dimensional setting with the use of momentum.
The main differences here are the following: $(i)$  the transport vector field $V$ is not divergence free, 
but however, we construct a symplectic correction of the torus diffeomorphism to still obtain the conjugation to
\eqref{figaro11}; $(ii)$ in the procedure we have to control the $C^{k}$ norms with respect to 
the parameters $\omega$ and $\mathtt{h}$. 
The construction of the diffeomorphism and the conjugation of \eqref{figaro10} to \eqref{figaro11}
are obtained in Lemma \ref{conju.tr}.

We then use such diffeomorphism to conjugate the whole operator $\mathcal{L}_2$
 in \eqref{elle2INTRO}
 to 
 \begin{equation}\label{elle3INTRO}
\begin{aligned}
\mathcal{L}_3
=
\omega\cdot\pa_{\vphi}\uno+\mathtt{m}_1\cdot\nabla
+\ii E\Omega(D)+\opw(A_3+B_3)+\mathcal{R}_3\,,
\end{aligned}
\end{equation}
where $\mathcal{R}_3$ is  smoothing remainder,
$B_3$ is  $2\times2$ matrix of symbols of order $0$ while 
\[
A_{3}:=A_{3}(\vphi,x,\x):=\ii \begin{pmatrix}
a_3(\vphi, x, \xi) & b_3(\vphi, x, \xi) \\
- b_3(\vphi, x, \xi) & - a_3(\vphi, x,  \xi)
\end{pmatrix}
\]
is a matrix of symbols of order $1/2$ satisfies the algebraic properties \eqref{cal A3 Hamiltoniano}.
This is the content of 
Proposition \ref{riduzione trasporto}.
The main difficulty here  which do not appear in \cite{BiFMT} (or in the one dimensional cases)
is that 
we must conjugate under the torus diffeomorphism 
the  general pseudo-differential operators 
$\opw(A_2+B_2)$ that are not homogeneous in the variable $\x$. As already remarked this happens due to the complicated structure of the Dirichlet-Neumann operator.
To obtain such conjugation is a challenging technical problem. We exploit a quantitative version of 
the Egorov Theorem to have sharp estimates on the symbols of the conjugate operators.
This part is developed, following some of the ideas introduced in 
\cite{FGtrave}-\cite{BFPT1}, in section \ref{sec:egothm}.

\vspace{0.9em}
\medskip
\emph{3) Block-decoupling of the order $1/2$.} 
In section \ref{sym.low.orderZERO} we symmetrize the operator $\mathcal{L}_{3}$ in \eqref{elle3INTRO}
at order $1/2$. More precisely, our aim is to diagonalize the matrix of symbols
\begin{equation}\label{matrix12INTRO}
\ii E\Omega(\x)+\chi_0(\x)A_{3}(\vphi,x,\x)\,,
\end{equation}
where $\chi_0(\x)$ is the cut-off function in \eqref{cut-offCHIZERO}, up to 
regularizing remainders\footnote{note that, in view of the definition of the cut-off $\chi_0$ the symbols in  $(1-\chi_0(\x))A_{3}(\vphi,x,\x)$ are of order $-M$ for any $M>0$.}.
The basic idea
%, which has been implemented also in \cite{FGtrave}, 
is to consider the pseudo-differential 
change  of coordinates $\opw(F(\vphi,x,\x))$ where 
$F(\vphi,x,\x)$ is a $2\times2$ matrix of symbols of order zero
\[
F(\vphi, x, \xi) :=  \begin{pmatrix}
f_{+}(\vphi, x, \xi) & f_{-}(\vphi, x, \xi) \\
f_{-}(\vphi, x, \xi) & f_{+} (\vphi, x, \xi)
\end{pmatrix}
\]
where the columns of $F$ are the eigenvectors of the matrix \eqref{matrix12INTRO}.
Such matrix can be diagonalized since it has distinct, real eigenvalues
\[
\begin{aligned}
 \lambda_\pm(\vphi, x, \xi) &  = \pm \ii  \lambda(\vphi, x, \xi)\,,
  \\
\lambda(\vphi,x,\x) & :=\sqrt{  \Big(\Omega(\xi)   
+ \chi_0(\xi)a_3(\vphi, x, \xi) \Big)^2   +  \chi_0(\xi)^2b_3(\vphi, x, \xi)^2}\,.
\end{aligned}
\]
Note that $\lambda-\Omega(\x)$ is a symbol of order $1/2$, small in size. 
This is the content of Lemma \ref{lemma autovalori decoupling ordine alto}.
The matrix of eigenvectors $F(\vphi,x,\x)$ is  constructed explicitly in Lemma \ref{lem:lambdino}
in such a way 
\[
F(\vphi, x, \xi)^{- 1} \Big[\ii E \Omega(\xi) + A_3(\vphi,x, \xi) \Big]F(\vphi, x, \xi)  \\
=
\begin{pmatrix}
\ii \lambda(\vphi, x, \xi) & 0 \\
0 & - \ii \lambda(\vphi, x, \xi)\,.
\end{pmatrix}
\]
Therefore the map $\opw(F(\vphi,x,\x))$ is the good candidate to be the correct map which 
diagonalizes at order $1/2$ the operator $\mathcal{L}_{3}$. 
Unfortunately such a  map is not symplectic. 
The crucial point of section \ref{sym.low.orderZERO} is to find a suitable symplectic correction.
This is done by taking the flow $\Phi_{\mathcal{F}}=\exp\{(\opw(G(\vphi,x,\x)))\}$ generated by 
a suitable $2\times2$ matrix of pseudo-differential operator $G$, see \eqref{caffe6}. 
Such matrix $G$
is chosen in such a way $\Phi_{\mathcal{F}}-\opw(F(\vphi,x,\x))$ is a pseudo-differential operator of order $-1$.
This is actually done in Lemma \ref{lemma M  exp P}.
The crucial difference here with the one dimensional case \cite{FGtrave}
is that the symbols $f_{\pm}$ are general symbols of order zero while in \cite{FGtrave}
they simple are functions independent of $\x\in \R^{d}$.
On the contrary, in our case also the matrix $G(\vphi,x,\x)$ has general symbols as entries.
As a consequence to obtain a precise structure of $\Phi_{\mathcal{F}}$ is not trivial. We tackle this problem
by using Lemma \ref{lemma.potenze.sharp} which provide a precise analysis of the structure 
of the exponential of a pseudo-differential operator 
with sharp bounds on the symbols.
We conclude, in Proposition \ref{prop operatore cal L3}, 
by providing the conjugate of the operator $\mathcal{L}_{3}$, obtaining
\[
\mathcal{L}_{4}:=\omega\cdot\partial_{\vphi} 
+ \opw(\ii \mathtt{m}_{1}\cdot \xi) \uno + \ii E \opw(\lambda) 
+ {\mathcal B}_4 + {\mathcal R}_4\,,
\]
for some  smoothing reminder $\mathcal{R}_4$ and 
where ${\mathcal B}_4$
is a 
$2\times2$ matrix of pseudo-differential operators of
order $0$.
Finally we remark that not only the operator $\mathcal{L}_{4}$ is Hamiltonian and momentum preserving, 
but also 
$\opw(\ii \mathtt{m}_{1}\cdot \xi) \uno$, $ \ii E \opw(\lambda) $
and ${\mathcal B}_4$ are so separately.

\vspace{0.9em}
\medskip
\emph{4) Block-decoupling lower orders.} 
In Proposition \ref{blockTotale} of section \ref{subsec:loweroffdiag}
we symmetrize $\mathcal{L}_{4}$ up to arbitrarily smoothing remainders.
This is done by conjugating, inductively,  the operator $\mathcal{L}_{4}$ 
under the flow of pseudo-differential operators of negative 
orders. This step is more standard and do not present much difficulties or differences with 
respect to the
one dimensional case. 
At the end of the procedure we get that $\mathcal{L}_{4}$ is transformed into
\begin{equation}\label{opL5INTRO}
\mathcal{L}_{5}:=\omega\cdot\partial_{\vphi} 
+ \opw(\ii \mathtt{m}_{1}\cdot \xi) \uno + \ii E \opw(\lambda) 
+ {\mathcal B}_5 + {\mathcal R}_5\,,
\end{equation}
where $\mathcal{R}_{5}$ is some smoothing remainders, while
\[
\mathcal{B}_5=\opw\big(B_{5}(\vphi,x,\x)\big)\,,\qquad B_{5}(\vphi,x,\x)=\ii \left(
\begin{matrix}
b_{5}(\vphi,x,\x) & 0 \\ 0 & - b_{5}(\vphi,x,-\x)
\end{matrix}
\right)\,.
\]
We remark that, thanks to the Hamiltonian structure and the momentum conservation 
 (that we preserved along the previous steps) we have that 
the operator $\mathcal{B}_{5}$ is Hamiltonian and momentum preserving.

%symbol $b_{5}$ satisfies the following properties which will be fundamental in the following:

%This property turns out to be fundamental for the next step.

\medskip
\emph{5) Reduction at order $1/2$.} We now shall describe the reduction to constant coefficients of the order $\frac12$. In order to highlight the main idea, we simply consider the "scalar" operator 
\[
{\mathcal P}_5 := \omega \cdot \partial_\vphi + \mathtt m_1 \cdot \nabla 
+  \ii \opw\Big(\lambda(\vphi, x, \xi) \Big)\,.
\]
First of all we remark that since the operator is Hamiltonian and momentum preserving, 
on the symbol $\lambda(\vphi,x,\x)$ we have the following crucial properties:
\begin{itemize}

\item $\lambda$ is real valued;

\item the conservation of momentum of the operator $\opw(\lambda)$
is equivalent to say (see Lemma \ref{lem:mompressimbolo}) 
that the symbol $\lambda$ satisfies the following:
there is a smooth function $ \T^{\nu}\times \R^{d}\ni(\Theta,\x)\mapsto \widetilde{\lambda}(\Theta,\x)$
such that (recall Def. \ref{def:quasitravelling} and \eqref{velocityvec2})
\[
\lambda(\varphi, x, \xi)=\widetilde{\lambda}(\varphi-\tV x, \xi)\,.
\]
By passing to the Fourier side in $(\vphi,x)$, i.e. by writing 
$\lambda=\sum_{\ell\in \Z^{\nu}}\sum_{j\in \Gamma^*}\widehat{\lambda}(\ell,j,\x)e^{\ii\ell\cdot\vphi+\ii j\cdot x}$
the condition above imply that
\[
{\rm if }\;\;\;\widehat{\lambda}(\ell,j,\x)\neq0\quad {\rm then}\quad 
\tV^{T}\ell + j =0
\]
\end{itemize}
and thus the Fourier expansion of $\lambda$ becomes 
\begin{equation}\label{espansioneLambda}
\lambda(\varphi, x, \xi) = \sum_{\begin{subarray}{c}
(\ell, j) \in \Z^\nu \times \Gamma^* \\
\tV^{T}\ell + j  = 0
\end{subarray}} \widehat{\lambda}(\ell, j, \xi) e^{\ii \ell \cdot \vphi} e^{\ii j \cdot x}\,.
\end{equation}

Now, 
in order to remove the dependence on time and space $(\vphi, x)$ from the symbol $\lambda(\vphi, x, \xi)$ which is of order $1/2$, we consider the time one flow map $\Phi_g$ generated by $\ii \opw(g(\vphi,x, \xi))$ 
for some unknown symbol $g$ of order $1/2$. 
By standard Lie expansion series (see Lemma \ref{flussi coniugi eccetera}), we get that 
$$
\begin{aligned}
\Phi_g^{- 1} {\mathcal P}_5 \Phi_g & =  
\omega \cdot \partial_\vphi + \mathtt m_1 \cdot \nabla + \ii\opw\Big( (\omega \cdot \partial_\vphi + \mathtt m_1 \cdot \nabla ) g(\vphi, x, \xi) + \lambda(\vphi, x, \xi) \Big)  + l.o.t. 
\end{aligned}
$$
where by $l.o.t.$ we denote pseudo-differential operators of order smaller or equal than zero. 
Our aim is to solve the equation 
$$
(\omega \cdot \partial_\vphi + \mathtt m_1 \cdot \nabla) g(\vphi, x, \xi) + \lambda(\vphi, x, \xi) = \langle \lambda \rangle_\vphi(x, \xi)
$$
where 
$$
\langle \lambda \rangle_\vphi(x, \xi) := \frac{1}{(2 \pi)^\nu}\int_{\T^\nu} \lambda(\vphi, x, \xi)\, d \vphi\,. 
$$
By imposing the diophantine condition 
$$
|\omega \cdot \ell + \mathtt m_1 \cdot j| \geq \frac{\gamma}{|\ell|^\tau}, \quad \forall \ell \in \Z^\nu \setminus \{ 0 \}, \quad j \in \Gamma^* 
$$
($0 < \gamma \ll 1$, $\tau \gg 0$), one has that the latter equation has a solution of the form 
$$
g(\vphi, x, \xi) = - \sum_{\begin{subarray}{c}
\ell \in \Z^\nu \setminus \{ 0 \} \\
j \in \Gamma^*
\end{subarray}} \dfrac{\widehat \lambda(\ell, j, \xi)}{\ii (\omega \cdot \ell + \mathtt m_1 \cdot j)} e^{\ii (\ell \cdot \vphi + j \cdot x)}\,.
$$
Hence, in general we can remove only the dependence on time $\vphi \in \T^\nu$, which would not be sufficient for our aims.
The key idea to overcome this problem and to eliminate the space-time dependence from the symbol at order $1/2$, 
is to exploit
%But now it enters the key point that exploits 
the conservation of momentum. 
Indeed, since the symbol
$\lambda$ is momentum preserving, one can compute the average in $\vphi$ of $\lambda$
by recalling \eqref{espansioneLambda} 
obtaining that
% then we look even for $g$ which is momentum preserving. 
%Moreover if one computes $\langle \lambda \rangle_\vphi(x, \xi)$ one obtains that 
\[
\begin{aligned}
\langle \lambda \rangle_\vphi(x, \xi) & = \int_{\T^\nu} \Big( \sum_{\begin{subarray}{c}
(\ell, j) \in \T^\nu \times \Gamma^* \\
\mathtt V^T \ell + j = 0
\end{subarray}} \widehat \lambda(\ell, j, \xi) e^{\ii (\vphi \cdot \ell + j \cdot x)}\Big) \, d \vphi \\
& = \sum_{\begin{subarray}{c}
(\ell, j) \in \T^\nu \times \Gamma^* \\
\mathtt V^T \ell + j = 0
\end{subarray}} \widehat \lambda(\ell, j, \xi) e^{\ii j \cdot x} \int_{\T^\nu} e^{\ii \vphi \cdot \ell } \, d \vphi  
%\\& 
= \widehat \lambda(0, 0, \xi)
\end{aligned}
\]
since $\ell = 0$ implies $\mathtt V^T \ell = 0$ and hence $j = 0$. 

Therefore, for a momentum preserving pseudo-differential operator, 
the time average is actually constant w.r.t. the variable $x$ 
(i.e. it is a Fourier multiplier $\widehat \lambda(0, 0, \xi)$). 
Hence, with this additional cancellation due to the momentum, we get that 
\[
g(\vphi, x, \xi) = - \sum_{\begin{subarray}{c}
\ell \in \Z^\nu \setminus \{ 0 \}\,,\, j \in. \Gamma^* \\
\mathtt V^T \ell + j = 0
\end{subarray}} \dfrac{\widehat \lambda(\ell, j, \xi)}{\ii (\omega \cdot \ell + \mathtt m_1 \cdot j)} 
e^{\ii (\ell \cdot \vphi + j \cdot x)}
\]
solves the  equation
$$
(\omega \cdot \partial_\vphi + \mathtt m_1 \cdot \nabla) g(\vphi, x, \xi) + \lambda(\vphi, x, \xi) =  \widehat\lambda(0, 0, \xi)\,. 
$$
In particular, this also implies that $g$ is itself momentum preserving.
What we discussed above is implemented in Lemma 
\ref{proposizione equazione omologica generale riduzione in ordini}.
To summarize (see  Proposition \ref{riduzione diagonal ordine 1 2}) we conjugate the operator
$\mathcal{L}_{5}$ in \eqref{opL5INTRO}
\begin{equation}\label{opL6INTRO}
\mathcal{L}_{6}:=\omega\cdot\partial_{\vphi} 
+ \opw(\ii \mathtt{m}_{1}\cdot \xi) \uno + \ii E \opw(\mathtt{m}_{6}(\x)) +
\ii \opw\left(
\begin{matrix}
b_{6}(\vphi,x,\x) & 0 \\ 0 & - b_{6}(\vphi,x,-\x)
\end{matrix}
\right)
 + {\mathcal R}_6\,,
\end{equation} 
where ${\mathcal R}_6$ is a smoothing remainder, 
$b_{6}(\vphi,x,\x)$ a real valued symbol  of order $0$, and where we defined the Fourier multiplier 
$\mathtt{m}_{6}(\x):=\widehat{\lambda}(0,0,\x)$.
Notice that $\mathtt{m}_{6}(\x)-\Omega(\x)$ is a symbol (independent of $x$) of order $1/2$
and small in size (see estimates \eqref{stima.m.6}).

\medskip
\emph{6) Reduction at order lower orders.}
In section \ref{sec:reductionlower} we iterate the procedure described in the step above, with the aim of reducing the operator $\mathcal{L}_{6}$
to a constant coefficients one up to smoothing remainders.  
The key idea is still contained in Lemma \ref{proposizione equazione omologica generale riduzione in ordini}.
At the end of the procedure, see Proposition \ref{riduzione diagonal ordini bassi},
we get that $\mathcal{L}_{6}$ in \eqref{opL6INTRO} becomes
\begin{equation}\label{opL7INTRO}
\mathcal{L}_{7}=\omega\cdot\partial_{\vphi} 
+ \opw(\ii \mathtt{m}_{1}\cdot \xi) \uno + \ii E \opw(\mathtt{m}_{7}(\x)) 
 + {\mathcal R}_7\,,
\end{equation} 
where $\mathcal{R}_{7}$ is a smoothing in space remainder satisfying tame estimates like
\eqref{stima tame cal R7 new},
and where $\mathtt{m}_{7}(\x)$ is a real valued, independent of 
$(\vphi,x)\in \mathbb{T}^{\nu}\times\T^{d}_{\Gamma}$ symbol  of order $1/2$.
In particular, it is close to $\Omega(\x)$ (see estimates \eqref{stima lambda 7}).

\medskip
\emph{7) KAM reducibility scheme.} 
In the previous steps from $1)$ to $6)$, we have obtained (taking also into account the $L^{2}$ projectors on the normal variables)
that the linearized operator $\mathcal{L}_{\omega}$ in \eqref{representation Lom} is transformed into
(see \eqref{operatoreintero} in Theorem \ref{red.Lomega.smooth.rem})
\[
\cL_7 
=  \omega \cdot \partial_\vphi {\mathbb I}_{\perp}+ \ii \cD_\bot + \cR_\bot 
\]
where $\mathcal{D}_{\perp}$
 is a diagonal operator of the form
 \[
 {\mathcal D}_{\perp} := \left(\begin{matrix}
\mathtt{D}_{\perp} & 0
\\  0 &-\overline{\mathtt{D}}_{\perp}
\end{matrix}
\right)\,,
\qquad \mathtt{D}_{\perp}:={\rm diag}_{j \in S_0^c} \mu_0 (j)\,,\qquad 
\mu_0 (j) :=  \mathtt{m}_{1} \cdot j + \mathtt{m}_{7}(j)\,, 
 \]
 where we defined $S_{0}^{c}:=\Gamma^{*}\setminus(S\cup\{0\})$, while 
 $\cR_\bot $ is a smoothing in space operator satisfying tame  estimates  as in 
 \eqref{tame riducibilita cal R bot iniziale}.
In section \ref{sec:KAMreducibility}
we implement a (more classical) reducibility KAM procedure in order 
to completely diagonalize $\cL_7 $.
Differently for the previous steps $1)-6)$, the aim is to construct a change of variables which 
reduces quadratically the size of the remainders at  each step of the procedure.
As it is well know, this requires some quantitive lower bounds 
on differences of eigenvalues, the so called second Melnikov conditions. 
A crucial problem, typical in PDEs in higher space dimension, 
is the presence of strong resonance phenomena. As a consequence one typically 
is able to only impose very ``weak'' Melnikov conditions, namely we require
\begin{equation}\label{secondeINTRO}
 |\omega \cdot \ell + \mu_0(j; \lambda, i) \pm  \mu_0(j'; \lambda, i)| 
\geq \frac{2  \,\gamma}{\langle \ell \rangle^\tau | j' |^\tau } \,, \qquad 
 		   \ell \neq0\,, \;\; {\mathtt V}^T\ell  + j \pm j' =0\,.
\end{equation}
By requiring the bounds above, one is able to solve the homological equation
\begin{equation}\label{omoINTRO}
\ii\omega\cdot\pa_{\vphi}\Psi_{0}(\vphi)+[ \ii \cD_\bot, \Psi_{0}(\vphi)]+\mathcal{R}_{\perp}(\vphi)=
\widehat{\mathcal{R}}_{\perp}(0)\,,
\qquad \widehat{\mathcal{R}}_{\perp}(0)=\int_{\T^{\nu}}\mathcal{R}_{\perp}
(\vphi) d\vphi\,,
\end{equation}
so that one has 
\[
e^{-\Psi_0}\mathcal{L}_{7}e^{\Psi_0}=
 \omega \cdot \partial_\vphi {\mathbb I}_{\perp}+ \ii \cD_\bot + \widehat{\mathcal{R}}_{\perp}(0)+ \cR_1\,,
\]
where $\mathcal{R}_1$ is still a smoothing in space operator, but with smaller size 
$\mathcal{R}_1\simeq \mathcal{R}_{\perp}^2$.
Again, in view of the conservation of momentum, one can infer that 
$\ii \cD_\bot + \widehat{\mathcal{R}}_{\perp}(0)$ is diagonal in the Fourier basis.
In order to make rigorous the reasoning above, one needs suitable estimates on the solution
$\Psi_0$ of \eqref{omoINTRO}, which, in the Fourier basis is represented by the matrix\footnote{here to simplify the notation we are assuming that the operator $\Psi_{0}$ acts on \emph{scalar} functions. Actually $\Psi_0$ is a $2\times2$ matrix of operators. For details we refer to \ref{shomo1}.}
coefficients  (see \eqref{matrrep})
\[
\widehat{\Psi}_{0}(\ell)_{j}^{j'}=\frac{\widehat{\mathcal{R}}_{\perp}(\ell)_{j}^{j'}}{ \ii (\omega \cdot \ell 
+ \mu_0(j)-  \mu_0(j'))}\,,
\qquad 
 		   \ell \neq0\,, \;\; {\mathtt V}^T\ell  + j - j' =0\,.
\]
It is  important to remark  that the  weak lower bounds \eqref{secondeINTRO} produces  an a priori
loss regularity in both time and space in estimating the operator $\Psi_0$.
Nevertheless, 
we are able to control such loss exploiting the special property of the remainder
$\mathcal{R}_{\perp}$, which are preserved along the KAM iteration.
First of all it is momentum preserving.
Secondly one has that, see \eqref{tame riducibilita cal R bot iniziale} and the discussion on the remainders of the Dirichlet-Neumann,
the operators 
\[
\langle D \rangle^{M} \partial_\vphi^\beta \cR_\bot \langle D \rangle^{-\mathtt{c}}\,,
\]
for $M\geq \mathtt{c}\gg1$  and $|\beta|\leq \beta_0$ with with $M,\beta_0$ arbitrarily large,
are bounded on the space-time Sobolev spaces $H^{s}(\T^{\nu}\times\T^{d}_{\Gamma})$, with ``tame'' estimates.

Therefore, (i) the loss of derivatives in time will be controlled by the regularity in 
$\vphi$ of $\mathcal{R}_{\perp}$, which essentially imply a strong decay in $|\ell|$ of the 
  coefficients $\mathcal{R}_{\perp}(\ell)_{j}^{j'}$;
(ii) the loss of derivatives in space will be controlled by the fact that 
the operator with coefficients $|j'|^{\tau}\mathcal{R}_{\perp}(\ell)_{j}^{j'}$
is still bounded.
Also here a key point is the use of the momentum condition which ensure
a constraint  among the indexes $j,j'$ and $\ell$.
This ingredients are used in the crucial Lemma \ref{stime eq omologica} to provide sharp estimates on the 
conjugating map.
The main non trivial, and non standard 
 difficulty, is to check that the map $e^{\Psi_0}$ is such that the new remainder obtained after the conjugation has
 still the same properties of $\cR_{\perp}$.
 This technical point is handled in the proof of estimates
  \eqref{stima cal R nu} inductively in Proposition \ref{iterazione riducibilita}.

  \smallskip
  \vspace{0.9em}
\noindent
$(iv)$ 
{\bf Non-resonance conditions.}
In section \ref{sec:nash} we rescale $(\eta,\psi)\mapsto (\e\eta,\e\psi)$
and we introduce on 
 the tangential sites the action angle variables $(\theta,I)\in \T^{\nu}\times \R^{\nu}$ in \eqref{definizione.v}. 
Therefore the system \eqref{eq:113} becomes the Hamiltonian system 
generated by the Hamiltonian function
\[
H_{\e}=
\overline{\omega}(\th)\cdot I+\frac{1}{2}({\bf \Omega} z, z)_{L^{2}}+\e P
\]
where $P$ is given in \eqref{P.epsilon}, $\bar{\omega}(\th)$ in 
\eqref{LinearFreqWW} 
and ${\bf \Omega}$ in \eqref{HamLinearReal}.
In order to look for quasi-periodic solutions, with frequencies $\omega \simeq \bar \omega(\mathtt h)$, 
of the  Hamiltonian system generated by $H_{\e}$
it is convenient (as done for instance in \cite{BM20, BBHM, FGtrave})
to use $\omega\in\R^{\nu}$ as parameters and introduce ``counterterms'' $\alpha,\mu\in \R^{\nu}$
in the family of Hamiltonians
\[
H_{\alpha,\mu}:=\alpha\cdot I+\frac{1}{2}({\bf \Omega} z, z)_{L^{2}}+\e P+\mu\cdot\theta\,,
\]
(see \eqref{HHalphaMumu} for details). Then in Theorem \ref{NMT}  we prove that, for $\e$ small enough,
there exist a map $ \alpha_\infty : \mathtt{\Omega} \times [\th_1,\th_2] \mapsto \R^\nu $ and 
a \emph{traveling} embedded torus (see \eqref{definition mathtt Omega}, \eqref{RTTT}) 
$i_{\infty}=i_{\infty}(\omega,\th,\varepsilon;\vphi)$
of the form
\[
i_{\infty} \,:\, \mathbb{T}^{\nu} \,\to \,\mathbb{T}^{\nu}\times\mathbb{R}^{\nu}\times H_{S}^{\perp}\,,
\qquad \vphi\mapsto i_{\infty}(\vphi)=(\theta_{\infty}(\vphi),I_{\infty}(\vphi),z_{\infty}(\vphi)) \,,
\]
such that, for any choices of the parameters $(\omega,\mathtt{h})$
belonging to a suitable set $\cG_{\infty} \subseteq  \mathtt \Omega \times [\th_1, \th_2] $,
$(i_\infty,\alpha_\infty, 0)(\omega,\th)$
is a zero of the nonlinear  operator (see \eqref{F_op})
\[
\cF(\omega,\th,\varepsilon;i,\alpha, \mu) := \omega\cdot\pa_\vphi i(\vphi) - X_{H_{\alpha, \mu}}(i(\vphi)) \,.
\]
Theorem \ref{NMT} is based on a Nash-Moser 
scheme in which, at each step, we construct an ``approximate'' invariant torus $i_{n}$\footnote{actually we construct $(i_n,\alpha_n, \mu_{n})(\omega,\th)$ in such a way 
$\mathcal{F}(i_n,\alpha_n, \mu_{n})\sim \delta_{n}$ for some $\delta_{n}\to0$ as $n\to+\infty$}
for the system.
The key point in the result, is not only to prove the convergence of the scheme, but, more important, is also to characterize the set of 
``good'' parameters $\mathcal{G}_{\infty}$ for which the convergence occurs and providing measure estimates
to conclude that the solution exists for ``most'' choices of the parameters.

The set $\mathcal{G}_{\infty}$ is defined  implicitly in the following way. 
As explained before in this introduction, 
in the normal form procedure for the linearized operator $\mathcal{L}=\mathcal{L}(i)$
 in \eqref{linWWINTRO}, 
 we need several non-resonance conditions in order to deal with small divisors, which of course depend on the step of the Nash-Moser. 
 More precisely, 
 given $\lambda=(\omega,\th)\in \mathtt{\Omega}\times[\th_1,\th_2]$ we denote by 
 $\lambda\mapsto \mu_{\infty}(j;\lambda, i)$
 the eigenvalues 
 of the linearized $\mathcal{L}(i)$ at a traveling embedding $i=i(\vphi,\lambda)$.
 We must require that for suitable choices of $\lambda$ one satisfies 
 $0$-th, $1$-st and $2$-nd Melnikov conditions,
 namely (recall \eqref{figaro11})
 \begin{equation}\label{Melnikov fine intro}
 \begin{aligned}
|(\omega-\tV\mathtt{m}_{1}(\lambda,i))\cdot\ell| 
				&\geq  4\gamma\jap{\ell}^{-\tau} \,,\qquad  \forall\,\ell\neq 0\,,
				\\
	| \oo\cdot \ell + \mu_{\infty}(j;\lambda,i)| 
&\ge 2\gamma \jap{\ell}^{-\tau}\,, 
\quad \ell\neq0\,,\;\; \tV^{T}\ell + j= 0\,,
\\
 |\omega \cdot \ell + \mu_\infty(j; \lambda, i) \pm  \mu_\infty(j'; \lambda, i)| 
&\geq \frac{2  \,\gamma}{\langle \ell \rangle^\tau | j' |^\tau } \,, \qquad 
 		   \ell \neq0\,, \;\; {\mathtt V}^T\ell  + j \pm j' =0\,.
\end{aligned} 
 \end{equation}
 In \eqref{Melnikov-invert} we denote by ${\bf \Lambda}_{\infty}^{\gamma}(i)$ 
 the set of $\lambda\in \mathtt{\Omega}\times[\th_1,\th_2]$ such that the conditions above holds
(see  \eqref{tDtCn}, \eqref{prime.di.Melnikov}, \eqref{Cantor set} for the precise definition of the set of good parameters).
We then define $\mathcal{G}_{\infty}$
as the intersection (see \eqref{defGinfty})
\[
\mathcal{G}_{\infty}:=\mathcal{G}_{\infty}(\gamma):=\big({\mathtt \Omega}\times[\th_1,\th_2]\big) \cap \bigcap_{n \geq 0} 
 {\bf \Lambda}_{\infty}^{\gamma_n}(\tilde \imath_{n}) \,,\qquad \gamma_n := \gamma(1 + 2^{- n})\,,
\]
where $\tilde \imath_{n}$ is the approximately invariant torus at each step of the Nash-Moser.
It is clear that, by definition, on the set $\mathcal{G}_{\infty}$ one is able to obtain the complete reducibility
of the linearized operator $\mathcal{L}_{\omega}(\tilde \imath_{n})$ for any $n\geq0$.
On the other hand, in this way we are imposing \emph{infinitely many} Melnikov condition
at \emph{each} step $n$ of the Nash-Moser. This means that it is not obvious that 
$\mathcal{G}_{\infty}$ is ``large'' in measure. This part is highly technical and it is done in all details in Section \ref{subsec:measest}. 
We focus our discussion on the 2nd Melnikov conditions ( third estimate in \eqref{Melnikov fine intro}), which are the most difficult to verify. The main difficulties related to our problem are the following.
\begin{itemize}
\item {The eigenvalues $\mu_\infty(j; \cdot)$ do not have an asymptotic expansion in powers of $\frac{1}{|j|}$. This is a typical problem in higher dimension and, as one can see in \cite{BBHM}, \cite{BFM2}, and \cite{FGtrave}, such an expansion has been deeply exploited, in one dimensional cases, to prove measure estimates.  In our case the eigenvalues $\mu_\infty(j)$ has the expansion}
$$
\mu_\infty(j) = \mathtt m_1 \cdot j + \mathtt m_7(j) + O(\e |j|^{- N}) \quad \text{for some} \quad N \gg 0\,. 
$$
The term $\mathtt m_7(j)$ satisfies $|\mathtt m_7(j) - \Omega(j)| \lesssim \e |j|^{\frac12}$ so it is a small correction of the dispersion relation with the same growth in $j$. 
On the other hand,  the crucial property is that $\mathtt m_7 \in S^{\frac12}$ is a Fourier multiplier of order $1/2$. In particular, we strongly use the fact that 
$|\nabla \mathtt{m}_{7}(j)|\lesssim \langle j \rangle^{-1/2}$
%$|(\mathtt m_7 - \Omega)(j) - (\mathtt m_7 - \Omega)(j')| \lesssim \e |j - j'|$ for any $j, j' \in \Gamma^*$, 
which is clearly consequence of the fact that $\mathtt m_7$ is a symbol of order $1/2$. 
This fact is crucially  used in Lemma \ref{lemma inclusione cantor}. 

\item As a byproduct of the Nash-Moser scheme, it turns out that the final 
frequencies of oscillations $\omega$ will be close to the unperturbed ones $\bar{\omega}(\th)=\left(
\omega_{\overline{\jmath}_{i}}(\th)
\right)_{i=1,\ldots,\nu}$, see \eqref{LinearFreqWW}.
The first fundamental step is to prove the Melnikov conditions at the unperturbed level, 
namely
\[
\begin{aligned}
&|\bar{\omega}(\th)\cdot\ell+\bar{\omega}_{j}(\th)-\bar{\omega}_{j'}(\th)|
\geq \frac{\gamma}{\langle \ell,j'\rangle^{\tau}}\,,
\qquad \ell\in \Z^{\nu}\setminus\{0\}\,,\;\;\;j,j'\in \Gamma^{*}\setminus S\,,
\\
&\mathtt{V}^{T}\ell+j-j'=0\,.
\end{aligned}
\]
This is the content of Section \ref{sec:degener}. A problem in verifying this bound, connected with the high dimensional setting,
is that there could be exact resonances between tangential sites and normal ones.
For instance it could happen that, for some $j,j'\in \Gamma^{*}\setminus S$, there exist
$i,k=1,\ldots,\nu$, such that $|j|=|\overline{\jmath}_{i}|$ and $|j'|=|\overline{\jmath}_{k}|$.
This fact is avoided  thanks to the irrationality of the lattice required in Hypothesis \ref{hyp:lattice}.
Such an assumption guarantees that (see Lemma \ref{lem:irrazionale})
if $|j|=|\overline{\jmath}_{i}|$ and $|j'|=|\overline{\jmath}_{k}|$ then $j=\pm \overline{\jmath}_{i}$ and $j'=\pm \overline{\jmath}_{k}$. Of course since $j,j'$ do no belong to the set $S$   one can only have the case 
$j=-\overline{\jmath}_{i}$ and $j'=- \overline{\jmath}_{k}$.
Therefore the small divisor takes the form
\[
\bar{\omega}(\th)\cdot\ell+\bar{\omega}_{j}(\th)-\bar{\omega}_{j'}(\th)=\bar{\omega}(\th)\cdot\big(\ell+ e_i-e_k\big)\,,
\]
where $\{e_{i}\}_{i=1,\ldots,\nu}$ denotes the standard basis of $\R^{\nu}$.
Then if the vector $\ell+ e_i-e_k\neq0$ the desired lower bounds follows by a Diophantine condition on 
$\bar{\omega}(\th)$.
If on the contrary one has $\ell+ e_i-e_k=0$ then the momentum condition is violated, indeed
one has
$\mathtt{V}^{T}\ell+j-j'= 2(\overline{\jmath}_{k}-\overline{\jmath}_{i})\neq0$, which holds true since, by assumption, all the tangential sites have different modulus.

\end{itemize}

\section{Functional setting}

\subsection{Hamiltonian structure of water waves and action angle coordinates.}
Recall the lattice $\Gamma$ introduced in \eqref{def:lattice}-\eqref{def:lattice2}.
We expand a function $u(x)$ in $L^{2}(\mathbb{T}_{\Gamma}^{d};\C)$
as 
\[
u(x)=\sum_{j\in\Gamma^{*}}\widehat{u}(j)e^{\ii j\cdot x}\,,
\qquad \widehat{u}(j)=\frac{1}{|\T^d_\Gamma|}\big(u, e^{\ii j\cdot x}\big)_{L^2_{\Gamma}}\,,
\qquad \forall \, j\in \Gamma^{*}\,,
\]
where
%\footnote{
%For simplicity in the following we shall use the following notation:
%$u^{+}_{j}=u_{j}=\widehat{u}(j)$, \;\; 
%$u^{-}_{j}=\ov{u}_{j}=\ov{\widehat{u}(j)}$.
%}
 $(\cdot,\cdot)_{L^2_{\Gamma}}$ is the standard Hermitian product
\begin{equation}\label{scalprod}
(u,v)_{L^2_{\Gamma}}:=\int_{\mathbb{T}^{d}_{\Gamma}} u(x)\ov{v(x)}dx\,,
\end{equation}
and $|\T^d_\Gamma|$ is the Lebesgue measure of $\T_{\Gamma}^{d}$.
We consider the scale of Sobolev spaces, $s\in\R$,
\begin{equation}\label{spaziosoloX}
H^{s}_{x}:=H^{s}(\mathbb{T}_{\Gamma}^{d};\C):=
\Big\{
u\in L^{2}(\mathbb{T}_{\Gamma}^{d};\C)\; : \; \|u\|^{2}_{H_{x}^s}
:=\sum_{j\in \Gamma^{*}} \langle j\rangle^{2s} |\widehat{u}(j)|^{2}<+\infty
\Big\}\,,
\end{equation}
where
$\langle j\rangle:=\sqrt{1+|j|^{2}}$ for $j\in \Gamma^*$.
\begin{rmk}
Notice that the space $H^{s}(\mathbb{T}_{\Gamma}^{d};\C)$ with $s>d/2$ is a Banach algebra.
\end{rmk}
The classical 
Craig-Sulem-Zakharov formulation in
\cite{Zak1}, \cite{CrSu} guarantees that equations \eqref{eq:113} can be seen as 
an infinite dimensional dynamical system on the phase space 
$H_0^{1}(\mathbb{T}_{\Gamma}^{d};\R)\times \dot{H}^{1}(\mathbb{T}_{\Gamma}^{d};\R) $
(see \eqref{phase-space}). Let us consider
 the symplectic form
\begin{equation}\label{symReal}
\mathcal{W}(W_1,W_2)
:=\int_{\mathbb{T}_{\Gamma}^{d}}W_1
\cdot J^{-1}W_2dx=
\int_{\mathbb{T}_{\Gamma}^{d}}-(\eta_1\psi_2-\psi_1\eta_2)dx\,,
\quad
J=\sm{0}{1}{-1}{0}\,,
\end{equation}
for 
\[
W_1=\vect{\eta_1}{\psi_1}\,, W_2=\vect{\eta_2}{\psi_2}\in 
H_0^1(\mathbb{T}_{\Gamma}^{d};\R)\times \dot{H}^{1}(\mathbb{T}_{\Gamma}^{d};\R)\,.
\]
%On the phase space 
%$H_0^{1}(\mathbb{T}_{\Gamma}^{d};\R)\times \dot{H}^{1}(\mathbb{T}_{\Gamma}^{d};\R) $
One can note that 
the vector field (see \eqref{HS})
\begin{equation*}%\label{realvecWW}
X(\eta,\psi)=X_{H}(\eta,\psi)=J \nabla H(\eta,\psi)
=\left(\begin{matrix}\nabla_{\psi}H(\eta,\psi) \\ 
-\nabla_{\eta} H(\eta,\psi) \end{matrix}\right)
\end{equation*}
 is the Hamiltonian 
vector field of $H$ in \eqref{HamiltonianWW} w.r.t. the symplectic form
\eqref{symReal}, i.e.
\[
-\mathcal{W}(X_{H}(\eta,\psi), h)=dH((\eta,\psi))[h]\,\quad 
h=\vect{\hat{\eta}}{\hat{\psi}}\,,
\]
while the Poisson  bracket between functions 
$F(\eta,\psi), H(\eta,\psi)$
%in \eqref{poissonBra}
are defined as
\begin{equation}\label{poissonBra22}
\{F,H\}=\mathcal{W}(X_{F},X_{H})
=\int_{\mathbb{T}_{\Gamma}^{d}}
\big(\nabla_{\eta}H\nabla_{\psi}F-\nabla_{\psi}H\nabla_{\eta}F\big)dx\,.
\end{equation}
Recalling \eqref{blackwhite}, \eqref{wild} 
we introduce the symplectic complex variables 
\begin{equation}\label{CVWW}
\begin{aligned}
& 
\left(\begin{matrix}\eta \\ \psi\end{matrix}\right) 
= \mathcal{C}^{-1} 
\left(\begin{matrix}
u \\ \bar u \end{matrix}\right) = 
\frac{1}{\sqrt{2}} 
\left(\begin{matrix}
\Omega^{1/2}(D)( u + \bar u ) \\
 - \ii  \Omega^{-1/2}(D)(  u - \bar u )
\end{matrix}\right)\,, 
\\& 
\begin{pmatrix}
u \\
\bar u \end{pmatrix} = \mathcal{C}\begin{pmatrix}
\eta \\
\psi
\end{pmatrix} = \frac{1}{\sqrt{2}}\begin{pmatrix}
\Omega^{-1/2}(D)\eta+\ii \Omega^{1/2}(D)\psi \\
\Omega^{-1/2}(D)\eta- \ii \Omega^{1/2}(D)\psi
\end{pmatrix}\,.
\end{aligned}
\end{equation}
The Poisson bracket in \eqref{poissonBra22} 
assumes the form 
\begin{equation}\label{poissonBraComp}
\{F,H\}=\frac{1}{\ii}
\int_{\mathbb{T}_{\Gamma}^{d}}(\nabla_{u}H\nabla_{\bar{u}}F-\nabla_{\bar{u}}H\nabla_{{u}}F)dx\,.
\end{equation}
Moreover,
in these coordinates, the vector field $X=X_{H}$ in \eqref{eq:113}  
assumes the form
(setting $H_{\mathbb{C}}:=H\circ\mathcal{C}^{-1}$)
\begin{equation}\label{compvecWW}
X^{\mathbb{C}}:=X_{H_{\C}}=
 \left(
\begin{matrix}
-\ii\pa_{\bar{u}}H_{\mathbb{C}} \\
\ii\pa_{u}H_{\mathbb{C}}
\end{matrix}
\right)\,.
%=
%\frac{1}{\sqrt{2 \pi}} \sum_{k\in\Z\setminus\{0\}} 
%  \left(\begin{matrix}
%- \ii  \pa_{ \ov{u_k}} H_{\C} \, e^{\ii k x}  \\
% \ii  \pa_{u_k} H_{\C}   \, e^{- \ii k x} 
%\end{matrix} 
%\right) \, ,
\end{equation}

\noindent 
\textbf{Splitting of the phase space.}
Recall $S$ in \eqref{TangentialSitesWW} and define 
$S_0^{c}:=\Gamma^{*}\setminus\big(S\cup\{0\}\big)$.
We decompose the phase 
space as follows. For any $j \in \Gamma^* \setminus \{ 0 \}$, we define 
\begin{equation}\label{definition Vj}
V_j := \Big\{ \begin{pmatrix}
\Omega(j)^{\frac12} \Big( \alpha_j \cos(j \cdot x) + \beta_j \sin(j \cdot x) \Big) \\
\Omega(j)^{- \frac12} \Big( \alpha_j \sin(j \cdot x)  - \beta_j \cos(j \cdot x)  \Big)
\end{pmatrix} : \alpha_j , \beta_j \in \R \Big\}
\end{equation}
and we denote by $\Pi_j$, $j \in \Gamma^*$ the corresponding orthogonal projectors. We define 
\[
H_S := \oplus_{j \in S} V_j\,,, \quad H_S^\bot := \oplus_{j \in S_0^{c}} V_j\,.
\]
Note that $V_j, H_S, H_S^\bot$ are symplectic, Lagrangian subspaces and we have that 
\begin{equation}\label{decomposition}
\begin{aligned}
 H_0^1(\mathbb{T}_{\Gamma}^{d};\R)\times 
 \dot{H}^1(\mathbb{T}_{\Gamma}^{d};\R) &:=H_S\oplus H_S^{\perp}\,.  
\end{aligned}
\end{equation}
We denote by $\mathbb{P}_S, \mathbb{P}_S^{\perp}$ 
the corresponding orthogonal projectors on $H_S, H_S^\bot$. 
The subspaces $H_S$ and $H_S^{\perp}$ are symplectic 
orthogonal respect to the $2$-form 
$\mathcal{W}$ (see \eqref{symReal}). 
On tangential sites we introduce the 
action angle variables $(\theta,I)\in \T^{\nu}\times \R^{\nu}$
as follows:
\begin{equation}\label{definizione.v}
\begin{aligned}
v\, :\, &\T^{\nu}\times \R^{\nu} \to H_{S}
\\&
\;\;(\theta,I)\mapsto v(\theta,I)=\sum_{j\in S} \begin{pmatrix}
\Omega(j)^{\frac12} \Big( \alpha_j \cos(j \cdot x) + \beta_j \sin(j \cdot x) \Big)
% \nonumber
\\
\Omega(j)^{- \frac12} \Big( \alpha_j \sin(j \cdot x)  - \beta_j \cos(j \cdot x)  \Big)
\end{pmatrix}
\end{aligned}
\end{equation}
with
\begin{equation*}
\begin{aligned}
\alpha_{j}&:=\sqrt{2(\zeta_{j}+I_{j})}\cos(\theta_j)\,,
\qquad 
\beta_{j}:=  \sqrt{2(\zeta_{j}+I_{j})}\sin(\theta_j)\,, 
\quad j\in S\,,
\end{aligned}
\end{equation*}
where
$\zeta_{j}>0$, $j=1,\ldots,\nu$ and the variables $I_{j}$
satisfy $|I_{j}|<\zeta_{j}$.
This induces a map
\begin{equation}\label{actionanglevar}
\begin{aligned}
A\, :\, \T^{\nu}\times \R^{\nu}\times H_{S}^{\perp}\;&\to \; H_{S}\oplus H_{S}^{\perp}
\\
(\theta,I,z)\;\;\;&\mapsto\;\; v(\theta,I)+z
\end{aligned}
\end{equation}
and hence 
\begin{equation}\label{v theta I vera}
\begin{aligned}
v (\theta, I) & = \sum_{j\in S} \begin{pmatrix}
\Omega(j)^{\frac12}  \Big( \alpha_j \cos(j \cdot x) + \beta_j \sin(j \cdot x) \Big) \\
\Omega(j)^{- \frac12} \Big( \alpha_j \sin(j \cdot x)  - \beta_j \cos(j \cdot x)  \Big)
\end{pmatrix} \\
& = \sum_{j\in S} \begin{pmatrix}
\Omega(j)^{\frac12} \sqrt{2(\zeta_{j}+I_{j})} \Big( \cos(\theta_j) \cos(j \cdot x) + \sin(\theta_j) \sin(j \cdot x) \Big) \\
\Omega(j)^{- \frac12} \sqrt{2(\zeta_{j}+I_{j})} \Big( \cos(\theta_j)\sin(j \cdot x)  - \sin(\theta_j)\cos(j \cdot x)  \Big)
\end{pmatrix} \\ 
& = \sum_{j\in S} \begin{pmatrix}
\Omega(j)^{\frac12} \sqrt{2(\zeta_{j}+I_{j})} \cos(\theta_j - j \cdot x)   \\
- \Omega(j)^{- \frac12} \sqrt{2(\zeta_{j}+I_{j})} \sin(\theta_j - j \cdot x) 
\end{pmatrix}
\end{aligned}
\end{equation}
which put in action angle variables the tangential sites
while leave unchanged the normal component $z$.
The symplectic $2$-form in \eqref{symReal} reads
\begin{equation}\label{symformAction}
\mathcal{W}:=\Big(\sum_{j\in S}d\theta_{j}\wedge dI_{j}\Big)\oplus\mathcal{W}_{|H_{S}^{\perp}}
=d\Lambda
\end{equation}
where $\Lambda$ is the Liouville $1$-form
\begin{equation}\label{liouville}
\Lambda_{\theta,I,z}[\widehat{\theta}, \widehat{I},\widehat{z}]:=
-\sum_{j\in S}I_{j}\widehat{\theta}_{j}-\frac{1}{2}(Jz,\widehat{z})_{L^{2}}\,.
\end{equation}
Given a Hamiltonian $K : \T^{\nu}\times\R^{\nu}\times H_{S}^{\perp}\to \R$,
the associated vector field, w.r.t. the symplectic form in \eqref{symformAction}
is
\[
X_{K}:=\big(\pa_{I}K, -\pa_{\theta}K, \mathbb{P}_{S}^{\perp} J\nabla_{z}K\big)\,.
\]

\vspace{0.9em}
\noindent
{\bf Tangential and normal subspaces in Complex variables.}
By means of the transformation ${\mathcal C}$ in \eqref{CVWW}, 
each two-dimensional subspace 
$V_j$, $j \in \Gamma^*\setminus\{0\}$ (see \eqref{definition Vj}) is isomorphic to 
\begin{equation}\label{def bf Hj}
{\bf H}_j := \Big\{ \begin{pmatrix}
z_j e^{\ii j \cdot x} \\
\overline z_j e^{- \ii j \cdot x}
\end{pmatrix} : z_j \in \C\Big\} \,.
\end{equation}
We denote by ${\bf \Pi}_j$ the orthogonal projectors on ${\bf H}_j$ and 
one easily sees that $\mathbb{P}_j = {\mathcal C}^{- 1} {\bf \Pi}_j {\mathcal C}$. 
We define 
\begin{equation}\label{decomp siti tangenziali coordinate complesse}
\begin{aligned}
{\bf H}_S & := \oplus_{j \in S} {\bf H}_j := \Big\{ \begin{pmatrix}
u \\
\overline u 
\end{pmatrix} : u(x) = \sum_{j \in S} \widehat u(j) e^{\ii j \cdot x}\Big\} \\
{\bf H}_S^\bot & := \oplus_{j \in S_0^{c}} {\bf H}_j = \Big\{ \begin{pmatrix}
u \\
\overline u 
\end{pmatrix} : u(x) = \sum_{j \in S_0^{c}} \widehat u(j) e^{\ii j \cdot x}\Big\}
\end{aligned}
\end{equation}
and we denote by ${\bf \Pi}_S, {\bf \Pi}_S^\bot$ the orthogonal projectors on ${\bf H}_S, {\bf H}_S^\bot$. One clearly has 
\begin{equation}\label{decomp siti tangenziali coordinate complesse2}
\mathbb{P}_S = {\mathcal C}^{- 1} {\bf \Pi}_S {\mathcal C}\,, \qquad  
\mathbb{P}_S^\bot = {\mathcal C}^{- 1} {\bf \Pi}_S^\bot {\mathcal C}\,. 
\end{equation}
Written in a more explicit way 
\begin{equation*}
\begin{aligned}
& {\bf \Pi}_S \begin{pmatrix}u (x) \\
 \overline u(x) 
 \end{pmatrix} = \begin{pmatrix}
\sum_{j \in S} \widehat u(j) e^{\ii j \cdot x} \\
\sum_{j \in S} \overline{\widehat u(j)} e^{- \ii j \cdot x}
\end{pmatrix}\,, 
\qquad
 {\bf \Pi}_S^\bot \begin{pmatrix}u (x) \\
 \overline u(x) 
 \end{pmatrix} = \begin{pmatrix}
\sum_{j \in S^c_0} \widehat u(j) e^{\ii j \cdot x} \\
\sum_{j \in S^c_0} \overline{\widehat u(j)} e^{- \ii j \cdot x}
\end{pmatrix}\,.
\end{aligned} 
\end{equation*}
We also define the projector $\Pi_{S},\Pi_{S}^{\perp}$ as 
\begin{equation}\label{projComplex}
(\Pi_{S}u)(x)=\sum_{j\in S}\widehat{u}(j)e^{\ii j\cdot x}\,,\qquad 
(\Pi_{S}^{\perp}u)(x)=\sum_{j\in S^c_0}\widehat{u}(j)e^{\ii j\cdot x}\,,
\end{equation}
so that ${\bf \Pi}_{S}(u,\bar{u})=(\Pi_{S}u, \Pi_{-S}\bar{u})=\big(\Pi_{S}u, \overline{\Pi_{S}u}\,\big)$.

\subsection{Spaces of Functions}\label{subsec:function spaces}

Let $\nu\geq1$. Along the paper we shall deal with families
\[
\mathbb{T}^{\nu}\ni\vphi\mapsto u(\vphi)\in H^{s}(\mathbb{T}_{\Gamma}^{d};\C)\,,
\]
with a certain regularity in the variables $\vphi\in \mathbb{T}^{\nu}$. 
Given $u\in L^{2}(\mathbb{T}^{\nu}; H^{s}(\T_{\Gamma}^{d};\C))$
we can expand
\begin{equation*}
u(\vphi,x)=
\sum_{\ell\in\Z^{\nu}}u(\ell;x)e^{\ii \ell\cdot\vphi}=
\sum_{\ell\in \Z^{\nu}, j\in \Gamma^{*}}\widehat{u}(\ell,j)e^{\ii (\ell\cdot\vphi+j\cdot x)}\,,
\qquad u(\ell;x)=\frac{1}{(2\pi)^\nu}
(u(\cdot,x), e^{\ii\ell\cdot\vphi})_{L^{2}}\,,
\end{equation*}
where 
$(\cdot,\cdot)_{L^{2}}$ is the standard $L^{2}$-Hermitian 
product of the straight torus. We introduce the scale of Hilbert spaces
\begin{equation}\label{SobolevSpazioAngoli}
\begin{aligned}
H^{s}&:=H^{s}(\mathbb{T}^{\nu}\times \mathbb{T}_{\Gamma}^{d};\C)
\\&:=
\Big\{ 
u(\vphi,x)\in L^{2}(\mathbb{T}^{\nu}\times \mathbb{T}_{\Gamma}^{d};\C)\;:\;
\|u\|_{s}^{2}:=\sum_{\ell\in \Z^{\nu}, j\in \Gamma^{*}}
\langle \ell,j\rangle^{2s}
|\widehat{u}(\ell,j)|^{2}
\Big\}\,,
\end{aligned}
\end{equation}
where $\langle \ell,j\rangle:=\sqrt{1+|\ell|^{2}+|j|^{2}}$.
To lighten the notation we shall write $\T_{*}^{\nu+d}$ instead of 
$\mathbb{T}^{\nu}\times \mathbb{T}_{\Gamma}^{d}$.

\noindent
The Sobolev norm $ \| \cdot \|_s $ defined in \eqref{SobolevSpazioAngoli} 
is equivalent to 
\begin{equation}\label{sobnormseparati}
\| u \|_s \simeq \| u \|_{H^s_\vphi L^2_x} + \| u \|_{L^2_\vphi H^s_x} \, . 
\end{equation}
\begin{rmk}
It is easy to check that 
\[
H^{s}(\mathbb{T}^{\nu}\times \mathbb{T}_{\Gamma}^{d};\C)\sim
\bigcap_{p_1+p_2=s}
H^{p_1}(\mathbb{T}^{\nu}; H^{p_2}(\mathbb{T}_{\Gamma}^{d};\C))\,.
\]
\end{rmk}
\noindent
In the following we shall consider the scale 
of Sobolev spaces $(H^s)_{s=0}^{\bar{s}}$ where
$$
\bar{s}\gg s_0 > \frac{\nu + d}{2}  \,.
$$
We remark that $H^{s_0}$ is continuously 
embedded in $L^\infty(\T^{\nu+d}_*, \C)$.

\subsubsection{Traveling functions}\label{sec:travelfunction}
Since the water waves Hamiltonian $H$ in \eqref{HamiltonianWW}
Poisson commutes with the momenta in \eqref{Hammomento} 
we have that the system \eqref{eq:113} is invariant under space translations.
Indeed, by defining the translation operator
\begin{equation}\label{def:vec.tau}
\tau_\vs : h(x) \mapsto h(x + \vs)\,,\;\;\vs\in \R^{d}\,,
\end{equation}
the Hamiltonian \eqref{HamiltonianWW} satisfies 
$H\circ\tau_{\vs}=H$ for any $\vs\in \R^{d}$.
In this paper we look for solutions which are traveling according 
to Definition \ref{def:quasitravelling}.
Recalling the $\nu\times d$-matrix $\tV:=(\tv_1,\ldots,\tv_{\nu})^{T}$,
with $\tv_1, \dots, \tv_\nu \in\Gamma^{*}$,
this is equivalent to look for embeddings
$u\in L^{2}(\T_{*}^{\nu+d};\R^{2})$
satisfying 
\begin{equation}\label{condtraembedd}
u(\vphi-\mathtt{V}\vs,x)=
( \tau_\vs \circ u)(\vphi,x)=u(\vphi,x+\vs)\,,  
\quad \forall \vphi \in \T^\nu\,,\; \vs\in \R^{d}\,, x\in\R^{d}\,.
\end{equation}
We define the subspace
\begin{equation}\label{SV}
S_{\tV}:=\{u(\varphi, x)\in L^2(\T_{*}^{\nu+d}, \C^2) : u(\varphi, x)=U(\varphi- \tV x)\,,\; U(\Theta)\in L^2(\T^\nu, \C^2)\}\,.
\end{equation}
With abuse of notation, we denote by $S_{\tV}$ also the subspace of $L^2(\T_{*}^{\nu+d}, \C)$ 
of scalar functions $u(\varphi, x)$ of the form 
$u(\varphi, x)=U(\varphi- \tV x),\, U(\Theta)\in L^2(\T^\nu, \C)$). 
We can easily note that
\begin{equation}\label{travelequi}
\eqref{condtraembedd} \qquad \Leftrightarrow \qquad u\in S_{\tV}
\qquad  \Leftrightarrow\qquad
\tV^{T}\pa_{\vphi}u+\nabla u=0\,,\;\;\forall(\vphi,x)\in \T_{*}^{\nu+d}\,.
\end{equation}

\smallskip
\noindent
{\bf Quasi-periodic traveling waves in action-angle-normal coordinates.}
We now discuss how the momentum preserving 
condition reads in the \emph{action-angles} coordinates $(\theta,I,z)$ introduced in 
\eqref{actionanglevar}.
More precisely,
we give the following Definition.
\begin{defn}{\bf (Traveling wave embedding).}\label{def:travelembedd}
Recalling \eqref{decomposition}, 
consider a torus embedding 
\begin{equation}\label{torusprototipo}
\begin{aligned}
i \,:\, \mathbb{T}^{\nu} \,&\to \,\mathbb{T}^{\nu}\times\mathbb{R}^{\nu}\times H_{S}^{\perp}\,,
\qquad \vphi\mapsto i(\vphi)=(\theta(\vphi),I(\vphi),z(\vphi,x)) 
\\
&\mathcal{I}(\vphi):=i(\vphi)-(\vphi,0,0)=((\Theta(\vphi),I(\vphi),z(\vphi,x)) )\,.
\end{aligned}
\end{equation}
We say that $\mathcal{I}(\vphi)$ (or $i(\vphi)$) is a \emph{traveling wave embedding}
if satisfies the conditions
\begin{equation}\label{cond:travelemb}
\begin{aligned}
&\Theta(\vphi)=\Theta(\vphi-\mathtt{V}\vs)\,,\quad I(\vphi)=I(\vphi-\mathtt{V}\vs)\,,
\\&
z(\vphi,x+\vs)=z(\vphi-\tV \vs,x)\,,\qquad \qquad
\forall(\vphi,x)\in \T_{*}^{\nu+d}\,,\vs\in\R^{d}\,.
\end{aligned}
\end{equation}
\end{defn}
We have the following simple result.
\begin{lemma}
Consider a torus embedding as in \eqref{torusprototipo}.
The following holds.

\noindent
$(i)$ The condition \eqref{cond:travelemb} holds if and only if
the embedding $i(\vphi)$ satisfies
\begin{equation}\label{eq:travelemb}
\tV^{T}\pa_{\vphi}\Theta(\vphi)=0\,,
\quad \tV^{T}\pa_{\vphi}I(\vphi)=0\,,
\quad
(\tV^{T}\pa_{\vphi}+\nabla)z(\vphi,x)=0\,.
\end{equation}
\noindent
$(ii)$ Define the function $u(\vphi,x):=(Ai)(\vphi)$ 
where $A$ is the map in \eqref{actionanglevar}.
Then $u$ is a traveling wave according to Def. \ref{def:quasitravelling} if and only if
$i(\vphi)$ is a traveling wave embedding according to Def. \ref{def:travelembedd}.
\end{lemma}
\begin{proof}
Item $(i)$ follows by a straightforward computation using \eqref{cond:travelemb}.
Now we note that 
conditions \eqref{eq:travelemb} and the definition of the map $A$
imply that the function 
$u$ satisfies the transports equations $(\tV^{T}\pa_{\vphi}+\nabla)u=0$.
Therefore item $(ii)$ follows using \eqref{travelequi}.
\end{proof}

\subsubsection{Whitney-Sobolev functions}\label{sec:Whitney-Sobolev}
In the paper we will use Sobolev norms for real or complex functions 
$u(\omega, \th, \vphi, x)$, $(\vphi,x) \in \T^\nu \times \T^d_\Gamma$, 
depending on parameters $(\omega,\th) \in F$ 
in a Lipschitz way together with their derivatives in the sense of Whitney, 
where $F$ is a closed subset of $\R^{\nu+1}$. 
We use the compact notation $\lambda := (\omega,\th)$ 
to collect the frequency $\om$ and the depth $\th$ into a parameter vector. 
Also recall that $\| \cdot \|_s$ denotes the  norm 
of the Sobolev space $H^s(\T^{\nu+d}_*, \C) = H^s_{(\vphi,x)}$
introduced in \eqref{SobolevSpazioAngoli}. 
We now define  the ``Whitney-Sobolev'' norm 
$\| \cdot \|_{s,F}^{k+1,\g}$.
\begin{defn} {\bf (Whitney-Sobolev functions)}\label{def:Lip F uniform}
Let $F$ be a closed subset of $\R^{\nu+1}$.  
Let $k \geq 0$ be an integer,  $\g \in (0,1]$, and $s \geq 0$. 
We say that a function $u : F \to H^s_{(\vphi,x)}$ 
belongs to $\Lip(k+1,F,s,\g)$ if there exist functions 
\[
u^{(j)} : F \to H^s_{(\vphi,x)}\,, \quad j \in \N^{\nu+1}\,, \ 0 \leq |j| \leq k\,, 
\]
with $u^{(0)} = u$, and a constant $M > 0$ 
such that, if 
$R_j(\lambda,\lambda_0) := R_j^{(u)}(\lambda,\lambda_0) $ 
is defined by 
\begin{equation*} %\label{16 Stein uniform}
u^{(j)}(\lambda) 
= \sum_{\ell \in \N^{\nu+1} : 
|j+\ell| \leq k} \frac{1}{\ell!} \, u^{(j+\ell)}(\lambda_0) \, (\lambda - \lambda_0)^\ell 
+ R_j(\lambda,\lambda_0), \quad \lambda, \lambda_0 \in F, 
\end{equation*} 
then 
\begin{equation} \label{17 Stein uniform}
\g^{|j|} \| u^{(j)}(\lambda) \|_s \leq M, \quad 
\g^{k+1} \| R_j(\lambda, \lambda_0) \|_s \leq M |\lambda - \lambda_0|^{k+1 - |j|} \quad 
\forall \lambda, \lambda_0 \in F, \ |j| \leq k.
\end{equation} 
An element of $\Lip(k+1,F,s,\g)$ 
is in fact the collection $\{ u^{(j)} : |j| \leq k \}$. 
The norm of $u \in \Lip(k+1,F,s,\g)$ is defined as
\begin{equation} \label{def norm Lip Stein uniform}
\| u \|_{s,F}^\kug 
:= \| u \|_s^\kug
:= \inf \{ M > 0 : \text{\eqref{17 Stein uniform} holds} \}.
\end{equation}
If $ F = \R^{\nu+1} $ by $ \Lip(k+1,\R^{\nu+1},s, \g) $  
we shall mean the 
space of the functions $ u = u^{(0)} $ for which
there exist $ u^{(j)} = \pa_{\lambda}^j u $, $ |j| \leq k $,  
satisfying  \eqref{17 Stein uniform}, 
with the same norm \eqref{def norm Lip Stein uniform}.
\end{defn}
We make some remarks.
\begin{enumerate}
\item 
If $ F = \R^{\nu+1} $, and $ u \in \Lip(k+1,F,s,\g) $
the $ u^{(j)} $, $ | j | \geq 1 $, are uniquely determined as the partial derivatives 
$ u^{(j)} = \pa_{\lambda}^j u $, $ | j | \leq k $,  of $ u = u^{(0)} $. 
Moreover all the derivatives $ \pa_{\lambda}^j u $, $ | j | = k $
are Lipschitz.  Since $ H^s $ is a Hilbert space we have that % the space of functions 
$ \Lip(k+1,\R^{\nu+1},s,\g) $ coincides with the Sobolev space
$ W^{k+1, \infty} (\R^{\nu + 1}, H^s) $. 
\item
The Whitney-Sobolev norm  of $ u $ in \eqref{def norm Lip Stein uniform} is 
equivalently given by 
\begin{equation*} %\label{def -equiv}
\| u \|_{s,F}^\kug := \| u \|_s^\kug
= \max_{|j| \leq k} 
\Big\{ \g^{|j|} \sup_{\lambda \in F}  \| u^{(j)}(\lambda) \|_s, \g^{k+1} 
\sup_{\lambda \neq \lambda_0}  
\frac{\| R_j(\lambda, \lambda_0) \|_s}{|\lambda - \lambda_0|^{k+1 - |j|}}  \Big\} \, .
\end{equation*}
\item
The exponent of $\g$ in \eqref{17 Stein uniform} 
gives the number of ``derivatives'' of $u$ that are involved 
in the Taylor expansion (taking into account that 
in the remainder there is one derivative more than in the Taylor polynomial);
on the other hand, the exponent of $|\lambda-\lambda_0|$ 
gives the order of the Taylor expansion of $u^{(j)}$ with respect to $\lambda$. 
This is the reason for the difference of $|j|$ between the two exponents. 
The factor $\g$ is normalized by the rescaling 
%$B.7$ in \cite{BBHM}.
\eqref{resca}.
\end{enumerate}
Theorem \ref{thm:WET} and \eqref{Wg} 
%Theorem $B.2$ and $(B.10)$ in \cite{BBHM}
provide an extension operator which associates 
to an element $u \in \Lip(k+1,F,s,\g)$ 
an extension $\tilde u \in \Lip(k+1,\R^{\nu+1},s,\g) $. 

\noindent
As already observed, %the space 
$\Lip(k+1,\R^{\nu+1},s,\g) $ 
coincides with $ W^{k+1, \infty} (\R^{\nu + 1}, H^s)$, 
with equivalence of the norms (see 
%$(B.9)$ in \cite{BBHM})
\eqref{0203.1})
\[
\| u \|_{s, F}^\kug  \sim_{\nu, k} 
\| \tilde u \|_{W^{k+1,\infty, \gamma}(\R^{\nu+1}, H^s)} := 
\sum_{|\alpha| \leq k+1} \g^{|\alpha|} 
\| \pa_\lambda^\alpha \tilde u \|_{L^\infty(\R^{\nu+1}, H^s)}  \,.
\]
By Lemma \ref{lemma:2702.1}, 
%By Lemma $B.3$ in \cite{BBHM}
the extension $\tilde u$ is independent of the Sobolev space $H^s$.
We can identify any element $u \in \Lip(k+1, F, s, \g)$
(which is a collection $u = \{ u^{(j)} : |j| \leq k\}$) 
with the equivalence class of functions 
$f \in W^{k+1, \infty}(\R^{\nu+1}, H^s) / \! \! \sim$  
with respect to the equivalence relation 
$f \sim g$ when $\pa_\lambda^j f(\lambda) = \pa_\lambda^j g(\lambda)$
for all $\lambda \in F$, for all $|j| \leq k+1$.

\noindent 
For a scalar function $f: \R^{\nu+1}\to \R$ independent of $(\varphi, x)$ the norm $\| \cdot \|_{s}^{k, \gamma}$ are independent of $s$ and we denote it simply $| \cdot |^{k, \gamma}$. 

\noindent
For any $N>0$, we introduce the smoothing operators
\begin{equation*}%\label{def:smoothings}
(\Pi_N u)(\vphi,x) := 
\sum_{\la \ell,j \ra \leq N} u_{\ell j} e^{\ii (\ell\cdot\vphi + j\cdot x)} \qquad
\Pi^\perp_N := {\rm Id} - \Pi_N\,.
\end{equation*}
We now collect some results proved in \cite{BBHM}. 
For completeness we refer the reader to
Appendix \ref{sec:U}.
\begin{lemma} {\bf (Smoothing).} \label{lemma:smoothing}
Consider the space $\Lip(k+1,F,s,\g)$ 
defined in Definition \ref{def:Lip F uniform}. 
The smoothing operators $\Pi_N, \Pi_N^\perp$ 
satisfy the estimates
\begin{align}
\| \Pi_N u \|_{s}^\kug 
& \leq N^\alpha \| u \|_{s-\alpha}^\kug\, , \quad 0 \leq \alpha \leq s\,, 
\label{p2-proi} \\
\| \Pi_N^\perp u \|_{s}^\kug 
& \leq N^{-\alpha} \| u \|_{s + \alpha}^\kug\, , \quad  \alpha \geq 0\,.
\label{p3-proi}
\end{align}
\end{lemma}

%\begin{proof}
%See \cite{BBHM}.
%%Appendix \ref{sec:U}.
%\end{proof}
\begin{lemma} {\bf (Interpolation)}\label{lemma:interpolation}
Consider the space $\Lip(k+1,F,s,\g)$ defined in Definition \ref{def:Lip F uniform}. 

\noindent
(i) Let $s_1<s_2$. Then for any $\theta\in(0,1)$ one has
\begin{equation}\label{2202.3}
\| u \|_s^{k+1, \gamma} \leq (\| u \|_{s_1}^{k+1, \gamma})^\theta 
(\| u \|_{s_2}^{k+1, \gamma})^{1 - \theta}\,, \quad s:=\theta s_1 + (1-\theta) s_2\,. 
\end{equation}

\noindent
(ii) Let $a_0, b_0 \geq0$ and $p,q>0$. For all $\epsilon>0$, there exists a constant $C(\epsilon):= C(\epsilon,p,q)>0$, which satisfies $C(1)<1$, such that
\begin{equation}\label{2202.2}
\| u \|^{k+1,\g}_{a_0 + p} \| v \|^{k+1,\g}_{b_0 + q} 
\leq \epsilon \| u \|^{k+1,\g}_{a_0 + p +q} \| v \|^{k+1,\g}_{b_0} 
+ C(\epsilon) \| u \|^{k+1,\g}_{a_0} \| v \|^{k+1,\g}_{b_0 + p +q}\ .
\end{equation}
\end{lemma}

%\begin{proof}
%See Appendix \ref{sec:U}.
%\end{proof}
\begin{lemma}{\bf (Product and composition).}
\label{lemma:LS norms}
Consider the space $\Lip(k+1,F,s,\g)$ defined in Definition \ref{def:Lip F uniform}.
For all $ s \geq s_0  := (\nu + d)/2 + 1$, we have
\begin{align}
\| uv \|_{s}^\kug
& \leq C(s, k) \| u \|_{s}^\kug \| v \|_{s_0}^\kug 
+ C(s_0, k) \| u \|_{s_0}^\kug \| v \|_{s}^\kug\,. 
\label{p1-pr}
\end{align}
Let $ \| \beta \|_{2s_0+1}^\kug \leq \delta (s_0, k) $ small enough. 
Then the composition operator 
\[
\cB : u \mapsto \cB u, \quad  
(\cB u)(\vphi,x) := u(\vphi, x + \beta (\vphi,x)) \, , 
\]
satisfies the following tame estimates: for all $ s \geq s_0$, 
\begin{equation}\label{pr-comp1}
\| \cB u \|_{s}^\kug \lesssim_{s, k} \| u \|_{s+k+1}^\kug 
+ \| \beta \|_{s}^\kug \| u \|_{s_0+k+2}^\kug \, .
\end{equation}
Let $ \| \beta \|_{2s_0+k+2}^\kug \leq \delta (s_0, k) $ small enough.  The function $ \breve \beta $ defined by 
the inverse diffeomorphism 
$ y = x + \beta (\vphi, x) $ if and only if $ x = y + \breve \beta ( \vphi, y ) $,  
satisfies 
\begin{equation}\label{p1-diffeo-inv}
\| \breve \beta \|_{s}^\kug \lesssim_{s, k}  \| \beta \|_{s+k+1}^\kug \, . 
\end{equation}
\end{lemma}
%\begin{proof} See Appendix \ref{sec:U}. % and Lemma 2.30 in \cite{BertiMontalto}.  
%\end{proof}
If  $ \omega$ belongs to the set of 
Diophantine  vectors $ \mathtt{DC}(\gamma, \tau) $
defined in \eqref{def:DCgt}, 
%where
%\begin{equation}\label{DC tau0 gamma0}
%\mathtt{DC}(\gamma, \tau) := \Big\{ \omega \in \mathtt \R^\nu :
% |\omega \cdot \ell| \geq \frac{\gamma}{|\ell|^{\tau}} \ \  
% \forall \ell \in \Z^\nu \setminus  \{ 0 \} \Big\} \, ,
%\end{equation}
the equation $\ompaph v = u$, where $u(\vphi,x)$ 
has zero average with respect to $ \vphi $, 
has the periodic solution 
\begin{equation}\label{def:ompaph}
(\omega\cdot \pa_\vphi )^{-1} u := 
\sum_{\ell \in \Z^{\nu} \setminus \{0\}, j \in \Gamma^{*}} 
\frac{ u_{\ell, j} }{\ii \omega\cdot \ell }e^{\ii (\ell \cdot \vphi + j x )} \,.
\end{equation}
For all $\omega\in \R^\nu$ we define its extension
\begin{equation} \label{def ompaph-1 ext}
(\ompaph)^{-1}_{ext} u(\vphi,x) := 
\sum_{(\ell, j) \in \Z^{\nu}\times\Gamma^*} 
\frac{\chi( \omega\cdot \ell \g^{-1} \langle \ell \rangle^{\t}) }{\ii \omega\cdot \ell}\,
u_{\ell,j} \, e^{\ii (\ell \cdot \vphi + jx)}\,,
\end{equation}
where $\chi \in \cC^\infty(\R,\R)$ is an even and positive cut-off function such that 
\begin{equation}\label{cut off simboli 1}
\chi(\xi) = \begin{cases}
0 & \quad \text{if } \quad |\xi| \leq \frac13 \\
1 & \quad \text{if} \quad  \ |\xi| \geq \frac23\,,
\end{cases} \qquad 
\partial_\xi \chi(\xi) > 0 \quad \forall \xi \in \Big(\frac13, \frac23 \Big) \, . 
\end{equation}
Note that $(\ompaph)^{-1}_{ext} u = (\ompaph)^{-1} u$
for all $\omega\in \mathtt{DC}(\gamma, \tau)$.

\begin{lemma} { \bf (Diophantine equation)}
\label{lemma:WD}
Let $\mathtt \Omega \subset \R^\nu$ be a bounded open set. For all $u \in W^{k+1,\infty,\g}(\mathtt \Omega \times [\mathtt h_1, \mathtt h_2], H^{s+\mu})$, we have that $(\omega\cdot \pa_\vphi )^{-1}_{ext} u \in W^{k+1,\infty,\g}(\mathtt \Omega \times [\mathtt h_1, \mathtt h_2], H^{s})$ and 
\begin{equation} \label{2802.2}
\| (\omega\cdot \pa_\vphi )^{-1}_{ext} u \|_{s}^\kug
\leq C(k) \g^{-1} \| u \|_{s+\mu}^\kug, 
\qquad \mu := k+1 +  \t(k+2)\,. 
\end{equation}
Moreover, for $F \subseteq \mathtt{DC}(\g,\t) \times [\mathtt h_1, \mathtt h_2]$, one has 
\begin{equation} \label{Diophantine-1}
\| (\omega\cdot \pa_\vphi )^{-1} u \|_{s,F}^\kug
\leq C(k) \g^{-1}\| u \|_{s+\mu, F}^\kug \,.  
\end{equation}
\end{lemma}
%\red{aggiungere qui lo statement del lemma B.4 in \cite{BBHM}.}
%\begin{proof} See Appendix \ref{sec:U}. 
%\end{proof}
Let $\mathtt{m}\in \R$ be a constant 
of the form $\mathtt{m}=\mathtt{m}(\lambda,i(\lambda))$, $\lambda\in \R^{\nu+1}$ and 
$i=i(\lambda)$ is a smooth traveling wave embedding 
$i$ as in Definition \ref{def:travelembedd} with $\lambda = (\omega, \mathtt h) \in \mathtt \Omega \times [\mathtt h_1, \mathtt h_2]$ where $\mathtt \Omega \subset \R^\nu$ is a bounded open set. 
Assume also that 
 $|\mathtt{m}|^{k_0,\gamma}\lesssim\e\ll1$ and that $\mathtt{m}$ depends on $i$ in a 
 Lipschitz way.  More precisely, we assume that 
 \begin{equation}\label{stimaemminodelta12}
 |\mathtt{m}(i_1)-\mathtt{m}({i_2})|\lesssim \e \|i_{1}-i_{2}\|_{s_0+\s}\,,
 \end{equation}
for some $\s>0$\,. 
Let us define the set
	\begin{equation}\label{tDtCzero}
			\begin{aligned}
				\tT\tC(\gamma,\tau) 
				& :=\tT\tC(\gamma,\tau; i)\\
				&:= \Big\{ (\omega,\th)\in\mathtt \Omega \times[\th_1,\th_2] \,:\, 
				|(\omega-\tV\mathtt{m})\cdot\ell| 
				\geq  \mathtt{c}_0\gamma\jap{\ell}^{-\tau} \,\ \forall\,\ell\neq 0 \Big\}\,,
			\end{aligned}
		\end{equation}
		for some $\mathtt{c}_0>1$.
For $\lambda=(\omega,\th)\in \tT\tC(\gamma,\tau) $ we
consider the operator $((\omega-\mathtt{V}\mathtt{m})\cdot\pa_{\vphi})^{-1}$
acting on traveling periodic functions as
\begin{equation}\label{def:ompaph2}
((\omega-\mathtt{V}\mathtt{m})\cdot\pa_{\vphi})^{-1} u := 
\sum_{\ell \in \Z^{\nu} \setminus \{0\}} 
\frac{ u_{\ell} }{\ii (\omega-\mathtt{V}\mathtt{m})\cdot \ell }e^{\ii \ell \cdot \vphi} \,.
\end{equation}
For all $\lambda\in \mathtt \Omega \times [\mathtt h_1, \mathtt h_2]$, we define its extension
\begin{equation} \label{def ompaph-1 ext2}
((\omega-\mathtt{V}\mathtt{m})\cdot\pa_{\vphi})^{-1}_{ext} u(\vphi) := 
\sum_{\ell \in \Z^{\nu}} 
\frac{\chi( (\omega-\mathtt{V}\mathtt{m})\cdot \ell \g^{-1} \langle \ell \rangle^{\t}) }{\ii 
(\omega-\mathtt{V}\mathtt{m})\cdot \ell}\,
u_{\ell} \, e^{\ii \ell \cdot \vphi }\,,
\end{equation}
where $\chi$ is as in \eqref{cut off simboli 1}.
We have the following.

\begin{lemma} { \bf (Diophantine equation 2)}
\label{lemma:WD2}
For all $u \in W^{k+1,\infty,\g}(\mathtt \Omega \times [\mathtt h_1, \mathtt h_2], H^{s+\mu}(\T^{\nu}))$, 
if $\e \gamma^{-k-1} \ll1$, 
we have that $((\omega-\mathtt{V}\mathtt{m})\cdot \pa_\vphi )^{-1}_{ext} u \in  W^{k+1,\infty,\g}(\mathtt \Omega \times [\mathtt h_1, \mathtt h_2], H^{s}(\T^{\nu}))$ and 
\begin{equation} \label{2802.2nuova}
\| ((\omega-\mathtt{V}\mathtt{m})\cdot \pa_\vphi )^{-1}_{ext} u \|_{s}^\kug
\leq C(k) \g^{-1} \| u \|_{s+\mu}^\kug, 
\qquad \mu := k+1 +  \t(k+2)\,. 
\end{equation}
The operator $((\omega-\mathtt{V}\mathtt{m})\cdot \pa_\vphi )^{-1}_{ext}$ concides 
with $((\omega-\mathtt{V}\mathtt{m})\cdot \pa_\vphi )^{-1}$
one the set $\tT\tC(\gamma,\tau) $.
Moreover we also have 
\begin{equation}\label{2802.2nuova2}
\begin{aligned}
\| ((\omega-\mathtt{V}\mathtt{m}(i_1))\cdot \pa_\vphi )^{-1}_{ext} u
&-
 ((\omega-\mathtt{V}\mathtt{m}(i_2))\cdot \pa_\vphi )^{-1}_{ext} u
 \|_{s,\R^{\nu+1}}
 \\&\lesssim\e\gamma^{-2}\|i_1-i_2\|_{s_0+\s}\|u\|_{s+\mu}\,,\qquad \mu=2\tau+1\,.
 \end{aligned}
\end{equation}
\end{lemma}

We finally state a standard Moser tame estimate for the nonlinear 
composition operator
$$
u(\vphi, x) \mapsto {\mathtt f}(u)(\vphi, x) := f(\vphi, x, u(\vphi, x)) \, . 
$$ 
Since the variables $ (\vphi, x) := y $ have the same role, 
we state it for a  generic Sobolev space  $ H^s (\T^d ) $.  
\begin{lemma}{\bf (Composition operator)} \label{Moser norme pesate}
Let $ f \in \cC^{\infty}( \T^\nu \times \T^d_\Gamma \times \R, \C )$ and $C_0 > 0$. 
Consider the space $\Lip(k+1,F,s,\g)$ given in Definition \ref{def:Lip F uniform}. 
If $u(\lambda) \in H^s(\T^\nu \times \T^d_\Gamma, \R)$, $\lambda \in F$ is a family of Sobolev functions
satisfying $\| u \|_{s_0,F}^{k+1, \gamma} \leq C_0$, then, for all $ s \geq s_0 := (d + \nu)/2 + 1$, 
\begin{equation} \label{0811.10}
\| {\mathtt f}(u) \|_{s, F}^{k+1, \gamma} 
\leq C(s, k, f, C_0 ) ( 1 + \| u \|_{s, F}^{k+1, \gamma}) \,.
\end{equation} 
The constant $C(s,k,f, C_0)$ depends on $s$, $k$ 
and linearly on $\| f \|_{\cC^m}$, 
where $m$ is an integer larger than $s+k+1$, 
and $B \subset \R$ is a bounded interval such that 
$u(\lambda, \vphi, x) \in B$ for all $\lambda \in F$, $(\vphi, x) \in  \T^\nu \times \T^d_\Gamma$,
for all $\| u \|_{s_0,F}^{k+1, \gamma} \leq C_0$.
\end{lemma}
%\begin{proof} See Appendix \ref{sec:U}. 
%\end{proof}

\subsubsection{Space of Functions on a Strip}\label{sec:funzionistrip}
In this subsection 
we introduce 
appropriate spaces (and norms) of functions defined on the strip. 
These spaces will be used in order to study the pseudo-differential structure of the Dirichlet-Neumann
operator $G(\eta)\psi$ in \eqref{eq:112a} in Section \ref{sec:DNsezione}.

\smallskip
We follow the notations in \cite{BeMasVent1}.
We define for any $s \in \N$, the space\footnote{
By expanding $u(y, x) = \sum_{\begin{subarray}{c}
\ell \in \Z^\nu \\
\xi \in \Gamma^*
\end{subarray}} \widehat u(y,  \xi) e^{\ii  x \cdot \xi}$ in Fourier series w.r.t. $x$, the latter norm can be also written as follows. For any $0 \leq k \leq s$
\begin{equation}\label{norma sobolev striscia seconda form}
\begin{aligned}
\| \partial_y^k u \|_{L^2([- 1, 0], H^{s - k})} & =  \Big( \int_{- 1}^0 \| \partial_y^k u(y, \cdot) \|_{H^{s - k}}^2\, d y \Big)^{\frac12}  
%\\& 
= \Big( \int_{- 1}^0 \sum_{ \xi} \langle  \xi \rangle^{2 (s - k)} |\partial_y^k \widehat u(y, \xi)|^2 \, d y \Big)^{\frac12} 
\\& 
= \Big( \sum_{\xi \in \Gamma^*} \langle  \xi \rangle^{2(s - k)} \int_{- 1}^0 |\partial_y^k \widehat u(y, \xi)|^2\, d y \Big)^{\frac12} 
%\\& 
= \Big( \sum_{\xi \in \Gamma^*} \langle  \xi \rangle^{2(s - k)}  \| \partial_y^k \widehat u(\cdot, \xi) \|_{L^2_y}^2  \Big)^{\frac12}\,.
\end{aligned}
\end{equation}
}

\begin{equation}\label{def sobolev striscia}
\begin{aligned}
 \cH^s  &:= \Big\{ u: [- 1, 0] \times \T^d_\Gamma \to \C :    \| u \|_{\cH^s}  < + \infty \Big\}\,,
\\
 \| u \|_{\cH^s} &:= \sum_{k = 0}^s \Big( \int_{- 1}^0 \| \partial_y^k u(y, \cdot) \|_{H^{s - k}_{x}}^2\, d y \Big)^{\frac12} 
= \sum_{k = 0}^s \| \partial_y^k u \|_{L^2_y H_{x}^{s - k}}\,, 
\\
 \| \partial_y^k u \|_{L^2_yH^{s - k}} &:= \| \partial_y^k u \|_{L^2([- 1, 0], H^{s - k})}\,. 
\end{aligned}
\end{equation}

\noindent
For any $s \geq 0$, we now define the space $\cO^s := H^s \big( \T^\nu \times \T^d_\Gamma \times [- 1, 0] ;\C\big)$
%$\cO^s$ as
%\begin{equation}\label{def cal Os}
%\cO^s := H^s \Big( \T^\nu \times \T^d_\Gamma \times [- 1, 0] \Big)\,.
%\end{equation} 
of the Sobolev functions $u : \T^\nu \times \T^d_\Gamma \times [- 1, 0] \to \C$.

\begin{lemma}\label{algebra striscia}
(i) Let $s_0 > \frac{\nu + d + 1}{2}$, $s \geq s_0$ and $u, v \in \cH^s$, resp. $u, v \in \cO^s$. 
Then $u v \in \cH^s$, resp. $u v \in \cO^s$  and
\begin{equation}\label{strippoperlalgebra}
\begin{aligned}
\| u v \|_{\cH^s} &\lesssim_s \| u \|_{\cH^s} \| v \|_{\cH^{s_0}} + \| u \|_{\cH^{s_0}} \| v \|_{\cH^s}\,,
\\
\| u v \|_{\cO^s} &\lesssim_s \| u \|_{\cO^s} \| v \|_{\cO^{s_0}} + \| u \|_{\cO^{s_0}} \| v \|_{\cO^s}\,.
\end{aligned}
\end{equation}
Moreover, for any $s \geq 0$, one has the equivalences
\begin{equation}\label{prop L2 vhi Hs Hs vphi cal 2 xy}
\begin{aligned}
\|u\|_{\cH^{s}}&\simeq_s\|u\|_{L^{2}_{y}H^{s}_{x}}+
\|u\|_{H^{s}_{y}L^{2}_{x}}
\qquad \quad 
\| u \|_{\cO^s} \simeq_s \| u \|_{L^2_\vphi \cH^s} + \| u \|_{H^s_\vphi \mathcal{H}^{0}}\,.
%\simeq\|u\|_{L^{2}_{\vphi}L^{2}_{y}H^{s}_{x}}+
%\|u\|_{L^{2}_{\vphi}H^{s}_{y}L^{2}_{x}}+\| u \|_{H^s_\vphi L^{2}_{y}L^{2}_x}\,. 
\end{aligned}
\end{equation}

$(ii)$ Let $s_1<s_2$, $\theta\in(0,1)$ and $u\in \mathcal{H}^{s_2}$ (resp. $u\in \mathcal{O}^{s_2}$). Then one has
\begin{equation}\label{strippointrippo}
\begin{aligned}
\|u\|_{\mathcal{H}^{s}}&\leq \|u\|_{\mathcal{H}^{s_1}}^{\theta}\|u\|_{\mathcal{H}^{s_2}}^{1-\theta}\,,
\qquad \quad 
\|u\|_{\mathcal{O}^{s}}\leq \|u\|_{\mathcal{O}^{s_1}}^{\theta}\|u\|_{\mathcal{O}^{s_2}}^{1-\theta}\,,
\end{aligned}
\end{equation}
for any $s=\theta s_1+(1-\theta)s_2$.
\end{lemma}

\begin{proof}
In the case of the half-line $y\in(-\infty,0]$ such properties follows 
 by interpolation properties of Sobolev spaces. For more details  one can see 
 \cite{BeMasVent1} and also Chapter $2$ in \cite{LionsMagenes1972}.
 In the case of the segment $[-1,0]$ one can reduce 
 to the previous case by providing an extension to the whole 
 half-line. In order to do this, following Theorem $2.1$ in Chapter $2$ (see in particular formula $(2.19)$)
 in \cite{LionsMagenes1972}, 
 one introduce $C^{\infty}([-1,0])$ funtions $\alpha$ and $\beta$ such that
 \[
 \alpha(y)+\beta(y)=1\,,\quad \forall\,y\in[-1,0]\,,
 \]
 and 
 $\alpha(y)\equiv0$ in a neighborhood of $0$ while  $\beta(y)\equiv0$ in a neighborhood of $-1$.
 Then one split
 \[
 u(y,\cdot)=\alpha(y)u(y,\cdot)+\beta(y)u(y,\cdot)\,,
 \]
and one considers the extensions $u_1$ defined on the half-line $[-1,+\infty)$ and $u_2$ defined on 
 the half-line $(-\infty,0]$ of the form
 \[
 \begin{aligned}
 u_1=\alpha(y)u\,,\;\;\; y\in[-1,0]\,,\qquad u_1:=0 \,,\,\;y>0\,,
 \\
  u_2=\beta(y)u\,,\;\;\; y\in[-1,0]\,,\qquad u_2:=0 \,,\,\;y<-1\,,
 \end{aligned}
 \]
 The desired properties for the function $u$ on $[-1,0]$ follow by the properties on the functions
 $u_1,u_2$ defined on the half-line.
%  one can see, for details, Theorem $2.1$ in Chapter $2$ (see in particular formula $(2.19)$)
% in \cite{LionsMagenes1972}.
% 
% We refer to \cite{LionsMagenes1972} for instance.
%We also refer the reader to \cite{BeMasVent1} for the case $y\in(-\infty,0]$.
\end{proof}
%\begin{rmk}\label{rmk:equiOHs}
%By standard properties one has that 
%\begin{equation}\label{prop L2 vhi Hs Hs vphi cal 2 xy}
%\| u \|_{\cO^s} \simeq \| u \|_{L^2_\vphi \cH^s} + \| u \|_{H^s_\vphi \mathcal{H}^{0}}
%\simeq\|u\|_{L^{2}_{\vphi}L^{2}_{y}H^{s}_{x}}+
%\|u\|_{L^{2}_{\vphi}H^{s}_{y}L^{2}_{x}}+\| u \|_{H^s_\vphi L^{2}_{y}L^{2}_x}\,. 
%\end{equation}
%\end{rmk}

\subsection{Linear operators}
We introduce general classes of linear operators. 
Given two Banach spaces $E,F$, we denote by 
$\mathcal{L}(E;F)$ the space of 
bounded linear operators from $E$ to $F$
endowed with the standard operator norm.
If $F=E$ we simply write $\cL(E)$.
We consider $\varphi-$dependent families of operators 
$A: \T^\nu \to \cL(L^2(\T_{\Gamma}^d, \R))$ (or $\cL(L^2(\T_{\Gamma}^d, \C))$), 
$\varphi\mapsto A(\varphi)$, acting on functions 
$u(x)\in L^2(\T_{\Gamma}^d, \R)$ (or $L^2(\T_{\Gamma}^d, \C)$). 
We also regard $A$ as an operator acting on functions 
$u(\varphi, x)\in L^2(\T_{*}^{\nu+d}, \R)$ (or $L^2(\T_{*}^{\nu+d}, \C)$) of the space-time. 
In other words, the action of $A\in \cL(L^2(\T_{*}^{\nu+d}, \R))$ is defined by
\[
(Au)(\varphi, x)= (A(\varphi)u(\varphi, \cdot))(x)\,.
\]
We represent a linear operator $\cQ$  acting on $L^2(\T_{*}^{\nu+d}, \R^2)$ 
by a matrix
\begin{equation}\label{covid}
\cQ:= \begin{pmatrix}
A && B\\
C && D
\end{pmatrix}
\end{equation}
where $A, B, C, D$  are operators 
acting on the scalar valued components 
$\eta, \psi\in L^2(\T_{*}^{\nu+d}, \R)$. 
The action of $A$ on periodic function is given by
\begin{equation}\label{linearoperator}
(Au)(\varphi, x)=
\sum_{\ell\in\Z^\nu, j\in\Gamma^*} 
\big(\sum_{\ell'\in\Z^\nu, j'\in\Gamma^*} \widehat{A}(\ell-\ell')_j^{j'}
\widehat{u}(\ell', j')\big) e^{\im(\ell\cdot\varphi+j\cdot x)}\,.
\end{equation}

We shall identify the operator $A$ with the matrix $( \widehat{A}(\ell)_j^{j'})_{\ell\in\Z^\nu, j,j'\in\Gamma^*}$. 
We say that $A$ is a \textit{real} operator 
if it maps real valued functions in real valued functions. 
For the matrix coefficients this means that
\begin{equation}\label{piscina}
\overline{\widehat{A}(\ell)_j^{j'}}=\widehat{A}(-\ell)_{-j}^{-j'}\,.
\end{equation}

It will be convenient to work with the complex variables $(u, \overline{u}) = \cC(\eta, \psi)$ 
introduced in \eqref{CVWW}.
Thus we study how linear operators $\cQ$ as in \eqref{covid} 
transform under the map $\cC$. We have
\begin{equation}\label{divano}
\begin{aligned}
\cT:=\cC\circ\cQ\circ\cC^{-1}:=\begin{pmatrix}
\cT_1 && \cT_2
\\
\overline{\cT_2} && \overline{\cT_1}
\end{pmatrix}\,, 
\end{aligned}
\end{equation}
where the \textit{conjugate} operator $\overline{\cT_{i}}$, $i=1,2$ is defined as:
\begin{equation}\label{coniugato}
\overline{\cT_{i}}[h]:=\overline{\cT_{i}[\overline{h}]}\,.
\end{equation}
In particular, recalling \eqref{covid} and \eqref{CVWW}, we note that 
\begin{equation}\label{coniugio complesse simmetrizzate}
\begin{aligned}
{\mathcal T}_{1} := & \frac12 \Big(  \Omega^{- \frac12} A \Omega^{\frac12} 
+ \Omega^{\frac12} D \Omega^{- \frac12} + \ii \Omega^{\frac12} C \Omega^{\frac12} 
- \ii \Omega^{- \frac12} B \Omega^{- \frac12} \Big)\,, 
\\
{\mathcal T}_{2} := & \frac12 \Big( \Omega^{- \frac12} A \Omega^{\frac12} 
- \Omega^{\frac12} D \Omega^{- 1} + \ii \Omega^{- \frac12} B \Omega^{- \frac12} 
+ \ii \Omega^{\frac12} C \Omega^{\frac12}   \Big)\,.
\end{aligned}
\end{equation}

\begin{defn}\label{operatorerealtoreal}
We say that $\cT(\varphi)$ is \emph{real-to-real} 
if it has the form \eqref{divano}. Equivalently 
if it preserves the subspace
\begin{equation}\label{calUUU}
\cU:=\{(u, v)\in L^2(\T_{*}^{\nu+d}, \C^2) : v=\overline{u}\}\,.
\end{equation}
\end{defn}
\noindent
By setting $u^+=u, u^-=\overline{u}$ for any $u\in L^2(\T_{*}^{\nu+d}, \C)$, 
the action of an operator $\cT$ of the form \eqref{divano} is given by
\begin{equation}\label{matricepm}
\Big(\cT \begin{bmatrix}
u\\ \overline{u}
\end{bmatrix} \Big)(\varphi, x)= \begin{pmatrix}
\cT_1u^+ + \cT_2u^-\\
\overline{\cT_2}u^+ + \overline{\cT_1}u^-
\end{pmatrix}\,.
\end{equation}
By formul\ae\, \eqref{linearoperator} and \eqref{matricepm} 
we shall identify the operator $\cT$ in \eqref{divano}
where the matrix coefficients of $\cT_1$, $\cT_{2}$ are given by 
\begin{equation}\label{matrrep}
\begin{aligned}
{\cT_1}(\vphi)_{j}^{j'}&:=\frac{1}{|\T^d_\Gamma|}
\int_{\T^{d}} \cT_{1} [e^{\im j'\cdot x} ] e^{-\ii j\cdot x}\, dx\,,
\qquad 
\widehat{\cT_1}(\ell)_{j}^{j'}:=\frac{1}{(2\pi)^\nu}\int_{\T^{\nu}}{\cT_1}(\vphi)_{j}^{j'}e^{-\ii \ell\cdot\vphi}d\vphi
\\
{\cT_2}(\vphi)_{j}^{j'}&:=\frac{1}{|\T^d_\Gamma|}
\int_{\T^{d}} \cT_{1} [e^{-\im j'\cdot x} ] e^{-\ii j\cdot x}\, dx\,,
\qquad 
\widehat{\cT_2}(\ell)_{j}^{j'}:=\frac{1}{(2\pi)^\nu}\int_{\T^{\nu}}{\cT_2}(\vphi)_{j}^{j'}e^{-\ii \ell\cdot\vphi}d\vphi\,.
\end{aligned}
\end{equation}

\subsubsection{Hamiltonian and Momentum preserving operators}
We give the following Definition.
\begin{defn}\label{operatoreHam}
Consider an operator
$\cT=\cT(\varphi)$ acting on $L^2(\T^{\nu+d}_{*}, \C^2)$ of the form \eqref{divano}.

\noindent
$(i)$ {\bf (Hamiltonian).} We say that $\cT(\varphi)$ is \emph{Hamiltonian} 
if it has the form 
\begin{equation}\label{ham1}
\cT=\im E \widetilde{\cT}\,,\quad E=\sm{1}{0}{0}{-1}\,,
\end{equation}
where $\widetilde{\cT}$ has the form \eqref{divano} and 
\begin{equation}\label{ham2}
(\widetilde{\cT}_1)^* = \widetilde{\cT}_1\,, 
\quad 
(\widetilde{\cT}_2)^* = \overline{\widetilde{\cT}_2}\,,
\end{equation}
where $(\widetilde{\cT}_i)^*$ 
denotes the adjoint operator with respect to the complex 
scalar product of $L^2(\T^d_{\Gamma}, \C)$;

\noindent
$(ii)$ {\bf (Momentum preserving).} We say that $\cT(\varphi)$ is 
\emph{momentum preserving}  
if 
 for any $\vs\in \R^d$ one has
\begin{equation}\label{vec.tau}
\cT(\vphi - {\mathtt V}\vs) \circ \tau_\vs = \tau_\vs \circ \cT(\vphi)\,;
\quad \forall \vphi \in \T^\nu
\end{equation}
(recall Definition \ref{def:vec.tau})
\end{defn}

\begin{rmk}{\bf (Self-adjointeness).}\label{rmk:selfop}
Notice that a real-to-real operator $\widetilde{\mathcal{T}}(\vphi)$ (i.e. of the form \eqref{divano})
satisfying the conditions \eqref{ham2} is self-adjoint 
with respect to the 
scalar product on $L^2(\T_{\Gamma}^{d}, \C^2)\cap\mathcal{U}$ (recall \eqref{calUUU}) 
\[
\langle U,V\rangle:={\rm Re}\int_{\T^{d}_{\Gamma}} u\cdot\bar{v} dx\,,
\qquad U=\vect{u}{\bar{u}}\,,\; 
V=\vect{v}{\bar{v}}\in L^2(\T^{d}_{\Gamma}, \C^2)\cap\mathcal{U}\,.
\]
\end{rmk}
In view of the definitions above we have the following.
\begin{lemma}{\bf (Real-to-real/Self-adjoint operators).}\label{lemma:hamfou}
Let $\cT$ be 
a real-to-real, self-adjoint operator, i.e. satisfying \eqref{divano} and \eqref{ham2} 
(recall Remark \ref{rmk:selfop}), 
and let $\widehat{\cT}_{i}(\ell)_{j}^{k}\,\; \ell\in\Z^\nu, j, k\in\Gamma^*, i=1,2$  
be the coefficients representing the operators
$\cT_i$, $i=1,2$ as in \eqref{matrrep}. Then
\begin{equation}
\overline{{\cT}_i}:= \Big( \widehat{(\overline{\cT_i})}(\ell)_{j}^{k} \Big)\,, 
\quad 
\widehat{(\overline{\cT_i})}(\ell)_{j}^{k}  = \overline{\widehat{\cT_{i}}(-\ell)_{-j}^{-k} }\,,
\end{equation}
\begin{equation}
(\cT_i)^*:=(\overline{\cT_i})^T:= \Big( \widehat{(\cT_{i})^*}(\ell)_{j}^{k} \Big)\,, 
\quad \widehat{(\cT_{i})^*}(\ell)_{j}^{k} 
 = \overline{\widehat{\cT_{i}}(-\ell)_{k}^{j} }\,.
\end{equation}
In particular, since $\cT$ is self-adjoint (see \eqref{ham2}) and real-to-real, one has
\begin{equation}
\widehat{\cT_1}(\ell)_{j}^{k}= \overline{\widehat{\cT_1}(-\ell)_{k}^{j} }\, ,
\quad 
\widehat{\cT_{2}}(\ell)_{-j}^{-k}= \widehat{((\cT_2)^T)}(\ell)_{k}^{j}= \widehat{\cT_{2}}(\ell)_{k}^{j}\,.
\end{equation}
%and 
%\begin{equation}
%(\cA)_{-\sigma, j}^{-\sigma', k}(\ell)= \overline{(\cA)_{\sigma, -j}^{\sigma', -k}(-\ell)}\,,
%\qquad
%(\cA)_{-\sigma, j}^{-\sigma', k}(\ell)= \overline{(\cA)_{\sigma, -j}^{\sigma', -k}(-\ell)}\,.
%\end{equation}
\end{lemma}
\begin{proof}
See Section $2$ in \cite{FGtrave}.
\end{proof}

\begin{rmk}{\bf (Momentum preserving).}\label{verofood}
It is easy to note that if $\mathcal{T}(\vphi)$ is momentum preserving 
  then it preserves the subspace $S_{\tV}$ in \eqref{SV}.
  Moreover,
by identifying the operator ${\cT}$ with its  matrix coefficients 
as in \eqref{matrrep}
one has that $\mathcal{T}(\vphi)$ is momentum preserving if and only
\begin{equation}
(\widehat{\widetilde{\cT_1}})(\ell)_{j}^{k}\neq 0\;\; \Rightarrow \;\;
\tV^{T}\ell + j - k =0\,,
\quad \forall\, \ell\in\Z^\nu, j, k\in\Gamma^*, \sigma, \sigma'=\pm\,,
\end{equation}
which is equivalent to 
\[
\tV^{T}\pa_{\vphi}\widetilde{\cT}(\vphi)
+\nabla \circ \widetilde{\cT}(\vphi)-\widetilde{\cT}(\vphi)\circ\nabla=0\,.
\]
\end{rmk}

\begin{rmk}\label{rmk:actionoftravelop}
According to $\eqref{vec.tau}$ an operator $\mathcal{T}(\vphi)$
acting on $\T^{\nu}\times\R^{\nu}\times H_{S}^{\perp}$
is momentum preserving 
if
\begin{equation}\label{condTTzeta}
\cT(\vphi - {\mathtt V}\vs) \circ D\tau_\vs = D\tau_\vs \circ \cT(\vphi)\,,
\qquad \forall\vs\in\R^{d}\,,
\end{equation}
where
\[
D\tau_{\vs}(\Theta,I,z)=(\Theta,I,\tau_{\vs}z)\,.
\]
If $\mathcal{T}(\vphi)$ satisfies \eqref{condTTzeta} and $i$  is a traveling
wave embedding (according to Def. \ref{def:travelembedd})
then $\mathcal{T}(\vphi)i(\vphi)$ is still a traveling wave embedding.
\end{rmk}

Consider a linear invertible map
\begin{equation}\label{portofluviale}
\T^{\nu}\ni\vphi\mapsto \Phi(\vphi)\in\mathcal{L}(H^{s}
\big(\mathbb{T}_{\Gamma}^{d};\C)\times H^{s}(\mathbb{T}_{\Gamma}^{d};\C)\big)\,,
\qquad s\geq0\,.
\end{equation}
We give the following definition.
\begin{defn}{\bf Symplectic maps}\label{def simplettica complesse}
We say that a map $\vphi \mapsto \Phi(\vphi)$ 
as in \eqref{portofluviale}
is symplectic if, for any $\vphi \in \T^\nu$, 
\[
\overline{\Phi(\vphi)}^* J \Phi(\vphi) = J\,, \quad \forall \vphi \in \T^\nu \,,\qquad 
J:=\begin{pmatrix}
0 && 1\\
-1 && 0
\end{pmatrix}\,.
\]
\end{defn}
 
\begin{lemma}
Consider an operator 
\[
\begin{aligned}
&\mathcal{L}(\vphi):=\omega\cdot\partial_{\vphi}-\mathcal{T}(\vphi)
\end{aligned}
\]
where $\mathcal{T}(\vphi)$ is real-to-real, Hamiltonian 
and momentum preserving according to 
Definition \ref{operatoreHam} and a map $\Phi(\vphi)$ as in \eqref{portofluviale}.
Assume that $\Phi(\vphi)$ is real-to-real, momentum preserving and symplectic
according to Definition \ref{def simplettica complesse}.
Then the operator 
\[
{\mathcal L}_+(\vphi) = \omega \cdot \partial_\vphi + {\mathcal T}_+(\vphi)\,, 
\quad 
{\mathcal T}_+(\vphi) := \Phi(\vphi)^{- 1} {\mathcal T}(\vphi) \Phi(\vphi) 
- \Phi(\vphi)^{- 1} \omega \cdot \partial_\vphi \Phi(\vphi)
\]
is still real-to-real, Hamiltonian and momentum preserving.
\end{lemma}
\begin{proof}
It follows straightforward by Definition 
\ref{def simplettica complesse} and Definition \ref{operatoreHam}.
\end{proof}

\subsubsection{Tame and modulo-tame operators}
In this subsection we introduce a special class of ``tame'' operators.
We first need the following Definition.
\begin{defn}\label{def:maj} 
Given a linear operator $ A $ as in \eqref{linearoperator}, we define the operator
\begin{enumerate}
\item $ | A | $  {\bf (majorant operator)} the linear operator, 
whose matrix elements are given by  $  | \widehat{A}(\ell - \ell')_j^{j'}| $,  
\item $ \Pi_N A $, $ N \in \N  $ {\bf (smoothed operator)}  
whose matrix elements are
\begin{equation}\label{proiettore-oper}
 \widehat{(\Pi_N A)}(\ell- \ell')^{j'}_j  := 
 \begin{cases}
 \widehat{A}(\ell- \ell')^{j'}_j \quad {\rm if} \quad   \langle \ell - \ell' \rangle \leq N \\
 0  \qquad \qquad \quad {\rm  otherwise} \, .
 \end{cases} 
\end{equation}
We also denote $ \Pi_N^\perp := {\rm Id} - \Pi_N $, 
\item $ \langle \pa_{\vphi} \rangle^{\tb}  A $, $ \tb \in \R $,   % the operator 
whose matrix elements are  
$  \langle \ell - \ell' \rangle^\tb \widehat{A}(\ell - \ell')_{j}^{j'} $,
\item $\partial_{\vphi_m} A(\vphi) = [\partial_{\vphi_m}, A] = 
\partial_{\vphi_m} \circ A - A \circ \partial_{\vphi_m}$ 
{\bf (differentiated operator)} 
whose matrix elements are $\ii (\ell_m - \ell_m') \widehat{A}(\ell - \ell')_{j}^{j'}$. 
\end{enumerate}
Similarly the commutator $ [\pa_{x_{i}}, A ]$ is represented by the matrix 
with entries $ \ii (j_{i} - j_{i}') \widehat{A} (\ell- \ell')_{j}^{j'} $
for $i=1,\ldots,d$.
\end{defn}
Fix $k_0\in\N$ and consider $\mathtt \Omega \subset \R^\nu$ a bounded, open set and 
\begin{equation}\label{def:setlambda0}
\lambda:=(\oo, \th)\in\mathtt{\Lambda}_{0}
\subseteq  \mathtt \Omega \times [\th_1, \th_2]\,.
\end{equation}
Let $A:=A(\lambda)$ be a linear operator $k_0$-times 
differentiable with respect to the parameter $\lambda\in\mathtt{\Lambda}_{0}$.

\begin{defn}{\bf ($\cD^{k_0}$-$\sigma$-tame operators).}\label{Dksigmatame} 
Let $0 < s_0< \bar{s} \le \infty$ be two real  numbers and $\sigma\ge 0$. 
We say that a linear operator  $A:=A(\lambda)$ is $\cD^{k_0}$-$\sigma$-tame w.r.t. 
a non-decreasing function of positive real numbers 
$\{\fM_A(\sigma, s)\}_{s=s_0}^{\bar{s}}$ 
if (recall \eqref{SobolevSpazioAngoli})
\[
\sup_{|k|\le k_0} \sup_{\lambda\in\mathtt{\Lambda}_{0}} 
\gamma^{|k|} \|\partial_\lambda^k A u\|_{s} 
\le  
\fM_A(\sigma, s)\|u\|_{s_0+\sigma} + \fM_A(\sigma, s_0)\|u\|_{s+\sigma}\,, 
\qquad u\in H^s\,,
\]
for any $s_0\le s\le \bar{s}$. 
We call $\fM_A(\sigma, s)$ a $k_{0}$\textit{-tame constant} 
for the operator $A$. When the index $\sigma$ is not relevant we
write $\fM_A(\sigma, s)=\fM_A(s)$.
\end{defn}
When the ``loss of derivatives"  $ \sigma $ is zero, 
we simply write $ \cD^{k_0} $-tame instead of $ \cD^{k_0} $-$ 0 $-tame. 
Also note that representing the operator $ A $ by its matrix elements 
$ \big(\widehat{A}(\ell - \ell')_j^{j'}  \big)_{\ell, \ell' \in \Z^\nu, j, j' \in \Gamma^*} $, 
we have, for all
$ |k| \leq k_0 $, $ j' \in \Gamma^*$, $ \ell' \in \Z^\nu $,  
\begin{equation}\label{tame-coeff}
\gamma^{2 |k|} {\mathop \sum}_{\ell , j} \langle \ell, j \rangle^{2 s} 
|\partial_\lambda^k \widehat{A}(\ell - \ell')_j^{j'}|^2 
\leq 
2 \big({\mathfrak M}_A(s_0) \big)^2 \langle \ell', j' \rangle^{2 (s+\s)} 
+ 2 \big({\mathfrak M}_A(s) \big)^2 \langle \ell', j' \rangle^{2 (s_0+\s)}\,. 
\end{equation}
%For convenience, for $\tb\in\N$ we define
%\begin{equation}\label{Mdritto}
%\mathbb{M}_{A} (s, \mathtt{b}):=
%\max_{0\leq |b|\leq \mathtt{b}} 
%\sup_{\substack{m_{1}+m_{2}=\rho-|b|\\ m_{1},m_{2}\geq0}}
%\fM_{\langle D\rangle^{m_{1}} \pa_{\vphi}^{b}A 
%\langle D\rangle^{m_{2}}}(0,s)\,.
%\end{equation}
%For properties of $ {\mathcal D}^{k_0} $-$ \s $-tame operators 
%we refer the reader to Appendix \ref{app:A}
%(see Lemmata 
%\ref{composizione operatori tame AB}-\ref{lemma operatore e funzioni dipendenti da parametro}).

We collects some properties of  $ {\mathcal D}^{k_0} $-$ \s $-modulo-tame operators.

\begin{lemma}\label{composizione operatori tame AB} {\bf (Composition).}
Let $ A, B $ be respectively $ {\mathcal D}^{k_0} $-$\sigma_A$-tame and 
$ {\mathcal D}^{k_0} $-$\sigma_B$-tame operators with tame 
constants respectively $ {\mathfrak M}_A (s) $ and $ {\mathfrak M}_B (s) $. 
Then the composition 
$ A \circ B $ is $ {\mathcal D}^{k_0} $-$(\sigma_A + \sigma_B)$-tame 
with a tame constant satisfying 
\[
 {\mathfrak M}_{A B} (s) \leq  
 C(k_0) \big( {\mathfrak M}_{A}(s) 
 {\mathfrak M}_{B} (s_0 + \sigma_A) + {\mathfrak M}_{A} (s_0) 
{\mathfrak M}_{B} (s + \sigma_A) \big)\,.
\]
The same estimate holds if $A,B$ are matrix operators as in 
\eqref{divano}.
\end{lemma}

\begin{proof}
See Lemma 2.20 in \cite{BM20}. 
\end{proof}

We now discuss the action of a $ {\mathcal D}^{k_0} $-$ \s $-tame operator  
$ A(\lambda) $ on a family of  Sobolev functions  $ u (\lambda) \in H^s $
for $\lambda$ belonging to the set $\mathtt{\Lambda}_0$ in \eqref{def:setlambda0}.

\begin{lemma}\label{lemma operatore e funzioni dipendenti da parametro}
{\bf (Action on $ H^s $)}
Let $ A := A(\lambda) $ be a $ {\mathcal D}^{k_0} $-$ \s $-tame operator. 
Then, $ \forall s_0 \leq s \leq \bar s $, for any family of Sobolev functions 
$ u := u(\lambda) \in H^{s+\s} $ 
which is $k_0$ times differentiable with respect to $ \lambda $,  
we have 
\[
\| A u \|_s^{k_0, \gamma} \lesssim_{k_0}  
{\mathfrak M}_A(s_0) \| u \|_{s + \sigma}^{k_0, \gamma} 
+ {\mathfrak M}_A(s) \| u \|_{s_0 + \sigma}^{k_0, \gamma}  \,.
\]
The same estimate holds if $A$ is a matrix operator as in 
\eqref{divano}.
\end{lemma}
\begin{proof}
See Lemma 2.22 in \cite{BM20}. 
\end{proof}

\begin{rmk}
Let $\mathtt{c}_1\geq \mathtt{c}\geq 0$. It is easy to note (using Def. \ref{Dksigmatame})
that if $A\langle D\rangle^{-\mathtt{c}}$ is  a  $\cD^{k_0}$-tame operator (with $\s=0$)
then also $A\langle D\rangle^{-\mathtt{c}_1}$ is so.
This simply follows from the fact that $\langle D\rangle^{-(\mathtt{c}_1-\mathtt{c})}$ is 
bounded on the spaces $H^{s}$.
\end{rmk}

We give the following definition.
\begin{defn}{\bf ($\cD^{k_0}$-$\sigma$-modulo-tame operators).}\label{Dk0modulotame}
Let $0<s_0< \bar{s} \le \infty$ 
be two real numbers and $\sigma\ge 0$. 
We say that a linear operator  
$A:=A(\lambda)$ is $\cD^{k_0}-\sigma$-modulo tame 
w.r.t. a non-decreasing function of positive real numbers 
$\{\fM_A^{\#}(\sigma, s)\}_{s=s_0}^{\bar{s}}$ if (recall \eqref{SobolevSpazioAngoli})
\[
\sup_{|k|\le k_0} \sup_{\lambda\in\mathtt{\Lambda}_{0}} 
\gamma^{|k|} \||\partial_\lambda^k A| u\|_{s} 
\le 
\fM_A^{\#}(\sigma, s)\|u\|_{s_0+\sigma} + \fM_A^{\#}(\sigma, s_0)\|u\|_{s+\sigma}\,, 
\qquad u\in H^s\,,
\]
for any $s_0\le s\le  \bar{s}$. 
We call $\fM_A^{\#}(\sigma, s)$ a $k_{0}$\textit{-modulo tame constant} 
for the operator $A$. When the index $\sigma$ is not relevant we
write $\fM_A^{\#}(\sigma, s)=\fM_A^{\#}(s)$.
\end{defn}

%\noindent
%For properties of $ {\mathcal D}^{k_0} $-$ \s $-modulo-tame operators 
%we refer  to Appendix \ref{app:A}
%(see Lemmata 
%\ref{A versus |A|}, \ref{interpolazione moduli parametri}).

\begin{rmk}
If $ A $, $B$ are $ {\mathcal D}^{k_0}$-modulo-tame operators, with
$ | \widehat A_j^{j'} (\ell) | \leq | \widehat B_j^{j'} (\ell) | $, then one has the estimate 
$ {\mathfrak M}^\sharp_{A}(s) \leq {\mathfrak M}^\sharp_{B}(s) $. 
\end{rmk}
We collect some results on $ {\mathcal D}^{k_0} $-$ \s $-modulo-tame operators.
\begin{lemma}\label{A versus |A|}
An operator $ A $ that is $ {\mathcal D}^{k_0}$-modulo-tame is also 
$ {\mathcal D}^{k_0}$-tame and
$ {\mathfrak M}_A (s) \leq  {\mathfrak M}_A^\sharp (s) $. 
The same holds if $A$ is a matrix operator as in \eqref{divano}.
\end{lemma}

\begin{proof}
See Lemma 2.24 in \cite{BM20}. 
\end{proof}

\begin{lemma} \label{interpolazione moduli parametri} {\bf (Sum and composition)}
Let $ A, B $ be $ {\mathcal D}^{k_0} $-modulo-tame 
operators with modulo-tame constants respectively 
$ {\mathfrak M}_A^\sharp(s) $ and $ {\mathfrak M}_B^\sharp(s) $. 
Then 
$ A+ B $ is $ {\mathcal D}^{k_0} $-modulo-tame 
with a modulo-tame constant satisfying 
\begin{equation}\label{modulo-tame-A+B}
{\mathfrak M}_{A + B}^\sharp (s) 
\leq 
{\mathfrak M}_A^\sharp (s)  + {\mathfrak M}_B^\sharp (s)  \,.
\end{equation}
The composed operator 
$  A  \circ B $ is $ {\mathcal D}^{k_0} $-modulo-tame 
with a modulo-tame constant satisfying 
\begin{equation}\label{modulo tame constant for composition}
 {\mathfrak M}_{A B}^\sharp (s) \leq  C(k_0) \big( {\mathfrak M}_{A}^\sharp(s) 
 {\mathfrak M}_{B}^\sharp (s_0) + {\mathfrak M}_{A}^\sharp (s_0) 
{\mathfrak M}_{B}^\sharp (s) \big)\,.
\end{equation}
Assume in addition that 
$ \langle \partial_{\vphi} \rangle^{\mathtt b} A $, 
$ \langle \partial_{\vphi } \rangle^{\mathtt b}  B $ 
(see Definition \ref{def:maj}) are $ {\mathcal D}^{k_0}$-modulo-tame 
with a modulo-tame constant respectively 
$ {\mathfrak M}_{\langle \partial_{\vphi} \rangle^{\mathtt b} A}^\sharp (s) $ 
and  
 $ {\mathfrak M}_{\langle \partial_{\vphi } \rangle^{\mathtt b} B}^\sharp (s) $. 
Then $ \langle \partial_{\vphi } \rangle^{\mathtt b} (A  B) $ is 
$ {\mathcal D}^{k_0}$-modulo-tame 
with a modulo-tame constant satisfying
\begin{equation}
\begin{aligned}\label{K cal A cal B}
{\mathfrak M}_{\langle \partial_{\vphi} \rangle^{\mathtt b} (A  B)}^\sharp (s) & \leq C(\tb)C(k_0)\Big(
{\mathfrak M}_{\langle \partial_{\vphi} \rangle^{\mathtt b} A}^\sharp (s) 
{\mathfrak M}_{B}^\sharp (s_0) + 
{\mathfrak M}_{\langle \partial_{\vphi} \rangle^{\mathtt b} A }^\sharp (s_0) 
{\mathfrak M}_{B}^\sharp (s) 
 \\ & 
 \qquad \qquad \qquad \quad 
 + {\mathfrak M}_{A}^\sharp (s) 
 {\mathfrak M}_{ \langle\partial_{\vphi} \rangle^{\mathtt b} B}^\sharp (s_0) 
+ {\mathfrak M}_{A}^\sharp (s_0) 
{\mathfrak M}_{ \langle \partial_{\vphi } \rangle^{\mathtt b} B}^\sharp (s)\Big)\,,
\end{aligned}
\end{equation}
for some constants $ C(k_0) , C( {\mathtt b} ) \geq 1 $. 
The same statement holds if $A$ and $B$ are matrix operators as in \eqref{divano}.
\end{lemma}
\begin{proof}
See Lemma 2.25 in \cite{BM20}. 
\end{proof}

Iterating \eqref{modulo tame constant for composition}-\eqref{K cal A cal B} 
we obtain that, for any $n \geq 2$, 
\begin{equation}\label{M Psi n}
{\mathfrak M}_{A^n}^\sharp (s) 
\leq 
\big(  2C(k_0){\mathfrak M}_{A}^\sharp (s_0) \big)^{n - 1} 
{\mathfrak M}_{A}^\sharp(s)\, , 
\end{equation}
\begin{equation}\label{K Psi n}
\begin{aligned}
{\mathfrak M}_{\langle \pa_\vphi \rangle^{ \mathtt b} A^n}^\sharp (s) 
& \leq  
(4C(\tb)C(k_0))^{n - 1} \Big( 
{\mathfrak M}^\sharp_{\langle \partial_\vphi \rangle^{\mathtt b} A}(s) 
\big[ {\mathfrak M}^\sharp_A(s_0) \big]^{n - 1} 
 \\&\quad+ {\mathfrak M}^\sharp_{\langle \partial_\vphi \rangle^{\mathtt b} A}(s_0) 
 {\mathfrak M}_A^\sharp(s) 
 \big[ {\mathfrak M}_A^\sharp(s_0) \big]^{n - 2} 
 \Big)\,.
 \end{aligned}
\end{equation}
As an application of \eqref{M Psi n}-\eqref{K Psi n}  
we obtain the following  result.

\begin{lemma}{\bf (Exponential map).}\label{serie di neumann per maggiorantia}
$(i)$ Let $ A $ and 
$ \langle \partial_{\vphi} \rangle^{\mathtt b}  A $ be 
${\mathcal D}^{k_0}$ modulo-tame operators and 
assume that 
${\mathfrak M}_{A}^\sharp : [s_0, \bar{s}] \to [0, + \infty)$ 
is a modulo-tame constant  satisfying
\begin{equation}\label{piccolezza neumann tamea}
  {\mathfrak M}_{A}^\sharp (s_0)  \leq 1 \, .
\end{equation}
Then the operators $ \Phi^{\pm 1} := {\rm exp}(\pm A) $, 
$ \Phi^{\pm 1} - {\rm Id} $ and 
$\langle \partial_\vphi \rangle^{\mathtt b}(\Phi^{\pm 1} - {\rm Id})$
 are ${\mathcal D}^{k_0}$ modulo tame operators 
 with modulo-tame constants satisfying, for any $s_0 \leq s \leq \bar{s}$, 
\begin{equation} \label{exp-MT}
\begin{aligned}
{\mathfrak M}_{ \Phi^{\pm 1} - {\rm Id}}^\sharp(s) & \lesssim  
{\mathfrak M}^\sharp_{A}(s)\,, 
\\
 {\mathfrak M}_{\langle \partial_\vphi \rangle^{\mathtt b} 
 (\Phi^{\pm 1} - {\rm Id})}^\sharp(s) &  \lesssim_{\mathtt b}  
 {\mathfrak M }_{\langle \partial_\vphi \rangle^{\mathtt b } A}^\sharp(s) 
 + {\mathfrak M}_A^\sharp(s)  
 {\mathfrak M }_{\langle \partial_\vphi \rangle^{\mathtt b } A}^\sharp(s_0) \,. 
\end{aligned}
\end{equation}
$(ii)$ Let $ A $,  $ \langle \partial_{\vphi} \rangle^{\mathtt b}  A $, 
$\langle D \rangle^{- m}A \langle D \rangle^{ m}$, 
$ \langle \partial_{\vphi} \rangle^{\mathtt b}  (\langle D \rangle^{- m}A \langle D \rangle^{ m} )$
 be ${\mathcal D}^{k_0}$ 
 modulo-tame operators and 
assume that 
${\mathfrak M}_{A}^\sharp : [s_0, \bar{s}] \to [0, + \infty)$ 
is a modulo-tame constant  satisfying
\begin{equation}\label{piccolezza neumann smooth}
  {\mathfrak M}_{\langle D \rangle^{- m}A \langle D \rangle^{ m}}^\sharp (s_0)  \leq 1 \, ,
\end{equation}
and set $ \Phi^{\pm 1} := {\rm exp}(\pm A) $.
Then the operators $\langle D \rangle^{- m} \Phi^{\pm 1} \langle D \rangle^{m}$, 
 $\langle \partial_\vphi \rangle^{\mathtt b} 
 \langle D \rangle^{- m} \Phi^{\pm 1} \langle D \rangle^{m}$
 are ${\mathcal D}^{k_0}$ modulo tame operators 
 with modulo-tame constants satisfying, for any $s_0 \leq s \leq \bar{s}$, 
\begin{equation} \label{exp-MT smooth}
\begin{aligned}
{\mathfrak M}_{ \langle D \rangle^{- m} 
\Phi^{\pm 1} \langle D \rangle^{m}}^\sharp(s) & \lesssim 
1 +     {\mathfrak M}_{\langle D \rangle^{- m}A \langle D \rangle^{ m}}^\sharp(s)\,, 
\\
 {\mathfrak M}_{\langle \partial_\vphi \rangle^{\mathtt b} (\langle D \rangle^{- m} 
 \Phi^{\pm 1} \langle D \rangle^{m})}^\sharp(s) &  \lesssim_{\mathtt b}  
 1+
 {\mathfrak M }_{ \langle \partial_{\vphi} \rangle^{\mathtt b}  
 (\langle D \rangle^{- m}A \langle D \rangle^{ m} )}^\sharp(s)  
 \\& \quad 
 + {\mathfrak M}_{\langle D \rangle^{- n}A \langle D \rangle^n}^\sharp(s)  
 {\mathfrak M }_{ \langle \partial_{\vphi} \rangle^{\mathtt b}  
 (\langle D \rangle^{- m}A \langle D \rangle^{ m} )}^\sharp(s_0) \,. 
\end{aligned}
\end{equation}
\end{lemma}

\begin{proof}
$(i)$ The item $(i)$ is proved in \cite{BMK21}. 
We observe that,
for any $n \in \N$, 
$\langle D \rangle^{- m} A^n \langle D \rangle^m 
= (\langle D \rangle^{- m} A \langle D \rangle^m)^n$, 
which implies that 
\[
\langle D \rangle^{- m}{\rm exp}(\pm A) \langle D \rangle^m 
= {\rm exp}\big( \pm \langle D \rangle^{- m} A \langle D \rangle^m \big)\,. 
\]
Therefore item $(ii)$ follows 
by the item $(i)$.
\end{proof}

\begin{lemma} \label{lemma:smoothing-tame} {\bf (Smoothing).}
Suppose that $ \langle \pa_{\vphi} \rangle^{\mathtt b} A $, 
$ {\mathtt b} \geq  0 $, is $ {\mathcal D}^{k_0} $-modulo-tame. 
Then 
the operator $ \Pi_N^\perp A $ (see Definition \ref{def:maj}) is 
$ {\mathcal D}^{k_0} $-modulo-tame with a modulo-tame constant satisfying 
\begin{equation*}%\label{proprieta tame proiettori moduli}
{\mathfrak M}_{\Pi_N^\perp A}^\sharp (s) \leq 
N^{- {\mathtt b} }{\mathfrak M}_{ \langle \pa_{\vphi} \rangle^{\mathtt b} A}^\sharp (s) \, ,
\quad
{\mathfrak M}_{\Pi_N^\perp A}^\sharp (s) \leq  {\mathfrak M}_{ A}^\sharp (s) \, . 
\end{equation*}
The same estimate holds when $A$ is a matrix operator of the form 
\eqref{divano}.
\end{lemma}
\begin{proof}
See Lemma 2.27 in \cite{BM20}. 
\end{proof}

In order to verify that an operator is modulo-tame, 
we shall use the following Lemma. 
Notice that the right hand side of \eqref{tame op modulo tame op lemma astratto} 
below contains tame constants
(not modulo-tame) of operators which control 
more space and time derivatives than 
$\langle D \rangle^\fm \langle \partial_{\vphi} \rangle^{\mathtt b} A \langle D \rangle^{-\mathtt c}$.  %\langle D \rangle^\fm$. 

\begin{lemma}\label{lem: Initialization astratto} Let  $\mathtt b, \fm, \mathfrak{c} \geq 0$, $p_0 := (\nu+d)/2 + 1$, $\beta_0 = 2p_0 + \mathtt b$  and assume that $A$ is a momentum preserving operator 
according to Definition \ref{operatoreHam}-$ii$. 
Assume also that the operators
$$\langle D\rangle^{\fm} A \langle D\rangle^{-\mathfrak{c}}\,, \quad
\langle D\rangle^{\fm} \pa_{\vphi_i}^{\beta_0 }A  \langle D\rangle^{-\mathfrak{c}}\,, \quad i=1,\ldots,\nu$$
are $\mathcal{D}^{k_0}$-tame operators according to Def. \ref{Dksigmatame} (with $s_0\rightsquigarrow p_0$ ).
Then  one has 
that
$\langle D \rangle^{\fm}  \langle \partial_{\vphi} \rangle^{\mathtt b} A \langle D \rangle^{-\mathfrak{c}} $ is $\mathcal{D}^{k_0}$- modulo tame and 
  \begin{equation}\label{tame op modulo tame op lemma astratto}
  \begin{aligned}
 {\mathfrak M}_{\langle D \rangle^{\fm}  A\langle D \rangle^{-\mathfrak{c}}  }^\sharp &( s)  \,,
 \, {\mathfrak M}_{\langle D \rangle^{\fm} \langle \partial_{\vphi} 
 \rangle^{\mathtt b} A\langle D \rangle^{-\mathfrak{c}}  }^\sharp ( s) 
 \\& \lesssim_{\mathtt b}  {\mathfrak M}_{ \langle D\rangle^{\fm } A \langle D\rangle^{-\mathfrak{c}} }(s) 
  %\\&
  + \max_{i =1, \ldots,  \nu} \big\{
{\mathfrak M}_{ \langle D\rangle^{\fm} \pa_{\vphi_i}^{\beta_0 }A  \langle D\rangle^{-\mathfrak{c}} }(s)  \big\}   \,. 
\end{aligned}
\end{equation}
\end{lemma}

\begin{proof}
We denote by $ {\mathbb N} (s, \mathtt b)  $ the right hand side in 
\eqref{tame op modulo tame op lemma astratto}.
For any $\alpha \in \N $, the matrix elements of the operator 
$\langle D \rangle^m\partial_{\vphi_i}^\alpha  A  \langle D \rangle^{m'}$
are  
$\ii^\alpha  \langle j \rangle^m (\ell_i - \ell'_i)^\alpha  \widehat{A}(\ell - \ell')_j^{j'}\langle j' \rangle^{m'} $. 
Then, by \eqref{tame-coeff} with $ \s = 0 $, applied to 
the operator
 $\langle D\rangle^{\fm} \pa_{\vphi_i}^{2p_0 +\tb}A  \langle D\rangle^{-\mathfrak{c}}$, 
 $i=1,\ldots,\nu$,
we get, using the inequality 
$ \langle \ell - \ell' \rangle^{2 (2p_0+{\mathtt b})}   \lesssim_{\mathtt b} 
1 + \max_{i = 1, \ldots, \nu}|\ell_i - \ell'_i|^{2 (2 p_0+{\mathtt b})} $, the bounds
\begin{equation}\label{carletto2}
\begin{aligned}
& \!\! \! \gamma^{2|k|}  {\mathop \sum}_{\ell, j} \langle \ell, j \rangle^{2 s} 
 \langle \ell - \ell' \rangle^{2 ({2 p_0+\tb})}  
  | \partial_\lambda^k \widehat{A}(\ell - \ell')_j^{j'}|^2 \langle j \rangle^{2\fm}  \langle j' \rangle^{-2\mathfrak{c}}  
 \\ 
& \qquad 
\lesssim_{\mathtt b}
 {\mathbb N}^2(p_0, {\mathtt b}) \langle \ell', j'\rangle^{2 s} + 
{\mathbb N}^2(s, {\mathtt b}) \langle \ell', j' \rangle^{2 p_0}  \, . 
\end{aligned}
\end{equation}
Since $A$ is a momentum preserving operator, we recall that 
(see Remark \ref{verofood})
\[
\widehat{A}(\ell - \ell')_j^{j'} \neq 0\quad  \Longrightarrow \quad\mathtt{V}^{T}(\ell - \ell') + j - j' = 0\,,
\]
implying that 
\begin{equation}\label{1309.1}
\langle j - j' \rangle \lesssim \langle \ell - \ell' \rangle  \quad \text{and hence} 
\quad \langle j - j' \rangle^{2 p_0} \lesssim \langle \ell - \ell' \rangle^{2 p_0} \,. 
\end{equation}
Therefore, for any $|k | \leq k_0 $, by the Cauchy-Schwartz inequality (using that 
$\sum_{\ell', j'} \langle \ell - \ell' \rangle^{- 2 p_0} \langle j - j' \rangle^{- 2 p_0} = C(p_0) 
< + \infty$, since $p_0 = \frac{\nu + d}{2} + 1$), 
we get
\begin{align}
\| | &\langle D \rangle^{\fm} \langle \partial_{\vphi} \rangle^{\mathtt{b}}
 \partial_\lambda^k A  \langle D \rangle^{-\mathfrak{c}}  | h \|_s^2  
 \nonumber
\\&\leq
{\mathop \sum}_{\ell, j} \langle \ell, j \rangle^{2 s} 
\Big( {\mathop \sum}_{\ell', j'}  \langle \ell - \ell' \rangle^{\mathtt b} 
 |\partial_\lambda^k   \widehat{A}(\ell - \ell')_j^{j'}| \langle j \rangle^{\fm}  \langle j' \rangle^{-\mathfrak{c}}   
 |h_{\ell', j'}| \Big)^2 
\nonumber 
\\
&= {\mathop \sum}_{\ell, j} \langle \ell, j \rangle^{2 s} 
\Big(  {\mathop \sum}_{\ell', j'} \langle \ell - \ell' \rangle^{p_0 + \mathtt b} \langle j - j' \rangle^{p_0}
|\partial_\lambda^k \widehat{A}(\ell - \ell')_j^{j'}| \langle j \rangle^{\fm} |h_{\ell', j'}| 
\frac{ \langle j' \rangle^{-\mathfrak{c}}  }{\langle \ell - \ell' \rangle^{p_0} \langle j - j' \rangle^{p_0} } \Big)^2 
\nonumber
\\
& \stackrel{\eqref{1309.1}}{\lesssim} {\mathop \sum}_{\ell, j} \langle \ell, j \rangle^{2 s} 
\Big(  {\mathop \sum}_{\ell', j'} \langle \ell - \ell' \rangle^{2p_0+ \mathtt b} 
|\partial_\lambda^k \widehat{A}(\ell - \ell')_j^{j'}| \langle j \rangle^{\fm} |h_{\ell', j'}|
\frac{ \langle j' \rangle^{-\mathfrak{c}}  }{\langle \ell - \ell' \rangle^{p_0} \langle j - j' \rangle^{p_0} } \Big)^2 
\nonumber 
\\
& \lesssim {\mathop \sum}_{\ell, j} \langle \ell, j \rangle^{2 s} 
{\mathop \sum}_{\ell', j'} \langle \ell - \ell' \rangle^{2 (2 p_0 + \mathtt b)}  
 |\partial_\lambda^k \widehat{A}(\ell - \ell')_j^{j'}|^2 |h_{\ell', j'}|^2 \langle j \rangle^{2 \fm } 
 \langle j' \rangle^{-2\mathfrak{c} } 
 \nonumber 
 \\
& \stackrel{ \eqref{carletto2}}{{\lesssim}_{\mathtt b}}  \gamma^{- 2 |k|}
{\mathop \sum}_{\ell', j'} |h_{\ell', j'}|^2 
\big( {\mathbb N}^2(p_0, \mathtt b) \langle \ell', j' \rangle^{2 s} 
+ {\mathbb N}^2(s, \mathtt b) \langle \ell', j' \rangle^{2 p_0}    \big)
\nonumber   
\\
 & \lesssim_{\mathtt b}  \gamma^{- 2|k|} 
 \big( {\mathbb N}^2(p_0, \mathtt b) \| h \|_{s}^2 + {\mathbb N}^2(s, \mathtt b) \| h \|_{p_0}^2 \big)\,.   
\nonumber
\end{align}
The claimed statement has then been proved.
\end{proof}

\subsection{Pseudo-differential operators}\label{sec:pseudodiffsimboli}

\subsubsection{Symbols and norms}
Following \cite{BM20} we give the following definitions.
\begin{defn}\label{standard} Let $u=\sum_{j\in\Gamma^*} u_je^{\im j\cdot x}$ 
and $m\in\R$. A linear operator $A$ defined by
\begin{equation}\label{standquanti}
(Au)(x):=\sum_{j\in\Gamma^*} a(x, j) \widehat{u}(j)e^{\im j\cdot x}
\end{equation}
is called \textit{pseudo-differential} of order $\le m$ if its symbol $a(x, j)$ 
is the restriction to $\T^d_{\Gamma}\times\Gamma^*$ of a function $a(x, \xi)$  
which is in $C^\infty(\T^{d}_{\Gamma}\times\R^{d})$ 
and satisfies the inequality
\begin{equation}\label{pseudo}
|\pa_{x}^\alpha \partial_\xi^\beta a(x,\xi)|\le C_{\al, \beta} \jap{\xi}^{m-|\beta|}\,, 
\quad \forall\, \al, \beta\in \N^d\,.
\end{equation}
We call $a(x,\xi)$ the symbol of the operator $A$, that we denote
\[
A=\op(a)=a(x, D)\,, \quad D:=D_x:=\frac{1}{i}\pa_{x}\,.
\]
\end{defn}

We denote by $S^m$ the class of 
all the symbols $a(x,\xi)$ satisfying \eqref{pseudo}, 
and by $OPS^m$ the set of pseudo-differential operators of order $m$. 
We set $OPS^{-\infty}:=\cap_{m\in\R} OPS^m$.

\noindent
We shall use the following notation, 
used also in \cite{AB15}. 
For any $m\in\R\setminus\{0\}$, we set 
\begin{equation}\label{opcutofffunct}
|D|^m:= \op(\chi(\xi)|\xi|^m)\,,
\end{equation}
where $\chi\in C^\infty (\R^d, \R)$ is a positive cut-off functions satisfying
\begin{equation}\label{cutofffunct}
\begin{aligned}
&\chi\in C^\infty (\R^d, \R)\,,\quad \chi(\x)=\chi(-\x)\,,\quad 0\leq \chi\leq 1\,,\quad \chi(0)=0\,,
\\&
{\rm supp}\big(1-\chi\big)\subset\big\{\x\in \R^{d}\,:\, |\x|\leq \frac{2}{3}\big\}\,.
\end{aligned}
\end{equation}
In the same spirit, for the operator $\Omega(D)$ defined in \eqref{wild}, we set
\begin{equation}\label{Omega D cut off}
\Omega (D) = {\rm Op}(\Omega(\xi)), \quad \Omega(\xi) = |\xi|^{\frac12} \tanh(\mathtt h |\xi|)^{\frac12} \chi(\xi)\,.
\end{equation}
Similarly, for any $\alpha\in\mathbb{N}^{d}$, 
we set 
\begin{equation}\label{opcutofffunct2}
\pa_{x}^{\alpha}:= \ii^{|\alpha|} \op(\x^{\alpha})\equiv  \ii^{|\alpha|} \op(\chi(\xi)\x^\alpha)\,.
\end{equation}
The identifications above do not affect the action of these kind of  pseudo-differential operators on periodic functions
since $\chi(\x)\equiv1$ for $|\x|\geq1$ and $\chi(0) = 0$.

%We also identify the Hilbert transform $\cH$, 
%acting on the $2 \pi$-periodic functions, defined by
%\begin{equation}\label{Hilbert-transf}
%\cH ( e^{\im j x} ) := - \im \, \sign(j) e^{\im jx} \,  
%\quad  \forall j \neq 0 \, ,  \quad \cH (1) := 0\,, 
%\end{equation}
%with the Fourier multiplier $\opw (- \im\, \sign(\xi) \chi(\xi) )$. 

We now introduce a norm (inspired to Metivier \cite{M08}, chapter 5) 
which controls the regularity in $( \varphi, x)$, 
and the decay in $\xi$, of the symbol $a(\varphi, x, \xi)\in S^m$, 
together with its derivatives 
$\partial_\xi^\beta a\in S^{m-|\beta|}, 0\le |\beta| \le \alpha$, 
in the Sobolev norm $\|\cdot\|_s$.

\begin{defn}\label{def:normsimbolo}
Let $A:=a(\varphi, x, D)\in OPS^m$ be a family of 
pseudo-differential operators with symbol $a(\varphi, x, \xi)\in S^m$, 
$m\in\R$. For $\alpha\in\N$, $s\ge 0$ we define the norm
\begin{equation}\label{norma pseudodiff}
\|A\|_{m,s,\alpha} := 
\max_{0\le |\beta|\le \alpha} \sup_{\xi\in\R^d} 
\|\partial_\xi^\beta a(\cdot, \cdot, \xi) \|_s \jap{\xi}^{-m+|\beta|}\,.
\end{equation}
We will use also the notation $\|a\|_{m,s,\alpha}=\|A\|_{m,s,\alpha}$.
\end{defn} 

\begin{rmk}{\bf (Monotonicity of the norm $\|\cdot\|_{m,s,\alpha}$).}\label{monpseudo} 
We remark that
\begin{align*}
m\le m' &\;\Rightarrow\; \|\cdot\|_{m',s,\alpha}\le \|\cdot\|_{m,s,\alpha}\,,
\\
s\le s' &\;\Rightarrow\; \|\cdot\|_{m,s,\alpha}\le \|\cdot\|_{m,s',\alpha}\,,
%\\
\qquad \alpha\le \alpha' \;\Rightarrow\; \|\cdot\|_{m,s,\alpha}\le \|\cdot\|_{m,s,\alpha'}\,.
\end{align*}
\end{rmk}

Along the paper we also consider families of pseudo-differential operators
\[
A(\lambda):=\op(a(\lambda,\varphi, x, \xi))
\]
which are $k_0$-times differentiable with respect to a parameter
$\lambda:=(\omega,\th)\in \mathtt{\Lambda}_0$ (see \eqref{def:setlambda0}),
where the regularity constant $k_0\in\N$ will be fixed.
Let us now define a norm that controls also the regularity on $\lambda$.
\begin{defn}\label{norma pesata}
Let $A(\lambda):=\op(a(\lambda,\varphi, x, \xi))\in OPS^m$ be a family of pseudo-differential operators with symbol $a(\lambda,\varphi, x, \xi)\in S^m$, $m\in\R$, which is $k_0$-times differentiable with respect to a parameter $\lambda$. For $\gamma\in (0,1), \alpha\in\N, s\ge0$, we define the weighted norm
$$
\|A\|_{m, s, \alpha}^{k_0, \gamma} :=
\sum_{|k|\le k_0} \gamma^{|k|} \sup_{\lambda\in\mathtt{\Lambda}_{0}} 
\|\partial_\lambda^k A(\lambda)\|_{m, s, \alpha}\,.
$$
\end{defn}

\begin{rmk}\label{rmk:delambdaop}
$\bullet$\;For $k_0, \gamma$ fixed, the same monotonicity rules of 
Remark \ref{monpseudo} hold.

\noindent
$\bullet$\; Note that 
\[
\partial_\lambda^k \op(a)=\op(\partial_\lambda^k a)\,, 
\quad \forall\,k\in\N^{\nu+1}, |k|\le k_0\,.
\]
\end{rmk}

%\noindent
%Pseudo-differential operators are tame operators according to 
%Definition \ref{Dksigmatame}. 
%We refer to Lemma \ref{lemma: action Sobolev}.

\noindent
We have  the following interpolation lemma. 
\begin{lemma}\label{interpolazione fine}
Let $ a_0, b_0 \geq 0$ and $ p,q >  0 $. 
Let $a(\lambda, \vphi, x, \xi)\in S^m$  and 
$u(\vphi, x)\in H^{b_0+q+p}( \T^{\nu+d}_*)$.
Then for all $\epsilon > 0 $ there exists a family of constant 
$ C(\epsilon) := C(\epsilon, p, q ) > 0 $, decreasing with $\epsilon$,
which satisfies $ C(1) < 1 $, such that 
\begin{equation*} 
\| a \|_{m, a_0 + p,0}^{k_0,\g} \| u \|_{b_0 + q}^{k_0, \g} \leq 
\epsilon \| a \|_{m, a_0 + p + q, 0}^{k_0,\g} \| u \|_{b_0}^{k_0, \g}
+ 
C(\epsilon)\| a \|_{m, a_0,0 }^{k_0, \g} \| u \|_{b_0 + p + q}^{k_0, \g} \,.
\end{equation*}
\end{lemma}

\begin{proof}
We prove  the lemma for the norm $\|\cdot\|_{m,s,\alpha}$ in \eqref{norma pseudodiff},
i.e. for $k_0=0$. The case $k_0\neq0$  follows similarly.
By classical interpolation estimates for the norm $\|\cdot\|_{s}$ (see for instance Lemma $6.1$ in \cite{BBM14})
we have
\begin{equation*}
\begin{aligned}
\| a(\cdot, \xi) \|_{a_0+p}\jap{\xi}^m &\leq 
\| a(\cdot, \xi)\|_{a_0}^\mu  \| a(\cdot, \xi) \|_{a_0+p+q}^{1-\mu}\jap{\xi}^m\,, 
\qquad   \mu := \frac{q}{p+q} \, ,
\\
\| u \|_{b_0+q} &\leq 
\| u \|_{b_0}^\eta  \| u \|_{b_0+p+q}^{1-\eta} \, ,
\qquad  \eta := \frac{p}{p+q} \, . 
\end{aligned}
\end{equation*}
Hence, using  that $ \eta + \mu = 1  $, we have  
\[
\| a(\cdot, \xi) \|_{a_0+p} \jap{\xi}^m \| u \|_{b_0+q}  
\leq 
( \| a(\cdot, \xi) \|_{a_0+p+q}  
\| u \|_{b_0} )^{\eta}  ( \| a(\cdot, \xi) 
\|_{a_0} \| u \|_{b_0+p+q} )^\mu\jap{\xi}^m   \, .  
\]
By the asymmetric Young inequality we get, for any $ \epsilon > 0  $, 
\[
\| a(\cdot, \xi) \|_{a_0+p}\jap{\xi}^m  \| u \|_{b_0+q}  \leq 
\epsilon \| a(\cdot, \xi) \|_{a_0+p+q}\jap{\xi}^m   \| u \|_{b_0} + C(\epsilon, p, q)  
 \| a(\cdot, \xi) \|_{a_0} \jap{\xi}^m  \| u \|_{b_0+p+q}\,,
\]
where $ C(\epsilon, p, q ) := \mu ( \eta \slash \epsilon )^{\frac{\eta}{\mu}} = \frac{q}{p + q}  
\big( \frac{p}{\epsilon (p + q)} \big)^{ p / q }  $. 
Note that for $ \epsilon = 1 $ the constant $ C(1, p, q) < 1 $. 
Passing to the $\sup$ on $\xi\in\R^d$ we obtain the desired  bound
for $k_0=0$.
The general case $k_0>0$ follows reasoning similarly
recalling Definition \ref{norma pesata}
and Remark \ref{rmk:delambdaop}.
\end{proof}

%\subsubsection{Weyl/standard quantization}

\vspace{0.5em}
\noindent
{\bf Weyl/standard quantization.}
The quantization of a symbol in $S^m$ given in \eqref{pseudo} 
is called Standard quantization. 
We shall compare this quantization with an \emph{equivalent} one: 
given a symbol $a(x,\xi)\in S^m$, we define its \emph{Weyl quantization} as 
\begin{equation}\label{Weyl}
\opw (a(x,\xi))[u]:= \sum_{j\in\Gamma^*} \Big( \sum_{k\in\Gamma^*} 
\widehat{a}\big(j-k, \frac{j+k}{2}\big)\widehat{u}(k)\Big)e^{\im j\cdot x}\,,
\qquad 
u=\sum_{k\in\Gamma^*} \widehat{u}(k) e^{\im k\cdot x}\,.
\end{equation}
Let us compare the definitions of pseudo-differential operator given 
thorough the Standard quantization (\ref{standard}) and the Weyl one. 
It is known (see paragraph 18.5 of \cite{M08}) 
that the two quantizations are equivalent
and that one can transform the symbols between 
different quantizations by using the following formulas:
\[
\opw(a)=\op(b) \,, \qquad \widehat{a}(k,j)=\widehat{b}(k, j-\frac{k}{2})\,.
\]
The next lemma guarantees that given $b(x,\xi)\in S^{m}$ also $a(x,\xi)$ is a symbol belonging to 
the same class.

\begin{lemma}\label{passaggiostweyl} 
Let $b(\lambda, \varphi, x, \xi)\in S^m$ be a $k_{0}$-differentiable symbol in $\lambda\in\mathtt{\Lambda}_{0}$. 
Consider the function 
\begin{equation}\label{simboaabb}
a(\lambda, \varphi, x, \xi)=\sum_{j\in \Gamma^{*}}e^{\ii j\cdot x}\widehat{a}(\varphi, j,\x)\,,
\qquad
\widehat{a}(\varphi, j,\x)= \widehat{b}(\varphi, j, \x-\frac{j}{2})\,,
\end{equation}
for any $(\vphi,j,\x)\in \T^{\nu}\times\Gamma^{*}\times\R^{d}$\,. Let $\sigma := (\nu + d)/2 + 1$.
Then the following holds:

\noindent
$(i)$ the symbol $a\in S^m$ satisfies 
 (recall \eqref{standquanti}, \eqref{Weyl})
\begin{equation}\label{aweylbstand}
\opw( a(\varphi, x, \xi))= \op(b(\varphi, x, \xi))\,,
\end{equation}
and the estimate
\begin{equation}\label{stima passaggio weyl-st}
\|a\|^{k_{0}, \g}_{m, s, \alpha}\lesssim_{m,s,\alpha} 
\|b\|^{k_{0}, \g}_{m, s+2\sigma + |m|+\alpha, \alpha}\,.
\end{equation}

\noindent
$(ii)$ For any  $N\in \N$ we have the expansion $a=a_{<N}+a_{\geq N}$ with
\begin{align}
a_{< N}(\varphi, x, \xi)&:=
\sum_{p=0}^{N-1}
\sum_{|k|=p} \frac{(-1)^{p}}{(2\ii)^{p}k!} \pa_{x}^{k}\pa_{\xi}^{k} b(\varphi, x, \xi)
\label{miracolo1}
\\ 
a_{\ge N}(\varphi, x, \xi)&:= \sum_{j\in\Gamma^{*}} e^{\ii j \cdot x}\sum_{|k|=N} 
\frac{(-1)^{|k|}j^k}{2^{|k|}}
\frac{|k|}{k!} \int_{0}^{1} (1-\tau)^{|k|-1} \pa_{\xi}^{k}\widehat{b}(\vphi, j, \xi-\frac{j\tau}{2})\, d\tau\,. 
	\label{miracolo2}
\end{align}
Moreover the symbols $a_{< N}, a_{\ge N}$ satisfies the following bounds
\begin{align}
\|a_{< N}(\varphi, x, \xi)\|_{m, s, \alpha}&\lesssim_{m,s, \alpha, N}\|b\|_{m, s+2\sigma+N, \alpha+N}\,,
\label{fili}
\\
\|a_{\ge N}(\varphi, x, \xi)\|_{m-N, s, \alpha}&\lesssim_{m,s, \alpha, N}\|b\|_{m, s+2\sigma+2N+|m|+\alpha, \alpha+N}\,. \label{pino}
%\quad \forall\, 0\le\alpha\le \alpha_{*}
\end{align}
\end{lemma}

\begin{proof}
\emph{Proof of item} $(i)$. 
Formula \eqref{aweylbstand} follows by a direct computation  
using the definition of the two quantizations in \eqref{standquanti} and \eqref{Weyl}
and the definition of the symbol $a$ in \eqref{simboaabb}.
In order to prove estimate \eqref{stima passaggio weyl-st}
we write, for $\x\in \R^{d}$,
\begin{align*}
\|\pa_{\x}^{\beta }a(\cdot, \xi) \|^2_{s} \jap{\xi}^{-2m+2|\beta|} &=
\sum_{j\in \Gamma^*, \ell\in \Z^\nu} \jap{\ell, j}^{2s} |\widehat{a}(\ell, j, \xi)|^2\jap{\xi}^{-2m+2|\beta|}
\\
&\stackrel{\eqref{simboaabb}}{=}\sum_{j\in \Gamma^*, \ell\in \Z^\nu} \jap{\ell, j}^{2s} 
|\widehat{b}(\ell, j, \xi-\frac{j}{2})|^2
\jap{\xi}^{-2m+2|\beta|}
%%\\&
=: S_1 + S_2
\end{align*}
where  we defined
\begin{equation}\label{alberoalto}
S_i:= \sum_{\ell\in\Z^{\nu},j\in\mathcal{B}_i(\x)}\jap{\ell, j}^{2s} 
|\widehat{b}(\ell, j, \xi-\frac{j}{2})|^2
\jap{\xi}^{-2m+2|\beta|}\,,
\end{equation}
and $\mathcal{B}_i$, $i=1,2$ are the sets
\begin{equation}\label{defcalBB1}
\mathcal{B}_1:=\{j\in\Gamma^{*}:\jap{j}\leq 2\langle\x\rangle\}\,,
\qquad  
\mathcal{B}_2:=\mathcal{B}_1^{c}\,.
\end{equation}
%A simple computation recalling the norms \eqref{norma pseudodiff} 
%and \eqref{SobolevSpazioAngoli} shows that, 
We claim that, for any $p\geq 0$, one has
\begin{equation}\label{decadimentocoeffFourierAngoli2}
|\widehat{\pa_{\x}^{\beta}b}(\ell, j, \xi)|
\lesssim_{p} 
\|b\|_{m, \sigma +p,|\beta|}\jap{\ell, j}^{-p} \jap{\xi}^{m-|\beta|}\,.
\end{equation}
Therefore, using \eqref{decadimentocoeffFourierAngoli2} with $p\rightsquigarrow s+\sigma$, we have
\[
\begin{aligned}
S_1&=\sum_{\ell\in\Z^{\nu},j\in\mathcal{B}_1(\x)}\jap{\ell, j}^{2s} 
|\widehat{b}(\ell, j, \xi-\frac{j}{2})|^2
\jap{\xi}^{-2m+2|\beta|}
\\&
\stackrel{\eqref{decadimentocoeffFourierAngoli2}}{\lesssim_{s}}
\|b\|_{m,s+2\sigma,\alpha}^2
\sum_{\ell\in\Z^{\nu},j\in\mathcal{B}_1(\x)}\langle\x-\frac{j}{2}\rangle^{2m-2|\beta|}
\jap{\xi}^{-2m+2|\beta|}\langle\ell,j\rangle^{-2\sigma}
\\&
\stackrel{\eqref{defcalBB1}}{\lesssim_{s,\alpha,m}}\|b\|_{m,s+2\sigma,\alpha}^2
\sum_{\ell\in\Z^{\nu},j\in \Gamma^*}
\langle\ell,j\rangle^{-2\sigma}\lesssim_{s,\alpha,m}\|b\|_{m,s+2\sigma,\alpha}^2\,,
\end{aligned}
\]
where we used the definition of $\mathcal{B}_1$ to get the simple estimate
$\langle\x-j/2\rangle\lesssim 2\langle\x\rangle$.
Moreover
\[
\begin{aligned}
S_2&
\stackrel{\eqref{alberoalto}, \eqref{decadimentocoeffFourierAngoli2}}{\lesssim_s}
\|b\|_{m,s+2\sigma+|m|+\alpha,\alpha}^2
\sum_{\ell\in\Z^{\nu},j\in\mathcal{B}_2(\x)}\langle\x-\frac{j}{2}\rangle^{2m-2|\beta|}
\jap{\xi}^{-2m+2|\beta|}\langle\ell,j\rangle^{-2(\sigma+|m|+|\beta|)}
\\&
\stackrel{\eqref{defcalBB1}}{\lesssim_{s,m,\alpha}}\|b\|_{m,s+2\sigma+|m|+\alpha,\alpha}^2\,,
\end{aligned}
\]
where we used that $j\in \mathcal{B}_2=\mathcal{B}_1^{c}$ so that 
$\langle\x\rangle\lesssim\langle j\rangle\lesssim\langle\ell,j\rangle$.
The discussion above implies the bound \eqref{stima passaggio weyl-st} with $k_0=0$.
 The general case for the norm $\|\cdot\|_{m,s,\alpha}^{k_0,\gamma}$ follows similarly 
recalling Remark \ref{rmk:delambdaop}.
%This implies the thesis of item $(i)$ of Lemma \ref{passaggiostweyl}.
%The proof of Lemma \ref{passaggioweylst} is similar and we omit it.

\medskip
\noindent
\emph{Proof of item $(ii)$}.  
The bound \eqref{fili} follows straightforward using the explicit formula \eqref{miracolo1}.
The estimate for the symbol  $a_{\ge N}$ in \eqref{miracolo2} is more delicate.
We note that
\begin{align*}
\|\pa_{\xi}^{\beta}a_{\ge N}(\cdot, \xi) &\|^2_{s} \jap{\xi}^{-2m+2|\beta|+2N}\lesssim S_{1} + S_{2}
\end{align*}
where
\begin{equation*}
S_{i}:= \sum_{j\in \cB_{i}(\xi), \ell\in\Z^{\nu}} \jap{\ell, j}^{2s+ 2N} \sup_{\tau\in[0,1]}\sum_{|k|=N} 
\left| \pa_{\xi}^{k+\beta}\widehat{b}(\vphi, j, \xi-\frac{j\tau}{2})\right|^{2} 
%\frac{\jap{\xi-\frac{j\tau}{2}}^{2m-2|\beta|-2N}}{\jap{\xi-\frac{j\tau}{2}}^{2m-2|\beta|-2N}}  
\jap{\xi}^{-2m+2|\beta|+2N}\,,
\end{equation*}
with $\cB_{i}(\xi)$ are defined in \eqref{defcalBB1}. Reasoning exactly as in the proof of item $(i)$, we obtain that
\begin{equation*}
S_{1} \lesssim_{m,s,\alpha, N} \|b\|_{m, s+2\sigma+N, \alpha+N}\,, 
\quad 
S_{2} 
\lesssim_{m,s,\alpha, N} 
\|b\|_{m, s+2\sigma+|m|+ 2 N+ \alpha, \alpha+N}\,,
\end{equation*}
which implies \eqref{pino} with $k_0=0$.
The general case for the norm $\|\cdot\|_{m,s,\alpha}^{k_0,\gamma}$ follows similarly 
recalling Remark \ref{rmk:delambdaop}.
%The proof of item $(ii)$ of Lemma \ref{passaggioweylst} follows as well.

\smallskip
\noindent
It remains to prove the claim \eqref{decadimentocoeffFourierAngoli2}.
Let $p,q_1,q_2\in \N$ such that $q_1+q_2\leq p$.
By a simple integration by parts, for $\ell\in\Z^{\nu}$, $k\in \Gamma^*$, one has
\begin{align*}
\ell_{j}^{q_1}k_i^{q_2} \widehat{b}(\ell, k, \xi) &= \frac{1}{(-\im)^{q_1+q_2}}\frac{1}{(2\pi)^\nu|\T^d_\Gamma|}
\int_{\T^{\nu}}\int_{\T^d_\Gamma} b(\vphi,x,\xi) \partial_{x_i}^{q_2} e^{-\im k\cdot x} 
\partial_{\vphi_j}^{q_1} e^{-\im \ell\cdot \vphi} \, dxd\vphi
\\&
= \frac{(-1)^{q_1+q_2}}{(-\im)^{q_1+q_2}}\frac{1}{(2\pi)^{\nu}|\T^d_\Gamma|}
\int_{\T^{\nu}}\int_{\T^d_\Gamma} \pa_{\vphi_j}^{q_1}\partial_{x_i}^{q_2} b(\vphi,x,\xi) 
e^{-\im k\cdot x}e^{-\im \ell\cdot \vphi} \, dxd\vphi\,,
\end{align*}
from which one deduces
\[
\langle \ell\,,\, k \rangle^{p}|\pa_{\x}^{\beta}\widehat{b}(\ell, k, \xi)|\lesssim_{p}
\|\pa_{\x}^{\beta}b(\cdot,\x)\|_{\sigma + p} 
\]
which implies \eqref{decadimentocoeffFourierAngoli2}.
\end{proof}

\noindent
The following lemma is the counterpart of Lemma \ref{passaggiostweyl}. 
Since the proof is essentially the same we omit it.
\begin{lemma}\label{passaggioweylst}  
Let $a(\lambda, \varphi, x, \xi)\in S^m$ be a $k_{0}$-differentiable symbol in $\lambda\in\mathtt{\Lambda}_{0}$. 
Consider the function 
\begin{equation}\label{simbobbaa}
b(\lambda, \varphi, x, \xi)=\sum_{j\in \Gamma^{*}}e^{\ii j\cdot x}\widehat{b}(\varphi, j,\x)\,,
\qquad
\widehat{b}(\varphi, j,\x)= \widehat{a}(\varphi, j, \x+\frac{j}{2})\,,
\end{equation}
for any $(\vphi,j,\x)\in \T^{\nu}\times\Gamma^{*}\times\R^{d}$\,. Let $\sigma := (\nu + d)/2 + 1$.
Then the following holds:

\noindent
$(i)$ the symbol $b\in S^m$ satisfies 
 (recall \eqref{standquanti}, \eqref{Weyl})
\[
\op(b(\varphi, x, \xi))=\opw( a(\varphi, x, \xi))\,,
\]
and the estimate
\begin{equation}\label{stima passaggio st-weyl}
\|b\|^{k_{0}, \g}_{m, s, \alpha}\lesssim_{m,s,\alpha} 
\|a\|^{k_{0}, \g}_{m, s+2\sigma + |m|+\alpha, \alpha} \,.
\end{equation}

\noindent
$(ii)$ For any  $N\in \N$ we have the expansion $b=b_{<N}+b_{\geq N}$ with
\begin{align*}
b_{< N}(\varphi, x, \xi)&:=
\sum_{p=0}^{N-1}
\sum_{|k|=p} \frac{1}{(2\ii)^{p}k!} \pa_{x}^{k}\pa_{\xi}^{k} a(\varphi, x, \xi)\,,
\\
b_{\ge N}(\varphi, x, \xi)&:=  
\sum_{j\in\Gamma^{*}} e^{\ii j \cdot x} 
\sum_{|k|=N}  \left(\frac{j^k}{2^{|k|}}\right)  
\frac{|k|}{k!} \int_{0}^{1} (1-\tau)^{|k|-1} \pa_{x}^{k}\pa_{\xi}^{k}a(\vphi, j, \xi+\frac{j\tau}{2})\, d\tau\,.
\end{align*}
Moreover the symbols $b_{< N}, b_{\ge N}$ satisfies the following bounds
\begin{align}
\|b_{< N}(\varphi, x, \xi)\|_{m, s, \alpha}&\lesssim_{m,s, \alpha, N}\|a\|_{m, s+2\sigma+N, \alpha+N}\,,
\label{forza1}
\\
\|b_{\ge N}(\varphi, x, \xi)\|_{m-N, s, \alpha}&\lesssim_{m,s, \alpha, N}\|a\|_{m, s+2\sigma+2N+|m|+\alpha, \alpha+N}\,.
\label{forza2}
\end{align}
\end{lemma}

\noindent
In the next lemma we show that a pseudo-differential operator is $\mathcal{D}^{k_0}$-tame.
\begin{lemma}{\bf (Action of pseudo-differential operators).}\label{lemma: action Sobolev}
Let $ A = a(\lambda, \vphi, x, D) \in OPS^0 $  be a 
family of pseudo-differential operators
that are ${k_0}$ times differentiable 
with respect to  $\lambda\in\mathtt{\Lambda}_{0}$. 
If $ \| a \|_{0, s, 0}^{k_0, \gamma} < + \infty $, 
$ s \geq s_0 > \frac{\nu + d}{2}$, then  $ A $ is ${\mathcal D}^{k_0}$-tame 
(see Def. \ref{Dksigmatame})
with a tame constant satisfying 
\begin{equation}\label{interpolazione parametri operatore funzioni}
{\mathfrak M}_A(s) \leq C(s)  \| a \|_{0, s, 0}^{k_0, \gamma}\,.
\end{equation}
As a consequence 
\begin{equation}\label{interpolazione parametri operatore funzioni (2)}
\| A h \|_s^{k_0, \gamma} \leq 
C(s_0, k_0) \| a \|_{0, s_0, 0}^{k_0, \gamma} 
\| h \|_{s}^{k_0, \gamma} 
+ C(s, k_0) \| a \|_{0, s, 0}^{k_0, \gamma} 
\| h \|_{s_0}^{k_0, \gamma}\,.
\end{equation}
The same statement holds if $ A $ 
is a matrix operator of the form \eqref{divano}.
\end{lemma}

\begin{proof}
See Lemma 2.21 in \cite{BM20} for the proof of 
\eqref{interpolazione parametri operatore funzioni}.
Then Lemma \ref{lemma operatore e funzioni dipendenti da parametro} 
implies \eqref{interpolazione parametri operatore funzioni (2)}. 
\end{proof}

\vspace{0.5em}
\noindent
{\bf Adjoint and conjugate operator.}
It is known that the $L^2$-adjoint of a pseudo-differential 
operator is pseudo-differential. 
%In the Weyl quantization it is easier, w.r.t. the standard quantization, 
%to determine the symbol of the adjoint operator. 
%In particular we have the following.

\vspace{0.5em}
\noindent
\emph{Standard quantization.} If $B=\op(b)$, with $b\in S^{m}$ then the $L^2$-adjoint is the pseudo-differential operator
\[
B^{*}=\op(b^{*})\quad {\rm with \;symbol}\quad b^{*}(x,\x)=\overline{\sum_{j\in\Gamma^*}\widehat{b}(j,\x-j)e^{\ii j\cdot x} }\,.
\]

\vspace{0.5em}
\noindent
\emph{Weyl quantization.}
If $A:=\opw(a(x, \xi))$, with $a(x, \xi)\in S^m$ then the operator $\overline{A}$ defined in \eqref{coniugato} 
and the $L^2$-adjoint have the form 
\begin{equation}\label{robertosilamenta}
\overline{A}:=\opw(\overline{a(x, -\xi)}), \qquad A^*:=\opw(\overline{a(x, \xi)})\,.
\end{equation}

\begin{rmk}\label{rmk:operatorerealtoreal}
In view of formula \eqref{robertosilamenta} one can see that
$A$ is a \emph{real} operator if and only if $\overline{a(x, -\xi)}=a(x,\x)$
while it is self-adjoint if and only if $\overline{a(x, \xi)}=a(x,\x)$.
Moreover, recalling Lemma \ref{passaggiostweyl}
one has that the   operator $\opw(a)$ is real-to-real if and only if $\op(b)$ is so.
\end{rmk}

\subsubsection{Compositions}
In this subsection we discuss the composition between pseudo-differential 
operators both in the Standard (see \eqref{standquanti}) and in the Weyl (see \eqref{Weyl}) quantizations.

\vspace{0.5em}
\noindent
\emph{Composition in Standard quantization.}
Let $c\in S^{m}$ and  $d\in S^{m'}$. Then the composition
\[
\op(c)\circ\op(d)=\op(c\#d)\,,
\]
is a pseudo-differential operator with symbol
\begin{equation}\label{composizioStandard}
c\#d(x,\x):=\sum_{j\in \Gamma^{*}}a(x,\x+j)\widehat{b}(j,\x)e^{\ii j\cdot x}
=
\sum_{j,j'\in\Gamma^*}\widehat{a}(j'-j,\x+j)\widehat{b}(j,\x)e^{\ii j'\cdot x}\,.
\end{equation}
The symbol $c\#d$ admits the following asymptotic expansion: for any $N\geq1$ one has
\begin{equation}\label{espansionecompostandard}
\begin{aligned}
c\#d(x,\x)&=c\#_{<N}d+\widetilde{r}_{N}\,,\qquad c\#_{<N}d\in S^{m+m'}\,,\;\; \widetilde{r}_{N}\in S^{m+m'-N}\,,
\\
c\#_{<N}d(x,\x)&:=
\sum_{n=0}^{N-1}c\#_{n}d(x,\x),
\qquad c\#_{n}d(x,\x):=\sum_{|\beta|=n}\frac{1}{\ii^{|\beta|}\beta!}(\pa_{\x}^{\beta}c)(x,\x)(\pa_{x}^{\beta}d)(x,\x)\,,
\\
c\#_{\geq N}d(x,\x)&:=\widetilde{r}_{N}(x,\x):=\widetilde{r}_{N,cd}(x,\x)
\\&:=\frac{1}{ \ii^{N}}\int_{0}^{1}(1-\tau)^{N-1}\sum_{|\beta|=N}
\sum_{j\in \Gamma^*}\frac{1}{\beta!}(\pa_{\x}^{\beta}c)(x,\x+\tau j)\widehat{\pa_{x}^{\beta}d}(j,\x)e^{\ii j\cdot x}d\tau\,.
\end{aligned}
\end{equation}
We recall the following important result.

\begin{lemma}{\bf (Composition in standard quantization).}\label{lemmacomposizioneSTANDARD} 
Consider two symbols $c=c(\lambda, \vphi, x, \xi)\in S^m$, $d=d(\lambda, \vphi, x, \xi)\in S^{m'}$, 
$m, m'\in\R$. Assume that $c$ and $d$ are $k_{0}$-times differentiable in 
$\lambda\in\mathtt{\Lambda}_{0}$. Then, the following hold.

\noindent
$(i)$ The symbols $c d$ and $c\# d$ in \eqref{composizioStandard} belongs to 
$S^{m+m'}$ and, for all $\al\in\N, s\ge s_0 > \frac{\nu + d}{2}$ 
\begin{equation}\label{stimasharpSTANDARD}
\begin{aligned}
\|c d\|_{m+m', s, \al}^{k_0, \gamma} \lesssim_{m, s,\alpha} 
&\|c\|_{m, s, \al}^{k_0, \gamma} \|d\|_{m', s_0, \al}^{k_0, \gamma} 
%\\&
+\|c\|_{m, s_0, \al}^{k_0, \gamma} \|d\|_{m', s, \al}^{k_0, \gamma}\,, \\
\|c\# d\|_{m+m', s, \al}^{k_0, \gamma} \lesssim_{m,  s, \alpha} 
&\|c\|_{m, s, \al}^{k_0, \gamma} \|d\|_{m', s_0+|m|+ \alpha, \al}^{k_0, \gamma} 
%\\&
+\|c\|_{m, s_0, \al}^{k_0, \gamma} \|d\|_{m', s+ |m| + \alpha, \al}^{k_0, \gamma}\,.
\end{aligned}
\end{equation}

\noindent
\noindent
$(ii)$  The symbol $c\#_{n} d$ in  
\eqref{espansionecompostandard} satisfies for all $\al\in\N, s\ge s_0 > \frac{\nu + d}{2}$, the estimate
\small{\begin{equation}\label{stimacancellettoesplicitoAlgrammo}
\|c\#_{n} d\|_{m+m'-n,s,p}^{\gamma}\lesssim_{m,s,p,n}
\sum_{\substack{\beta_1,\beta_2\in\N \\ \beta_1+\beta_2=p}}
\|c\|_{m,s,\beta_1+n}^{k_0, \gamma}
\|d\|_{m',s_0+n,\beta_2}^{k_0, \gamma}
+\|c\|_{m,s_0,\beta_1+n}^{k_0, \gamma}
\|d\|_{m',s+n,\beta_2}^{k_0, \gamma}\,.
\end{equation}}

\noindent
$(iii)$ For any integer $N\ge 1$, $\al\in\N, s\ge s_0 > \frac{\nu + d}{2}$, the symbol $\widetilde{r}_{N}$ in \eqref{espansionecompostandard} satisfies 
\begin{equation}\label{stima:resto composizioneSTANDARD}
\begin{aligned}
\|\widetilde{r}_N\|^{k_0, \gamma}_{m+m'-N, s, \al} 
\lesssim_{m, N, \al} 
&\frac{1}{N!}\big( C(s)\|c\|^{k_0, \gamma}_{m, s, N+\al} \|d\|^{k_0, \gamma}_{m', s_0+2N+|m|+\al, \al} 
\\&\qquad\qquad
+ C(s_0)\|c\|^{k_0, \gamma}_{m, s_0, N+\al} \|d\|^{k_0, \gamma}_{m', s+2N+|m|+\alpha, \al}\big)  \,,
\end{aligned}
\end{equation}
for some $C(s)>0$.

\noindent
$(iv)$ If the  symbol $d(\x)$ is a Fourier multiplier in $S^{m'}$ 
then 
 $\|c\#d\|_{m+m', s, \alpha}^{k_0, \gamma}\lesssim_{m, s,\alpha} \|c\|_{m,  \alpha}^{k_0, \gamma}$.
\end{lemma}
\begin{proof}
Items $(i), (iii), (iv)$ can be proved as  Lemma $2.13$ in \cite{BM20}. Item $(ii)$ follows using the explicit formula 
in \eqref{espansionecompostandard} and recalling Def. \ref{def:normsimbolo}.
\end{proof}
Consider two symbols $c\in S^{m}$, $d\in S^{m'}$. 
By \eqref{composizioStandard} the commutator between two pseudo-differential operators
is pseudo-differential and has the form
\begin{equation}\label{Moyalstandard}
\big[\op(c), \op(d)\big]=\op(c\star d)\,,
\qquad c\star d:=c\#d-d\#c\in S^{m+m'-1}\,.
\end{equation}
We have the following.
\begin{lemma}\label{lemmaMoyalStandard}
Under the assumptions of Lemma \ref{lemmacomposizioneSTANDARD} the following holds.
The symbol $c\star d \in S^{m + m' - 1}$ in \eqref{Moyalstandard} satisfies for any $\alpha \in \N$, $s \geq s_0 > \frac{\nu + d}{2}$
\begin{equation}\label{stimaMOYALSTANDARD}
\begin{aligned}
\|c\star d\|_{m+m' - 1, s, \al}^{k_0, \gamma} \lesssim_{m, m', s,\alpha} 
&\|c\|_{m, s+2+|m'|+\alpha, \al+1}^{k_0, \gamma} \|d\|_{m', s_0+2+|m|+ \alpha, \al+1}^{k_0, \gamma} 
\\&
+\|c\|_{m, s_0+2+|m'|+\alpha, \al+1}^{k_0, \gamma} \|d\|_{m', s+ 2+|m| + \alpha, \al+1}^{k_0, \gamma}\,.
\end{aligned}
\end{equation}
\end{lemma}
\begin{proof}
See Lemma $2.15$ in \cite{BM20}.
\end{proof}

\vspace{0.5em}
\noindent
\emph{Composition in Weyl quantization.}
If $ A = \opw(a( x, \xi))$, $ B = \opw(b(x, \xi)) $, $ m , m' \in \R $,  
are pseudo-differential operators with symbols $ a \in S^m $, $ b \in S^{m'} $ 
then the composition operator 
$ A B := A \circ B = \opw(a\#^Wb (x, \xi)) $ is a pseudo-differential operator
with  symbol 
\begin{equation}\label{formula composizione pseudo}
a\#^W b (x, \xi) =
\sum_{k\in \Gamma^*} \Big(\sum_{j\in \Gamma^*} 
\widehat{a}(k - j, \xi + \frac{j}{2}) 
\widehat{b}(j, \xi- \frac{k-j}{2})\Big) e^{\ii k\cdot x}\,.
\end{equation}
Moreover, let
\[
\sigma(D_x, D_\xi, D_y, D_\eta):=D_\xi D_y-D_x D_\eta
\]
where $D_x=-\im \pa_{x}$ and $D_\xi, D_y, D_\eta$ are similarly defined, 
then the symbol $a\#^W b$ has the following asymptotic expansion:

\begin{defn}{\bf (Asymptotic expansion of composition symbol).}\label{def: espansione simbolo composizione} 
Let $N\ge0$. Consider symbols $a(x, \xi)\in S^m $ and $b(x, \xi)\in S^{m'} $. We define
\begin{equation}\label{espanzionecomposizioneesplicita}
(a\#_N^W b)(x, \xi):= \sum_{q=0}^{N-1} \frac{1}{q!}  
\sum_{k=0}^{q} \sum_{\substack{|\al|=k\\ |\beta|=q-k}} 
\frac{(-1)^{q-k}}{(2\ii)^{q}} \left(\begin{matrix}q \\ k \end{matrix}\right)
\pa_{\xi}^{\al}\pa_{x}^{\beta}a(x, \xi) \pa_{\xi}^{\beta}\pa_{x}^{\alpha} b(x, \xi)
\end{equation}
modulo symbols in $S^{m+m'-N}$.  
\end{defn}
\begin{rmk}
Note that
\begin{equation}\label{PoissonSimbolo}
(a\#_N^W b)(x, \xi)= ab + \frac{1}{2\im} \{a, b\}\,,
\qquad
\{a, b\}:=\sum_{k=1}^{d}(\partial_{\xi_k} a \partial_{x_k} b - \partial_{x_k} a\partial_{\xi_k} b)
\end{equation}
up to a symbol in $S^{m+m'-2}$, where $\{\cdot,\cdot\}$
denotes the Poisson bracket.
For compactness we shall also write
\begin{equation}\label{restoerreenne}
r_{N}(a,b):=a\#^W b-a\#_N^W b\,.
\end{equation}
\end{rmk}
The counterpart of Lemma \ref{lemmacomposizioneSTANDARD} in the Weyl quantization is the following.

\begin{lemma}{\bf (Composition in Weyl quantization).}\label{lemmacomposizione} 
Let $a=a(\lambda, \vphi, x, \xi)\in S^m$, $b=b(\lambda, \vphi, x, \xi)\in S^{m'}$, 
$m, m'\in\R$. Assume that $a$ and $b$ are $k_{0}$-times differentiable in $\lambda\in\mathtt{\Lambda}_{0}$. Then the following holds.

\noindent
$(i)$ The symbol $a\#^Wb$ in \eqref{formula composizione pseudo} belongs to 
$S^{m+m'}$ and, for all $\al\in\N, s\ge s_0 > \frac{\nu + d}{2}$,  
%there exists 
%$\sigma=\sigma(s_0, |m|, |m'|, \alpha)=\red{2(s_0+|m|+|m'|) + 3\alpha}$ such that 
\begin{equation}\label{stimasharp}
\begin{aligned}
\|a\#^Wb\|_{m+m', s, \al}^{k_0, \gamma} \lesssim_{m,m', s,\alpha} 
&\|a\|_{m, s+\s, \al}^{k_0, \gamma} \|b\|_{m', s_0+\s , \al}^{k_0, \gamma} \\
&+\|a\|_{m, s_0+\s , \al}^{k_0, \gamma} \|b\|_{m', s+ \s , \al}^{k_0, \gamma}
\end{aligned}
\end{equation}
for some constant $\sigma = \sigma(m, m', \alpha)> 0$ large enough. 

\noindent
$(ii)$ For any integer $N\ge 1$ , the symbol $a\#^Wb$ admits the expansion (recalling \eqref{espanzionecomposizioneesplicita}
\begin{equation}\label{expansion symbol}
(a\#^W b)(x, \xi)=(a\#_N^W b)(x, \xi)+ r_N
\end{equation}
where 
$a\#_N^W b\in S^{m+m'}$ and (see \eqref{restoerreenne}) 
$r_N:=r_{N}(a,b)\in S^{m+m'-N}$.
Moreover there exist a positive constant  $\mu=\mu(N, m, m', \al) > 0$
such that for any $\alpha \in \N$, $s > s_0 > \frac{\nu + d}{2}$, 
\begin{equation}\label{stima:sharptroncato}
\begin{aligned}
\| a\#_N^W b\|^{k_{0}, \g}_{m+m', s, \alpha}
&\lesssim_{m, m', s,\alpha, N} 
\sum_{\alpha_{1} + \alpha_{2}= \alpha}
\|a\|^{k_{0}, \g}_{m, s+N, \alpha_{1}+N} \|b\|^{k_{0}, \g}_{m', s_{0}+N, \alpha_{2}+N}
\\&\qquad \qquad 
+\sum_{\alpha_{1} + \alpha_{2}= \alpha}\|a\|^{k_{0}, \g}_{m, s_{0} +N, \alpha_{1}+N} 
\|b\|^{k_{0}, \g}_{m', s+N, \alpha_{2}+N}\,,
\end{aligned}
\end{equation}

\begin{equation}\label{stima:resto composizione}
\begin{aligned}
\|r_N\|^{k_0, \gamma}_{m+m'-N, s, \al} 
\lesssim_{m, s,\al, N} 
&\frac{1}{N!}\big( C(s)\|a\|^{k_0, \gamma}_{m, s +\mu, N+\al} \|b\|^{k_0, \gamma}_{m', s_0+\mu, \al} 
\\&\qquad\qquad
+ C(s_0)\|a\|^{k_0, \gamma}_{m, s_0+\mu, N+\al} \|b\|^{k_0, \gamma}_{m', s+\mu, \al}\big)  \,.
\end{aligned}
\end{equation}

\noindent
$(iii)$ Assume that the symbol $b(\x)$ is a Fourier multiplier in $S^{m'}$ 
(i.e. it does not depend on $(\vphi,x)\in \T^{\nu+d}_*$).
Then for any $\alpha \in \N$, $s \geq s_0 > \frac{\nu + d}{2}$, we have
\begin{equation}\label{lemma composizione multiplier}
\|a\#^{W}b\|_{m+m', s, \al}^{k_0, \gamma}\lesssim_{m, m', s,\al} \|a\|_{m,s+\sigma,\al}^{k_0, \gamma}
\end{equation}
for some constant $\sigma =\sigma(m, m') > 0$. 
\end{lemma}

\begin{proof}
Given symbols $a\in S^{m}$, $b\in S^{m'}$, we consider 
 the functions 
$c, d\in C^{\infty}(\T^{\nu+d}_* \times \R^d, \C)$ 
such that
\begin{equation}\label{alberobasso}
\widehat{a}(\varphi, k,j)= \widehat{c}(\varphi, k, j-\frac{k}{2})\,, 
\qquad 
\widehat{b}(\varphi, k,j)= \widehat{d}(\varphi, k, j-\frac{k}{2})\,.
\end{equation}
Note that $c\in S^{m}$ and  $d\in S^{m'}$ by Lemma \ref{passaggioweylst}.
Moreover, by an explicit computation using \eqref{alberobasso}
and formul\ae\, \eqref{formula composizione pseudo}, \eqref{composizioStandard},  one deduces that 
\[
\opw\big((a\#^W b) (\vphi, x, \xi)\big) =\op\big( (c\#d)(\vphi,x,\x)\big)\,.
\]
We now apply the estimates \eqref{stima passaggio weyl-st} 
and \eqref{stima passaggio st-weyl} and the estimates 
on the composition of pseudo-differential operators in Standard 
quantization  by Lemma \ref{lemmacomposizioneSTANDARD}. We then get \eqref{stimasharp}.
%We obtain
%\begin{equation*}
%\begin{aligned}
%\|a\#^W b\|^{k_0, \gamma}_{m+ m', s, \alpha} &\stackrel{\eqref{stima passaggio weyl-st}}{\lesssim_{s, s_0}}
%\|c\# d\|^{k_0, \gamma}_{m+m', s+2s_0 + |m| + |m'|+\alpha, \alpha} 
%\\&
%\stackrel{\eqref{stimasharpSTANDARD}}{\lesssim_{s, s_0}}
%\|c\|^{k_0, \gamma}_{m, s+2s_0 + |m| + |m'|+\alpha, \alpha}  
%\|d\|^{k_0, \gamma}_{m', s_0 + |m|+  \alpha, \alpha}
%\\&\qquad\qquad
%+\|c\|^{k_0, \gamma}_{m, s_0, \alpha}  
%\|d\|^{k_0, \gamma}_{m', s+ 2s_0 + 2|m| + |m'|+ 2\alpha, \alpha}
% \\&
% \stackrel{\eqref{stima passaggio st-weyl}}{\lesssim_{s, s_0}}
% \|a\|^{k_0, \gamma}_{m, s+3s_0 + 2|m| + |m'|+2\alpha, \alpha}  
% \|b\|^{k_0, \gamma}_{m', 2s_0 + |m|+|m'|+  2\alpha, \alpha}
% \\&\qquad\qquad
% +\|a\|^{k_0, \gamma}_{m, 2s_0 + |m|+\alpha, \alpha}  
% \|b\|^{k_0, \gamma}_{m', s+ 3s_0 + 2|m| + 2|m'|+ 3\alpha, \alpha}\,,
%\end{aligned}
%\end{equation*}
%which implies 
%with $\sigma= 2(s_0+|m|+ |m'|)+\alpha$.

\noindent
Concerning the item $(ii)$, the estimate  \eqref{stima:sharptroncato} follows 
recalling formula  \eqref{espanzionecomposizioneesplicita} for $a\#^{W}_{N}b$, 
and by an explicit computation distributing the $\alpha$ derivatives in $x$ using Leibniz. 
%To prove \eqref{stima:resto composizione} we reason as follows.

\noindent
The estimates \eqref{stima:resto composizione} and 
\eqref{lemma composizione multiplier} follow 
immediately applying \eqref{stima passaggio weyl-st} 
and \eqref{stima passaggio st-weyl} in Lemmata 
\ref{passaggiostweyl}-\ref{passaggioweylst} and 
the composition Lemma \ref{lemmacomposizioneSTANDARD}.
%$2.13$ in \cite{BM20}.
\end{proof}

\begin{rmk}
The norm $\|\cdot\|_{m, s, \al}^{k_{0}, \g}$ in Def. \ref{norma pesata} is closed 
under composition and satisfies ``tame'' bounds w.r.t the parameters $m,s$ 
(see \eqref{stimasharpSTANDARD} and \eqref{stimasharp}).  
Regarding the third parameter $\al$ (which controls the derivatives w.r.t. the variable $\xi$) 
the composition is not tame.\footnote{If one does not need to apply Egorov theorem 
to control a change of variables generated by a pseudo-differential operator, 
then this lack of tameness is not a problem.} However, if we restrict our 
attention to $\#^W_N$ then we have tame bounds also w.r.t. $\al$ (see \eqref{stima:sharptroncato}). 
\end{rmk}

We have the following result which provides a connection between pseudo-differential 
and $\cD^{k_0} $-tame operators.
\begin{lemma}\label{constantitamesimbolo} 
Let $A=\opw(a(\lambda, \vphi, x, \xi))$ with $a\in S^{m}$, $m\in\R$, 
be a family of pseudo-differential operators that are ${k_0}$ 
times differentiable with respect to   
$\lambda\in\mathtt{\Lambda}_{0}$. 
Let $\beta_0\in\N$. For any $m_1+m_2=m$, 
the operator $\jap{D}^{-m_1}\pa_\vphi^\beta \opw(a)\jap{D}^{-m_2}$ 
is $\mathcal{D}^{k_0}-$ tame (see Def. \ref{Dksigmatame}) 
with tame constant satisfying, for $s\geq s_0>(\nu+d)/2$,
\begin{align*}
\fM_{\jap{D}^{-m_1}\opw(a)\jap{D}^{-m_2}} (s) &\lesssim_{s, m_1, m_2} 
\|a\|^{\g, k_{0}}_{m, s+\sigma, 0}\, 
\\
\fM_{\jap{D}^{-m_1}\pa_\vphi^\beta\opw(a)\jap{D}^{-m_2}} (s) &\lesssim_{s, m_1, m_2} 
\|a\|^{\g, k_{0}}_{m, s+\sigma+\beta_0, 0} \qquad |\beta|\le \beta_0\,,
\end{align*}
for some constant $\s=\s(|m|)$. 
%\red{ègiusta la seconda con quella scelta di $m_{1}, m_{2}$??}
\end{lemma}
\begin{proof}
It follows by Lemmata \ref{lemmacomposizione}, \ref{lemma: action Sobolev} 
and \ref{lemma operatore e funzioni dipendenti da parametro}.
\end{proof}

We now provide estimates on commutators between pseudo-differential operators in Weyl quantization.
Given two symbols $a\in S^{m}$, $b\in S^{m'}$,
by \eqref{formula composizione pseudo}, we have that  the commutator between two pseudo-differential operators
is pseudo-differential and has the form
\begin{equation}\label{formula commutatore}
\big[\opw(a), \opw(b)\big]=\opw(a\star^W b)\,,
\qquad a\star^W b:=a\#^Wb-b\#^Wa\in S^{m+m'-1}\,.
\end{equation}

\begin{lemma}{\bf (Commutator).}\label{commutatorsimbol}
Let $a=a(\lambda, \vphi, x, \xi)\in S^m$, $b=b(\lambda, \vphi, x, \xi)\in S^{m'}$, 
$m, m'\in\R$, be a family of symbols that are ${k_0}$ times differentiable with respect to   $\lambda\in\mathtt{\Lambda}_{0}$. 
The following hold.

\noindent
$(i)$ The symbol $a\star^Wb$ in \eqref{formula commutatore} belongs to $S^{m+m'-1}$ and there exists a constant $\sigma=\sigma(m, m', \alpha)$ such that for all $\alpha\in\N$, $s\ge s_0 > \frac{\nu + d}{2}$ 
\begin{equation}\label{stimasharp comutatore}
\begin{aligned}
\|a\star^Wb\|_{m+m'-1, s, \al}^{k_0, \gamma} 
&\lesssim_{m, m', s,\al} 
\|a\|_{m, s+\sigma, \al+1}^{k_0, \gamma} \|b\|_{m', s_0+\sigma, \al+1}^{k_0, \gamma} 
\\&\qquad \qquad\qquad
+\|a\|_{m, s_0+\sigma, \al+1}^{k_0, \gamma} \|b\|_{m', s+ \sigma, \al+1}^{k_0, \gamma}\,.
\end{aligned}
\end{equation}

\noindent
$(ii)$ The symbol  $ a \star^{W}  b$ admits the expansion 
%\begin{equation}\label{Expansion Moyal bracket}
%a \star b = - \ii \{ a, b \} + {\mathtt r_{\mathtt 2}} (a, b )\,,
%\end{equation}
\begin{equation}\label{espstar}
a\star^{W} b =a\#^{W}b-b\#^{W}a = -\ii \{a,b\}+
\sum_{|\beta|=2}^{N-1} 
(a\#^{W}_{|\beta|}b-b\#^{W}_{|\beta|}a) + \tr_N\,,
\end{equation}
where 
%\begin{equation}\label{espstar2}
%\{a, b\}:= \partial_\xi a\partial_x b - \partial_x a\partial_\xi b\,,
%\qquad
% \tr_N :=r_{N, ab}-r_{N, ba}\,.
%\end{equation}
where $\{\cdot,\cdot\}$
is the Poisson bracket in \eqref{PoissonSimbolo},
and, recalling 
\eqref{expansion symbol} and \eqref{restoerreenne},  
\[
\tr_N:=\tr_N(a,b) :=r_{N}(a,b)-r_{N}(b,a) \, . 
\]  
Moreover, the Poisson bracket $ \{ a, b \}  \in S^{m + m' - 1} $
and the remainder ${\mathtt r}_{N} (a, b ) \in S^{m+m'-N}$
satisfy for any $\alpha \in \N$, $s \geq s_0 > \frac{\nu + d}{2} $, the estimates

\begin{align}
\| \{ a, b \} \|_{m+m' -1, s, \alpha}^{k_0, \gamma} 
&  \leq_{m, m',s, \alpha} 
\| a \|_{m, s + 1, \alpha + 1}^{k_0, \gamma} 
\| b  \|_{m', s_0 + 1, \alpha + 1}^{k_0, \gamma} \nonumber
\\& \ \ \qquad \qquad\qquad\qquad
+  \| b  \|_{m, s_0 + 1, \alpha + 1}^{k_0, \gamma} 
\| a \|_{m', s + 1, \alpha + 1}^{k_0, \gamma}  \, ,\label{stima Poisson}
\\
%\| {\mathtt r_{\mathtt 2}} (a,b)\|_{m+m' -3, s, \alpha}^{k_0, \gamma} 
%&  \leq_{s,\alpha} 
%\| a \|_{m, s + \s, \alpha + 1}^{k_0, \gamma} 
%\| b  \|_{m', s_0 + \s, \alpha + 1}^{k_0, \gamma} \nonumber
%\\& \ \ \qquad\qquad \qquad\qquad\qquad
%+  \| b  \|_{m, s_0 + \s, \alpha + 1}^{k_0, \gamma} 
%\| a \|_{m', s + \s, \alpha + 1}^{k_0, \gamma}  \, .
\|\tr_N\|^{k_0, \gamma}_{m+m'-N, s, \alpha} 
&\lesssim_{m, m',s, N, \alpha} 
\|a\|^{k_0, \gamma}_{m, s+ \s, \alpha+N} 
\|b\|^{k_0, \gamma}_{m', s_0+\s, \alpha+N} \nonumber
\\
&\qquad
+ \|a\|^{k_0, \gamma}_{m, s_0+\s, \alpha+N} 
\|b\|^{k_0, \gamma}_{m', s+\s, \alpha+N} \,, \label{stima resto commutatore}
\end{align}
for some $\sigma = \sigma(N, \alpha, m, m') > 0$. 
\end{lemma}

\begin{proof}
Recalling formul\ae\, \eqref{formula composizione pseudo}, \eqref{composizioStandard}, \eqref{simbobbaa},
\eqref{formula commutatore}
and \eqref{Moyalstandard} one can prove that there are
$c\in S^{m}$, $d\in S^{m'}$ such that
$\opw(a\star^W b)= \op (c\star d)$. Therefore, the bound  \eqref{stimasharp comutatore} 
follows using Lemmata  \ref{passaggiostweyl}, \ref{passaggioweylst} and \ref{lemmaMoyalStandard}.
Estimate \eqref{stima Poisson} follows directly using 
the explicit formula \eqref{PoissonSimbolo}.
Estimate \eqref{stima resto commutatore} follows
using \eqref{stima:resto composizione}.
\end{proof}

\noindent
Iterating estimate \eqref{stimasharp comutatore}, 
given $A \in OPS^{m}$ and $B \in OPS^{m'}$, 
we  get  estimates of the operators ${\rm Ad}^n_{A}(B)$, 
$ n \in \N $, defined inductively by 
\begin{equation}\label{liberitutti3}
{\rm Ad}_{A}(B) := [A, B]\,, \qquad {\rm Ad}_A^{n + 1}(B) := [A, {\rm Ad}_A^n (B)] \, , \  n \in \N \,.
\end{equation}
\begin{lemma}\label{lemma Ad Pseudo diff}

Let $a=a(\lambda, \vphi, x, \xi)\in S^m$, $b=b(\lambda, \vphi, x, \xi)\in S^{m'}$, 
$m, m'\in\R$, be a family of symbols that are ${k_0}$-times 
differentiable with respect to   
$\lambda\in\mathtt{\Lambda}_{0} $. 
Consider $A=\opw(a)$ and $B=\opw(b)$. Then ${\rm Ad}^n_A(B) = {\rm Op}^W(p_n(a, b))$ with $p_n(a, b) \in S^{n m + m' - n}$ and for any 
$n , \alpha \in \N$, $s \geq s_0 > \frac{\nu + d}{2}$,  one has that  
\begin{equation} \label{stima Ad pseudo diff}
\begin{aligned}
\| p_n(a, b) \|_{n m + m' - n, s, \alpha}^{k_0, \gamma} 
& \lesssim_{m, m', s, \alpha,n} 
(\| a\|^{k_0, \gamma}_{m, s_0 + \mathtt c_n, \alpha + n})^n 
\| b \|_{m', s + \mathtt c_n, \alpha + n}^{k_0, \gamma}
\\
& \quad 
+ (\| a\|^{k_0, \gamma}_{m, s_0 + \mathtt c_n, \alpha + n})^{n - 1} 
\| a\|^{k_0, \gamma}_{m, s + \mathtt c_n, \alpha + n}
\| b \|_{m', s_0 + \mathtt c_n, \alpha + n}^{k_0, \gamma}\,,
\end{aligned}
\end{equation}
for some constant $\mathtt c_n=\mathtt c_n(m, m', \alpha, n) > 0$ large enough.
\end{lemma}

\begin{proof}
Estimate \eqref{stima Ad pseudo diff} follows by applying iteratively 
\eqref{stimasharp comutatore}.
\end{proof}

We conclude this subsection with a technical result 
which will be used in subsection \ref{sym.low.orderZERO}.

\begin{lemma}\label{lemma.potenze.sharp}
Let be $a\in S^{0}$ a symbol ${k_0}$ 
times differentiable with respect to   
$\lambda\in\mathtt{\Lambda}_{0}$ and let $\sigma := \frac{\nu + d}{2} + 1$. 
Then the following holds.
For any $k\in \N$
\[
\opw(a)^{k} = \opw( a^k ) + \opw (b_{k})\,,
\]
where $b\in S^{-1}$. Moreover, for any $s\ge s_{0} > \frac{\nu + d}{2}$ and $\alpha\in \N$ one has that
\begin{equation}\label{penna}
\begin{aligned}
\| a^{k} \|^{k_{0}, \gamma}_{0, s, \alpha} 
&\lesssim_{s, \alpha}  
\Big( \tC_{1} \|a\|^{k_{0}, \gamma}_{0, s_{0}, \alpha}\Big)^{k-1} \|a\|_{0, s , \alpha}^{k_{0}, \gamma}\,,
\\
\|b_{k} \|^{k_{0}, \gamma}_{-1, s, \alpha} 
&\lesssim_{s, \alpha} 
\Big( \tC_{1} 
\|a\|^{k_{0}, \gamma}_{0, s_{0}+\alpha, \alpha}\Big)^{k-1} \|a\|^{k_{0}, \g}_{0, s+ 4\sigma  + 3 + 3\alpha, \alpha+1}\,,
\end{aligned}
\end{equation}
for some $\tC_1=\tC_1(s,\alpha)>0$.
Finally there exists a constant $\delta(s, \alpha)$ small enough such 
that if for any $s\ge \s$ and $\alpha\in \N$  one has
\begin{equation}\label{piccolissimome} 
\| a \|^{k_{0}, \gamma}_{0, s_0 + 4\sigma + 3+4\alpha , \alpha+1} \le \delta(s, \alpha)\,,
\end{equation}
then 
\[
\exp(\opw(a))=\opw\big(e^a\big)
%\sum_{k={0}}^{\infty} \frac{\opw(a)^{k}}{k!} = \sum_{k={0}}^{\infty} \frac{\opw(a^{k})}{k!} 
+ \opw(b)
\]
where 
\begin{align}
\| \opw(b)\|^{k_{0}, \gamma}_{-1, s, \alpha} 
\lesssim_{s, \alpha} \|a\|^{k_{0}, \g}_{0, s+ 4\sigma + 3 + 4\alpha, \alpha+1}\,. \label{tine}
\end{align}
 \end{lemma}
 
 \begin{proof}
First of all we apply Lemma \ref{passaggioweylst}$-(ii)$ (with $N=1$) 
and we obtain that
\[
\opw(a)=\op(a) + \op(\widetilde{a})\,,\quad  \widetilde{a}\in S^{-1}
\]
with the following estimate
\begin{equation}\label{margherita2}
\begin{aligned}
\|\widetilde{a}\|^{k_{0}, \gamma}_{-1, s, \alpha}
&\lesssim_{s, \alpha}
\|a\|^{k_{0}, \gamma}_{0, s+2\sigma+2+\alpha, \alpha+1}\,.
\end{aligned}
\end{equation}
Moreover, we note that
\begin{equation}\label{margherita}
\begin{aligned}
\opw(a)^{k} &= \Big( \op(a) + \op(\widetilde{a}) \Big)^{k}= 
\sum_{j=0}^{k} \left(\begin{matrix}k\\ j \end{matrix}\right)
  \op(a)^{k-j} \circ\op(\widetilde{a})^{j}
  \\&
  = \op(a)^{k} + \sum_{j=1}^{k} \left(\begin{matrix}k\\ j \end{matrix}\right)  \op(a)^{k-j}\circ \op(\widetilde{a})^{j}\,.
\end{aligned}
\end{equation}
We prove by induction on $k\in\N$ that the operator $\op(a)^{k}$ satisfies the following:
\begin{equation}\label{margherita3}
\op(a)^{k}=\op(a^{k}) +\op(\tr_{k})\,, \quad \tr_{k}\in S^{-1}\,,
\end{equation}
where $a^{k}$ and $\tr_{k}$ satisfy, for any $s\ge s_{0}, \alpha\in\N$, the bounds
\begin{align}
\| a^{k}\|^{k_{0}, \gamma}_{0, s, \alpha} 
&\leq \tC_{1}\Big( \tC_{1}
\|a\|^{k_{0}, \gamma}_{0, s_{0}, \alpha}\Big)^{k-1} \|a\|_{0, s, \alpha}^{k_{0}, \gamma} \,,
 \label{madrid1}
\\
\| \tr_{k}\|^{k_{0}, \gamma}_{-1, s, \alpha} 
&\leq
\tC_{1} \Big( \tC_{1} \|a\|^{k_{0}, \gamma}_{0, s_{0}+\alpha+2, \alpha+1}\Big)^{k-1} 
\|a\|_{0, s+ \alpha+2, \alpha+1}^{k_{0}, \gamma} \,,\label{madrid2}
\end{align}
for some constant $\tC_{1}=\tC_{1}(s,\alpha)>0$.

\noindent
The case $k=1$ is trivial.

\vspace{0.5em}
\noindent
\emph{ Case $k=2$}.
Clearly the symbol $a^{2}$ satisfies \eqref{madrid1}.
%\eqref{penna}. 
By using Lemma \ref{lemmacomposizioneSTANDARD}
(see estimate \eqref{stima:resto composizioneSTANDARD} with $N=1$)
we have that 
\begin{equation*}
\begin{aligned}
 \|\mathtt{r}_{2} \|^{k_{0}, \gamma}_{-1, s, \alpha}
&\lesssim _{s, \alpha}
 \|a\|^{{k_{0}, \gamma}}_{0, s_{0}, \alpha+1}\|a\|^{{k_{0}, \gamma}}_{0, s+\alpha+2, \alpha}
 +  \|a\|^{{k_{0}, \gamma}}_{0, s, \alpha+1}\|a\|^{{k_{0}, \gamma}}_{0, s_{0}+\alpha+2, \alpha}
 \\
 &\lesssim _{s, \alpha}
 \|a\|^{{k_{0}, \gamma}}_{0, s_{0}+\alpha + 2, \alpha+1}\|a\|^{{k_{0}, \gamma}}_{0, s+\alpha+2, \alpha+1}
 \,.
\end{aligned}
\end{equation*}

\vspace{0.5em}
\noindent
\emph{ Case $k\geq3$}.
If we now suppose the claim true for $k\ge2$ and let us prove it for $k+1$.
Recalling the notation in \eqref{espansionecompostandard}, we write
\begin{equation*}
\begin{aligned}
\op(a)^{k+1}
&=
 \op(a)^{k}\circ \op(a)
 \stackrel{\eqref{margherita3}}{=}   \big(\op(a^{k}) + \op(\tr_{k})\big) \circ \op(a) 
=\op( a^{k}\#a) + \op( \tr_{k} \# a)
\\
&=
\op(a^{k+1}) + \op(a^{k} \#_{\geq1} a) + \op(\tr_{k} \# a)\,.
\end{aligned}
\end{equation*}
Clearly the symbol $a^{k+1}$ satisfies \eqref{madrid1} with $k \rightsquigarrow k+1$.
By estimate \eqref{stima:resto composizioneSTANDARD} with $N=1$ in 
Lemma \ref{lemmacomposizioneSTANDARD}
one has, for some constant $C=C(s,\alpha)>0$,
\begin{equation*}
\begin{aligned}
\| a^{k}\#_{\geq1} a \|^{k_{0}, \gamma}_{0,s,\alpha}
&\leq
C  \|a^{k}\|^{{k_{0}, \gamma}}_{0, s, \alpha+1}\|a\|^{{k_{0}, \gamma}}_{0, s_{0}+\alpha+2, \alpha}
 + C \|a^{k}\|^{{k_{0}, \gamma}}_{0, s, \alpha+1}\|a\|^{{k_{0}, \gamma}}_{0, s_{0}+\alpha+2, \alpha}
 \\&
\stackrel{\eqref{madrid1}}{\leq}
C \tC_1( \tC_1 \|a\|^{{k_{0}, \gamma}}_{0, s_{0}, \alpha+1})^{k-1}\|a\|^{{k_{0}, \gamma}}_{0, s, \alpha+1}
\|a\|^{{k_{0}, \gamma}}_{0, s_{0}+\alpha+2, \alpha}
\\ &
+ C (\tC_1\|a^{k}\|^{{k_{0}, \gamma}}_{0, s_{0}, \alpha+1})^{k}
\|a\|^{{k_{0}, \gamma}}_{0, s+\alpha+2, \alpha}
\\&
\leq 2C (\tC_1\|a^{k}\|^{{k_{0}, \gamma}}_{0, s_{0}+\alpha+2, \alpha+1})^{k}
\|a\|^{{k_{0}, \gamma}}_{0, s+\alpha+2, \alpha+1}\,.
\end{aligned}
\end{equation*}
%so the estimate holds provided that we take $\tC_{1}(s, \alpha)\ge 2C(s, \alpha)$. 
Reasoning as above, using \eqref{stimasharpSTANDARD} and the inductive assumption \eqref{madrid2}, 
we have also
\begin{equation*}
\begin{aligned}
\| \tr_{k}\# a\|^{k_{0}, \gamma}_{-1,s,\alpha} 
&\leq
C\| \tr_{k} \|^{k_{0}, \gamma}_{-1, s, \alpha}\| a \|^{k_{0}, \gamma}_{0, s_{0}+\alpha+1, \alpha}
+ C\| \tr_{k} \|^{k_{0}, \gamma}_{-1, s_{0}, \alpha}\| a \|^{k_{0}, \gamma}_{0, s+\alpha+1, \alpha}
\\&
\leq 
C \tC_{1} \Big( \tC_{1} \|a\|^{k_{0}, \gamma}_{0, s_{0}+\alpha + 2, \alpha+1}\Big)^{k-1} 
\|a\|_{0, s+ \alpha+2, \alpha+1}^{k_{0}, \gamma}
\| a \|^{k_{0}, \gamma}_{0, s_{0}+\alpha+1, \alpha}
\\&
+  C \Big( \tC_{1} \|a\|^{k_{0}, \gamma}_{0, s_{0}+\alpha+2, \alpha+1}\Big)^{k} 
\| a \|^{k_{0}, \gamma}_{0, s+\alpha+1, \alpha}
\\&
\leq 2C (\tC_1\|a^{k}\|^{{k_{0}, \gamma}}_{0, s_{0}+\alpha+2, \alpha+1})^{k}
\|a\|^{{k_{0}, \gamma}}_{0, s+\alpha+2, \alpha+1}\,.
\end{aligned}
\end{equation*}
Therefore, by collecting the estimates above, 
setting $\mathtt{r}_{k+1}:=a^{k} \#_{\geq1} a+\mathtt{r}_{k}\#a$ and provided that $\tC_1\geq 4 C$,  we deduce \eqref{madrid2} 
with $k\rightsquigarrow k+1$.

Let us now study the second summand in \eqref{margherita}.
By \eqref{stimasharpSTANDARD} in Lemma \ref{lemmacomposizioneSTANDARD}
we first note that there are constants $\mathtt{C}_1:=\mathtt{C}_1(s,\alpha)\geq 
\widetilde{\mathtt{C}}_1:=\widetilde{\mathtt{C}}_1(s,\alpha)>0$ such that 
(recall $j\geq1$)\footnote{Here we also use that $\widetilde{a}\in S^{-1}$ 
and so 
$\|\widetilde{a}\|^{k_{0}, \gamma}_{0,s,\alpha}
\lesssim \|\widetilde{a}\|^{k_{0}, \gamma}_{-1,s,\alpha}$.
}
\[
\begin{aligned}
\|\op(a)^{k-j}\|^{k_{0}, \gamma}_{0,s,\alpha}&\leq \mathtt{C}_{1}
\big(\mathtt{C}_{1}\|a\|^{k_{0}, \gamma}_{0,s_0+\alpha,\alpha}\big)^{k-j-1}
\|a\|^{k_{0}, \gamma}_{0,s+\alpha,\alpha}\,,
\\
\|\op(\widetilde{a})^{j-1}\|^{k_{0}, \gamma}_{0,s,\alpha}
%&\leq 
%\widetilde{\mathtt{C}}_1
%\big(\widetilde{\mathtt{C}}_1\|\widetilde{a}\|^{k_{0}, \gamma}_{0,s_0+\alpha,\alpha}\big)^{\max\{j-2,0\}}
%\|\widetilde{a}\|^{k_{0}, \gamma}_{0,s+\alpha,\alpha}
%\\
&\leq  
\widetilde{\mathtt{C}}_1
\big(\widetilde{\mathtt{C}}_1
\|\widetilde{a}\|^{k_{0}, \gamma}_{-1,s_0+\alpha,\alpha}
\big)^{\max\{j-2,0\}}
\|\widetilde{a}\|^{k_{0}, \gamma}_{-1,s+\alpha,\alpha}\,.
\end{aligned}
\]
In particular, using again the bound \eqref{stimasharpSTANDARD} with $m=0$ and $m'=-1$, we also get
\[
\begin{aligned}
\|\op(\widetilde{a})^{j}\|^{k_{0}, \gamma}_{-1,s,\alpha}
&\lesssim_{s} 
\|\op(\widetilde{a})^{j-1}\|^{k_{0}, \gamma}_{0,s,\alpha}
\|\op(\widetilde{a})\|^{k_{0}, \gamma}_{-1,s_0+\alpha,\alpha}
+
\|\op(\widetilde{a})^{j-1}\|^{k_{0}, \gamma}_{0,s_0,\alpha}
\|\op(\widetilde{a})\|^{k_{0}, \gamma}_{-1,s+\alpha,\alpha}
\\&
\lesssim_{s}
\widetilde{\mathtt{C}}_1
\big(\widetilde{\mathtt{C}}_1
\|\widetilde{a}\|^{k_{0}, \gamma}_{-1,s_0+\alpha,\alpha}
\big)^{j-1}
\|\widetilde{a}\|^{k_{0}, \gamma}_{-1,s+\alpha,\alpha}
\\&\stackrel{\eqref{margherita2}}{\leq}
\mathtt{C}_{1}\big(\mathtt{C}_{1}\|{a}\|^{k_{0}, \gamma}_{0,s_0 + 2 \sigma+2+2\alpha,\alpha+1}\big)^{j-1}
\|{a}\|^{k_{0}, \gamma}_{0,s+2\sigma+2+2\alpha,\alpha+1}\,,
\end{aligned}
\]
provided $\tC_1=\tC_1(s,\alpha)>0$ is large enough.
By the estimates above, using again  Lemma \ref{lemmacomposizioneSTANDARD}
and the monotonicity of the norms (see Remark \ref{monpseudo}),
 we deduce
\begin{equation}\label{margherita5}
\begin{aligned}
\|\op(a)^{k-j}\circ\op(\widetilde{a})^{j}\|^{k_{0}, \gamma}_{-1,s,\alpha}
&\lesssim_{s}\|\op(a)^{k-j}\|^{k_{0}, \gamma}_{0,s,\alpha}\|\op(\widetilde{a})^{j}\|^{k_{0}, \gamma}_{-1,s_0+\alpha,\alpha}
\\&+\|\op(a)^{k-j}\|^{k_{0}, \gamma}_{0,s_0,\alpha}\|\op(\widetilde{a})^{j}\|_{-1,s+\alpha,\alpha}
\\&\leq
\mathtt{C}\big(\mathtt{C}\|a\|^{k_{0}, \gamma}_{0,s_0+ 2 \sigma +2+3\alpha,\alpha+1}\big)^{k}
\|a\|^{k_{0}, \gamma}_{0,s+2\sigma+2+3\alpha,\alpha+1}\,,
\end{aligned}
\end{equation}
for some $\mathtt{C}=\mathtt{C}(s,\alpha)\geq \mathtt{C}_1$.
Therefore, recalling \eqref{margherita}-\eqref{margherita3},  we obtained that
\begin{equation}\label{margherita11}
\begin{aligned}
\opw(a)^{k}
&= \op(a^{k}) + \op(\tc_{k})\\
\op(\tc_{k})&:= \op(\tr_{k}) + \sum_{j=1}^{k} \left(\begin{matrix}k\\ j \end{matrix}\right)  \op(a)^{k-j}\circ \op(\widetilde{a})^{j}
\end{aligned}
\end{equation}
with (see \eqref{madrid2}, \eqref{margherita5})
\begin{equation}\label{margherita10}
\| \tc_{k}\|^{k_{0}, \g}_{-1, s, \alpha} \leq k!
\mathtt{C}\big(\mathtt{C}\|a\|^{k_{0}, \gamma}_{0,s_0 + 2 \sigma +2+3\alpha,\alpha+1}\big)^{k}
\|a\|^{k_{0}, \gamma}_{0,s+2\sigma+2+3\alpha,\alpha+1}\,.
\end{equation}
We now apply Lemma \ref{passaggiostweyl}-$(ii)$ to the symbol $a^{k}$. Using the expansions 
\eqref{miracolo1}-\eqref{miracolo2} with $N=1$ we obtain
\[
\op(a^{k})=\opw(a^k)+\opw(\mathtt{d}_{k})\,,\qquad \mathtt{d}_{k}\in S^{-1}\,,
\]
with $\mathtt{d}_k$ satisfying (recall \eqref{pino} and \eqref{madrid1})
\begin{equation}\label{margherita20}
\begin{aligned}
\|\mathtt{d}_{k}\|^{k_0,\gamma}_{-1, s, \alpha}&\lesssim_{s, \alpha}
\|a^k\|^{k_0,\gamma}_{0, s+2s_{0}+2+\alpha, \alpha+1}
\\&\stackrel{\eqref{madrid1}}{\leq}
 \tC \Big( \tC \|a\|^{k_{0}, \gamma}_{0, s_{0} + 2 \sigma +2+\alpha, \alpha+1}\Big)^{k-1} 
 \|a\|_{0, s+2\sigma+2+\alpha, \alpha+1}^{k_{0}, \gamma} \,,
 \end{aligned}
\end{equation}
for some $\tC=\mathtt{C}(s,\alpha)\gg1$.
Reasoning similarly, using Lemma \ref{passaggiostweyl}-(i), for the symbol $\mathtt{c}_{k}$ we get
\[
\op(\mathtt{c}_{k})=\opw(\widetilde{\mathtt{c}}_{k})\,,\quad \widetilde{\mathtt{c}}_{k}\in S^{-1}\,,
\]
where
\begin{equation}\label{margherita21}
\begin{aligned}
\|\widetilde{\mathtt{c}}_{k}\|^{k_0,\gamma}_{-1,s,\alpha}
&\stackrel{\eqref{stima passaggio weyl-st}}{\lesssim_{s,\alpha}}
\|\mathtt{c}_{k}\|^{k_{0}, \g}_{-1, s+2\sigma + 1+\alpha, \alpha}
\\&\stackrel{\eqref{margherita10}}{\leq}k!
\mathtt{C}\big(\mathtt{C}\|a\|^{k_{0}, \gamma}_{0,s_0 + 4\sigma+3+4\alpha,\alpha+1}\big)^{k}
\|a\|^{k_{0}, \gamma}_{0,s+4\sigma+3+4\alpha,\alpha+1}\,.
\end{aligned}
\end{equation}
As a consequence, \eqref{margherita11} becomes
\[
\opw(a)^{k} = \opw(a^{k}) + \opw(\widetilde{\mathtt{c}}_{k}+\mathtt{d}_{k})\,,
\]
which implies
\[
\sum_{k={0}}^{\infty} \frac{\opw(a)^{k}}{k!} = \sum_{k={0}}^{\infty} \frac{\opw(a^{k})}{k!} + \opw(b)\,,
\qquad {b}:=\sum_{k\geq0}\frac{1}{k!} (\widetilde{\mathtt{c}}_{k}+\mathtt{d}_{k})\,.
\]
The estimate \eqref{tine} on the symbol $b$ defined above follows by
\eqref{margherita20}, \eqref{margherita21} and recalling the smallness condition 
\eqref{piccolissimome}.
This concludes the proof.
 \end{proof}

\vspace{0.5em}
\noindent
{\bf Notation.} 
Let $a,b,c,d$ be symbols in the class $S^m$, 
$m\in\R$, $k_0$-times differentiable w.r.t. $\lambda \in \mathtt \Lambda_0$. We shall write 
\begin{equation}\label{matricidisimboli}
A:=A(\varphi, x, \xi):=
 \begin{pmatrix}
a && b\\
c && d
\end{pmatrix}%\in S^m\otimes \cM_2(\C)
\end{equation} 
to denote $2\times 2$ matrices of symbols.
With a slight abuse of notation we still write $A\in S^{m}$.
We define $\op(A)$ in the natural way as the $2\times2$ matrix 
whose entries are the pseudo-differential operators associated to the symbols $a,b,c,d$.

With abuse of notation we write $\|A\|_{m,s, \alpha}^{k_{0}, \g}$ to denote 
$\max\{ \|f\|_{m, s, \alpha}^{k_{0}, \g}, f=a,b,c,d\}$.
 %and we write $A\in S^m$ when $A\in S^m\otimes \cM_2(\C)$.

\subsubsection{Algebraic properties of symbols}
We provide some characterizations of momentum preserving 
pseudo-differential operators in terms of symbols. 

\begin{defn}\label{def:mompressimbo}
Consider a symbol $a(\varphi, x, \xi)\in S^m$, $m\in\R$. We say that $a(\vphi,x,\x)$ is
\emph{momentum preserving} if the operator $\opw(a) $ is momentum preserving 
according to Definition \ref{operatoreHam}.
\end{defn}

\begin{lemma}{\bf (Momentum preserving symbols).} \label{lem:mompressimbolo}
Consider a symbol $a(\varphi, x, \xi)\in S^m$, $m\in\R$. 
The following conditions are equivalent.

\begin{enumerate}%[(a)]
\item $a(\vphi,x,\x)$ is momentum  preserving according to Definition \ref{def:mompressimbo}.
\item There is $\widetilde{a}(\Theta,\x)$, $\Theta\in\T^\nu, \xi\in\R^d$ such that
\begin{equation}\label{invariantsymbol}
a(\varphi, x, \xi)=\widetilde{a}(\varphi-\tV x, \xi)\,.
\end{equation}
\item for any $\vs\in \R^{d}$ one has
\[
a(\vphi-\tV \vs,x,\x)=a(\vphi,x+\vs,\x)\,,
\qquad \forall (\vphi,x)\in \T^{\nu+d}_{*}\,,\;\; \x\in \Gamma^{*}\,.
\]
\end{enumerate}
\end{lemma}

\begin{proof}
It follows straightforward by applying Definition \ref{operatoreHam} to $A=\opw(a)$.
\end{proof}

\begin{lemma}{\bf (Real-to-real/Self-adjoint matrices of symbols).}
\label{lemma hamiltoniano simboli}
Consider the matrix of operators (recall \eqref{matricidisimboli})
\[
L=\opw(A(\varphi, x, \xi)), \qquad A\in S^m \,.
\]
We have that $L$ is {real-to-real} according to Definition \ref{operatoreHam}
%\ref{realtoreal} 
if and only if the matrix $A$ has the form
\begin{equation}\label{matricedisimboli}
A(\varphi, x, \xi):= \begin{pmatrix}
a(\varphi, x, \xi) && b(\varphi, x, \xi)\\
\overline{b(\varphi, x, -\xi)} && \overline{a(\varphi, x, -\xi)}
\end{pmatrix}\,,
\end{equation}
and $L$ is {self-adjoint}, i.e. satisfies \eqref{ham2}, if and only if the matrix of symbols $A$ satisfies
\begin{equation}\label{matrixSymself}
\overline{a(\varphi, x, \xi)}=  a(\varphi, x, \xi), \qquad b(\varphi, x, -\xi)=b(\varphi, x, \xi)\,.
\end{equation}
Finally $L$ is {Hamiltonian}, i.e. satisfies \eqref{ham1} if and only if the matrix $A$ satisfies
\begin{equation*}
\overline{a(\varphi, x, \xi)}= - a(\varphi, x, \xi), \qquad b(\varphi, x, -\xi)=b(\varphi, x, \xi)\,.
\end{equation*}
\end{lemma}
\begin{proof}
 It follows straightforward by applying Definitions \ref{operatorerealtoreal},  \ref{operatoreHam} and \eqref{robertosilamenta}
\end{proof}

\begin{rmk}\label{rmk:algsimboli}
We remark that, in view of Lemma \ref{lem:mompressimbolo}, 
if $a$ and $b$ are two momentum preserving symbols, then also
the symbols $a\#^{W}b$, $a\#^{W}_{N}b$ 
in \eqref{formula composizione pseudo} and 
\eqref{espanzionecomposizioneesplicita} respectively are momentum preserving symbols.
The same holds for $a\# b$, $a\#_{N}b$  in \eqref{composizioStandard} and \eqref{espansionecompostandard}.
\end{rmk}

\smallskip
{\bf Notation.}
We will also a dependence on a parameter 
$\T^{\nu}\ni\vphi\to i(\vphi)\in \T^{\nu}\times \R^{\nu} \times H_{S}^{\perp}$ 
of which we want to control the regularity. 
In particular, we will ask to be Lipschitz in the parameter $i$ 
and therefore it is natural to introduce the following quantity.\\
Given the map
$i \mapsto g(i)$ (see \eqref{torusprototipo})
where $g$ is an operator (or a map or a scalar function), we define 
\begin{equation}\label{deltaunodue}
\Delta_{12} g := g(i_{2}) - g(i_{1})
\end{equation}

\begin{rmk}\label{rmk:deltaunodueop}
 Note that  for a pseudo-differential operator one has
$\Delta_{12}\op(a)=\op(\Delta_{12}a)$.
\end{rmk}

\section{Egorov Theory and conjugations}\label{sec:Egorov}
We consider the flow map of
\begin{equation}\label{flussosimboloGenerico}
\begin{cases}
\partial_\tau \Phi_a^\tau ( \vphi) = \ii {\rm Op}^W\big( a(\tau;\vphi, x, \xi) \big) \circ \Phi_a^\tau ( \vphi) 
\\
\Phi_a^0( \vphi) = {\rm Id}\,. 
\end{cases}\qquad \vphi\in \T^{\nu}\,,\;\; x\in \T^d_{\Gamma}\,, \quad \tau \in [- 1, 1]
\end{equation}
where $a(\vphi,x,\x)$ 
is a symbol assumed to have one of the following forms:
\begin{align}
a(\tau;\vphi,x,\x)&:=\breve{b}(\tau; \varphi, x)\cdot\x\,,\qquad 
\breve{b}(\tau; \varphi, x):=(\uno+\tau \nabla_x\breve\beta(\varphi, x))^{-T}\breve\beta(\varphi, x)\,,
\label{egogenerator}
\\
 \breve\beta(\varphi, x)&\in H^s(\T_{*}^{\nu+d};\R^d)\,,
\qquad \uno:=\sm{1}{0}{0}{1}\,,
\nonumber
\\
a(\tau;\vphi,x,\x)&\equiv a(\vphi,x,\x)\in S^{m}\,,\qquad m<1\,,\qquad {\rm real\; valued}\,.
\label{egogenerator2}
\end{align}
 The symbol  $a := a(i)$ may depend also on the “approximate” torus $i(\vphi)$
 (recall Def.  \ref{def:travelembedd}).
First of all we state a result which guarantees the  
well-posedness of the flow  \eqref{flussosimboloGenerico}. 

\begin{lemma}\label{lemma:buonaposFlussi}{\bf (Well-posedness of flows)}
Let $a(\vphi, x, \xi)\in S^{m}$, $m\le1$ 
with $a=a(\lambda, i(\lambda))$, $k_{0}$-differentiable 
in $\lambda\in\mathtt{\Lambda}_{0}$ 
and Lipschitz in the variable $i$. % and fix $p\in\N$.
Consider the flow $\Phi_a$ in \eqref{flussosimboloGenerico}.
%The following holds.
For any $m_1, m_2 \in \Z$  such that $m_1 + m_2 \ge (\beta_0 + k_{0}) m $,
 there exist positive constants 
$\sigma_{1} = \sigma_{1}(\beta_0, m_1, m_2)\geq \sigma_{2} = \sigma_{2}(\beta_0, m_1, m_2)$ 
such that the following holds.
%with for $p+\s_{2}\le s_{0} + \s_{1}$ and 
For any $\bar{s}\geq s_0+\s_1$, $s_0 > \frac{d + \nu}{2}$
there exists $\delta=\delta(\beta_{0}, m_{1}, m_{2},\bar{s})$ 
such that if 
\[
\| a \|_{m, s_0 + \sigma_{1}, \sigma_{1}}^{k_0, \gamma}\le \delta \,, 
\]
then for any $\beta \in \N^\nu$, $|\beta| \leq \beta_0$, 
the linear operator 
$\langle D \rangle^{-m_1} \partial_\vphi^\beta \Phi_a^\tau \langle D \rangle^{-m_2}$ 
is ${\mathcal D}^{k_0}$-tame 
with tame constant satisfying, for $ s_0\leq s\leq \bar{s}-\s_1$ and $p+\s_2\leq s_0+\s_1$, 
\[
{\mathfrak M}_{\langle D \rangle^{ -m_1} \partial_\vphi^\beta 
\Phi_a^\tau \langle D \rangle^{ -m_2}}(s) 
\lesssim_{\beta_0, m_1, m_2, s} 
1 + \| a \|_{m, s + \sigma_{1}, \sigma_{1}}^{k_0, \gamma}\,, 
\]
\[
\sup_{\tau\in[0,1]}
%sup_{|\beta|\leq \beta_{0}}
\sup_{m_1+m_2=m(\beta_{0}+1)}
\|\jap{D}^{-m_1} \pa_\vphi^\beta  \Delta_{12} \Phi_a^\tau \jap{D}^{-m_2}\|_{\mathcal{L}(H^p;H^p)}\lesssim_{\beta_0, m_1, m_2, p} \|\Delta_{12} a\|^{ k_{0}, \g}_{m, p+\s_{2}, \s_{2}}
\]
\end{lemma}
\begin{proof}
%It follows by applying directly Propositions 
%\ref{Teorema totale partial vphi beta k D beta k Phi} and \ref{lemma:tame derivate flusso ancora diverso derivate}.
See Appendix $A$ in \cite{BM20}.
\end{proof}

\begin{lemma}\label{lem:momentoflusso}
Under the assumptions of Lemma \ref{lemma:buonaposFlussi} the following holds.\\
$(i)$ The flow $\Phi_a^{\tau}$ in  \eqref{flussosimboloGenerico}
is symplectic according to Def. 
\ref{def simplettica complesse}.\\
%\smallskip
\noindent
$(ii)$ If  the symbol $a$ in \eqref{egogenerator}-\eqref{egogenerator2} is 
\emph{momentum preserving} according to Def. \ref{def:mompressimbo}, then
the flow map $\Phi_a^{\tau}$ satisfies  \eqref{vec.tau},
i.e. is momentum preserving.
\end{lemma}

\begin{proof}
Item $(i)$ follows since the right hand side of  \eqref{flussosimboloGenerico}
is \emph{Hamiltonian} according to Def. 
\ref{operatoreHam}.

Let us check item $(ii)$ by assuming that (see  Lemma \ref{lem:mompressimbolo})
\begin{equation}\label{cappotto3}
a(\vphi-\tV \vs,x,\x)=a(\vphi,x+\vs,\x)\,,
\qquad \forall (\vphi,x)\in \T^{\nu+d}_{*}\,,\;\; \x\in \Gamma^{*}\,,
\end{equation}
 for any $\vs\in \R^d$. 
In order to prove that $\Phi_{a}^{\tau}(\vphi)$
 is momentum preserving
 we need to show that \eqref{vec.tau} is satisfied.
 For any fixed $\vs\in \R^d$ and any $\tau\in[- 1,1]$, we define
 \[
 L^{\tau}:=\Phi_{a}^{\tau}(\vphi-\mathtt{V}\vs)\circ \tau_{\vs} -\tau_{\vs}\circ\Phi_a^{\tau}\,.
 \]
 First of all we note that $L^{0}\equiv0$, since $\Phi_{a}^{0}$ is the identity.
 Moreover, using \eqref{flussosimboloGenerico}, one has
  \[
 \begin{aligned}
 \pa_{\tau}L^{\tau}&=
 \ii \opw\big(a(\tau;\vphi-\mathtt{V}\vs,x,\x)\big)\Phi_{a}^{\tau}(\vphi-\mathtt{V}\vs)\circ \tau_{\vs}
 \\&
 -\ii\opw(a(\tau;\vphi,x+\vs,\x))\tau_{\vs}\circ \Phi_a^{\tau}
 \\&
 \stackrel{\eqref{cappotto3}}{=}
 \ii\opw(a(\tau;\vphi-\mathtt{V}\vs,x,\x))[L^{\tau}]\,,
 \end{aligned}
 \]
 form which we deduce that $L^{\tau}\equiv0$ for any $\tau\in [- 1,1]$.
 This implies the thesis.
\end{proof}

\begin{rmk}\label{rmk:algflussoEgo}
$(i)$ Note that for $m\le0$ Lemma \eqref{lemma:buonaposFlussi} 
holds for $m_{1} + m_{2}=0$.
%We remark that the flow $\Phi_a^{\tau}$ in  \eqref{flussosimboloGenerico}
%is symplectic according to Def. 
%\ref{def simplettica complesse} since the r.h.s. of \ref{flussosimboloGenerico}
%is \emph{Hamiltonian} according to Def. 
%\ref{operatoreHam}.
%Notice also that if $a\in S_{\tV}$ in \eqref{SV}
%the flow map $\Phi_a^{\tau}$ satisfies  \eqref{vec.tau},
%i.e. is momentum preserving.

\noindent
$(ii)$ In the case that $m\leq 0$, 
the well-posedness results holds true also in the case of the flow ${\bf \Psi}^{\tau}$
\begin{equation}\label{flussoMatrice}
\pa_{\tau}{\bf \Psi}^{\tau}=\ii E \opw(M(\vphi,x,\x)){\bf \Psi}^{\tau}\,,\;\;\; {\bf \Psi}^{0}=\id\,,
\end{equation}
where $M\in S^{m}$ is a \emph{matrix} of symbols in the class $S^{m}$ of the form
\eqref{matricedisimboli}-\eqref{matrixSymself}. In this case the flow ${\bf \Psi}^{\tau}$ is also symplectic.
If in addition the entries of the matrix $M$ are momentum preserving symbols (see \eqref{invariantsymbol})
then ${\bf \Psi}^{\tau}$ is momentum preserving as well.
\end{rmk}

In the  rest of this section we collect some results about  how 
pseudo-differential  and ``tame'' operators 
transform under the flow 
map in \eqref{flussosimboloGenerico}.
In particular, in subsection \ref{sec:conjrules} we study the case of generator as in \eqref{egogenerator2},
while in subsection 
\ref{sec:egothm} we study the case \eqref{egogenerator} proving a quantitative version of Egorov Theorem.

\subsection{Conjugation rules under pseudo-differential flows}\label{sec:conjrules}

We start by studying how ``tame'' operators (see Def. \ref{Dksigmatame}), which are smoothing in space,
are conjugated under the flow
\eqref{flussosimboloGenerico}.

 \begin{lemma}{\bf (Conjugation of smoothing operators).}\label{conj.NEWsmoothresto}
Let $a=a(\vphi, x, \xi)\in S^{m}$, 
$ m\leq 1$ 
with  $a=a(\lambda, i(\lambda))$  
 $k_{0}$-times 
differentiable in $\lambda\in\mathtt{\Lambda}_{0}$
and Lipschitz in the variable $i$.
Assume also that $a$ has the form \eqref{egogenerator}, for $m=1$, or \eqref{egogenerator2} for $m<1$
and let  
$\Phi_{a}$ be the 1-time flow of \eqref{flussosimboloGenerico}.
Fix $\beta_0>0$, 
$\tc\geq m(\beta_0 + k_{0})$
%$\alpha\geq0$ and $N\in \Z$ with \red{$3+\alpha(k_0+\beta_0)\leq N$}
and $M\geq \tc$.
Consider an operator $\cR$ with the following properties:
\begin{itemize}
\item[(a)] %For any $M\leq N-\alpha(\beta_0+k_0)-2$ and 
For any
$\beta\in\Z^{\nu}$ with $|\beta|\le \beta_{0}$,
the operator
$\partial_{\vphi}^{\beta}\jap{D}^{M+m(\beta_{0}+k_{0})} \cR \jap{D}^{m(\beta_0+k_{0})-\tc}$ is $\cD^{k_{0}}$-tame
(see Def. \ref{Dksigmatame}) with tame constants bounded for $ \frac{d + \nu}{2} < s_0\leq s\leq \bar{s}$ for some $\bar{s}\gg1$
large;

\item[(b)] $\cR=\cR(i)$ depends  in a Lipschitz way on a  parameter $i$.  For  %$M \le N-\beta_{0}-2$,
any $\beta\in\Z^{\nu}$ with $|\beta|\le \beta_{0}$, and some  $s_0\leq p\leq \bar{s}$
one has
$\jap{D}^{M+m\beta_{0}} \partial_{\vphi}^{\beta}\Delta_{12}\cR  \jap{D}^{ m\beta_0-\tc}\in\cL(H^{p}, H^{p})$.

\end{itemize}

\noindent
Then there exist constants $\s_{1}=\s_{1}(\beta_{0}, M,  k_{0}, \tc)\geq \s_{2}=\s_{2}(\beta_{0}, M, k_{0}, \tc)$ 
such that if $\bar{s}\geq s_0+\s_1$
there is $\delta=\delta(M,\bar{s})$ such that if
\begin{equation}\label{tazza.rob}
\|a\|^{k_{0} , \g}_{m, s_{0}+\s_{1}, \s_{1}}\leq \delta\,,
\end{equation}
then the operator  $\cQ:= \Phi_{a}^{-1} \cR \Phi_{a}$ satisfies the following.

\begin{itemize}
\item
For 
% \red{$M'\le N-(\alpha+m)( \beta_{0} + k_{0})-2$} and 
any $\beta\in\Z^{\nu}$ 
 with $|\beta|\le \beta_{0}$, one has that the operator 
 $\jap{D}^{M} \partial_{\vphi}^{\beta} \cQ \jap{D}^{-\tc}$ is 
 $\cD^{k_{0}}$-tame and satisfies, for any $s_0 \leq s \leq \bar{s}-\s_1$, the estimate
\begin{equation}\label{stima:NEWconiugatoresto1}
\begin{aligned}
{\mathfrak M}_{\jap{D}^{ M} \partial_{\vphi}^{\beta} \cQ \jap{D}^{-\tc}}(s) 
\lesssim_{m, s, M, \beta_0} 
{\mathbb N}_{\mathcal R}(s, \beta_0) 
+ \| a \|_{ m, s + \sigma_{1},  \sigma_{1}}^{k_0, \gamma}  
{\mathbb N}_{\mathcal R}(s_0, \beta_0)\,,
\end{aligned}
\end{equation}
where 
${\mathbb N}_{\mathcal R}(s, \beta_0) := 
{\rm max}\big\{ {\mathfrak M}_{ \langle D \rangle^{M+m(\beta_0 + k_{0}) }\partial_\vphi^\beta
{\mathcal R} \langle D \rangle^{m(\beta_0 + k_{0} )- \tc}}(s) : |\beta| \leq \beta_0  \big\}$;

\item for 
% $M'\le N- (\alpha+m)\beta_{0}-3$ and 
any 
$|\beta|\le \beta_{0}$ we have that, for any $p+\s_2\leq s_0+\s_1$,
\begin{equation}\label{stima:NEWdelta12coniugatoresto1}
\begin{aligned}
&\|\jap{D}^{M} \pa_\vphi^\beta  \Delta_{12} \cQ \jap{D}^{-\tc }\|_{\mathcal{L}(H^p;H^p)}
\\
&\qquad\qquad
\lesssim_{m, p,  \beta_0,M} 
\|\Delta_{12} a\|_{m, p+\s_{2}, \s_{2}}
\sup_{j=1,2}\|\jap{D}^{M+m(\beta_0 + 1)} 
\pa_\vphi^\beta  \cR(i_{j}) \jap{D}^{m(\beta_0 + 1 )- \tc}\|_{\mathcal{L}(H^p;H^p)}
\\&\qquad\qquad\qquad\quad
+ (1+ \sup_{j=1,2}\| a(i_{j})\|^{ k_{0}, \g}_{m, p+\s_{2}, \s_{2}}) 
\|\jap{D}^{M+m\beta_0} \pa_\vphi^\beta  \Delta_{12} \cR \jap{D}^{m\beta_0 - \tc}\|_{\mathcal{L}(H^p;H^p)}\,.
\end{aligned}
\end{equation}

\end{itemize}
If $m\leq 0$ the statements holds with $m=0$.
\end{lemma}
\begin{proof}
The proof is similar to the one of Lemma $B.10$ in \cite{FGP19}.
By Leibniz rule we have that
\begin{equation*}
\begin{aligned}
 \jap{D}^{M}  \pa_{\vphi}^\beta \Phi_{a}^{-1} \cR \Phi_{a}&\jap{D}^{-\tc} 
 \\&= \sum_{\beta_1+\beta_2+\beta_2=\beta} C_{\beta_{1}, \beta_{2},\beta_{3}}
 \jap{D}^{M}\pa_{\vphi}^{\beta_1}\Phi_{a}^{-1}\pa_{\vphi}^{\beta_2}\cR\pa_{\vphi}^{\beta_3}\Phi_{a}
 \jap{D}^{-\tc}
 \\&=
 \sum_{\beta_1+\beta_2+\beta_2=\beta} C_{\beta_{1}, \beta_{2},\beta_{3}}
 A_{\beta_1}\circ A_{\beta_2}\circ A_{\beta_3}\,,
\end{aligned}
\end{equation*}
where we defined
and note that
\begin{equation*}
\begin{aligned}
A_{\beta_1}&:= \jap{D}^{M}\pa_{\vphi}^{\beta_1}\Phi_{a}^{-1}\jap{D}^{-M-m(\beta_0+k_0)}\,,
\qquad
A_{\beta_2}:=\jap{D}^{M+m(\beta_0+k_0)}\pa_{\vphi}^{\beta_2}\cR\jap{D}^{-\tc + m(\beta_0+k_0)}\,,
\\
A_{\beta_3}&:= \jap{D}^{\tc-m(\beta_0+k_0)}\pa_{\vphi}^{\beta_3}\Phi_{a}\jap{D}^{-\tc}\,.
\end{aligned}
\end{equation*}
%First of all note that, by the choice of $M'$, one has $M'+m(\beta_0+k_0)=M\leq N-\alpha(\beta_0+k_0)-2$. Then
The operator $A_{\beta_2}$ is controlled thanks to the assumption $(a)$ on $\mathcal{R}$.
Moreover, we apply  Lemma \ref{lemma:buonaposFlussi} to $A_{\beta_1}$, 
with  $m_1\rightsquigarrow -M$, $m_2\rightsquigarrow M+m(\beta_0+k_0)$ (which verify the assumption $m_1+m_2\geq m(\beta_0+k_0)$). Similarly we shall use Lemma \ref{lemma:buonaposFlussi}
on $A_{\beta_3}$ with $m_1\rightsquigarrow -\tc+m(\beta_0+k_0)$
and $m_2\rightsquigarrow \tc$
(which again verify the assumption $m_1+m_2\geq m(\beta_0+k_0)$). Then estimate 
\eqref{stima:NEWconiugatoresto1} on the operator $A_{\beta_1}\circ A_{\beta_2}\circ A_{\beta_3}$ follows by composition Lemma \ref{composizione operatori tame AB} and the smallness condition \ref{tazza.rob}.
The estimate \eqref{stima:NEWdelta12coniugatoresto1} follows reasoning as above using assumption $(b)$
and the estimates on the Lipschitz variation of $\Phi_a$ in Lemma \ref{lemma:buonaposFlussi}.
This concludes the proof.
\end{proof}

\begin{rmk}
Notice that in the case $0<m\leq 1$, the conjugate operator $\mathcal{Q}$ in Lemma \ref{conj.NEWsmoothresto}
is less regular than the original one. Indeed, $M<M+m\beta_0$ 
and moreover we have more loss of derivatives coming from
the time derivatives $\pa_{\vphi}$.
On the contrary, when the generator of the flow $\Phi_a$ in \eqref{flussosimboloGenerico} is bounded (case $m\leq 0$)
we have no extra loss in the conjugation procedure. 
 We will only apply Lemma \ref{conj.NEWsmoothresto} with 
$0<m\leq 1$ in subsections \ref{sec:almoststraightening} and \ref{sec:riduzione diagonal ordine 1 2}
a \emph{finite} number of times. 
In performing the infinitely many change of coordinates in the KAM procedure of section 
\ref{sec:KAMreducibility} we will instead use always Lemma \ref{conj.NEWsmoothresto} with $m\leq 0$.
\end{rmk}

\begin{rmk}\label{rmk:restimatrici}
We remark that the result of Lemma \ref{conj.NEWsmoothresto} holds (following word by word its proof)
in the following cases:

$(i)$ for the composition $\Phi_a\circ \mathcal{R}\circ \Phi_{g}$ where $\Phi_{g}$ is the flow (as in \eqref{flussosimboloGenerico}) 
generated by 
$\opw(g)$, with $g\in S^{m}$ a symbol satisfying the same assumptions of the symbol $a$.
This property will be used in sub-section \ref{sec:riduzione diagonal ordine 1 2} to deal with the off-diagonal term of a matrix 
of smoothing operators.

$(ii)$ for  the composition ${\bf \Psi}^{-1}\circ \mathcal{Q}\circ {\bf \Psi}$ 
where $\mathcal{Q}$ is a $2\times 2$
matrix of operators satisfying the properties of $\mathcal{R} $ in Lemma \ref{conj.NEWsmoothresto}
and where ${\bf \Psi}$ is the time one flow of \eqref{flussoMatrice}
generated by a matrix of symbols $M$ in the class $S^{m}$ with $m\leq 0$.
This property will be used in sub-section \ref{subsec:loweroffdiag}.
\end{rmk}

We now study the conjugate of pseudo-differential operators
under the flow in \eqref{flussosimboloGenerico}. More precisely, we consider operators of the form
\begin{equation}\label{operatoreProvaconiugio}
\mathcal{L}=\omega\cdot\pa_{\vphi}+\mathtt{m}\cdot \nabla+\opw(b)\,,
\end{equation}
where  
\begin{itemize}
\item  $b\in S^{n}$, 
$ n\in \R$ 
with  $b=b(\lambda, i(\lambda))$  
 $k_{0}$-times 
differentiable in $\lambda = (\omega, \mathtt h)\in\mathtt{\Lambda}_{0}$
and Lipschitz in the variable $i$. 

\item $\mathtt{m}=\mathtt{m}(\lambda,i(\lambda))\in \R^d$ independent of  
$(\vphi,x,\x)\in \T^{\nu}\times\T_{\Gamma}^{d}\times\R^{d}$, $k_0$-times differentiable 
 in $\lambda = (\omega, \mathtt h)\in\mathtt{\Lambda}_{0}$
and Lipschitz in the variable $i$, 
with $|\mathtt{m}|^{k_0,\gamma}\lesssim1$.
\end{itemize}

One has the following.

\begin{lemma}\label{flussi coniugi eccetera}
Let $a=a(\vphi, x, \xi)\in S^{m}$, 
$ m< 1$ 
with  $a=a(\lambda, i(\lambda))$  
 $k_{0}$-times 
differentiable in $\lambda\in\mathtt{\Lambda}_{0}$
and Lipschitz in the variable $i$ and real valued (see \eqref{egogenerator2}). Let 
$\Phi_{a}$ be the 1-time flow of \eqref{flussosimboloGenerico}
and consider an operator $\mathcal{L}$ as in \eqref{operatoreProvaconiugio}.
Fix $\beta_0 \in \N, \tc>0, M>\tc$ with
$M-\tc- m (\beta_0+k_{0})\geq 3$ and  $\al_{*}\in \N$, $s_0 > \frac{\nu + d}{2}$. 
There exist constants
 $\s_{1}=\s_{1}(\beta_{0}, M, \tc,  k_{0},\alpha_*)\geq \s_{2}=\s_{2}(\beta_{0}, M, \tc,\al_{*})$ 
such that for any $\bar{s}\geq s_0+\s_1$ the following holds. 
There exists $\delta=\delta(M, \tc,\beta_{0}, n, m,\bar{s},\alpha_*)$ 
such that if 
\begin{equation}\label{liberitutti4}
\| a \|_{m, s_0 + \sigma_{1},\alpha+ \sigma_{1}}^{k_0, \gamma}\le \delta \,, 
\end{equation}
then 
\[
\begin{aligned}
\Phi_a^{-1} \circ \,\mathcal{L} \circ \Phi_a=\omega\cdot\pa_{\vphi}+\mathtt{m}\cdot\nabla
&+\ii\opw \big( \omega \cdot \partial_\vphi a + \mathtt m \cdot \nabla a\big)
\\&-\ii\big[\opw(a),\opw(b)\big]
+\opw(c_M)+\mathcal{R}_{M}\,,
\end{aligned}
\]
where the symbol $c_{M}$ and the remainder ${\mathcal R}_M$ satisfy the following:
\\
\noindent
{\bf (Symbol).}
$c_M(\vphi, x, \xi)\in S^{n + 2m - 2}$  
satisfies, for 
 $s\ge s_{0}$, $0 \le \alpha \le \al_{*}$ and $p+\s_2\leq s_0+\s_1$,
\begin{equation}\label{stimaCNN}
\begin{aligned}
\| c_M \|_{n + 2m - 2, s, \alpha}^{k_0, \gamma}  
&\lesssim_{n, m, s, \alpha, M} 
\| a \|_{m, s + \sigma_{1}, \alpha + \sigma_{1}}^{k_0, \gamma}  
\| b \|_{n, s_0 + \sigma_{1}, \alpha + \sigma_{1}}^{k_0, \gamma}   
\\&\qquad\qquad\quad\qquad\qquad
+ \| a \|_{m, s_0 + \sigma_{1}, \alpha + \sigma_{1}}^{k_0, \gamma} 
\| b \|_{n, s + \sigma_{1}, \alpha + \sigma_{1}}^{k_0, \gamma}\,,
\\
\| \Delta_{12} c_M \|_{n + 2m - 2, p, \alpha}^{k_0, \gamma} 
 &\lesssim_{n, m, p, \alpha, M} 
 \sup_{j=1,2}\| a(i_{j}) \|_{m, p+ \sigma_{2}, \alpha + \sigma_{2}}^{k_0, \gamma} 
(\| \Delta_{12} b \|_{n, p + \sigma_{2}, \alpha + \sigma_{2}}+| \Delta_{12} \mathtt{m} |)
 \\&+
\| \Delta_{12} a \|_{m, p + \sigma_{2}, \alpha + \sigma_{2}} 
\sup_{j=1,2}(\| b(i_{j}) \|_{n, p+ \sigma_{2}, \alpha + \sigma_{2}}^{k_0, \gamma}+|\mathtt{m}(i_{j})|^{k_0,\gamma})\,.
\end{aligned}
\end{equation}

\noindent
{\bf (Remainder).}
 For
%{$M \leq N- m (\beta_0+k_{0})-2$ and $|\beta| \leq \beta_0$ da cambiare}
$|\beta| \leq \beta_0$ 
the operator 
$\langle D \rangle^{M} \partial_\vphi^\beta {\mathcal R}_M \langle D \rangle^{- \tc}$ 
%$\langle D \rangle^{M} \partial_\vphi^\beta {\mathcal R}_N \langle D \rangle^{-m(\beta_0+k_0)-2}$ 
is ${\mathcal D}^{k_0}$-tame with tame constant satisfying for every $s_0 \leq s \leq \bar{s}-\s_1$
\begin{equation}\label{restoconiugazione4}
\begin{aligned}
{\mathfrak M}_{\langle D \rangle^{M} \partial_\vphi^\beta 
{\mathcal R}_M \langle D \rangle^{-\tc}}(s) 
&\lesssim_{n, m, s, M}   
\| a \|_{m, s +{\sigma}_{1},  {\sigma}_{1}}^{k_0, \gamma} 
\| b \|_{n, s_0 +{\sigma}_{1},  {\sigma}_{1}}^{k_0, \gamma}  
\\&\qquad\qquad\qquad+ \| a \|_{m, s_0 + {\sigma}_1, {\sigma}_{1}}^{k_0, \gamma}  
\| b \|_{n, s + {\sigma}_{1},  {\sigma}_{1}}^{k_0, \gamma}\,,
\end{aligned}
\end{equation}
and for any
%,  for $M'\le M-\tc - m\beta_{0}-3$ \red{check} and any 
$|\beta|\le \beta_{0}$ we have that, for any $p+\s_2\leq s_0+\s_1$,
\[
\begin{aligned}
%\|\jap{D}^{M'} \pa_\vphi^\beta  \Delta_{12} &\mathcal{R}_{N} 
%\jap{D}^{-(\alpha+m)\beta_0  -2 }\|_{\mathcal{L}(H^p;H^p)}
\|\jap{D}^{M} \pa_\vphi^\beta  \Delta_{12} &\mathcal{R}_{M} 
\jap{D}^{-\tc }\|_{\mathcal{L}(H^p;H^p)}
\\&
\lesssim_{ \beta_0, m,p,M} 
 \sup_{i=1,2}\| a(i_{j}) \|_{m, p+ \sigma_{2}, \alpha + \sigma_{2}}^{k_0, \gamma} 
(\| \Delta_{12} b \|_{n, p + \sigma_{2}, \alpha + \sigma_{2}}+| \Delta_{12} \mathtt{m} |)
 \\&\qquad\qquad+
\| \Delta_{12} a \|_{m, p + \sigma_{2}, \alpha + \sigma_{2}} 
\sup_{j=1,2}(\| b(i_{j}) \|_{n, p+ \sigma_{2}, \alpha + \sigma_{2}}^{k_0, \gamma}
+|\mathtt{m}(i_{j})|^{k_0,\gamma})\,.
\end{aligned}
\]
If $m\leq 0$ the statements holds with $m=0$.
Finally, if the symbols $a,b$ are momentum preserving (see Def. \ref{def:mompressimbo}), then also
$c_{M}$ and ${\mathcal R}_M$ are momentum preserving.
\end{lemma}

\begin{proof}
We shall write
\begin{equation}\label{fogliodicarta}
\Phi_a^{- 1} \circ \mathcal{L} \circ \Phi_a=
\underbrace{\Phi_a^{- 1} \circ \omega\cdot\pa_{\vphi} \circ \Phi_a}_{I}
+\underbrace{\Phi_a^{- 1} \circ (\mathtt{m}\cdot\nabla) \circ \Phi_a}_{II}+
\underbrace{\Phi_a^{- 1} \circ \opw(b) \circ \Phi_a}_{III}\,.
\end{equation}
We study each term separately.

\noindent
$(i)$ Let us study the operator $III$. Recalling \eqref{liberitutti3} and using
the usual Lie expansion up to order $L$ we get
\begin{align}
III&=\opw(b)-\ii\big[\opw(a),\opw(b)\big]
+\sum_{q=2}^{L}\frac{(-1)^q}{q!}{\rm ad}_{\ii\opw(a) }^{q}[\opw(b)]
\nonumber%\label{Lieuffa2}
 \\&+  \frac{(-1)^{L+1}}{L!}
  \int_{0}^{1}  (1- \theta)^{L} {\Phi}_{a}^{\theta} 
  {\rm ad}_{\ii\opw(a) }^{L+1}[\opw(b) ] 
  ({ \Phi}_{a}^{\theta})^{-1} 
d \theta\,.
\label{Lieuffa3}
\end{align}
We choose $L$  in such a way that  
$  (L+1)(1-m) - n \geq M- \tc + 2m(k_{0} + \beta_{0})$.  
In view of Lemma \ref{lemma Ad Pseudo diff} we have that 
\[
\sum_{q=2}^{L}\frac{(-1)^q}{q!}{\rm ad}_{\ii\opw(a) }^{q}[\opw(b)]=:
\opw(c_{N}^{(1)})\,,\qquad c_{N}^{(1)}\in S^{n+{2}m-2}\,.
\]
In particular, by estimate \eqref{stima Ad pseudo diff},
we deduce that $c_{N}^{(1)}$ satisfies \eqref{stimaCNN}.
The estimates on the Lipschitz variation follow similarly recalling Remark \ref{rmk:deltaunodueop}.
Reasoning in the same way,
 we deduce that 
\[
Q:={\rm ad}_{\ii\opw(a) }^{L+1}[\opw(b) ] \in OPS^{-N}\,, \quad N:= M - \tc + 2m(\beta_{0} + k_{0})\,,
\]
with symbol satisfying an estimate like \eqref{stimaCNN} with $n+m-2\rightsquigarrow -N$.
Moreover by Lemma   \ref{constantitamesimbolo}, applied on $Q$, 
we have that for any $|\beta|\leq \beta_0$, $\langle D\rangle^{M+m( \beta_{0}+k_{0})}\pa_{\vphi}^{\beta}Q\langle D\rangle^{-\tc+m( \beta_{0}+k_{0})}$ is 
$\mathcal{D}^{k_0}$-tame 
with tame constant bounded form above by the r.h.s. of \eqref{stimaCNN}.
In conclusion, by 
Lemma \ref{conj.NEWsmoothresto} %(applied with $\alpha=0$)
we obtain that the term in \eqref{Lieuffa3} is a remainder satisfying \eqref{restoconiugazione4}.

\smallskip
\noindent
$(ii)$ Recalling \eqref{fogliodicarta} and using again Lie series expansion we shall write
\begin{align*}
I+II&=\omega\cdot{\pa_{\vphi}}+ \opw(\ii \omega\cdot{\pa_{\vphi}} a )
+ \opw( \ii \mathtt{m}\cdot\nabla a ) 
\\&
+\sum_{k=2}^{L}\frac{(-1)^{k}}{k!}{\rm ad}^{k-1}_{\opw(\ii a)}
[\opw(\ii\omega\cdot\pa_{\vphi} a +\ii\mathtt{m}\cdot\nabla a)]
\\
& +\frac{(-1)^{L+1}}{L!}\int_{0}^{1} (1-\theta)^{L} \Phi_{a}^{\theta}
\Big({\rm ad}^{L}_{\opw(\ii a)}[\opw(\ii\omega\cdot\pa_{\vphi}  a+\ii\mathtt{m}\cdot\nabla a)] \Big)(\Phi_{a}^{\theta})^{-1} d\theta\,.
\end{align*}
Hence one gets the thesis reasoning exactly as done for the term $III$.

Finally, if $a$ and $b$ are momentum preserving, by Remark \ref{rmk:algsimboli} 
one can deduce that also the final symbol $c_{M}$
is momentum preserving and purely imaginary.
\end{proof}

\begin{rmk}
We remark that we shall apply Lemma \ref{flussi coniugi eccetera}
both in the case $\mathtt{m}\equiv0$ or ${b}\equiv0$.
\end{rmk}

\subsection{A quantitative Egorov Theorem}\label{sec:egothm}
In this subsection we consider the flow \eqref{flussosimboloGenerico} in the case 
$m=1$.
More precisely, we consider 
the flow of the equation 
\begin{equation}\label{gianduiotto}
\begin{cases}
\partial_\tau \breve\cA^\tau = \opw (\im \breve{b}(\tau; \varphi, x)\cdot\xi)\breve\cA^\tau 
\\
\breve\cA^0= {\rm Id}
\end{cases}
\end{equation}
with $\breve{b}(\tau; \varphi, x)$ defined in \eqref{egogenerator}.
By an explicit computation is easy to check that 
\begin{align}
\partial_\tau \breve\cA^\tau u = \breve{b}(\tau; \varphi, x) \pa_x \breve\cA^\tau u 
&
+ \frac{1}{2} {\rm div}(\breve{b}(\tau; \varphi, x))\breve\cA^\tau u=:\op(c(\tau; \varphi, x, \xi))
\breve\cA^\tau u\,, \nonumber
\\
c(\tau; \varphi, x, \xi)&:= \im \breve{b}(\tau; \varphi, x)\cdot\xi + \frac{1}{2} {\rm div}(\breve{b}(\tau; \varphi, x))\,.\label{zucchero}
\end{align}
In particular, we have explicit expression for the flow of \eqref{gianduiotto}. Indeed one has
\begin{align}\label{flussocorretto}
\breve\cA^{\tau} h(\varphi, x)&:= \cM^\tau \circ T^{\tau}_{\breve\beta} h(\varphi, x) \,, 
\quad &&
\varphi\in\T^\nu, x\in\T_{\Gamma}^d\,,
\\
(\breve\cA^{\tau})^{-1} h(\varphi, y)&:= \widetilde{\cM^\tau} \circ T^{\tau}_{{\beta}}h(\varphi, x)  
 \,,  &&\varphi\in\T^\nu, y\in\T_{\Gamma}^d\,,\label{belgio}
\end{align}
for $\tau\in [0,1]$,
where 
\begin{align*}
\cM^\tau h(\varphi, x) &= \sqrt{\det(\uno+\tau\nabla_x \breve\beta(\varphi, x))} h(\varphi, x)\,,
\\
\widetilde{\cM^\tau} h(\varphi, x) &= \sqrt{\det(\uno+\tau\nabla_y {\beta}(\varphi, y))} h(\varphi, x)
\end{align*}
and 
\begin{align}
T_{\breve\beta}^{\tau} h(\varphi, x) &= h(\varphi, x+ \tau\breve\beta(\varphi, x))\label{def:taubeta}
\qquad 
T_{{\beta}}^{\tau} h(\varphi, x) = h(\varphi, x+ \tau{\beta}(\varphi, x))\,,
\end{align} 
where $\breve\beta$ is some smooth vector-valued function and 
${\beta}$ is such that\footnote{We remark that, in view of \eqref{diffeoinv} we can write \eqref{belgio} as
\begin{equation}\label{oralegale}
\breve{\mathcal{A}}^{-1}=T_{{\beta}}\circ \frac{1}{\sqrt{\det(\uno+\tau\nabla_x \breve\beta(\varphi, x))}\,.
}
\end{equation}}
\begin{equation}\label{diffeoinv}
x\mapsto y=x+ \tau\breve\beta(\varphi, x) \iff y\mapsto x= y+{\beta}(\tau; \varphi, y)\,, \quad \tau\in[0,1]\,.
\end{equation}
It is easy to check that the map \eqref{flussocorretto} 
is the flow of the equation \eqref{gianduiotto}. 
The well-posedness of $\breve{\mathcal{A}}^{\tau}$ is guaranteed by Lemma \ref{lemma:buonaposFlussi}.

We now discuss the conjugation of a pseudo-differential operator
under the flow \eqref{gianduiotto}.

\noindent
{\bf Notation.} Consider an integer $n\in \N$.
To simplify the notation for now on we shall write, 
$\Sigma^{*}_{n}$ the sum over indexes $k_1,k_2,k_3\in\N$ such that 
$k_1<n$, $k_1+k_2+k_3=n$ and $k_1+k_2\geq 1$.

\begin{thm}{\bf (Egorov Theorem).}\label{quantitativeegorov}
Fix $\tb \in \N$, $m\in \R$, $\tc\geq \tb + k_{0}+2$ and $M\geq \tc$, $s_0 := \frac{\nu + d}{2} + 1$.
%with $M- \tc- (\tb+k_{0})-2\geq m$ \red{check}.
Let $w(\varphi, x, \xi)\in S^m$ with 
$w=w(\lambda, i(\lambda))$, $k_0$-times differentiable  in 
$\lambda\in\mathtt{\Lambda}_{0}$ 
and Lipschitz in the variable $i$. 
Let $\breve\cA^{\tau}$, $\tau\in[0,1]$ be the flow of the system \eqref{gianduiotto}. 
There exist $\sigma_1:=\sigma_1(m, M, \tc,\tb,k_0)\geq \s_{2}:=\s_{2}(m, M, \tc,\tb)$  
for any $\bar{s}\geq s_0+\s_1$ the following holds.
There exists $\delta:=\delta(m, M, \tc,\tb,k_0,\bar{s})$ such that, if
\begin{equation}\label{buf}
\|\breve\beta\|^{k_0, \gamma}_{s_0+\sigma_1} \leq\delta\,,  
\end{equation}
then 
\begin{equation}\label{simbotransportato}
\breve\cA^\tau \opw (w(\varphi, x, \xi))(\breve\cA^\tau)^{-1}_{|\tau=1}= \opw(q(\varphi, x, \xi)) + R\,,
\end{equation}
where $q\in S^m$ and in particular 
\begin{equation}\label{formaEsplic}
q(\varphi, x, \xi)=q_0(\varphi, x,\xi)+q_1(\varphi, x,\xi)\,,
\end{equation}
where $q_0\in S^{m}$, $q_1\in S^{m-1}$ and  (recall \eqref{diffeoinv})
\begin{equation}\label{formaEsplic2}
q_0(\varphi, x,\xi)=
w\big(\vphi,x+\breve\beta(\varphi, x), (1+\nabla_y{\beta}(1, \varphi, y))^{T}_{|y=x+\breve\beta(\varphi, x)}\xi  \big)\,.
\end{equation}
Moreover,  $s\ge s_{0}$, $\alpha\in\N$
%$0 \le \alpha \le \al_{*}$ 
and $p+\s_2\leq s_0+\s_1$, one has that
\begin{equation}\label{parlare}
\|q_i\|_{m-i,s,\alpha}^{k_{0}, \g}
\lesssim_{m,s, \alpha, M}
(1- i) \|w\|_{m,s,\alpha+\sigma_1}^{k_{0}, \g} + \sum_{s}^* \|w\|_{m,k_1,\alpha+k_2+\sigma_1}^{k_{0}, \g} \|\breve\beta\|^{k_{0}, \g}_{k_3+\sigma_1}\,,
\end{equation}
\begin{multline}\label{troppo}
\|\Delta_{12}q_{i}\|_{m-i,p,\alpha}
 \lesssim_{m,p, \alpha, M} 
(1-i)\big(\sup_{j=1,2}\|w(i_{j})\|_{m,p+1,\alpha+\sigma_2} \|\Delta_{12} \breve\beta\|_{p+1} 
+ \|\Delta_{12} w\|_{m, p, \alpha+\sigma_2}\big)
\\
 + \sum_{p+1}^* \sup_{j=1,2}\|w(i_{j})\|_{m,k_1,\alpha+k_2+\sigma_2} 
 \sup_{j=1,2}\|\breve\beta(i_{j})\|_{k_3+\sigma_2} 
 \|\Delta_{12}\breve\beta \|_{s_0+1} 
 \\
 +  \sum_{p}^* 
 \|\Delta_{12}w\|_{m,k_1,\alpha+k_2+\sigma_2} 
 \sup_{j=1,2}\|\breve\beta(i_{j})\|_{k_3+\sigma_2}\,,
\end{multline}
for $i=0,1$.
Furthermore, the operator $R$ satisfies the following: 
 
 \noindent
$\bullet$ For %any $M\leq N-(\tb+k_0)-2$ \red{check} and 
any $b\in\N^{\nu}$ with
$|b|\leq \tb$, 
the operator $\jap{D}^{M} \pa_\vphi^{b} R \jap{D}^{-\tc}$ 
is $\cD^{k_0}$-tame with tame constants satisfying, for $s_0\leq s\leq \bar{s}-\s_1$,
%(recall \eqref{Mdritto})
\begin{equation}\label{francia1}
\begin{aligned}
%\mathbb{M}_R(s,\tb) 
&\mathfrak{M}_{\jap{D}^{M} \pa_\vphi^{b} R \jap{D}^{-\tc}}(s)
\lesssim_{m,s, M} 
%\|w\|_{m,s+\rho,\sigma_1}^{k_0, \gamma} + 
\sum_{s+\s_1}^{*} \|w\|_{m,k_1,\al+k_2+\sigma_1}^{k_0, \gamma} 
\|\breve\beta\|_{k_3+\sigma_1}^{k_0, \gamma}\,.
\end{aligned}
\end{equation}

\noindent
$\bullet$ %For $M'\le N- \tb-2$ \red{check} and any 
For any $|b|\le \tb$ we have that, for any $p+\s_2\leq s_0+\s_1$,
\begin{equation}\label{francia2}
\begin{aligned}
\|\jap{D}^{M} \pa_\vphi^b &\Delta_{12} R \jap{D}^{-\tc}\|_{\mathcal{L}(H^p;H^p)}
\lesssim_{m, p, M}
 \sum_{p+M-\tc}^* \|\Delta_{12}w\|_{m,k_1, k_2+\sigma_2}
 \sup_{j=1,2} \|\breve\beta(i_{j})\|_{k_3+\sigma_2}
\\
& + \sum_{p+M-\tc}^* \ \sup_{j=1,2}\|w(i_{j})\|_{m,k_1,k_2+\sigma_2} 
\sup_{j=1,2}\|\breve\beta(i_{j})\|_{k_3+\sigma_2} 
  \|\Delta_{12}\breve\beta \|_{s_0+\sigma_2} \,.
\end{aligned}
\end{equation} 
%for $0\le \tb\le \rho-3\,.$
Finally if $\breve\beta$ in \eqref{egogenerator} and $w$ \eqref{simbotransportato} 
are momentum preserving,
then
the final symbol 
$q$ in \eqref{formaEsplic} and the remainder $R$ are momentum preserving.
\end{thm}

In the proof we use  the following lemma proved 
in the Appendix of \cite[Lemma A.7]{FGP19}.
%Given a square matrix $A$, we denote by $A^T$ the transpose matrix and if $A$ 
%is invertible we denote by $A^{- T}$ the inverse of the transpose matrix. 

\begin{lemma}\label{Lemmino}
Let $s_0 := \frac{\nu + d}{2} + 1$, $\breve\beta \in C^\infty (\T^{\nu+d}_{*}; \R^{d})$ satisfy 
$\| \breve\beta \|^{k_0, \gamma}_{2s_0+2}< 1$. 
Then, for any symbol $w\in S^m$, 
\begin{equation}
		A w:=w\Big( x+\breve\beta(x), ({\rm Id} + \nabla_x \breve\beta(x))^{-T}\xi \Big)
	\end{equation}
	is a symbol in $ S^m$ satisfying, for any $ s \geq s_0 $,    
	\begin{equation}\label{stima}
		\| A w \|^{k_0, \gamma}_{m, s, \al}\lesssim_{m,s,\al} \| w \|^{k_0, \gamma}_{m, s, \al}
		+ \sum_{s}^{*}
		\| w \|^{k_0, \gamma}_{m, k_1, \al+k_2} \| \breve\beta \|^{k_0, \gamma}_{k_3+s_0+2} \, .
	\end{equation}
	% for some $ C:=C(s,p)>0$.
	 For $s=s_0$  we have the rougher estimate 
$	\| A w \|^{k_0, \gamma}_{m, s_0,\al}\lesssim \| w \|^{k_0, \gamma}_{m, s_0, \al+s_0} $. 
\end{lemma}
We are now in position to prove the main result of this section.

\begin{proof}[{\bf Proof of Theorem \ref{quantitativeegorov}}]
First of all, in view of \eqref{flussocorretto} and \eqref{oralegale}, we shall write
\[
\begin{aligned}
\breve\cA^\tau \opw (w(\varphi, x, \xi))(\breve\cA^\tau)^{-1}&=
{\cM^\tau} \circ \big(T^{\tau}_{\breve{\beta}}\circ\opw (w(\varphi, x, \xi))\circ T^{\tau}_{{\beta}}\big)\circ ({\cM^\tau} )^{-1}
\\&
={\cM^\tau} \circ \big(T^{\tau}_{\breve{\beta}}\circ\opw (w(\varphi, x, \xi))\circ (T^{\tau}_{\breve{\beta}})^{-1}\big)\circ ({\cM^\tau} )^{-1}\,.
\end{aligned}
\]
Since $\cM^{\tau}$ is just a multiplication operator,
it is easy to check that the thesis follows by proving the theorem only for the conjugate
\[
T^{\tau}_{\breve{\beta}}\circ\opw (w(\varphi, x, \xi))\circ (T^{\tau}_{\breve{\beta}})^{-1}\,.
\]
Moreover, by Lemma \ref{passaggioweylst}, we know that there are symbols 
$b_{<n}\in S^{m}$ and $b_{\geq n}\in S^{m-n}$
satisfying \eqref{forza1}-\eqref{forza2} with $a\rightsquigarrow w$, $N\rightsquigarrow n$
such that $\opw(w)=\op(b_{<n})+\op(b_{\geq_{n}})$.
Therefore
\[
\begin{aligned}
T^{\tau}_{\breve{\beta}}\circ\opw (w(\varphi, x, \xi))\circ (T^{\tau}_{\breve{\beta}})^{-1}&=
T^{\tau}_{\breve{\beta}}\circ\op (b_{<n}(\varphi, x, \xi))\circ (T^{\tau}_{\breve{\beta}})^{-1}+\mathcal{Q}\,,
\\
\mathcal{Q}&:=T^{\tau}_{\breve{\beta}}\circ\op(b_{\geq n}(\varphi, x, \xi))\circ (T^{\tau}_{\breve{\beta}})^{-1}\,.
\end{aligned}
\]
By applying estimates \eqref{forza2}, Lemma \ref{constantitamesimbolo} 
with $n$ sufficiently large w.r.t. $M,\tb,k_0$,
combined with Lemma \ref{conj.NEWsmoothresto},
we deduce that $\mathcal{Q}$ satisfies the estimates \eqref{francia1}-\eqref{francia2}.
We reduced to study only the conjugate of $\op(b_{<n})$ under the map $T_{\breve\beta}^{\tau}$,
which can be done studying the conjugate of $\op(w)$ for any symbol $w$.

We shall follow some arguments which are similar to the ones of Theorem $3.4$ of \cite{FGP19} (see also Thm. 6.1 of \cite{BFPT1}).
Let us consider for $\tau\in[0,1]$ the composition operators
\begin{equation}\label{ignobel}
T_{\breve\beta}^{\tau}h(x):=h(\varphi, x+\tau\breve\beta(x))\,, 
\quad 
(T_{\breve\beta}^{\tau})^{-1}h(\varphi, y):=
h(y+{\beta}(\tau; y))\, , 
\end{equation}
\noindent
For any $\tau\in[0,1]$, the conjugated operator
	\begin{equation}\label{opPtautau}
	P^{\tau}:=T_{\breve\beta}^{\tau} \circ {\rm Op}(w)\circ  (T_{\breve\beta}^{\tau})^{-1}
	\end{equation}
	solves  the Heisenberg equation
	\begin{equation}\label{ars}
	\partial_{\tau} P^{\tau}=[{X}^{\tau}, P^{\tau}]\,,  \quad P^{0}={\rm Op}(w) \, , 
	\end{equation}
	where\footnote{By using the definition of $T_{\breve\beta}^{\tau}$ in \eqref{ignobel}
	and \eqref{def:Calphagenerator}
	one can easily check that
	\begin{equation}\label{flussodiffeo}
      \partial_{\tau}T_{\breve\beta}^{\tau}=X^{\tau}T_{\breve\beta}^{\tau}\,,\quad 
     T_{\breve\beta}^{0}={\rm Id}\,.
	\end{equation}
	} 
	\begin{equation}\label{def:Calphagenerator}
		{X}^{\tau}:=\breve{b}(\tau;\vphi,x) \cdot \nabla_{x}={\rm Op}(\chi)\,,\qquad 
		\chi :=\chi(\tau;\vphi,x,\xi):=\ii \breve{b}(\tau;\vphi,x)\cdot\xi\,,
	\end{equation}
	and $\breve{b}(\tau;\vphi,x)$ is \eqref{egogenerator}.
Notice that
\begin{equation}\label{normachichi}
\begin{aligned}
\|\chi\|_{1,s,\al}^{k_0,\gamma}&\lesssim_{s}\|\breve{b}\|_{s}^{k_0,\gamma}
\lesssim_{s}\|\breve\beta\|_{s+1}^{k_0,\gamma}\,,\qquad \forall \,\alpha\geq0\,,
\\
\|\Delta_{12}\chi\|_{1,p,\al}&\lesssim_{p}\|\Delta_{12}\breve\beta\|_{p+1}\,, \qquad p+1\leq s_0+\s_1\,.
\end{aligned}
\end{equation}
Let us fix 
\begin{equation}\label{sceltarho}
\rho:=M -\tc +4 \big( k_0+\tb\big)+1\,.
\end{equation}
	We look for an approximate solution of \eqref{ars} 
	of the form % as a pseudo-differential operator %${\rm Op}(q)$ with 
	\begin{equation}\label{balconata0}
		Q^{\tau}:={\rm Op}(q(\tau;x,\xi))\,,\quad 
		q=q(\tau;x,\xi)=\sum_{k=0}^{m+\rho-1} q_{m-k}(\tau;x, \xi)\,,
	\end{equation}
	where $q_{m-k}$ are  symbols in $S^{m-k}$ to be determined iteratively	so that
\begin{equation}\label{approssimo}
	\partial_{\tau} Q^{\tau}=[{X}^{\tau}, Q^{\tau}] + \mathcal{M}^\tau\,, \quad Q^{0}={\rm Op}(w) \, , 
\end{equation}
	where $\mathcal{M}^\tau={\rm Op}(\mathtt{r}_{-\rho}(\tau;x,\xi))$ with symbol $\mathtt{r}_{-\rho}\in S^{-\rho}$.
	Passing to the symbols we obtain 
	(recall \eqref{espansionecompostandard}, \eqref{balconata0} and \eqref{Moyalstandard})
	\begin{equation}\label{probapproxsimboloq}
\left\{\begin{aligned}
			&\partial_{\tau}q(\tau;x,\xi)=\chi(\tau;x,\xi)\star q(\tau;x,\xi) + \mathtt{r}_{-\rho}(\tau;x,\xi)
			\\
			&q(0;x,\xi)=w(x,\xi)\,.
		\end{aligned}\right.
	\end{equation}
	where the unknowns are now $q(\tau;x,\xi),\mathtt{r}_{-\rho}(\tau;x,\xi)$.

	We  expand $\chi(\tau;x,\xi)\star q(\tau;x,\xi)$ into a sum of symbols with decreasing orders.
	Using that $\chi$ is linear in $\xi\in\R^{d}$ (see \eqref{def:Calphagenerator}), 
	\eqref{espansionecompostandard} %, \eqref{remainder cal RN composizione},
	and by the expansion \eqref{balconata0} (together with the ansatz that $q_{m-k}\in S^{m-k}$)
	we note that %(using the linearity of $\partial_{x}$ and of $\#$)
	\[
	\begin{aligned}
	\chi\star q&=\chi\#q-q\#\chi=\chi q+\frac{1}{i}(\nabla_{\xi}\chi)\cdot(\nabla_{x}q)-q\#\chi
	\\&
	=
	\sum_{k=0}^{m+\rho-1} \frac{1}{i}\{\chi,q_{m-k}\}
	-\sum_{k=0}^{m+\rho-1}
	\sum_{n=2}^{m-k+\rho}q_{m-k}\#_{n}\chi-\sum_{k=0}^{m+\rho-1}q_{m-k}\#_{\geq m-k+\rho+1}\chi\,.
	\end{aligned}
	\]
	
	By rearranging the sums we can write
\begin{equation}\label{civilwar}
\begin{aligned}
	\chi\star q&=\underbrace{-i\{\chi,q_{m}\}}_{{\rm ord}\; m}
	-\overbrace{\sum_{k=1}^{m+\rho-1}\underbrace{\big(-i\{\chi,q_{m-k}\}+r_{m-k}\big)}_{{\rm ord}\; m-k} }^{{\rm orders\, from}\;-\rho+1 \;{\rm to}\;  m-1}
	-\underbrace{\mathtt{r}_{-\rho}}_{{\rm ord}\;  -\rho}
\end{aligned}
\end{equation}
where we defined, denoting $\mathtt{w}=\mathtt{w}(k,h):=k-h+1$ and recalling again \eqref{espansionecompostandard},
\begin{equation}\label{sperobene}
\begin{aligned}
r_{m-k}&:=
-\sum_{h=0}^{k-1}
q_{m-h}\#_{\mathtt{w}}\chi
%\frac{1}{h! i^{h}} (\partial_{\xi}^\mathtt{w} q_{m-h})(\partial_x^\mathtt{w} \chi)
\in S^{(m-h)+1-(k-h+1)}\equiv S^{m-k}\,,
\end{aligned}
\end{equation}
\begin{equation}\label{sperobeneResto}
\begin{aligned}
\mathtt{r}_{-\rho}&:=
\sum_{k=0}^{m+\rho-1}r_{m-k+1+\rho}(q_{m-k},\chi)\in S^{m-k+1-(m-k+1+\rho)}\equiv S^{-\rho}\,.
\end{aligned}
\end{equation}
We have reduced the problem to  finding symbols $q_{m-k}\in S^{m-k}$, $0\leq k\leq m+\rho-1$  which solve for $k=0$
(recall the form of $\chi$ in \eqref{def:Calphagenerator})
	\begin{equation}\label{ordm}
		\left\{\begin{aligned}
			&\partial_{\tau}q_{m}(\tau;x,\xi)= \{\breve{b}(\tau;x)\xi, q_{m}(\tau;x,\xi)\}
			\\
			&q_{m}(0;x,\xi)=w(x,\xi)\,,
		\end{aligned}\right.
	\end{equation}
while for $1\le k\le m+\rho-1$ 
	\begin{equation}\label{ordmmenok}
	\left\{\begin{aligned}
		&\partial_{\tau}q_{m-k}(\tau;x,\xi)= \{\breve{b}(\tau;x) \cdot \xi, q_{m-k}(\tau;x,\xi)\}+r_{m-k}(\tau;x,\xi)
		\\
		&q_{m-k}(0;x,\xi)=0\,.
	\end{aligned}\right.
\end{equation}
We note that the symbols $r_{m-k}$ in \eqref{sperobene} with $1\leq k\leq m+\rho-1$ depend only
on $q_{m-h}$ with $0\leq h<k$. This means that equations \eqref{ordmmenok} can be solved iteratively.

\vspace{0.5em}
\noindent
{\bf Order $m$.} To solve \eqref{ordm}, 	we consider the solutions of the Hamiltonian system
	\begin{equation}\label{charsys}
		\left\{\begin{aligned}
			&\frac{d}{ds}x(s)=-\breve{b}(s,x(s))\\
			&\frac{d}{ds}\xi(s)= \big(\nabla_{x} \breve{b}(s,x(s))\big)^{T}\xi
		\end{aligned}\right.\qquad (x(0),\xi(0))=(x_0,\xi_0)\in \T_{\Gamma}^{d}\times \R^{d}\,.
	\end{equation}
	One can note that if $q_m$ is a solution of \eqref{ordm}, then 
	it is constant when evaluated along the flow of \eqref{charsys}. In other words
	setting
	$g(\tau):=q_{m}(\tau;x(\tau),\xi(\tau))$ one has that 
	$	\frac{d}{d\tau}g(\tau)=0 $ implies that $ g(\tau)=g(0) $, for any $ \tau\in[0,1] $. 
	Let us denote by $\gamma^{\tau_0,\tau}(x,\xi)$ the solution 
	of the characteristic system \eqref{charsys}
	with initial condition $\gamma^{\tau_0,\tau_0}=(x,\xi)$\footnote{
		In other words $(x(\tau),\xi(\tau))=\gamma^{0,\tau}(x_0,\xi_0)$
		and the inverse flow is given by $\gamma^{\tau,0}(x,\xi)=(x_0,\xi_0)$.
	}.
	Then the equation \eqref{ordm} has the solution
	\begin{equation}\label{ars3}
		q_m(\tau;x, \xi)=w(\gamma^{\tau, 0}(x, \xi))
	\end{equation}
	where $\gamma^{\tau, 0}(x, \xi)$ has the explicit form (recall \eqref{egogenerator})
	\begin{equation}\label{ars40}
	\begin{aligned}
		\gamma^{\tau, 0}(x, \xi)&=\big(f(\tau;x), g(\tau;x)\xi\big)\,,
		\\
		f(\tau;x)&:=x+\tau\breve\beta(x)\,, \quad 
		g(\tau;x):= \left( \mathbb{I} + \tau \nabla_{x} \breve\beta(x) \right)^{-T}\,.
		\end{aligned}
	\end{equation}
	For an explicit computation of the solution above, see \cite{CMT24}.\\
	Hence, by Lemma \ref{Lemmino} applied to $q_{m}(\tau;x,\xi)$ in \eqref{ars3} 
	we have that bound \eqref{parlare} with $i=0$ holds true.
	The bound \eqref{troppo} with $i=0$ on the Lipschitz variation follows 
	using the explicit formula \eqref{ars3}, 
	Leibniz  rule, the tame estimates of Lemma \ref{lemma:LS norms} 
	and using again Lemma \ref{Lemmino}.
	
	\vspace{0.5em}
	\noindent
	{\bf Order $m-k$ with $1\leq k\leq m+\rho-1$.} 
	We now assume inductively that 
	we have already found the appropriate solutions $q_{m-h}\in S^{m-h}$ 
	of \eqref{ordmmenok}
	with $0\leq h<k$, for some $k\ge 1$. Moreover, we assume that 
	\begin{align}
		\| q_{m-h} \|^{k_0,\gamma}_{m-h, s, \al}&\lesssim_{m, s, \al, \rho}  
		\sum_{s}^{*} 
		\| w \|^{k_0,\gamma}_{m, k_1, k_2+\sigma_{h}+\al} 
		\| \breve\beta \|^{k_0,\gamma}_{k_3+\sigma_{h}}\,,\quad 1\leq h<k\,,\label{zeppelinIndut}
		\\
\|\Delta_{12}q_{i}\|_{m-i,p,\alpha} 
&\lesssim_{m,s, \alpha, \rho} 
 \sum_{p}^* \|\Delta_{12}w\|_{m,k_1,\alpha+k_2+\sigma_h}
 \sup_{j=1,2} \|\breve\beta(i_{j})\|_{k_3+\sigma_h}
 \nonumber
\\&+
 \sum_{p+1}^*  \sup_{j=1,2}\|w(i_{j})\|_{m,k_1,\alpha+k_2+\sigma_h} 
 \sup_{j=1,2} \|\breve\beta(i_{j})\|_{k_3+\sigma_h} \|\Delta_{12}\breve\beta \|_{s_0+1} \,,
 \label{zeppelinIndutDelta12}
		\end{align}
	for some non decreasing sequence of parameters $\sigma_{h}$ 
	depending only on $\nu$, $d$, $|m|,\rho$. % and \textcolor{red}{$p\leq p_{*}$}.
	Notice that in \eqref{zeppelinIndut} the Sobolev  norm of  $\breve\beta$ 
	does not contain  the parameter $\al$. This is due to the fact that $\chi$ is linear in $\xi$.
	We now construct the solution $q_{m-k}$ of \eqref{ordmmenok} 
	satisfying estimate \eqref{zeppelinIndut} with $h=k$.
	Using Lemma \ref{lemmacomposizioneSTANDARD} (see \eqref{stimacancellettoesplicitoAlgrammo}) 
	we deduce 
(recall also that $\chi$ is linear in $\xi$ and $\mathtt{w}:=k-h+1$)\footnote{We are applying estimate 
\eqref{stimacancellettoesplicitoAlgrammo}
with $a\rightsquigarrow q_{m-h}\in S^{m-h}$, 
$b\rightsquigarrow \chi\in S^{1}$, $n\rightsquigarrow\mathtt{w}=k-h+1$.}
\begin{equation}\label{stimaResti}
\begin{aligned}
\|r_{m-k}\|_{m-k,s,\al}^{k_0,\gamma}
&\stackrel{\eqref{sperobene},\eqref{normachichi}}{\lesssim_{m,s,\al,\rho}}
\|q_{m}\|_{m,s,\al+k+1}^{k_0,\gamma}
\|\breve\beta\|_{s_0+k+2}^{k_0,\gamma}
+
\|q_m\|_{m,s_0,\al+k+1}^{k_0,\gamma}
\|\breve\beta\|_{s+k+2}^{k_0,\gamma}
\\&+
\sum_{h=1}^{k-1}
\|q_{m-h}\|_{m-h,s,\al+\mathtt{w}}^{k_0,\gamma}
\|\breve\beta\|_{s_0+1+\mathtt{w}}^{k_0,\gamma}
+
\|q_{m-h}\|_{m-h,s_0,\al+\mathtt{w}}^{k_0,\gamma}
\|\breve\beta\|_{s+1+\mathtt{w}}^{k_0,\gamma}\, .
\end{aligned}
\end{equation} 
Let us consider first the second  summand, 
which is the most complicated: by  the inductive assumption \eqref{zeppelinIndut} with $1 \leq h\leq k-1$, 
we deduce %(recall that $\mathtt{w}=k-h+1$)
\[
\begin{aligned}
\sum_{h=1}^{k-1}
\|q_{m-h}\|_{m-h,s,\al+\mathtt{w}}^{k_0,\gamma}
&\|\breve\beta\|_{s_0+1+\mathtt{w}}^{k_0,\gamma}
+
\|q_{m-h}\|_{m-h,s_0,\al+\mathtt{w}}^{k_0,\gamma}
\|\breve\beta\|_{s+1+\mathtt{w}}^{k_0,\gamma}
\\&\lesssim_{m, s, \al,\rho} \sum_{h=0}^{k-1}
\sum_{s}^{*} 
\| w \|^{k_0,\gamma}_{m, k_1, k_2+\sigma_{h}+\al+\mathtt{w}} 
\| \breve\beta \|^{k_0,\gamma}_{k_3+\sigma_{h}+\mathtt{w}}
\\ &\qquad\quad
+
\sum_{h=0}^{k-1} \  \sum_{1\le k_1+k_2 \le s_0}
\| w \|^{k_0,\gamma}_{m, k_1, k_2+\sigma_{h}+\al+\mathtt{w}} 
\|\breve\beta\|_{s-s_0+1+\mathtt{w} +s_0 }^{k_0,\gamma} 
\\&	
\lesssim_{m, s, \al, \rho} 
\sum_{s}^{*} 
\| w \|^{k_0,\gamma}_{m, k_1, k_2+\widehat{\sigma}_{k-1}+\al} 
\|\breve \beta \|^{k_0,\gamma}_{k_3+\widehat{\sigma}_{k-1}} \, . 
\end{aligned}
\]
In the second and third line we have used  the smallness assumption
\eqref{buf} 
(assuming $\sigma\geq {\sigma}_{k-1}+1+ k %\ge \s_h + \mathtt{w} +1
$), and in the last line we have set $\widehat\sigma_{k-1}:= \sigma_{k-1} +k +2 +s_0$.
In dealing with the first summand  in \eqref{stimaResti}
we proceed in the same way, only we substitute \eqref{parlare} instead of \eqref{zeppelinIndut}.
We conclude that 
\begin{equation}
	\label{verdone1}
\|r_{m-k}\|_{m-k,s,\al}^{k_0,\gamma}
{\lesssim}_{m, s, \al, \rho} 
\sum_{s}^{*} 
\| w \|^{k_0,\gamma}_{m, k_1, k_2+\widehat{\sigma}_{k-1}+\al} 
\|\breve \beta \|^{k_0,\gamma}_{k_3+\widehat{\sigma}_{k-1}}\,.
\end{equation}
Using again \eqref{sperobene} we note that
\[
\Delta_{12}r_{m-k}=-\sum_{h=0}^{k-1}
\Big[(\Delta_{12} q_{m-h})\#_{\mathtt{w}}\chi(i_2)+(q_{m-h}(i_1))\#_{\mathtt{w}}\Delta_{12}\chi
%\frac{1}{h! i^{h}} \Big[(\partial_{\xi}^\mathtt{w}\Delta_{12} q_{m-h})(\partial_x^\mathtt{w} \chi)
%+(\partial_{\xi}^\mathtt{w}q_{m-h})(\partial_x^\mathtt{w} \Delta_{12} \chi)
\Big]\,.
\]
Reasoning as in \eqref{stimaResti}-\eqref{verdone1}, using \eqref{normachichi} and the inductive assumption 
\eqref{zeppelinIndut}-\eqref{zeppelinIndutDelta12}, one can check that
\begin{equation}\label{verdone1Delta12}
\begin{aligned}
\|\Delta_{12}r_{m-k}\|_{m-i,p,\alpha} 
&\lesssim_{m,p, \alpha, \rho} 
\sum_{p}^* \|\Delta_{12}w\|_{m,k_1,\alpha+k_2+\widehat{\sigma}_{k-1}} 
\sup_{j=1,2} \|\breve\beta(i_{j})\|_{k_3+\widehat{\sigma}_{k-1}}
\\&+
 \sum_{p+1}^* \sup_{j=1,2}\|w(i_{j})\|_{m,k_1,\alpha+k_2+\widehat{\sigma}_{k-1}} 
 \sup_{j=1,2}\|\breve\beta(i_{j})\|_{k_3+\widehat{\sigma}_{k-1}} 
 \|\Delta_{12}\breve\beta \|_{s_0+1} \,.
\end{aligned}
\end{equation}
	We now reason similarly to  $q_m$, and solve \eqref{ordmmenok} by variation of constants.
	We define $f_{m-k}(\tau):=q_{m-k}(\tau;x(\tau),\xi(\tau))$
	where $x(\tau),\xi(\tau)$ are the solution of the Hamiltonian system \eqref{charsys}.
	One has that, if $q_{m-k}$ solves \eqref{ordmmenok}, then
	\[
	\frac{d}{d\tau}f_{m-k}(\tau)=r_{m-k}(\tau;x(\tau),\xi(\tau))\quad \Rightarrow
	\quad
	f_{m-k}(\tau)=\int_{0}^{\tau} r_{m-k}(\sigma ; x(\sigma),\xi(\sigma))d\sigma\,,
	\]
	where we used that $f(0)=q_{m-k}(0,x(0),\xi(0))=0$.
	Therefore the solution of \eqref{ordmmenok} is
	\begin{equation*}%\label{gnomo}
		q_{m-k}(\tau;x, \xi)=
		\int_0^{\tau} r_{m-k}(\gamma^{0, \sigma} \gamma^{\tau, 0}(x, \xi)) \,d\sigma \, .
	\end{equation*}
	We observe also that 
	\begin{equation*}
		\gamma^{0, \sigma} \gamma^{\tau, 0}(x, \xi)=(\tilde{f}, \tilde{g}\,\xi)\,,\qquad \sigma,\tau\in[0,1]\,,
	\end{equation*}
	with
	\begin{equation*}
		\tilde{f}(\sigma, \tau,x):=x+\tau \breve\beta(x)+{\beta}(\sigma, x+\tau \breve\beta(x))\,, 
		\qquad 
		\tilde{g}(\sigma, \tau,x):= \left( \nabla_{x} f (\sigma, \tau, x) \right)^{-T} \, . 
	\end{equation*}
	Thus if $\tilde{A} r:=r(\tilde{f},\xi\cdot \tilde{g})$ we have (recall that $\tau\in [0, 1]$)
	\begin{equation*}
		\begin{aligned}
			\| q_{m-k} \|^{k_0,\gamma}_{m-k, s, \al}
			&\lesssim_{m,s, \al,\rho} 
			\| \tilde{A} r_{m-k} \|^{k_0,\gamma}_{m-k, s, \al}\,,
			\\
			\| q_{m-k} \|^{k_0,\gamma}_{m-k, s_0, \al}
			&\lesssim_{m,s_0,\al,\rho} 
			\| \tilde{A} r_{m-k} \|^{k_0,\gamma}_{m-k, s_0, \al}
			\lesssim_{m,s_0,\al,\rho} 
			\| r_{m-k} \|^{k_0,\gamma}_{m-k, s_0, \al+s_0}\,,
		\end{aligned}
	\end{equation*}
	and by Lemma \ref{Lemmino} with  
	$\breve\beta\rightsquigarrow  \tau \breve\beta(x)+{\beta}(\sigma, x+\tau \breve\beta(x))$\footnote{Notice that 
	the function $t(\tau,\sigma,x):=\tau \breve\beta(x)+{\beta}(\sigma, x+\tau \breve\beta(x))$
	satisfies, using that $y+{\beta} $ is the inverse diffeomorphism of $x+\tau\breve\beta$,
	%(recall for instance  \textcolor{red}{stima sul diffeo inverso}) ,
	the estimate	$ \| t(\tau,\sigma)\|_{s}^{k_0,\gamma}\lesssim_s  \|  \breve{\beta}\|_{s+s_0}^{k_0,\gamma} $
	uniformly in $\tau,\sigma\in[0,1]$.
	}
	\begin{equation}\label{interpol}
	\begin{aligned}
		\| q_{m-k} \|^{k_0,\gamma}_{m-k, s, \al}&\lesssim_{m, s, \alpha,\rho} 
		\| r_{m-k} \|^{k_0,\gamma}_{m-k, s, \al}
		+\sum_{s}^{*} 
		\| r_{m-k} \|^{k_0,\gamma}_{m-k, k_1, \al+k_2} \| \breve\beta \|^{k_0,\gamma}_{k_3+2s_0+2}\,,
		\end{aligned}
		\end{equation}
\begin{equation}\label{interpolDelta12}	
\begin{aligned}	
		\|\Delta_{12}q_{m-k}\|_{m-i,p,\alpha} 
&\lesssim_{m,p, \alpha, \rho} 
\|\Delta_{12} r_{m-k} \|_{m-k, s, \al}
		\\&+\sum_{p}^{*} 
		\| \Delta_{12} r_{m-k} \|_{m-k, k_1, \al+k_2} \sup_{j=1,2}\|\breve\beta(i_{j}) \|_{k_3+2s_0+2}
		\\&
		+ \sup_{j=1,2}\| r_{m-k}(i_{j}) \|^{\g, \cO}_{m, p+1, \al}
		+ \sum_{p+1}^{*}
		\\&+ \sup_{j=1,2}\|r_{m-k} (i_{j})\|^{\g, \cO}_{m, k_1, \al+k_2} 
		 \sup_{j=1,2}\| \breve\beta(i_{j}) \|^{\g, \cO}_{k_3+s_0+2} 
		\|\Delta_{12}\breve\beta\|_{s_0+1}\,.
 \end{aligned}
	\end{equation}
	It remains to prove that actually the symbol
	$q_{m-k}$ satisfies the bound \eqref{zeppelinIndut} with $h=k$ and for some 
	new $\sigma_{k}\geq \sigma_{k-1}$ depending only on $\nu$, $d$, $|m|$ and $\rho$.
By substituting \eqref{verdone1} in the estimate \eqref{interpol}
we get 
\begin{equation}\label{verdone2}
\begin{aligned}
\| q_{m-k} \|^{k_0,\gamma}_{m-k, s, \alpha}&\lesssim_{m,s, \al,\rho} 
\sum_{s}^{*} 
		\| w \|^{k_0,\gamma}_{m, k_1, k_2+\widehat{\sigma}_{k-1}+\alpha} 
		\| \breve\beta \|^{k_0,\gamma}_{k_3+\widehat{\sigma}_{k-1}}
		\\&+\sum_{s}^{*} 
		\Big(
	\sum_{k_1}^{*} 
		\| w \|^{k_0,\gamma}_{m, k_1', k_2'+\widehat{\sigma}_{k-1}+\al+k_2} 
		\| \breve\beta \|^{k_0,\gamma}_{k_3'+\widehat{\sigma}_{k-1}}
			\Big) \| \breve\beta \|^{k_0,\gamma}_{k_3+2s_0+2}\,,
\end{aligned}
\end{equation}
where the sums run over indexes satisfying $k_1+k_2+k_3=s$ and $k_1'+k_2'+k_3'=k_1$.
By classical interpolation estimates on Sobolev spaces (see Lemma \ref{interpolazione fine})\footnote{applied if $k_3'\neq0$, otherwise we do nothing}, using \eqref{buf}, we deduce that \eqref{verdone2} 
becomes (by renaming the indexes in the sum)
\[
\begin{aligned}
\| q_{m-k} \|^{k_0, \gamma}_{m-k, s, p}&\lesssim_{m,s, p,\rho} 
			\sum_{s}^{*} \| w \|^{k_0,\gamma}_{m, k_1, k_2+{\sigma}_{k}+p} 
		\| \breve\beta \|^{k_0,\gamma}_{k_3+{\sigma}_{k}}\,,
\end{aligned}
\]
for some $\sigma_{k}\geq \widehat{\sigma}_{k-1}+2s_0+2$ depending only on $\nu$, $|m|$, $\rho$.
This is the estimate \eqref{zeppelinIndut} with $h=k$.
The estimate \eqref{zeppelinIndutDelta12} with $h=k$ follows reasoning as above recalling 
\eqref{interpolDelta12} and \eqref{verdone1Delta12}.

Let us now estimate the remainder term in \eqref{sperobeneResto}.
Similarly to \eqref{stimaResti}, 
we have\footnote{Here we use  \eqref{stima:resto composizioneSTANDARD} 
with $a\rightsquigarrow q_{m-k}\in S^{m-k}$, $b\rightsquigarrow\chi\in S^{1}$, $N\rightsquigarrow m-k+1+\rho$.
Notice also that $ 2<N\leq |m|+1+\rho$.}
\begin{equation}\label{stimaRestiRHO}
\begin{aligned}
\|\mathtt{r}_{-\rho}\|_{-\rho,s,\al}^{k_0,\gamma}
&
%\stackrel{\eqref{sperobeneResto}}{\lesssim}
\stackrel{\eqref{sperobeneResto},\eqref{stima:resto composizioneSTANDARD}}{\lesssim_{m,s,\al,\rho}}
\sum_{k=0}^{m+\rho-1}
\|q_{m-k}\|^{k_0,\gamma}_{m-k, s, |m|+1+\rho+\alpha}
 \|\chi\|^{k_0,\gamma}_{1, s_0+3|m|+3+3\rho+\alpha, \alpha} 
\\&
\qquad\quad 
+ \|q_{m-k}\|^{k_0,\gamma}_{m-k, s_0, |m|+1+\rho+\alpha} 
\|\chi\|^{k_0,\gamma}_{1, s+3|m|+3+3\rho+\alpha, \alpha}
\\&
\stackrel{\eqref{normachichi}}{\lesssim_{m,s,\alpha,\rho}}
\sum_{k=0}^{m+\rho-1}
\|q_{m-k}\|^{k_0,\gamma}_{m-k, s, |m|+1+\rho+\alpha} \|\breve\beta\|^{k_0,\gamma}_{s_0+3|m|+4+3\rho+\alpha} 
\\&
\qquad\quad 
+ \|q_{m-k}\|^{k_0,\gamma}_{m-k, s_0, |m|+1+\rho+\alpha} \|\breve\beta\|^{k_0,\gamma}_{s+3|m|+4+3\rho+\alpha}\,.
\end{aligned}
\end{equation} 
Substituting the bound \eqref{zeppelinIndut}, we get
\begin{equation}\label{stimaerrinomenorho}
\begin{aligned}
	\|\mathtt{r}_{-\rho}\|_{-\rho,s,\alpha}^{k_0,\gamma}
\lesssim_{m,s,\alpha,\rho}\sum_{s}^{*} \| w \|^{k_0,\gamma}_{m, k_1, k_2+\widehat{\sigma}+\alpha} 
\| \breve\beta \|^{k_0,\gamma}_{k_3+\widehat{\sigma}+\alpha}\,,
\end{aligned}
\end{equation}
for some $\widehat{\sigma}\geq \sigma_{m+\rho-1}+3|m|+4+3\rho$ 
(recall that $\sigma_{k}$ is a non-decreasing sequence in $k$)
and provided that \eqref{buf} holds for some $\sigma\geq \widehat{\sigma}$.
In the last inequality we also used the interpolation estimate in Lemma \ref{interpolazione fine}
following the reasoning used for \eqref{verdone2}.
	\noindent
	To estimate the Lipschitz variation we first note (recall \eqref{sperobeneResto})
	\[
	\Delta_{12}\mathtt{r}_{-\rho}=
	\sum_{k=0}^{m+\rho-1}r_{m-k+1+\rho}(\Delta_{12}q_{m-k},\chi(i_2))
	+r_{m-k+1+\rho}(q_{m-k}(i_1),\Delta_{12}\chi)\,.
	\]
	Then, reasoning as done for $\mathtt{r}_{-\rho}$ one can check that
	\begin{equation}\label{stimaerrinomenorhoDelta12}
	\begin{aligned}
	\|\Delta_{12}\mathtt{r}_{-\rho}\|_{-\rho,p,\alpha} 
&\lesssim_{m,s, \alpha, \rho} 
  \sum_{p}^* \|\Delta_{12}w\|_{m,k_1,\alpha+k_2+\widehat{\sigma}} 
  \sup_{j=1,2}\|\breve\beta(i_j)\|_{k_3+\widehat{\sigma}}
\\&+
 \sum_{p+1}^*   \sup_{j=1,2}\|w(i_j)\|_{m,k_1,\alpha+k_2+\widehat{\sigma}}
   \sup_{j=1,2} \|\breve\beta(i_j)\|_{k_3+\widehat{\sigma}} 
 \|\Delta_{12}\breve\beta \|_{s_0+1} \,,
\end{aligned}
	\end{equation}
	for some $\widehat{\sigma}$ (possibly smaller than the one in \eqref{stimaerrinomenorho}).
	
	\noindent
We now conclude the construction of the solution of \eqref{ars}.	
We set (recall \eqref{balconata0})
	\begin{equation}\label{opPtau}
		P^{\tau}=Q^{\tau}+R^{\tau}\,, \qquad Q^\tau={\rm Op}(q)\in OPS^m\,.
	\end{equation}
	By construction, we have that (recall \eqref{balconata0},  \eqref{probapproxsimboloq}) 
	\[
	\partial_{\tau} Q^{\tau}=[{X}^{\tau}, Q^{\tau}] + \mathcal{M}^\tau\,, \quad  \mathcal{M}^\tau:={\rm Op}(\mathtt{r}_{-\rho})  \, , 
	\]
so that the remainder $R^{\tau}$ is the solution of the problem
	\begin{equation}\label{equazioneRRtauesplicita}
	\partial_{\tau} R^{\tau}=[{X}^{\tau}, R^{\tau}]-\mathcal{M}^{\tau}\,,\qquad R^{0}=0\,.
	\end{equation}
	We set $V^{\tau}:=(T_{\breve\beta}^{\tau})^{-1}\circ R^{\tau}\circ T_{\breve\beta}^{\tau}$. 
	It is easy to check, using 
	 \eqref{equazioneRRtauesplicita}, \eqref{flussodiffeo}, that $V^{\tau}$ has to solve the problem
	 \[
	 \partial_{\tau} V^{\tau} = -(T_{\breve\beta}^{\tau})^{-1}\circ\mathcal{M}^{\tau}\circ  T_{\breve\beta}^{\tau}\,,
	 \qquad V^0=0\,,
	 \]
	 	whose solution is 
	\begin{equation*}%\label{blacksabbath}
		V^{\tau}=-
		\int_0^{\tau} 
		(T_{\breve\beta}^{t})^{-1} \circ
		\mathcal{M}^{t}\circ T_{\breve\beta}^{t}\,dt \,.
	\end{equation*}
	We then deduce that 
	\begin{equation}\label{restofinalissimo} 
		R^{\tau}=-
		\int_0^{\tau} T_{\breve\beta}^{\tau}\circ (T_{\breve\beta}^{t})^{-1} \circ
		\mathcal{M}^t \circ T_{\breve\beta}^{t} \circ(T_{\breve\beta}^{\tau})^{-1}\,dt\,.
	\end{equation}
	It remains to prove the bounds \eqref{francia1}-\eqref{francia2}. 
	We first recall that $\mathcal{M}^{\tau}$ is a pseudo-differential operator of order $-\rho$, with $\rho$ given in 
	\eqref{sceltarho}. Then, using \eqref{stimaerrinomenorho} and 
	Lemma \ref{constantitamesimbolo}
	one can prove that 
	the operator $\mathcal{M}^{\tau}$ satisfies,%, for $M\leq N-m(\tb+k_0)-2$, and
	for $s_0\leq s\leq \bar{s}-\widehat{\s}$\footnote{Actually $\mathcal{M}^{\tau}$ 
	satisfies even better bounds w.r.t. \eqref{stimesuMop}. However, \eqref{stimesuMop} is sufficient for our aims.},
	\begin{equation}\label{stimesuMop}
	\begin{aligned}
&\mathfrak{M}_{\jap{D}^{M+ 2(\tb + k_{0})} \pa_\vphi^{b} \mathcal{M}^{\tau} \jap{D}^{-\tc + 2(\tb + k_{0})}}(s)
\lesssim_{m,s, M} 
%\|w\|_{m,s+\rho,\sigma_1}^{k_0, \gamma} + 
\sum_{s+\s_1}^{*} \|w\|_{m,k_1,\al+k_2+\widehat{\sigma}}^{k_0, \gamma} \|\breve\beta\|_{k_3+\widehat{\sigma}}^{k_0, \gamma}\,,
\end{aligned}
	\end{equation}
	uniformly in $\tau\in[0,1]$.
	Using also Lemma \ref{lemma: action Sobolev} and estimate \eqref{stimaerrinomenorhoDelta12},
	we deduce
	\[
	\begin{aligned}
	&\sup_{\tau\in [0,1]}\|\jap{D}^{-\tc + 2(\tb + k_{0})} \pa_\vphi^b \Delta_{12} \mathcal{M}^{\tau} \jap{D}^{-\tc + 2(\tb + k_{0})}\|_{\mathcal{L}(H^p;H^p)}
\\&\lesssim_{m, p, M}
  \sum_{p+M-\tc}^* 
  \sup_{j=1,2}\|w(i_{j})\|_{m,k_1,k_2+\sigma_2} 
  \sup_{j=1,2}\|\breve\beta(i_{j})\|_{k_3+\sigma_2} 
  \|\Delta_{12}\breve\beta \|_{s_0+\sigma_2} 
\\&
+  \sum_{p+M-\tc}^* \|\Delta_{12}w\|_{m,k_1, k_2+\sigma_2} 
\sup_{j=1,2}\|\breve\beta(i_{j})\|_{k_3+\sigma_2}\,.
 \end{aligned}
	\]
	The latter bound implies that $\mathcal{M}^{\tau}$ satisfies the assumptions of 
	Lemma \ref{conj.NEWsmoothresto} %with $\alpha=0$. 
	Since  $T_{\breve\beta}^{\tau}$ is the flow of \eqref{flussodiffeo}
and satisfies the estimates provided by Lemma \ref{lemma:buonaposFlussi}, we have that the reminder
$R^{\tau}$ in \eqref{restofinalissimo} satisfies the desired estimates \eqref{francia1}-\eqref{francia2}
using 	Lemma \ref{conj.NEWsmoothresto}.
By formul\ae\,\eqref{opPtau}, \eqref{opPtautau}, \eqref{balconata0}
we have obtained that
\[
P^{1}=T_{\breve\beta}^{\tau} \circ {\rm Op}(w)\circ  (T_{\breve\beta}^{\tau})^{-1}_{|\tau=1}
=\op(q_0+\widetilde{q}_1)+R\,,
\]
with $q_0$ as in \eqref{formaEsplic2}, $R$ satisfying \eqref{francia1}-\eqref{francia2} and 
$\widetilde{q}_1=\sum_{k=0}^{m+\rho-1} q_{m-k}$ satisfying bounds like \eqref{parlare}-\eqref{troppo} with $i=1$.
The thesis \eqref{simbotransportato} follows 
by applying  Lemma \ref{passaggiostweyl} to pass to the Weyl quantization.
	 
 Finally assume that $\breve\beta$ in \eqref{egogenerator} and $w$ \eqref{simbotransportato} 
are momentum preserving so that the map $\breve{\mathcal{A}}^{\tau}$ is momentum preserving
by Lemma \ref{lem:momentoflusso}.
In view of Remark \ref{rmk:algsimboli},
%following word by word the proof of Theorem $3.4$ in \cite{FGP19},
one can check
that
the final symbol 
$q$ in \eqref{formaEsplic} is momentum preserving.
The remainder $R$ is  momentum preserving by difference. This concludes the proof.
\end{proof}
 
\section{Tame estimates for the Dirichlet-Neumann operator}\label{sec:DNsezione}
We now present some fundamental properties of the Dirichlet-Neumann operator $G(\eta)\psi$ 
defined in \eqref{eq:112a} 
that are used in the paper. 
There is a huge literature about it for which we refer for instance to the work of Alazard-Delort \cite{ADel1}
and the book of Lannes \cite{Lan2}, and references therein. 
We remark that for our purposes it is sufficient to work in the class of smooth $C^{\infty}$ profiles 
$\eta(x)$ because at each step of the Nash-Moser iteration we perform a $C^{\infty}$-regularization.
The map $(\eta,\psi)\to G(\eta)\psi$ is linear with respect to $\psi$
and nonlinear with respect to the profile $\eta$.
The derivative with respect to $\eta$ (which is called ``shape derivative'')
is given by the formula (see for instance \cite{Lan1}) %, \cite{Lan2})
\begin{equation}\label{shapeDer}
G'(\eta)[\hat{\eta}]\psi=\lim_{\epsilon\to0} \frac{\big(G(\eta+\epsilon\hat{\eta})\psi-
G(\eta)\psi
\big)}{\varepsilon}=-G(\eta)(B\hat{\eta})-\div\big(V\hat{\eta}\big)\,,
\end{equation}
where we denoted 
the horizontal and vertical components 
of the velocity field at the free interface by
\begin{align} 
\label{def:V}
& V =  V (\eta, \psi) :=  (\nabla \underline{\phi}) (x, \eta) = 
\nabla\psi - \nabla\eta B \, , 
\\
\label{form-of-B}
& B =  B(\eta, \psi) := (\pa_y \underline{\phi}) (x, \eta) =  
\frac{G(\eta) \psi + \nabla\eta\cdot\nabla\psi}{ 1 + |\nabla\eta|^2} \, .
\end{align}
We are interested in study the pseudo-differential structure 
of $G(\eta)$. 
Such structure has been widely investigated, 
we refer, for instance, to 
\cite{ADel1}, \cite{Lan2}, \cite{BM20}
and reference therein.
However, in this paper we need very precise and tame estimates 
of $G(\eta)$ posed in high space dimension.
This is the content of Theorem \ref{lemma totale dirichlet neumann} below.
In order to state our main result we recall that (see \eqref{eq:112a})
\begin{equation*}
 G(\eta)\psi := \sqrt{1+|\nabla\eta|^2}(\partial_n\underline{\phi})(x,\eta(x))
 = (\partial_z\uphi -\nabla\eta\cdot \nabla\uphi)(x,\eta(x))
\end{equation*}
where $\uphi=\uphi(x,z)$ is the solution of the elliptic problem\footnote{we defined $\Delta_{z,x}:=
\pa_{zz}+\Delta$.}
\begin{equation}\label{elliptic1}
\left\{\begin{aligned}
\Delta_{z,x}\uphi&=0\,,\qquad x\in\T^{d}_\Gamma\,,\;\; - \mathtt h <z<\eta(x)\,,\\
\uphi(x,\eta(x))&=\psi(x)\,, \quad x \in \T^d_\Gamma\,, \\
(\pa_{z}\uphi)(x,-\mathtt{h})&=0\,, \quad x \in \T^d_\Gamma\,. 
\end{aligned}\right.
\end{equation}
Our aim is to study the pseudo-differential structure of $G(\eta)[\psi]$
when the functions $\eta(t,x),\psi(t,x)$ satisfy the following assumption.
\begin{hypo}{\bf (Assumptions on the free surface $\eta$).}\label{hypo:eta}
Given (see \eqref{def:setlambda0}) we consider 
\[
\mathtt \Omega \times [\mathtt h_1, \mathtt h_2] \ni\lambda:=(\omega,\mathtt{h})\mapsto (\eta(\lambda),\psi(\lambda))\in 
C^{\infty}(\T^{\nu}\times\T_{\Gamma}^{d};\R^{2})\,,
\]
where $\mathtt \Omega \subset \R^\nu$ is a bounded iopen set, satisfying the following:

\noindent
$\bullet$ {\bf (Reality).} The functions $\eta,\psi$ are real valued;

\noindent
$\bullet$ {\bf (Whitney-Sobolev regularity).}  $\eta,\psi$ depend on the 
parameters $\lambda=(\omega,\mathtt{h})$  in 
 a Lipschitz way
together with their derivatives in the sense of Whitney
and they have finite $\|\cdot\|_{s}^{k_0,\gamma}$ norms (see Section \eqref{sec:Whitney-Sobolev});

%\noindent
%$\bullet$ {\bf (Dependence on the embeddings).}  $\eta,\psi$ depends on the 
%``approximate'' torus $i(\vphi)$ (see \eqref{torusprototipo})
%in a Lipschitz way;

\noindent
$\bullet$ {\bf (Traveling wave).}  $\eta,\psi$  are traveling wave in the sense of
\eqref{def:vec.tau}-\eqref{condtraembedd}.
\end{hypo}
With abuse of notation, differently form \eqref{deltaunodue},
 in this section and in the next section \ref{sec:pseudoDN}, we denote $\Delta_{12}$ the variation of symbols, operators and functions,
with respect the profile $\eta$ instead of the embedding $i$, namely for any function $f (\eta)$ depending on $\eta$, we write
$$
\Delta_{12} f := f(\eta_1) - f(\eta_2)\,. 
$$
Moreover in this section and in the next section \ref{sec:pseudoDN}, we fix the minimal regularity $s_0$ as
\begin{equation}\label{s0 Dirichlet Neumann}
s_0 := \frac{\nu + d + 1}{2} + 1\,. 
\end{equation}
Our main result is the following.
\begin{thm}\label{lemma totale dirichlet neumann}
Let $\beta_0, k_0 \in \N$, $M\in \mathbb{N}$, $\mathtt{c}\geq \beta_0+k_0+2$
and $M\geq \mathtt{c}$.
%There exist $\s=\s(\beta_0,k_0,M)$ such that
%Then there exist $\bar N \equiv \bar N(\beta_0, k_0) \gg 0$ such that for any $N\geq \bar{N}$
There exist
$\sigma_1 =\sigma_1(\beta_0, k_0, M)\geq \s_2=\sigma_2(\beta_0, M) \gg 0$ 
such that the following holds.
%For 
For any $\bar{s}\geq s_0+\s_1$ there exists $\delta=\delta(\beta_0,k_0,M,\bar{s})$ such that if 
\begin{equation}\label{smalleta}
\| \eta \|_{s_0 + \sigma_1}^{k_0, \gamma}\leq \delta\,,
\end{equation}
and satisfies  Hypothesis \ref{hypo:eta},
then  the Dirichlet-Neumann operator $G(\eta)$  in \eqref{eq:112a}
admits the expansion
\begin{equation}\label{pseudoespansione}
G(\eta) =  |D| \tanh(\mathtt h |D|) +   {\mathcal A}_G+ {\mathcal B}_G + {\mathcal R}_G
\end{equation}
where ${\mathcal A}_G = {\rm Op}^W\big( a_G(\vphi, x, \xi; \lambda) \big) \in OPS^1$, 
${\mathcal B}_G = {\rm Op}^W\big( b_G(\vphi, x, \xi; \lambda) \big) \in OPS^0$ 
and the remainder ${\mathcal R}_G$ is a regularizing remainder satisfying the following.

The pseudo-differential operators ${\mathcal A}_G, {\mathcal B}_G$ satisfy the estimates, for any $\alpha\geq0$,
\begin{equation}\label{stimefinalisimboliDN}
\begin{aligned}
\| {\mathcal A}_G \|_{1, s, \alpha}^{k_0,\gamma}\,,\, \| {\mathcal B}_G \|_{0, s, \alpha}^{k_0, \gamma} 
&\lesssim_{s,\alpha, M}  \| \eta \|_{s + \sigma_1}^{k_0, \gamma} \,, \qquad  s_0\leq s\leq \bar{s}-\s_1\,,
\\
 \| \Delta_{12}{\mathcal A}_G \|_{1, p, \alpha}\,,\, \| \Delta_{12}{\mathcal B}_G\|_{1, p, \alpha}
 &\lesssim_{p,\alpha, M}\| \eta_1-\eta_2 \|_{p + \sigma_2}\,, \;\;\;p+\s_2\leq s_0+\s_1\,,
\end{aligned}
\end{equation}
The symbols $a_G, b_G$ are real valued and the symbol $a_G$ satisfies the symmetry condition
\begin{equation}\label{contagion1}
a_G(\vphi, x, \xi) = a_G(\vphi, x, - \xi)\,. 
\end{equation}
The smoothing remainder ${\mathcal R}_G$ satisfies the following property. 
For any  $|\beta| \leq \beta_0$, 
one has that $ \langle D \rangle^{M} \partial_\vphi^\beta{\mathcal R}_G \langle D \rangle^{- \mathtt{c}}$ 
is a ${\mathcal D}^{k_0}$-tame operator with tame constant satisfying the bounds
\begin{equation}\label{stimefinaliRR}
\begin{aligned}
{\mathfrak M}_{ \langle D \rangle^{M} \partial_\vphi^\beta {\mathcal R}_G \langle D \rangle^{-  \mathtt{c}}}(s) 
&\lesssim_{\bar{s}, M}
\| \eta \|_{s + \sigma_1}^{k_0, \gamma}\,, 
\quad  s_0 \leq s \leq \bar{s}-\s_1\,,
\\
\|\langle D \rangle^{M} 
 \partial_\vphi^\beta\Delta_{12}
\cR_G
 \langle D \rangle^{- \mathtt{c}} \|_{\mathcal{L}({H}^{p};{H}^{p})}
&\lesssim_{\bar{s}, M} \|\eta_{1}-\eta_{2}\|_{p+\s_2}\,,\;\;\;p+\s_2\leq s_0+\s_1\,.
\end{aligned}
\end{equation}
Furthermore, $G(\eta)$, ${\mathcal A}_G$, ${\mathcal B}_G$ and ${\mathcal R}_G$ are real-to-real 
and momentum preserving operators. 
\end{thm}

An important consequence of the result above is the following.

\begin{lemma}\label{stime tame dirichlet neumann}
There exist $\s=\s(s_{0}, k_{0})>0$ such that, 
for any $s\geq s_0$, there is $\delta(s, k_0) > 0$ such that if  
 $\| \eta\|_{s_0+\s}^{k_0, \gamma} \leq \delta$,  then 
 \begin{align}
  \| \big(G(\eta) - |D| \tanh(\th |D|) \big) \psi \|_s^{k_0, \gamma} 
  & \lesssim_{s, k_0}  \| \eta \|^{k_0, \gamma}_{s+\s} \| \psi\|_{s_0+\s}^{k_0, \gamma} 
   + \| \eta \|^{k_0, \gamma}_{s_0+\s} \| \psi\|_{s+\s}^{k_0, \gamma}\,, \label{stima tame dirichlet neumann}
\\
  \|  G'(\eta)[\widehat \eta] \psi \|_s^{k_0, \gamma}  
  &\lesssim_{s, k_0} \| \psi\|_{s + \s}^{k_0, \gamma} \| \widehat \eta\|_{s_0 + \s }^{k_0, \gamma} 
  + \| \psi\|_{s_0 + \s }^{k_0, \gamma} \| \widehat \eta\|_{s + \s }^{k_0, \gamma}  \nonumber
  \\
  &\quad \quad + \| \eta\|_{s+\s}^{k_0, \gamma} \| \widehat \eta\|_{s_0  + \s }^{k_0, \gamma} 
  \| \psi\|_{s_0 + \s }^{k_0, \gamma}\,, \label{stima tame derivata dirichlet neumann}
  \\
  \| G''(\eta)[\widehat \eta, \widehat \eta] \psi \|_s^{k_0, \gamma}  
  & \lesssim_{s, k_0} 
  \big(\| \psi\|_{s + \s}^{k_0, \gamma}+\| \eta\|_{s + \s}^{k_0, \gamma} \| \psi\|_{s_0 + \s}^{k_0, \gamma}  \big)
  \big( \| \widehat \eta\|_{s_0 + \s}^{k_0, \gamma} \big)^2  \nonumber
  \\
  &\quad \quad + \| \psi\|_{s_0 + \s}^{k_0, \gamma} \| \widehat \eta\|_{s + \s}^{k_0, \gamma} 
  \| \widehat \eta\|_{s_0 + \s}^{k_0, \gamma}  \,.  \label{stima tame derivata seconda dirichlet neumann}
 \end{align}
   \end{lemma}
\begin{proof}
Given $s\geq s_0$,
we apply Theorem \ref{lemma totale dirichlet neumann} taking $\bar{s}\rightsquigarrow s+\s$ 
%in such a way $\bar{s}-\s\gg s$
for some $\s>0$ larger than $\s_1$ of Theorem \ref{lemma totale dirichlet neumann}. 
Therefore,
the estimate \eqref{stima tame dirichlet neumann} follows by the 
formula \eqref{pseudoespansione}, Lemmata \ref{lemma: action Sobolev}, \ref{constantitamesimbolo}, \ref{lemma operatore e funzioni dipendenti da parametro} and estimates \eqref{stimefinalisimboliDN} and \eqref{stimefinaliRR}.
The estimate \eqref{stima tame derivata dirichlet neumann} follows by the shape derivative formula \eqref{shapeDer}, applying  \eqref{stima tame dirichlet neumann}, \eqref{p1-pr} and the fact that the functions
$ B , V $ defined in \eqref{def:V}-\eqref{form-of-B} satisfy  
\[
\| B\|_s^{k_0, \gamma},\, \| V\|_{s}^{k_0, \gamma} 
\leq_s \| \psi\|_{s + \s}^{k_0, \gamma} + \| \eta\|_{s + \s}^{k_0, \gamma} \| \psi\|_{s_0 + \s}^{k_0, \gamma}\,.
\]
The estimate \eqref{stima tame derivata seconda dirichlet neumann} 
follows by differentiating the shape derivative formula \eqref{shapeDer} 
and by applying the same kind of arguments. 
\end{proof}

The proof of Theorem \ref{lemma totale dirichlet neumann} involves several arguments.
First of all, in Section \ref{sec:tameLaplace} we shall provide a priori tame estimates
on the solution $\underline{\phi}$
of the elliptic problem \eqref{elliptic1}.
This part is somehow more classical. However, we shall provide suitable \emph{tame} estimates
in the spaces $H^{s}$ in \eqref{SobolevSpazioAngoli}
taking into account both the regularity in $\vphi\in \T^{\nu}$ and in $\lambda\in \R^\nu \times [\mathtt h_1, \mathtt h_2]$.
We partially follow the strategy and notations in \cite{BeMasVent1}.
Then in Section \ref{sec:pseudoDN} we provide the pseudo-differential expansion \eqref{pseudoespansione}
and conclude the proof of Theorem \ref{lemma totale dirichlet neumann}.

\begin{rmk}
We make some comments about the constants $\s_1,\s_2$ appearing in Theorem 
\ref{lemma totale dirichlet neumann}, and which 
represent the \emph{loss of derivatives} accumulated
along the procedure. 
For our aims the crucial point is to ensure that $\s_1$ is independent of the Sobolev index $s$.
It turns out that it depends only of the parameters $k_0, \beta_0,M$
which are  \emph{fixed arbitrarily} at the beginning of the procedure. 
In our procedure we need a smallness condition on the \emph{low} 
norm $\|\eta\|_{s_0+\s_1}$
in \eqref{smalleta}. The smallness parameter $\delta$ will depends also on the ``high'' 
Sobolev index $\bar{s}$ (fixed arbitrarily large). This allow us to provide estimates on symbols and operators
in the expansion of the Dirichlet-Neumann operator
for any Sobolev regularity $s$ in the range $s_0\leq s\leq \bar{s}-\s_1$ (see \eqref{stimefinalisimboliDN}-\eqref{stimefinaliRR}) which is arbitrarily large by taking $\bar{s}$ larger.
We also need 
 to estimate the variation with respect to the shape of the water $\eta(\vphi,x)$ in a low norm
 $\|\cdot\|_{p}$ \ 
for all Sobolev indexes $p$ such that
\begin{equation}\label{lownormasigma2}
p+\sigma_2 \leq s_0 +\s_1 \,, \quad \text{ for some } \ 0<\sigma_2:=\sigma_2(M,\beta_0)\leq \s_1\,. 
\end{equation}
In general, the \emph{loss} of derivatives $\s_2$ is smaller than $\s_1$
and it will not depends on the parameter $k_0$, since we do not need to control the Lipschitz variations of operators 
taking into account their regularity in the parameters $\lambda$.
 Condition \eqref{lownormasigma2} guarantees that we are below the threshold of regularity for which the 
 function $\eta$ is small, this simplify a lot the estimates for the Lipschitz variation.
 Secondly, \eqref{lownormasigma2} can be fulfilled just by taking $\s_1$ larger.
\end{rmk}

\noindent
{\bf (Notation).} Sometimes, when it does not create confusion, we shall omit the dependence 
of $\eta$ on the parameters $\lambda, \vphi$.

\subsection{The elliptic problem on the straight strip}\label{sec:tameLaplace}

Following \cite{Lan1},
in order study the problem
\eqref{elliptic1}  (which is defined on $-\mathtt{h}\leq z\leq \eta(x)$),
we introduce 
new coordinates where  \eqref{elliptic1}
transforms into a variable coefficients elliptic problem defined on a \emph{straight}
strip.

%\subsection{Straightening the free surface}\label{stripdritta}

More precisely, consider the new vertical variable $y \in [- 1, 0]$
\begin{equation}\label{change1}
z=\eta(x)+y(\mathtt h+\eta(x))\,,\quad -1\leq y\leq 0\,,
\end{equation}
so that the closure of $\Omega_{t}$ in \eqref{Deta} becomes
\begin{equation*}
\mathcal{S}_1=\{(x,y)\in \T^{d}_\Gamma \times \R\; :\; -1\leq y\leq 0\}\,.
\end{equation*}
Consider the function
\begin{equation}\label{func:phi}
\phi(x,y):=\uphi(x,\eta(x)+y( \mathtt h +\eta(x)))\,.
\end{equation}
We have the following.
\begin{lemma}\label{lem:firststright}
Assume Hyp. \ref{hypo:eta} and let $s\geq s_0$. Assume also that 
$\uphi(x,z)$ solves the problem \eqref{elliptic1}.  
There are functions $\alpha_{1}\in C^{k}((-1,0);(H^{s})^d)$, 
$\alpha_{i}\in C^{k}((-1,0);H^{s})$ for $i=2,3$, $k\in\N$, such that
the function
$\phi(x,y)$ in
\eqref{func:phi} solves the problem
\begin{equation}\label{elliptic2}
\left\{\begin{aligned}
\widetilde{\mathcal{L}}\phi&=0\,,\qquad x\in\T^{d}_{\Gamma}\,,\;\; -1<y<0\,,\\
\phi(x,0)&=\psi(x)\,,\\
(\pa_{y}\phi)(x,-1)&=0\,,
\end{aligned}\right.
\end{equation}
where 
\begin{equation}\label{op:tildeL}
\widetilde{\mathcal{L}}:=\widetilde{\mathcal{L}}(\eta;y):=
\pa_{yy}+(\alpha_1(x,y) \cdot \nabla+\alpha_2(x,y)\big)\pa_{y}+\mathtt{h}^2\alpha_{3}(x,y)\Delta\,.
\end{equation}
%and $\alpha_1, \alpha_2, \alpha_3$ satisfies the following estimates. 
Moreover, there is $\s=\s(k_0)\gg0$ such that for any
$s \geq s_0$, $n \in \N$, $|k| \leq k_0$ there is $\delta=\delta(s,k_0)\ll1$ such that 
 if \eqref{smalleta} holds with $\s_1\geq \s$,
%$\| \eta \|_{s_0 + \sigma}^{k_0, \gamma} \ll 1$, 
%then for any $s \geq s_0$, $n \in \N$, $|k| \leq k_0$, 
then
\begin{equation}\label{stimealphai}
\begin{aligned}
 \| \partial_\lambda^k \partial_y^n \alpha_1 \|_{{\mathcal H}^s}\,,\, 
\| \partial_\lambda^k  \partial_y^n \alpha_2 \|_{{\mathcal H}^s}\,,\,  
\| \partial_\lambda^k \partial_y^n (\alpha_3 - 1) \|_{{\mathcal H}^s} 
&\lesssim_{s, n, k} \gamma^{- |k|} \| \eta \|_{s + \sigma}^{k_0, \gamma}\,, 
\\
 \| \partial_y^n \alpha_1 \|_s^{k_0, \gamma}\,,\, \| \partial_y^n \alpha_2 \|_s^{k_0, \gamma}\,,\,  
 \| \partial_y^n (\alpha_3 - 1) \|_s^{k_0, \gamma} 
 &\lesssim_{s, n} \| \eta \|_{s + \sigma}^{k_0, \gamma}\,,
\end{aligned}
\end{equation}
%One also has that
\begin{equation}\label{stimealphaiDelta12}
\| \partial_y^k \Delta_{12}\alpha_1 \|_p\,,\, \| \partial_y^k \Delta_{12}\alpha_2 \|_p\,,\, 
\| \partial_y^k \Delta_{12}\alpha_3  \|_p \lesssim_{p, k} \| \eta_1-\eta_2 \|_{p + \widetilde{\sigma}}\,,
\quad p\geq s_0\,.
\end{equation}
for some $0<\widetilde{\s}\leq \s$ and $p+\widetilde{\s}\leq s_0+\s$. 
The estimates \eqref{stimealphai}-\eqref{stimealphaiDelta12} are uniform in $y\in[-1,0]$.
\end{lemma}
\begin{proof}
%We refer to Appendix \ref{proof:lem:firststright}.
First of all by \eqref{change1} and \eqref{func:phi} we note
\begin{equation}\label{nuovecondini}
\begin{aligned}
&\uphi(x,z)=\phi\Big(x, \frac{z-\eta(x)}{\mathtt{h}+\eta(x)}\Big)\,,
\\
&\phi(x,0)=\uphi(x,\eta)=\psi\,,\qquad
(\pa_{y}\phi)(x,-1)=(\mathtt h +\eta)(\pa_{z}\uphi)(x,- \mathtt h)=0\,.
\end{aligned}
\end{equation}
%and 
%\[
%\uphi(x,z)=\phi\Big(x, \frac{z-\eta(x)}{\mathtt{h}+\eta(x)}\Big)\,.
%\]
%Then we have
%\begin{align*}
%\pa_{z}\uphi(x,z)&=\frac{1}{\mathtt{h}+\eta(x)}\pa_{y}\phi(x,y)\,,
%\qquad
%\pa_{zz}\uphi(x,z)=\frac{1}{(\mathtt{h}+\eta(x))^2}\pa_{yy}\phi(x,y)\,,
%\\
%\nabla\uphi(x,z)&=\nabla\phi(x,y)+(\pa_y\phi)(x,y)\nabla\Big(\frac{z-\eta(x)}{\mathtt{h}+\eta(x)}\Big)=
%\nabla\phi(x,y)-\frac{(1+y)\nabla\eta}{\mathtt{h}+\eta(x)}(\pa_y\phi)(x,y)\,,
%\\
%\Delta\uphi(x,z)&=\Delta\phi(x,y)-\frac{2(1+y)\nabla\eta}{\mathtt{h}+\eta(x)}(\pa_y\nabla\phi)(x,y)
%+\frac{(1+y)^2|\nabla\eta|^2}{(\mathtt{h}+\eta(x))^2}(\pa_{yy}\phi)(x,y)
%\\&-\frac{(1+y)\Delta\eta}{\mathtt{h}+\eta(x)}(\pa_y\phi)(x,y)
%+\frac{(1+y)|\nabla\eta|^2}{(\mathtt{h}+\eta(x))^2}(\pa_y\phi)(x,y)\,.
%\end{align*}
%Recalling \eqref{elliptic1} we have that the function $\phi(x,y)$ in \eqref{func:phi}
%solves 
%\begin{align*}
%0=\pa_{zz}\uphi+\Delta\uphi&=
%\Big(\frac{1+(1+y)^2|\nabla\eta|^2}{(\mathtt{h}+\eta(x))^2}\Big)\pa_{yy}\phi+\Delta\phi
%\\&+\Big[-\frac{2(1+y)}{\mathtt{h}+\eta(x)}\nabla\eta\cdot\nabla-\frac{(1+y)\Delta\eta}{\mathtt{h}+\eta(x)}
%+\frac{(1+y)|\nabla\eta|^2}{(\mathtt{h}+\eta(x))^2}\Big]\pa_{y}\phi
%\end{align*}
%with initial condition \eqref{nuovecondini}. This is equivalent to
Therefore, by an explicit computation, using also \eqref{elliptic1}, one can check that 
$\phi(x,y)$ in \eqref{func:phi} solves 
\[
\begin{aligned}
0&=\pa_{yy}\phi+(\alpha_1(x,y) \cdot \nabla+\alpha_2(x,y)\big)\pa_{y}\phi+\mathtt{h}^{2} \alpha_{3}(x,y)\Delta\phi\,,
\end{aligned}
\]
with  conditions \eqref{nuovecondini} and 
where we defined
\begin{align*}
\alpha_{3}(x,y)&:=\frac{\big(1+\frac{\eta}{\mathtt{h}}\big)^2}{1+(1+y)^2|\nabla\eta|^2}\,,
\qquad
\alpha_1(x,y):=-\frac{2(1+y)(\mathtt{h}+\eta)}{1+(1+y)^2|\nabla\eta|^2}\nabla\eta\,,
\\
\alpha_{2}(x,y)&:=(1+y)\frac{-(\mathtt{h}+\eta)\Delta\eta+2|\nabla\eta|^2}{1+(1+y)^2|\nabla\eta|^2}\,.
\end{align*}
The discussion above, together with \eqref{nuovecondini}, implies
\eqref{elliptic2}-\eqref{op:tildeL}.
By an explicit computation one deduces that 
the functions $\alpha_{i}$ are smooth in $y\in [-1,0]$ and satisfies 
the estimates
\eqref{stimealphai} for some $\s>0$
provided that the norm $\| \eta \|_{s_0 + \sigma}^{k_0, \gamma}$
is small enough w.r.t. constants depending only on $s$ and $\mathtt{h}\in[\mathtt{h_1},\mathtt{h}_2]$.
The estimates \eqref{stimealphaiDelta12} on the Lipschitz variation in $\eta$
of the functions $\alpha_i$ follows similarly by using their explicit definition.
\end{proof}
We also prove the following lemma that we shall use later in a systematic way. 
\begin{lemma}\label{lemma stima tame prodotto con gli alpha i}
Assume that $\alpha(\vphi, y, x; \lambda )$ satisfies for any $s \geq 0$, $\beta \in \N^\nu$, 
$k \in \N^{\nu + 1}$, $|k| \leq k_0$, $n \in \N$ the estimate 
\[
\sup_{y \in [- 1, 0]} \sup_{\vphi \in \T^\nu} 
\| \partial_\lambda^k \partial_\vphi^\beta \partial_y^n \alpha \|_{H^s_x} 
\lesssim_{s, n, k, \beta} \gamma^{- |k|} \| \eta \|_{s + |\beta| + \sigma}^{k_0, \gamma}\,,
\]
for some $\sigma > 0$ independent of $n$ and $\beta$ 
and let $s \geq 0$, $u \in {\mathcal H}^s$. Then 
\[
\| (\partial_\lambda^k \partial_\vphi^\beta \alpha) u \|_{{\mathcal H}^s} 
\lesssim_{s, k , \beta} \gamma^{- |k|} \Big( \| \eta \|_{s + |\beta| 
+ \sigma'}^{k_0, \gamma} \| u \|_{L^2_{x, y}} 
+ \| \eta \|_{s_0 + |\beta| + \sigma'}^{k_0, \gamma} \| u \|_{{\mathcal H}^s} \Big)\,.
\]
for some $\s'\geq \s$.
\end{lemma}
\begin{proof}
By \eqref{prop L2 vhi Hs Hs vphi cal 2 xy}
we use that $\| \partial_\lambda^k \partial_\vphi^\beta \alpha u \|_{{\mathcal H}^s} \simeq \| \partial_\lambda^k \partial_\vphi^\beta \alpha u \|_{L^2_y H^s_x} + \| \partial_\lambda^k \partial_\vphi^\beta \alpha u \|_{H^s_y L^2_x}$ and we estimate separately these two terms. 
%since $s \geq s_0 > d/2$, 
By using the "mixed" interpolation estimate of $H^s(\T^d_\Gamma)$ 
\[
\| u v \|_{H^s_x} \lesssim \| u \|_{L^\infty_x} \| v \|_{H^s_x} + \| u \|_{C^s_x} \| v \|_{L^2_x}\,,
\]
one gets that 
\begin{equation}\label{stima beta per u a}
\begin{aligned}
\| \partial_\lambda^k &\partial_\vphi^\beta \alpha u \|_{L^2_y H^s_x}^2 
 = \int_{- 1}^0 
\| \partial_\lambda^k \partial_\vphi^\beta \alpha(y, \cdot) u(y, \cdot) \|_{H^s_x}^2 \, d y  
\\& 
\lesssim_s 
\int_{- 1}^0 \| \partial_\lambda^k \partial_\vphi^\beta \alpha(y, \cdot) \| _{L^\infty_x}^2 
\|u(y, \cdot) \|_{H^s_x}^2 \, d y 
%\\&\qquad
+ \int_{- 1}^0 \| \partial_\lambda^k \partial_\vphi^\beta \alpha(y, \cdot)\|_{C^s_x}^2 
\| u(y, \cdot) \|_{L^2_x}^2 \, d y 
\\& 
\lesssim_s 
\| \partial_\lambda^k \partial_\vphi^\beta \alpha \|_{L^\infty_y H^{s_0}_x}^2 
\| u \|_{L^2_y H^s_x}^2 
+ \| \partial_\lambda^k \partial_\vphi^\beta \alpha \|_{L^\infty_y H^{s + s_0}_x}^2 
\| u \|_{L^2_y L^2_x}^2 
\\& 
{\lesssim_{s, k, \beta}} \
\gamma^{- 2 |k|}  (\| \eta \|_{s_0 +|\beta| + \sigma}^{k_0, \gamma})^2 
\| u \|_{{\mathcal H}^s}^2 
+  \gamma^{- 2 |k|} (\| \eta \|_{s +s_0+ |\beta| + \sigma}^{k_0, \gamma})^2 
\| u \|_{L^2_{x, y}}^2\,. 
\end{aligned} 
\end{equation}
Moreover 
\begin{equation}\label{stima beta per u b}
\begin{aligned}
\| (\partial_\lambda^k \partial_\vphi^\beta \alpha) u \|_{H^s_y L^2_x} 
& \lesssim_s 
\sum_{s_1 + s_2 \leq s} \| \partial_y^{s_1} (\partial_\lambda^k \partial_\vphi^\beta \alpha) \partial_y^{s_2} u \|_{L^2_y L^2_x}
\\& 
\lesssim_s \sum_{s_1 + s_2 \leq s} \| \partial_y^{s_1} (\partial_\lambda^k \partial_\vphi^\beta \alpha) \|_{L^\infty_y L^\infty_x} \| \partial_y^{s_2} u \|_{L^2_y L^2_x}  
\\& 
\lesssim_s  \sum_{s_1 + s_2 \leq s} \| \partial_y^{s_1} (\partial_\lambda^k \partial_\vphi^\beta \alpha) \|_{L^\infty_y H^{s_0}_x} \| \partial_y^{s_2} u\|_{L^2_y L^2_x} 
\\&
{\lesssim_{s, k, \beta}} \
  \gamma^{- |k|} \| \eta \|_{s_0 +|\beta| + \sigma}^{k_0, \gamma} \sum_{s_2 \leq s} \| \partial_y^{s_2} u\|_{L^2_y L^2_x} 
\\& 
{\lesssim_{s, k, \beta}} \
\gamma^{- |k|} \| \eta \|_{s_0 +|\beta| + \sigma}^{k_0, \gamma} \| u \|_{{\mathcal H}^s}\,.
\end{aligned}
\end{equation}
The thesis  follows by combining \eqref{stima beta per u a}, 
\eqref{stima beta per u b}. 
\end{proof}

\begin{rmk}
Notice that formula \eqref{eq:112a}
becomes
\begin{equation}\label{eq:112abis}
 G(\eta)\psi = 
\frac{1 +|\nabla\eta|^2}{\mathtt h +\eta}(\pa_{y}\phi)(x,0)-\nabla\eta\cdot\nabla \psi \,.
\end{equation}
\end{rmk}
In the following we study the solution of \eqref{elliptic2} with a perturbative argument.

\medskip
\noindent
{\bf Perturbative setting.} 
We notice that at the flat surface $\eta(x)\equiv 0$ the problem \eqref{elliptic1} reduces to
\begin{equation}\label{elliptic0}
\left\{\begin{aligned}
\Delta_{z,x}\underline \vphi_0&=0\,,\qquad x\in\T^{d}_\Gamma\,,\;\; -\mathtt{h}<z<0\,,\\
\underline \vphi_0(x,0)&=\psi(x)\,,\\
(\pa_{z} \underline \vphi_0)(x,- \mathtt h)&=0\,.
\end{aligned}\right.
\end{equation}
By an explicit computation one can check that the solution of the problem above is given by
\[
\underline \vphi_0(x,z)=\frac{\cosh((\mathtt h +z)|D|)}{\cosh(\mathtt h |D|)}[\psi]\,.
\]
For our aim it is more convenient to study the problem \eqref{elliptic0}
on a \emph{fixed} strip $[-1,0]$.
We have that the function $\vphi_0(x, y) = \underline \vphi_0(x, \mathtt h y)$ solves the problem
\begin{equation}\label{problema vphi 0 riscalato laplace}
\begin{aligned}
& \begin{cases}
 \partial_{yy} \vphi_0 + \mathtt h^{2} \Delta \vphi_0 = 0, \quad x \in \T^d_\Gamma, \quad - 1 < y < 0 \\
\vphi_0(x, 0) = \psi(x)\,,  \\
\partial_y \vphi_0(x, - 1) = 0\,. 
\end{cases}
\end{aligned}
\end{equation}
A direct calculation shows that $\vphi_0(x, y)$ is given by
\begin{equation}\label{def vphi 0}
\vphi_0(x, y):=\mathcal{L}_0[\psi]\,,\qquad \mathcal{L}_0[\cdot]=\mathcal{L}_0(\th, |D|)[\cdot]:=
\frac{\cosh((1 +y)\th|D|)}{\cosh(\mathtt h |D|)}[\cdot]\,,\quad y\in[-1,0]\,.
\end{equation}
We have the following.
\begin{lemma}\label{lemma stima cal L0 sol omogenea laplace}
For any $k \in \N^{\nu + 1}$, $s \geq 0$, $\psi \in H^{s + |k|}_x$, one has that 
\begin{equation}\label{stima cal L0 nel lemma cal Hs}
\|  \partial_\lambda^k\cL_0 [ \psi] \|_{\cH^s} \lesssim_{s,k} \| \psi \|_{H^{s + |k|}_x}\,. 
\end{equation}
Furthermore, if $\psi \in H^s$, then 
\begin{equation}\label{stima cal L0 nel lemma cal Os}
\| \langle D \rangle^{- |k|} \partial_\lambda^k\cL_0 [ \psi] \|_{\cO^s} = \| \partial_\lambda^k\cL_0 [\langle D \rangle^{- |k|} \psi] \|_{\cO^s} \lesssim \| \psi \|_{s}\,. 
\end{equation}
\end{lemma}
\begin{proof}
We prove the estimate \eqref{stima cal L0 nel lemma cal Hs}. One has that 
\[
\begin{aligned}
{\mathcal L}_0 \psi & = {\rm Op}(\Lambda(y, \xi; \mathtt h))[\psi] =  
\sum_{\xi \in \Gamma^*} \Lambda (y, \xi; \mathtt h) \widehat \psi(\xi) e^{\ii x \cdot \xi} \,,\;\;
\quad 
%\text{where} \\
 \Lambda(y, \xi ; \mathtt h)  = \dfrac{\cosh((1 +y)\th|\xi|)}{\cosh(\mathtt h |\xi|)}\,. 
\end{aligned}
\]
Since $\mathcal{L}_0$ does not depend on $\omega$, instead of computing the derivatives $\pa_{\lambda}$
it is enough to control derivatives in the parameter $\th$.
A direct calculation shows that for any $n, k \in \N$
\begin{equation}\label{stima cal L0 simbolo}
|\partial_{\mathtt h}^k \partial_y^n \Lambda(y, \xi; \mathtt h)| 
\lesssim_{k, n}  
\langle \xi \rangle^{n + k}\,, 
\qquad 
\forall y \in [- 1, 0]\,, \mathtt h \in [\mathtt h_1, \mathtt h_2]\,, \xi \in \Gamma^*\,. 
\end{equation}
Note that 
$\partial_{\mathtt h}^k {\mathcal L}_0 = {\rm Op}(\partial_{\mathtt h}^k \Lambda(y, \xi; \mathtt h))$ and that
\[
\| \partial_{\mathtt h}^k {\mathcal L}_0 [\psi] \|_{{\mathcal H}^s}  
\simeq 
\| \partial_{\mathtt h}^k {\mathcal L}_0 [\psi] \|_{L^2_y H^s_x} 
+ 
\| \partial_{\mathtt h}^k {\mathcal L}_0 [\psi]\|_{H^s_y L^2_x}\,.
\]
Then one deduces
\[
\begin{aligned}
\| \partial_{\mathtt h}^k {\mathcal L}_0 [\psi] \|_{L^2_y H^s_x}^2 
& =  
\int_{- 1}^0 \| \partial_{\mathtt h}^k {\mathcal L}_0 [\psi] \|_{H^s_x}^2\, d y 
\\& 
= \int_{- 1}^0 \sum_{\xi \in \Gamma^*} \langle \xi \rangle^{2 s} 
|\partial_{\mathtt h}^k \Lambda(y, \xi; \mathtt h)|^2 |\widehat \psi(\xi)|^2\, d y 
\\& 
\stackrel{\eqref{stima cal L0 simbolo}}{\lesssim_k} 
\int_{- 1}^0 \sum_{\xi \in \Gamma^*} \langle \xi \rangle^{2 (s + k)} 
|\widehat \psi(\xi)|^2\, d y  
\lesssim_k \| \psi \|_{H^{s + k}_x}^2
\end{aligned}
\]
and 
\[
\begin{aligned}
\| \partial_{\mathtt h}^k {\mathcal L}_0 [\psi]\|_{H^s_y L^2_x}^2 
& \lesssim_s \sum_{0 \leq \alpha \leq s} 
\| \partial_{\mathtt h}^k \partial_y^\alpha {\mathcal L}_0 [\psi] \|_{L^2_y L^2_x}^2 
\\& 
\lesssim_s \sum_{0 \leq \alpha \leq s} \int_{- 1}^0 \sum_{\xi \in \Gamma^*} 
|\partial_{\mathtt h}^k \partial_y^\alpha \Lambda(y, \xi; \mathtt h)|^2 |\widehat \psi(\xi)|^2\, d y 
\\& 
\stackrel{\eqref{stima cal L0 simbolo}}{\lesssim_{s, k}} \sum_{0 \leq \alpha \leq s} 
\int_{- 1}^0 \sum_{\xi \in \Gamma^*}   
\langle \xi \rangle^{2(\alpha + k)}|\widehat \psi(\xi)|^2\, d y \lesssim_{s, k} \| \psi \|_{H^{s + k}_x}^2
\end{aligned}
\]
and so the  claimed bound \eqref{stima cal L0 nel lemma cal Hs}  follows. 
The bound \eqref{stima cal L0 nel lemma cal Os} follows by a direct easy computation using that 
the operator $\cL_0$ does not depend on $\vphi \in \T^\nu$, namely $\partial_\vphi^\beta \cL_0 = 0$ for any $\beta \in \N^\nu$. 
\end{proof}

In the next lemma we state an estimate of the Green operator associated to the elliptic problem
\begin{equation}\label{freeProbg}
\begin{cases}
L u  = g, \quad L  \equiv L(\mathtt h) := \partial_{yy} + \mathtt h^{2} \Delta  \\
u(0, x) = 0 \\
\partial_y u(- 1, x) = 0
\end{cases}
\end{equation}
in the case where $g \equiv g(x, y) \in {\mathcal H}^s$. We call such an operator $L^{- 1}$.

\begin{lemma}\label{lemma laplace forzato dip da vphi}
For any $s \geq 0$ and any $k\in \mathbb{N}$, one has that the unique solution $u = L^{- 1} g$ 
of the problem \eqref{freeProbg} satisfies the elliptic estimate 
\[
\| \partial_{\mathtt h}^k u \|_{\cH^{s + 2}} \lesssim_{s, k} \| g \|_{\cH^s}\,.
\] 
\end{lemma}

\begin{proof}
For $k = 0$ it follows by standard elliptic regularity. Then one argues by induction. 
%\red{METTERE DIMOSTRAZIONE FORSE}
%For completeness we refer to Appendix \ref{lemmaElliptic regularity}.
\end{proof}

\begin{rmk}
By Lemma \ref{lemma laplace forzato dip da vphi} one deduces 
$\| u \|_{H^s_\vphi \cH^2} \lesssim \| g \|_{H^s_\vphi L^2_{x, y}}$ and $\| u \|_{\cO^s} \lesssim_s \| g \|_{\cO^s}$.
\end{rmk}

\noindent
The solution $\phi(x,y)$\footnote{Actually  $\phi(x,y)\equiv\phi(\vphi;x,y)$ for $\vphi\in \T^{\nu}$ since the functions
$\alpha_{i}$, $i=1,2,3$ depend on $\vphi$ through the embedding $\eta(x)=\eta(\vphi,x)$ (see assumption \ref{hypo:eta}.}
of the problem \eqref{elliptic2}  can be written as 
\begin{equation}\label{solphiphi}
\phi (x, y) = \vphi_0(x, y) + u(x, y)\,,
\end{equation}
where $\vphi_0(x,y)$ solves \eqref{problema vphi 0 riscalato laplace}
while 
$u(x, y)$ solves the problem
\begin{equation}\label{prob omogeneop forzato}
\begin{cases}
\partial_y^2 u + \th^2 \Delta u + \widetilde{F}[u]= f \\
u(x, 0) = 0 \\
\partial_y u(x,- 1) = 0
\end{cases}
\end{equation}
where we defined the differential operator 
\begin{equation}\label{def op F vphi}
\widetilde{F}:=\widetilde{F}(\vphi)[\cdot ] := 
\big( \alpha_1(\vphi, x, y) \cdot \nabla + \alpha_2( \vphi, x, y) \big) \partial_y   +\th^2 (\alpha_3(\vphi, x, y) -1)\Delta  \,,
\end{equation}
and the ``error''
\begin{equation}\label{def: errore effino}
\begin{aligned}
& f := - \widetilde{F}[\vphi_0] =-\widetilde{\mathtt{L}}[\psi]\,,
%- \Big( \alpha_1( x, y) \cdot \nabla + \alpha_2( x, y) \Big) \partial_y \vphi_0  - (\alpha_3(\, x, y) -1)\Delta\vphi_0 \,.
\end{aligned}
\end{equation}
%
%\begin{rmk}
%Again we remark that  the operator $\widetilde{F}$ is actually a family of operators
%$\mathbb{T}^{\nu}\ni\vphi\mapsto \widetilde{F}[\cdot]=\widetilde{F}(\vphi)[\cdot]=\widetilde{F}(\vphi;x,y)[\cdot]$.
%\end{rmk}
where recalling \eqref{def vphi 0} we defined
\begin{equation}\label{operatoLLLL}
\widetilde{\mathtt{L}}:=\widetilde{\mathtt{L}}(\vphi,y):=\widetilde{F}\circ\mathcal{L}_0\,.
\end{equation}

At least at formal level, one can note that $u$ solves \eqref{prob omogeneop forzato} if and only if 
\begin{equation}\label{propagKK}
u = \cK(\vphi, y)[\psi] \quad \text{where} \quad \cK(\vphi, y) :=  
\Big( {\rm Id} + L^{- 1} \circ \widetilde{F}(\vphi, y) \Big)^{- 1} \circ L^{- 1} 
\circ \widetilde{\mathtt{L}}(\vphi,y)[\cdot]\,,
\end{equation}
where $L^{-1}$ is the Green operator for the problem \eqref{freeProbg}.
The main result of this section is to prove the following estimates on ${\mathcal{K}}(\vphi,y)$ in \eqref{propagKK}, which is actually the Green operator of the elliptic problem \eqref{prob omogeneop forzato}. 
\begin{prop}\label{stima cal Kn equazione di laplace}
Let $k \in \N^{\nu + 1}$, $|k| \leq k_0$. Then there exists 
$\sigma = \sigma(k_0)  \gg 0$ 
large enough such that for any $s\geq  2$, $\beta\in \N^{\nu}$ 
there exists $\delta=\delta(s,\beta,k_0)>0$ small 
such that, 
if $\| \eta \|_{s_0 + \sigma}^{k_0, \gamma} \leq \delta$,
%then for any $s \geq 2$, 
%$\beta \in \N^\nu$ there exists $C(s) \gg 0$ large enough such that 
then the operator $\cK(\vphi, y)$ 
(recall \eqref{propagKK}) satisfies 
\begin{equation}\label{bottle10}
\|  \partial_\lambda^k \partial_\vphi^\beta  \cK(\vphi, y)  [\psi] \|_{\cH^s} 
\lesssim_{s, k, \beta}  
\gamma^{- |k|} \Big( \| \eta\|_{s_0 + |\beta| + \sigma}^{k_0, \gamma} 
\| \psi \|_{H^{s + |k|}_x} + \| \eta \|_{s + |\beta| + \sigma}^{k_0, \gamma} \| \psi \|_{H^{ |k| + 2}_x} \Big)
\end{equation}
and for any $m = 0,1,2$, 
\begin{equation}\label{bottle11}
\begin{aligned}
\| \partial_{x, y}^m \partial_\lambda^k \partial_\vphi^\beta   &\cK(\vphi, y)  [\psi] \|_{\cH^s} 
\\&\lesssim_{s,\beta,k} \gamma^{- |k|} \Big(   \| \eta\|_{s_0 + |\beta| + \sigma}^{k_0, \gamma} \| \psi \|_{H^{s + |k| + m}_x} 
+ \| \eta \|_{s + m + |\beta| + \sigma}^{k_0, \gamma} \| \psi \|_{H^{ |k| + 2}_x} \Big)\,.
\end{aligned}
\end{equation}
Moreover, one has
\begin{equation}\label{bottle12}
\| \partial_\vphi^\beta\Delta_{12} \cK(\vphi, y)\|_{\mathcal{L}(\mathcal{H}^{p+2+|k|};\mathcal{H}^{p})}
\lesssim_{p,k,\beta} \|\eta_{1}-\eta_{2}\|_{p+\widetilde{\s}+|\beta|}\,,\;\;\;s_0+2\leq p\leq s_0 + \s-\widetilde{\s}\,.
\end{equation}
for some $\widetilde{\s}\leq \s$. % and for $p+\widetilde{\s}\leq s_0+\s$, $p\geq s_0$.
\end{prop}
The proof of the Proposition above involves several arguments. In subsection \ref{sec:stimeforzante Elle tilde}
we first give estimates on the forcing term $f$ in \eqref{def: errore effino}.
Then in subsection \ref{sec:stimegreen} we conclude the proof of Proposition \ref{stima cal Kn equazione di laplace}.

\subsubsection{Estimates on the forcing term f}\label{sec:stimeforzante Elle tilde}
In this section we prove the following.
\begin{prop}\label{stime striscia mathtt L vphi y O s}
 Let $k \in \N^{\nu + 1}$, $|k| \leq k_0$
 and consider the operator $\widetilde{\mathtt{L}}=\widetilde{\mathtt{L}}(\vphi,y)$ in \eqref{operatoLLLL}. 
There exists $\sigma \equiv \sigma(k_0)  \gg 0$ such that, for any $\beta \in \N^\nu$, $s\geq 0$, 
there exists $\delta=\delta(s,\beta,k_0)>0$ small 
such that, 
if $\| \eta \|_{s_0 + \sigma}^{k_0, \gamma} \leq \delta$,
then
one has the estimates
\begin{align}
\!\!\! \| \partial_\lambda^k \partial_\vphi^\beta \widetilde{\mathtt{L}}(\vphi, y) [\psi] \|_{\cH^s} 
&\lesssim_{s,k, \beta} 
\gamma^{- |k|} \Big(  \| \eta \|_{s_0 + |\beta| + \sigma }^{k_0,\gamma} \| \psi \|_{H^{s + |k| +  2}_x} 
+ \| \eta \|_{s + |\beta| + \sigma}^{k_0,\gamma} \| \psi \|_{H^{ |k| + 2}_x} \Big)\,,\label{bottle1}
\\
\!\!\!\| \partial_\vphi^\beta\Delta_{12} \widetilde{\mathtt{L}}(\vphi, y)&\|_{\mathcal{L}(\mathcal{H}^{p+2};\mathcal{H}^{p})}
 \lesssim_{p,\beta} 
\|\eta_{1}-\eta_{2}\|_{p+\widetilde{\s}+|\beta|}\,,\;\;\;0\leq p\leq s_0+\s-\widetilde{\s}\,,\label{bottle2}
\end{align}
and 
\begin{align}
\| \partial_\lambda^k \partial_\vphi^\beta \widetilde{\mathtt{L}}(\vphi, y) \langle D \rangle^{- |k| - 2}  [\psi] \|_{\cO^s} 
&\lesssim_{s, k, \beta } 
\gamma^{- |k|} \Big(  \| \eta \|_{s_0 + |\beta| + \sigma }^{k_0,\gamma} \| \psi \|_{s } 
+ \| \eta \|_{s + |\beta| + \sigma}^{k_0,\gamma} \| \psi \|_{0 } \Big)\,,\label{bottle1A}
\\
\| \partial_\vphi^\beta\Delta_{12} \widetilde{\mathtt{L}}(\vphi, y) \langle D \rangle^{- 2}\|_{\mathcal{L}(\mathcal{O}^{p};\mathcal{O}^{p})}
& \lesssim_{p,\beta} 
\|\eta_{1}-\eta_{2}\|_{p+\widetilde{\s}+|\beta|}\,,%\;\;\;p\geq 0\,.
\;\;\;0\leq p\leq s_0+\s-\widetilde{\s}\,,
\label{bottle2A}
\end{align}
for some $\widetilde{\s}\leq \s$. % and for $p+\widetilde{\s}\leq s_0+\s$, $p\geq0$.
\end{prop}
In order to prove the result above we first need some preliminary results.
In particular, we provide estimates on
the forcing term  in \eqref{def op F vphi}.

\begin{lemma}\label{F eta vphi 0 L eta psi}
 Let $\widetilde F$ be the linear operator defined in \eqref{def op F vphi}. Let $k \in \N^{\nu + 1}$, $|k| \leq k_0$. Then there exists 
$\sigma = \sigma(k_0)  \gg 0$ 
large enough such that for any $s\geq  0$, $\beta\in \N^{\nu}$ there exists $\delta=\delta(s,\beta,k_0)>0$ small 
such that, 
if 
$\| \eta \|_{s_0 + \sigma}^{k_0, \gamma} \leq \delta$, 
%$\| \eta \|_{s_0 + \sigma}^{k_0, \gamma} \leq \delta$, 
%then for any $s \geq 2$, 
%$\beta \in \N^\nu$ there exists $C(s) \gg 0$ large enough such that 
then the operator $\widetilde F(\vphi, y)$ 
%(recall \eqref{propagKK}) 
satisfies 
\begin{equation}\label{bottle10 tilde F}
\begin{aligned}
\|  \partial_\lambda^k \partial_\vphi^\beta  \widetilde F(\vphi, y)  [u] \|_{\cH^s} 
& \lesssim_{s,k,\beta}  
\gamma^{- |k|} \Big( \| \eta\|_{s_0 + |\beta| + \sigma}^{k_0, \gamma} 
\| u \|_{{\mathcal H}^{s + 2}} + \| \eta \|_{s + |\beta| + \sigma}^{k_0, \gamma} \| u \|_{{\mathcal H}^{ 2}} \Big)\,, \\
\| \partial_\vphi^\beta\Delta_{12} \widetilde{F}(\vphi, y)\|_{\mathcal{L}(\mathcal{H}^{p+2};\mathcal{H}^{p})} & 
\lesssim_{p,\beta} 
\|\eta_{1}-\eta_{2}\|_{p+\widetilde{\s}+|\beta|}\,,%\;\;\;p\geq 0\,.
\;\;\;0\leq p\leq s_0+\s-\widetilde{\s}\,,
\end{aligned}
\end{equation}
for some $\widetilde{\s}\leq \s$. %and for $p+\widetilde{\s}\leq s_0+\s$, $p\geq0$.
\end{lemma}

\begin{proof}
By the definition of $\widetilde F$ in \eqref{def op F vphi}, one has that 
$$
\begin{aligned}
& \widetilde F = {\mathcal R}_1 + {\mathcal R}_2 + {\mathcal R}_3\,, \\
& {\mathcal R}_1 := \alpha_1 \cdot \nabla \partial_y \,, \quad {\mathcal R}_2 := \alpha_2 \partial_y \,, \quad {\mathcal R}_3 := \th^2(\alpha_3 - 1) \Delta\,. 
\end{aligned}
$$
We estimate ${\mathcal R}_1$. The operators ${\mathcal R}_2$ and ${\mathcal R}_3$ can be estimated similarly\,. One has that 
$$
\partial_\lambda^k \partial_\vphi^\beta {\mathcal R}_1 [u] = \big(\partial_\lambda^k \partial_\vphi^\beta \alpha_1 \big) \cdot \nabla \partial_y \,.
$$
By the estimates \eqref{stimealphai}, one has that there exists $\sigma \gg 0$ large enough such that for any $n \in \N$, $s \geq s_0$, 
\begin{equation}\label{stime alpha 1 partial lambda vphi}
\begin{aligned}
\| \partial_\vphi^\beta \partial_\lambda^k  \partial_y^n \alpha_1 \|_{H^s_x} & \lesssim  \|  \partial_\lambda^k  \partial_y^n \alpha_1 \|_{H^{s_0 + |\beta|}_\vphi H^s_x} \lesssim  \|  \partial_\lambda^k  \partial_y^n \alpha_1 \|_{s + s_0 + |\beta|} \\
& \lesssim_{s, \beta, k, n} \gamma^{- |k|}\| \eta \|_{s + |\beta| + \sigma}^{k_0, \gamma}
\end{aligned}
\end{equation}
for some $\sigma \gg 0$ large enough indpendent of $n$ and $\beta$. Hence 
$$
\begin{aligned}
\|\partial_\lambda^k \partial_\vphi^\beta {\mathcal R}_1 [u] \|_{{\mathcal H}^s} & 
\stackrel{\eqref{stime alpha 1 partial lambda vphi}\,,\,Lemma \ref{lemma stima tame prodotto con gli alpha i}}{\lesssim_{s, k,\beta}} 
 \gamma^{- |k|}\| \eta \|_{s + |\beta| + \sigma}^{k_0, \gamma} \|  \nabla \partial_y u \|_{L^2_{x, y}} + \gamma^{- |k|}\| \eta \|_{s_0 + |\beta| + \sigma}^{k_0, \gamma} \|  \nabla \partial_y u \|_{{\mathcal H}^{s}} \\
& \lesssim_{s, k,\beta} 
\gamma^{- |k|}\| \eta \|_{s + |\beta| + \sigma}^{k_0, \gamma} \|  u \|_{{\mathcal H}^{ 2}} + \gamma^{- |k|}\| \eta \|_{s_0 + |\beta| + \sigma}^{k_0, \gamma} \|  u \|_{{\mathcal H}^{s + 2}}\,. 
\end{aligned}
$$
The estimate for $\partial_\vphi^\beta \Delta_{12} {\mathcal R}_1$ can be proved similarly. 
\end{proof}

\begin{proof}[{\bf Proof of Proposition \ref{stime striscia mathtt L vphi y O s}}]
By formula \eqref{operatoLLLL} one gets directly that
%By Lemma \ref{F eta vphi 0 L eta psi}, one gets directly that  
\[
\begin{aligned}
\partial_\lambda^k \partial_\vphi^\beta \widetilde{\mathtt{L}}(\vphi, y) & = 
\sum_{k_1 +k_2=k}
 (\partial_\lambda^{k_1} \partial_\vphi^\beta \widetilde F) \circ \partial_\lambda^{k_2} {\mathcal L}_0\,.
\end{aligned}
\]
The claimed bounds \eqref{bottle1}, \eqref{bottle2} then easily follow by Lemmata \ref{lemma stima cal L0 sol omogenea laplace}, \ref{F eta vphi 0 L eta psi}. 

\noindent
We shall prove the estimate \eqref{bottle1A}. 

To shorten notations in the proof we write $\| \eta \|_s$ instead of $\| \eta \|_s^{k_0, \gamma}$. 
By Lemma \ref{algebra striscia},
%Remark \ref{rmk:equiOHs}, 
it 
is enough to provide an estimate of 
\begin{align}
\|  \partial_\lambda^k \partial_\vphi^\beta \widetilde{\mathtt{L}}(\vphi, y)
\langle D \rangle^{- |k| - 2}[\psi] \|_{\cO^s}
&  \simeq 
\|  \partial_\lambda^k \partial_\vphi^\beta \widetilde{\mathtt{L}}(\vphi, y)
\langle D \rangle^{ - |k| - 2}[\psi] \|_{L^2_\vphi \cH^s}  
\label{pizzetta1}
\\
& \quad 
+ \|   \partial_\lambda^k \partial_\vphi^\beta \widetilde{\mathtt{L}}(\vphi, y)
\langle D \rangle^{ - |k| - 2}[\psi]  \|_{H^s_\vphi L^2_{x, y}} \,.
\label{pizzetta2}
\end{align}

\medskip

\noindent
{\sc Estimate of  \eqref{pizzetta1}.
%$\|   \partial_\lambda^k \partial_\vphi^\beta \widetilde{\mathtt{L}}(\vphi, y)
%\langle D \rangle^{ - |k| - 2}[\psi] \|_{L^2_\vphi \cH^s} $.
}  
By using the estimate \eqref{bottle1} 
one gets, for any $s \geq 0 $, that
\begin{equation}\label{bla bla elliptic forcing 4A}
\begin{aligned}
&\| \partial_\lambda^k \partial_\vphi^\beta \widetilde{\mathtt{L}}(\vphi, y)
\langle D \rangle^{ - |k| - 2}[\psi]   \|_{L^2_\vphi \cH^{s}}  
\\& 
\lesssim_{s,k,\beta} \gamma^{- |k|} \Big(  \| \eta \|_{s_0 + |\beta| + \sigma}  \| \langle D \rangle^{ - |k| - 2}\psi \|_{L^2_\vphi H^{s + |k| + 2}_x}  
%\\& \qquad 
+ \| \eta \|_{s  + |\beta|  + \sigma}  \| \langle D \rangle^{ - |k| - 2}\psi \|_{L^2_\vphi H^{ |k| + 2}_x}\Big) 
%\\& 
%\lesssim_s \gamma^{- |k|} \Big(  \| \eta \|_{s_0 + |\beta| + \sigma}  \| \psi \|_{L^2_\vphi H^{s }_x}  
%%\\& \qquad 
%+ \| \eta \|_{s  + |\beta|  + \sigma}  \| \psi \|_{L^2_\vphi L^2_x}\Big) 
\\
& \stackrel{\eqref{spaziosoloX},\eqref{sobnormseparati}}{\lesssim_{s,k,\beta}} 
\gamma^{- |k|} \Big(  \| \eta \|_{s_0 + |\beta| + \sigma}  \| \psi \|_{s}   
+ \| \eta \|_{s  + |\beta|  + \sigma}  \| \psi \|_{0}\Big)\,. 
\end{aligned}
\end{equation}

\smallskip

\noindent
{\sc Estimate of \eqref{pizzetta2}.
%$\|  \partial_\lambda^k \partial_\vphi^\beta \widetilde{\mathtt{L}}(\vphi, y)
%\langle D \rangle^{ - |k| - 2}[\psi]  \|_{H^s_\vphi L^2_{x, y}} $.
} 
We first note that 
\begin{equation}\label{bla bla elliptic forcing 0A}
\begin{aligned}
&\|  \partial_\lambda^k \partial_\vphi^\beta \widetilde{\mathtt{L}}(\vphi, y)\langle D \rangle^{ - |k| - 2}[\psi]  \|_{H^s_\vphi L^2_{x, y}} 
 \lesssim_s 
\sum_{|\alpha| \leq s} \Big\| \partial_\vphi^\alpha \Big(   \partial_\lambda^k \partial_\vphi^\beta \widetilde{\mathtt{L}}(\vphi, y)\langle D \rangle^{ - |k| - 2}[\psi]  \Big) 
\Big\|_{L^2_\vphi L^2_{x, y}} 
\\& \qquad\quad
\lesssim_s 
\sum_{|\alpha| \leq s} \sum_{\alpha_1 + \alpha_2 = \alpha} 
\Big\| \partial_\lambda^k \partial_\vphi^{\alpha_1 + \beta}  \widetilde{\mathtt{L}}(\vphi, y) [\langle D \rangle^{ - |k| - 2}\partial_\vphi^{\alpha_2} \psi] \Big\|_{L^2_\vphi L^2_{x, y}}\,.
\end{aligned}
\end{equation}
Secondly, applying again the estimate \eqref{bottle1}, 
for any $|\alpha| \leq s$, $\alpha_1 + \alpha_2 = \alpha$, one gets that
\begin{equation}\label{bla bla elliptic forcing 1A}
\begin{aligned}
 \Big\| \partial_\lambda^k \partial_\vphi^{\alpha_1 + \beta}  &\widetilde{\mathtt{L}}(\vphi, y) [\langle D \rangle^{ - |k| - 2}\partial_\vphi^{\alpha_2} \psi] \Big\|_{L^2_\vphi L^2_{x, y}}
%& 
%\lesssim C(s)^n \gamma^{- |k|} \Big( \| \eta \|_{s_0 + \sigma}^{n - 1} \| \eta \|_{s_0 + |\alpha_1| + |\beta| + \sigma}  
%\| \partial_\vphi^{\alpha_2} f \|_{L^2_\vphi L^2_{x, y}}  
%\\& \qquad 
%+ \| \eta \|_{|\alpha_1| + |\beta|  + \sigma} \| \eta \|_{s_0 + \sigma}^{n - 1} \| \partial_\vphi^{\alpha_2} f \|_{L^2_\vphi L^2_{x, y}} \Big) \\
 \\
 & \lesssim_{k,\beta}  \gamma^{- |k|}   \| \eta \|_{s_0 + |\alpha_1| + |\beta| + \sigma}  \|  \langle D \rangle^{ - |k| - 2} \psi \|_{H^{|\alpha_2|}_\vphi H^{ |k| + 2}_x} \\
 & \lesssim_{k, \beta} \gamma^{- |k|}   \| \eta \|_{s_0 + |\alpha_1| + |\beta| + \sigma}  \|   \psi \|_{H^{|\alpha_2|}_\vphi L^2_x}\,.
\end{aligned}
\end{equation}
By interpolation\footnote{Here we are using interpolation estimates of Lemma \ref{lemma:interpolation}. 
We remark that such estimates holds also on the space $H_{\vphi}^{s}$ of functions of the $\vphi\in \T^{\nu}$ variables.} and by the Young inequality\footnote{recall that for $a,b\geq0$ one has
\begin{equation}\label{centoYoung}
ab\leq \frac{a^{p}}{p}+\frac{b^q}{q}\,,\quad p,q>1\;\;\; {\rm with}\;\;\; \frac{1}{p}+\frac{1}{q}=1\,.
\end{equation}} 
with exponents $|\alpha|/|\alpha_1|$ and $|\alpha|/|\alpha_2|$,
%$a b \lesssim a^{\frac{|\alpha|}{|\alpha_1|}} + b^{\frac{|\alpha|}{|\alpha_2|}}$, 
one gets
\begin{equation}\label{bla bla elliptic forcing 2A}
\begin{aligned}
\| \eta \|_{s_0 + |\alpha_1| + |\beta| + \sigma} & \|  \psi \|_{H^{|\alpha_2|}_\vphi L^2_{x}} 
\\& \lesssim 
\Big( \| \eta \|_{s_0 + |\beta| + \sigma} \| \psi \|_{H^{|\alpha|}_\vphi L^2_x} \Big)^{\frac{|\alpha_2|}{|\alpha|}} 
\Big( \| \eta \|_{s_0 + |\alpha| + |\beta| + \sigma} \| \psi \|_{L^2_\vphi L^2_x} \Big)^{\frac{|\alpha_1|}{|\alpha|}} 
\\& 
\lesssim 
\| \eta \|_{s_0 + |\beta| + \sigma} \| \psi \|_{H^{|\alpha|}_\vphi L^2_x}  
+ \| \eta \|_{s_0 + |\alpha| + |\beta| + \sigma} \| \psi \|_{L^2_\vphi L^2_x} 
\\& 
\lesssim \| \eta \|_{s_0 + |\beta| + \sigma} \| \psi \|_{s}  + \| \eta \|_{s +  |\beta| + \sigma} \| \psi \|_{0}\,,
\end{aligned}
\end{equation}
for some $\sigma \gg 0$ large enough. 
Hence, by \eqref{bla bla elliptic forcing 0A}, \eqref{bla bla elliptic forcing 1A}, \eqref{bla bla elliptic forcing 2A}, one gets that for any $s \geq s_0$
\begin{equation}\label{bla bla elliptic forcing 3A}
\begin{aligned}
 \| \partial_\lambda^k \partial_\vphi^{ \beta}  &\widetilde{\mathtt{L}}(\vphi, y) [\langle D \rangle^{ - |k| - 2}\psi]  
 \|_{H^s_\vphi L^2_{x, y}} 
\\
&\lesssim_{s, k ,\beta} 
 \gamma^{- |k|}   \Big(  \| \eta \|_{s_0 + |\beta| + \sigma} \| \psi \|_{s}  
+ \| \eta \|_{s +  |\beta| + \sigma} \| \psi \|_{0} \Big)\,.
\end{aligned}
\end{equation}
The estimates  \eqref{bla bla elliptic forcing 4A}, \eqref{bla bla elliptic forcing 3A} 
imply the claimed bound \eqref{bottle1A}. The bound \eqref{bottle2A} can be proved by similar arguments, using also the estimate \eqref{bottle2}. 
\end{proof}

\subsubsection{Space-time Green's function estimates 1}\label{sec:stimegreen}

In this section we  provide tame estimates for the operator $\cK(\vphi, y)$ 
and we conclude the proof of 
Proposition \ref{stima cal Kn equazione di laplace}.
First of all,
using the Neumann series, we shall write the operator $\cK(\vphi,y)$ in \eqref{propagKK} as
(recall \eqref{def op F vphi}, \eqref{operatoLLLL}) 
\begin{equation}\label{propagatore cal K cal Kn vphi}
\begin{aligned}
\cK(\vphi, y) & = \sum_{n \geq 0} (- 1)^n \cK_n(\vphi, y)\,,  \\
\cK_n(\vphi, y) & :=\cA_n(\vphi, y) \circ \widetilde{\mathtt{L}}(\vphi, y)\,, \quad \forall n \geq 0\,,
%\\ A(\vphi, y )&:= L^{- 1} \circ \widetilde{F}(\vphi, y)\,.
\end{aligned}
\end{equation}
where,
for $n \in \N$ we defined 
 the operators 
\begin{equation}\label{AAnn}
\cA_n(\vphi, y) := 
(L^{- 1} \circ \widetilde{F}(\vphi, y))^n \circ L^{- 1}\,,\;\;n\geq1\,,\quad \cA_0(\vphi,y):= L^{- 1}\,.
\end{equation}

We need the following lemma.
\begin{lemma}\label{Corollario stima L inv F 0}
 Under the assumptions of Proposition \ref{stima cal Kn equazione di laplace}
 we have that the operators  $\cA_n(\vphi, y)$ in \eqref{AAnn} satisfy, for $n\geq1$, $s \geq 2$, $\beta\in\N^{\nu}$,
\[
\begin{aligned}
\| \partial_\lambda^k \partial_\vphi^\beta& \cA_n(\vphi, y) [ f] \|_{\cH^s} 
\\&\leq 
\gamma^{- |k|}C(s)^n\Big( (\| \eta \|_{s_0 + \sigma}^{k_0,\gamma})^{n - 1} \| \eta \|_{s_0 + |\beta| + \sigma}^{k_0,\gamma}  
\| f \|_{\cH^{s - 2}} 
+ \| \eta \|_{s  + |\beta|  + \sigma}^{k_0,\gamma} (\| \eta \|_{s_0 + \sigma}^{k_0,\gamma})^{n - 1} \| f \|_{L^2_{x, y}}\Big)\,,
\end{aligned}
\]
for some $C(s)=C(s,\beta,k)>0$, 
and
\begin{equation*}
\| \partial_\vphi^\beta\Delta_{12} \mathcal{A}_{n}(\vphi, y)\|_{\mathcal{L}(\mathcal{H}^{p-2};\mathcal{H}^{p})}
\leq C(p)^n (\| \eta_1 \|_{s_0 + \widetilde{\s}+|\beta|}
+\| \eta_2 \|_{s_0 +\widetilde{\s}+|\beta|})^{n - 1}\|\eta_{1}-\eta_{2}\|_{p+\widetilde{\s}+|\beta|}\,,%\;\;\;p\geq 2\,.
\end{equation*}
for some $C(p)=C(p,\beta)>0$ and 
some $\widetilde{\s}\leq \s$ and for $p+\widetilde{\s}\leq s_0+\s$, $p\geq2$.
\end{lemma}
\begin{proof}
We start by proving that the operators
\[
A_{ n}(\vphi, y ) :=  (A(\vphi, y))^n \,,\quad n\geq1\,,
\qquad A(\vphi,y):=L^{- 1} \circ \widetilde{F}(\vphi, y)\,,
\]

satisfy  for any $n\geq1$ and any $s\geq  2$,
\begin{align}
\| & \partial_\lambda^k \partial_\vphi^\beta A_{ n}(\vphi, y)[f] \|_{\cH^s} \nonumber
\\&\leq 
\gamma^{- |k|}C(s)^n \Big( (\| \eta \|_{s_0 + \sigma}^{k_0,\gamma})^{n - 1} \| \eta \|_{s_0 + |\beta| + \sigma}^{k_0,\gamma} \| f \|_{\cH^s} 
+ \| \eta \|_{s + |\beta| + \sigma}^{k_0,\gamma} (\| \eta \|_{s_0 + \sigma}^{k_0,\gamma})^{n - 1} \| f \|_{\cH^{ 2}} \Big)\,,
\label{abruzzo3}
\\
\|& \partial_\vphi^\beta\Delta_{12} A_{n}(\vphi, y)\|_{\mathcal{L}(\mathcal{H}^{p};\mathcal{H}^{p})}\nonumber
\\&\leq C(p)^n (\| \eta_1 \|_{s_0 + \widetilde{\sigma}+|\beta|}
+\| \eta_2 \|_{s_0 + \widetilde{\sigma}+|\beta|})^{n - 1}\|\eta_{1}-\eta_{2}\|_{p+\widetilde{\s}+|\beta|}\,,
\;\;\;p\geq  2\,,
\label{AnDelta12}
\end{align}
for some $\s\gg1$ depending on $k_0$ and some $\widetilde{\s}\leq \s$ independent of $k_0$.

\noindent
To shorten notations in the proof, we write $\| \eta \|_s$ instead of $\| \eta \|_s^{k_0, \gamma}$. 
We argue by induction on $n$.

\smallskip
\noindent
{\sc Case $n=1$.} The claimed bound for $n = 1$ follows by a direct application of Lemmata \ref{lemma laplace forzato dip da vphi}, \ref{F eta vphi 0 L eta psi}

\smallskip
\noindent
{\sc Induction step.} 
%We argue by induction on $n$. For $n = 1$ the statement follows by the item $(i)$. 
Now assume that the claimed statement holds for $n \geq 1$ and let us prove it for $n + 1$. We first note that
\begin{equation*}
\begin{aligned}
  \partial_\lambda^k \partial_\vphi^\beta  A(\vphi, y)^{n + 1} &= 
  \partial_\lambda^k \partial_\vphi^\beta \Big( A(\vphi, y)^n A(\vphi, y) \Big)    
  \\&=  \sum_{\begin{subarray}{c}
\beta_1 + \beta_2 =   \beta \\
k_1 + k_2 = k 
\end{subarray}} C(\beta_1, \beta_2, k_1, k_2) \partial_\lambda^{k_1}\partial_\vphi^{\beta_1} A(\vphi,y)^n \circ \partial_\lambda^{k_2}\partial_\vphi^{\beta_2} A(\vphi, y) \,,
\end{aligned}
\end{equation*}
for some constants $C(\beta_1, \beta_2, k_1, k_2)>0$.
By the item $(i)$ and by the inductive hypothesis, 
for any $\beta_1, \beta_2 \in \N^\nu$, $k_1, k_2 \in \N^{\nu + 1}$, $|k| \leq k_0$ with $\beta_1 + \beta_2 = \beta$, $k_1 + k_2 = k$, $s \geq  2$,
one obtains the estimate 
\begin{equation}\label{partial vphi beta 1 beta 2 A n A}
\begin{aligned}
\!\!\! \| \partial_\lambda^{k_1}&\partial_\vphi^{\beta_1} A(\vphi,y)^n \circ  
 \partial_\lambda^{k_2}\partial_\vphi^{\beta_2} A(\vphi, y) [f] \|_{\cH^s}  
 \\&  
 \stackrel{\eqref{strippoperlalgebra},\eqref{abruzzo3}_n}{\leq} \gamma^{- |k_1|}C(s)^n 
 \| \eta \|_{s_0 + \sigma}^{n - 1} \| \eta \|_{s_0 + |\beta_1| + \sigma} 
 \| \partial_\lambda^{k_2}\partial_\vphi^{\beta_2} A(\vphi, y) [f] \|_{\cH^s} 
 \\&\qquad 
 + \gamma^{- |k_1|}C(s)^n 
 \| \eta \|_{s + |\beta_1| + \sigma} \| \eta \|_{s_0 + \sigma}^{n - 1} 
 \| \partial_\lambda^{k_2} \partial_\vphi^{\beta_2} A(\vphi, y) [f] \|_{\cH^{ 2}} 
 \\& 
 \stackrel{\eqref{abruzzo3}_1}{\lesssim_s} \gamma^{- |k|} C(s)^n   
 \| \eta \|_{s_0 + \sigma}^{n - 1} \| \eta \|_{s_0 + |\beta_1| + \sigma} \| \eta \|_{s_0 + |\beta_2| + \sigma} \| f \|_{\cH^s}   
 \\&  
 + \gamma^{- |k|}C(s)^n \| \eta \|_{s_0 + \sigma}^{n - 1} \Big(   \| \eta \|_{s_0 + |\beta_1| + \sigma} \| \eta \|_{s + |\beta_2| + \sigma} +   \| \eta \|_{s + |\beta_1| + \sigma} \| \eta \|_{s_0 + |\beta_2| + \sigma} \Big) \| f \|_{\cH^{ 2}}
 \\&
   \leq \gamma^{- |k|} C(s)^{n + 1} \| \eta \|_{s_0 + \sigma}^{n }\big(  \| \eta \|_{s_0 + |\beta| + \sigma}  \| f \|_{\cH^s} 
 + \| \eta \|_{s + |\beta| + \sigma} \| f \|_{\cH^{ 2}} \big) \,, 
\end{aligned}
\end{equation}
by taking $C(s)=C(s, k , \beta) \gg 0$ large enough, where in the last inequality we used the 
interpolation estimate \eqref{2202.2}\footnote{In particular, we used \eqref{2202.2} with 
$a_0=b_0=s_0+\s$ (or $b_0=s+\s$)
and $p=|\beta_1|, q=|\beta_2|$ to obtain 
\[
\begin{aligned}
 \| \eta \|_{s_0 + |\beta_1| + \sigma} \| \eta \|_{s_0 + |\beta_2| + \sigma}&\lesssim 
 \| \eta \|_{s_0 + \sigma} \| \eta \|_{s_0 + |\beta| + \sigma}
 \\
 \| \eta \|_{s_0 + |\beta_1| + \sigma} \| \eta \|_{s + |\beta_2| + \sigma}&\lesssim 
 \| \eta \|_{s_0 + \sigma} \| \eta \|_{s + |\beta| + \sigma} \,.
 \end{aligned}
\]
 }.
The estimate above is \eqref{abruzzo3} at the step $n + 1$, then the thesis follows. 
In order to estimate the variation in $\eta$ one can note that
\begin{equation*}
\begin{aligned}
\partial_\vphi^\beta \Delta_{12}  A(\vphi, y)^{n + 1} 
%&= 
%\partial_\vphi^\beta \Big( \Delta_{12}A(\vphi, y)^n A(\vphi, y) +A(\vphi, y)^n \Delta_{12}A(\vphi, y) \Big)    
%  \\
  &=  
  \sum_{\begin{subarray}{c}
\beta_1 + \beta_2 =   \beta 
\end{subarray}} C(\beta_1, \beta_2) \partial_\vphi^{\beta_1} \Delta_{12}A(\vphi,y)^n \circ \partial_\vphi^{\beta_2} A(\vphi, y) 
 \\&+  
  \sum_{\begin{subarray}{c}
\beta_1 + \beta_2 =   \beta 
\end{subarray}} C(\beta_1, \beta_2) \partial_\vphi^{\beta_1} A(\vphi,y)^n \circ \partial_\vphi^{\beta_2} \Delta_{12}A(\vphi, y) \,.
\end{aligned}
\end{equation*}
Therefore the bound \eqref{AnDelta12} with $n+1$ 
follows by the inductive assumption 
and \eqref{abruzzo3}.
We are now in position to conclude the proof of the lemma.
Recalling \eqref{AAnn} one has that 
\[
\partial_\lambda^k \partial_\vphi^\beta \cA_n(\vphi, y) 
=  \sum_{k_1 + k_2 = k} \partial_\lambda^{k_1} \partial_\vphi^\beta A(\vphi, y)^n \circ \partial_\lambda^{k_2} L^{- 1}\,.
\] 
By Lemma  \ref{lemma laplace forzato dip da vphi} 
one has $\| \partial_\lambda^{k_2} L^{- 1} f \|_{\cH^s} \lesssim_{s, k} \| f \|_{\cH^{s - 2}}$ for any $s\geq2$.
Therefore, the thesis follows using   \eqref{abruzzo3}. 
For the variation in $\eta$ one can reason similarly using \eqref{AnDelta12}.
\end{proof}

In view of Lemma \ref{Corollario stima L inv F 0} we shall conclude the proof of the main proposition.
\begin{proof}[{\bf Proof of Proposition \ref{stima cal Kn equazione di laplace}}]
To shorten notations in the proof we write $\| \eta \|_s$ instead of $\| \eta \|_s^{k_0, \gamma}$. 

\smallskip
\noindent
For any $n \geq 0$, we first estimate $\cK_n(\vphi, y)$. % in \eqref{propagatore cal K cal Kn vphi}. 
Recalling \eqref{AAnn} and \eqref{propagatore cal K cal Kn vphi} we shall write
\[
\begin{aligned}
\cK_n(\vphi, y) &= \cA_n(\vphi, y) \circ \widetilde{\mathtt{L}}(y, \vphi)\,,  
\\
 \partial_\lambda^k \partial_\vphi^\beta \cK_n(\vphi, y) &= 
 \sum_{\begin{subarray}{c}
\beta_1 + \beta_2 = \beta  \\
k_1 + k_2 = k
\end{subarray}} 
C(\beta_1, \beta_2, k_1, k_2) \partial_\lambda^{k_1}\partial_\vphi^{\beta_1}  
\cA_n(\vphi, y) \circ \partial_\lambda^{k_2}\partial_\vphi^{\beta_2}  \widetilde{\mathtt{L}}(\vphi, y)\,.
\end{aligned}
\]
Therefore, in order to prove the claimed statement, it is enough to estimate, 
for any 
$\beta_1, \beta_2 \in \N^\nu$, $k_1, k_2 \in \N^{\nu + 1}$, $\beta_1 + \beta_2 = \beta$, $k_1 + k_2 = k$, 
 the linear operator 
\begin{equation}\label{bottle15}
\begin{aligned}
\cR(\vphi, y ) & := \partial_\lambda^{k_1}\partial_\vphi^{\beta_1} \cA_n(\vphi, y) 
\circ \partial_\lambda^{k_2} \partial_\vphi^{\beta_2} \widetilde{\mathtt{L}}(y, \vphi) \,. 
\end{aligned}
\end{equation}
Then  by applying Lemma \ref{stime striscia mathtt L vphi y O s} and 
Lemma \ref{Corollario stima L inv F 0} (recall that $C(s)\equiv C(s,\beta,k)$) 
one gets for any $s \geq   2$, 
\begin{equation}\label{stima cal R per stima Kn Neumann}
\begin{aligned}
\| \cR(\vphi, y)[\psi] \|_{\cH^s} & \leq
\gamma^{- |k_1|}C(s)^n\Big( \| \eta \|_{s_0 + \sigma}^{n - 1} \| \eta \|_{s_0 + |\beta_1| + \sigma}  
\| \partial_\lambda^{k_2}\partial_\vphi^{\beta_2} \widetilde{\mathtt{L}}(\vphi, y)[\psi] \|_{\cH^{s - 2}}  
\\& 
\qquad 
+ \| \eta \|_{s  + |\beta_1|  + \sigma} \| \eta \|_{s_0 + \sigma}^{n - 1} 
\| \partial_\lambda^{k_2} \partial_\vphi^{\beta_2} \widetilde{\mathtt{L}}(\vphi, y)[\psi]  \|_{L^2_{x, y}}\Big) 
\\& 
\stackrel{\eqref{bottle1}}{\lesssim_{s,\beta,k}} 
\gamma^{- |k|}C(s)^n \| \eta \|_{s_0 + \sigma}^{n - 1} \| \eta \|_{s_0 + |\beta_1| + \sigma} 
\| \eta \|_{s_0 + |\beta_2| + \sigma} \| \psi \|_{H^{s + |k|}_x}  
\\& 
\qquad 
+ \gamma^{- |k|}C(s)^n \| \eta \|_{s_0 + \sigma}^{n - 1} 
\| \eta \|_{s + |\beta_1| + \sigma} \| \eta \|_{s_0 + |\beta_2| + \sigma} \| \psi \|_{H^{ |k| + 2}_x}
\\& 
\qquad 
+ \gamma^{- |k|}C(s)^n \| \eta \|_{s_0 + \sigma}^{n - 1}  
\| \eta \|_{s_0 + |\beta_1| + \sigma} \| \eta \|_{s + |\beta_2| + \sigma}  \| \psi \|_{H^{|k| + 2}_x}\,. 
\end{aligned}
\end{equation}
%Arguing as in \eqref{trick assurdi 0}-\eqref{trick assurdi 2}, 
%in the proof of Lemma \ref{Corollario stima L inv F 0}, one can show that 
%\begin{equation*}
%\begin{aligned}
% \| \eta \|_{s_0 + |\beta_1| + \sigma} \| \eta \|_{s_0 + |\beta_2| + \sigma} 
% &\lesssim 
% \| \eta \|_{s_0 + \sigma} \| \eta \|_{s_0 + |\beta| + \sigma}, 
% \\
% \| \eta \|_{s + |\beta_1| + \sigma} \| \eta \|_{s_0 + |\beta_2| + \sigma} 
%&\lesssim 
%\| \eta \|_{s_0 + \sigma} \| \eta \|_{s + |\beta| + \sigma}  \,,
%\\  
%\| \eta \|_{s_0 + |\beta_1| + \sigma} \| \eta \|_{s + |\beta_2| + \sigma} 
%&\lesssim \| \eta \|_{s_0 + \sigma} \| \eta \|_{s + |\beta| + \sigma}\,,
%\end{aligned}
%\end{equation*}
%which, together with the estimate \eqref{stima cal R per stima Kn Neumann}, 
By the interpolation estimate \eqref{2202.2} (arguing as in \eqref{partial vphi beta 1 beta 2 A n A}
in the proof of Lemma \ref{Corollario stima L inv F 0})
and by \eqref{stima cal R per stima Kn Neumann}
one deduces that
the operator $\cK_n$ satisfies  
\[
\begin{aligned}
\| \partial_\lambda^k \partial_\vphi^\beta \cK_n(\vphi, y)[\psi] \|_{\cH^s} 
& \lesssim_{s,\beta,k}
\gamma^{- |k|}C(s)^n \| \eta \|_{s_0 + \sigma}^n \Big( \| \eta \|_{s_0 + |\beta| + \sigma} \| \psi \|_{H^{s + |k|}_x} 
+   \| \eta \|_{s + |\beta| + \sigma} \| \psi \|_{H^{  |k| + 2}_x} \Big)\,.
\end{aligned}
\]
Since by \eqref{propagatore cal K cal Kn vphi} one has 
$\cK = \sum_{n \geq 0} (- 1)^n \cK_n$, the bound \eqref{bottle10} follows by the latter estimate taking 
$C(s)\| \eta \|_{s_0 + \sigma} \ll1$ small enough
 in such a way that $\sum_{n \geq 0} (C(s) \| \eta \|_{s_0 + \sigma})^n \leq 2$.
 
\smallskip
\noindent
The bound \eqref{bottle11} follows by \eqref{bottle10} and the estimate, 
for $m = 0,1,2$, 
\[
\| \partial_{x, y}^m \partial_\lambda^k \partial_\vphi^\beta   \cK(\vphi, y)  [\psi] \|_{\cH^s} 
\lesssim
\| \partial_\lambda^k \partial_\vphi^\beta   \cK(\vphi, y)  [\psi] \|_{\cH^{s + m}}.
\] 
Finally by \eqref{bottle15} we shall write
\[
\Delta_{12}\cR(\vphi, y ) := \partial_\lambda^{k_1}\partial_\vphi^{\beta_1} \Delta_{12}\cA_n(\vphi) 
\circ \partial_\lambda^{k_2} \partial_\vphi^{\beta_2} \widetilde{\mathtt{L}}(y, \vphi) 
+
\partial_\lambda^{k_1}\partial_\vphi^{\beta_1} \cA_n(\vphi) 
\circ \Delta_{12}\partial_\lambda^{k_2} \partial_\vphi^{\beta_2} \widetilde{\mathtt{L}}(y, \vphi) \,.
\]
Therefore, the bound \eqref{bottle12} follows 
reasoning as above recalling also
\eqref{bottle2}
\end{proof}

We need the following result, which is a corollary of 
Proposition \ref{stima cal Kn equazione di laplace}.
%section \ref{sec:stimegreen}.
\begin{cor}\label{stima of di green forzante qualunque}
Let  $D := \partial_y$ or $D := \langle D \rangle$
and define 
\begin{equation}\label{ascensore1}
 \mathcal{B}(\vphi, y) :=  \sum_{n \geq 0} (-1)^{n}D \circ \cA_n(\vphi, y)\,,
 \end{equation}
 where $\cA_n(\vphi, y)$ is in \eqref{AAnn}.
 For any $|k|\leq k_0$
there exists $\sigma \equiv \sigma(k_0) \gg 0$ large enough such that,
for any $s\geq 0$, $\beta \in \N^\nu$
there exists $\delta=\delta(s,\beta,k_0)>0$ small 
such that, 
if 
$\| \eta \|_{s_0 +|\beta|+ \sigma}^{k_0, \gamma} \leq \delta$, 
% if $\| \eta \|_{s_0 + |\beta| + \sigma}^{k_0, \gamma} \lesssim_{\beta,s}1$, 
 then
\begin{align}
\| \partial_\lambda^k \partial_\vphi^\beta \mathcal{B}(\vphi, y) [ f] \|_{\cO^{s}} 
&\lesssim_{s,\beta,k} \gamma^{- |k|}
\Big( \| f \|_{\cO^{s }} + \| \eta \|_{s + |\beta| + \sigma}^{k_0, \gamma} \| f \|_{\cO^{0}}\Big)\,,\label{corocoro}
\\
\| \partial_\vphi^\beta\Delta_{12}\mathcal{B}(\vphi, y)\|_{\mathcal{L}(\mathcal{O}^{p };\mathcal{O}^{p})}
&\lesssim_{p,\beta} \|\eta_{1}-\eta_{2}\|_{p+\widetilde{\s}+|\beta|}\,,%\;\;\;p\geq 0 \,.
\;\;\;0\leq p\leq s_0+\s-\widetilde{\s}\,,
\label{corocoroDelta12}
\end{align}
for some $\widetilde{\s}\leq \s$. %and for $p+\widetilde{\s}\leq s_0+\s$, $p\geq0$.
\end{cor}
\begin{proof}
To shorten notations in the proof we write $\| \eta \|_s$ instead of $\| \eta \|_s^{k_0, \gamma}$. 
By Lemma \ref{algebra striscia}, 
%Remark \ref{rmk:equiOHs}, 
it 
is enough to provide for any $n \geq 0$ an estimate of 
\begin{equation*}%\label{abruzzo10}
\begin{aligned}
\| D \circ \partial_\lambda^k \partial_\vphi^\beta \cA_n(\vphi, y) [f] \|_{\cO^s}&  \simeq 
\| D \circ \partial_\lambda^k \partial_\vphi^\beta \cA_n(\vphi, y) [f]  \|_{L^2_\vphi \cH^s} 
+ \| D \circ \partial_\lambda^k \partial_\vphi^\beta \cA_n(\vphi, y) [f]  \|_{H^s_\vphi L^2_{x, y}} 
\\&  
\lesssim 
\underbrace{\|  \partial_\lambda^k \partial_\vphi^\beta \cA_n(\vphi, y) [f]  \|_{L^2_\vphi \cH^{s + 1}} }_{I}
+\underbrace{ \|  \partial_\lambda^k \partial_\vphi^\beta \cA_n(\vphi, y) [f]  \|_{H^s_\vphi \cH^1}}_{II} \,. 
\end{aligned}
\end{equation*}
For $n= 0$, since $\cA_0 = L^{- 1}$ one immediately has (recall Lemma \ref{lemma laplace forzato dip da vphi}) that 
\begin{equation}\label{bla bla elliptic forcing 5}
\| D \circ \partial_\lambda^k L^{- 1} f \|_{\cO^s} \lesssim \| f \|_{\cO^s}, \quad \forall s \geq 0.
\end{equation}
We now prove the estimate for $\cA_n(\vphi, y)$ with $n \geq 1$,  using the bounds 
obtained in Lemma \ref{Corollario stima L inv F 0}. 

\medskip

\noindent
{\sc Estimate of $I$.}  By using the estimate of Lemma \ref{Corollario stima L inv F 0} %\eqref{abruzzo3} with Lemma \ref{lemma laplace forzato dip da vphi} 
one gets, for any $s \geq s_0 $, that
\begin{equation}\label{bla bla elliptic forcing 4}
\begin{aligned}
\| \partial_\lambda^k \partial_\vphi^\beta &\cA_n(\vphi, y) [f]  \|_{L^2_\vphi \cH^{s + 1}}  
\leq 
\| \partial_\lambda^k \partial_\vphi^\beta \cA_n(\vphi, y) [f]  \|_{L^2_\vphi \cH^{s + 2}} 
\\& 
\leq \gamma^{- |k|}C(s)^n\Big( \| \eta \|_{s_0 + \sigma}^{n - 1} \| \eta \|_{s_0 + |\beta| + \sigma}  \| f \|_{L^2_\vphi\cH^{s}} 
+ \| \eta \|_{s  + |\beta|  + \sigma} \| \eta \|_{s_0 + \sigma}^{n - 1} \| f \|_{L^2_\vphi L^2_{x, y}}\Big) 
\\& 
\stackrel{\eqref{prop L2 vhi Hs Hs vphi cal 2 xy}}{\leq} \gamma^{- |k|}
C(s)^n\Big( \| \eta \|_{s_0 + \sigma}^{n - 1} \| \eta \|_{s_0 + |\beta| + \sigma}  \| f \|_{\cO^s} + \| \eta \|_{s  + |\beta|  + \sigma} \| \eta \|_{s_0 + \sigma}^{n - 1} \| f \|_{\cO^{0}}\Big)\,.
\end{aligned}
\end{equation}

\smallskip

\noindent
{\sc Estimate of $II$.} We first note that 
\begin{equation}\label{bla bla elliptic forcing 0}
\begin{aligned}
\| \partial_\lambda^k \partial_\vphi^\beta \cA_n(\vphi, y) [f]  \|_{H^s_\vphi \cH^1} 
%& \lesssim_s 
%\sum_{|\alpha| \leq s} \Big\| \partial_\vphi^\alpha \Big( \partial_\lambda^k \partial_\vphi^\beta \cA_n(\vphi, y) [f] \Big) 
%\Big\|_{L^2_\vphi \cH^1} 
%\\& 
%\lesssim_s 
%\sum_{|\alpha| \leq s} \sum_{\alpha_1 + \alpha_2 = \alpha} 
%\Big\| \partial_\lambda^k \partial_\vphi^{\alpha_1 + \beta} \cA_n(\vphi, y) [\partial_\vphi^{\alpha_2} f] 
%\Big\|_{L^2_\vphi \cH^1} \\
& \lesssim_s 
\sum_{|\alpha| \leq s} \sum_{\alpha_1 + \alpha_2 = \alpha} 
\Big\| \partial_\lambda^k \partial_\vphi^{\alpha_1 + \beta} \cA_n(\vphi, y) [\partial_\vphi^{\alpha_2} f] \Big\|_{L^2_\vphi \cH^2}\,.
\end{aligned}
\end{equation}
Secondly, applying again the estimate of Lemma \ref{Corollario stima L inv F 0}, 
for any $|\alpha| \leq s$, $\alpha_1 + \alpha_2 = \alpha$ one gets that
\begin{equation}\label{bla bla elliptic forcing 1}
\begin{aligned}
\| \partial_\lambda^k \partial_\vphi^{\alpha_1 + \beta} \cA_n(\vphi, y) [\partial_\vphi^{\alpha_2} f] \|_{L^2_\vphi \cH^{ 2}} 
%& 
%\lesssim C(s)^n \gamma^{- |k|} \Big( \| \eta \|_{s_0 + \sigma}^{n - 1} \| \eta \|_{s_0 + |\alpha_1| + |\beta| + \sigma}  
%\| \partial_\vphi^{\alpha_2} f \|_{L^2_\vphi L^2_{x, y}}  
%\\& \qquad 
%+ \| \eta \|_{|\alpha_1| + |\beta|  + \sigma} \| \eta \|_{s_0 + \sigma}^{n - 1} \| \partial_\vphi^{\alpha_2} f \|_{L^2_\vphi L^2_{x, y}} \Big) \\
& \lesssim_{\beta,k} C(s)^n \gamma^{- |k|}   \| \eta \|_{s_0 + \sigma}^{n - 1}
 \| \eta \|_{s_0 + |\alpha_1| + |\beta| + \sigma}  
\|  f \|_{H^{|\alpha_2|}_\vphi L^2_{x, y}}\,. 
\end{aligned}
\end{equation}
By interpolation (see \eqref{2202.3}) and by \eqref{centoYoung} with $p\rightsquigarrow |\alpha|/|\alpha_1|$,
$q\rightsquigarrow |\alpha|/|\alpha_2|$, one gets
%$a b \lesssim a^{\frac{|\alpha|}{|\alpha_1|}} + b^{\frac{|\alpha|}{|\alpha_2|}}$, 
%one obtains that 
\begin{equation}\label{bla bla elliptic forcing 2}
\begin{aligned}
\| \eta \|_{s_0 + |\alpha_1| + |\beta| + \sigma} & \|  f \|_{H^{|\alpha_2|}_\vphi L^2_{x, y}} 
\\& \lesssim 
\Big( \| \eta \|_{s_0 + |\beta| + \sigma} \| f \|_{H^{|\alpha|}_\vphi L^2_{x, y}} \Big)^{\frac{|\alpha_2|}{|\alpha|}} 
\Big( \| \eta \|_{s_0 + |\alpha| + |\beta| + \sigma} \| f \|_{L^2_\vphi L^2_{x, y}} \Big)^{\frac{|\alpha_1|}{|\alpha|}} 
\\& 
\lesssim_{\alpha}
\| \eta \|_{s_0 + |\beta| + \sigma} \| f \|_{H^{|\alpha|}_\vphi L^2_{x, y}}  
+ \| \eta \|_{s_0 + |\alpha| + |\beta| + \sigma} \| f \|_{L^2_\vphi L^2_{x, y}} 
\\& 
\stackrel{\eqref{prop L2 vhi Hs Hs vphi cal 2 xy}}{\lesssim_{s}} 
\| \eta \|_{s_0 + |\beta| + \sigma'} \| f \|_{\cO^{s }}  + \| \eta \|_{s +  |\beta| + \sigma'} 
\| f \|_{\cO^{0}}\,,
\end{aligned}
\end{equation}
for some $\sigma'\geq \sigma \gg 0$ large enough. 
Hence, by \eqref{bla bla elliptic forcing 0}, \eqref{bla bla elliptic forcing 1}, \eqref{bla bla elliptic forcing 2}, one gets that for any $s \geq 0$
\begin{equation}\label{bla bla elliptic forcing 3}
\begin{aligned}
\| \partial_\lambda^k \partial_\vphi^\beta &\cA_n(\vphi, y) [f]  \|_{H^s_\vphi \cH^1} 
\\&\lesssim_{s,\beta,k} 
C(s)^n \gamma^{- |k|}   \| \eta \|_{s_0 + \sigma'}^{n - 1} \Big(  \| \eta \|_{s_0 + |\beta| + \sigma'} 
\| f \|_{\cO^{s}}  
+ \| \eta \|_{s +  |\beta| + \sigma'} \| f \|_{\cO^{0}} \Big)\,.
\end{aligned}
\end{equation}
The estimates \eqref{bla bla elliptic forcing 5}, \eqref{bla bla elliptic forcing 4}, \eqref{bla bla elliptic forcing 3} 
imply the claimed bound \eqref{corocoro}, using that $\mathcal{B} = \sum_{n \geq 0} D \circ \cA_n$ 
and using that $\sum_{n \geq 1} C(s,\beta,k)^n \| \eta \|^n_{s_0 + \sigma'} \leq 2$ 
(since $C(s,\beta,k)\| \eta \|_{s_0 + |\beta| + \sigma'} \ll 1$) where $C(s,\beta,k)>0$ is some fixed constant large.
The estimate \eqref{corocoroDelta12} can be obtained reasoning similarly
using the estimates on $\Delta_{12}\cA_n(\vphi, y)$ in Lemma \ref{Corollario stima L inv F 0}.
\end{proof}

\subsection{The non-homogeneous problem}\label{sec:Problemastriscianuova}
The \eqref{elliptic2} satisfied by the function $\phi$ is a linear ODE in $\pa_{y}$ with variable coefficients in $y$. 
It is technically convenient, for the constructions to follow, to have instead 
a solution of a linear constant coefficients ODE. 
For this reason we introduce in \eqref{func:tildevphi} below the function 
$\widetilde{\Phi}$ that will satisfy the system \eqref{elliptic5}.

Differently from section
\ref{sec:tameLaplace}, we 
consider the new vertical variable
\begin{equation}\label{change1bis}
z=\eta(x)+ \mathtt h y\,,\quad - 1  -\th^{-1}\eta(x)\leq y\leq 0\,,
\end{equation}
so that the closure of $\Omega_{t}$ in \eqref{Deta} becomes
\begin{equation}\label{Deta2bis}
\mathcal{S}_2=\{(x,y)\in \T^{d}_\Gamma \times \R\; :\; - 1  -\th^{-1}\eta(x) \leq y\leq 0\}\,.
\end{equation}
Consider the function
\begin{equation}\label{func:phi2}
\Phi(x,y):=\uphi(x,\eta(x)+ \mathtt h y)\stackrel{\eqref{func:phi}}{=}\phi\Big(x,\frac{\mathtt h y}{\mathtt h +\eta(x)} \Big)\,.
\end{equation}
We have the following.
\begin{lemma}\label{lemma laplace storto}
If $\uphi(x,z)$ solves the problem \eqref{elliptic1} then the function
$\Phi(x,y)$ in
\eqref{func:phi2} solves 
\begin{equation}\label{elliptic2bis}
\left\{\begin{aligned}
{\mathcal{L}}\Phi&=0\,,\qquad x\in\T^{d}_\Gamma\,,\;\; - 1  -\th^{-1}\eta(x)<y<0\,,\\
\Phi(x,0)&=\psi(x)\,, \quad x \in \T^d_\Gamma\\
(\pa_{y}\Phi)(x,- 1 -\th^{-1}\eta(x))&=0\,, \quad x \in \T^d_\Gamma
\end{aligned}\right.
\end{equation}
where 
\begin{equation}\label{op:L}
{\mathcal{L}}:={\mathcal{L}}(\eta):= \partial_{yy} + \mathtt h^{2} \Delta + \beta_1(\vphi,x) \cdot \nabla \partial_y 
+ \beta_2(\vphi,x) \partial_y + \beta_3(\vphi,x) \Delta\,,
\end{equation}
where 
\begin{equation}\label{op laplace raddrizzata buona}
\begin{aligned}
\beta_1  := - \mathtt h \dfrac{2 \nabla \eta}{1 + |\nabla \eta|^2}, \quad \beta_2 := -\mathtt h \dfrac{\Delta \eta}{1 + |\nabla \eta|^2}\,, \quad \beta_3 := - \mathtt h^{ 2} \dfrac{|\nabla \eta|^2}{1 + |\nabla \eta|^2}\,. 
\end{aligned}
\end{equation}
Moreover, there exists $\sigma \gg 0$ large enough such that for any
 $s\geq s_0$ there is $\delta=\delta(s,k_0)>0$ such that,  if 
\eqref{smalleta} holds, 
%$\| \eta \|_{s_0 + \sigma}^{k_0, \gamma} \ll1 $, 
then 
\begin{equation}\label{stime beta i coefficienti operatore giusto}
\begin{aligned}
\| \beta_i \|_s^{k_0, \gamma} &\lesssim_s \| \eta \|_{s + \sigma}^{k_0, \gamma}\,, 
%\quad \forall s \geq s_0\,,
\quad
\| \Delta_{12}\beta_i \|_p\lesssim_{p}\| \eta_1-\eta_2 \|_{p + \widetilde{\sigma}}\,,
\qquad i=1,2,3\,,
%\qquad p\geq s_0\,,
\end{aligned}
\end{equation}
for some $\widetilde{\s}\leq \s$ and for $p+\widetilde{\s}\leq s_0+\s$.
%for $i=1,2,3$.
\end{lemma}

\begin{proof}[{\bf Proof of Lemma \ref{lemma laplace storto}}]\label{proof of lemma laplace storto}
By \eqref{func:phi2} and \eqref{change1bis} we shall write
\[
\begin{aligned}
&\underline{\phi}(x,z)=\Phi\Big(x,\frac{z-\eta(x)}{\th}\Big)\,,
\\
&\Phi(x,0)=\uphi(x,\eta)=\psi\,,\qquad
(\pa_{y}\Phi)(x,-1-\th^{-1}\eta)=\th(\pa_{z}\uphi)(x,-1)=0\,.
\end{aligned}
\]
Therefore, the thesis follows by an explicit computation using \eqref{elliptic1}.
\end{proof}

\begin{rmk}
Notice that formula \eqref{eq:112a} (see also \eqref{eq:112abis})
becomes
\begin{equation}\label{eq:112aTRIS}
 G(\eta)\psi = 
 \th^{-1}
(1+|\nabla\eta|^2)(\pa_{y}\Phi)(x,0)-\nabla\eta\cdot\nabla\Phi(x,0)\,.
\end{equation}
\end{rmk}

An important point is that the operator $\mathcal{L}$ has coefficients depending only on the  variable
$x\in \mathbb{T}^{d}_{\Gamma}$ (and on $\vphi\in \mathbb{T}^{\nu}$ as parameters) but not on the 
variable $y\in[-1,0]$ differently from $\widetilde{\mathcal{L}}$ in \eqref{op:tildeL}. 
In addition to this, it can be also written as composition of two operators of the first order.
This property will be fundamental in section \ref{sec:decoupling}. 
%, we refer the reader to Appendix \ref{sec:divo1}.
For the moment we shall not exploit this fact.
The price to pay to have coefficients independent of $y\in[-1,0]$ is that the problem
\eqref{elliptic2bis} is posed on a strip with variable bottom.

\noindent
With the aim of still working on the \emph{straight} strip $[-1,0]$ we
consider an even  cut-off function $\chi\in C_c^{\infty}(\R;\R)$
\begin{equation}\label{cut-off}
\chi(y)=\left\{\begin{aligned}
&1\,,\quad |y|\leq 1/4\,,\\
&0\,,\quad |y|\geq1/2\,.
\end{aligned}\right.
\end{equation}
and introduce the function
\begin{equation}\label{func:tildevphi}
\widetilde{\Phi}(x,y):=\chi(y)\Phi(x,y)+(1-\chi(y))\vphi_0(x,y)
\end{equation}
where $\vphi_0$ and $\Phi$ are respectively the solutions of the problems \eqref{problema vphi 0 riscalato laplace}  
and  \eqref{elliptic2bis}. 
A straightforward computation\footnote{Indeed,
using \eqref{func:tildevphi}, we have $\vphi(x,0)=\Phi(x,0)=\psi$, 
$(\pa_{y}\vphi)(x,-1)=0$ and 
\begin{equation*}
\begin{aligned}
\mathcal{L}\widetilde{\Phi}&=\chi(y)\mathcal{L}(\Phi-\vphi_0)+[\mathcal{L},\chi(y)](\Phi-\vphi_0)
+\mathcal{L}\vphi_0
\\&
\stackrel{\eqref{elliptic2bis}}{=}[\mathcal{L},\chi(y)](\Phi-\vphi_0)+(1-\chi(y))\mathcal{L}\vphi_0\,.
\end{aligned}
\end{equation*}
Then \eqref{elliptic5}-\eqref{func:gg} follow recalling \eqref{problema vphi 0 riscalato laplace} and \eqref{op:L}.
} 
show that 
the function $\widetilde{\Phi}$ solves the elliptic (non-homogeneous) problem
\begin{equation}\label{elliptic5}
\left\{\begin{aligned}
{\mathcal{L}} \widetilde{\Phi}&=g\,,\qquad x\in\T^{d}_\Gamma\,,\;\; -1<y<0\,,\\
\widetilde{\Phi}(x,0)&=\psi(x)\,,\\
(\pa_{y}\widetilde{\Phi})(x,-1)&=0\,,
\end{aligned}\right.
\end{equation}
$\mathcal{L}$ is in \eqref{op:L} and the forcing term ${g}$ is defined as
\begin{equation}\label{func:gg}
\begin{aligned}
{g}:={g}(x, y)&:={g}(\vphi,x,y)
\\&: =
\Big[ \partial_{yy} + \beta_1(\vphi,x) \cdot \nabla \partial_y + \beta_2(\vphi,x) \partial_y, \chi(y) \Big] (\Phi - \vphi_0) 
\\
& \qquad  + (1 - \chi(y))\Big( \beta_1(\vphi,x) \cdot \nabla \partial_y + \beta_2(\vphi,x) \partial_y + \beta_3(\vphi,x) \Delta \Big) \vphi_0\,.
\end{aligned}
\end{equation}

\medskip
\noindent
{\bf Perturbative setting.} 
The solution $\widetilde{\Phi}$ of the problem \eqref{elliptic5}  can be written as 
\begin{equation}\label{vphi vphi0 u zero cond bordo}
\widetilde{\Phi} (x, y) = \vphi_0(x, y) + v(x, y)
\end{equation}
where $\vphi_0(x,y)$ solves \eqref{problema vphi 0 riscalato laplace}
while the ``new unknown''  
$v(x, y)$ solves the problem
\begin{equation}\label{elliptic4tris}
\left\{\begin{aligned}
{\mathcal{L}}v&=g+h\,,\qquad x\in\T^{d}_\Gamma\,,\;\; -1<y<0\,,\\
v(x,0)&=0\,, \\
(\pa_{y}v)(x,-1)&=0\,, 
\end{aligned}\right.
\end{equation}
with  $g$ as in \eqref{func:gg} and the new forcing term $h$ is  given by 
\begin{equation}\label{func:hh}
h:=h(x, y):=h(\vphi,x,y):= -F(\vphi)[\vphi_0]\,,
\end{equation}
where $F(\vphi)$ is the second order operator defined as (recall \eqref{op:L})
\begin{equation}\label{secondFF}
F(\vphi)[\cdot]:= \beta_1(\vphi,x) \cdot \nabla \partial_y + \beta_2(\vphi,x) \partial_y + \beta_3(\vphi,x) \Delta\,.
\end{equation}
\begin{rmk}\label{rmk:DNconV}
Notice that by \eqref{vphi vphi0 u zero cond bordo}, one has that 
\[
\partial_y \widetilde{\Phi}_{| y = 0} = (\partial_y \vphi_0)_{| y = 0} 
+ \partial_y v_{| y = 0} = |D| \tanh(\mathtt h |D|) \psi + \partial_y v_{| y = 0}\,,
\] 
hence the Dirichlet-Neumann operator $G(\eta)$ takes the form, recall \eqref{eq:112aTRIS},
\begin{equation}\label{eq:112aquatuor}
\begin{aligned}
 G(\eta)\psi 
& = \th^{-1}(1+|\nabla\eta|^2)(\pa_{y}\widetilde{\Phi})(x,0)-\nabla\eta\cdot\nabla \psi(x) 
 \\
& = \th^{-1}(1+|\nabla\eta|^2)) |D| \tanh(\mathtt h |D|) \psi +\th^{-1}(1+|\nabla\eta|^2)) \partial_y v_{| y = 0} -\nabla\eta\cdot\nabla \psi(x)\,.
\end{aligned}
\end{equation}
\end{rmk}
\noindent
The aim of this section is to prove the following.
\begin{prop}\label{lemma tame laplace partial y u langle D rangle u}
Let $v$ be the solution of \eqref{elliptic4tris}. Then the following holds. 
\begin{equation}\label{soluzioneV}
v = \cF_0(\vphi, y)[\psi] \,, \quad \partial_y v 
= \cF_1(\vphi, y)[\psi] \,, \quad \langle D \rangle v = \cF_2(\vphi, y)[\psi]
\end{equation}
where the linear operators $\cF_i(\vphi, y)$, $i=0,1,2$ satisfy the following properties. 
For any $|k|\leq k_0$ there exists $\s=\s(k_0)\gg0$ such that, 
for any $s\geq s_0$ (recall \eqref{s0 Dirichlet Neumann}), $\beta\in \N^{\nu}$, $i=0,1,2$,
there is $\delta=\delta(s,k_0,\beta)$ such that,
if $\| \eta \|_{s_0 + |\beta| + \sigma}^{k_0, \gamma} \leq \delta$,  
then 
one has
\begin{align}
\| \partial_\lambda^k \partial_\vphi^\beta \cF_i(\vphi, y)[\langle D \rangle^{- |\beta|  -k_0 - 2} \psi] \|_{\cO^s} 
&\lesssim_{s,\beta,k}
\gamma^{- |k|} \Big( \| \eta \|_{s_0 + |\beta| + \sigma}^{k_0, \gamma} \| \psi \|_{s} 
+ \| \eta \|_{s + |\beta| + \sigma}^{k_0, \gamma} \| \psi \|_{s_0} \Big)\,,\label{propagF}
\\
\| \partial_\vphi^\beta\Delta_{12}\cF_i(\vphi, y)\langle D\rangle^{ -|\beta|-3}\|_{\mathcal{L}(\mathcal{O}^{p};\mathcal{O}^{p})}
&\lesssim_{p, \beta} \|\eta_{1}-\eta_{2}\|_{p+\widetilde{\s}+|\beta|}\,,%\;\;\;p\geq s_0\,.
\label{propagFDelta12}
\end{align}
for some $\widetilde{\s}\leq \s$ and for $p+\widetilde{\s}\leq s_0+\s$, $p\geq s_0$.
\end{prop}

The proof of the proposition above is similar to the one of Corollary \ref{stima of di green forzante qualunque}.
A difference is in the coefficients of $\mathcal{L}$ in \eqref{op:L}, that in this case does not depend on $y\in[-1,0]$, differently from
$\widetilde{\mathcal{L}}$ in \eqref{op:tildeL}.
However, in this case the most important difficulty is to give 
a precise description of the forcing terms $g$ and $h$ in \eqref{func:gg}, \eqref{func:hh}.
This is the content of the following subsections.

\subsubsection{Estimates on the forcing term h}\label{sec:stimeforzanteHH Elle vero}
Recalling \eqref{secondFF},
\eqref{def vphi 0} we shall define the operator
\begin{equation}\label{def:opLLTT}
\mathtt{L}[\cdot]=\mathtt{L}(\vphi,y)[\cdot]=F(\vphi)\circ\mathcal{L}_0[\cdot]\,,
\end{equation}
so that, by \eqref{func:hh}, we have
\begin{equation}\label{termineH}
h=\mathtt{L}(\vphi,y)[\psi]\,.
\end{equation}
The following proposition can be proved as Proposition \ref{stime striscia mathtt L vphi y O s} (using the estimates of Lemma \ref{lemma laplace storto}).
\begin{prop}\label{stime mathttL su Os}
 Let  $k \in \N^{\nu + 1}$, $|k| \leq k_0$
 and consider the operator ${\mathtt{L}}={\mathtt{L}}(\vphi,y)$ in \eqref{def:opLLTT}. 
There exists $\sigma \equiv \sigma(k_0)  \gg 0$ such that, for any $\beta\in\N^{\nu}$, $s\geq 0$, 
there exists $\delta=\delta(s,\beta,k_0)>0$ small 
such that, 
if $\| \eta \|_{s_0 + \sigma}^{k_0, \gamma} \leq \delta$,
%$\| \eta \|_{s_0 + \sigma}^{k_0, \gamma} \ll  1$, 
then
one has the estimates
\begin{align}
\| \partial_\lambda^k \partial_\vphi^\beta \mathtt L(\vphi, y)[\langle D \rangle^{ - |k| - 2} \psi] \|_{\cO^s} 
&\lesssim_{s,\beta,k}
\gamma^{- |k|} \Big( \| \eta \|_{s_0 + |\beta| + \sigma}^{k_0, \gamma} \| \psi \|_s + \| \eta \|_{s + |\beta| + \sigma}^{k_0, \gamma} \| \psi \|_{0} \Big)\,,
\label{stima cal Os op mathtt L beta k}
\\
\| \partial_\vphi^\beta\Delta_{12}{\mathtt{L}}(\vphi, y)\langle D\rangle^{ -2}\|_{\mathcal{L}(\mathcal{O}^{p};\mathcal{O}^{p})}
&\lesssim_{p, \beta} \|\eta_{1}-\eta_{2}\|_{p+\widetilde{\s}+|\beta|}\,,%\;\;\;p\geq 0\,.
\;\;\;0\leq p\leq s_0+\s-\widetilde{\s}\,,
\label{stima cal Os op mathtt L beta k Delta12}
\end{align}
for some $\widetilde{\s}\leq \s$. %and for $p+\widetilde{\s}\leq s_0+\s$, $p\geq0$.
\end{prop}

\subsubsection{Estimates on the forcing term g}\label{sec:stimeforzanteGG Elle vero}
In this section we provide estimates on the forcing term $g$ in \eqref{func:gg}.
In order to describe precisely the forcing $g$ we proceed as follows.
Consider  the linear operator 
\begin{equation}\label{cambio di variabile striscia}
\cA (\vphi) : u(x, y) \mapsto u(x, \alpha(\vphi, x) y), \quad \alpha(\vphi, x) := \frac{\mathtt h}{\mathtt h + \eta(\vphi, x)}
\end{equation}
with $\eta\in C^\infty$ satisfying Hypothesis \ref{hypo:eta} and $\| \eta \|_{s_0 + \sigma}^{k_0, \gamma} \ll 1$ small enough, for some $\sigma \gg 0$. Note that by the Sobolev embedding $\| \eta \|_{L^\infty} \lesssim \| \eta \|_{s_0} \ll 1$ and since $\mathtt h_1 \leq \mathtt h \leq \mathtt h_2$, for $\| \eta \|_{s_0 } \ll \mathtt h_1$, the function $\alpha$ is always well defined. With this notation, the solution $\Phi$ of \eqref{elliptic2bis} 
can be written in terms of the solution $\phi$ in \eqref{solphiphi}
of the problem \eqref{elliptic2}. More precisely, one has
\begin{equation}\label{strutturaPhiPhi}
\Phi(x, y) = \cA(\vphi)[\phi]\stackrel{\eqref{solphiphi}}{=}\mathcal{A}(\vphi)[\vphi_0]+\mathcal{A}(\vphi)[u] 
\stackrel{\eqref{def vphi 0}, \eqref{propagKK}}{=}\cA(\vphi) \circ \cL_0[\psi] + \cA(\vphi) \circ \cK(\vphi, y)[\psi]\,.
\end{equation}
As a consequence, we shall write
\begin{equation}\label{strutturaGG}
\begin{aligned}
g&=M(\vphi,y)[\psi]
\\
M(\vphi,y)&:=
\Big[ \partial_{yy} + \beta_1(\vphi,x) \cdot \nabla \partial_y + \beta_2(\vphi,x) \partial_y, \chi(y) \Big] 
\circ \Big(\cA(\vphi) \circ \big(\cL_0+\cK(\vphi, y)\big) -\mathcal{L}_0\Big) 
\\
&
+ (1 - \chi(y))\Big( \beta_1(\vphi,x) \cdot \nabla \partial_y + \beta_2(\vphi,x) \partial_y + \beta_3(\vphi,x) \Delta \Big) \circ\mathcal{L}_0\,.
\end{aligned}
\end{equation}

\begin{rmk}[{\bf Lipschitz variation}]\label{rmk:variazioneAA}
 We also need to control the Lipschitz variation of $\mathcal{A}(\vphi)$ 
w.r.t. the embedding $\eta(\vphi,x)$. 
Given  some functions $\eta_{i}$ satisfying the smallness
$ \| \eta_{i} \|_{s_0  + \sigma}^{k_0, \gamma} \ll 1 $, $i=1,2$, for some $\s\gg0$, we define,
for $\tau\in[0,1]$, the maps
\[
\cB^{\tau} (\vphi) : u(x, y) \mapsto u(x, \widetilde{\alpha}(\tau;\vphi, x) y)\,,
\qquad
\cA_{i} (\vphi) : u(x, y) \mapsto u(x,{\alpha}_{i}(\vphi, x) y)\,,\;i=1,2\,,
\]
where
\begin{equation}\label{def:alphatilde}
\widetilde{\alpha}=\widetilde{\alpha}(\tau;x,\vphi):=\alpha_2+\tau(\alpha_1-\alpha_2)\,,
\qquad \alpha_{i}:=\frac{\mathtt{h}}{\mathtt{h}+\eta_{i}(\vphi,x)}\,,\;\;i=1,2\,.
\end{equation}
With this notation we shall write
we note that
\begin{equation}\label{eq:variazioneAA}
\Delta_{12}\mathcal{A}(\vphi))[\cdot]=
\int_{0}^{1}(\alpha_1-\alpha_2)\mathcal{B}^{\tau}\circ\pa_{y}[\cdot]d\tau\,.
\end{equation}
\end{rmk}

The main result of this section is the following.
\begin{prop}\label{prop:combo}
Let 
$k \in \N^{\nu + 1}$, $|k| \leq k_0$. There exists $\sigma \equiv \sigma(k_0) \gg 0$ large enough such that,
for any $s\geq s_0$, $\beta\in\N^{\nu}$,
there exists $\delta=\delta(s,\beta,k_0)>0$ small 
such that, 
if $\| \eta \|_{s_0 +|\beta|+\sigma}^{k_0, \gamma} \leq \delta$,
 %if $\eta$ satisfies \eqref{smalleta},
 %$\| \eta \|_{s_0 + |\beta| + \sigma}^{k_0, \gamma} \lesssim_{s,\beta,k_0} 1$, 
 then the linear operator $M(\vphi, y)$ in \eqref{strutturaGG}  satisfies 
{\small \begin{align}
\| \partial_\lambda^k \partial_\vphi^\beta M(\vphi, y)[\langle D \rangle^{- |\beta| - k_0 - 2} \psi] 
&\|_{\cO^s} \lesssim_{s, k,\beta}
\gamma^{- |k|}  \Big(  \| \eta \|_{s_0 + |\beta| + \sigma}^{k_0, \gamma} \| \psi \|_s 
+ \| \eta \|_{s + |\beta| + \sigma}^{k_0, \gamma} \| \psi \|_{s_0} \Big)\,,\label{albero1}
\\
\| \partial_\vphi^\beta\Delta_{12}{M}(\vphi, y)
&\langle D\rangle^{-|\beta|-3}\|_{\mathcal{L}(\mathcal{O}^{p};\mathcal{O}^{p})}
\lesssim_{p, \beta} \|\eta_{1}-\eta_{2}\|_{p+\widetilde{\s}+|\beta|}\,,%\;\;\;p\geq s_0\,.
\label{albero2}
\end{align}}
for some $\widetilde{\s}\leq \s$ and for $p+\widetilde{\s}\leq s_0+\s$, $p\geq s_0$.
Moreover $M(\vphi, y) \equiv 0$ near the boundary $y = 0$. 
\end{prop}

In order to prove the result above, we need to provide estimates on the operator $\mathcal{A}(\vphi)$ 
in \eqref{cambio di variabile striscia} (see Lemmata 
\ref{lemma cambio di variabile striscia strano}, \ref{stime cal A vphi - Id ca,bio variabile strano striscia}),
and on the compositions $\mathcal{A}(\vphi)\circ\mathcal{K}(\vphi,y)$ 
(see Lemmata \ref{lemma stima cal A partial x y m cal K vphi y}, \ref{corollario derivata cal A vphi cal K vphi})
and $\mathcal{A}(\vphi)\circ\mathcal{L}_0$ 
(see Lemma \ref{corollario derivata cal A vphi - Id cal L0}).

\subsubsection{Estimates on the change of variable on the strip}\label{sec:stimeA}
We consider a cut-off function
  \begin{equation}\label{cutoffZETA}
  \begin{aligned}
  &\zeta\in C^{\infty}_{c}((-1,0);\R)\,, \quad {\rm supp}(\zeta) \subset (- 1 + \epsilon, - \epsilon)\,,
  \qquad {\rm for\; some }\quad \epsilon>0\,,
  %\\&\red{\zeta\equiv1\;\;\; y\in (-1+ 2\epsilon,- 2\epsilon)}\qquad {\rm for\; some }\quad \epsilon>0\,,
 \end{aligned}
 \end{equation}
we define 
the operator
\begin{equation}\label{calAlcaA}
\widetilde{\mathcal{A}}(\vphi):=\zeta(y)\cA (\vphi)\,.
\end{equation}

\begin{lemma}{\bf (Estimates on $\widetilde{\mathcal{A}}$)}\label{lemma cambio di variabile striscia strano}
Let $k \in \N^{\nu + 1}$, $|k| \leq k_0$ and $0<\epsilon\ll1$. 
Then there exists $\sigma \equiv \sigma(k_0)   \gg 0$ 
large enough such that
the linear operator $\widetilde\cA(\vphi)$ in \eqref{calAlcaA} 
satisfies the following bounds.
%
%if $\| \eta \|_{s_0 + \sigma}^{k_0, \gamma} \ll 1$, 
%then the linear operator $\widetilde\cA(\vphi)$ in \eqref{calAlcaA} 
%satisfies the following bounds

\medskip
 
\noindent
$(i)$ There is $\delta=\delta(k_0,\epsilon)$ such that,
if $\| \eta \|_{s_0 + \sigma}^{k_0, \gamma}\leq \delta$, then
\begin{equation}\label{stime A tilde norma bassissima}
\begin{aligned}
& \| \widetilde\cA(\vphi) [u] \|_{L^2_{x, y}} \lesssim \| u \|_{L^2_{x, y}}, \quad u \in L^2([- 1, 0] \times \T^d_\Gamma)\,,  \\
& \| \widetilde\cA(\vphi) [u] \|_{{\mathcal H}^1} \lesssim \| u \|_{{\mathcal H}^1}, \quad u \in {\mathcal H}^1\,. 
\end{aligned}
\end{equation}

\noindent
$(ii)$ For any $s \geq 1$
there is $\delta=\delta(k_0,\epsilon,s)$ such that,
if $\| \eta \|_{s_0 + \sigma}^{k_0, \gamma}\leq \delta$,
%if $\| \eta \|_{s_0 + \sigma}^{k_0, \gamma} \ll 1$, 
then
\begin{align}
\| \widetilde\cA(\vphi) [u] \|_{\cH^s}  &\lesssim_s \| u \|_{\cH^s} + \| \eta \|_{s + \sigma} \| u \|_{{\mathcal H}^1}\,,
\label{lagrima1}
\\
 \| \partial_\lambda^k \widetilde\cA(\vphi)[u] \|_{L^2_{x, y}} 
&\lesssim_{k}
\gamma^{- |k|}  \| u \|_{\cH^{|k|}} \,,
\nonumber
\\
 \| \partial_\vphi^\beta \partial_\lambda^k \widetilde\cA(\vphi)[u] \|_{L^2_{x, y}} 
&\lesssim_{s,k}
\gamma^{- |k|} \Big(  \| u \|_{\cH^{s + |k|}} 
+ \| \eta \|_{s + \sigma}^{k_0, \gamma} \| u \|_{{\mathcal H}^{1 + |k|}} \Big)\,, \quad |\beta| = s\,,
\label{lagrima2}
\end{align}
for any $u\in \mathcal{H}^{s}$ and 
\begin{align}
\| \partial_\vphi^\beta\Delta_{12}\widetilde\cA(\vphi)\|_{\mathcal{L}(\mathcal{H}^{p+1};L^2_{x, y})}
&\lesssim_{p} \|\eta_{1}-\eta_{2}\|_{p+\widetilde{\s}}\,,%\;\;\;p\geq s_0\,,
\;\;|\beta|=p\,,
\label{lagrima2Delta12}
\end{align}
for some $\widetilde{\s}\leq \s$ and for $p+\widetilde{\s}\leq s_0+\s$, $p\geq s_0$.

\noindent
$(iii)$ For any $\beta \in \N^\nu$, any $s \geq s_0$, 
there is $\delta=\delta(k_0,\epsilon,s)$ such that,
if $\| \eta \|_{s_0 + |\beta|+\sigma}^{k_0, \gamma}\leq \delta$ for some $\s\gg1 $
 depending on $\beta$,
 %$\| \eta \|_{s_0 + |\beta| + \sigma}^{k_0, \gamma} \lesssim_{s,\beta,k_0} 1$, 
 then 
\begin{align}
\| \partial_\lambda^k \partial_\vphi^\beta \widetilde\cA(\vphi) [u] \|_{\cH^s} 
&\lesssim_{s, k, \beta}
\gamma^{- |k|} \Big(  \| u \|_{\cH^{s + |\beta| + |k|}} 
+ \| \eta \|_{s + |\beta| + \sigma}^{k_0, \gamma} \| u \|_{\cH^{s_0 + |\beta| + |k|}} \Big)\,, \label{lagrima3}
\\
\| \partial_\vphi^\beta\Delta_{12}\widetilde\cA(\vphi)\|_{\mathcal{L}(\mathcal{H}^{p+|\beta|+1};\mathcal{H}^{p})}
&\lesssim_{p, \beta} \|\eta_{1}-\eta_{2}\|_{p+\widetilde{\s}+|\beta|}\,,%\;\;\;p\geq s_0\,,
\label{lagrima4}
\end{align}
for any $u\in \mathcal{H}^{s+|\beta|+|k|}$, 
for some $\widetilde{\s}\leq \s$ and for $p+\widetilde{\s}\leq s_0+\s$, $p\geq s_0$.

\end{lemma}
\begin{proof}
First of all we recall that, 
by the standard Moser composition Lemma (see \cite{moser66}) one has 
\[
\| \alpha - 1 \|_s^{k_0, \gamma} \,,\, \| \alpha^{-  1}  - 1\|_s^{k_0, \gamma} \lesssim_s 
 \| \eta \|_{s + \sigma}^{k_0, \gamma}, \quad \forall s \geq s_0\,,
\]
and, by the Sobolev embedding,
\begin{equation}\label{stima funzione cambio di variabile striscia}
\begin{aligned}
%& \| \alpha - 1 \|_s^{k_0, \gamma} \,,\, \| \alpha^{-  1}  - 1\|_s^{k_0, \gamma} \lesssim_s 
% \| \eta \|_{s + \sigma}^{k_0, \gamma}, \quad \forall s \geq s_0 \\
% & \text{and by the Sobolev embedding} \quad 
 \| \alpha^{\pm 1} - 1 \|_{{\mathcal C}^s}^{k_0, \gamma}
 \lesssim \| \alpha^{\pm 1} - 1 \|_{s + s_0}^{k_0, \gamma}  \lesssim_s 
 \| \eta \|_{s + \sigma}^{k_0, \gamma}, \quad \forall s \geq 0\,. \\
% & \text{Moreover, if} \quad 
% \| \eta \|_{s_0 + \sigma}^{k_0, \gamma} \ll 1, \quad \text{then} \\
% & \| \alpha^{\pm 1} \|_{{\mathcal C}^1}^{k_0, \gamma} \lesssim \| \alpha^{\pm 1} \|_{s_0 + 1}^{k_0, \gamma} \lesssim 1\,.  
 \end{aligned}
\end{equation}
Moreover, recalling the smallness condition on $\eta$, one gets
\[
\| \alpha^{\pm 1} \|_{{\mathcal C}^1}^{k_0, \gamma} 
\lesssim \| \alpha^{\pm 1} \|_{s_0 + 1}^{k_0, \gamma} \lesssim 1\,.  
\]
%Similarly, by \eqref{smalleta} with $\delta$ small enough,
Secondly, using the smallness of $\| \eta \|_{s_0 + \sigma}^{k_0, \gamma}$ w.r.t. to $\epsilon$ 
one can conclude that 
if $y\in (-1+\epsilon,-\epsilon)$
then $\alpha(\vphi,x)y\in [-1,0]$.
%\red{Rifare con calma faa di bruno con norma equiv con la funzione cutoff}

\medskip

\noindent
{\sc Proof of $(i)$.} We prove the first estimate in \eqref{stime A tilde norma bassissima}. One has that 
\[
\begin{aligned}
\| \widetilde{\mathcal A}(\vphi)[u]\|_{L^2_{x, y}}^2  & = \int_{- 1}^0 \int_{\T^d_\Gamma} |\zeta(y)|^2 |u (x, \alpha(x) y)|^2 \, d x dy\,.
\end{aligned}
\]
By making the change of variable $z = \alpha(x) y$, one gets that 
\[
\begin{aligned}
\| \widetilde{\mathcal A}(\vphi)[u]\|_{L^2_{x, y}}^2  
& = \int_{- 1}^0 \int_{\T^d_\Gamma} |\zeta(z/\alpha(x))|^2 |u (x, z)|^2 |\alpha(x)^{- 1}| \, d x dz\,. 
\\
& \stackrel{\eqref{cutoffZETA}, \eqref{stima funzione cambio di variabile striscia}}{\lesssim} 
\int_{- 1}^0 \int_{\T^d_\Gamma}  |u (x, z)|^2  \, d x dz \lesssim \| u \|_{L^2_{x, y}}^2\,. 
\end{aligned}
\]
Now let us prove the second estimate in \eqref{stime A tilde norma bassissima}. 
Arguing as in the previous estimate, one obtains that 
\[
\begin{aligned}
\| \nabla_{y, x} &(\widetilde{\mathcal A}(\vphi)[u]) \|_{L^{2}_{x, y}} 
\\& \lesssim 
\| \widetilde{\mathcal A}(\vphi)[\nabla_x u]  \|_{L^2_{x, y}} +
\| \widetilde{\mathcal A}(\vphi)[\pa_y u] \nabla_x \alpha \|_{L^2_{x, y}} 
%\\&
+ \| \zeta'(\cdot) {\mathcal A}(\vphi)[u] \|_{L^2_{x, y}} 
+ \| \widetilde{\mathcal A}(\vphi)[\partial_y u]  \alpha \|_{L^2_{x, y}}
\\&
\lesssim\| \nabla_x u \|_{L^2_{x, y}}+
\| u \|_{L^2_{x, y}}
+\| \pa_y u \|_{L^2_{x, y}} \| \nabla_x \alpha \|_{L^\infty_{x, y}} 
+\| \partial_y u \|_{L^2_{x, y}} \| \alpha \|_{L^\infty_{x, y}} 
\stackrel{\eqref{stima funzione cambio di variabile striscia}}{\lesssim} \| u \|_{{\mathcal H}^1} \,.
\end{aligned}
\]
The claimed statement then follows. 
%and arguing as in the previous estimate, one obtains that 
%\[
%\begin{aligned}
%\| \widetilde{\mathcal A}(\vphi)[\pa_y u] \nabla_x \alpha \|_{L^2_{x, y}} 
%& \lesssim 
%\| \pa_y u \|_{L^2_{x, y}} \| \nabla_x \alpha \|_{L^\infty_{x, y}} 
%\stackrel{\eqref{stima funzione cambio di variabile striscia}}{\lesssim} \| u \|_{{\mathcal H}^1} 
%\\
%\| \widetilde{\mathcal A}(\vphi)[\nabla_x u] \|_{L^2_{x, y}} 
%+\| \zeta'(\cdot) {\mathcal A}(\vphi)[u] \|_{L^2_{x, y}} 
%& \lesssim \| \nabla_x u \|_{L^2_{x, y}}+
%\| u \|_{L^2_{x, y}} \lesssim \| u \|_{{\mathcal H}^1} \\
% \| \widetilde{\mathcal A}(\vphi)[\partial_y u]  \alpha \|_{L^2_{x, y}} 
% & \lesssim \| \partial_y u \|_{L^2_{x, y}} \| \alpha \|_{L^\infty_{x, y}} 
% \stackrel{\eqref{stima funzione cambio di variabile striscia}}{\lesssim} \| u \|_{{\mathcal H}^1}\,. 
%\end{aligned}
%\]
%The claimed statement then follows. 
%

\medskip

\noindent
{\sc Proof of $(ii)$.} 

\medskip
\noindent
{\bf Proof of \eqref{lagrima1}}
The case $s=0$ follows trivially by \eqref{stime A tilde norma bassissima} 
recalling the definition \eqref{def sobolev striscia}.
For any $s \geq 1$, in view of \eqref{prop L2 vhi Hs Hs vphi cal 2 xy}, we estimate 
\begin{equation}\label{trenoRoma1}
\| \widetilde{\mathcal A}(\vphi)[u] \|_{{\mathcal H}^s} 
\simeq_s
\underbrace{\| \widetilde{\mathcal A}(\vphi)[u] \|_{L^2_y H^s_x}}_{I} + 
\underbrace{\| \widetilde{\mathcal A}(\vphi)[u] \|_{H^s_y L^2_x}}_{II} 
\end{equation}
and we estimate separately the two terms above. 

\smallskip
\noindent
{\it Estimate of $II$ % $\| \widetilde{\mathcal A}(\vphi)[u] \|_{H^s_y L^2_x}$
in \eqref{trenoRoma1}.} For any $1 \leq m \leq s$, one has that (recall \eqref{cambio di variabile striscia}, \eqref{calAlcaA})
\[
\begin{aligned}
\partial_y^m \big( \widetilde{\mathcal A}(\vphi)[u] \big) 
%& = \sum_{m_1 + m_2 = m} C(m_1, m_2) \partial_y^{m_1} \zeta \partial_y^{m_2}\big( {\mathcal A}(\vphi)[u] \big) 
%\\
%& = \sum_{m_1 + m_2 = m} C(m_1, m_2) \partial_y^{m_1} \zeta \partial_y^{m_2}\big( {\mathcal A}(\vphi)[u] \big) 
%\\
&  = \sum_{m_1 + m_2 = m} C(m_1, m_2) \partial_y^{m_1} \zeta   {\mathcal A}(\vphi)[\partial_y^{m_2} u] \alpha(\vphi, x)^{m_2}\,,
\end{aligned}
\]
for some constants $C(m_1, m_2) > 0$. Hence 
\[
\begin{aligned}
\| \partial_y^m \big( \widetilde{\mathcal A}(\vphi)[u] \big) \|_{L^2_{x, y} } 
& 
\lesssim_m  \sum_{m_1 + m_2 = m} \| \partial_y^{m_1} \zeta   {\mathcal A}(\vphi)[\partial_y^{m_2} u]  \|_{L^2_{x, y}} \| \alpha \|_{L^\infty}^{m_2}   
%\\
%& \stackrel{ \eqref{stima funzione cambio di variabile striscia}}{\lesssim_m} \sum_{m_1 + m_2 = m} \| \partial_y^{m_1} \zeta   {\mathcal A}(\vphi)[\partial_y^{m_2} u]  \|_{L^2_{x, y}}  
\\
& \stackrel{\eqref{stime A tilde norma bassissima}, \eqref{stima funzione cambio di variabile striscia}}{\lesssim_m} 
\sum_{0 \leq m_2 \leq m} \| \partial_y^{m_2} u \|_{L^2_{x, y}} \lesssim_m \| u \|_{{\mathcal H}^m}\,.
\end{aligned}
\]
The latter estimate implies that for $s \geq 1$, 
\begin{equation}\label{stima Hs y L2 x A tilde}
\| \widetilde{\mathcal A}(\vphi)[u] \|_{H^s_y L^2_x} \lesssim_s \| u \|_{{\mathcal H}^s}\,. 
\end{equation}

\medskip

\noindent
{\it Estimate of $I$
%$\| \widetilde{\mathcal A}(\vphi)[u] \|_{L^2_y H^s_x}$, $s \geq 1$
in \eqref{trenoRoma1}.} 
We shall prove that for any $s \geq 1$, 
$u \in {\mathcal H}^s$, and 
if $\| \eta \|_{s_0 + \sigma} \ll 1$, $\sigma \gg 1$, 
then 
\begin{equation}\label{stima A tilde L2y Hs x}
\|  \widetilde{\mathcal A}(\vphi)[u] \|_{L^2_y H^s_x} \lesssim_s \| u \|_{{\mathcal H}^s} + \| \eta \|_{s + \sigma} \| u \|_{{\mathcal H}^1}\,. 
\end{equation}
We argue by induction. For $s = 1$ the estimate follows by the item $(i)$. 
Now assume that \eqref{stima A tilde L2y Hs x} holds for some 
$s \geq 1$ and let us prove it for $s + 1$. 
One has that if $u \in {\mathcal H}^{s + 1}$, then 
\begin{equation}\label{stima A tilde L2y Hs xA}
\begin{aligned}
\|  \widetilde{\mathcal A}(\vphi)[u] &\|_{L^2_y H^{s + 1}_x}  
\simeq 
\|  \widetilde{\mathcal A}(\vphi)[u] \|_{L^2_y L^2_x} 
+ \|  \nabla_x (\widetilde{\mathcal A}(\vphi)[u]) \|_{L^2_y H^{s}_x}  
%\\& 
%\lesssim 
%\| u \|_{L^2_{x, y}} + \| \widetilde{\mathcal A}(\vphi)[\pa_yu] \nabla_x \alpha \|_{L^2_y H^{s}_x} 
%+ \| \widetilde{\mathcal A}(\vphi)[\nabla_xu]  \|_{L^2_y H^{s}_x} 
\\
& \lesssim \| u \|_{L^2_{x, y}}  
+ \| \widetilde{\mathcal A}(\vphi)[\nabla_xu]  \|_{L^2_y H^{s}_x}
\\&\qquad+ \| \widetilde{\mathcal A}(\vphi)[\pa_y u] \|_{L^2_y H^s_x} \| \nabla_x \alpha \|_{L^\infty} 
  +  \| \widetilde{\mathcal A}(\vphi)[\pa_y u] \|_{L^2_y L^2_x} \| \nabla_x \alpha \|_{{\mathcal C}^s} \\
& \stackrel{\eqref{stima funzione cambio di variabile striscia}}{\lesssim_{s}}  \| u \|_{{\mathcal H}^s} + \| \alpha \|_{{\mathcal C}^{s + 1}} \| u \|_{{\mathcal H}^1} 
+ \| \widetilde{\mathcal A}(\vphi)[\pa_yu] \|_{L^2_y H^s_x}+ \| \widetilde{\mathcal A}(\vphi)[\nabla_xu]  \|_{L^2_y H^{s}_x}  \,.
\end{aligned}
\end{equation}
By the induction assumption, for $m\in \N^{\d+1}$ with $|m|\leq 1$, one has that 
\begin{equation}\label{stima A tilde L2y Hs xB}
\begin{aligned}
\| \widetilde{\mathcal A}(\vphi)[\pa_{y,x}^{m} u] \|_{L^2_y H^s_x}  & \lesssim_s 
\| \pa_{y,x}^m u \|_{{\mathcal H}^s} + \| \eta \|_{s + \sigma} \| \pa_{y,x}^{m} u\|_{{\mathcal H}^1} \\
& \lesssim_s  \|  u \|_{{\mathcal H}^{s + 1}} + \| \eta \|_{s + \sigma} \|  u\|_{{\mathcal H}^2}\,. 
\end{aligned}
\end{equation}
By interpolation estimates \eqref{2202.2} and \eqref{strippointrippo},  and by  Young \eqref{centoYoung},
%with $p\rightsquigarrow s$, $q\rightsquigarrow s/(s-1)$, 
one gets
%inequality $a b \lesssim a^{s} + b^{\frac{s}{s - 1}}$, one obtains that 
\begin{equation}\label{stima A tilde L2y Hs xC}
\begin{aligned}
\| \eta \|_{s + \sigma} \|  u\|_{{\mathcal H}^2} 
%& \lesssim 
%\| \eta \|_{\sigma + 1}^{\frac{1}{s}} \| \eta \|_{s + 1 + \sigma}^{\frac{s - 1}{s}}  
%\| u \|_{{\mathcal H}^{ 1}}^{\frac{s - 1}{s}} \| u \|_{{\mathcal H}^{s + 1}}^{\frac{1}{s}} 
%\\
& 
\lesssim 
\Big( \| \eta \|_{\sigma + 1} \| u \|_{{\mathcal H}^{s + 1}} \Big)^{\frac{1}{s}} 
\Big(\| \eta \|_{s + 1 + \sigma}  \| u \|_{{\mathcal H}^{ 1}}  \Big)^{\frac{s - 1}{s}} \\
& \lesssim_{s} \| \eta \|_{\sigma + 1} \| u \|_{{\mathcal H}^{s + 1}} 
+ \| \eta \|_{s + 1 + \sigma}  \| u \|_{{\mathcal H}^{ 1}} \,.
%\\& 
%\stackrel{\| \eta \|_{\sigma + 1} \ll 1}{\lesssim} 
%\| u \|_{{\mathcal H}^{s + 1}} + \| \eta \|_{s + 1 + \sigma}  \| u \|_{{\mathcal H}^{ 1}}\,. 
\end{aligned}
\end{equation}
By combining 
\eqref{stima A tilde L2y Hs xA}, \eqref{stima A tilde L2y Hs xB}, \eqref{stima A tilde L2y Hs xC} 
and using the smallness condition on $\eta$, one deduces 
the claimed bound \eqref{stima A tilde L2y Hs x} for $s + 1$. 
In conclusion estimate \eqref{lagrima2} follows by \eqref{stima Hs y L2 x A tilde}, \eqref{stima A tilde L2y Hs x}. 

%\noindent
%Recalling \eqref{def sobolev striscia}, in order to get the estimate \eqref{lagrima1},
%we have to estimate, for $0\leq k\leq s$ the term 
%\[
%\begin{aligned}
%\int_{-1+\delta}^{-\delta}\|\pa_{y}^{k}\zeta(y)\mathcal{A}(\vphi)u\|_{H_{x}^{s-k}}dy&\lesssim_{k}
%%\pa_{y}^{k}\zeta(y)u(x,\alpha(\vphi,x)y)
%\sum_{k_1+k_2=k}\int_{-1+\delta}^{-\delta}\|\pa_{y}^{k_1}\zeta(y)\pa_{y}^{k_2}\mathcal{A}(\vphi)u\|_{H_{x}^{s-k}}dy
%\\&\lesssim_{k}C\sum_{0\leq k\leq s}\int_{-1+\delta}^{-\delta}\|\pa_{y}^{k}\mathcal{A}(\vphi)u\|_{H_{x}^{s-k}}dy
%\end{aligned}
%\]
%where $C>0$ is some constant depending on $\|\zeta\|_{W_{y}^{k,\infty}}$. 
%Note that (recall \eqref{cambio di variabile striscia})
%\[
%\pa_{y}^{k}\mathcal{A}(\vphi)u=(\partial_y^k u)(x,\alpha(\vphi,x)y)(\alpha(\vphi,x))^k\,.
%\]
%The bound \eqref{lagrima1} then follows by 
%classical composition estimate on Sobolev spaces (see \cite{moser66}) {\color{red} CITAZIONE SBAGLIATA! CORREGGERE!}.
%
%***************
% is again a classical composition estimates which follows 
%by using Faa di Bruno Formula. For more details of similar estimates we refer (for instance) to 
%\cite{moser66} and \cite{baldi2013}.

\medskip

\noindent
{\bf Proof of \eqref{lagrima2}.}
We argue by induction on $s \geq 0$. 
For $s=|\beta|=0$ we have, by 
Faa di Bruno formula, 
%\[
%\begin{aligned}
%\partial_\lambda^k  \widetilde\cA(\vphi)[u]
%%&= 
%%\zeta(y)\partial_\lambda^k  \cA(\vphi)[u]
%% \\&
%% \simeq_{s, \beta, k} 
% =
% \sum_{1 \leq m \leq  |k|} \sum_{\begin{subarray}{c}
%k_1 + \ldots + k_m = k
%\end{subarray}}  
%\widetilde\cA(\vphi)[\partial_y^m u] y^m 
%\partial_\lambda^{k_1} \alpha \ldots \partial_\lambda^{k_m} \alpha\,.
%\end{aligned}
%\]
%Therefore
one gets
\[
\begin{aligned}
\| \partial_\lambda^k
\widetilde{\mathcal A}(\vphi)[u] \|_{L^2_{x, y}} 
&\lesssim_{k}
 \sum_{1 \leq m \leq  |k|} \sum_{\begin{subarray}{c}
k_1 + \ldots + k_m = k
\end{subarray}}  
\|\widetilde\cA(\vphi)[\partial_y^m u] y^m 
\partial_\lambda^{k_1} \alpha \ldots \partial_\lambda^{k_m} \alpha\|_{L^2_{x, y}} 
 \\
 &\lesssim_{k}
 \| \widetilde\cA(\vphi)[\partial_y^m u]  \|_{L^2_{x, y}} 
\| \partial_\lambda^{k_1} \alpha \ldots \partial_\lambda^{k_m} \alpha \|_{L^{\infty}_{x,y}} 
\stackrel{\eqref{stime A tilde norma bassissima},\eqref{stima funzione cambio di variabile striscia}}{\lesssim_{k}}
\|u\|_{\mathcal{H}^{|k|}}
\end{aligned}
\]
where we also used the smallness condition on the low norms of $\eta$. The latter bound implies 
\eqref{lagrima2} for $s=0$.
%For $s = 1$, it follows by a direct computation. 
We now prove the induction step, hence let us assume that \eqref{lagrima2} holds for any $\beta \in \N^\nu$, $1 \leq |\beta| \leq s$. Let $\beta \in \N^\nu$, $|\beta| = s$, then for any $i = 1, \ldots, \nu$ one has that  
\[
\begin{aligned}
\| \partial_\lambda^k \partial_\vphi^\beta \partial_{\vphi_i}  
\widetilde{\mathcal A}(\vphi)[u] \|_{L^2_{x, y}} 
& \simeq \| \partial_\lambda^k \partial_\vphi^\beta   
(\widetilde{\mathcal A}(\vphi)[\partial_y u] \partial_{\vphi_i} \alpha  ) \|_{L^2_{x, y}} 
%\\
%& \lesssim_s \sum_{\begin{subarray}{c}
%\beta_1 + \beta_2 = \beta \\
%k_1 + k_2 = k
%\end{subarray}} \| (\partial_\lambda^{k_1} \partial_\vphi^{\beta_1}   \widetilde{\mathcal A}(\vphi)[\partial_y u]) (\partial_\lambda^{k_2} \partial_\vphi^{\beta_2}\partial_{\vphi_i} \alpha  ) \|_{L^2_{x, y}} 
\\& 
\lesssim_{s,k} \sum_{\begin{subarray}{c}
\beta_1 + \beta_2 = \beta \\
k_1 + k_2 = k
\end{subarray}} \| \partial_\lambda^{k_1} \partial_\vphi^{\beta_1}   \widetilde{\mathcal A}(\vphi)[\partial_y u] \|_{L^2_{x, y}}  \|\partial_\lambda^{k_2} \partial_\vphi^{\beta_2}\partial_{\vphi_i} \alpha   \|_{L^\infty_{x, y}} \\
& \lesssim_{s,k}  \gamma^{- |k|}\sum_{\begin{subarray}{c}
\beta_1 + \beta_2 = \beta \\
k_1 + k_2 = k
\end{subarray}} \big(\| u \|_{{\mathcal H}^{|\beta_1| + |k_1| + 1}} + \| \eta \|_{|\beta_1| + \sigma} \| u \|_{{\mathcal H}^{|k_1| + 1}} \big)   \|  \alpha   \|_{{\mathcal C}^{|\beta_2| + 1}}\,.
\end{aligned}
\]
By using that 
$\ \|  \alpha   \|_{{\mathcal C}^{|\beta_2| + 1}} \lesssim 
1 + \| \eta \|_{|\beta_2| + 1 + \sigma}$, $\sigma \gg 0$ 
(see \eqref{stima funzione cambio di variabile striscia}), one obtains 
\begin{equation}\label{cavolfiore 0}
\begin{aligned}
\| \partial_\lambda^k &\partial_\vphi^\beta \partial_{\vphi_i} 
 \widetilde{\mathcal A}(\vphi)[u] \|_{L^2_{x, y}}  
 \lesssim_{s, k} \gamma^{- |k|}\sum_{\begin{subarray}{c}
\beta_1 + \beta_2 = \beta \\
k_1 + k_2 = k
\end{subarray}} \| u \|_{{\mathcal H}^{|\beta_1| + |k_1| + 1}} + \| \eta \|_{|\beta_1| + \sigma} \| u \|_{{\mathcal H}^{|k_1| + 1}}  
 \\& 
 + \gamma^{- |k|} \sum_{\begin{subarray}{c}
\beta_1 + \beta_2 = \beta \\
k_1 + k_2 = k
\end{subarray}} \| u \|_{{\mathcal H}^{|\beta_1| + |k_1| + 1}} \| \eta \|_{|\beta_2| + 1 + \sigma} + \| \eta \|_{|\beta_1| + \sigma} \| \eta \|_{|\beta_2| + 1 + \sigma} \| u \|_{{\mathcal H}^{|k_1| + 1}}\,.
\end{aligned}
\end{equation}
By interpolation (see \eqref{2202.2} and \eqref{strippointrippo}), one obtains that 
\[
\begin{aligned}
& \| u \|_{{\mathcal H}^{|\beta_1| + |k_1| + 1}}  \lesssim 
\| u \|_{{\mathcal H}^{|\beta| + |k_1| + 1}}^{\frac{|\beta_1|}{|\beta|}} 
\| u \|_{{\mathcal H}^{ |k_1| + 1}}^{\frac{|\beta_2|}{|\beta|}} 
\lesssim 
\| u \|_{{\mathcal H}^{|\beta| + |k| + 1}}^{\frac{|\beta_1|}{|\beta|}} \| u \|_{{\mathcal H}^{ |k| + 1}}^{\frac{|\beta_2|}{|\beta|}} \,, \\
&   \| \eta \|_{|\beta_1| + \sigma} \lesssim  \| \eta \|_{ \sigma}^{\frac{|\beta_2|}{|\beta|}} \| \eta \|_{|\beta| + \sigma}^{\frac{|\beta_1|}{|\beta|}}\,, \qquad   \| \eta \|_{|\beta_2| + 1 + \sigma} \lesssim \| \eta \|_{ \sigma + 1}^{\frac{|\beta_1|}{|\beta|}}  \| \eta \|_{|\beta| + 1 + \sigma}^{\frac{|\beta_2|}{|\beta|}}\,.
\end{aligned}
\]
Therefore by using the Young inequality \eqref{centoYoung}, recalling $|\beta|=s$ and  
%$ab \lesssim a^{\frac{|\beta|}{|\beta_1|}} + b^{\frac{|\beta|}{|\beta_2|}}$ and 
$\| \eta \|_{\sigma + 1} \ll 1$, one gets
\[
\begin{aligned}
 \| u \|_{{\mathcal H}^{|\beta_1| + |k_1| + 1}} \| \eta \|_{|\beta_2| + 1 + \sigma} & \lesssim \
 \Big( \| u \|_{{\mathcal H}^{|\beta| + |k| + 1}} \| \eta \|_{ \sigma + 1} \Big)^{\frac{|\beta_1|}{|\beta|}} \Big( \| \eta \|_{|\beta| + 1 + \sigma} \| u \|_{{\mathcal H}^{ |k| + 1}} \Big)^{\frac{|\beta_2|}{|\beta|}}   \\
%& \lesssim  \| u \|_{{\mathcal H}^{|\beta| + |k| + 1}}
%+ \| \eta \|_{|\beta| + 1 + \sigma} \| u \|_{{\mathcal H}^{ |k| + 1}}  
%\\
& \lesssim_{s} \| u \|_{{\mathcal H}^{s + 1+ |k| }}+ \| \eta \|_{s  + 1 + \sigma} \| u \|_{{\mathcal H}^{ |k| + 1}} \,,
%\quad \text{and} 
\\
\| \eta \|_{|\beta_1| + \sigma} \| \eta \|_{|\beta_2| + 1 + \sigma} \| u \|_{{\mathcal H}^{|k_1| + 1}} & \lesssim \| \eta \|_{ \sigma}^{\frac{|\beta_2|}{|\beta|}} \| \eta \|_{|\beta| + \sigma}^{\frac{|\beta_1|}{|\beta|}} \| \eta \|_{ \sigma + 1}^{\frac{|\beta_1|}{|\beta|}}  \| \eta \|_{|\beta| + 1 + \sigma}^{\frac{|\beta_2|}{|\beta|}}  \| u \|_{{\mathcal H}^{|k_1| + 1}} \\
& \lesssim_{s}  \| \eta \|_{s + 1 + \sigma}  \| u \|_{{\mathcal H}^{|k| + 1}}\,.
\end{aligned}
\]
The latter chain of inequalities, together with \eqref{cavolfiore 0} imply the claimed bound \eqref{lagrima2} for $s + 1$. 

\noindent
The estimate \eqref{lagrima2Delta12} follows by similar arguments. %reasoning as in the proof of \eqref{lagrima2Delta12}.

%\noindent
%{\bf Proof of \eqref{lagrima2Delta12}.} First of all we note that
%the function $\widetilde{\alpha}$ in \eqref{def:alphatilde} satisfies the bounds
%\[
%\| \widetilde{\alpha} - 1 \|_s^{k_0, \gamma} \,,\, \| \widetilde{\alpha}^{-  1}  - 1\|_s^{k_0, \gamma} \lesssim_s 
%\sup_{i=1,2} \| \eta_{i} \|_{s + \sigma}^{k_0, \gamma}, \quad \forall s \geq s_0\,,
%\] 
%uniformly in $\tau\in[0,1]$ when $\alpha_{i}(\vphi,x):=\mathtt{h}(\mathtt h + \eta_{i}(\vphi, x))^{-1}$, $i=1,2$, 
%for some functions $\eta_{i}$ satisfying 
%$ \| \eta_{i} \|_{s_0 + |\beta| + \sigma}^{k_0, \gamma} \lesssim 1$, $i=1,2$.
%Therefore estimate \eqref{lagrima2Delta12} follows by exploiting formula \eqref{eq:variazioneAA},
%using tame estimates on the product,
%and the fact that $\mathcal{A}_1, \mathcal{B}^{\tau}$ satisfies estimates like \eqref{lagrima2} uniformly in $\tau\in[0,1]$.

\smallskip

\noindent
{\sc Proof of $(iii)$.} Let $\beta \in \N^\nu$, $k \in \N^{\nu + 1}$, $|k| \leq k_0$. 
By the Faa di Bruno Formula one gets that 
\[
\begin{aligned}
\|\partial_\lambda^k \partial_\vphi^\beta \widetilde\cA(\vphi)[u]\|_{\mathcal{H}^{s}}&
\lesssim 
\|\zeta(y)\partial_\lambda^k \partial_\vphi^\beta \cA(\vphi)[u]\|_{\mathcal{H}^{s}}
 \\&\lesssim_{s, \beta, k} \sum_{1 \leq m \leq |\beta| + |k|} \sum_{\begin{subarray}{c}
\beta_1 + \ldots + \beta_m = |\beta| \\
k_1 + \ldots + k_m = k
\end{subarray}}  
\|\widetilde\cA(\vphi)[\partial_y^m u] y^m 
\partial_\lambda^{k_1}\partial_\vphi^{\beta_1} \alpha 
\ldots \partial_\lambda^{k_m}\partial_\vphi^{\beta_m} \alpha\|_{\mathcal{H}^s}\,.
\end{aligned}
\]
Then for any $\beta_1, \ldots, \beta_m$, $k_1, \ldots, k_m$ with $\beta_1 + \ldots + \beta_m = \beta$, $k_1 + \ldots + k_m = k$, 
and recalling Lemma \ref{algebra striscia},
one gets that for any $s \geq s_0$, 
\begin{equation}\label{stima cambio variabile x y cal A brutto 0}
\begin{aligned}
 \| \widetilde\cA(\vphi)[\partial_y^m u] y^m& \partial_\lambda^{k_1}\partial_\vphi^{\beta_1} \alpha \ldots 
\partial_\lambda^{k_m}\partial_\vphi^{\beta_m} \alpha \|_{\cH^s} 
\qquad\qquad\\& 
\lesssim_s \| \widetilde\cA(\vphi)[\partial_y^m u]  \|_{\cH^s} 
\| \partial_\lambda^{k_1} \partial_\vphi^{\beta_1} \alpha \ldots \partial_\lambda^{k_m} \partial_\vphi^{\beta_m} \alpha \|_{s_0 + \sigma} 
\\&\qquad \qquad+ \| \widetilde\cA(\vphi)[\partial_y^m u]  \|_{\cH^{s_0}} 
\| \partial_\lambda^{k_1} \partial_\vphi^{\beta_1} \alpha \ldots \partial_\lambda^{k_m} \partial_\vphi^{\beta_m} \alpha \|_{s + \sigma}\,. 
\end{aligned}
\end{equation}
By the tame product estimate, classical interpolation and   
using that by \eqref{stima funzione cambio di variabile striscia} $\| \alpha \|_{s_0 + |\beta| + \sigma}^{k_0, \gamma} 
\lesssim 1 + \| \eta \|_{s_0 + |\beta| + \sigma}^{k_0, \gamma} \lesssim 1$, one has that for any $s \geq s_0,$
$$
\begin{aligned}
&\| \partial_\lambda^{k_1}\partial_\vphi^{\beta_1} \alpha 
\ldots \partial_\lambda^{k_m}\partial_\vphi^{\beta_m} \alpha \|_{s + \sigma} 
\\& \lesssim_{s, \beta} 
\sum_{i = 1}^m \| \partial_\lambda^{k_i}\alpha \|_{s + |\beta_i| + \sigma} 
\prod_{j \neq i} \| \partial_\lambda^{k_j} \alpha \|_{s_0 + |\beta_j| + \sigma}  
%\\& 
\lesssim_{s, k, \beta}  
\sum_{i = 1}^m \gamma^{- |k_i|}\| \alpha \|_{s + |\beta_i| + \sigma}^{k_0, \gamma} 
\prod_{j \neq i} \gamma^{- |k_j|}\|  \alpha \|_{s_0 + |\beta_j| + \sigma}^{k_0, \gamma}  
\\& 
\stackrel{|k_1| + \ldots + |k_m| = |k|}{\lesssim_{s, k, \beta}} 
\gamma^{- |k|}  \| \alpha \|_{s + |\beta| + \sigma}^{k_0, \gamma}  \big( \| \alpha \|_{s_0 + |\beta| + \sigma}^{k_0, \gamma}\big)^{m - 1} 
%\\& 
\stackrel{\| \alpha \|_{s_0 + |\beta| + \sigma}^{k_0, \gamma} \lesssim 1}{\lesssim_{s,k, \beta}}
 \gamma^{- |k|}\big( 1 + \| \eta \|_{s + |\beta| + \sigma}^{k_0, \gamma}\big)\,. 
\end{aligned}
$$
Hence \eqref{lagrima3} follows by \eqref{stima cambio variabile x y cal A brutto 0} and by using the item $(i)$. 

\noindent
The estimate \eqref{lagrima4} follows by similar arguments. %reasoning as in the proof of \eqref{lagrima2Delta12}.
\end{proof}

In the next lemma we shall provide some tame estimates for the operator $\cA(\vphi) - {\rm Id}$. 
\begin{lemma}{\bf (Estimates on $\mathcal{A}(\vphi)-{\rm Id}$)}\label{stime cal A vphi - Id ca,bio variabile strano striscia}
Let $k \in \N^{\nu + 1}$, $m \in \N^{d + 1}$, $|k| \leq k_0$, $|m| \leq 1$. 
Then there exist $\sigma \equiv \sigma(k_0)  \geq \widetilde{\s} \gg 0$ large enough, with $\widetilde{\s}$ independent of $k_0$, such that 
the linear operator $\partial_{x, y}^m \circ \big( \cA(\vphi) - {\rm Id} \big)$ 
satisfies the following bounds.

\noindent
$(i)$ For any $s \geq 0$, there is $\delta=\delta(s,k_0)$ such that,  if 
$\| \eta \|_{s_0 + \sigma}^{k_0, \gamma} \leq\delta$, 
then 
\[
\begin{aligned}
\Big\|\zeta(\cdot) \partial_\vphi^s \partial_\lambda^k \partial_{x, y}^m\big( \cA(\vphi) - {\rm Id} \big)[u] \Big\|_{L^2_{x, y}} 
&\lesssim_{s,k} 
\gamma^{- |k|} \Big( \| \eta \|_{s_0 + \sigma}^{k_0, \gamma}  
\| u \|_{\cH^{s + |k| + 2}} + \| \eta \|_{s + \sigma}^{k_0, \gamma} \| u \|_{L^2_{x, y}} \Big)\,,
\\
\| \zeta(\cdot)\partial_\vphi^p\Delta_{12} 
\partial_{x, y}^m\big( \cA(\vphi) - {\rm Id} \big)\|_{\mathcal{L}(\mathcal{H}^{p+3};L^2_{x, y})}
&\lesssim_{p} \|\eta_{1}-\eta_{2}\|_{p+\widetilde{\s}+|\beta|}\,,
\;\;\;s_0\leq p\leq s_0+\s-\widetilde{\s}\,,
\end{aligned}
\]
for any $u \in \cH^{s + |k| + 2}$.

\noindent
$(ii)$ For any $\beta \in \N^\nu$, any $s \geq s_0$, $\delta=\delta(s,\beta,k_0)$ such that,
if $\| \eta \|_{s_0 + |\beta| + \sigma}^{k_0, \gamma} \leq \delta$, then 
\begin{align*}
\Big\|\zeta(\cdot) \partial_\lambda^k \partial_\vphi^\beta \partial_{x, y}^m&\big(  \cA(\vphi) - {\rm Id} \big)[u] \Big\|_{\cH^s} 
\\&\lesssim_{s, k, \beta} 
\gamma^{- |k|} \big(  \| \eta \|_{s_0 + |\beta| + \sigma}^{k_0, \gamma} 
\| u \|_{\cH^{s + |\beta| + |k| + 2}} + \| \eta \|_{s + |\beta| + \sigma}^{k_0, \gamma} 
\| u \|_{\cH^{s_0 + |\beta| + |k| + 2}} \big)\,,
\\
\| \zeta(\cdot)\partial_\vphi^\beta\partial_{x, y}^m\Delta_{12}
&\big(  \cA(\vphi) - {\rm Id} \big)\|_{\mathcal{L}(\mathcal{H}^{p+|\beta|+3};\mathcal{H}^{p})}
\lesssim_{p, \beta} \|\eta_{1}-\eta_{2}\|_{p+\widetilde{\s}+|\beta|}\,,\;\;\;s_0\leq p\leq s_0+\s-\widetilde{\s}\,,
\end{align*}
for any $u \in \cH^{s + |\beta|+|k| + 2}$.
%and for some $\widetilde{\s}\leq \s$ and for $p+\widetilde{\s}\leq s_0+\s$, $p\geq s_0$.
\end{lemma}

\begin{proof}
We prove the lemma for $\zeta(\cdot)\partial_y \circ \big( \cA(\vphi) - {\rm Id} \big)$. 
The bounds for the operator $\zeta(\cdot)\partial_{x_i} \circ \big( \cA(\vphi) - {\rm Id} \big)$ can be done similarly. 
By the mean value theorem 
$$
\begin{aligned}
\partial_y \circ \big( \cA(\vphi) - {\rm Id} \big)[u] 
%& = \alpha(\vphi, x) \big(  \cA(\vphi) - {\rm Id} \big)[\partial_y u] 
%\\
& 
=   \int_0^1 \cA_\theta(\vphi)[\partial_{yy} u] g(\vphi,x)(1+\theta g(\vphi,x))\, d \theta\,, 
\quad g(\vphi, x) := %\alpha(\vphi, x)
\big( \alpha(\vphi, x) - 1 \big)
\end{aligned}
$$
where the linear operator $\cA_\theta(\vphi)$ (for $\theta \in [0, 1]$) is defined by 
$$
\cA_\theta(\vphi) [u] := u\Big(x, y\big( 1 + \theta g(\vphi, x) \big) \Big)\,.
$$
By the estimates \eqref{stima funzione cambio di variabile striscia} and by Lemma \ref{algebra striscia}, 
one gets that if $\| \eta \|_{s_0 + \sigma}^{k_0, \gamma} \ll 1$, then 
$\| g \|_s^{k_0, \gamma} \lesssim_s \| \eta \|_{s + \sigma}^{k_0, \gamma}$ for any $s \geq s_0$. 
Then, arguing as in the proof of Lemma \ref{lemma cambio di variabile striscia strano}, 
one obtains that $\cA_\theta(\vphi)$ satisfies the estimates for the items $(i)$ and $(ii)$ uniformly in $\theta \in [0, 1]$ 
and hence one obtains the claimed bounds. 
\end{proof}

\vspace{0.5em}
\noindent
{\bf Composition with the Green's function of the flat surface problem}.
Let $m \in \N^{d + 1}$ with $|m| \leq 1$ and 
let us consider the linear operator (recall \eqref{cambio di variabile striscia} and \eqref{def vphi 0})
\[
\cQ_m(\vphi, y) : =  \zeta(\cdot) \partial_{x, y}^m \circ \big( \cA(\vphi) - {\rm Id}\big) \circ \cL_0\,,
\] 
We have the following.
\begin{lemma}\label{corollario derivata cal A vphi - Id cal L0}
For any $k \in \N^{\nu + 1}$, $|k| \leq k_0$ there exists 
$\sigma \equiv \sigma(k_0) \gg 0$ large enough such that the following holds.
For any $\beta \in \N^\nu$, $s\geq s_0$, there is $\delta=\delta(s,\beta,k_0)$ such that 
if $\| \eta \|_{s_0 + |\beta| + \sigma}^{k_0, \gamma} \leq \delta$, 
then one has
{\small \begin{align}
\| \partial_\lambda^k \partial_\vphi^\beta \cQ_m(\vphi, y)[\langle D \rangle^{- |\beta| - k_0 - 2} \psi] \|_{\cO^s} 
&\lesssim_{s, k, \beta} 
\gamma^{- |k|} 
\big(  \| \eta \|_{s_0 + |\beta| + \sigma}^{k_0, \gamma} \| \psi \|_s 
+ \| \eta \|_{s + |\beta| + \sigma}^{k_0, \gamma} \| \psi \|_{s_0} \big)\,,\label{boundQQm}
\\
\| \partial_\vphi^\beta\Delta_{12}\cQ_m(\vphi, y)\langle D\rangle^{-|\beta|-3}\|_{\mathcal{L}(\mathcal{O}^{p};\mathcal{O}^{p})}
&\lesssim_{p, \beta} \|\eta_{1}-\eta_{2}\|_{p+\widetilde{\s}+|\beta|}\,,%\;\;\;p\geq s_0\,. 
\label{boundQQmDelta12}
\end{align}}
for some $\widetilde{\s}\leq \s$ and for $p+\widetilde{\s}\leq s_0+\s$, $p\geq s_0$.
\end{lemma}
\begin{proof}
We start by proving \eqref{boundQQm}.
Let $\beta \in \N^\nu$, $k \in \N^{\nu + 1}$, $|k| \leq k_0$ and let $\cA_m(\vphi) := \zeta(\cdot) \partial_{x, y}^m \circ \big(\cA(\vphi) - {\rm Id} \big)$. Since $\cL_0$ is $\vphi$-independent, we deduce
\[
\partial_\lambda^k \partial_\vphi^\beta \cQ_m(\vphi, y) = \sum_{\begin{subarray}{c}
k_1 + k_2 = k
\end{subarray}} \partial_\vphi^\beta \partial_\lambda^{k_1}\cA_m(\vphi) \circ \partial_\lambda^{k_2} \cL_0 \,.
\]
Hence, in order to obtain the bound \eqref{boundQQm}, it suffices to estimate the operator
\[
\mathtt{Q}(\vphi,y):=\partial_\vphi^\beta \partial_\lambda^{k_1}\cA_m(\vphi) 
\circ \partial_\lambda^{k_2} \cL_0 \circ 
\langle D \rangle^{- |\beta| - k_0 - 2}\,.
\]
Moreover, recalling Lemma \ref{algebra striscia},
%Remark \ref{rmk:equiOHs}, 
we have
\begin{equation}\label{bottle21}
\|\mathtt{Q}(\vphi,y)[\psi]\|_{\mathcal{O}^s}\simeq_{s} \underbrace{\|\mathtt{Q}(\vphi,y)[\psi]\|_{L^2_\vphi \cH^s}}_{I}+
\underbrace{\|\mathtt{Q}(\vphi,y)[\psi]\|_{H^s_\vphi L^2_{x, y}}}_{II}\,.
\end{equation}

\smallskip

\noindent
{\sc Estimate of $I$ in \eqref{bottle21}.}
%$\| \partial_\vphi^\beta \partial_\lambda^{k_1}\cA_m(\vphi) \circ \partial_\lambda^{k_2} \cL_0  [\langle D \rangle^{- |\beta| - k_0 - 2} \psi] \|_{L^2_\vphi \cH^s}$.} 
By Lemma \ref{stime cal A vphi - Id ca,bio variabile strano striscia}-$(ii)$ and Lemma
 \ref{lemma stima cal L0 sol omogenea laplace}, one gets 
\begin{equation*}
\begin{aligned}
\|\mathtt{Q}(\vphi,y)[\psi]\|_{L^2_\vphi \cH^s}&\lesssim_{s,\beta, k}
 \gamma^{- |k_1|} \Big(  \| \eta \|_{s_0 + |\beta| + \sigma}^{k_0, \gamma} 
 \| \partial_\lambda^{k_2} \cL_0  [\langle D \rangle^{- |\beta| - k_0 - 2} \psi]  \|_{L^2_\vphi \cH^{s + |\beta| + |k_1| + 2}}  
 \\& 
 \quad\qquad 
 + \| \eta \|_{s + |\beta| + \sigma}^{k_0, \gamma} \| \partial_\lambda^{k_2} \cL_0  [\langle D \rangle^{- |\beta| - k_0 - 2} \psi]  
 \|_{L^2_\vphi \cH^{s_0 + |\beta| + |k_1| + 2}} \Big) 
 \\
& \lesssim_{s,\beta,k}
 \gamma^{- |k|} \Big(\| \eta \|_{s_0 + |\beta| + \sigma}^{k_0, \gamma} 
\|   \langle D \rangle^{- |\beta| - k_0 - 2} \psi  \|_{L^2_\vphi {H}^{s + |\beta| + |k| + 2}_x} 
\\& 
\quad\qquad + \| \eta \|_{s + |\beta| + \sigma}^{k_0, \gamma} 
\| \langle D \rangle^{- |\beta| - k_0 - 2} \psi  \|_{L^2_\vphi {H}^{s_0 + |\beta| + |k| + 2}_x}  \Big) 
\\& 
\lesssim_{s,\beta,k} 
\gamma^{- |k|} \Big(\| \eta \|_{s_0 + |\beta| + \sigma}^{k_0, \gamma} \|  \psi  \|_{L^2_\vphi {H}^{s}_x} 
+ \| \eta \|_{s + |\beta| + \sigma}^{k_0, \gamma} \|  \psi  \|_{L^2_\vphi {H}^{s_0 }_x}  \Big) 
\\ & 
\lesssim_{s,\beta,k}
 \gamma^{- |k|} \Big(\| \eta \|_{s_0 + |\beta| + \sigma}^{k_0, \gamma} \|  \psi  \|_{s} 
+ \| \eta \|_{s + |\beta| + \sigma}^{k_0, \gamma} \|  \psi  \|_{s_0}  \Big)\,. 
\end{aligned}
\end{equation*}

\smallskip

\noindent
{\sc Estimate of $II$ in \eqref{bottle21}.} In order to estimate $\|\mathtt{Q}(\vphi,y)[\psi]\|_{H^s_\vphi L^2_{x, y}}$, we shall control, for any $\alpha \in \N^\nu$, $|\alpha| \leq s$,  the terms
\begin{equation}\label{bottle20}
\begin{aligned}
\|\pa_{\vphi}^{\alpha}\mathtt{Q}&(\vphi,y)[\psi]\|_{L^{2}_\vphi L^2_{x, y}}
\\& \lesssim_s 
 \sum_{\alpha_1 + \alpha_2 = \alpha} \Big\| \partial_\vphi^{\beta + \alpha_1} 
 \partial_\lambda^{k_1}\cA_m(\vphi) \circ \partial_\lambda^{k_2} \cL_0  
 [\langle D \rangle^{- |\beta| - k_0 - 2} \partial_\vphi^{\alpha_2} \psi] \Big\|_{L^2_\vphi L^2_{x, y}}\,.
\end{aligned}
\end{equation}
By applying Lemma \ref{stime cal A vphi - Id ca,bio variabile strano striscia}-$(i)$ and Lemma \ref{lemma stima cal L0 sol omogenea laplace}, 
one then gets 
\begin{equation*}
\begin{aligned}
\Big\| \partial_\vphi^{\beta + \alpha_1} \partial_\lambda^{k_1}\cA_m(\vphi) &\circ \partial_\lambda^{k_2} \cL_0  
[\langle D \rangle^{- |\beta| - k_0 - 2} \partial_\vphi^{\alpha_2} \psi] \Big\|_{L^2_\vphi L^2_{x, y}} 
\\&  
\lesssim_{s,\beta,k} 
\gamma^{- |k_1|}  \| \eta \|_{s_0 + \sigma}^{k_0, \gamma}  
\| \partial_\lambda^{k_2} \cL_0  [\langle D \rangle^{- |\beta| - k_0 - 2} 
\partial_\vphi^{\alpha_2} \psi]  \|_{L^2_\vphi \cH^{|\beta| + |\alpha_1| + |k_1| + 2}}  
\\& 
\qquad\quad + \gamma^{- |k_1|}  \| \eta \|_{|\beta| + |\alpha_1| + \sigma}^{k_0, \gamma}  \| \partial_\lambda^{k_2} \cL_0  
[\langle D \rangle^{- |\beta| - k_0 - 2} \partial_\vphi^{\alpha_2} \psi]  \|_{L^2_\vphi L^2_{x, y}} 
\\& 
\lesssim_{s,\beta,k} 
\gamma^{- |k_1|}  \| \eta \|_{s_0 + \sigma}^{k_0, \gamma}  
\| \langle D \rangle^{- |\beta| - k_0 - 2} \partial_\vphi^{\alpha_2} \psi \|_{L^2_\vphi {H}^{|\beta| + |\alpha_1| + |k| + 2}_x}  
\\& 
\qquad \quad + \gamma^{- |k_1|}  \| \eta \|_{|\beta| + |\alpha_1| + \sigma}^{k_0, \gamma}  
\| \langle D \rangle^{- |\beta| - k_0 - 2} \partial_\vphi^{\alpha_2} \psi  \|_{L^2_\vphi H^{|k_2|}_x} 
\\& 
\lesssim_{s,\beta,k}
 \gamma^{- |k_1|}  \| \eta \|_{s_0 + \sigma}^{k_0, \gamma}  
\|  \partial_\vphi^{\alpha_2} \psi \|_{L^2_\vphi {H}^{ |\alpha_1|}_x}  
 + \gamma^{- |k_1|}  \| \eta \|_{|\beta| + |\alpha_1| + \sigma}^{k_0, \gamma}  
 \|  \partial_\vphi^{\alpha_2} \psi  \|_{L^2_\vphi L^2_x} 
 \\& 
 \lesssim_{s,\beta,k} 
 \gamma^{- |k|}  \Big( \| \eta \|_{s_0 + \sigma}^{k_0, \gamma}  \|   \psi \|_{s}  
 +  \| \eta \|_{|\beta| + |\alpha_1| + \sigma}^{k_0, \gamma} \| \psi \|_{|\alpha_2|} \Big)
 \\& 
\lesssim_{s,\beta,k} \gamma^{- |k|}  \Big( \| \eta \|_{s_0 +|\beta| + \sigma}^{k_0, \gamma}  
\|   \psi \|_{s}  +  \| \eta \|_{s + |\beta| + \sigma}^{k_0, \gamma} \| \psi \|_{0} \Big)\,,
\end{aligned}
\end{equation*}
where in the last inequality we used the interpolation
\eqref{2202.2}\footnote{here we used \eqref{2202.2} with 
$a_0=|\beta|+\s$, $b_0=2$, $p=|\alpha_1|$ and $q=|\alpha_2|$ to obtain
\[
\| \eta \|^{k_0,\gamma}_{|\beta| + |\alpha_1| + \sigma} \| \psi \|_{|\alpha_2|}  
\lesssim_{\alpha,\beta}
 \| \eta \|^{k_0,\gamma}_{|\beta| + \sigma} \| \psi \|_{|\alpha|} 
 + \| \eta \|^{k_0,\gamma}_{|\alpha| + |\beta| + \sigma} \| \psi \|_0
\]}.
%
%By interpolation and by the Young inequality $a b \lesssim a^{\frac{|\alpha|}{|\alpha_1|}} + b^{\frac{|\alpha|}{|\alpha_2|}}$, one gets that 
%\[
%\begin{aligned}
%\| \eta \|_{|\beta| + |\alpha_1| + \sigma} \| \psi \|_{|\alpha_2|}  & \lesssim 
%\Big( \| \eta \|_{|\beta| + \sigma} \| \psi \|_{|\alpha|} \Big)^{\frac{|\alpha_2|}{|\alpha|}} 
%\Big(\| \eta \|_{|\alpha| + |\beta| + \sigma} \| \psi \|_0\Big)^{\frac{|\alpha_1|}{|\alpha|}}  
%\\& 
%\lesssim \| \eta \|_{|\beta| + \sigma} \| \psi \|_{|\alpha|} + \| \eta \|_{|\alpha| + |\beta| + \sigma} \| \psi \|_0\,,
%\end{aligned}
%\]
%and hence the latter two chains of inequalities imply that for any $\alpha \in \N^\nu$, $|\alpha| \leq s$ 
%and $\alpha_1 + \alpha_2 = \alpha$
%$$
%\begin{aligned}
% \Big\| \partial_\vphi^{\beta + \alpha_1} \partial_\lambda^{k_1}\cA_m(\vphi) 
%&\circ \partial_\lambda^{k_2} \cL_0  [\langle D \rangle^{- |\beta| - k_0 - 2} \partial_\vphi^{\alpha_2} \psi] \Big\|_{L^2_\vphi L^2_{x, y}} 
%\\& 
%\lesssim_s \gamma^{- |k|}  \Big( \| \eta \|_{s_0 +|\beta| + \sigma}^{k_0, \gamma}  
%\|   \psi \|_{s}  +  \| \eta \|_{s + |\beta| + \sigma}^{k_0, \gamma} \| \psi \|_{0} \Big)\,.
%\end{aligned}
%$$
This implies, recalling \eqref{bottle20},  that also $\|\mathtt{Q}(\vphi,y)[\psi]\|_{H^s_\vphi L^2_{x, y}}$
satisfies the same bound. 
By the discussion above, and using \eqref{bottle21}, we deduce \eqref{boundQQm}.

\noindent
To prove \eqref{boundQQmDelta12} we first note that 
\[
 \partial_\vphi^\beta\Delta_{12}\cQ_m(\vphi, y)\langle D\rangle^{-|\beta|-3}= 
 \partial_\vphi^\beta\Delta_{12}\cA_m(\vphi)\circ\mathcal{L}_0\circ\langle D\rangle^{-|\beta|-3}\,.
\]
Therefore, \eqref{boundQQmDelta12} follows reasoning as above end using estimates on the Lipschitz variation given in 
Lemma \ref{stime cal A vphi - Id ca,bio variabile strano striscia}.
The proof  is then concluded. 
\end{proof}

\vspace{0.5em}
\noindent
{\bf Composition with the full Green operator.}
Our aim is to provide estimates on the composition operators
\begin{equation}\label{op fondamentale stime tame resto DN}
\cT_m(\vphi, y) :=\zeta(y) \cA(\vphi) \circ \partial_{x, y}^m \cK(\vphi, y)\,,\qquad m \in \N^{d + 1}\,,\;\; |m| \leq 1\,.
\end{equation}
\begin{equation}\label{operatorRRmm}
%\mathcal{A}(\vphi)\circ\mathcal{K}(\vphi,y)\,,\qquad 
%\cT_m(\vphi, y) := \cA(\vphi) \circ \partial_{x, y}^m \circ\cK(\vphi, y)\,,
%\qquad 
\cR_m(\vphi, y) : =  \zeta(y)\partial_{x, y}^m \circ \cA(\vphi) \circ \cK(\vphi, y)\,,
\quad  m\in \mathbb{N}^{d+1}\,,\;\;|m|\leq 1\,,
\end{equation}
where $\mathcal{A}(\vphi)$ is in \eqref{cambio di variabile striscia} and $\mathcal{K}(\vphi,y)$ is in \eqref{propagKK}.

The first technical result we need is the following.

\begin{lemma}\label{lemma stima cal A partial x y m cal K vphi y}
Let $m \in \N^{d + 1}$, $|m|\leq 1$, $k \in \N^{\nu + 1}$, $|k| \leq k_0$. 
There exists $\sigma \equiv \sigma(k_0)   \gg 0$ large enough such that
the following holds.
For any $\beta \in \N^\nu$, any $s\geq s_0$, there is $\delta=\delta(s,\beta,k_0)$ such that 
 if  
$\| \eta \|_{s_0 + |\beta| +  \sigma}^{k_0, \gamma} \leq \delta$, then the linear operator  
$\partial_\lambda^k \partial_\vphi^\beta \cT_m(\vphi, y)$ satisfies
the estimates
\begin{align}
\!\!\!\| \partial_\lambda^k \partial_\vphi^\beta \cT_m(\vphi, y) [ \langle D \rangle^{- |\beta| - k_0 - 2}\psi] \|_{\cO^s}  
&\lesssim_{s,k, \beta}  
\gamma^{- |k|} \Big( \| \eta \|_{s_0 + |\beta| + \sigma}^{k_0, \gamma} \| \psi \|_{s} 
+ \| \eta \|_{s + |\beta| + \sigma}^{k_0, \gamma} \| \psi \|_{s_0} \Big)\,,\label{boundTTm}
\\
\| \partial_\vphi^\beta\Delta_{12}\cT_m(\vphi, y)\langle D\rangle^{-|\beta|-3}&\|_{\mathcal{L}(\mathcal{O}^{p};\mathcal{O}^{p})}
\lesssim_{p,\beta} \|\eta_{1}-\eta_{2}\|_{p+\widetilde{\s}+|\beta|}\,,%\;\;\;p\geq s_0\,. 
\label{boundTTmDelta12}
\end{align}
for some $\widetilde{\s}\leq \s$ and for $p+\widetilde{\s}\leq s_0+\s$, $p\geq s_0$.
\end{lemma}
%\begin{proof}
%See Appendix \ref{app:green2}.
%\end{proof}
The  result is quite technical. We first provide estimates on the 
operator $\cT_m(\vphi, y)$ in \eqref{op fondamentale stime tame resto DN} 
in the spaces $\mathcal{H}^{s}$ (see \eqref{def sobolev striscia}). 
Then we conclude the estimates in $\mathcal{O}^{s}$. 
We have the following.
\begin{lemma}\label{lemma stima cal Tm Hs vphi Hs x}
Let $k \in \N^{\nu + 1}$, $|k| \leq k_0$. 
Then there exist $\sigma \equiv \sigma(k_0)\geq \widetilde{\s}   \gg 0$ large enough, with 
$\widetilde{\s}$ independent of $k_0$,
 such that the linear operator 
$\cT_m(\vphi, y)$ satisfies the following estimates. 

\noindent
$(i)$ For any $s \geq s_0$, for any $\beta \in \N^\nu$,
there is $\delta=\delta(s,\beta,k_0)$ such that  if  $\| \eta \|_{s_0 + |\beta| + \sigma}^{k_0, \gamma} \leq\delta$, 
for any $\psi \in H^s_x$, 
\begin{align*}
\| \partial_\lambda^k \partial_\vphi^\beta& \cT_m(\vphi, y) [\psi] \|_{\cH^s}  
\\&\lesssim_{s,k, \beta} 
\gamma^{- |k|} \Big(  \| \eta \|_{s_0 + |\beta| + \sigma}^{k_0, \gamma}\| \psi \|_{H^{s +  |\beta| + |k| + 1}_x} 
+ \| \eta \|_{s +  |\beta| + \sigma}^{k_0, \gamma} \| \psi \|_{H^{s_0  + |\beta| + |k|+ 1}_x} \Big)\,,  
\\
\| \partial_\vphi^\beta\Delta_{12}&\cT_{m}(\vphi,y)\|_{\mathcal{L}(\mathcal{H}^{p+|\beta|+2};\mathcal{H}^{p})}
\lesssim_{p, \beta} \|\eta_{1}-\eta_{2}\|_{p+\widetilde{\s}+|\beta|}\,,\;\;\;s_0\leq p\leq s_0+\s-\widetilde{\s}\,.
\end{align*}

\noindent
$(ii)$ For any $s \geq 0$ there is $\delta=\delta(s,k_0)$ such that if 
$\| \eta \|_{s_0 +  \sigma}^{k_0, \gamma} \leq \delta$, then
\begin{align*}
\| \partial_\lambda^k \partial_\vphi^\beta \cT_m(\vphi, y) [u] \|_{L^2_{x, y}} 
&\lesssim_{s,k} 
\gamma^{- |k|} \Big(  \| \eta \|_\sigma^{k_0, \gamma} \| \psi \|_{H^{s + |k| + 1 }_x} 
+ \| \eta \|_{s +  \sigma}^{k_0, \gamma} \| \psi \|_{H^{|k| + 2}_x} \Big)\,,\;\;\;|\beta|=s\,,
\\
\| \partial_\vphi^\beta\Delta_{12}\cT_{m}(\vphi,y)\|_{\mathcal{L}(\mathcal{H}^{p+3};L^2_{x, y})}
&\lesssim_{p, \beta}
 \|\eta_{1}-\eta_{2}\|_{p+\widetilde{\s}+|\beta|}\,,\;\;\;
s_0\leq p\leq s_0+\s-\widetilde{\s}\,.\label{lagrima2Delta12}
\end{align*}
\end{lemma}
\begin{proof}
To shorten notations we write $\| \eta \|_s$ instead of $\| \eta \|_s^{k_0, \gamma}$. 

\noindent
{\sc Proof of $(i)$.} It follows by the product rule and by applying  Proposition \ref{stima cal Kn equazione di laplace} 
and Lemma \ref{lemma cambio di variabile striscia strano}-$(iii)$.

\noindent
{\sc Proof of $(ii)$.} By the product rule one gets that 
\[
\begin{aligned}
 \| \partial_\lambda^k \partial_\vphi^\beta \cT_m(\vphi, y) \|_{L^2_{x, y}} 
 &\lesssim_{s, k} 
 \sum_{\begin{subarray}{c}
\beta_1 + \beta_2 = \beta \\
k_1 + k_2 = k 
\end{subarray}}  
\| \cT_{\beta_1, \beta_2, k_1, k_2}(\vphi, y) \|_{L^2_{x, y}} 
\\
\cT(\vphi, y) \equiv \cT_{\beta_1, \beta_2, k_1, k_2}(\vphi, y) &:= 
\partial_\lambda^{k_1}\partial_\vphi^{\beta_1} \cA(\vphi) 
\circ  \partial_\lambda^{k_2}\partial_\vphi^{\beta_2} \partial_{x, y}^m \cK(\vphi, y)\,,
\end{aligned}
\]
hence, it suffices to estimate $\| \cT(\vphi, y)[u] \|_{L^2_{x, y}}$ for any $\beta_1 + \beta_2 = \beta$, 
$|\beta|=s$, $k_1 + k_2 = k$. 
%By Proposition \ref{stima cal Kn equazione di laplace} and 
By Lemma \ref{lemma cambio di variabile striscia strano}-$(ii)$, one has that 
\begin{equation*}%\label{prima stima cal T s1 s2}
\begin{aligned}
 &\|  \partial_\lambda^{k_1}\partial_\vphi^{\beta_1} \cA(\vphi) \circ  \partial_\lambda^{k_2} 
 \partial_\vphi^{\beta_2} \partial_{x, y}^m \cK(\vphi, y)[\psi]\|_{L^2_{x, y}}   
 \\&\stackrel{\eqref{lagrima2}}{\lesssim_{s,k}}
  \gamma^{- |k_1|} 
 \Big(  \| \partial_\lambda^{k_2}\partial_\vphi^{\beta_2} \partial_{x, y}^m \cK(\vphi, y)[\psi] 
 \|_{\cH^{|\beta_1| + |k_1|}} 
% \\&\qquad
 + \| \eta \|_{|\beta_1| + \sigma} \| \partial_\lambda^{k_2} \partial_\vphi^{\beta_2} 
 \partial_{x, y}^m \cK(\vphi, y)[\psi] \|_{%{L^2_{x, y}} 
 \mathcal{H}^{1+|k_1|}
 } \Big)\,.
\end{aligned}
\end{equation*}
Moreover, using that $m\leq1$, $|\beta_1|+|\beta_2|=s$ and Proposition \ref{stima cal Kn equazione di laplace}
we deduce
%Then, since $m \leq 1$, one has 
\begin{equation*}%\label{s1 s2 m K nella proof 0}
\begin{aligned}
\| \partial_\lambda^{k_2} \partial_\vphi^{\beta_2} \partial_{x, y}^m& \cK(\vphi, y)[\psi] \|_{\cH^{|\beta_1| + |k_1|}} 
% \leq 
%\|  \partial_\lambda^{k_2} \partial_\vphi^{s_2}\cK(\vphi, y)[\psi]\|_{\cH^{s_1 + |k_1| + 1}}  
%\\& 
\stackrel{\eqref{bottle10}}{\lesssim_{s,k}} \gamma^{- |k_2|} \Big(  \| \eta \|_{|\beta_2| + \sigma} 
\| \psi \|_{{H}^{|\beta_1| + |k| + 1}_x} 
+ \| \eta \|_{s+  \sigma} \| \psi \|_{H^{|k| + 2}_x} \Big)\,,
\end{aligned}
\end{equation*}
and 
\begin{equation*}%\label{s1 s2 m K nella proof 1}
\begin{aligned}
\| \partial_\lambda^{k_2}\partial_\vphi^{\beta_2} \partial_{x, y}^m \cK(\vphi, y)[\psi] \|_{
%L^{2}_{x, y}
 \mathcal{H}^{1+|k_1|}
} 
& \lesssim_{s,k} \gamma^{- |k_2|} \| \eta \|_{|\beta_2| + \sigma} \| \psi \|_{H^{|k_1|+|k_2| + 2}_x}\,.
\end{aligned}
\end{equation*}
%Hence \eqref{prima stima cal T s1 s2}, \eqref{s1 s2 m K nella proof 0}, 
%\eqref{s1 s2 m K nella proof 1} imply that 
By combining the estimates above, and setting $|\beta_1|=s_1$, $|\beta_2|=s_2$, we get
\begin{equation}\label{seconda stima cal T s1 s2}
\begin{aligned}
&\|  \partial_\lambda^{k_1}\partial_\vphi^{s_1} \cA(\vphi) \circ  
\partial_\lambda^{k_2}\partial_\vphi^{s_2} \partial_{x, y}^m \cK(\vphi, y)[\psi]\|_{L^2_{x, y}} 
%\\& 
%\lesssim_s \gamma^{- |k|} \Big(  \|  \partial_\lambda^{k_2}\partial_\vphi^{s_2} \partial_{x, y}^m \cK(\vphi, y)[\psi] \|_{\cH^{s_1 + |k_1|}}  
%%\\& \quad 
%+ \| \eta \|_{s_1 + \sigma} \| \partial_\lambda^{k_2} \partial_\vphi^{s_2} \partial_{x, y}^m \cK(\vphi, y)[\psi] \|_{{L^2_{x, y}}} \Big) 
\\& 
\lesssim_{s,k} \gamma^{- |k|} 
\Big(  \| \eta \|_{s_2 + \sigma} \| \psi \|_{\cH^{s_1 + |k|+ 1}_x} + \| \eta \|_{s+  \sigma} \| \psi \|_{H^{ |k| + 2}_x} 
%\\& \qquad  \qquad\qquad
+  \| \eta \|_{s_1 + \sigma} \| \eta \|_{s_2 + \sigma} \| \psi \|_{H^{|k| + 2}_x} \Big)\,.
\end{aligned}
\end{equation}
By interpolation (recall \eqref{2202.3}), one has the inequalities
\[
\begin{aligned}
 \| \eta \|_{s_1 + \sigma} &\lesssim \| \eta \|_\sigma^{\frac{s_2}{s}} \| \eta \|_{s + \sigma}^{\frac{s_1}{s}}\,, 
\qquad 
\| \eta \|_{s_2 + \sigma} \lesssim \| \eta \|_\sigma^{\frac{s_1}{s}} \| \eta \|_{s + \sigma}^{\frac{s_2}{s}}\,, 
\\
\| \psi \|_{H^{s_1 + |k|+ 1}_x}& \lesssim \| \psi \|_{H^{|k| + 1}_x}^{\frac{s_2}{s}} 
\| \psi \|_{H^{s + |k|+ 1}_x}^{\frac{s_1}{s}} \,,
\end{aligned}
\]
and hence, using also \eqref{centoYoung} with $p=s/s_2, q=s/s_1$, 
%Hence (use also the Young inequality $a b \lesssim a^{\frac{s}{s_1}} + b^{\frac{s}{s_2}}$)
one gets
\begin{equation*}
\begin{aligned}
 \| \eta \|_{s_1 + \sigma} \| \eta \|_{s_2 + \sigma}
 &  \lesssim  
 \| \eta \|_\sigma^{\frac{s_2}{s}} \| \eta \|_{s + \sigma}^{\frac{s_1}{s}}  
 \| \eta \|_\sigma^{\frac{s_1}{s}} \| \eta \|_{s + \sigma}^{\frac{s_2}{s}} 
 \stackrel{s_1 + s_2 = s}{\lesssim} 
 \| \eta \|_\sigma \| \eta \|_{s + \sigma} \stackrel{\| \eta \|_\sigma \ll 1}{\lesssim } 
 \| \eta \|_{s + \sigma}\,, 
 \\
 \| \eta \|_{s_2 + \sigma} \| \psi \|_{\cH^{s_1 + |k| + 1}_x} 
 & \lesssim 
 \Big( \| \eta \|_{s + \sigma}  \| \psi \|_{H^{|k| + 1}_x} \Big)^{\frac{s_2}{s}} 
 \Big( \| \eta \|_\sigma \| \psi \|_{H^{s +|k| + 1}_x} \Big) ^{\frac{s_1}{s}}  
 \\& 
 \lesssim_{s} \| \eta \|_{s + \sigma}  \| \psi \|_{H^{|k| + 1}_x} + \| \eta \|_\sigma \| \psi \|_{H^{s + |k|+ 1}_x}\,.
%& \stackrel{\| \eta \|_\sigma \ll 1}{\lesssim} \| \psi \|_{H^{s + 2}_x} +  \| \eta \|_{s + \sigma}  \| \psi \|_{H^2_x}\,. 
\end{aligned}
\end{equation*}
Therefore, by \eqref{seconda stima cal T s1 s2}, one obtains the estimate 
\begin{equation}\label{terza stima cal T s1 s2}
\begin{aligned}
\|  \partial_\lambda^{k_1}\partial_\vphi^{s_1} \cA(\vphi) \circ  \partial_\lambda^{k_2}\partial_\vphi^{s_2} 
&\partial_{x, y}^m \cK(\vphi, y)[\psi]\|_{L^2_{x, y}} 
\\ &
\lesssim_{s,k} \gamma^{- |k|} \Big(  \| \eta \|_\sigma \| \psi \|_{H^{s + |k|+ 1}_x} +  \| \eta \|_{s + \sigma}  
\| \psi \|_{H^{|k| + 2}_x} \Big)\,,
\end{aligned}
\end{equation}
which implies the claimed bound.
The bound on the Lipschitz variation follows similarly.
\end{proof}

\begin{proof}[{\bf Proof of Lemma \ref{lemma stima cal A partial x y m cal K vphi y}}]
To shorten notations in the proof we write $\| \eta \|_s$ instead of $\| \eta \|_s^{k_0, \gamma}$. 
In order to obtain the bound \eqref{boundTTm}, it suffices to estimate the operator
\[
\mathtt{T}(\vphi,y):=\partial_\lambda^k \partial_\vphi^\beta \cT_m(\vphi, y)\circ \langle D \rangle^{- |\beta| - k_0 - 2}\,.
\]
Moreover, given $\psi \in H^s$, $s \geq s_0$ and  recalling Lemma \ref{algebra striscia},
%Remark \ref{rmk:equiOHs}, 
we shall estimate
\begin{equation}\label{bottleTTm}
\|\mathtt{T}(\vphi,y)[\psi]\|_{\mathcal{O}^s}\simeq_{s} \underbrace{\|\mathtt{T}(\vphi,y)[\psi]\|_{L^2_\vphi \cH^s}}_{I}+
\underbrace{\|\mathtt{T}(\vphi,y)[\psi]\|_{H^s_\vphi L^2_{x, y}}}_{II}\,.
\end{equation}
We estimate the two terms separately. 

\noindent
{\sc Estimate of $I$.} By Lemma \ref{lemma stima cal Tm Hs vphi Hs x}-$(i)$ one has that 
 \begin{equation*}%\label{partial vphi beta cal Tm L2 vphi cal Hs}
\begin{aligned}
\| \partial_\lambda^k \partial_\vphi^\beta &\cT_m(\vphi, y) [ \langle D \rangle^{- |\beta|  - k_0- 2}\psi] \|_{L^2_\vphi \cH^s } 
\\&
\lesssim_{s,k, \beta} 
\gamma^{- |k|} 
\Big(   \| \eta \|_{s_0 + |\beta| + \sigma}\| \langle D \rangle^{- |\beta| - k_0 - 2}\psi \|_{L^2_\vphi H^{s +  |\beta| + |k| + 1}_x}  
\\&\qquad\qquad+ \| \eta \|_{s +  |\beta| + \sigma} \| \langle D \rangle^{- |\beta| - k_0 - 2}\psi  \|_{L^2_\vphi H^{s_0  + |\beta| + |k| + 1}_x} \Big) 
\\& 
\lesssim_{s,k, \beta}
 \gamma^{- |k|} \Big( \| \eta \|_{s_0 + |\beta| + \sigma}\| \psi \|_{L^2_\vphi H^{s }_x} 
+ \| \eta \|_{s +  |\beta| + \sigma} \| \psi  \|_{L^2_\vphi H^{s_0 }_x} \Big)
\\& \lesssim_{s,k, \beta}
 \gamma^{- |k|} \Big( \| \eta \|_{s_0 + |\beta| + \sigma}\| \psi \|_{s} + \| \eta \|_{s +  |\beta| + \sigma} \| \psi  \|_{s_0} \Big)\,. 
\end{aligned}
\end{equation*}

\medskip

\noindent
{\sc Estimate of $II$.} We need to estimate for any $\alpha \in \N^\nu$, $|\alpha| \leq s$, 
\begin{equation}\label{partial alpha beta D beta 2 psi 0}
\begin{aligned}
 \Big\| \partial_\vphi^\alpha \Big( \partial_\lambda^k &\partial_\vphi^\beta \cT_m(\vphi, y) 
 [ \langle D \rangle^{- |\beta| - k_0 - 2}\psi]  \Big) \Big\|_{L^2_\vphi L^2_{x, y}} 
 \\& 
  \lesssim_\alpha 
  \sum_{\alpha_1 + \alpha_2 = \alpha} 
  \Big\|  \partial_\lambda^k \partial_\vphi^{\beta + \alpha_1} \cT_m(\vphi, y) [ \langle D \rangle^{- |\beta| - k_0 - 2} 
  \partial_\vphi^{\alpha_2} \psi]  \Big\|_{L^2_\vphi L^2_{x, y}}\,.
\end{aligned}
\end{equation}
Then, by applying Lemma \ref{lemma stima cal Tm Hs vphi Hs x}-$(ii)$, for any $\alpha_1, \alpha_2 \in \N^\nu$ with $\alpha_1 + \alpha_2 = \alpha,$ one has that 
{\small \begin{equation*}%\label{partial alpha beta D beta 2 psi 1}
\begin{aligned}
\Big\| &\partial_\lambda^k \partial_\vphi^{\beta + \alpha_1} \cT_m(\vphi, y) 
[ \langle D \rangle^{- |\beta| - k_0 - 2} \partial_\vphi^{\alpha_2} \psi]  \Big\|_{L^2_\vphi L^2_{x, y}} 
\\& 
\lesssim_{\alpha, k , \beta} 
\gamma^{- |k|} 
\Big( \| \eta \|_\sigma \| \langle D \rangle^{- |\beta|  - k_0 - 2} \partial_\vphi^{\alpha_2} 
\psi\|_{L^2_\vphi H^{|\alpha_1| + |\beta| + |k| + 1 }_x}  
+ \| \eta \|_{|\beta| + |\alpha_1| +  \sigma} \| \langle D \rangle^{- |\beta| - k_0 - 2} \partial_\vphi^{\alpha_2} \psi \|_{L^2_\vphi H^{|k| + 2}_x} \Big) 
\\& 
\lesssim_{\alpha, k , \beta}
 \gamma^{- |k|} \Big(  \| \eta \|_\sigma \| \partial_\vphi^{\alpha_2} \psi \|_{L^2_\vphi H^{|\alpha_1|}_x} 
%\\& \quad 
+ \| \eta \|_{|\beta| + |\alpha_1| +  \sigma} \|  \partial_\vphi^{\alpha_2} \psi \|_{L^2_\vphi L^2_x} \Big) 
\\& 
\lesssim_{s,k,  \beta} 
\gamma^{- |k|} \Big(  \| \eta \|_\sigma \| \psi \|_{s} + \| \eta \|_{|\beta| + |\alpha_1| +  \sigma} \| \psi \|_{|\alpha_2|} \Big)\,. 
\end{aligned}
\end{equation*}}
By interpolation we have
\begin{equation*}
\begin{aligned}
 & \| \eta \|_{|\beta| + |\alpha_1| +  \sigma} 
 \lesssim 
 \| \eta \|_{|\beta| + \sigma}^{\frac{|\alpha_2|}{|\alpha|}} \| \eta \|_{|\beta| + |\alpha| + \sigma}^{\frac{|\alpha_1|}{|\alpha|}}\,, 
 \qquad 
 \| \psi \|_{|\alpha_2|} 
 \lesssim 
 \| \psi \|_0^{\frac{|\alpha_1|}{|\alpha|}} \| \psi \|_{|\alpha|}^{\frac{|\alpha_2|}{|\alpha|}}\,,
\end{aligned}
\end{equation*}
and hence by the Young inequality $a b \lesssim a^{\frac{|\alpha|}{|\alpha_1|}} + b^{\frac{|\alpha|}{|\alpha_2|}}$,
we deduce
\begin{equation*}
\begin{aligned}
\| \eta \|_{|\beta| + |\alpha_1| +  \sigma} \| \psi \|_{|\alpha_2|} & \lesssim \Big( \| \eta \|_{|\beta| + \sigma} \| \psi \|_{|\alpha|} \Big)^{\frac{|\alpha_2|}{|\alpha|}} \Big( \| \eta \|_{|\beta| + |\alpha| + \sigma} \| \psi \|_0 \Big)^{\frac{|\alpha_1|}{|\alpha|}}  \\
& \lesssim_{\alpha} \| \eta \|_{|\beta| + \sigma} \| \psi \|_{|\alpha|} + \| \eta \|_{|\beta| + |\alpha| + \sigma} \| \psi \|_0 \\
& \lesssim_{s} \| \eta \|_{|\beta| + \sigma} \| \psi \|_{s} + \| \eta \|_{s + |\beta| +  \sigma} \| \psi \|_0\,.
\end{aligned}
\end{equation*}
In conclusion, recalling \eqref{partial alpha beta D beta 2 psi 0}, %\eqref{partial alpha beta D beta 2 psi 1}, 
one obtains the bound
\[
\| \partial_\lambda^k \partial_\vphi^\beta \cT_m(\vphi, y) [ \langle D \rangle^{- |\beta| - k_0 - 2}\psi] \|_{H^s_\vphi L^2_{x, y}} 
\lesssim_{s, k, \beta} 
\gamma^{- |k|} \Big(  \| \eta \|_{|\beta| + \sigma} \| \psi \|_{s} + \| \eta \|_{s + |\beta| +  \sigma} \| \psi \|_0 \Big)\,. 
\]
The bound \eqref{boundTTm} follows by the discussion above recalling \eqref{bottleTTm}.
%Finally, the latter bound together with \eqref{partial vphi beta cal Tm L2 vphi cal Hs} imply the claimed statement. 
The Lipschitz variation \eqref{boundTTmDelta12} follows reasoning similarly.
 \end{proof}

We are now in position to prove the following Lemma.
\begin{lemma}\label{corollario derivata cal A vphi cal K vphi}
Let $m \in \N^{d + 1}$ with $|m| \leq 1$ and let us consider the linear operator $\cR_m(\vphi, y)$
in \eqref{operatorRRmm}.
For any  $k \in \N^{\nu + 1}$, $|k| \leq k_0$ there exists $\sigma \equiv \sigma(k_0) \gg 0$ 
large enough such that the following holds.
For any $\beta \in \N^\nu$, $s\geq s_0$, there is $\delta=\delta(s,\beta,k_0)$ such that,
if $\| \eta \|_{s_0 + |\beta| + \sigma}^{k_0, \gamma} \leq \delta$, 
then 
{\small \begin{align}
\| \partial_\lambda^k \partial_\vphi^\beta \cR_m(\vphi, y)[\langle D \rangle^{- |\beta| - k_0 - 2} \psi] \|_{\cO^s} 
&\lesssim_{s,k , \beta} 
\gamma^{- |k|} \Big(  \| \eta \|_{s_0 + |\beta| + \sigma}^{k_0, \gamma} \| \psi \|_s 
+ \| \eta \|_{s + |\beta| + \sigma}^{k_0, \gamma} \| \psi \|_{s_0} \Big)\,,\label{stimasuRRm}
\\
\| \partial_\vphi^\beta\Delta_{12}\cR_m(\vphi, y)\langle D\rangle^{-|\beta|-3}\|_{\mathcal{L}(\mathcal{O}^{p};\mathcal{O}^{p})}
&\lesssim_{p, \beta} \|\eta_{1}-\eta_{2}\|_{p+\widetilde{\s}+|\beta|}\,,%\;\;\;p\geq s_0\,.
 \label{boundRRmDelta12}
\end{align}}
for some $\widetilde{\s}\leq \s$ and for $p+\widetilde{\s}\leq s_0+\s$, $p\geq s_0$.
\end{lemma}

\begin{proof}%[{\bf Proof of Lemma \ref{corollario derivata cal A vphi cal K vphi}}]
We prove the claimed bounds for the linear operator 
$\cR(\vphi, y) =\zeta(y) \partial_y \circ \cA(\vphi) \circ \cK(\vphi, y)$. 
The bounds for $\zeta(y)\partial_{x_i} \circ \cA(\vphi) \circ \cK(\vphi, y)$, $i = 1, \ldots, d$ 
can be proved by similar arguments.  
A direct calculation, recalling \eqref{cambio di variabile striscia}, shows that 
\[
\cR(\vphi, y) =  \alpha(\vphi, x) \circ \cT(\vphi, y), \quad \cT(\vphi, y) :=  
\zeta(y)\cA(\vphi) \circ \partial_y \cK(\vphi, y)\,. %\quad \alpha(\vphi, x) := \frac{\mathtt h}{\mathtt h + \eta(\vphi, x)}
\]
Moreover, for any $s \geq s_0$, one has that 
\[
\begin{aligned}
\| \partial_\lambda^k \partial_\vphi^\beta \cR(\vphi, y) &[\langle D \rangle^{- |\beta| - k_0 - 2} \psi] \|_{\mathcal{O}^s}
\\&\lesssim_{s,\beta,k}
\sum_{\begin{subarray}{c}
\beta_1 + \beta_2 = \beta \\
k_1 + k_2 = k
\end{subarray}} \| \partial_\lambda^{k_1} \partial_\vphi^{\beta_1}\alpha(\vphi, x) 
\circ \partial_\lambda^{k_2} \partial_\vphi^{\beta_2}\cT(\vphi, y) 
[\langle D \rangle^{- |\beta| - k_0 - 2} \psi] \|_{\cO^s}\,.
\end{aligned}
\]
By the estimate \eqref{stima funzione cambio di variabile striscia} and Lemma \ref{algebra striscia}, 
using the smallness condition on $\eta$,  one has that for any $s \geq s_0$, 
\[
\begin{aligned}
\| \partial_\lambda^{k_1} \partial_\vphi^{\beta_1} \alpha(\vphi, x) \circ \partial_\lambda^{k_2} \partial_\vphi^{\beta_2}&\cT(\vphi, y) 
[\langle D \rangle^{- |\beta| - k_0 - 2} \psi] \|_{\cO^s}  
\\&
\lesssim_s 
\| \partial_\lambda^{k_1} \partial_\vphi^{\beta_1}\alpha(\vphi, x)\|_s 
\| \partial_\lambda^{k_2} \partial_\vphi^{\beta_2}\cT(\vphi, y) [\langle D \rangle^{- |\beta| - k_0 - 2} \psi] \|_{\mathcal{O}^{s_0}}  
\\& 
\qquad\quad 
+  \| \partial_\lambda^{k_1} \partial_\vphi^{\beta_1}\alpha(\vphi, x)\|_{s_0} 
\| \partial_\lambda^{k_2} \partial_\vphi^{\beta_2}\cT(\vphi, y) [\langle D \rangle^{- |\beta| - k_0 - 2} \psi] \|_{\mathcal{O}^{s}} 
\\& 
\lesssim_{s,k, \beta} \gamma^{- |k_1|} \Big( \| \eta \|_{s + |\beta| + \sigma}^{k_0, \gamma} 
\| \partial_\lambda^{k_2} \partial_\vphi^{\beta_2}\cT(\vphi, y) [\langle D \rangle^{- |\beta| - k_0 - 2} \psi] \|_{s_0} 
\\& \qquad\quad 
+ \| \partial_\lambda^{k_2} \partial_\vphi^{\beta_2}\cT(\vphi, y) [\langle D \rangle^{- |\beta| - k_0 - 2} \psi] \|_{\mathcal{O}^{s}}\Big)\,.
\end{aligned}
\]
The  bound \eqref{stimasuRRm} follows combining the estimate \eqref{boundTTm} in  
Lemma \ref{lemma stima cal A partial x y m cal K vphi y}. 
The bound \eqref{boundRRmDelta12} on the Lipschitz variation follows similarly using \eqref{boundTTmDelta12}.
\end{proof}

\begin{rmk}
We remark that the  estimates \eqref{boundTTm}-\eqref{boundTTmDelta12} (and consequently also in 
\eqref{stimasuRRm}-\eqref{boundRRmDelta12}) are not sharp.
One can improve such bounds by substituting the operators $\langle D\rangle^{-|\beta|-k_0-2}$,
$\langle D\rangle^{-|\beta|-3}$
with $\langle D\rangle^{-|\beta|-k_0-1}$,
$\langle D\rangle^{-|\beta|-2}$ respectively. 
The proof of the bounds with a loss of $|\beta|+k_0+2$ is slightly simpler. 
Moreover, this simplification does not affect 
the main result in Proposition \ref{prop:combo} for which  we must  consider 
contributions coming from the operators $\cQ_m(\vphi, y)$, whose bounds are given in Lemma
 \ref{corollario derivata cal A vphi - Id cal L0}, and which,  on the contrary, cannot be improved.
\end{rmk}

\subsubsection{Estimates on g: concluded.} We now conclude the proof of  the main result of the section.

\begin{proof}[{\bf Proof of Proposition \ref{prop:combo}}]
Recalling \eqref{strutturaGG} (see also \eqref{cambio di variabile striscia}, \eqref{propagKK} and \eqref{def vphi 0})
we shall write
\[
\begin{aligned}
M(\vphi,y)&= 
\Big( 2 \chi'(y) \partial_y + \chi''(y) + \chi'(y) \beta_1\cdot \nabla + \chi'(y) \beta_2  \Big) 
\circ \Big(\cA(\vphi) \circ \big(\cL_0+\cK(\vphi, y)\big) -\mathcal{L}_0\Big)  
\\& 
\qquad +  (1 - \chi(y))\Big( \beta_1 \cdot \nabla \partial_y + \beta_2 \partial_y + \beta_3 \Delta \Big) \cL_0
\\&
=: M_1(\vphi, y) + M_2(\vphi, y) + M_3(\vphi, y) + M_4(\vphi, y)\,,
\end{aligned}
\]
where we defined 
\[
\begin{aligned} 
 M_1(\vphi, y) &:= 2 \chi'(y) \partial_y \Big( \cA(\vphi) - {\rm Id}\Big) \circ \cL_0  
 + 2 \chi'(y) \partial_y \cA(\vphi) \circ \cK(\vphi, y)\,, 
 \\
 M_2(\vphi, y) &:= \chi'(y) \beta_1 \cdot \nabla \Big( \cA(\vphi) - {\rm Id}\Big) \circ \cL_0 
 + \chi'(y) \beta_1 \cdot \nabla \cA(\vphi) \circ \cK(\vphi, y) \,, 
 \\ 
 M_3(\vphi, y)& := \Big(\chi''(y) + \chi'(y) \beta_2 \Big) \Big( \cA(\vphi) - {\rm Id}\Big) \circ \cL_0 
 + \Big(\chi''(y) + \chi'(y) \beta_2 \Big) \cA(\vphi) \circ \cK(\vphi, y)\,, 
 \\ M_4(\vphi, y) &:= (1 - \chi(y))\Big( \beta_1 \cdot \nabla \partial_y + \beta_2 \partial_y 
 + \beta_3\Delta \Big) \cL_0\,.
\end{aligned}
\]
First of all, since $\chi \equiv 1$ near $y = 0$, we note that
 $\chi', \chi'', 1 - \chi$ satisfy the conditions in \eqref{cutoffZETA}.
 In particular, they
 vanishes near $y = 0$ and 
 hence also $M = M_1 + M_2 + M_3 + M_4$ 
 vanishes in a neighborhood of $y = 0$. 
 Moreover, the operators $M_{i}$, $i=1,2,3,4$, satisfy the bounds \eqref{albero1}-\eqref{albero2}
%Then one has that $M_1, M_2, M_3, M_4$ satisfies the claimed bound 
by using the estimates \eqref{stime beta i coefficienti operatore giusto} on $\beta_1, \beta_2, \beta_3$ 
and by applying Lemmata \ref{corollario derivata cal A vphi cal K vphi}, 
\ref{lemma stima cal L0 sol omogenea laplace} and \ref{corollario derivata cal A vphi - Id cal L0}.
\end{proof}

\subsubsection{Space-time Green's function estimates 2}\label{sec:stimegreen2}
We are now in position to prove the main Proposition \ref{lemma tame laplace partial y u langle D rangle u} 
of this section.
\begin{proof}[{\bf Proof of Proposition \ref{lemma tame laplace partial y u langle D rangle u}}]
Recalling \eqref{op:L} and \eqref{secondFF}, at least at formal level, we observe that $v$ solves \eqref{elliptic4tris} if and only if 
\begin{equation*}
v = 
\big( {\rm Id} + L^{- 1} \circ {F}(\vphi) \big)^{- 1} \circ L^{- 1} [g+h]\,,
\end{equation*}
where $L^{-1}$ is the Green operator for the problem \eqref{freeProbg}, $g$ is in \eqref{func:gg} and $h$ in \eqref{func:hh}.
Furthermore, in view of \eqref{strutturaGG} and \eqref{termineH}
we shall write
\begin{equation}\label{definizionesolV}
v=\cF(\vphi,y)\circ M_1(\vphi,y)[\psi]\,,
\qquad \begin{aligned} M_1(\vphi, y) &:=  \mathtt L(\vphi, y)+ M(\vphi, y)\,,
\\
\cF(\vphi,y)&:=\big( {\rm Id} + L^{- 1} \circ {F}(\vphi) \big)^{- 1} \circ L^{- 1} \,.
\end{aligned}
\end{equation}
Therefore we get formula \eqref{soluzioneV}
by setting
\[
\mathcal{F}_1(\vphi,y)[\cdot]:=\pa_{y}\cF(\vphi,y)\circ M_1(\vphi,y)\,,
\qquad \mathcal{F}_2(\vphi,x)[\cdot]:=\langle D\rangle\cF(\vphi,y)\circ M_1(\vphi,y)\,.
\]
Now, by Propositions \ref{stime mathttL su Os} and \ref{prop:combo}, 
the operator $M_1(\vphi, y)$ satisfies for any 
$\beta \in \N^\nu$, $k \in \N^{\nu + 1}$, $|k| \leq k_0$, 
for any $s \geq s_0$, 
the estimates
{\small \begin{align}
\| \partial_\lambda^k \partial_\vphi^\beta M_1(\vphi, y)[\langle D \rangle^{- |\beta| - k_0 - 2} \psi] \|_{\cO^s}
& \lesssim_{s,k, \beta} 
 \gamma^{- |k|} \Big( \| \eta \|_{s_0 + |\beta| + \sigma}^{k_0, \gamma} \| \psi \|_{s} 
 + \| \eta \|_{s + |\beta| + \sigma}^{k_0, \gamma} \| \psi \|_{s_0} \Big)\,,\label{barzelletta assurda 1}
 \\
\| \partial_\vphi^\beta\Delta_{12}{M}_1(\vphi, y)\langle D\rangle^{-|\beta|-3}\|_{\mathcal{L}(\mathcal{O}^{p};\mathcal{O}^{p})}
&\lesssim_{p, \beta} \|\eta_{1}-\eta_{2}\|_{p+\widetilde{\s}+|\beta|}\,,%\;\;\;p\geq s_0\,.
\label{barzelletta assurda 1Delta12}
\end{align}}
for some $\s\geq \widetilde{\s}$ and $s_0\leq p\leq s_0+\s-\widetilde{\s}$,
provided that
 $\| \eta \|_{s_0 + |\beta| + \sigma}^{k_0, \gamma} \ll_{s,\beta}1$.
So, to get the result, we need to bounds the operators $\pa_{y}\mathcal{F}(\vphi,y)$, $\langle D\rangle \mathcal{F}(\vphi,y)$.
We just give the proof for the first one. The other is similar.

We note that the operator $\pa_{y}\mathcal{F}(\vphi,y)$ (see  \eqref{definizionesolV})
has the same structure of the  operator in
\eqref{ascensore1}. The only difference is in the definition of $\mathcal{A}_{n}(\vphi,y)$ in \eqref{AAnn}
where $\widetilde{F}(\vphi,y)$ (see \eqref{def op F vphi}) must be replaced by $F(\vphi)$ in \eqref{secondFF}.
Now, the latter operators have coefficients $\beta_i$, independent of $y$ 
and satisfying \eqref{stime beta i coefficienti operatore giusto} which are the same estimates 
satisfied by the coefficients $\widetilde{F}(\vphi,y)$ (see \eqref{stimealphai}-\eqref{stimealphaiDelta12}).
So we shall apply 
 Corollary \ref{stima of di green forzante qualunque}
to obtain
%$$
%\partial_y u = \cF(\vphi, y)[f] \quad \text{where} \quad f = M_1(\vphi, y), \quad M_1(\vphi, y) :=  \mathtt L(\vphi, y)+ M(\vphi, y)
%$$
%see \eqref{termnotoG2} and 
that the operator $\pa_{y}\cF(\vphi, y)$ satisfies for any $\beta \in \N^\nu$, $k \in \N^{\nu + 1}$, $|k| \leq k_0$, 
%if $\| \eta \|_{s_0 + |\beta| + \sigma}^{k_0, \gamma} \ll 1$, 
for any $s \geq 0$, the estimates
\begin{align}
\| \partial_\lambda^k \partial_\vphi^\beta \pa_{y}\cF(\vphi, y)[f] \|_{\cO^s} 
&\lesssim_{s,k, \beta}
\gamma^{- |k|} 
\Big( \| f \|_{\cO^s} + \| \eta \|_{s + |\beta| + \sigma}^{k_0, \gamma} \| f \|_{\cO^0} \Big)\,,
\label{barzelletta assurda 0}
\\
\| \partial_\vphi^\beta\Delta_{12}\pa_{y}\cF(\vphi, y)\|_{\mathcal{L}(\mathcal{O}^{p};\mathcal{O}^{p})}
&\lesssim_{p, \beta} \|\eta_{1}-\eta_{2}\|_{p+\widetilde{\s}+|\beta|}\,,%\;\;\;p\geq s_0\,.
\label{barzelletta assurda 0Delta12}
\end{align}
for some $\s\geq \widetilde{\s}$ and $s_0\leq p\leq s_0+\s-\widetilde{\s}$,
provided that
 $\| \eta \|_{s_0 + |\beta| + \sigma}^{k_0, \gamma} \ll_{s,\beta}1$.
Therefore, using that $\partial_y u = \cF_1(\vphi, y)[\psi] = \pa_{y}\cF(\vphi, y) \circ M_1(\vphi, y)[\psi]$
and that
\[
\partial_\lambda^k \partial_\vphi^\beta\pa_{y} \cF_1(\vphi, y) = \sum_{\begin{subarray}{c}
\beta_1 + \beta_2 = \beta \\
k_1 + k_2 = k
\end{subarray}} \partial_\lambda^{k_1} \partial_\vphi^{\beta_1} \pa_{y}\cF(\vphi, y) 
\circ \partial_\lambda^{k_2} \partial_\vphi^{\beta_2} M_1(\vphi, y)\,,
\] 
one obtains \eqref{propagF} 
by applying the estimates \eqref{barzelletta assurda 0}, \eqref{barzelletta assurda 1}.
The Lipschitz variation \eqref{propagFDelta12} follows similarly using \eqref{barzelletta assurda 0Delta12}, \eqref{barzelletta assurda 1Delta12}.
\end{proof}

\section{Pseudo-differential expansion of the Dirichlet-Neumann operator}\label{sec:pseudoDN}
In Sections \ref{sec:tameLaplace} and \ref{sec:Problemastriscianuova} we provide tame estimates for the solution $\phi$ in \eqref{lem:firststright} of the problem
\eqref{elliptic2} from which one deduces estimates on the solution of the original problem \eqref{elliptic1}.
In this section we shall also show that actually the 
Dirichlet-Neumann operator \eqref{eq:112a} (see also \eqref{eq:112abis})
has a pseudo-differential structure. 
\subsection{Decoupling and reduction to the zero boundary condition}\label{sec:decoupling}

In view of Remark \ref{rmk:DNconV}, in order to understand the structure of the Dirichlet-Neumann operator
we shall study the properties of the solution  $\widetilde{\Phi}$ 
(see \eqref{vphi vphi0 u zero cond bordo}) of the problem
\eqref{elliptic5}, close to the boundary $y=0$.

First of all, we shall perform a full decoupling into a forward and 
a backward elliptic evolution equations of the problem \eqref{elliptic5}.
We have the following.
\begin{lemma}\label{divo1}
Fix any $\beta_0,k_0\in\N$ and constants $M,\mathtt{c}$ such that $M\geq\mathtt{c} \geq \beta_0+k_0+2$.
There are symbols $a,A\in S^{1}$ 
such that
\begin{equation}\label{splittingfinale}
(\pa_{y}+\op(a))\circ(\pa_{y}-\op(A))[\cdot]=
\mathcal{L}[\cdot]+\mathcal{R}_{\mathcal{L}}(\eta)[\cdot]\,,
\end{equation}
where $\mathcal{L}$ is in \eqref{op:L}, and $\mathcal{R}_{\mathcal{L}}(\eta)$ is a smoothing operator such that 
the following holds. The linear operators ${\rm Op}(a)$, ${\rm Op}(A)$, ${\mathcal R}_{{\mathcal L}}$ 
are real-to-real operators, and the symbols $a,A$ are momentum preserving.
There exists $\s=\s(M,\beta_0,k_0)\gg0$ such that, 
for any $\bar{s}\geq s_0 + \sigma$ (recall \eqref{s0 Dirichlet Neumann}), $\beta\in \N^{\nu}$, $|\beta|\leq \beta_0$,
there is $\delta=\delta(s,k_0,\beta_0)$ such that,
 if
then smallness condition \eqref{smalleta} holds,
%$\| \eta \|_{s_0 +\sigma} \lesssim_{\bar{s},\beta_0,k_0} 1$,
then

\noindent
$(i)$ {\bf (Asymptotic).}  One has that the symbols $a,A\in S^1$ have the form
\[
\begin{aligned}
a(\vphi,x,\x)&=a_{1}(\vphi,x,\x)+a_{\leq0}(\vphi,x,\x)\,,
\\
A(\vphi,x,\x)&=A_{1}(\vphi,x,\x)+A_{\leq0}(\vphi,x,\x)\,,
\end{aligned}
\]
where $a_{\leq0},A_{\leq0}\in S^{0}$ and 
(recall the definition of the cut-off function $\chi$ in \eqref{cutofffunct})
\begin{equation}\label{simboliordine1}
\begin{aligned}
a_{1}(x,\x)&:= - \ii b \mathtt h \nabla\eta\cdot\x \chi(\xi)
+  \mathtt h \chi(\xi) \sqrt{b|\x|^{2}-(b\nabla\eta\cdot\x)^{2}}\,,
\\
A_{1}(x,\x)&:=\ii b \mathtt h \nabla\eta\cdot\x \chi(\xi)
+ \mathtt h \chi(\xi)\sqrt{b|\x|^{2}-(b\nabla\eta\cdot\x)^{2}}\,,\qquad b:=b(\eta;x):=\frac{1}{1+|\nabla\eta|^{2}}\,,
\end{aligned}
\end{equation}
and satisfy for any $s_0\leq s \leq \bar{s}-\s$, $\alpha \in \N$,
\begin{align}
& \| a_{1}-\th|\x| \chi(\xi)\|_{1, s, \alpha}^{k_0, \gamma}\,,\,  \| A_{1}-\th|\x| \chi(\xi) \|_{1, s, \alpha}^{k_0, \gamma} \,, 
 \| a_{\leq0} \|_{0, s, \alpha}^{k_0, \gamma}\,,\, \| A_{\leq0} \|_{0, s, \alpha}^{k_0, \gamma} 
\lesssim_{s, \alpha, M} \| \eta \|_{s + \sigma}^{k_0, \gamma}\,,\label{a vicino |xi|}
\\
 &\| \Delta_{12}a \|_{1, p, \alpha}\,,\, \| \Delta_{12}A \|_{1, p, \alpha}
 \lesssim_{p,\alpha, M}
 \| \eta_1-\eta_2 \|_{p + \sigma}\,, \;\;\;s_0\leq p\leq s_0+\s-\widetilde{\s}\,,
 \label{a vicino |xi|Delta12}
\end{align}
for some $\widetilde{\s}\leq \s$.
%for any $\alpha\in\N$ .

\noindent
$(ii)$ {\bf (Remainder).} The smoothing operator  $\mathcal{R}_{\mathcal{L}}(\eta)$
satisfies the following.
For any $|\beta|\leq \beta_0$, %$M\leq N-(\beta_0+k_0)-2$
one has that $\partial_\vphi^\beta \langle D \rangle^{M} {\mathcal R}_{\mathcal{L}}
 \langle D \rangle^{-\mathtt{c}}$ 
is a ${\mathcal D}^{k_0}$-tame operator and, for $\bar{s}\geq s_0+\s$,  it satisfies 
\begin{equation}\label{stima tame resto splitting}
\begin{aligned}
{\mathfrak M}_{\partial_\vphi^\beta \langle D \rangle^{M} {\mathcal R}_{\mathcal{L}} 
\langle D \rangle^{-\mathtt{c}}}(s) 
&\lesssim_{ s, M} 
\|\eta\|_{s+\s}^{k_0,\gamma}\,, \quad \forall s_0 \leq s \leq \bar{s}-\s\,,
\\
\|  \langle D \rangle^{M} 
\partial_\vphi^\beta \Delta_{12}{\mathcal R}_{\mathcal{L}}\langle D \rangle^{-\mathtt{c}} 
\|_{\mathcal{L}(H^{p} , H^{p})}
& \lesssim_{p, M}
 \|\eta_1-\eta_2\|_{p+\widetilde{\s}}\,,
\;\;\;s_0\leq p\leq s_0+\s-\widetilde{\s}\,,
\end{aligned}
\end{equation}
for some $\widetilde{\s}\leq \s$. %and for $p+\widetilde{\s}\leq s_0+\s$.
\end{lemma}
\begin{proof} 
To simplify notations in this proof we write $\| \cdot \|$ instead of $\| \cdot \|^{k_0, \gamma}$. Consider the function $b$ in \eqref{simboliordine1}. Then, recalling \eqref{op:L}-\eqref{op laplace raddrizzata buona},
we rewrite the operator $\mathcal{L}(\eta)$ 
as 
\begin{equation}\label{operatorLL}
\mathcal{L}[\cdot]:=\mathcal{L}(\eta)[\cdot]:=
\pa_{yy}-\op(\ii 2 \mathtt h b\nabla\eta \cdot \x \chi(\xi)+ \mathtt h b\Delta\eta)\pa_{y}
-\op(\mathtt h^2 b |\x|^2\chi(\xi)^2)\,.
\end{equation}
Recall that the cut off function $\chi$ is given in \eqref{cutofffunct} and that by \eqref{opcutofffunct}, \eqref{opcutofffunct2}, we use the identifications 
$$
\Delta = - |D|^2 = - |D| \circ |D| = - {\rm Op}(|\xi| \chi(\xi)) \circ {\rm Op}(|\xi|\chi(\xi)) =  - {\rm Op}(|\xi|^2 \chi(\xi)^2)
$$ and $\nabla \eta \cdot \nabla = {\rm Op}(\ii \nabla \eta \cdot \xi \chi(\xi))$. 
Our purpose is to determine two symbols $a,A\in S^{1}$, such that 
\[
(\pa_{y}+\op(a))\circ(\pa_{y}-\op(A))=\pa_{yy}+\op(a-A) \partial_y -\op(a)\circ\op(A) = {\mathcal L} + {\mathcal R}_{\mathcal L}(\eta)
\]
where ${\mathcal R}_{\mathcal L}(\eta)$ is a suitable smoothing operator. 
In order to achieve this purpose, we will solve (see \eqref{operatorLL})
 \begin{equation}\label{sommaprodotto}
\begin{aligned}
a-A&=-\ii 2 \mathtt h b\nabla\eta \cdot \x \chi(\xi)- \mathtt h b\Delta\eta\,,
\\
a\#A&=\th^2 b|\x|^{2} \chi(\xi)^2 + \text{smoothing operator of order}\,\, \rho\,,
\end{aligned}
\end{equation}
for an arbitrary $\rho \in \N$. This parameter $\rho$ will be chosen suitably later, large with $M$.
 In order to solve \eqref{sommaprodotto}, 
we look for symbols of the form
\begin{equation}\label{laghetto3}
a=\sum_{j=0}^{\rho}a_{1-j}\,,
\quad A=\sum_{j=0}^{\rho}A_{1-j}\,,\qquad
a_{1-j},A_{1-j}\in S^{1-j}\,,
\end{equation}
for some $\rho\gg1$ to be chosen later.
Note that, recalling  \eqref{espansionecompostandard},
\[
\op(a)\circ \op(A)=\sum_{j,k=0}^{\rho}\op(a_{1-j})\circ \op(A_{1-k})=\op(c_{\rho})+R_{\rho}\,,
\]
where we defined 
\begin{equation}\label{laghetto2}
c_{\rho}:=\sum_{j,k=0}^{\rho}\sum_{n=0}^{1-j-k+\rho}a_{1-j}\#_{n}A_{1-k}\,,
\qquad
R_{\rho}:=\sum_{j,k=0}^{\rho}a_{1-j}\#_{\geq 2-j-k+\rho}A_{1-k}\,.
\end{equation}
In this way 
\[
(\pa_{y}+\op(a))\circ(\pa_{y}-\op(A))=\pa_{yy}+\op(a-A) \partial_y -\op(c_{\rho})-R_{\rho}
\]

and 
 \begin{equation}\label{sommaprodotto2}
\begin{aligned}
c_{\rho}\stackrel{\eqref{laghetto2}}{=}
\sum_{j,k=0}^{\rho}
\sum_{n=0}^{1-j-k+\rho} a_{1-j}\#_{n}A_{1-k}&=\th^2b|\x|^{2} \chi(\xi)^2 + r_\infty\,,
\\
a-A&=-2\ii \th b\nabla\eta\cdot\x \chi(\xi)-\th b\Delta\eta\,,
\end{aligned}
\end{equation}
where $r_\infty \in S^{- \infty}$ is a smoothing pseudo-differential operator which has to be determined. In this way, \eqref{splittingfinale} holds true by setting $\mathcal{R}_{\mathcal{L}}:=-R_{\rho} - {\rm Op}(r_\infty)$. 
We can solve the equations above iteratively. 

\smallskip
\noindent
{\bf Equation at order $1$.} Equations \eqref{sommaprodotto2} for the highest order symbols $a_1$, $A_1$ 
have the form
\[
\begin{aligned}
a_{1} A_{1}&=\th^2 b|\x|^2\chi(\xi)^2\,,
\qquad \quad
a_1-A_1&=-2\ii \th b\nabla\eta\cdot\x \chi(\xi)\,,
\end{aligned}
\]
whose solutions are given by 
\begin{equation}\label{AAord1}
\begin{aligned}
a_1&=-\ii \th b\nabla\eta\cdot\x \chi(\xi)+\th \chi(\xi)\sqrt{b|\x|^2-(b\nabla\eta\cdot\x)^2}\,,
\\
A_1&=\ii \th b\nabla\eta\cdot\x \chi(\xi)+\th \chi(\xi)\sqrt{b|\x|^2-(b\nabla\eta\cdot\x)^2}\,.
\end{aligned}
\end{equation}
Note that 
\begin{equation}\label{reality A1 a1 PARAMETRIX DN}
a_1(\vphi, x, \xi) = \overline{a_1(\vphi, x, - \xi)} \qquad \text{and} \qquad A_1(\vphi, x, \xi) = \overline{A_1(\vphi, x, - \xi)}\,,
\end{equation} implying that ${\rm Op}(a_1)$ and ${\rm Op}(A_1)$ are real-to-real operators. 
The bounds \eqref{a vicino |xi|}-\eqref{a vicino |xi|Delta12} follow by an explicit computation recalling 
the definition of $b$ in \eqref{simboliordine1}.
Notice that 
\[
%{\rm Re}(a_1) = {\rm Re}(A_1)=
a_1+A_1=2\th \chi(\xi)\sqrt{b|\x|^2-(b\nabla\eta\cdot\x)}
\]
is an \emph{elliptic} symbol for $\| \eta \|_{L^\infty} \ll 1$ small enough, namely 
\begin{equation}\label{prop ellittico simboli a1 + A1}
a_1(\vphi, x, \xi) + A_1(\vphi, x, \xi) \gtrsim |\xi| \gtrsim 1, \quad \forall |\xi| \geq 1, \quad (\vphi, x) \in \T^\nu \times \T^d_\Gamma\,. 
\end{equation}

\smallskip
\noindent
{\bf Equation at order  $0$.} Equations \eqref{sommaprodotto2} for the sub principal symbols $a_0$, $A_0$ reads
\[
\begin{aligned}
A_0&=a_0+b\Delta\eta
\\
a_{0} A_{1}+a_{1}A_0&=-a_{1}\#_{1}A_1+S^{-\infty}
 \stackrel{\eqref{espansionecompostandard}}{=}\ii  (\nabla_{\x}a_1) \cdot (\nabla_{x}A_1) + S^{- \infty}\,.
\end{aligned}
\]
  By substituting the expression of $A_0$ in the second equation, one obtains then the following equation for $a_0$ 
  \begin{equation}\label{eq per a0 splitting laplace}
  \begin{aligned}
  (A_1 + a_1)a_{0}   & =  f_0 + S^{- \infty} \,, \qquad  f_0  := \ii  (\nabla_{\x}a_1) \cdot (\nabla_{x}A_1)  - a_1 b \Delta \eta\,. 
  \end{aligned}
  \end{equation}
  By the identity \eqref{reality A1 a1 PARAMETRIX DN}, one gest that $f_0$ satisfies 
    \begin{equation}\label{reality f0 lemma splitting DN}
  f_0(\vphi, x, \xi) = \overline{f_0(\vphi, x, - \xi)}\,. 
  \end{equation}
  By the bounds \eqref{a vicino |xi|}-\eqref{a vicino |xi|Delta12}, the definition of $b$ in 
  \eqref{simboliordine1} and the smallness \eqref{smalleta}, 
  one can show that 
  %for $\sigma \gg 0$ large enough, $\| \eta \|_{s_0 + \sigma}^{k_0, \gamma} \ll 1$ small enough, 
  \begin{equation}\label{bound f0 lemma splitting laplace}
  f_0 \in S^1 \quad \text{and} \quad  \| f_0 \|_{1, s, \alpha} 
  \lesssim_{s, M}
   \| \eta \|_{s + \sigma}, \quad \forall s \geq s_0\,. 
  \end{equation}
  We now solve the equation \eqref{eq per a0 splitting laplace}. Let us introduce a cut-off function $\chi_0 \in C^\infty(\R^d, \R)$ with the following properties: 
 $$
\begin{aligned}
& 0 \leq \chi_0\leq 1, \quad \chi_0(\xi) = \chi_0(- \xi)\,, \\
&  \chi_0(\xi) = 0 \quad \text{if} \quad |\xi| \leq 1\,, \\
& \chi_0(\xi) = 1 \quad \text{if} \quad |\xi| \geq 2\,. 
\end{aligned}
$$
We split the symbol $f_0 $ as 
$$
f_0 =  \chi_0(\xi) f_0  + (1 - \chi_0(\xi)) f_0\,.
$$
Clearly $r_\infty^{(0)} := (1 - \chi_0(\xi)) f_0 \in S^{- \infty}$ since $1 - \chi_0 \equiv 0$ for $|\xi| \geq 2$ and by \eqref{bound f0 lemma splitting laplace} 
  \begin{equation}\label{bound f0 lemma splitting laplaceA}
   \| r_\infty^{(0)}  \|_{- \rho , s, \alpha} \lesssim_{s,\al, \rho,M} \| \eta \|_{s + \sigma}\,, \quad \forall s \geq s_0\,, \quad \forall \rho \,, \,\alpha \in \N\,.  
  \end{equation}
Note that for $|\xi| \geq 1$ (which is the support of $\chi_0$), the symbol $a_1 + A_1$ satisfies the bound \eqref{prop ellittico simboli a1 + A1}, therefore, by using the estimates \eqref{a vicino |xi|} for $a_1, A_1$, one gets that 
\begin{equation}\label{prop simboli a_1 + A_1 inv}
\dfrac{\chi_0}{a_1 + A_1} \in S^{- 1}, \quad \Big\|\dfrac{\chi_0}{a_1 + A_1}\Big\|_{- 1, s, \alpha} 
\lesssim_{s, \alpha, M}
 1 + \| \eta \|_{s + \sigma}, \quad \forall s \geq s_0, \quad \forall \alpha \in \N\,. 
\end{equation}
Hence, we can solve the equation \eqref{eq per a0 splitting laplace}, by solving 
$$
(A_1 + a_1)a_{0}    = \chi_0(\xi) f_0\,. 
$$
We then define
\begin{equation}\label{AAord0}
\begin{aligned}
a_{0}&: =\dfrac{\chi_0(\xi) f_0}{a_1+A_1}\,,
\qquad \quad
A_{0}: =a_0+b\Delta\eta
\end{aligned}
\end{equation}
and hence 
\[
\begin{aligned}
A_0&=a_0+b\Delta\eta
\\
a_{0} A_{1}+a_{1}A_0&= r_\infty^{(0)}\,.
\end{aligned}
\]
By the properties \eqref{reality A1 a1 PARAMETRIX DN}, \eqref{reality f0 lemma splitting DN}, one gets that $a_0$ and $A_0$ satisfies 
\begin{equation}\label{reality A0 a0 PARAMETRIX DN}
a_0(\vphi, x, \xi) = \overline{a_0(\vphi, x, - \xi)}, \quad A_0(\vphi, x, \xi) = \overline{A_0(\vphi, x, - \xi)}\,,
\end{equation}
implying that ${\rm Op}(a_0)$ and ${\rm Op}(A_0)$ are real-to-real operators. One has that $a_0 , A_0 \in S^0$, since $f_0 \in S^1$ and $\frac{\chi_0}{a_1 + A_1} \in S^{- 1}$ 
and by using the bounds \eqref{bound f0 lemma splitting laplace},  \eqref{prop simboli a_1 + A_1 inv} 
(using also \eqref{smalleta})
%(using that $\| \eta \|_{s_0 + \sigma} \ll 1$) 
together with Lemma \ref{lemmacomposizioneSTANDARD} one gets that $a_0, A_0$ satisfies the estimates
\begin{equation}\label{stime a0 A0 nel lemma splitting laplace}
\| a_0 \|_{0, s, \alpha}\,,\, \| A_0 \|_{0, s, \alpha} 
\lesssim_{s, \alpha, M}
 \| \eta \|_{s + \sigma}, \quad \forall s \geq s_0, \quad \forall \alpha \in \N\,. 
\end{equation} 

\smallskip
\noindent
{\bf Equation at order $\leq -1$.} For fixed $2\leq q\leq \rho$, the equations for the symbols $a_{1-q},A_{1-q}$, have the form
\begin{equation}\label{simbolor2menorho}
\begin{aligned}
A_{1-q}&=a_{1-q}
\\
a_{1-q} A_{1}+a_{1}A_{1-q}&=r_{2-q} + S^{- \infty}\,,
\qquad 
r_{2-q}:=-\sum_{\substack{ j+k+n=q\\ 0\leq j,k\leq q-1}}a_{1-j}\#_{n}A_{1-k}\in S^{2-q}\,,
\end{aligned}
\end{equation}
which can be solved iteratively since the term $r_{2-q}$ depends only on $a_{1-j}, A_{1-k}$ with $0 \leq j,k\leq q-1$, which have been already computed at the previous steps.
Then we reason inductively on $q$. The base of the induction has been already proved in 
\eqref{AAord1} and \eqref{AAord0}.
Assume now that 
$a_{1-j}, A_{1-k}$ with $j,k\leq q-1$ satisfy the bounds \eqref{a vicino |xi|}-\eqref{a vicino |xi|Delta12}. and they satisfy the reality condition
\begin{equation}\label{reality induttiva parametrix DN}
a_{1-j}(\vphi, x, \xi) = \overline{a_{1-j}(\vphi, x, - \xi) }, \quad A_{1-j}(\vphi, x, \xi) = \overline{A_{1-j}(\vphi, x, - \xi) }, 
\end{equation}
(i.e. ${\rm Op}(a_{1-j})$, ${\rm Op}(A_{1-j})$, $0 \leq j \leq q - 1$ are real-to-real operators). 
By \eqref{simbolor2menorho} and Lemma \ref{lemmacomposizioneSTANDARD} 
(see \eqref{stimacancellettoesplicitoAlgrammo}) we deduce that 
\begin{equation}\label{stima r q - 2 lemma splitting DN}
\|r_{2-q}\|_{2-q,s,\alpha} \lesssim_{s, \alpha, M} \| \eta \|_{s + \sigma} \,,
 %\| \Delta_{12}r_{2-q} \|_{2-q, p, \alpha}
% &
% \lesssim_{p,\alpha}\| \eta_1-\eta_2 \|_{p + \sigma}^{k_0, \gamma}\,,
\end{equation}
for some $\s \gg 0$ large enough. Moreover since \eqref{reality induttiva parametrix DN} holds, one easily verify that 
$$
(a_{1-j}\#_{n}A_{1-k})(\vphi, x, \xi) = \overline{(a_{1-j}\#_{n}A_{1-k})(\vphi, x, - \xi)}
$$  
which implies that 
\begin{equation}\label{reality r q - 2 parametrix DN}
r_{q - 2}(\vphi, x, \xi) = \overline{r_{q - 2}(\vphi,x, - \xi)}\,.
\end{equation}
One can immediately check that $r_\infty^{(1 - q)} := (1 - \chi_0(\xi)) r_{2 - q} \in S^{- \infty}$ 
(since $1 - \chi_0 \equiv 0$ for $|\xi| \geq 2$) and by \eqref{stima r q - 2 lemma splitting DN}
\begin{equation}\label{stima R q - 2 infty splitting DN}
\| r_\infty^{(1 - q)} \|_{- \rho, s, \alpha} 
\lesssim_{\rho, s, \alpha, M}
 \| \eta \|_{s + \sigma}, \quad \forall s \geq s_0, \quad \forall \rho, \alpha \in \N\,. 
\end{equation}
We define 
\[
a_{1-q}=A_{1-q}=\frac{\chi_0(\xi) r_{2-q}}{a_1+A_1}\,,
\]
and hence we solve \eqref{simbolor2menorho} (with the smoothing symbol $R_\infty^{(2 - q)}$). Hence, $a_{1 - q} = A_{1 - q} \in S^{1 - q}$ and by the estimates \eqref{prop simboli a_1 + A_1 inv}, \eqref{stima r q - 2 lemma splitting DN}, using also Lemma \ref{lemmacomposizioneSTANDARD} , one gets that 
\begin{equation}\label{stime a A 1 - q splitting DN}
\| a_{1 - q} \|_{1 - q, s, \alpha}\,,\, \| A_{1 - q} \|_{1 - q, s, \alpha}
\lesssim_{q, s, \alpha, M} 
\| \eta \|_{s + \sigma}, \quad \forall s \geq s_0, \quad \forall  \alpha \in \N\,. 
\end{equation} 
The discussion above proves the claim \eqref{sommaprodotto2} 
on $c_{\rho}$ with $r_\infty := \sum_{q = 1}^\rho r_\infty^{(1 - q)}$. 
By \eqref{reality A1 a1 PARAMETRIX DN}, \eqref{reality r q - 2 parametrix DN} 
one immediately deduces that $a_{1 - q}$ and $A_{1 - q}$ satisfy the condition
 \eqref{reality induttiva parametrix DN}. Hence, $a = \sum_{j = 0}^\rho a_{1 - j}$ and $A = \sum_{j = 0}^\rho A_{1 - j}$ satisfy $a(\vphi, x, \xi) = \overline{a(\vphi, x, - \xi)}$ and $A(\vphi, x, \xi) = \overline{A(\vphi, x, - \xi)}$ and hence ${\rm Op}(a)$, ${\rm Op}(A)$ are real-to-real operators. 
 By the estimates \eqref{stima R q - 2 infty splitting DN} for any $ 1 \leq q \leq \rho$, one gets that $r_\infty$ satisfies 
\begin{equation}\label{stima r infty finale splitting DN}
\| r_\infty \|_{- \rho, s , \alpha} 
\lesssim_{\rho, s, \alpha, M}
 \| \eta \|_{s + \sigma}, \quad \forall s \geq s_0, \quad \forall \alpha \in \N\,. 
\end{equation}
Moreover, the remainder $R_{\rho}$ in \eqref{laghetto2}, is pseudo-differential and, using the bounds
\eqref{a vicino |xi|}-\eqref{a vicino |xi|Delta12} on $a_{1-j},A_{1-k}$, the estimates \eqref{bound f0 lemma splitting laplaceA}, \eqref{stima R q - 2 infty splitting DN} and Lemma \ref{lemmacomposizioneSTANDARD}
(see \eqref{stima:resto composizioneSTANDARD}) one can deduce that 
$R_{\rho}$ satisfies 
\begin{equation}\label{stima R rho finale splitting DN}
\|R_{\rho}\|_{-\rho,s,\alpha} 
\lesssim_{\rho, s, \alpha, M} 
\| \eta \|_{s + \sigma}\,, \quad \forall s \geq s_0, \quad \forall \alpha \in \N\,. 
%\\
% \| \Delta_{12}R_{\rho} \|_{-\rho, p, \alpha}
% &
% \lesssim_{p,\alpha}\| \eta_1-\eta_2 \|_{p + \sigma}^{k_0, \gamma}\,.
\end{equation}
Therefore, by applying the estimates \eqref{stima r infty finale splitting DN}, \eqref{stima R rho finale splitting DN},   Lemma \ref{constantitamesimbolo} to ${\mathcal R}_{\mathcal L}(\eta) = R_{\rho} + {\rm Op}(r_\infty)$, taking $\rho$ sufficiently large w.r.t. $M,\beta_{0},k_0$
one deduce that ${\mathcal R}_{\mathcal L}(\eta)$ is $\mathcal{D}^{k_0}$-tame and fulfils estimates
\eqref{stima tame resto splitting}. Moreover, since ${\mathcal L}$, ${\rm Op}(a)$, ${\rm Op}(A)$ are real-to-real operators, then ${\mathcal R}_{\mathcal L}(\eta)$ is a real-to-real operator too. 
In particular, in view of the explicit construction of $a,A$, recalling 
Remark
\ref{rmk:algsimboli} and that $\eta$ is a traveling wave, one deduce that $a,A$ are mometum preserving symbols.
The estimates for $\Delta_{12}a, \Delta_{12} A$ and $\Delta_{12} {\mathcal R}_{\mathcal L}(\eta)$ follow by similar arguments. The proof of this lemma is then concluded.
\end{proof}

In view of Lemma \ref{divo1}, we shall rewrite problem \eqref{elliptic5} as
\begin{equation}\label{elliptic44}
\left\{\begin{aligned}
(\pa_{y} + \op(a))\circ(\pa_{y}-\op(A))\widetilde{\Phi}&= g+\mathcal{R}_{\mathcal{L}}[\widetilde{\Phi}] \,,
\qquad x\in\T^{d}_\Gamma\,,\;\; -1<y<0\,,\\
\widetilde{\Phi}(x,0)&=\psi\,,\\ 
(\pa_{y}\widetilde{\Phi})(x,-1)&=0\,,
\end{aligned}\right.
\end{equation}
where $g$ is given  in \eqref{func:gg}.
We recall again that we need to provide a precise structure of the normal derivative 
$\pa_{y}\widetilde{\Phi}$ close to the boundary $y\sim0$. To do this, we introduce 
the function
\begin{equation}\label{funz:ww}
w:=\chi(y)\big(\pa_{y}-{\rm Op}(A)\big)\widetilde{\Phi}\,,
\end{equation}
where $\chi$ is the cut-off function  given in \eqref{cut-off}.
A direct computation implies the following simple lemma.
\begin{lemma}\label{lem:strutturaDNconw}
The Dirichlet-Neumann operator in \eqref{eq:112a}-\eqref{eq:112aquatuor} has the form
\begin{equation}\label{eq:112TOTALE}
\begin{aligned}
 G(\eta)\psi 
 & = \th^{-1}(1+|\nabla\eta|^2))\op(A) [\psi] +\th^{-1}(1+|\nabla\eta|^2)) w_{| y = 0} -\nabla\eta\cdot\nabla \psi(x)
\end{aligned}
\end{equation}
where $A$ is the symbol given by Lemma \ref{divo1} and $w$ is given in \eqref{funz:ww}.
\end{lemma}
\begin{proof}
By \eqref{cut-off} we have $\chi(y)=1$ near $y = 0$, so that by \eqref{funz:ww} we deduce
\[
\big(\pa_{y}-{\rm Op}(A) \big)\widetilde{\Phi}=w \quad \Rightarrow\quad 
(\pa_{y}\widetilde{\Phi})(x,0)={({\rm Op}(A) \widetilde{\Phi})_{|y=0}}+w(x,0)\,.
\]
for $y$ near zero.
Recalling that $\widetilde{\Phi}$ solves the problem \eqref{elliptic44} we have 
$\widetilde{\Phi}(x,0)\equiv0$, and therefore 
\[
(\pa_{y}\widetilde{\Phi})(x,0)={\rm Op}(A)[\psi]+w(x,0)\,.
\]
Then, formula \eqref{eq:112TOTALE} follows recalling \eqref{eq:112aquatuor}.
\end{proof}

Here and in the following we shall obtain
a formula like
\begin{equation}\label{claimsuWW}
w(x,0)={\rm \; smoothing \;remainder}\,.
\end{equation}
By \eqref{elliptic44}-\eqref{funz:ww}, recalling also \eqref{cut-off},  the function $w$ solves the (parabolic) problem
\begin{equation}\label{Cauchy problem backward calore}
\left\{\begin{aligned}
&(\pa_{y} + \op(a))w= F\,,
\\ 
&w(x,-1)=0\,,
\end{aligned}\right.
\end{equation}
where the forcing $F(x,y)$ is given by 
\begin{equation}\label{forcingFF}
F(x, y)  =\chi(y)( g(x, y)+\mathcal{R}_{\mathcal{L}}[\widetilde{\Phi}])
+\chi'(y) \big(\pa_{y}- {\rm Op}(A) \big)\widetilde{\Phi}\,.
\end{equation}
%In order to achieve this aim, we construct a parametrix 
%for the equation \eqref{Cauchy problem backward calore}. 
The solution of the equation \eqref{Cauchy problem backward calore} has the form 
\begin{equation}\label{soluzione parametrica dirichlet Neumann}
w(\cdot, y) = \int_{- 1}^y\cU(\vphi, y - z) F(\cdot, z)\, d z\,,\qquad \forall\,y \in [- 1, 0]\,,
\end{equation}
where the operator $\cU(\vphi, \tau)$, $\vphi \in \T^\nu$ 
is the flow of the (homogeneous) parabolic equation 
\begin{equation}\label{calore striscia}
\begin{cases}
\partial_{\tau} \cU(\vphi, \tau) +  {\rm Op}\big(a(\vphi, x, \xi) \big) \cU(\vphi, \tau) = 0\,, 
\quad  0 \leq \tau \leq 1 \\
\cU(\vphi, 0) = {\rm Id} \,. 
\end{cases}
\end{equation}
%\eqref{calore striscia}.
We now provide some basic properties of $w$ in \eqref{soluzione parametrica dirichlet Neumann}.

\smallskip
\noindent
{\bf The forcing term $F$.} In view of Section \ref{sec:Problemastriscianuova}, 
we shall rewrite the forcing term $F$ in 
\eqref{forcingFF} in a more convenient way, by rewriting  it  as a linear operator acting 
on the function $\psi$.

First recall that $g$ is given  in \eqref{func:gg}, while $\mathcal{R}_{\mathcal{L}}$
is the remainder given by Lemma \ref{divo1}.
Secondly, by \eqref{vphi vphi0 u zero cond bordo}, we shall write\footnote{we recall that the operator $\mathcal{L}_0$ does not depend on $\vphi\in\T^{\nu}$ but only on the variable $y\in[-1,0]$.}
\[
\widetilde{\Phi}=\vphi_0+v\stackrel{\eqref{def vphi 0}, \eqref{soluzioneV}}{=}
\mathcal{L}_0[\psi]+\cF_0(\vphi, y)[\psi]\,.
\]
Using that 
\[
(\pa_{y}-\op(A))v
=\pa_{y}v-\big({\rm Op}(A)\langle D\rangle^{-1}\big)\langle D\rangle v\stackrel{\eqref{soluzioneV}}{=}
\mathcal{F}_1(\vphi,y)[\psi]-\big({\rm Op}(A)\langle D\rangle^{-1}\big)\circ\mathcal{F}_2(\vphi,y)[\psi]\,.
\]
and 
\[
(\pa_{y}-\op(A))\vphi_0=(\pa_{y}-\op(A))\mathcal{L}_0[\psi]=
(\pa_{y}-\th|D|)\mathcal{L}_0[\psi]-\op(A-\th|\x|)\mathcal{L}_0[\psi]\,,
\]
and, recalling  \eqref{strutturaGG}, we shall write
\begin{equation}\label{F x y espansione forzante calore bla}
\begin{aligned}
F(x, y)
& \stackrel{\eqref{forcingFF}}{=} \chi'(y) \mathfrak{Z}(y)[\psi] + \Gamma (\vphi, y)[\psi]
+S_{\mathcal{L}} (\vphi, y)[\psi]\,, 
\end{aligned}
\end{equation}
where
\begin{equation}\label{eq:opZZ}
\mathfrak{Z}(y)[\cdot]:=(\pa_{y}-\th|D|)\mathcal{L}_0[\cdot]\stackrel{\eqref{def vphi 0}}{=}
\th |D|\Big(\frac{\sinh((1+y)\th|D|)}{\cosh(\th|D|)}-\frac{\cosh((1+y)\th|D|)}{\cosh(\th|D|)}\Big)[\cdot]
\end{equation}
is independent of $\eta$, whereas 
\begin{equation}\label{restiGammaSelle}
\begin{aligned}
\Gamma (\vphi, y)[\cdot] & := \chi(y) M(\vphi, y)[\cdot] + \chi'(y) \mathcal{F}_1(\vphi,y)[\cdot]
-\chi'(y)\big({\rm Op}(A)\langle D\rangle^{-1}\big)\circ \mathcal{F}_2(\vphi,y)[\cdot]
\\&\quad-\chi'(y)\op(A-\th|\x|)\mathcal{L}_0[\cdot]
\\
S_{\mathcal{L}} (\vphi, y)[\cdot]& := \chi(y) \mathcal{R}_{\mathcal{L}}\circ\big(\mathcal{L}_0
+\cF_0(\vphi, y)\big)[\cdot] 
\end{aligned}
\end{equation}
where $\mathcal{F}_i$, $i=0,1,2$ are given by Proposition 
\ref{lemma tame laplace partial y u langle D rangle u}.
By an explicit computation we also note that\footnote{One has
\[
\begin{aligned}
\sinh((1+y)\th|D|)&-\cosh((1+y)\th|D|)
\\&
=\sinh(\th|D|)\Big(\cosh(y\th|D|)-\sinh(y\th|D|)\Big)
-\cosh(\th|D|)\Big(\cosh(y\th|D|)-\sinh(y\th|D|)\Big)
\\&=\big(\sinh(\th|D|)-\cosh(\th|D|)\big)e^{-y\th|D|}\,.
\end{aligned}
\]}
\begin{equation}\label{eq:opZZ2}
\mathfrak{Z}(y)[\cdot]:=e^{-y\th|D|} \th|D|\big(\tanh(\th|D|)-1\big)[\cdot]
\end{equation}
We have the following.

\begin{lemma}\label{lemma tame F Gamma vphi y che palle}
The linear operators $\Gamma (\vphi, y), S_{\mathcal{L}} (\vphi, y)$ 
in \eqref{F x y espansione forzante calore bla} 
satisfy the following properties. 
Fix any $\beta_0,k_0\in\N$ and  $M\geq\mathtt{c} \geq \beta_0+k_0+2$.
%For any  $k \in \N^{\nu + 1}$, $|k| \leq k_0$,
There exists $\sigma \equiv \sigma (k_0,M) \gg 0$
large enough such that
for any $s\geq s_0$
%, and any $\beta \in \N^\nu$, 
there is $\delta=\delta(s,\beta_0,k_0,M)$ such that the following holds.
 If $\| \eta \|_{s_0 + |\beta| + \sigma}^{k_0, \gamma} \leq\delta$, then, for any $|k|\leq k_0$ and 
 $|\beta|\leq \beta_0$,
one has 
\[
\begin{aligned}
\| \partial_\lambda^k \partial_\vphi^\beta \Gamma (\vphi, y) [ \langle D \rangle^{- \tc} \psi] \|_{\cO^s} 
&
\lesssim_{s,M} \gamma^{- |k|} \Big(  \| \eta \|_{s_0 + |\beta| + \sigma}^{k_0, \gamma}\| \psi \|_s 
+ \| \eta \|_{s + |\beta| + \sigma}^{k_0, \gamma} \| \psi \|_{s_0} \Big) \,,
\\
\| \partial_\vphi^\beta\Delta_{12}\Gamma(\vphi, y)\langle D\rangle^{-\tc}
\|_{\mathcal{L}(\mathcal{O}^{p};\mathcal{O}^{p})}
&\lesssim_{p,M} \|\eta_{1}-\eta_{2}\|_{p+\widetilde{\s}+|\beta|}\,,
%\;\;\;p\geq s_0\,,
\end{aligned}
\]
and 
\[
\begin{aligned}
\| \langle D\rangle^{M}\partial_\lambda^k \partial_\vphi^\beta \mathcal{S}_{\mathcal{L}} (\vphi, y) 
[ \langle D \rangle^{- \mathtt{c}} \psi] \|_{\cO^s} 
&
\lesssim_{s,M} \gamma^{- |k|} \Big(  \| \eta \|_{s_0 + |\beta| + \sigma}^{k_0, \gamma}\| \psi \|_s 
+ \| \eta \|_{s + |\beta| + \sigma}^{k_0, \gamma} \| \psi \|_{s_0} \Big) \,,
\\
\| \langle D\rangle^{M}\partial_\vphi^\beta\Delta_{12}\mathcal{S}_{\mathcal{L}} (\vphi, y)
\langle D\rangle^{-\mathtt{c}}\|_{\mathcal{L}(\mathcal{O}^{p};\mathcal{O}^{p})}
&\lesssim_{p,M} \|\eta_{1}-\eta_{2}\|_{p+\widetilde{\s}+|\beta|}\,,
%\;\;\;p\geq s_0\,,
\end{aligned}
\]
for some $\widetilde{\s}\leq \s$ and for $p+\widetilde{\s}\leq s_0+\s$.
Moreover, $\Gamma (\vphi, y) \equiv 0$ vanishes identically near the boundary $y = 0$. 
\end{lemma}

\begin{proof}
First of all note that if $\mathtt{c}\geq \beta_0+k_0+2$ then one also have $\mathtt{c}\geq \beta_0+3$
since $k_0\geq 1$. Moreover, it is easy to note that, if an operator $M$ is such that $M\langle D\rangle^{-c_1}$
is bounded, then also $M\langle D\rangle^{-c_2}$ with $c_2\geq c_1$ is bounded.
Therefore, 
the estimates for $\Gamma$ follow by applying 
Proposition \ref{prop:combo} to estimate $M(\vphi,y)$, 
Proposition \ref{lemma tame laplace partial y u langle D rangle u}
to estimate $\mathcal{F}_i$, $i=1,2$, the bounds  \eqref{a vicino |xi|}-\eqref{a vicino |xi|Delta12} 
(and action on Sobolev of pseudo-differential operators ) to estimate ${\rm Op}(A)$, Lemma 
\ref{lemma stima cal L0 sol omogenea laplace} to estimate $\mathcal{L}_0$ and using 
that $\chi'$ is identically zero near $y = 0$. 
The estimates on $\mathcal{S}_{\mathcal{L}}$ follows combining \eqref{stima tame resto splitting} 
in Lemma \ref{divo1}
with Proposition \ref{lemma tame laplace partial y u langle D rangle u} and 
Lemma \ref{lemma stima cal L0 sol omogenea laplace}.
\end{proof}

\begin{rmk}
Notice that the operators $\Gamma(\vphi,y)$, $\mathcal{S}_{\mathcal{L}}$
are \emph{small} with $\eta$. The only term which is not zero at the flat surface $\eta=0$
is in the operator $\mathfrak{Z}$ in \eqref{eq:opZZ}.
\end{rmk}

\smallskip
\noindent
{\bf The parabolic flow $\mathcal{U}$.} The flow of the parabolic equation
\eqref{calore striscia} satisfies the following tame estimates.

\begin{lemma} \label{stima derivate lambda vphi flusso parabolico}
Fix $ \beta_0, k_0 \in \N$. 
For any $m_1, m_2\in\R$ with $m_1+m_2\geq \beta_0+k_0$ there exists
$ \sigma \equiv \sigma(\beta_0, k_0, m_1, m_2) > 0$ such that the following holds.
For any $\bar{s}\geq s_0+\s$ there exists $\delta=\delta(\bar{s}, \beta_0,k_0)>0$ such that,
if \eqref{smalleta} holds.
%$\|\eta\|_{s_0+\s}^{k_0,\gamma}\leq \delta$, 
then
\begin{align}
\sup_{\tau \in [ 0,1]} \|  \cU(\vphi, \tau) u_0 \|_{s}^{k_0, \gamma}
& \lesssim_s  \| u_0 \|_{s}^{k_0, \gamma} 
+ \| \eta \|^{k_0, \gamma}_{s +\sigma } \| u_0\|_{s_0}^{k_0, \gamma}\,, 
\quad \forall s_0\leq s\leq \bar{s}-\s  \,.
\label{stima flusso parabolico}
\end{align}
Moreover,
for any $\beta \in \N^\nu$ with $|\beta| \leq \beta_0$, 
the operator 
$\langle D \rangle^{- m_1}  \partial_\vphi^\beta \cU(\vphi,\tau)   
\langle D \rangle^{- m_2}$ is $\cD^{k_0}$-tame operator satisfying 
\[
\begin{aligned}
\sup_{\tau \in [ 0,1]}{\mathfrak M}_{\langle D \rangle^{- m_1}  \partial_\vphi^\beta \cU(\cdot,\tau)   
\langle D \rangle^{- m_2}}(s) 
&\lesssim_{\bar{s}, m_{1}, m_{2}} 
1 + \| \eta \|_{s + \sigma}^{k_0, \gamma}\,, 
\qquad \forall s_0 \leq s \leq \bar{s}-\s\,,
\\
\sup_{\tau \in [ 0,1]}\|\Delta_{12}\langle D \rangle^{- m_1}  \partial_\vphi^\beta \cU(\cdot,\tau)   
\langle D \rangle^{- m_2}\|_{\mathcal{L}(H^{p};H^{p})}
&\lesssim_{\bar{s}, m_{1}, m_{2}}
 \|\eta_{1}-\eta_{2}\|_{p+\widetilde{\s}}\,,
\end{aligned}
\]
for some $\widetilde{\s}\leq \s$ and for $p+\widetilde{\s}\leq s_0+\s$.
\end{lemma}
\begin{proof}
See Appendix \ref{sec:backwardflow}.
\end{proof}

In view of formula \eqref{soluzione parametrica dirichlet Neumann} and Lemma
\ref{lem:strutturaDNconw}, the properties of the Dirichlet-Neumann operator
can be deduced studying the function
\begin{equation}\label{soluzione parametrica dirichlet Neumann2}
w(x, 0) =\int_{- 1}^0 \cU(\vphi, -z) F(x, z) [\psi]\, d z\,.
\end{equation}
However, the estimates given by  Lemmata \ref{lemma tame F Gamma vphi y che palle} 
and \ref{stima derivate lambda vphi flusso parabolico} are still not enough to prove the claim 
\ref{claimsuWW}.
In the next section we will provide a more explicit structure on $w(x,0)$.

\subsection{Pseudo-differential expansion and tame estimates}\label{sec:exppseudo}
%\section{The parabolic evolution equation}\label{sec:parabolicproblem}

\noindent
By \eqref{F x y espansione forzante calore bla}, \eqref{soluzione parametrica dirichlet Neumann}, 
\eqref{soluzione parametrica dirichlet Neumann2},
in order to compute $w(x, 0)$, we make the splitting 
\[
w(x, 0) = \cW(\vphi)[\psi] + \cS(\vphi)[\psi]\,,
\]
where 
\begin{equation}\label{def parte dominante parte trascurata w x 0}
\begin{aligned}
 \cW(\vphi)[\psi] & : = \int_{- 1}^0 \cU(\vphi, -y) \chi'(y) \mathfrak{Z}(y)[\psi] \, d y\,, \\
\cS(\vphi)[\psi]& := \int_{- 1}^0 \cU(\vphi, -y) 
\big(\Gamma(\vphi, y)+S_{\mathcal{L}} (\vphi, y)\big)[\psi]\, d y\,,
\end{aligned}
\end{equation}
where  $\Gamma(\vphi,y), S_{\mathcal{L}} (\vphi, y)$ are given in \eqref{restiGammaSelle}.
In this section we prove the following proposition.
\begin{prop}[{\bf Expansion of $w_{| y = 0}$}]\label{lemma espansione pseudo-diff w x 0}
Fix $ \beta_0, k_0 \in \N$, 
$M\in \mathbb{N}$, $\mathtt{c}\geq s_0+\beta_0+k_0+2$
and $M\geq \mathtt{c}$.
Then 
\begin{equation}\label{espansionewzero}
w(x, 0) = \cW(\vphi)[\psi] + \cS(\vphi)[\psi] =\th|D|\big(\tanh(\th|D|)-1\big)[\psi] + \cR_M(\vphi)[\psi]\,,
\end{equation}
where  the remainder $\cR_M(\vphi)$ satisfies the following properties. 
There exists $\sigma \equiv \sigma(M, \beta_0, k_0) \gg 0$ large enough
such that, for any $\bar{s}\geq s_0+\s$ there is  
$\delta \equiv \delta(\bar{s},\beta_0,k_0,M) \ll 1$ small enough such that if 
\eqref{smalleta} holds,
%$\| \eta \|_{s_0 + \sigma}^{k_0, \gamma} \leq\delta$,  
 then for any $\beta \in \N^\nu$, $|\beta| \leq \beta_0$, 
 %$0 \leq M \leq N - \beta_0 - k_0-2$, 
 the operator 
 $\langle D \rangle^M \partial_\vphi^\beta \cR_M(\vphi) \langle D \rangle^{- \mathtt{c}}$ 
 is $\cD^{k_0}$-tame with tame constant satisfying the estimate 
\begin{equation}\label{stimefinalew00}
\begin{aligned}
{\mathfrak M}_{\langle D \rangle^M \partial_\vphi^\beta \cR_M \langle D \rangle^{- \mathtt{c}}}(s) 
&\lesssim_{\bar{s},M}
\| \eta \|_{s + \sigma}^{k_0, \gamma}, \quad \forall s_0 \leq s \leq \bar{s}-\s\,,
\\
\| \partial_\vphi^\beta\Delta_{12}
\langle D \rangle^{M} 
\cR_M
 \langle D \rangle^{-\mathtt{c}} \|_{\mathcal{L}({H}^{p};{H}^{p})}
&\lesssim_{\bar{s},M} 
\|\eta_{1}-\eta_{2}\|_{p+\widetilde{\s}}\,,
\end{aligned}
\end{equation}
for some $\widetilde{\s}\leq \s$ and for $p+\widetilde{\s}\leq s_0+\s$.
\end{prop}

Note that 
$\Gamma(\vphi,y)$ in \eqref{def parte dominante parte trascurata w x 0}
vanishes close to $y = 0$ (hence it gives smoothing contributions), 
whereas $\mathcal{S}_{\mathcal{L}}[\cdot]$ is smoothing by 
the properties of $\mathcal{R}_{\mathcal{L}}$. 
By the estimates in Lemma
\ref{stima derivate lambda vphi flusso parabolico}
we then expect that the operator $\cS(\vphi)$ 
in \eqref{def parte dominante parte trascurata w x 0}
gives 
a smoothing contribution.
We also expect that $\mathcal{W}$ gives a smoothing contribution, with the difference that
it will not be small in $\eta$.

In order to prove these claims we first construct a parametrix of 
the problem \eqref{calore striscia}, showing that 
$\mathcal{U}(\vphi,y)$ has actually a pseudo-differential structure.
This is the content of
subsection \ref{sec:parametrixEVO}.
Then in subsections \ref{sec:analisiW} and \ref{sec:analisiS}
we analyze the structure of the operators  $ \cW(\vphi)$ and $ \cS(\vphi)$ in 
\eqref{def parte dominante parte trascurata w x 0}
and we conclude the proof of Proposition \ref{lemma espansione pseudo-diff w x 0}.
We first start by studying, in subsection \ref{subsec:Poissonsymbols}, some properties of
a special class of symbols.

\subsubsection{Poisson symbols on the strip}\label{subsec:Poissonsymbols}
We start by proving a technical result concerning symbols 
of the form $e^{y a(\vphi, x, \xi)}$, $y < 0$ and $a \in S^1$ satisfying\footnote{Recall that one can also assume 
that it is \emph{traveling} (see Def. \ref{def:mompressimbo}).} 
\begin{equation}\label{propsimboloa}
\begin{aligned}
 a(\vphi, x,  \xi ) &= \mathtt h|\xi| \chi(\xi) + r(\vphi, x, \xi)\,,  \quad r \in S^1 \,, 
\\
 \| r \|_{1, s, \alpha}^{k_0, \gamma}
 & \lesssim_{s, \alpha}
  \| \eta \|_{s + \sigma}^{k_0, \gamma}\,, 
\quad \forall s \geq s_0\,,\forall\, \alpha\in \N\,,
\\
 \| \Delta_{12}r \|_{1, p, \alpha} &\lesssim_{p,\alpha}\| \eta_1-\eta_{2} \|_{p + \sigma}\,.
 \;\;\;p\geq s_0\,,\forall\, \alpha\in \N\,.
\end{aligned}
\end{equation}
\begin{rmk}\label{rmk:simbolopiccolo}
Notice that, by Lemma \ref{divo1}, the symbol $a$ in \eqref{calore striscia} satisfies \eqref{propsimboloa}.
We also remark that the symbols $\pa_{x}^{\ell}a$, for $|\ell|\geq1$ are small in $\eta$ and satisfy estimates
like the one of $r$ in \eqref{propsimboloa}.
\end{rmk}

We prove the following.
\begin{lemma}\label{prop exp simbolo}
Let $a \in S^1$ as in \eqref{propsimboloa}. There is $\s\gg1$ such that, for any $s\geq s_0$, $n\in \N$ the following holds. There is $\delta=\delta(s,n)\ll1$ such that,
if $\| \eta \|_{s_0 + \sigma}^{k_0,\gamma} \leq\delta$, then

\noindent
$(i)$ one has that 
\begin{align}
\sup_{y \leq 0\,,\, \xi \in \R^d}\| y^n \langle \xi \rangle^n e^{y a(\cdot, \xi)}\|^{k_0,\gamma}_{s} 
&\lesssim_{s,n} 1 + \| \eta \|^{k_0,\gamma}_{s + \sigma}\,, 
%\quad \forall s \geq s_0\,,
\label{caramella1}
\\
 \sup_{y\leq 0,\x\in \R^{d}}\| \Delta_{12} y^n \langle \xi \rangle^n e^{y a(\cdot, \xi)}\|_{p} 
 &\lesssim_{p,n}\| \eta_1-\eta_{2} \|_{p + \widetilde{\sigma}}\,, \;\;\;s_0\leq p\leq s_0=\s-\widetilde{\s}\,,
 \label{caramella1Delta12}
\end{align}
for some $\widetilde{\s}\leq \s$.

\noindent
$(ii)$ For any $y < 0$
the symbol $y^n e^{y a}$ belongs to $S^{- n}$ and,   $\alpha\in \N$,
\begin{align}
\sup_{y \leq 0} \| y^n e^{y a} \|^{k_0,\gamma}_{- n , s, \alpha} 
&\lesssim_{s, \alpha,n}
 1 + \| \eta \|^{k_0,\gamma}_{s + \sigma}\,, 
%\qquad s \geq s_0\,,
\label{caramella2}
\\
 \sup_{y\leq 0}\| \Delta_{12} y^n  e^{y a(\cdot, \xi)}\|_{-n+1, p, \alpha} 
 &\lesssim_{p,\alpha,n}
 \| \eta_1-\eta_{2} \|_{p + \widetilde{\sigma}}\,, \;\;\;s_0\leq p\leq s_0+\s-\widetilde{\s}\,,
% \alpha\in \N\,.
 \label{caramella2Delta12}
\end{align}
for some $\widetilde{\s}\leq \s$.

\noindent
$(iii)$ For any $\ell \in \N^d$, $|\ell| \geq 1$ there is $\s_1\gg1$ (depending also on $\alpha$),
such that for any $s\geq s_0$,  there exists $\delta=\delta(s,\ell,k_0)$ such that the following holds.
If $\|\eta\|_{s_0+\s_1}^{k_0,\gamma}\leq\delta$ then
one has that 
$\partial_x^\ell(e^{y a})$ has the following expansion:
%There exists $Q_\alpha \gg 0$ such that 
\[
\partial_x^\ell(e^{y a}) = \sum_{q = 1}^{|\ell|} y^q e^{y a} a_q \quad \text{where} \quad a_q \in S^{q}\,,
\]
and the symbols $a_q$ satisfies the estimates, for $\alpha\in\N$,
\begin{align*}
\| a_q \|^{k_0,\gamma}_{q, s, \alpha} 
&\lesssim_{q, s, \alpha} 
 \| \eta \|^{k_0,\gamma}_{s + \sigma}\,, %\qquad \forall s \geq s_0\,,
\\
\| \Delta_{12} a_q\|_{q, p, \alpha} 
 &\lesssim_{q,p,\alpha}
 \| \eta_1-\eta_{2} \|_{p + \widetilde{\sigma}}\,, \;\;\;s_0\leq p\leq s_0+\s-\widetilde{\s}\,,
 %\;\forall\alpha\in\N\,.
\end{align*}
for some  $\widetilde{\s}\leq \s$.

\end{lemma}
\begin{proof}
We prove the lemma for the norms $\|\cdot\|_{m,s,\alpha}$. The bounds 
for  $\|\cdot\|^{k_0,\gamma}_{m,s,\alpha}$ follow similarly.

\smallskip
\noindent
{\sc Proof of $(i)$.} For convenience, we denote by 
$\vartheta = (\vphi, x) \in \T^\nu \times \T^d_\Gamma$. 
We actually provide a bound of $\| y^n \langle \xi \rangle^n e^{y a} \|_{\cC^s_\vartheta}$, 
$y <0$. For any $\ell \in \N^{\nu + d}$, $|\ell| \leq s$, 
by the Faa di Bruno formula, one has 
\begin{equation}\label{formula simbolo exp 0}
 y^n \langle \xi \rangle^n \partial_\vartheta^\ell (e^{y a} ) 
 =  
 y^n \langle \xi \rangle^n \sum_{1 \leq m \leq |\ell|} 
 \sum_{\ell_1 + \ldots + \ell_m = \ell} C(\ell_1, \ldots, \ell_m) e^{y a} y^m 
 (\partial_\vartheta^{\ell_1} a) \ldots (\partial_\vartheta^{\ell_m} a)\,.
\end{equation}
We first note that for any $1 \leq m \leq |\ell|$, $\ell_1 + \ldots + \ell_m = \ell$, we have 
\[
|y^{m + n} \langle \xi \rangle^{m + n} e^{y a}| \leq C(m, n), \quad \forall y  <0\,, 
\qquad 
(\vphi, x, \xi) \in \T^\nu \times \T^d_\Gamma \times \R^d\,, \;n\in\N\,.
\]
To prove the bound above we reason as follows.
In the case $|\x|\leq 1$  the bound is trivial. For $|\x|>1$ the symbol $a$ is bounded from below by $\mathtt{h}|\x|/2$
(by taking $\| \eta \|_{s_0 + \sigma}^{k_0, \gamma}\ll1$ small enough).
Therefore, we have
\[
 |y^{m + n} \langle \xi \rangle^{m + n} e^{y a}|=|\frac{\langle \xi \rangle^{m + n}}{a^{m+n}}  (ya)^{m + n} e^{y a}|
 \lesssim |(ya)^{m + n} e^{y a}|\lesssim_{m,n} \sup_{z\leq0} |z|^{m+n} e^{z}\,.
\]
Hence the bound follows.
Then one gets 
\begin{equation}\label{stima simbolo exp 1}
\begin{aligned}
\| y^{n + m} \langle \xi \rangle^n e^{y a} 
(\partial_\vartheta^{\ell_1} a) \ldots (\partial_\vartheta^{\ell_m} a)  \|_{\cC^0_\vartheta} 
& \lesssim_{m,n} 
\langle \xi \rangle^{- m} 
\| (\partial_\vartheta^{\ell_1} a) \ldots (\partial_\vartheta^{\ell_m} a)  \|_{\cC^0_\vartheta} 
\\& 
\lesssim_{m,n} 
\langle \xi \rangle^{- m} 
\| a \|_{\cC^{|\ell_1|}_\vartheta} \ldots \| a \|_{\cC^{|\ell_m|}_\vartheta}\,.
\end{aligned}
\end{equation}
By standard interpolation inequality one gets that, 
for any $i = 1, \ldots, m$, 
\[
\| a \|_{\cC_\vartheta^{|\ell_i|}} 
\lesssim
\| a \|_{\cC^0_\vartheta}^{1 - \frac{|\ell_i|}{|\ell|}} 
\| a \|_{\cC^{|\ell|}_\vartheta}^{\frac{|\ell_i|}{|\ell|}}\,.
\] 
Therefore, using  \eqref{stima simbolo exp 1},  $|\ell_1 | + \ldots + |\ell_m| = |\ell|$ and 
the Sobolev embedding $\| \cdot \|_{\cC^k_\vartheta} \lesssim \| \cdot \|_{k + s_0}$,
one obtains that
\begin{equation}\label{stima simbolo exp 2}
\begin{aligned}
\| y^{n + m} \langle \xi \rangle^n e^{y a} 
(\partial_\vartheta^{\ell_1} a) &\ldots (\partial_\vartheta^{\ell_m} a)  \|_{\cC^0_\vartheta} 
  \lesssim 
\langle \xi \rangle^{- m} \prod_{i = 1}^m 
\| a \|_{\cC^0_\vartheta}^{1 - \frac{|\ell_i|}{|\ell|}} 
\| a \|_{\cC^{|\ell|}_\vartheta}^{\frac{|\ell_i|}{|\ell|}}  
\\& 
\lesssim_{m,n}
\langle \xi \rangle^{- m} \| a \|_{\cC^0_\vartheta}^{m  - 1} \| a \|_{\cC^{|\ell|}_\vartheta} 
\\& 
\lesssim_{m,n}
\langle \xi \rangle^{- m} \| a (\cdot , \xi) \|_{s_0}^{m - 1} \| a (\cdot , \xi) \|_{|\ell| + s_0} 
\\& 
\lesssim_{m,n,\ell} 
\Big( \sup_{\xi} \langle \xi \rangle^{- 1} \| a (\cdot, \xi) \|_{s_0} \Big)^{m - 1} 
\Big( \sup_{\xi} \langle \xi \rangle^{- 1} \| a (\cdot, \xi) \|_{|\ell| + s_0} \Big) 
\\& 
\lesssim_{m,n,s} \| a \|_{1, s_0, 0}^{m - 1} \| a \|_{1 , s + s_0, 0} 
\stackrel{\eqref{a vicino |xi|}, \| \eta \|_{s_0 + \sigma } \ll 1}{\lesssim_{s,\alpha, M}} 
1 + \| \eta \|_{s + s_0 + \sigma}\,.
\end{aligned}
\end{equation}
Hence, by \eqref{formula simbolo exp 0}, \eqref{stima simbolo exp 2} 
and using that $\| \cdot \|_{s} \lesssim \| \cdot \|_{\cC^s}$, one obtains \eqref{caramella1}.
%$$
%\|  y^n \langle \xi \rangle^n e^{y a (\cdot, \xi)}  \|_s \lesssim_s 1 + \| \eta \|_{s + s_0 + \sigma}, \quad \forall s \geq s_0
%$$
%and uniformly in $y \leq 0$ and $\xi \in \R^d$. 

\noindent
In order to prove \eqref{caramella1Delta12} we first note that
\begin{equation}\label{expaDelta}
\begin{aligned}
 \Delta_{12} y^n \langle \xi \rangle^n e^{y a(\cdot, \xi)}&=
 y^n \langle \xi \rangle^ne^{y a^{(1)}}\Big(1-e^{y (a^{(2)}-a^{(1)})}\Big)
 \\&=
  y^n \langle \xi \rangle^ne^{y a^{(1)}}\int_{1}^{0}y(a^{(2)}-a^{(1)})e^{yt (a^{(2)}-a^{(1)})}dt\,,
 \end{aligned}
\end{equation}
where we set $a^{(i)}:=a(\eta_i;x,\x)$, $i=1,2$. Then the claimed bound follows by the tame estimates
on the norm $\|\cdot\|_{s}$, and using the bound on the Lipschitz variation in \eqref{propsimboloa}.

\smallskip
\noindent
{\sc Proof of $(ii)$.} In order to prove the bound \eqref{caramella2}, 
for any $\ell \in \N^d$, 
we need to estimate the norm $\| \partial_\xi^\ell (e^{y a(\cdot, \xi)}) \|_s$. 
By the Faa di Bruno formula one has that 
\begin{equation}\label{stima simbolo exp 3a}
y^n \partial_\xi^\ell (e^{y a}) 
= 
\sum_{1 \leq m \leq |\ell|} \sum_{\ell_1 + \ldots + \ell_m = \ell} 
C(\ell_1,\dots, \ell_m) e^{y a} y^{m + n} (\partial_\xi^{\ell_1} a) \ldots (\partial_\xi^{\ell_m} a)\,.
\end{equation}
Note that for any $i = 1, \ldots, m$, one has that 
\begin{equation}\label{stima simbolo exp 3}
\begin{aligned}
\| \partial_\xi^{\ell_i} a(\cdot, \xi) \|_s  & \lesssim 
\langle \xi \rangle^{1 - |\ell_i|} \| a \|_{1, s, \ell_i} 
\stackrel{|\ell_i| \leq |\ell|}{\lesssim} \langle \xi \rangle^{1 - |\ell_i|} \| a \|_{1, s, |\ell|}  
\stackrel{\eqref{propsimboloa}}{\lesssim_{s}} 
\langle \xi \rangle^{1 - |\ell_i|} \big(1 + \| \eta \|_{s + \sigma} \big)\,, 
\end{aligned}
\end{equation}
which, specialized for $s=s_0$, implies
\[
\| \partial_\xi^{\ell_i} a(\cdot, \xi) \|_{s_0}  
\lesssim 
\langle \xi \rangle^{1 - |\ell_i|} (1 + \| \eta \|_{s_0 + \sigma}) 
\stackrel{\| \eta \|_{s_0 + \sigma} \ll 1}{\lesssim} \langle \xi \rangle^{1 - |\ell_i|}\,. 
\] 
By using \eqref{caramella1} in  item $(i)$, and the smallness condition on $\eta$, one has that 
\begin{equation}\label{stima simbolo exp 4}
\begin{aligned}
 \| e^{y a(\cdot, \xi)} y^{m + n}\|_s  
 &\lesssim_{s,m,n} 
 \langle \xi \rangle^{- m - n} (1 + \| \eta \|_{s + \sigma})\,, 
 \quad\qquad \forall y \leq 0 \,, \qquad \forall \xi \in \R^d\,, 
 \\
  \| e^{y a(\cdot, \xi)} y^{m + n}\|_{s_0} 
  & \lesssim_{m,n}
  \langle \xi \rangle^{- m - n} (1 + \| \eta \|_{s_0 + \sigma}) 
{\lesssim} 
  \langle \xi \rangle^{- m - n}\,, \qquad \forall y \leq 0\, , \quad \forall \xi \in \R^d\,. 
\end{aligned}
\end{equation}
Thus, the estimates \eqref{stima simbolo exp 3}, \eqref{stima simbolo exp 4}, 
together with an iterative application of the interpolation estimate, 
imply that for any $1 \leq m \leq |\ell|$, $\ell_1 + \ldots + \ell_m = \ell$ with $|\ell|\leq \alpha$, 
\begin{equation}\label{stima simbolo exp 5}
\begin{aligned}
 \| e^{y a(\cdot, \xi)} y^{m + n} 
 (\partial_\xi^{\ell_1} a(\cdot, \xi)) \ldots (\partial_\xi^{\ell_m} a (\cdot, \xi)) \|_s  
 &\lesssim_{s,m,n}  
 \langle \xi \rangle^{- m - n } \prod_{i = 1}^m 
 \langle \xi \rangle^{1 - |\ell_i|} \big(1 + \| \eta \|_{s + \sigma} \big) 
 \\& 
 \stackrel{|\ell| = |\ell_1| + \ldots + |\ell_m|}{\lesssim_{s,m,n,\alpha}} 
 \langle \xi \rangle^{- n - |\ell|} (1 + \| \eta \|_{s + \sigma})\,,
\end{aligned}
\end{equation}
uniformly in $y \leq 0$ and $\xi \in \R^d$. 
The latter estimate, together with the formula \eqref{stima simbolo exp 3a}, 
implies that 
\[
\sup_{y \leq 0\,,\, \xi \in \R^d} \langle \xi \rangle^{n + |\ell|} \| \partial_\xi^\ell(e^{y a(\cdot, \xi)}) \|_{s}
\lesssim_{s,\alpha} 1 + \| \eta \|_{s + \sigma}\,,
\]
provided that $\| \eta \|_{s_0 + \sigma} $ is small enough and therefore \eqref{caramella2} follows.  
Estimate \eqref{caramella2Delta12} follows recalling formula \eqref{expaDelta} and reasoning as above.

\smallskip
\noindent
{\sc Proof of $(iii)$.} The arguments to prove the item $(iii)$ 
are similar to the ones used to prove $(ii)$
and recalling Remark \ref{rmk:simbolopiccolo} on the smallness of the symbols.
\end{proof}

\subsubsection{Parametrix for the parabolic evolution equation}\label{sec:parametrixEVO}
In this subsection we study some properties of the parabolic 
equation \eqref{calore striscia}
%\begin{equation}\label{equazione calore striscia}
%\partial_y w + {\rm Op}\big(a(\vphi, x, \xi) \big) w=0\,, \quad y \in [- 1, 0]\,.
%\end{equation}
where $a\in S^{1}$ is some 
symbol
satisfying  \eqref{propsimboloa}.
%We provide a parametrix for \eqref{calore striscia}.
\begin{prop}\label{proposizione totale parametrica calore}
Let $n_1, n_2 , \beta_0, k_0 \in \N$ and $N \geq n_1 + n_2 + \beta_0 + k_0$. 
There exist $\sigma \equiv \sigma(N,\beta_0, k_0)\geq \widetilde{\s}=\widetilde{\s}(N,\beta_0) \gg 0$ 
such that for any 
$\bar{s} > s_0 + \sigma$, there is $\delta \equiv \delta(\bar{s},N,\beta_0,k_0) \ll 1$ such that,
if $\| \eta \|_{s_0 + \sigma}^{k_0, \gamma} \leq \delta$, then
 the following holds.

\noindent
\emph{(i)}
the parabolic flow $\cU(-y)$, $y \in [- 1, 0]$ admits the expansion
\begin{equation}\label{formula parametrica proposizione generale}
\cU(\vphi, -y) = {\rm Op}\big( g(y, \vphi, x, \xi) \big) + \cH_N(\vphi, -y)
\end{equation}
where $g(y, \cdot) \in S^0$ and $\cH_N$ is a smoothing operator. 
More precisely, the symbol $g$ is of the form 
\begin{equation}\label{forma finale g}
\begin{aligned}
 g &= e^{y a} (1 + f ) \qquad \text{where} \qquad  
 f(y, \vphi, x, \xi) = \sum_{q = 2}^{2N} y^q f_q(\vphi, x, \xi), \quad f_q \in S^{q - 1}\,, 
\end{aligned}
\end{equation}
and the symbols $f_{q}$ satisfy the estimates
\begin{equation}\label{forma finale gbis}
\begin{aligned}
 \| f_q \|_{q - 1, s, \alpha}^{k_0, \gamma} &\lesssim_{ s, \alpha,N}  \| \eta \|_{s + \sigma}^{k_0, \gamma}\,, 
 %\quad \forall s \geq s_0\,,
 \\
 \|  \Delta_{12}f_q \|_{q-1,p,\alpha} &\lesssim_{ p, \alpha, N}
 \| \eta_1-\eta_2 \|_{p + \widetilde{\sigma}}\,,\qquad s_0\leq p\leq s_0+\s-\widetilde{\s}\,.
\end{aligned}
\end{equation}
The smoothing remainder $\cH_N(-y)$ satisfies the following property. For any $|\beta| \leq \beta_0$, the linear operator $\langle D \rangle^{n_1} \partial_\vphi^\beta \cH_N(\vphi, -y) \langle D \rangle^{n_2}$ is a $\cD^{k_0}$-tame operator with tame constant satisfying 
\begin{align}
\sup_{y \in [- 1, 0]} {\mathfrak M}_{\langle D \rangle^{n_1} 
\partial_\vphi^\beta \cH_N(-y) \langle D \rangle^{n_2}}(s) 
&\lesssim_{\bar{s}, N}  \| \eta \|_{s + \sigma}^{k_0, \gamma}\,, 
\qquad \forall s_0 \leq s \leq \bar{s}-\s\,,\label{stima resto parametrica caloreB}
\\
\| \partial_\vphi^\beta\Delta_{12}
\langle D \rangle^{n_1} 
 \cH_N(-y)
 \langle D \rangle^{n_2} \|_{\mathcal{L}({H}^{p};{H}^{p})}
&\lesssim_{\bar{s},N} \|\eta_{1}-\eta_{2}\|_{p+\s}\,,\;\;\;s_0\leq p\leq s_0+\s-\widetilde{\s}\,.
\label{stima resto parametrica caloreBDelta12}
\end{align}

\noindent
\emph{(ii)}
One has that  $g(y, \cdot) \in S^0$ for any $y \in [- 1, 0]$ 
and it satisfies the estimates, for $s_0\leq s\leq \bar{s}-\s$,
\begin{align}
\sup_{y \in [- 1, 0]} \| g(y, \cdot) \|_{0, s, \alpha}^{k_0, \gamma} 
&\lesssim_{s, \alpha, N} 
1 + \| \eta \|_{s + \sigma}^{k_0, \gamma}\,, 
\quad \sup_{y \in [- 1, 0]} \| g(y, \cdot) - e^{y a} \|_{- 1, s, \alpha}^{k_0, \gamma} 
\lesssim_{s, \alpha, N} 
\| \eta \|_{s + \sigma}^{k_0, \gamma}\,,\label{caramella10}
\\
 \|  \Delta_{12}g(y,\cdot) \|_{1,p,\alpha} 
 &\lesssim_{p, N}
 \| \eta_1-\eta_2 \|_{p + \widetilde{\sigma}}\,,\qquad 
 s_0\leq p\leq s_0+\s-\widetilde{\s}\,,\label{caramella10Delta12}
\end{align}
for any $\alpha \in \N$.
\end{prop}
The rest of the section is devoted to the proof of the latter proposition. 

%\begin{rmk}
%Notice that by formul\ae\,, \eqref{forma finale gbis}-\eqref{stima resto parametrica caloreB} the symbol $f$
%and the remainder $\mathcal{H}_{N}$ are \emph{small} in $\eta$.
%This is coherent with the fact that the equation \eqref{equazione calore striscia}
%(see also \eqref{calore striscia})
%reduces to 
%\[
%(\pa_{y}+\th|D|)w=0\,,
%\]
%whose solution is given by 
%\end{rmk}
\noindent
To shorten notations we often omit to write  the dependence on $\vphi$ from the symbols arising in our construction. 
First of all we note the following fact. If we look for a solution of 
\eqref{calore striscia} of the form
\[
\mathcal{U}(\tau)=\op(g(-\tau,x,\x))\,,\quad \tau\in[0,1]\,,
\]
for some symbol $g(y,x,\x)$, $y\in[-1,0]$, we have that, at the level of symbols, equation 
\eqref{calore striscia} reads
\begin{equation}\label{calore striscia simbolo}
(\pa_{y}g)=a(x,\x)\#g(y,x,\x)\,,\quad y\in[-1,0]\,, \qquad g(0,x,\x)=1\,. 
\end{equation}
For any integer $N \geq 1$, we construct an approximate solution of \eqref{calore striscia simbolo}
of the form 
\begin{equation}\label{propostaGn}
g(y, x, \xi) := \sum_{n = 0}^N g_{- n}( y,  x, \xi)\,, \qquad g_{- n}(y, \cdot) \in S^{-n}\,, \quad n = 0, \ldots, N\,.
\end{equation}
By the composition formula (recall \eqref{espansionecompostandard}), 
one has that 
\begin{equation}\label{comp Op a cal UN}
\begin{aligned}
a\#g&=\sum_{n=0}^{N}\sum_{q=0}^{N-n}a\#_{q} g_{-n}+r_{N}=\sum_{n=0}^{N}\big(a g_{-n}+r_{1-n}\big)+r_{N}\,,
\\
r_{1-n}&:=\sum_{\substack{k+q=n \\ 0\leq k\leq n-1 }}a\#_{q}g_{-k}\,,
\qquad
r_{N}:=\sum_{n=0}^{N}a\#_{\geq N-n+1}g_{-n}\,.
\end{aligned}
\end{equation}
Note that, since $a \in S^1$, $g_{- k} \in S^{- k}$, 
one has that 
\begin{equation*}%\label{grado forzante parametric heat}
\partial_\xi^\alpha a \partial_x^\alpha g_{- k} \in S^{1 - |\alpha| - k}\qquad  
\Longrightarrow \qquad 
r_{1 - n} \in S^{1 - n}\,,\quad {\rm and}\quad r_{N}\in S^{-N}\,.
\end{equation*}
Our aim is to find $g$ as in \eqref{propostaGn} in such a way that
\begin{equation}\label{parametrica calore approssimata}
\pa_{y}g=a\# g-r_{N}\,,\qquad g(0,x,\x)=1\,.
\end{equation}
%
%The ``approximate parabolic flow" $\cU_N(y)$ must solves the problem
%\begin{equation}\label{parametrica calore approssimata}
%\begin{cases}
%\partial_y \cU_N(y) = {\rm Op}(a) \cU_N(y) + \cR_N(y)\,, \\
%\cU_N(0) = {\rm Id}\,. 
%\end{cases}
%\end{equation}
Hence the iterative equations for the symbols $g_{- n}$, $n = 0, \ldots, N$ are given by 
\begin{equation}\label{eq g0 parametrica calore}
\begin{cases}
\partial_y g_0(y, x, \xi) = a(  x, \xi) g_0(y,  x, \xi)\,, \\
g(0,  x, \xi) = 1\,,
\end{cases}
\end{equation}
and for $n = 1, \ldots, N$
\begin{equation}\label{eq g - n parametrica calore}
\begin{cases}
\partial_y g_{- n}(y, x, \xi) = a( x, \xi) g_{- n}(y,  x, \xi) + r_{1 - n}(y,  x, \xi)\,, \\
g_{- n}(0,  x, \xi) = 0\,.
\end{cases}
\end{equation}
Note that the symbol $r_{1 - n}$ depends only  on $g_0, g_{- 1}, \ldots, g_{1 - n}$ 
hence the latter equations can be solved iteratively. 

\medskip

\noindent
{\sc Solving \eqref{eq g0 parametrica calore}.} By an explicit computation one has that 
\begin{equation}\label{espressione g0 parametrica calore}
g_0 = e^{y a}, \quad y \leq 0
\end{equation}
and hence by Lemma \ref{prop exp simbolo}, one has that $g_0 \in S^0$. 

\medskip

\noindent
{\sc Soving \eqref{eq g - n parametrica calore}.} 
Now assume to have determined $g_{- 1} \in S^{- 1}, g_{- 2} \in S^{- 2}, \ldots, g_{1 - n} \in S^{1 - n}$. 
Let us determine $g_{- n} \in S^{- n}$. 
We make the following ansatz. 
For any $ k = 0, \ldots , n-1$ %there exists $Q_k \gg 0$ such that 
the symbol 
$g_{- k}$ 
has the form 
\begin{equation}\label{forma induttiva simboli g - k}
g_{- k}(y, \vphi, x, \xi) = \sum_{q = 2}^{2k} y^q e^{y a(\vphi, x, \xi)} f_q^{(k)}(\vphi, x, \xi)\,, 
\qquad f_q^{(k)} \in S^{q - k}, \quad q = 2, \ldots, 2k\,.
\end{equation}
Moreover, we assume that there exist 
$\sigma_{k} \geq \widetilde{\s}_{k}\gg 0$ such that if 
$\| \eta \|_{s_0 + \sigma_k} \ll  1$ (the smallness is  w.r.t. constants depending on $\s$ ), 
then for any $k = 0, \ldots, n-1$, $q = 2, \ldots, 2k$ 
\begin{align}
\| f_q^{(k)} \|_{q - k, s, \alpha}^{k_0, \gamma} &\lesssim_{s, k, \alpha,q} 
 \| \eta \|_{s + \sigma_k}^{k_0, \gamma}\,, %\qquad \forall s \geq s_0\,,
 \label{stime coefficienti g - k}
 \\
 \|  \Delta_{12}f_q^{(k)} \|_{q-k,p,\alpha} &\lesssim_{p,q,\alpha}
 \| \eta_1-\eta_2 \|_{p + \widetilde{\sigma}_{k}}\,,\qquad s_0\leq p\leq s_0+\s_{k}-\widetilde{\s}_{k}\,.
 \label{stime coefficienti g - kDelta12}
\end{align}
We want to show that also $g_{- n}$ has the same form with coefficients satisfying similar estimates. 
\begin{lemma}\label{lemma r 1 - n parametrica calore}
For $n=0,\ldots, N$ there exist $\sigma_n\geq \widetilde{\s}_{m} \gg 0$ 
large enough such that if $\| \eta \|_{s_0 + \sigma_n} \lesssim_{s} 1$, the following holds. 
%There is $Q_n \gg 0$ large enough such that t
The symbol 
$r_{1 - n} \in S^{1 - n}$ has the form 
\[
r_{1 - n} = \sum_{q = 1}^{2n-1} y^q e^{y a} v_q^{(n)}\,, \qquad  v_q^{(n)}=v_q^{(n)}(x,\x) \in S^{1 + q - n}\,.
\]
Moreover, for any $q = 1, \ldots, 2n-1$ one has that 
\begin{equation}\label{succodiarancia1}
\begin{aligned}
\| v_q^{(n)} \|_{1 + q - n, s, \alpha}^{k_0,\gamma}
&\lesssim_{n, s, \alpha,q} 
 \| \eta \|_{s + \sigma_n}^{k_0,\gamma}\,, 
 %\qquad \forall s \geq s_0\,,
\\
\| \Delta_{12} v^{(n)}_q\|_{1+q-n, p, \alpha} 
 &\lesssim_{n,p,\alpha},q\| \eta_1-\eta_{2} \|_{p + \widetilde{\sigma}_n}\,, 
 \;\;\;s_0\leq p\leq s_0+\s_{n}-\widetilde{\s}_{n}\,.
\end{aligned}
\end{equation}
\end{lemma}
\begin{proof}
By recalling the formula \eqref{comp Op a cal UN}, the induction hypothesis 
\eqref{forma induttiva simboli g - k} 
we shall write
\[
\begin{aligned}
r_{1-n}&=\sum_{\substack{k+q=n \\ 0\leq k\leq n-1 }}\sum_{j=2}^{2k}a\#_q 
\Big( y^j e^{y a} f_j^{(k)}\Big)
\\&
\stackrel{\eqref{espansionecompostandard}}{=}
\sum_{\substack{k+q=n \\ 0\leq k\leq n-1 }}\sum_{j=2}^{2k}\sum_{|\alpha|=q}
\frac{1}{\ii^{|\alpha|}\alpha!}(\pa_{\x}^\alpha a)\pa_{x}^{\alpha}\Big( y^j e^{y a} f_j^{(k)}\Big)\,.
\end{aligned}
\]
Therefore, the bounds \eqref{succodiarancia1}
 follow by the product rule, Lemma \ref{prop exp simbolo}-$(iii)$, 
the tame estimates for the product of symbols and the induction hypotheses 
\eqref{forma induttiva simboli g - k}, \eqref{stime coefficienti g - k}.  
\end{proof}
We now can solve the equation \eqref{eq g - n parametrica calore} 
and analyze the structure of the symbol $g_{- n}$. 
By variation of constants one has that 
\[
g_{- n} = e^{y a} \gamma_{- n}\,,
\]
where the symbol $\gamma_{- n}$ solves the problem
\[
\begin{cases}
\partial_y \gamma_{- n}(y, x, \xi) = e^{- y a(y, x, \xi)} r_{1 - n}(y, x, \xi) \\
\gamma_{- n}(0, x, \xi) = 0\,.
\end{cases}
\]
By Lemma \ref{lemma r 1 - n parametrica calore}, 
one obtains that 
\[
e^{- y a(y, x, \xi)} r_{1 - n}(y, x, \xi) =\sum_{q = 1}^{2n-1} y^q  v_q^{(n)}\,, 
\qquad \text{where} \qquad v_q^{(n)} \in S^{1 + q - n}\,,
\]
implying that 
\[
\gamma_{- n}(y, x, \xi) = \sum_{q = 1}^{2n-1} \dfrac{y^{q + 1}}{q + 1}  v_q^{(n)}\,, 
\qquad v_q^{(n)} \in S^{1 + q - n}\,.
\]
Therefore, we conclude that
\[
g_{- n} 
= \sum_{q = 1}^{2n-1} e^{y a} \dfrac{y^{q + 1}}{q + 1}  v_q^{(n)} 
= \sum_{k = 2}^{2n} y^k e^{y a} f_k^{(n)}\,,
\]
where 
\[
f_k^{(n)} := \frac{v_{k - 1}^{(n)}}{k} \in S^{k - n}, \quad k = 2, \ldots, 2n\,. 
\]
This means (using again Lemma \ref{prop exp simbolo}) that
$g_{- n} \in S^{- n}$ has the form \eqref{forma induttiva simboli g - k}, while 
the estimates \eqref{stime coefficienti g - k} at the step $n$ follows by the estimates of Lemma \ref{lemma r 1 - n parametrica calore}. 

\noindent
From the properties \eqref{forma induttiva simboli g - k}, \eqref{stime coefficienti g - k} 
and \eqref{stime coefficienti g - kDelta12} for 
$k = 1, \ldots, N$, one can easily deduce that the symbol 
$g = g_0 + \ldots + g_{- N}$ satisfies \eqref{forma finale g} and the bounds \eqref{forma finale gbis}.
The properties \eqref{forma finale g}-\eqref{forma finale gbis},  together with Lemma \ref{prop exp simbolo} and the product estimates of symbols (see Lemma \ref{lemmacomposizioneSTANDARD})
imply the bounds \eqref{caramella10}-\eqref{caramella10Delta12}. 
This proves item $(ii)$ of the proposition. 
Finally, the bounds \eqref{caramella10}-\eqref{caramella10Delta12}, together with \eqref{propsimboloa}  
and the composition Lemma \ref{lemmacomposizioneSTANDARD}
imply that the remainder $\cR_N(y) = {\rm Op}(r_N(y, x, \xi)) \in OPS^{- N}$ 
in \eqref{comp Op a cal UN} satisfies the estimate
\begin{equation}\label{stima resto parametrica calore}
\begin{aligned}
\| r_N \|_{- N, s, \alpha}& \lesssim_{N, s, \alpha} \| \eta \|_{s + \sigma}^{k_0,\gamma}\,, 
\quad s \geq s_0\,,
\\
\| \Delta_{12}r_N \|_{- N, p, \alpha} &\lesssim_{N, p, \alpha}  \| \eta_1-\eta_2 \|_{p + \widetilde{\sigma}}\,, 
%\quad \forall p \geq s_0\,.
\end{aligned}
\end{equation}
for some $\widetilde{\s}\leq \s$ as long as $p+\widetilde{\s}\leq s_0+\s$.
We now set $\mathcal{U}_{N}(\tau):=\op(g(-\tau;x,\x))$, $\tau\in[0,1]$ with $g$ 
solving \eqref{parametrica calore approssimata}. We hence note that 
\[
\pa_{\tau}\mathcal{U}_{N}(\tau)=-\op(a)\mathcal{U}_{N}(\tau)-\mathcal{R}_{N}(-\tau)\,,\qquad \mathcal{U}_{N}(0)=1\,,
\]
for $\tau\in[0,1]$.
We are now in position to  
estimate the difference between $\cU(\tau)$ in \eqref{calore striscia} and $\cU_N(\tau)$, 
see \eqref{calore striscia}, \eqref{parametrica calore approssimata}. 
We write an equation for $\cH_N(\tau) := \cU(\tau) - \cU_N(\tau)$. 
Clearly one has 
\begin{equation}\label{parametrica vs parametrica approssimata calore}
\begin{cases}
\partial_y \cH_N(\tau) =- {\rm Op}(a) \cH_N(\tau) + \cR_N(-\tau) \\
\cH_N(0) = 0
\end{cases}
\end{equation}
and by Duhamel formula one obtains that 
\begin{equation}\label{lemma resto duhamel parametrica calore}
\cH_N(\tau) =  \int_0^\tau \cU(\tau - z) \cR_N(-z)\, d z\,, \quad \tau \in [0,1]\,. 
\end{equation}
By Proposition \ref{stima derivate lambda vphi flusso parabolico}, 
by the estimate \eqref{stima resto parametrica calore} 
and by applying Lemmata \ref{constantitamesimbolo} and \ref{composizione operatori tame AB},
one can deduce the estimates 
\eqref{stima resto parametrica caloreB}-\eqref{stima resto parametrica caloreBDelta12}
for the remainder $\cH_N(\tau)$. 
The proof of the Proposition \ref{proposizione totale parametrica calore} is then concluded.

\subsubsection{Analysis of the leading terms}\label{sec:analisiW}

In this section we analyse the structure of the term $\cW(\vphi)$ 
in \eqref{def parte dominante parte trascurata w x 0}, by exploiting the result on the flow 
$\mathcal{U}(\vphi,y)$ given in subsection \ref{sec:parametrixEVO}.
In the following, to simplify the notation we shall omit the dependence on $\vphi$ 
(and on $(y;x,\x)$)
in the symbols and operators.
More precisely, we prove the following Lemma.

\begin{lemma}\label{lem:espandoww}
Under the assumptions of Proposition \ref{lemma espansione pseudo-diff w x 0} the following holds.
One has that the operator $ \cW(\vphi)$ in \eqref{def parte dominante parte trascurata w x 0} admits the expansion
\begin{equation}\label{formulaespdiW}
 \cW(\vphi)[\psi] =\th|D|\big(\tanh(\th|D|)-1\big)[\psi] + \widetilde{\cQ}(\vphi)[\psi]
\end{equation}
where $\widetilde{\mathcal{Q}}$ is a smoothing remainder satisfying the estimates 
\eqref{stimefinalew00}.
\end{lemma}
\begin{proof}
In view of the 
expansion of $\mathcal{U}(\vphi,y)$ in
\eqref{formula parametrica proposizione generale}-\eqref{forma finale g}, we shall write, for some $N\gg1$,
\begin{equation}\label{espansioneNuova}
\begin{aligned}
\cU(\vphi, -y)&={\rm Op}\big( g(y, \vphi, x, \xi) \big) + \cH_N(\vphi, y)
\\&=\op(e^{y a})+\op(e^{y a}f)+\cH_N(\vphi, y)
\\&=\op(e^{y\th|\x|\chi(\x)})+\op\big(\mathtt{r}_1+e^{ya}f\big)+\cH_N(\vphi, y)\,,
\end{aligned}
\end{equation}
where, recalling the expansion of $a$ in \eqref{propsimboloa}, we defined 
\begin{equation}\label{simboRR1}
\mathtt{r}_1:=\mathtt{r}_1(y;\vphi,x,\x):=\int_{0}^{1}yr e^{y\th|\x|\chi(\x)+\tau yr}d\tau\,.
\end{equation}
The parameter $N$ in Proposition \ref{lemma espansione pseudo-diff w x 0} 
will be chosen later large with $M$.
We now define
\begin{equation}\label{definizionerestoQN}
\begin{aligned}
\widetilde{\mathcal{Q}}(\vphi)[\psi]&=\mathcal{Q}_{N}^{(1)}(\vphi)[\psi]+\mathcal{Q}_{N}^{(2)}(\vphi)[\psi]\,,
\\
\mathcal{Q}_{N}^{(1)}(\vphi)[\psi]&:=
 \int_{- 1}^0  \op\big(\mathtt{r}_1+e^{y a}f\big) \chi'(y) \mathfrak{Z}(y)[\psi] \, d y\,,
 \\
\mathcal{Q}_{N}^{(2)}(\vphi)[\psi]&:=
 \int_{- 1}^0 
  \cH_N(\vphi, y) \chi'(y) \mathfrak{Z}(y)[\psi] \, d y\,,
 \end{aligned}
\end{equation}
so that,  by \eqref{def parte dominante parte trascurata w x 0} and the expansion \eqref{espansioneNuova}, 
we get
\begin{equation*}
\begin{aligned}
 \cW(\vphi)[\psi] &  = \int_{- 1}^0 \cU(\vphi,- y) \chi'(y) \mathfrak{Z}(y)[\psi] \, d y
=
\underbrace{ \int_{- 1}^0 \op(e^{y\th|\x|\chi(\x)}) \chi'(y) \mathfrak{Z}(y)[\psi] \, d y}_{=: I}+
 \mathcal{Q}_{N}(\vphi)[\psi]\,.
 \end{aligned}
 \end{equation*}
 We now give an explicit expression for the term $I$.
 By formul\ae\, \eqref{eq:opZZ}-\eqref{eq:opZZ2} we deduce that
 \[
 \begin{aligned}
 I=\int_{- 1}^0 \op(e^{y\th|\x|\chi(\x)}) \chi'(y) \mathfrak{Z}(y)[\psi] \, d y
 &=\int_{- 1}^0 e^{y\th|D|} \chi'(y) e^{-y\th|D|} \th|D|\big(\tanh(\th|D|)-1\big)[\psi] \, d y
 \\&
 =\th|D|\big(\tanh(\th|D|)-1\big)[\psi]\int_{-1}^{0}\chi'(y)\,dy
 \\&
 =\th|D|\big(\tanh(\th|D|)-1\big)[\psi]\,.
 \end{aligned}
 \]
 Therefore, formula \eqref{formulaespdiW} follows with the remainder $\mathcal{Q}_{N}$ defined in 
 \eqref{definizionerestoQN}.
 It remains to prove the bounds \eqref{stimefinalew00} on $\mathcal{Q}_{N}$.

  \smallskip
 \noindent
 \emph{Estimates on $\mathcal{Q}_{N}^{(1)}$.}  
First of all recalling \eqref{eq:opZZ2}, we note that
 for any $\x\in \Gamma^{*}$, we have
 \begin{equation}\label{simbolomeno infinito}
 \tanh(\th|\x|\chi(\x))-1=-2\frac{e^{-2\th|\x|\chi(\x)}}{1+e^{-2\th|\x|\chi(\x)}}\in S^{-\infty}\,.
 %\qquad \Rightarrow\qquad 
 %\th|\x|( \tanh(\th|\x|)-1)\in S^{-\infty}\,.
 \end{equation}
 Therefore, we shall write
  \begin{equation}\label{simbolomeno infinito2}
 \mathfrak{Z}(y)[\cdot]=\op( \mathtt{c}(y,\x))\,,\quad
 \mathtt{c}(y,\x):=e^{-2\th|\x|\chi(\x)-y\th|\x|\chi(\x)} \frac{(-2)\th|\x|\chi(\x)}{1+e^{-2\th|\x|\chi(\x)}}
 \quad \forall\,y\in[-1,0]\,.
 \end{equation}
 Recalling \eqref{propsimboloa} and \eqref{definizionerestoQN}  we have
 \[
 \op\big(\mathtt{r}_1+e^{y a}f\big)  \mathfrak{Z}(y)=\op(\mathtt{d}(y,x,\x))\,,
 \]
 where
 \[
 \begin{aligned}
 \mathtt{d}(y,x,\x)&=\mathtt{r_1}\mathtt{c}+e^{ya}f\mathtt{c}
 \\&=\Big(
 \int_{0}^{1}yr e^{-2\th|\x|\chi(\x)+\tau yr}d\tau+e^{-2\th |\x|\chi(\x)+yr}f\Big)\frac{(-2)\th|\x|\chi(\x)}{1+e^{-2\th|\x|\chi(\x)}}\,.
 \end{aligned}
 \]
 Notice that
 \[
| {\rm Re}(2\th |\x|\chi(\x)-yr)|\geq \kappa|\x|\,, \quad |\xi| \geq 1
 \]
 for some small $\kappa>0$. Then one can reason exactly as in the proof of Lemma \ref{prop exp simbolo}.
 Therefore, using the estimates on the action of pseudo-differential operator, 
 to obtain that $\mathtt{d}\in S^{-N}$, for any $N\gg1$, and 
 \begin{align}
\| \mathtt{d} \|^{k_0,\gamma}_{- N , s, \alpha} 
&\lesssim_{N, s, \alpha} \| \eta \|^{k_0,\gamma}_{s + \sigma}\,, 
\qquad s \geq s_0\,,\label{caramella200}
\\
 \sup_{y\leq 0}\| \Delta_{12} \mathtt{d}\|_{-N, p, \alpha} 
 &\lesssim_{N, p,\alpha}
 \| \eta_1-\eta_{2} \|_{p + \widetilde{\sigma}}\,, \;\;\;s_0\leq p\leq s_0+\s-\widetilde{\s}\,,
 \label{caramella200Delta12}
\end{align}
for some large enough $\s\geq \widetilde{\s}\gg1$ (depending on $N$).
Then the  estimates \eqref{stimefinalew00} on $\mathcal{Q}_{N}^{(1)}$ in 
\eqref{definizionerestoQN} follow combining
\eqref{caramella200}-\eqref{caramella200Delta12} with Lemma \ref{constantitamesimbolo}
and taking $N=M-\mathtt{c}$ and $m_1=-M$, $m_2=\mathtt{c}$.
 %large enough with respect to $M$ and $\mathtt{c}$.

\smallskip
 \noindent
 \emph{Estimates on $\mathcal{Q}_{N}^{(2)}$.} 
 First of all (recall \eqref{eq:opZZ2}) by \eqref{simbolomeno infinito2}  and an explicit computation
 we deduce that $\mathtt{c}\in S^{1}$ uniformly in $y\in[-1,0]$.
 Therefore, the estimates \eqref{stimefinalew00} on $\mathcal{Q}_{N}^{(2)}$ follow by using
 \eqref{stima resto parametrica caloreB}-\eqref{stima resto parametrica caloreBDelta12},
 taking $n_1=M$, $n_2=0$, $N=M+\beta_0+k_0$.
 %and  Lemma \ref{constantitamesimbolo}
% , by taking $N$ large with $M$.
 This concludes the proof.
 \end{proof}

\subsubsection{Analysis of the remainders
%$\cS(\vphi)$ in \eqref{def parte dominante parte trascurata w x 0}.
}
\label{sec:analisiS}
In this section we analyze the structure of the term $\cS(\vphi)$ 
in \eqref{def parte dominante parte trascurata w x 0}, by exploiting the result on the flow 
$\mathcal{U}(\vphi,y)$ given in subsection \ref{sec:parametrixEVO}.

We recall that by the properties stated in Lemma \ref{lemma tame F Gamma vphi y che palle}, one has that $\Gamma$ vanishes identically near $y = 0$, hence there exists $- 1 < \bar y < 0$ such that $\Gamma(\vphi, y) = 0$ for any $y \in [\bar y, 0]$.
Therefore, $\cS(\vphi)$ in  \eqref{def parte dominante parte trascurata w x 0}
takes the form
\[
\cS(\vphi)[\psi] = \int_{- 1}^{\bar{y}} \cU(\vphi,- y) 
\Gamma(\vphi, y)[\psi] dy +  \int_{- 1}^{0}\cU(\vphi,- y) S_{\mathcal{L}} (\vphi, y)[\psi]\, d y\,,
\]
where $\Gamma(\vphi,y), S_{\mathcal{L}} (\vphi, y)$ are given in \eqref{restiGammaSelle}.
Moreover, 
according to the expansion of $\mathcal{U}(\vphi,-y)$ in 
Proposition \ref{proposizione totale parametrica calore} 
(see \eqref{formula parametrica proposizione generale}-\eqref{forma finale g})
one can expand $\cS(\vphi)$ in the following way: 
\begin{equation}\label{formulaespdiSS}
\cS(\vphi)  = \sum_{q = 1}^{2N} \cS_q(\vphi) + \cS_N(\vphi)+\mathcal{S}_{\mathcal{L}}(\vphi)
\end{equation}
where
\begin{align}
\cS_1(\vphi) & :=  \int_{- 1}^{\bar y} {\rm Op}\big(e^{y a(\vphi, x, \xi)} \big) \circ \Gamma(\vphi, y)\, d y\,, 
\label{opS1}
\\
\cS_q(\vphi) & := 
\int_{- 1}^{\bar y} {\rm Op}\big(e^{y a(\vphi, x, \xi)} y^q f_q(\vphi, x, \xi) \big) \circ \Gamma(\vphi, y)\, d y\,, 
\qquad q = 2, \ldots, 2N\,, 
\label{opSq}
\\
\cS_N(\vphi) & := \int_{- 1}^{\bar y} \cH_N(\vphi, -y)\circ \Gamma(\vphi, y)\, d y\,,
\label{opSN}
\\
\cS_{\mathcal{L}}(\vphi) &:=\int_{- 1}^{0}\cU(\vphi,- y) S_{\mathcal{L}} (\vphi, y)\, d y\,,
\label{opSL}
\end{align}
with $N:=M+\beta_0+k_0+1$.
%The parameter $N$ will be fixed later large with respect to $M$.
In the sequel, we shall make a repeated use of the following elementary fact. 
By the Cauchy Schwartz inequality if $f \in L^2([- 1, 0], H^s_{\vphi, x})$, 
one has that for any $\alpha \geq 0$
\begin{equation}\label{stima con cauchy schwatz integrale striscia}
\Big\| \int_{- 1}^{\bar y} \frac{1}{y^\alpha} f(y)\, d y\Big\|_s 
\lesssim \| f \|_{L^2_y H^s_{\vphi, x}}\,. 
\end{equation}

We prove the following.
\begin{lemma}\label{lemma stima resto cal S DN}
Fix any $\beta_0,k_0\in\N$ and  $M\geq\mathtt{c} \geq \beta_0+k_0+2$.
%Let $M \gg 0$, $\beta_0, k_0 > 0$.
Then there exists $\sigma \equiv \sigma(\beta_0, k_0, M) \gg 0$ large enough,
such that for any $\bar{s}>s_0+\s$ there is 
 $\delta \equiv \delta(\bar{s},\beta_0,k_0,M) \ll 1$ small enough such that if 
\begin{equation}\label{smallnessSS}
\| \eta \|_{s_0 + \sigma}^{k_0, \gamma} \leq \delta\,,
\end{equation} 
then for any $\beta \in \N^\nu$, $|\beta| \leq \beta_0$, 
the linear operator 
$\langle D \rangle^M \partial_\vphi^\beta \cS(\vphi) \langle D \rangle^{- \tc}$ 
(see \eqref{formulaespdiSS})
is $\cD^{k_0}$-tame operator with tame constant satisfying the bounds
\begin{equation}\label{stimacalSS}
\begin{aligned}
{\mathfrak M}_{\langle D \rangle^M \partial_\vphi^\beta \cS(\vphi) 
\langle D \rangle^{- \tc}}(s) 
&\lesssim_{\bar{s}, M}
 \| \eta \|_{s + \sigma}^{k_0, \gamma}\,, \qquad \;s_0 \leq s \leq \bar{s}-\s\,,
\\
\| \partial_\vphi^\beta\Delta_{12}
\langle D \rangle^{M} 
\cS(\vphi)
 \langle D \rangle^{-\tc} \|_{\mathcal{L}({H}^{p};{H}^{p})}
&\lesssim_{\bar{s},M}
 \|\eta_{1}-\eta_{2}\|_{p+\widetilde{\s}}\,,\;\;\;s_0\leq p\leq s_0+\s-\widetilde{\s}\,,
\end{aligned}
\end{equation}
for some $\widetilde{\s}\leq \s$.
\end{lemma}

\begin{proof}
We consider  separately each operator  $\cS_{\mathcal{L}}$, 
$\cS_1$, $\cS_p$, $p = 2, \ldots, 2N$ and $\cS_N$ 
appearing in \eqref{opS1}-\eqref{opSL}, and we show that they satisfy 
the bounds \eqref{stimacalSS}.

\medskip

\noindent
{\bf Estimate of  \eqref{opSL}.}
We recall that (see \eqref{restiGammaSelle}) has the form
\[
S_{\mathcal{L}} (\vphi, y)[\cdot] := \chi(y) \mathcal{R}_{\mathcal{L}}\circ\big(\mathcal{L}_0
+\cF_0(\vphi, y)\big)[\cdot] 
\]
where $\mathcal{R}_{\mathcal{L}}$ is the smoothing 
operator given by Lemma  \ref{divo1}, $\mathcal{L}_0$ is in \eqref{def vphi 0} and 
$\mathcal{F}_0$ is the operator appearing in \eqref{soluzioneV}.
Then it is easy ti check that  the operator  $\mathcal{S}_{\mathcal{L}}$ in \eqref{opSL}
satisfies estimates like \eqref{stimacalSS} combining 
Lemma \ref{stima derivate lambda vphi flusso parabolico} (to estimate $\mathcal{U}(\vphi,-y)$)
with Lemma  \ref{divo1}, Lemma \ref{lemma stima cal L0 sol omogenea laplace} and Proposition
\ref{lemma tame laplace partial y u langle D rangle u}.

\medskip

\noindent
{\bf { Estimate of  \eqref{opS1}.}}
Let $M \gg 0$ large enough. We rewrite \eqref{opS1} as
\begin{equation*}
\cS_1(\vphi) = \int_{- 1}^{\bar y} \frac{1}{y^M}  {\rm Op}\big(y^M e^{y a(\vphi, x, \xi)} \big) 
\circ \Gamma(\vphi, y)\, d y\,.
\end{equation*}
In view of the smallness condition  \eqref{smallnessSS}, we shall apply 
 Lemma \ref{prop exp simbolo}-$(ii)$, which implies the bounds
\begin{equation*}
\begin{aligned}
\| y^M e^{y a} \|_{- M, s, \alpha}^{k_0, \gamma} 
&\lesssim_{M, s, \alpha} 
1 + \| \eta \|_{s + \sigma}^{k_0, \gamma}\,, \quad \forall s_0\leq s \leq \bar{s}-\s\,,
\\
 \sup_{y\leq 0}\| \Delta_{12} y^M  e^{y a}\|_{-M, p, \alpha} 
 &\lesssim_{M, p,\alpha}\| \eta_1-\eta_{2} \|_{p + \widetilde{\sigma}}\,, \;\;\;s_0\leq p\leq s_0+\s-\widetilde{\s}\,,
 \end{aligned}
\end{equation*}
for some $\widetilde{\s}\leq \s$.
In particular, for $|k_1|\leq k_0$, $|\beta_1|\leq\beta_0$, using the estimate above,
Lemma \ref{lemmacomposizioneSTANDARD} and again the smallness condition on $\eta$ one gets
\begin{equation}\label{stima y M e y a stima cal S 1}
\begin{aligned}
\sup_{y \in [- 1, 0]} \Big\| \langle D \rangle^M \circ 
{\rm Op}\big(y^M \partial_\lambda^{k_1} \partial_\vphi^{\beta_1}e^{y a(\vphi, x, \xi)} \big) \Big\|_{0, s, 0} 
& \lesssim_{s, M}
 \gamma^{- |k_1|} \big( 1 + \| \eta \|_{s + \sigma}^{k_0, \gamma}\big)\,,
%\quad \forall s \geq s_0\,.
\end{aligned}
\end{equation}
for some $\s \equiv \sigma(\beta_0, k_0, M)\gg1$ and $s_0\leq s\leq \bar{s}-\s$.
Now, taking $k \in \N^{\nu + 1}$, $\beta \in \N^{\nu}$, $|\beta| \leq \beta_0$, 
and recalling $\tc\geq \beta_0+k_0+2$,
we have
\begin{equation}\label{D M - beta - k0 - 2 cal S}
\begin{aligned}
& \langle D \rangle^M \partial_\lambda^k \partial_\vphi^\beta \cS_1(\vphi) 
\langle D \rangle^{- \tc}  = \sum_{\begin{subarray}{c}
\beta_1 + \beta_2 = \beta \\
k_1 + k_2 = k
\end{subarray}} C(\beta_1, \beta_2, k_1, k_2) \cR_{\beta_1, \beta_2, k_1, k_2}(\vphi) \,, \\
& \cR(\vphi)  \equiv \cR_{\beta_1, \beta_2, k_1, k_2}(\vphi) \\
&  :=  \int_{- 1}^{\bar y} \frac{1}{y^M}  \langle D \rangle^M \circ {\rm Op}\Big(y^M \partial_\lambda^{k_1} \partial_\vphi^{\beta_1}e^{y a(\vphi, x, \xi)} \Big) \circ \partial_\lambda^{k_2} \partial_\vphi^{\beta_2}\Gamma(\vphi, y)\, d y \circ \langle D \rangle^{- \tc}\,,
\end{aligned}
\end{equation}
and hence, to get the result,  for any $k_1, k_2, \beta_1, \beta_2$, $k_1 + k_2 = k$, 
$\beta_1 + \beta_2 = \beta$ it suffices to estimate the linear operator $\cR(\vphi)$. 
By the inequality \eqref{stima con cauchy schwatz integrale striscia}
and Lemma \ref{lemma: action Sobolev},
one has that 
\begin{equation*}
\begin{aligned}
\| \cR \psi \|_s 
& \lesssim 
\Big\| \langle D \rangle^M \circ 
{\rm Op}\big(y^M \partial_\lambda^{k_1} \partial_\vphi^{\beta_1}e^{y a(\vphi, x, \xi)} \big) 
\circ \partial_\lambda^{k_2} \partial_\vphi^{\beta_2}\Gamma(\vphi, y)\,  
[ \langle D \rangle^{-\tc} \psi] \Big\|_{L^2_y H^s_{\vphi, x}} 
\\& 
\stackrel{\eqref{stima y M e y a stima cal S 1}}{\lesssim_{s, M}} 
\gamma^{- |k_1|}   \| \partial_\lambda^{k_2} \partial_\vphi^{\beta_2}\Gamma(\vphi, y)\,  
[ \langle D \rangle^{- |\beta| - k_0 - 2} \psi] \|_{L^2_y H^s_{\vphi, x}} 
\\&\qquad+\gamma^{-|k_1|} \| \eta \|_{s + \sigma}^{k_0, \gamma}
 \| \partial_\lambda^{k_2} \partial_\vphi^{\beta_2}\Gamma(\vphi, y)\, 
  [ \langle D \rangle^{-\tc} \psi] \|_{L^2_y H^{s_0}_{\vphi, x}}  
\\& 
%\lesssim_s 
%\gamma^{- |k_1|}  \| \partial_\lambda^{k_2} \partial_\vphi^{\beta_2}\Gamma(\vphi, y)\,  
%[ \langle D \rangle^{- |\beta| - k_0 - 2} \psi] \|_{\cO^s} 
%\\&\qquad+ \gamma^{- |k_1|}\| \eta \|_{s + \sigma}^{k_0, \gamma} \| \partial_\lambda^{k_2} \partial_\vphi^{\beta_2}\Gamma(\vphi, y)\,  
%[ \langle D \rangle^{-\tc} \psi] \|_{\cO^{s_0}} 
%\\& 
\stackrel{\eqref{prop L2 vhi Hs Hs vphi cal 2 xy},{\rm Lemma}\,\ref{lemma tame F Gamma vphi y che palle}}{\lesssim_{s, M}} 
\gamma^{- |k|} \Big( \| \eta \|_{s_0 + \sigma}^{k_0, \gamma}  \| \psi \|_s 
+ \| \eta \|_{s + \sigma}^{k_0, \gamma} \| \psi \|_{s_0}\Big)\,. 
\end{aligned}
\end{equation*}
This latter estimate  implies the first bound in \eqref{stimacalSS} for the operator $\mathcal{S}_1$. 
The Lipschitz bound follows similarly.

\medskip

\noindent
{\bf Estimate of \eqref{opSq}.}
Let $M \gg 0$ large enough. By recalling the definition of $\cS_q$ in \eqref{opSq}, we write 
\begin{equation}\label{prima formula cal Sp vphi}
\cS_q(\vphi) = \int_{- 1}^{\bar y} \frac{1}{y^M}  
{\rm Op}\big(y^{M + q}e^{y a(\vphi, x, \xi)} f_q(\vphi, x, \xi) \big) \circ \Gamma(\vphi, y)\, d y\,.
\end{equation}
Using the smallness condition \eqref{smallnessSS} on $\eta$ one can 
apply Lemma \ref{prop exp simbolo}-$(ii)$. Then by  estimates \eqref{forma finale g}-\eqref{forma finale gbis} 
on 
$f_q$ and using Lemma \ref{lemmacomposizioneSTANDARD},
%one has that if $\| \eta \|_{s_0 + \sigma}^{k_0,\gamma} \ll 1$ ($\sigma \gg 0$ large enough), 
one has that 
\begin{equation*}%\label{stima y M e y a stima cal Sp}
\begin{aligned}
\| y^{M + q} e^{y a} f_q \|_{- M- 1 , s, \alpha}^{k_0, \gamma} 
&\lesssim_{s, \alpha, M}  \| \eta \|_{s + \sigma}^{k_0, \gamma}\,, \qquad \forall s_0\leq s \leq \bar{s}-\s\,,
 \\
 \|  \Delta_{12}y^{M + q} e^{y a}f_q \|_{-M-1,p,\alpha} &\lesssim_{p,M,\alpha}
 \| \eta_1-\eta_2 \|_{p + \widetilde{\sigma}}\,,\qquad s_0\leq p\leq s_0+\s-\widetilde{\s}\,,
 \end{aligned}
\end{equation*}
for some $\widetilde{\s}\leq \s$.
In particular, for $|k_1|\leq k_0$, $|\beta_1|\leq\beta_0$, using the estimate above,
%By the estimate \eqref{stima y M e y a stima cal S} 
%and using 
Lemma \ref{lemmacomposizioneSTANDARD} and the smallness condition on $\eta$ one gets
%one gets that there exists $\sigma \equiv \sigma(\beta_0, k_0, M) \gg 0$ 
%large enough such that for $s\geq s_0$, if $\| \eta \|_{s_0 + \sigma}^{k_0, \gamma} \ll_{s} 1$, then 
\begin{equation}\label{stima y M e y a stima cal S p}
\begin{aligned}
\sup_{y \in [- 1, 0]} \Big\| \langle D \rangle^{M + 1} 
\circ 
{\rm Op}\big(y^{M + q} \partial_\lambda^{k_1} \partial_\vphi^{\beta_1} \big( e^{y a} f_q \big) \big) \Big\|_{0, s, 0} 
& \lesssim_{s, M} 
\gamma^{- |k_1|} \big( 1 + \| \eta \|_{s + \sigma}^{k_0, \gamma}\big)\,,
\end{aligned}
\end{equation}
for some $\s \equiv \sigma(\beta_0, k_0, M)\gg1$ and $s_0\leq s\leq \bar{s}-\s$.
Now, taking $k \in \N^{\nu + 1}$, $\beta \in \N^{\nu}$, $|\beta| \leq \beta_0$, one has that 
\begin{equation}\label{D M - beta - k0 - 2 cal Sp}
\begin{aligned}
& \langle D \rangle^{M + 1} \partial_\lambda^k \partial_\vphi^\beta \cS_q(\vphi) 
\langle D \rangle^{- \tc}  
= \sum_{\begin{subarray}{c}
\beta_1 + \beta_2 = \beta \\
k_1 + k_2 = k
\end{subarray}} C(\beta_1, \beta_2, k_1, k_2) \cR_{\beta_1, \beta_2, k_1, k_2}(\vphi) \,, 
\\& 
\cR(\vphi)  \equiv \cR_{\beta_1, \beta_2, k_1, k_2}(\vphi) 
\\&  
:=  \int_{- 1}^{\bar y} \frac{1}{y^M}  \langle D \rangle^{M + 1} \circ 
{\rm Op}\big(y^{M + q} \partial_\lambda^{k_1} \partial_\vphi^{\beta_1}\big(e^{y a} f_q \big) \big) 
\circ 
\partial_\lambda^{k_2} \partial_\vphi^{\beta_2}\Gamma(\vphi, y)\, d y 
\circ \langle D \rangle^{- \tc}\,,
\end{aligned}
\end{equation}
and hence, to get the result, for any $k_1, k_2, \beta_1, \beta_2$, $k_1 + k_2 = k$, 
$\beta_1 + \beta_2 = \beta$ it suffices to estimate the linear operator $\cR(\vphi)$. 
By using again the inequality \eqref{stima con cauchy schwatz integrale striscia}, 
and Lemma \ref{lemma: action Sobolev}, and the smallness condition on $\eta$, one gets
\[
\begin{aligned}
\| \cR \psi \|_s & \lesssim \Big\| \langle D \rangle^{M + 1} 
\circ 
{\rm Op}\big(y^{M + q} \partial_\lambda^{k_1} \partial_\vphi^{\beta_1} \big( e^{y a} f_q \big) \big) 
\circ 
\partial_\lambda^{k_2} \partial_\vphi^{\beta_2}\Gamma(\vphi, y)\,  
[ \langle D \rangle^{-\tc} \psi] \Big\|_{L^2_y H^s_{\vphi, x}} 
\\& 
\stackrel{\eqref{stima y M e y a stima cal S p}}{\lesssim_{s, M}} 
\gamma^{- |k_1|}   \| \partial_\lambda^{k_2} \partial_\vphi^{\beta_2}\Gamma(\vphi, y)\,  
[ \langle D \rangle^{-\tc} \psi] \|_{L^2_y H^s_{\vphi, x}} 
\\&\qquad
+ \gamma^{- |k_1|}\| \eta \|_{s + \sigma}^{k_0, \gamma} \| \partial_\lambda^{k_2} \partial_\vphi^{\beta_2}\Gamma(\vphi, y)\,  
[ \langle D \rangle^{- \tc} \psi] \|_{L^2_y H^{s_0}_{\vphi, x}} 
\\& 
%\lesssim_s 
%\gamma^{- |k_1|}  \| \partial_\lambda^{k_2} \partial_\vphi^{\beta_2}\Gamma(\vphi, y)\,  
%[ \langle D \rangle^{- \tc} \psi] \|_{\cO^s} 
%\\&\qquad 
%+ \gamma^{- |k_1|}\| \eta \|_{s + \sigma}^{k_0, \gamma} \| \partial_\lambda^{k_2} \partial_\vphi^{\beta_2}\Gamma(\vphi, y)\,  [ \langle D \rangle^{- \tc} \psi] \|_{\cO^{s_0}} 
%\\& 
\stackrel{\eqref{prop L2 vhi Hs Hs vphi cal 2 xy},{\rm Lemma}\,\ref{lemma tame F Gamma vphi y che palle}}{\lesssim_{s, M}} 
\gamma^{- |k|} \Big( \| \eta \|_{s_0 + \sigma}^{k_0, \gamma}  \| \psi \|_s 
+ \| \eta \|_{s + \sigma}^{k_0, \gamma} \| \psi \|_{s_0}\Big)\,. 
\end{aligned}
\]
This latter estimate, together with \eqref{D M - beta - k0 - 2 cal Sp},  implies 
the first bound in 
\eqref{stimacalSS}.
The Lipschitz bound follows similarly.

\medskip

\noindent
{\bf  Estimate of \eqref{opSN}.}
By recalling the definition of $\cS_N$ in \eqref{opSN}, and taking 
$k \in \N^{\nu + 1}$, $\beta \in \N^{\nu}$, $|\beta| \leq \beta_0$, 
and $\tc\geq\beta_0+k_0+2$,
%$0 \leq M \leq N - \beta_0 - k_0$, 
one has that 
\begin{equation}\label{D M - beta - k0 - 2 cal SN}
\begin{aligned}
& \langle D \rangle^{M } \partial_\lambda^k \partial_\vphi^\beta \cS_N(\vphi) \langle D \rangle^{- \tc}  = \sum_{\begin{subarray}{c}
\beta_1 + \beta_2 = \beta \\
k_1 + k_2 = k
\end{subarray}} C(\beta_1, \beta_2, k_1, k_2) \cR_{\beta_1, \beta_2, k_1, k_2}(\vphi) \,, 
\end{aligned}
\end{equation}
where
\[
\begin{aligned}
 &\cR(\vphi)  \equiv \cR_{\beta_1, \beta_2, k_1, k_2}(\vphi) 
 \\
 & :=  \int_{- 1}^{\bar y}   \langle D \rangle^{M } \circ \partial_\lambda^{k_1} \partial_\vphi^{\beta_1} \cH_N(\vphi,- y) \circ \partial_\lambda^{k_2} \partial_\vphi^{\beta_2}\Gamma(\vphi, y)\, d y \circ \langle D \rangle^{- \tc}\,,
\end{aligned}
\]
and hence, to get the result, for any $k_1, k_2, \beta_1, \beta_2$, $k_1 + k_2 = k$, 
$\beta_1 + \beta_2 = \beta$ it suffices to estimate the linear operator $\cR(\vphi)$. 
To do this we reason as follows. First of all, by the smallness condition on $\eta$ with $\s\gg1$ large enough, we shall apply 
 estimate \eqref{stima resto parametrica caloreB} 
in Proposition \ref{proposizione totale parametrica calore} 
with $n_1=M$, $n_2=0$, $N:=M+\beta_0+k_0+1$,
%(with $N,n_1,n_2$ large enough w.r.t. $M$),
in order to get that  
%by taking 
%$\sigma \equiv \sigma(\beta_0, k_0, M) \gg 0$ large enough such that 
%if $\| \eta \|_{s_0 + \sigma}^{k_0, \gamma} \ll 1$, one gets that 
$\langle D \rangle^M \partial_\vphi^{\beta} \cH_N$ is a $\cD^{k_0}$-tame 
operator with tame constant satisfying 
\begin{equation}\label{stimaclHHNN}
\begin{aligned}
\sup_{y \in [- 1, 0]} {\mathfrak M}_{\langle D \rangle^M \partial_\vphi^{\beta} \cH_N}(s) 
& \lesssim_{\bar{s}, M}
  \| \eta \|_{s + \sigma}^{k_0, \gamma}\,, \quad \forall s_0 \leq s \leq \bar{s}-\s\,,
\\
\| \partial_\vphi^\beta\Delta_{12}
\langle D \rangle^{M} 
 \cH_N(y)
 \|_{\mathcal{L}({H}^{p};{H}^{p})}
&\lesssim_{\bar{s},M} \|\eta_{1}-\eta_{2}\|_{p+\widetilde{\s}}\,,\;\;\;s_0\leq p\leq s_0+\s-\widetilde{\s}\,,
\end{aligned}
\end{equation}
for some $\widetilde{\s}\leq \s$.
%Therefore, the estimate above, together with 
%the inequality \eqref{stima con cauchy schwatz integrale striscia},  implies 
%%one has that for any $s \geq s_0$,
%\begin{equation}\label{stima cal R per cal Sp 0bis}
%\begin{aligned}
% \| \cR \psi \|_s  &\lesssim 
%\Big\| \langle D \rangle^{M } \circ \partial_\lambda^{k_1} \partial_\vphi^{\beta_1} \cH_N(\vphi, -y) 
%\circ \partial_\lambda^{k_2} \partial_\vphi^{\beta_2}\Gamma(\vphi, y) 
%[ \langle D \rangle^{-\tc} \psi] \Big\|_{L^2_y H^s_{\vphi, x}}\,.
%\end{aligned}
%\end{equation}
%which implies (together with \eqref{stima cal R per cal Sp 0bis}) 
%that  
Therefore, by the inequality \eqref{stima con cauchy schwatz integrale striscia} and Lemma \ref{lemma operatore e funzioni dipendenti da parametro},
one gets
\[
\begin{aligned}
\| \cR \psi \|_s 
&\lesssim 
\Big\| \langle D \rangle^{M } \circ \partial_\lambda^{k_1} \partial_\vphi^{\beta_1} \cH_N(\vphi, -y) 
\circ \partial_\lambda^{k_2} \partial_\vphi^{\beta_2}\Gamma(\vphi, y) 
[ \langle D \rangle^{-\tc} \psi] \Big\|_{L^2_y H^s_{\vphi, x}}
\\&  \stackrel{\eqref{stimaclHHNN}}{\lesssim_{\bar{s}, M}} 
\gamma^{- |k_1|}  \| \partial_\lambda^{k_2} \partial_\vphi^{\beta_2}\Gamma(\vphi, y)\,  
[ \langle D \rangle^{- \tc} \psi] \|_{L^2_y H^s_{\vphi, x}} 
\\&\qquad 
+ \gamma^{- |k_1|} \| \eta \|_{s + \sigma}^{k_0, \gamma} \| \partial_\lambda^{k_2} \partial_\vphi^{\beta_2}\Gamma(\vphi, y)\,  
[ \langle D \rangle^{-\tc} \psi] \|_{L^2_y H^{s_0}_{\vphi, x}} 
\\& 
%\lesssim_s 
%\gamma^{- |k_1|}  \| \partial_\lambda^{k_2} \partial_\vphi^{\beta_2}\Gamma(\vphi, y)\,  
%[ \langle D \rangle^{- \tc} \psi] \|_{\cO^s} 
%%\\&\qquad 
%+ \gamma^{- |k_1|} \| \eta \|_{s + \sigma}^{k_0, \gamma} \| \partial_\lambda^{k_2} \partial_\vphi^{\beta_2}\Gamma(\vphi, y)\,  
%[ \langle D \rangle^{- \tc} \psi] \|_{\cO^{s_0}} 
%\\& 
\stackrel{\eqref{prop L2 vhi Hs Hs vphi cal 2 xy},{\rm Lemma}\,\ref{lemma tame F Gamma vphi y che palle}}{\lesssim_{\bar{s}, M}} 
\gamma^{- |k|} \Big( \| \eta \|_{s_0 + \sigma}^{k_0, \gamma}  \| \psi \|_s 
+ \| \eta \|_{s + \sigma}^{k_0, \gamma} \| \psi \|_{s_0}\Big)\,. 
\end{aligned}
\]
This latter estimate, together with \eqref{D M - beta - k0 - 2 cal SN}, 
 implies \eqref{stimacalSS} on the operator $\mathcal{S}_{N}$. 
The bound on the Lipschitz variation follows similarly.

\noindent
By the discussion above one has that the operator $\mathcal{S}$ in \eqref{formulaespdiSS}
satisfies the bounds \eqref{stimacalSS}.
This concludes the proof.
\end{proof}

We are in position to conclude the proof of our main result.

\begin{proof}[{\bf Proof of Proposition \ref{lemma espansione pseudo-diff w x 0}}]
By formula
\eqref{def parte dominante parte trascurata w x 0}, and recalling the expansions \eqref{formulaespdiW} and 
\eqref{formulaespdiSS}, one gets
\[
 \cW(\vphi) =\th|D|\big(\tanh(\th|D|)-1\big)[\psi] + \cR_M(\vphi)[\psi]
\]
which is formula \eqref{espansionewzero}, where we set
\[
\cR_M(\vphi):=\widetilde{\cQ}+ \sum_{q = 1}^{2N} \cS_q(\vphi) 
+ \cS_N(\vphi)+\mathcal{S}_{\mathcal{L}}(\vphi)\,,
\]
with $N=M+\beta_0+k_0+1$.
The estimates \eqref{stimefinalew00} on $\cR_M(\vphi)$ follow by Lemma \ref{lem:espandoww} 
for the analysis of $\cQ(\vphi)$,
and  Lemma \ref{lemma stima resto cal S DN} for the estimates on $\cS(\vphi)$. 
\end{proof}

\subsection{Proof of Theorem \ref{lemma totale dirichlet neumann}}
By Lemma \ref{lem:strutturaDNconw} (see \eqref{eq:112TOTALE})
the Dirichlet-Neumann operator in \eqref{eq:112a} has the form
\begin{equation*}
\begin{aligned}
 G(\eta)\psi 
 & = \th^{-1}(1+|\nabla\eta|^2)\op(A) [\psi] +\th^{-1}(1+|\nabla\eta|^2) w_{| y = 0} 
 -\op(\ii\nabla\eta\cdot\x\chi(\x))\psi(x)
\end{aligned}
\end{equation*}
where $A$ is the symbol given by Lemma \ref{divo1} and $w$ is given in \eqref{funz:ww}. 
By the expansion given by Lemma \ref{divo1} (see \eqref{simboliordine1}), 
and recalling that $(1+|\nabla\eta|^2)b=1$, 
we note that
\begin{equation}\label{contagion3}
\begin{aligned}
&\th^{-1}(1+|\nabla\eta|^2)A=\ii   \nabla\eta\cdot\x\chi(\x)+b^{-1} \chi(\x)\sqrt{b|\x|^{2}-(b\nabla\eta\cdot\x)^{2}}+
b^{-1}A_{\leq0}
\\&=|\x|\chi(\x)
+\underbrace{|\x|\chi(\x)\big(\sqrt{1+|\nabla\eta|^2-|\x|^{-2}(\nabla\eta\cdot\x)^2}-1\big)}_{=:\mathfrak{m}_1\in S^1}
+\underbrace{b^{-1}A_{\leq0}}_{=:\mathfrak{m}_0\in S^0}
+\ii  \nabla\eta\cdot\x\chi(\x)\,.
\end{aligned}
\end{equation}
Since the symbol $A$ is real and momentum preserving, 
we deduce that $\mathfrak{m}_1,\mathfrak{m}_0$ are real (see 
Remark \ref{rmk:operatorerealtoreal}) and momentum preserving.
On the other hand, by
Proposition \ref{lemma espansione pseudo-diff w x 0} we have
\begin{equation*}
\begin{aligned}
 G(\eta)\psi 
 & =\th^{-1} b^{-1}\op(A) [\psi] +\th^{-1}b^{-1}\big(\th|D|\big(\tanh(\th|D|)-1\big)[\psi] -\cR_N(\vphi)[\psi]\big)
 -\nabla\eta\cdot\nabla \psi(x)\,.
\end{aligned}
\end{equation*}
Then using that
\[
\th^{-1}\op\big(\ii  \mathtt h \nabla\eta\cdot\x\chi(\x)\big)\psi-\nabla\eta\cdot\nabla\psi\equiv0\,,
\]
one deduces that
\begin{equation}\label{contagion2}
 G(\eta)\psi 
 =|D|\tanh(\th|D|)[\psi]+\op\big(\mathfrak{m}_1+\mathfrak{m}_0\big)\psi+\mathfrak{R}_0(\vphi)\psi\,,
\end{equation}
where
\[
\mathfrak{R}_0(\vphi):=|\nabla\eta|^2\big(|D|\big(\tanh(\th|D|)-1\big)
-\th^{-1}b^{-1}\cR_N(\vphi)\,.
\]
The remainder $\mathfrak{R}_0$ satisfies \eqref{stimefinaliRR} by \eqref{stimefinalew00}
and using that $(\tanh(\th|D|)-1)\in S^{-\infty}$. 
Moreover, since $G(\eta)[\cdot], \op(\mathfrak{m}_1+\mathfrak{m}_0)$
are real-to-real, one has that also $\mathfrak{R}_0$ is so (i.e. its matrix coefficients  satisfy \eqref{piscina}).
In order to obtain \eqref{pseudoespansione}
it remains to prove that 
\[
\op\big(\mathfrak{m}_1+\mathfrak{m}_0\big)=\opw\big(a_{G}+b_{G}\big)\,,
\]
for some real valued symbols $a_{G}\in S^{1}$, $b_{G}\in S^{0}$ satisfying the symmetry condition \eqref{contagion1}.
To do so, we first need some preliminary result provided in the following subsection.

\subsubsection{Spectrally localized symbols}
Let us consider a cut-off function $\chi \in C^\infty_c (\R, \R)$ such that 
\begin{equation}\label{cut off bony}
\begin{aligned}
& 0 \leq \tilde{\chi} \leq 1, \quad \tilde{\chi}(s) = \tilde{\chi}(- s), \quad \forall s \in \R\, \\ 
& \tilde{\chi}(s) = 1, \quad \forall s \in [- 1, 1], \quad \chi(s) = 0, \quad \forall s \in \R \setminus [- 2, 2]\,. 
\end{aligned}
\end{equation}
Given a symbol $a (\vphi, x, \xi) \in S^m$, we define the symbol $a_{\tilde{\chi}}$ as 
\begin{equation}\label{def a chi fourier}
a_{\tilde{\chi}}(\vphi, x, \xi) := \sum_{j \in \Gamma^*} \tilde{\chi}\Big( \frac{|j|}{\epsilon\langle \xi \rangle} \Big) \widehat a(\vphi, j, \xi) e^{\ii j \cdot x} = \sum_{\begin{subarray}{c}
\ell \in \Z^\nu \\
j \in \Gamma^*
\end{subarray}} \tilde{\chi}\Big( \frac{|j|}{\epsilon\langle \xi \rangle} \Big) \widehat a(\ell, j, \xi) e^{\ii (\ell \cdot \vphi + j \cdot x)}
\end{equation}
where we fix $0 < \epsilon \ll 1$ small enough. 
We then define  its \emph{Bony} quantization as 
\begin{equation}\label{bonyquantization}
{\rm Op}^B(a) := {\rm Op} (a_{\tilde{\chi}}) \quad \text{and} \quad {\rm Op}^{BW}(a) 
:= {\rm Op}^W(a_{\tilde{\chi}})\,. 
\end{equation}
The following lemma holds.
\begin{lemma}\label{lemma proprieta bony weil}
Consider a symbol $a \in S^m$, $m\in\R$.

\noindent
$(i)$ The operator ${\mathcal A} := {\rm Op}^{BW}(a)$ is {\bf spectrally localized}, namely there 
are two constants $0 < \mathtt c_1 < \mathtt c_2$  such that 
\[
\widehat{\mathcal A}(\ell)_j^{j'} \neq 0 \quad \Longrightarrow\quad \mathtt c_1 |j| \leq |j'| \leq \mathtt c_2 |j|\,. 
\]

\noindent
$(ii)$ For any $N \in \N$, $a - a_{\tilde{\chi}} \in S^{- N}$ and 
\[
\| a - a_\chi \|_{ - N, s, 0}^{k_0, \gamma} \lesssim_{s, N,m} 
\| a \|_{m, s + |m| + N, 0}^{k_0, \gamma}, \quad \forall s \geq 0\,. 
\]

\noindent
$(iii)$ The operator ${\mathcal R} := {\rm Op}^W(a - a_{\tilde{\chi}}) = {\rm Op}^{W}(a) - {\rm Op}^{BW}(a)$ 
satisfies the following properties. 
For any $\beta \in \N^\nu$, $N \in \N$,  the operator 
$\langle D \rangle^{N} \partial_\vphi^\beta {\mathcal R}$ is 
${\mathcal D}^{k_0}$-tame with tame constant satisfying 
$$
{\mathfrak M}_{\langle D \rangle^N \partial_\vphi^\beta {\mathcal R}}(s) 
\lesssim_{\bar{s},N} \| a \|_{- N, s + |\beta| + N + |m|, 0}^{k_0, \gamma}, \quad \forall s \geq s_0\,. 
$$
The latter estimate clearly implies that for any $\mathtt c \geq 0$
$$
{\mathfrak M}_{\langle D \rangle^N \partial_\vphi^\beta {\mathcal R} \langle D \rangle^{- \mathtt c}}(s) 
\lesssim_{N}
 {\mathfrak M}_{\langle D \rangle^N \partial_\vphi^\beta {\mathcal R}}(s) 
\lesssim_{\bar{s}, N} \| a \|_{- N, s + |\beta| + N, 0}^{k_0, \gamma}\,, \quad \forall s \geq s_0 \,. 
$$
\end{lemma}
\begin{proof}
{\sc Proof of $(i)$.}
One has (see \eqref{Weyl}, \eqref{def a chi fourier}) 
\[
{\mathcal A}(\vphi)[u](x) =  
\sum_{j\in\Gamma^*} \Big( \sum_{j'\in\Gamma^*} \tilde{\chi}\Big( \frac{|j - j'|}{ \epsilon \langle \frac{j+ j'}{2}\rangle} \Big)
\widehat{a}\big(\vphi, j- j', \frac{j+ j'}{2}\big)\widehat u(j')\Big)e^{\im j\cdot x}\,.
\]
The latter formula implies that for $\ell \in \Z^\nu$, $j, j' \in \Gamma^*$, 
\[
\widehat{\mathcal A}(\ell)_j^{j'} 
=  
\tilde{\chi}\Big( \frac{|j - j'|}{\epsilon\langle \frac{j+ j'}{2}\rangle} \Big) \widehat a\big(\ell, j - j', \frac{j + j'}{2} \big)\,.
\]
By the properties of the cut-off function $\tilde{\chi}$ in \eqref{cut off bony}, 
one immediately gets that 
\[
\widehat{\mathcal A}(\ell)_j^{j'} \neq 0 \qquad \Longrightarrow \qquad 
|j - j'| \leq 2 \epsilon \langle \frac{j + j'}{2} \rangle \leq C \epsilon (|j| + |j'|)\,,
\]
for some constant $C > 0$. By the latter inequality, and taking $\epsilon\ll1$ small enough, one deduces
the equivalence $|j|\sim |j|'$.
Then item $(i)$ follows.
%The proof of the item $(i)$ is then concluded. 

\noindent
{\sc Proof of $(ii)$.} By \eqref{cut off bony}, one has that 
$$
1 - \tilde{\chi}\Big( \frac{|j|}{\epsilon \langle \xi \rangle} \Big) \neq 0 \quad \text{if and only if} 
\quad \langle \xi \rangle \lesssim |j|,
$$
hence, given $N \in \N$, $k \in \N^{\nu + 1}$, $|k| \leq k_0$ one has that 
$$
\begin{aligned}
\langle \xi \rangle^{2 N} \| \partial_\lambda^k ( a - a_{\tilde{\chi}}) (\cdot, \xi) \|_{s}^2 & = \langle \xi \rangle^{2 N} \sum_{\begin{subarray}{c}
\ell \in \Z^\nu \\
j \in \Gamma^*
\end{subarray}}  \langle \ell, j \rangle^{2 s}\Big| 1 - \tilde{\chi}\Big( \frac{|j |}{\epsilon\langle \xi \rangle} \Big) \Big|^2 |\widehat{\partial_\lambda^k a}(\ell, j, \xi)|^2 \\
& \lesssim \langle \xi \rangle^{2 (N + |m|)} \sum_{|j| \gtrsim \langle \xi \rangle}  \langle \ell, j \rangle^{2 s} |\widehat{\partial_\lambda^k a}(\ell, j, \xi)|^2 \langle \xi \rangle^{- 2 m} \\
& \stackrel{\langle \xi \rangle^{2 (N + |m|)} \lesssim |j|^{2 (N + |m|)}}{\lesssim_{N,m}} \sum_{\ell, j}  \langle \ell, j \rangle^{2 (s + N + |m|)} |\widehat{\partial_\lambda^k a}(\ell, j, \xi)|^2 \langle \xi \rangle^{- 2 m}  \\
& \lesssim_{N,m} \| \partial_\lambda^k a \|_{m, s + N + |m|, 0} \lesssim_{N,m} 
\gamma^{- |k|} \| a \|_{m, s + N + |m|, 0}^{k_0, \gamma}\,.
\end{aligned}
$$
The latter chain of inequalities imply the claimed bound. 

\noindent
{\sc Proof of $(iii)$.}
It follows by $(i)$, using also Lemma \ref{constantitamesimbolo}. 
\end{proof}

\subsubsection{Conclusions}\label{self-adjoint-pezzi-DN}
Consider the expansion \eqref{contagion2} with symbols $\mathfrak{m}_1$, $\mathfrak{m}_0$
given in \eqref{contagion3}.
We split
$$
{\rm Op}(\mathfrak{m}_1+\mathfrak{m}_0) = {\rm Op}^{BW}(\mathfrak{m}_1) + {\rm Op}^{BW}(\mathfrak{q}_0) + \mathfrak{R}_1(\vphi)
$$
where, by Lemma \ref{divo1}, Lemma \ref{passaggiostweyl} and Lemma \ref{lemma proprieta bony weil} $(ii)$, $(iii)$, $\mathfrak{q}_0 \in S^0$ satisfies the estimates 
\begin{equation}\label{stima mathfrak q0}
\begin{aligned}
\| \mathfrak{q}_0 \|_{0, s, \alpha}^{k_0, \gamma} 
&\lesssim_{s, \alpha} \| \eta \|_{s + \sigma}^{k_0, \gamma}, \quad s \geq s_0, \quad \alpha \in \N\,,
\\
 \| \Delta_{12}\mathfrak{q}_0 \|_{0, p, \alpha}
 &\lesssim_{p,\alpha}\| \eta_1-\eta_2 \|_{p + \sigma}\,, \;\;\;p\geq s_0\,,
 \end{aligned}
\end{equation}
for $\sigma \gg 0$ large enough and $\| \eta \|_{s_0 + \sigma}^{k_0, \gamma} \lesssim_{s} 1$ 
and the remainder $\mathfrak{R}_1$ satisfies the 
estimates \eqref{stimefinaliRR}. 
In view of Remark \ref{rmk:operatorerealtoreal} we also have that 
$\opbw(\mathfrak{m}_1)$, $\opbw(\mathfrak{q}_0)$ are still real-to-real, and so it is the remainder 
$\mathfrak{R}_1$. Moreover the symbol $\mathfrak{q}_0$ is momentum preserving 
since $\mathfrak{m}_1,\mathfrak{m}_0$ are so.
We then set ${\mathfrak R} := {\mathfrak R}_0 + {\mathfrak R}_1$, and we get that 
\begin{equation}\label{sistemi dinamici merda 10}
G(\eta) = |D| \tanh(\mathtt h |D|) + {\rm Op}^{BW}(\mathfrak m_1 + \mathfrak q_0) + {\mathfrak R}(\vphi)
\end{equation}
where the remainder ${\mathfrak R}$ satisfies \eqref{stimefinaliRR}.
Since $G(\eta) = G(\eta)^*$ (see for instance \cite{Lan2}) 
and $\mathfrak m_1(\vphi , x, \xi)$ is real, one has that 
\begin{equation}\label{sistemi dinamici merda 11}
{\rm Op}^{BW}(\mathfrak q_0) + \mathfrak{R}(\vphi) = {\rm Op}^{BW}(\overline{\mathfrak q}_0) + \mathfrak{R}(\vphi)^*
\end{equation}
and hence we set
\begin{equation}\label{sistemi dinamici merda 12}
{\mathfrak Q}(\vphi) := \mathfrak{R}(\vphi) - \mathfrak{R}(\vphi)^* = {\rm Op}^{BW}(\overline{\mathfrak q}_0 - {\mathfrak q}_0)\,. 
\end{equation}
By Lemma \ref{lemma proprieta bony weil}-$(i)$,  the operator ${\mathfrak Q}(\vphi) = {\rm Op}^{BW}(\overline{\mathfrak q}_0 - {\mathfrak q}_0)$ is {\bf spectrally localized}, namely 
\begin{equation}\label{sistemi dinamici merda - 1}
\widehat{\mathfrak Q}(\ell)_j^{j'} = \widehat{\mathfrak R}(\ell)_j^{j'} 
-  \overline{\widehat{\mathfrak R}(- \ell)_{j'}^{j}} \neq 0 \quad \Longrightarrow \quad j \sim j'
\end{equation}
where $j \sim j'$ means that there are constants $0 < \mathtt c_1 < \mathtt c_2$ such that $\mathtt c_1| j| \leq |j'| \leq \mathtt c_2 |j|$. Since ${\mathfrak R}$, ${\mathfrak R}^*$ are real-to-real, one has that
\begin{equation}\label{reality frak R}
\overline{\widehat{\mathfrak R}(- \ell)_{j'}^{j}} = \widehat{\mathfrak R}( \ell)_{- j'}^{- j}\,.
\end{equation}
Moreover, since ${\mathfrak R}(\vphi)$ satisfies \eqref{stimefinaliRR}, one gets that for any 
$|k| \leq k_0$, and taking $\mathtt{c}\geq\beta_0+k_0+2$,
\[
\begin{aligned}
 \sum_{\ell, j} \langle \ell , j \rangle^{2 s}\langle j \rangle^{2 M}
 & \langle \ell - \ell' \rangle^{2 \beta_0}  |\partial_\lambda^k 
 \widehat{\mathfrak R}(\ell - \ell')_j^{j'} |^2 \langle j' \rangle^{- 2 \mathtt{c}} 
 \\& 
 \lesssim_{\bar{s},M} 
 \langle \ell', j' \rangle^{2 s} (\| \eta \|_{s_0 + \sigma}^{k_0, \gamma})^2 + \langle \ell', j' \rangle^{2 s_0} 
 (\| \eta \|_{s + \sigma}^{k_0, \gamma})^2\,. 
\end{aligned}
\]
The latter inequality implies that for $M \gg 0$, $\beta_0 \gg 0$ large enough, $s \geq s_0$, for any $\ell, \ell' \in \Z^\nu$, $j, j' \in \Gamma^*$
\begin{equation}\label{sistemi dinamici merda 0}
\begin{aligned}
 \gamma^{2 |k|}\langle \ell , j \rangle^{2 s}\langle j \rangle^{2 M} \langle \ell - \ell' \rangle^{2 \beta_0} 
& |\partial_\lambda^k \widehat{\mathfrak R}(\ell - \ell')_j^{j'} |^2 \langle j' \rangle^{- 2\mathtt{c}}  
\\& 
\lesssim_{\bar{s},M}
\langle \ell', j' \rangle^{2 s} (\| \eta \|_{s_0 + \sigma}^{k_0, \gamma})^2 + \langle \ell', j' \rangle^{2 s_0} 
(\| \eta \|_{s + \sigma}^{k_0, \gamma})^2
\end{aligned}
\end{equation}
and hence for $j \sim j'$, 
\begin{equation}\label{sistemi dinamici merda 1}
\begin{aligned}
 \gamma^{2 |k|} &\langle \ell , j \rangle^{2 s}\langle j \rangle^{2 M } \langle \ell - \ell' \rangle^{2 \beta_0 } 
 |\partial_\lambda^k \widehat{\mathfrak R}(\ell '- \ell)_{j'}^{j} |^2 \langle j' \rangle^{- 2 \mathtt{c}} 
\\& 
\stackrel{j \sim j', \eqref{reality frak R}}{\lesssim_{\bar{s},M}}  
\gamma^{2 |k|} \langle \ell , j' \rangle^{2 s}\langle j' \rangle^{2 M } \langle \ell - \ell' \rangle^{2 \beta_0 }  
|\partial_\lambda^k \widehat{\mathfrak R}(\ell - \ell')_{- j'}^{- j} |^2 \langle j \rangle^{- 2 \mathtt{c}} 
\\& 
%\lesssim_s 
%\langle \ell', j \rangle^{2 s} (\| \eta \|_{s_0 + \sigma}^{k_0, \gamma})^2 + \langle \ell', j\rangle^{2 s_0} 
%(\| \eta \|_{s + \sigma}^{k_0, \gamma})^2 
%\\& 
\lesssim_{\bar{s},M}
%\stackrel{j \sim j'}{\lesssim_s} 
\langle \ell', j' \rangle^{2 s} (\| \eta \|_{s_0 + \sigma}^{k_0, \gamma})^2 + \langle \ell', j' \rangle^{2 s_0} 
(\| \eta \|_{s + \sigma}^{k_0, \gamma})^2\,. 
\end{aligned}
\end{equation}
Thus, by \eqref{sistemi dinamici merda - 1}, \eqref{sistemi dinamici merda 0}, \eqref{sistemi dinamici merda 1}, 
one obtains that the operator ${\mathfrak Q}$ satisfies for 
$M \gg 0$, $\beta_0 \gg 0$ large enough, $s \geq s_0$, $|k| \leq k_0$, 
for any $\ell, \ell' \in \Z^\nu$, $j, j' \in \Gamma^*$
\begin{equation}\label{sistemi dinamici merda 2}
\begin{aligned}
 \gamma^{2 |k|} \langle \ell , j \rangle^{2 s}\langle j \rangle^{2 M} &\langle \ell - \ell' \rangle^{2 \beta_0}  
 |\partial_\lambda^k \widehat{\mathfrak Q}(\ell - \ell')_j^{j'} |^2 \langle j' \rangle^{- 2 \mathtt{c}}  
 \\& 
 \lesssim_{\bar{s},M}
 \langle \ell', j' \rangle^{ 2s} (\| \eta \|_{s_0 + \sigma}^{k_0, \gamma})^2 + \langle \ell', j' \rangle^{ 2 s_0} 
 (\| \eta \|_{s + \sigma}^{k_0, \gamma})^2\,.
\end{aligned}
\end{equation}
Then, we estimate for any $|\beta| \leq \beta_0 - 4 s_0$, $|k| \leq k_0$, the operator 
$ A := \langle D \rangle^{M } \partial_\vphi^\beta \partial_\lambda^k {\mathfrak Q}(\vphi) 
\langle D \rangle^{- \mathtt{c}}$. 
By a repeated use of the Cauchy Schwartz inequality, using that 
\[
\sum_{\ell', j'} \frac{1}{\langle \ell - \ell' \rangle^{2 s_0} \langle j - j'' \rangle^{2 s_0}}\,,\,\;\;
\sum_{\ell, j} \frac{1}{\langle \ell - \ell' \rangle^{2 s_0} \langle j - j'' \rangle^{2 s_0}} \leq C_0\,,
\] 
and since by the conservation of momentum $|j - j'| \lesssim |\ell - \ell'|$, one gets that 
\[
\begin{aligned}
&\| A h  \|_s^2    = 
\sum_{\ell, j} \langle \ell, j \rangle^{2 s} \Big| \sum_{\ell', j'} \langle j \rangle^{ M} \im^{|\beta|}( \ell - \ell' )^{\beta}  
\partial_\lambda^k \widehat{\mathfrak Q}(\ell - \ell')_j^{j'} \langle j' \rangle^{- \mathtt{c}}  \widehat h(\ell', j')  \Big|^2  
\\& 
\leq 
\sum_{\ell, j} \langle \ell, j \rangle^{2 s} \Big( \sum_{\ell', j'} \langle j \rangle^{ M}| \ell - \ell' |^{|\beta|}  
|\partial_\lambda^k \widehat{\mathfrak Q}(\ell - \ell')_j^{j'}| \langle j' \rangle^{- \mathtt{c}}  |\widehat h(\ell', j')|  \Big)^2 
\\& 
\lesssim 
\sum_{\ell, j} \langle \ell, j \rangle^{2 s}  \sum_{\ell', j'} \langle j \rangle^{ 2 M} 
\langle  \ell - \ell'  \rangle^{2 (|\beta| + s_0)} 
\langle j - j' \rangle^{2 s_0}  
|\partial_\lambda^k \widehat{\mathfrak Q}(\ell - \ell')_j^{j'}|^2 \langle j' \rangle^{-  2\mathtt{c}}  
|\widehat h(\ell', j')|^2 
%\\& 
%\stackrel{|j - j'| \lesssim |\ell - \ell'|}{\lesssim} 
%\sum_{\ell, j} \langle \ell, j \rangle^{2 s}  \sum_{\ell', j'} \langle j \rangle^{ 2 M} 
%\langle  \ell - \ell'  \rangle^{2 (|\beta| + 2 s_0)}   
%|\partial_\lambda^k \widehat{\mathfrak Q}(\ell - \ell')_j^{j'}|^2 \langle j' \rangle^{-  2(\beta_0 + k_0 + 2)}  |\widehat h(\ell', j')|^2 
\\& 
\stackrel{|j - j'| \lesssim |\ell - \ell'|}{\lesssim}
\sum_{\ell', j'} |\widehat h(\ell', j')|^2 \sum_{\ell, j} \dfrac{ \langle \ell, j \rangle^{2 s} \langle j \rangle^{ 2 M} 
\langle  \ell - \ell'  \rangle^{2 (|\beta| + 4 s_0)}   
|\partial_\lambda^k \widehat{\mathfrak Q}(\ell - \ell')_j^{j'}|^2 
\langle j' \rangle^{-  2\mathtt{c}}}{\langle \ell - \ell' \rangle^{2 s_0} \langle j - j' \rangle^{2 s_0}}  
%\\& 
%\stackrel{|\beta| + 4 s_0 \leq \beta_0 }{\lesssim_s} 
%\sum_{\ell', j'} |\widehat h(\ell', j')|^2 \sum_{\ell, j} \dfrac{ \langle \ell, j \rangle^{2 s} \langle j \rangle^{ 2 M}
%\langle  \ell - \ell'  \rangle^{2\beta_0}   
%|\partial_\lambda^k \widehat{\mathfrak Q}(\ell - \ell')_j^{j'}|^2 
%\langle j' \rangle^{-  2(\beta_0 + k_0 + 2)}}{\langle \ell - \ell' \rangle^{2 s_0} \langle j - j' \rangle^{2 s_0}}  
\\& 
\stackrel{|\beta| + 4 s_0 \leq \beta_0 , \eqref{sistemi dinamici merda 2}}{\lesssim_{\bar{s},M}} 
\gamma^{- 2 |k|}\sum_{\ell', j'} |\widehat h(\ell', j')|^2 \Big( \langle \ell', j' \rangle^{2 s} 
(\| \eta \|_{s_0 + \sigma}^{k_0, \gamma})^2 + \langle \ell', j' \rangle^{2 s_0} 
(\| \eta \|_{s + \sigma}^{k_0, \gamma})^2 \Big)  
\\& 
\lesssim_{\bar{s},M} \gamma^{- 2 |k|} \Big((\| \eta \|_{s_0 + \sigma}^{k_0, \gamma})^2 \| h \|_s^2 
+ (\| \eta \|_{s + \sigma}^{k_0, \gamma})^2 \| h \|_{s_0}^2 \Big)\,.
\end{aligned}
\]
Then, the latter chain of inequalities imply that for any $|\beta| \leq \beta_0 - 4 s_0$
\[
\mathfrak{M}_{\langle D \rangle^{M } \partial_\vphi^\beta  {\mathfrak Q}(\vphi)
 \langle D \rangle^{- \mathtt{c}}}(s) 
\lesssim_{\bar s,M} \| \eta \|_{s + \sigma}^{k_0, \gamma}\,, 
\quad \forall s_0 \leq s \leq \bar{s}-\s\,,
\]
and hence, using that ${\mathfrak R}^* = {\mathfrak R} + {\mathfrak R}^* - {\mathfrak R} 
= {\mathfrak R} - {\mathfrak Q} $ and since ${\mathfrak R}$ 
satisfies \eqref{stimefinaliRR}, one has that 
for any $|\beta| \leq \beta_0 - 4 s_0$
\begin{equation}\label{sistemi dinamici merda 4}
\mathfrak{M}_{\langle D \rangle^{M } \partial_\vphi^\beta {\mathfrak R}(\vphi)^* 
\langle D \rangle^{- \mathtt{c}}}(s) \lesssim_{\bar s,M} \| \eta \|_{s + \sigma}^{k_0, \gamma}\,, 
\quad \forall s_0 \leq s \leq \bar{s}-\s\,. 
\end{equation}
The estimates on the Lipschitz variation follow similarly.
By \eqref{sistemi dinamici merda 10}, we write 
\[
\begin{aligned}
G(\eta) & = |D| \tanh(\mathtt h |D|) + {\rm Op}^{BW}(\mathfrak m_1) + {\rm Op}^{BW}(\mathfrak q_0) + {\mathfrak R}(\vphi) 
\\& 
=  
|D| \tanh(\mathtt h |D|) + {\rm Op}^{BW}(\mathfrak m_1) 
+ {\rm Op}^{BW}\Big(\dfrac{\mathfrak q_0 + \overline{\mathfrak q}_0}{2} \Big) 
+ {\rm Op}^{BW}\Big(\dfrac{\mathfrak q_0 - \overline{\mathfrak q}_0}{2} \Big)+ {\mathfrak R}(\vphi) 
\\& 
\stackrel{\eqref{sistemi dinamici merda 11}}{=} 
|D| \tanh(\mathtt h |D|) + {\rm Op}^{BW}(\mathfrak m_1) 
+ {\rm Op}^{BW}\Big(\dfrac{\mathfrak q_0 + \overline{\mathfrak q}_0}{2} \Big) 
+ \dfrac{{\mathfrak R}(\vphi)^* - {\mathfrak R}(\vphi)}{2}+ {\mathfrak R}(\vphi) 
\\& 
%=  
%|D| \tanh(\mathtt h |D|) + {\rm Op}^{BW}(\mathfrak m_1) 
%+ {\rm Op}^{BW}\Big(\dfrac{\mathfrak q_0 + \overline{\mathfrak q}_0}{2} \Big) 
%+ \dfrac{{\mathfrak R}(\vphi)^* + {\mathfrak R}(\vphi)}{2} 
%\\& 
= 
|D| \tanh(\mathtt h |D|) + {\rm Op}^{W}(\mathfrak m_1) 
+ {\rm Op}^{W}\Big(\dfrac{\mathfrak q_0 + \overline{\mathfrak q}_0}{2} \Big)  
\\& \quad 
+ {\rm Op}^{BW}(\mathfrak m_1 + \dfrac{\mathfrak q_0 + \overline{\mathfrak q}_0}{2})  
- {\rm Op}^{W}\Big(\mathfrak m_1 + \dfrac{\mathfrak q_0 + \overline{\mathfrak q}_0}{2} \Big) 
+ \dfrac{{\mathfrak R}(\vphi)^* + {\mathfrak R}(\vphi)}{2}\,.
\end{aligned}
\]
Formula \eqref{pseudoespansione} follows by setting 
and
\[
a_G:=\mathfrak{m}_1=|\x|\chi(\x)\Big(\sqrt{1+|\nabla\eta|^2-\frac{(\nabla\eta\cdot\x)^2}{|\x|^2}}-1\Big)\,, \quad b_{G}:= \dfrac{\mathfrak{q}_0 + \overline{\mathfrak{q}_0}}{2}
\]
which are real symbols and
$$
{\mathcal R}_G(\vphi) :={\rm Op}^{BW}(a_G + b_G)  - {\rm Op}^{W}(a_G + b_G ) 
+  \dfrac{{\mathfrak R}(\vphi)^* + {\mathfrak R}(\vphi)}{2}\,. 
$$
Note that $a_G(\vphi, x, \xi) = a_G(\vphi, x ,- \xi)$. 
The estimates \eqref{stimefinalisimboliDN} follows by explicit computations on  the symbol $a_G$
and by the estimate \eqref{stima mathfrak q0} used 
to estimate $b_{G}$. The estimates on ${\mathcal R}_G$ follow by 
\eqref{stimefinaliRR} on $\mathfrak{R}$ and  \eqref{sistemi dinamici merda 4} on $\mathfrak{R}^*$, 
Lemma \ref{lemma proprieta bony weil} $(ii)$, $(iii)$ and by the estimates 
\eqref{stimefinalisimboliDN} on $a_G$ and $b_G$.  By construction the symbols 
$a_G$ and $b_G$ are momentum preserving. Since also the operator $G(\eta)$ is momentum preserving 
(see Lemma $B.3$ in \cite{FGtrave}) we deduce that also ${\mathcal R}_G$ is momentum preserving.
 This concludes the proof of Theorem \ref{lemma totale dirichlet neumann}.

\section{Degenerate KAM theory}\label{sec:degener}
In this section we analyze the transversality properties of the linear unperturbed frequencies which are fundamental in order to verify the non-resonance conditions that we need. 
We use the degenerate KAM theory, used already in dimension one in \cite{BM20,BBHM,BFM}
%\cite{BBMdeg,BM20,BBHM,BFM} 
in order to analyze the linear frequencies $\omega_{j}(\th)$ defined in \eqref{LinearFreqWW}
of the pure gravity water waves in finite depth.

The key result of this section is the transversality Proposition \ref{lem:transversality}.
We first need some preliminary results.
In order to obtain non-resonance conditions
we need to introduce the concept of non-degeneracy. 
We give the following definition

\begin{defn}
We say that an analytic function $f:[\th_1, \th_2] \to \R^n$ is non-degenerate if for any $c\in\R^n\setminus\{0\}$ the function $f\cdot c$ is not identically zero.
\end{defn} 

\begin{lemma}{\bf (Non-degeneracy).}\label{nondeglemma}
Consider the set of tangential sites $S = \{ \overline \jmath_1, \ldots, \overline \jmath_\nu \}$ 
in \eqref{TangentialSitesWW}
and the frequency vector $\overline{\omega}$ in \eqref{LinearFreqWW}.
Then the following assertions hold true.
\begin{enumerate}
\item The function $\bar{\oo}:[\th_1, \th_2] \to \R^\nu$ is non degenerate
\item Let $j\in\Gamma^* \setminus (S \cup \{ 0 \})$ such that $|j|\neq |\overline{\jmath}_{i}|$ for all $1\le i\le \nu$. 
Then the map  $ \mathtt h \in [\mathtt h_1, \mathtt h_2] \to (\bar{\oo}(\th), \oo_j(\th)) \in \R^{\nu + 1}$ is non degenerate.
\item Let $j, j'\in\Gamma^* \setminus (S \cup \{ 0 \})$ such that $|j|, |j'|\neq |\overline{\jmath}_{i}|$ 
for all $1\le i\le \nu$ and $|j|\neq |j'|$. 

Then the map $ \mathtt h \in [\mathtt h_1, \mathtt h_2] \to (\bar{\oo}(\th), \oo_j(\th), \oo_{j'}(\th)) \in \R^{\nu + 2}$ is non degenerate.
\end{enumerate}
\end{lemma}
\begin{proof}
It follows reasoning as in Lemma $3.2$ in \cite{BBHM}.
\end{proof}
As we already mentioned in the introduction, in order to verify transversality conditions at the unperturbed level (Propositions \ref{lem:transversality}  below), we need the irrationality of the lattice. The key property of the irrational lattice that we will use is the content of the following lemma. 
\begin{lemma}\label{lem:irrazionale} 
Assume Hypothesis \ref{hyp:lattice} and 
let $j, k \in\Gamma^*$ such that $|j|=|k|$. Then one has $j=\pm k$.
\end{lemma}
\begin{proof}
First of all we note that any vector $\xi\in\Gamma^*$ 
can be written as $\xi=  \sum_{i=1}^{d}\widehat{\x}_i {\bf K}_i$
with $\widehat{\x}=(\widehat{\x}_1,\ldots,\widehat{\x}_d)\in\Z^{d}$.
Moreover, setting
\[
\begin{aligned}
A&:=\big(A_{p}^{q}\big)_{p,q=1,\ldots,d}\,,
\qquad
A_{p}^{q}:=  {\bf K}_{p}\cdot {\bf K}_{q}\,.
\end{aligned}
\]
We note that  $|\xi|^2= A\widehat{\x}\cdot \widehat{\x}$.
Hypothesis \ref{hyp:lattice} implies that 
the vector $\vec{P}=\big(A_p^q\big)_{1\leq p\leq q\leq d}\in \R^{d(d+1)/2}$ is irrational. 
Now consider $j,k\in \Gamma^{*}$ with
\[
j= \sum_{i=1}^{d}\widehat{j}_i {\bf K}_i\,,
\qquad
k= \sum_{i=1}^{d}\widehat{k}_i {\bf K}_i\,,
\qquad |j|=|k|\,.
\]
We have that
\[
\begin{aligned}
0&=|j|^{2}-|k|^{2}=A\widehat{j}\cdot \widehat{j}-A\widehat{k}\cdot \widehat{k}
=\sum_{i=1}^{d}A_{i}^{i}(\widehat{j}_i^2-\widehat{k}_i^2)
+\sum_{i=1}^{d}\sum_{i'=1}^{i-1}2A_{i}^{i'}(\widehat{j}_i\widehat{j}_{i'}
-\widehat{k}_i\widehat{k}_{i'})\,.
\end{aligned}
\]
By the irrationality of the vector $\vec{P}$ we deduce that the equation above
implies that one must have
\[
\widehat{j}_{i}^2-\widehat{k}_{i}^2=0\,,
\qquad
\widehat{j}_i \widehat{j}_{i'}
-\widehat{k}_i\widehat{k}_{i'}=0\,,\quad \forall \,1\leq i'< i\leq d\,,
\]
which holds if and only if $\widehat{j}=\pm\widehat{k}$. Then the thesis follows.
\end{proof}

\begin{prop}{\bf (Transversality).}\label{lem:transversality}
Assume Hypothesis \ref{hyp:lattice}, consider the set $S$ in \eqref{TangentialSitesWW},
the frequency vector $\bar{\oo}(\th)$ in \eqref{LinearFreqWW}
and set
$\omega_{j}(\th):=\sqrt{|j|\tanh(\th|j|)}$ for $j\in \Gamma^{*}\setminus S$.
There exist $m_0\in\N$ and $\rho_{0}>0$ such that, 
for any $\th\in [\th_{1}, \th_{2}]$ the following hold:
\begin{enumerate}
\item For all $\ell\in\Z^\nu\setminus\{0\}$ one has
\begin{equation}\label{imp1}
\max_{0\le k\le m_{0}} \{|\partial^{k}_\th \bar{\oo}(\th)\cdot\ell |\} \ge \rho_0 \jap{\ell}\,.
\end{equation} 
\item Let $j\in\Gamma^*\setminus S$ and 
$\ell\in\Z^\nu\setminus\{0\}$ such that $\tV^{T}\ell+ j= 0$. Then 
\begin{equation}\label{imp2}
\max_{0\le k\le m_0}\{|\partial^{k}_\th (\bar{\oo}(\th)\cdot \ell \pm \oo_j(\th)) |\} \ge \rho_0  \jap{\ell}\,.
\end{equation} 
\item Let $j, j'\in\Gamma^*\setminus S$ and 
$\ell\in\Z^\nu\setminus\{0\}$ such that $\tV^{T}\ell+ j - j'= 0$. Then 
\begin{equation}\label{imp3}
\max_{0\le k\le m_0}\{|\partial^{k}_\th (\bar{\oo}(\th)\cdot \ell + \oo_j(\th)- \oo_{j'}(\th))  | \}
\ge \rho_0  \jap{\ell}\,.
\end{equation}
\item Let $j, j'\in\Gamma^*\setminus S$ and 
$\ell\in\Z^\nu\setminus\{0\}$ such that $\tV^{T}\ell+ j + j'= 0$. Then 
\begin{equation}\label{imp4}
\max_{0\le k\le m_0}\{|\partial^{k}_\th (\bar{\oo}(\th)\cdot \ell + \oo_j(\th) + \oo_{j'}(\th))  | \}
\ge \rho_0  \jap{\ell}\,.
\end{equation}
\end{enumerate}
\end{prop}

\begin{proof}
It will be useful in the following to note  that
\begin{equation}\label{asympomeghino}
\omega_{j}(\th) = |j|^\frac{1}{2} ( 1+ r_j(\th))\,,\quad j\in\Gamma^{*}\setminus S\,,
\end{equation}
where $r_j$ is analytic in $\th$ and it satisfies the following estimates: 
 \begin{equation}\label{stima resti freq tanh}
 \partial_{\mathtt h}^n r_j(\th)\lesssim_n e^{-c |j|}\,, 
 \quad 
 \forall \mathtt h \in [\mathtt h_1, \mathtt h_2]\,, \quad \forall n \geq 0\,,
 \end{equation} for some constant $c > 0$.
 
\noindent 
We prove separately \eqref{imp1}-\eqref{imp4}. 

\smallskip
\noindent
{\bf Proof of \eqref{imp1}.}
By contradiction, assume that
for any $n\in\N$ there exist $\th_n\in [\th_{1}, \th_{2}] $ 
and $\ell_n \in\Z^\nu\setminus\{0\}$ such that
\begin{equation}\label{eq:0_abs_m}
\Big| \partial_\th^k \overline{\omega}(\th_{n}) \cdot \frac{\ell_n}{\jap{\ell_n}}  \Big| 
< \frac{1}{\jap{n}} \,, 
\quad \forall\,0\leq k \leq n \, .
\end{equation}
The sequences $(\th_n)_{n\in\N}\subset [\th_{1}, \th_{2}]$ 
and $(\ell_n/\jap{\ell_n})_{n\in\N}\subset \R^\nu\setminus\{0\}$ 
are both bounded. 
By compactness, up to subsequences,
$\th_n\to \bar\th\in [\th_{1}, \th_{2}]$ and $\ell_n/\jap{\ell_n}\rightarrow\bar c\neq 0$. 
Therefore, for any $k\in\N$, passing to the limit
for $n\rightarrow + \infty$ in  \eqref{eq:0_abs_m}, 
we get 
$\partial_\th^k \overline{\omega}(\bar\th)\cdot \bar c = 0$.
By the analyticity of $ \overline{\omega}(\th)$, we deduce 
that the function $ \th \mapsto \overline{\omega}(\th) \cdot \bar c$ 
is identically zero on $[\th_{1}, \th_{2}]$, 
which contradicts Lemma \ref{nondeglemma}-$(1)$, since $ \bar c \neq 0 $.

\smallskip
\noindent
{\bf Proof of \eqref{imp2}.} 
By contradiction, assume that for any $n\in \N$ 
there exist $\ell_n, j_n, \th_n$ such that 
$\tV^{T}\ell_n+j_n=0$ and such that for any $m\le n$ we have that
\begin{equation}\label{4.5}
\left|\partial_\th^m (\bar{\oo} (\th_n) \cdot \frac{\ell_n}{\jap{\ell_n}} 
+ \oo_{j_{n}}(\th_n)) \frac{1}{\jap{\ell_n}}\right|< \frac{1}{\jap{n}}\,.
\end{equation}
We distinguish two cases, namely $\ell_n$ bounded or not bounded.

\noindent
\textbf{Case $\ell_n$ bounded:}\\
From the conservation of momentum also $j_n$ is bounded and $\th_n\in[\th_1, \th_2]$. Therefore, up to subsequences  $\ell_n\to \bar{\ell}, j_n\to \bar{j}, \th_n\to \bar{\th}$.   By passing to the limit in \eqref{4.5}, for any $m\in\N$,
\begin{equation}\label{4.5limite}
|\partial_\th^m (\bar{\oo} (\bar{\th}) \cdot \bar{\ell} + \oo_{\bar{j}}(\bar{\th})) |=0\,.
\end{equation}
Moreover, since $j_n\in \Gamma^*\setminus S$ 
then also the limit $\bar{j}\in \Gamma^*\setminus S$.
If $|\bar{j}|\neq |\bar{\jmath}_i|$ for every 
$1\le i\le \nu$ then \eqref{4.5limite} 
gives a contradiction by Lemma \ref{nondeglemma}-$(2)$.
Then we consider the case in which $|\bar{j}|= |\bar{\jmath}_i|$ 
for a certain $1\le i \le \nu$.
By  Lemma \ref{lem:irrazionale} this is possible only if 
$\bar{j}=\pm \bar{\jmath}_i$.
Since $\bar{j}\in \Gamma^*\setminus S$, then
$\bar{j}= -\bar{\jmath}_i$. Note that the momentum preserving condition implies 
\[
\tV^{T}\bar{\ell} + \bar{j}=0\;\;
\Longleftrightarrow\;\;
\tV^{T}\bar{\ell} - \bar{\jmath}_i=0\;\;
\Longleftrightarrow\;\;
\tV^{T}(\bar{\ell}- e_i)=0
\]

The relation \eqref{4.5limite} becomes 
\begin{equation}\label{4.7}
|\partial_\th^m \bar{\oo} (\bar{\th}) \cdot (\bar{\ell}+ e_i )) |=0\,, \quad \forall m \in \N\,. 
\end{equation}
We actually claim that  $\bar{\ell}+ e_i  \neq 0$. Indeed, if by contradiction $\bar{\ell}+ e_i  =  0$, then the momentum condition implies that  
\[
0 = \tV^{T}\bar{\ell}+\bar{j} \stackrel{\bar \ell = - e_i}{=} -2\bar{\jmath}_i \,,
\]
and so this is a contradiction since all the tangential sites $\bar{\jmath}_i \neq 0$ for any $1 \leq i \leq \nu$. Hence, since $\bar{\ell}+ e_i  \neq 0$, \eqref{4.7} implies that the function $\bar{\oo} (\bar{\th}) \cdot (\bar{\ell}+ e_i) $ 
is identically zero and this contradict Lemma \ref{nondeglemma}-$(1)$.
%If $\ell\neq e_i$ then the function $\bar{\oo} (\bar{\th}) \cdot \bar{\ell}+ e_i $ 
%is identically zero and this contradict Lemma \ref{nondeglemma}-$(1)$. 

\noindent
\textbf{Case $\ell_n$ unbounded:}\\
Up to subsequences, $|\ell_n|\to +\infty$ but $\frac{\ell_n}{\jap{\ell_n}}$ is bounded 
and it converges at $\bar{c}$. 
By the conservation of momentum, 
since $\tV^{T}\ell_n + j_n= 0$ 
then we have that $|j_n|\lesssim |\ell_n|$. 
Hence using 
\[
|\oo_{j_n}|\lesssim |j_{n}|^{1/2}\lesssim |\ell_n|^{1/2}
\]
we get
\begin{equation*}
\frac{\oo_{j_n}}{\jap{\ell_n}}\lesssim \frac{\jap{\ell_n}^{1/2}}{\jap{\ell_n}}
\lesssim \frac{1}{\jap{\ell_n}^{1/2}}\to 0\,.
\end{equation*}
Passing to the limit in \eqref{4.5} we obtain
\begin{equation}
\partial_\th^m (\bar{\oo}(\bar{\th}) \cdot \bar{c})=0\,,
\end{equation}
for any $m\in \N$. This implies that the function 
$\bar{\oo}(\bar{\th}) \cdot \bar{c}$ vanishes in a 
neighborhood of $\th$ and 
this contradicts  Lemma \ref{nondeglemma}-$(1)$.

\smallskip
\noindent
{\bf Proof of \eqref{imp3}.}
By contradiction, assume that for any $n\in\N$ there exist  
$\ell_n, j_n, j'_n, \th_n$ such that $\ell_n \in \Z^\nu \setminus \{ 0 \}$, $j_n, j_n' \in \Gamma^* \setminus S$, $\tV^{T}\ell_n+j_n-j'_n=0$ 
and such that for any $m\le n$ we have that
\begin{equation}\label{4.10}
\left|\partial_\th^m \Big(\bar{\oo} (\th_n) \cdot \frac{\ell_n}{\jap{\ell_n}} 
+ \oo_{j_n}(\th_n) \frac{1}{\jap{\ell_n}}-  \oo_{j_n'}(\th_n) \frac{1}{\jap{\ell_n}} \Big)\right|< \frac{1}{\jap{n}}\,.
\end{equation}
Again, we distinguish two cases. 

\noindent
\textbf{Case $\ell_n$ bounded:}\\
Up to subsequences $\ell_n\to \bar{\ell}$. 
By the conservation of momentum the difference $j_n - j'_n$ is bounded then either both are bounded or both are unbounded.\\
In the first case, assuming that up to subsequences $j_n\to \bar{j}$ and $j'_n\to \bar{j}'$,  
then passing to limit in \eqref{4.10} we obtain
\begin{equation}\label{4.10limite}
\begin{aligned}
|\partial_\th^m (\bar{\oo} (\th) \cdot \bar{\ell} &+ 
\oo_{\bar{j}}(\bar{\th} ) -  \oo_{\bar{j}'}(\bar{\th}) )|=0\,,\qquad \forall m\in\N\\
&\tV^{T}\bar{\ell} + \bar{j}- \bar{j}'=0\,.
\end{aligned}
\end{equation}
If $|\bar{j}|= |\bar{j}'| $  then  the condition \eqref{4.10limite} implies 
$\partial_\th^m (\bar{\oo} (\th)\cdot \bar{\ell}) =0$ 
for every $m\in\N$ and this implies that the function 
$\bar{\oo} (\th)\cdot \bar{\ell}$ is identically zero 
and this contradicts  Lemma \ref{nondeglemma}-$(1)$.

\noindent
Now let us suppose that  $|\bar{j}|\neq |\bar{j}'| $; we discuss three cases:

\noindent
\textbf{Case 1}: $|\bar{j}|, |\bar{j}'| \neq |\bar{\jmath}_i|$ for all $1\le i\le \nu$. 
In this case \eqref{4.10limite} implies 
$\bar{\oo} (\th) \cdot \bar{\ell} + \oo_{\bar{j}}(\bar{\th} ) 
-  \oo_{\bar{j}'}(\bar{\th})$ 
is identically zero and this contradicts Lemma \ref{nondeglemma}-$(3)$.

\noindent
\textbf{Case 2}: $|\bar{j}|= |\bar{\jmath}_k|$ for some $1\le k\le \nu$ 
and $|\bar{j}'| \neq |\bar{\jmath}_i|$ for all $1\le i\le \nu$. 
Since $\bar{j}\in \Gamma^*\setminus S$ 
then $|\bar{j}|= |\bar{\jmath}_k|$ implies, by Lemma \ref{lem:irrazionale}, 
that $\bar{j}=-\bar{\jmath}_k$. 
In this case we have that condition \eqref{4.10limite} reads
\[
\partial_\th^m\big(\bar{\oo}(\bar{\th} )\cdot (\bar{\ell} + e_k) -\oo_{\bar{j}'} (\bar{\th} )\big)=0
\]
for any $m$. This implies that 
$\bar{\oo}(\bar{\th} )\cdot (\bar{\ell} + e_k) -\oo_{\bar{j}'} (\bar{\th} )$ 
is identically zero in a neighborhood of $\bar{\th}$ 
and this contradicts Lemma \ref{nondeglemma}-$(2)$.

\noindent
\textbf{Case 3:} $|\bar{j}|= |\bar{\jmath}_k|$ and $|\bar{j}'| = |\bar{\jmath}_i|$  
for some $1\le k< i \le \nu$. 
Since $\bar{j}, \bar{j'}\in\Gamma^* \setminus S$ 
then Lemma \ref{lem:irrazionale} implies 
$ \bar{j}= - \bar{\jmath}_k$ and $ \bar{j}'= - \bar{\jmath}_i$ .
In this case we have that 
\[
\partial_\th^m(\bar{\oo}(\bar{\th} )\cdot (\bar{\ell} + e_k- e_i) (\bar{\th} ))=0
\]
for any $m$. We claim that $\bar{\ell} + e_k- e_i\neq 0$. Indeed, if by contradiction $\bar{\ell} + e_k- e_i =  0$, then, by the conservation of momentum $\tV^{T}\bar{\ell}+ \bar{j} - \bar{j}' =0$, we have that 
\begin{equation*}
\begin{aligned}
0&=\tV^{T}\bar{\ell}+ \bar{j} - \bar{j}' = \tV^{T}\bar{\ell} - \bar{\jmath}_k +  \bar{\jmath}_i  = \tV^{T}(e_i - e_k) - \bar{\jmath}_k +  \bar{\jmath}_i  = 2(\bar{\jmath}_i  - \bar{\jmath}_k )
\end{aligned}
\end{equation*}
and this is an absurd since $\bar{\jmath}_i \neq \bar{\jmath}_k$
by the choice of the tangential sites (see \eqref{TangentialSitesWW}).
Thus, the fact that  $\bar{\ell} + e_k- e_i\neq 0$ implies that 
$\bar{\oo}(\bar{\th} )\cdot (\bar{\ell} + e_k- e_i)$ 
is identically zero in a neighborhood of $\bar{\th}$ 
and this contradicts Lemma \ref{nondeglemma}-$(1)$.

Let us now suppose that $j_{n}, j'_{n}$ are both unbounded.
Using the expansion \eqref{asympomeghino} and the estimate
\eqref{stima resti freq tanh}
we note that
\begin{equation}\label{diff frequenze}
\begin{aligned}
|\partial_{\th}^{m}\left(\oo_{j_n}(\th) - \oo_{j'_n}(\th)\right)|&\leq
|\sqrt{|j_n|}- \sqrt{|j'_n|}| + |\partial_{\th}^{m}r_{j_n}(\th)| + |\partial_{\th}^{m}r_{j'_n}(\th)|
\\&\leq
\frac{||j_n|-|j'_n||}{\sqrt{|j_n|} + \sqrt{|j'_n|}}+ |\partial_{\th}^{m}r_{j_n}(\th)| + |\partial_{\th}^{m}r_{j'_n}(\th)|
\\
& \stackrel{}{\lesssim}
\frac{|j_n - j'_n|}{\sqrt{|j_n|}+\sqrt{|j'_n|}} +  |\partial_{\th}^{m}r_{j_n}(\th)| + |\partial_{\th}^{m}r_{j'_n}(\th)|\,,
\end{aligned}
\end{equation}
for $n\in\N$ large enough. Then, recalling that $|j_n-j'_n|\lesssim |\ell_n|$,
we obtain
\begin{equation}\label{diff frequenze 0}
\begin{aligned}
 \frac{|\partial_{\th}^{m} \left( \oo_{j_n}(\th) - \oo_{j'_n}(\th)\right)|}{\langle\ell_{n}\rangle} & \lesssim
\frac{1}{\langle\ell_{n}\rangle} \frac{|j_n - j'_n|}{\sqrt{|j_n|}+\sqrt{|j'_n|}} + \frac{1}{\langle\ell_{n}\rangle} ( |\partial_{\th}^{m}r_{j_n}(\th)| + |\partial_{\th}^{m}r_{j'_n}(\th)| ) \\
&{\lesssim} \frac{1}{\sqrt{|j_n|}+\sqrt{|j'_n|}} +  |\partial_{\th}^{m}r_{j_n}(\th)| + |\partial_{\th}^{m}r_{j'_n}(\th)|   \to 0
\end{aligned}
\end{equation}
as $n\to\infty$ (recall that by \eqref{stima resti freq tanh}, $|\partial_{\th}^{m}r_{j_n}(\th)| \,,\, |\partial_{\th}^{m}r_{j'_n}(\th) \to 0$ since $j_n, j_n' \to + \infty$ ).

\noindent
\textbf{Case $\ell_n$ unbounded:}\\
Up to subsequences, $|\ell_n|\to +\infty$ but $\frac{\ell_n}{\jap{\ell_n}}$ 
is bounded and it converges at $\bar{c}$. 
By the conservation of momentum, since 
$\tV^{T}\ell_n + j_n - j'_n= 0$,  we have that $|j_n-j'_n|\lesssim |\ell_n|$.

We now discuss three cases.

\noindent
{\bf Case 1:} $|j_n|, |j'_n|\le C$ for some constant $C>0$.
Recalling the asymptotic in \eqref{asympomeghino}, we get
\[
|\partial_{\mathtt h}^m(\oo_{j_n} - \oo_{j'_n})| \langle\ell_{n}\rangle^{-1}\lesssim \langle\ell_{n}\rangle^{-1}  \to0
\]
 as $n\to\infty$.
Then, passing to the limit in \eqref{4.10}, 
we get that $\partial_{\th}^n (\overline{\oo}(\overline{\th}) \cdot \bar{c})=0$
for any $m\in \N$, which implies that 
$\overline{\oo}(\overline{\th}) \cdot \bar{c}$
 is identically zero in a neighborhood of $\bar{h}$ and this contradicts 
 Lemma \ref{lem:irrazionale}-$(1)$.

\noindent
{\bf Case 2:}  $|j_n|\to +\infty, |j'_n|\le C$. By the momentum condition
we have $|j_n-j'_n|\lesssim |\ell_n|$ and hence $|j_n|\lesssim \jap{\ell_n}$.
Therefore, recalling \eqref{asympomeghino}, we deduce
that
\[
\frac{|\partial_{\th}^{m}\oo_{j_n}|+ |\partial_{\th}^{m}\oo_{j'_n}|}{\jap{\ell_n}} 
\lesssim 
\frac{|j_n|^{\frac{1}{2}}+C}{\jap{\ell_n}}
\lesssim \frac{\jap{\ell_n}^{\frac{1}{2}}}{\jap{\ell_n}}\to 0\,,
\]
(up to subsequences) as $n\to\infty$. Therefore,
we pass to the limit in \eqref{4.10} and we conclude as done in the previous case.

\noindent
{\bf Case 3:} $|j_n|, |j'_n|\to \infty$.
Arguing as in \eqref{diff frequenze}, \eqref{diff frequenze 0}, we have again that 
\[
\frac{|\partial_{\th}^{m}\left(\oo_{j_n}(\th) - \oo_{j'_n}(\th)\right) |}{\langle\ell_{n}\rangle}\lesssim
\frac{1}{\sqrt{|j_n|}+\sqrt{|j'_n|}} +  |\partial_{\th}^{m}r_{j_n}(\th)| + |\partial_{\th}^{m}r_{j'_n}(\th)|   \to 0
\]
as $n\to\infty$. Again we pass to the limit in \eqref{4.10} 
and we conclude as done in the previous cases.

\noindent
{\bf Proof of \eqref{imp4}.} It follows as \eqref{imp3} (it is actually simpler) hence we omit it.
\end{proof}

%\section{Non-linear functional setting}\label{sec:nonlinsetting}

\section{Nash-Moser theorem of hypothetical conjugation}\label{sec:nash}
For $0<\e\ll1$ we consider the rescaling 
$(\eta,\psi)\mapsto (\e\eta,\e\psi)$.
The equations
\eqref{eq:113}-\eqref{HS} trasform into
\begin{equation}\label{WW-riscalato}
\pa_t u = J {\bf \Omega} u + \e X_{P_\e} (u)  
\end{equation}
where $ J {\bf \Omega}  $ is the linearized Hamiltonian vector  field in \eqref{eq:113linear} and
\begin{equation}\label{campo hamiltoniano X P epsilon}
\begin{aligned}
& X_{P_\e} (u) 
 %:= X_{P_\e}(\kappa, u)  
 :=   
 \left(
 \begin{matrix} 
 \e^{- 1}(G(\e \eta)- G(0))\psi \\
- \frac{1}{2} |\nabla\psi|^2  
+\frac{1}{2} \frac{ \big(G(\e \eta)\psi + \e \nabla\eta\cdot\nabla\psi \bigr)^2}{1+ \e^2 |\nabla\eta|^2}
\end{matrix} 
\right)
\end{aligned}
\end{equation}
and the Hamiltonian $H$ in \eqref{HamiltonianWW} 
% (with $g=1$ and depth $\th$)
transforms into
\[
\mathcal{H}_{\e}(\eta,\psi)=\e^{-2}H(\e\eta,\e\psi)=H_{L}(\eta,\psi)+\e P_{\e}(\eta,\psi)
\]
where $H_{L}$ is the quadratic Hamiltonian in \eqref{HamLinearReal}
and
\begin{equation}\label{P.epsilon}
P_{\e}(\eta,\psi):=
\frac{1}{\e}\int_{\mathbb{T}_{\Gamma}^{2}}\psi\, \big(G(\e\eta )-G(0)\big)\psi\, dx\,.
\end{equation}
Recalling also the action angle variables introduced in \eqref{actionanglevar}
we define
\[
H_{\e}:=\mathcal{H}_{\e}\circ A=\e^{-2}H\circ\e A
\]
which reads, up to a constant,
\begin{equation}\label{hamsez4}
H_{\e}=\mathcal{N}+\e P\,,\qquad \mathcal{N}:=H_{L}\circ A=
\overline{\omega}(\th)\cdot I+\frac{1}{2}({\bf \Omega} z, z)_{L^{2}}\,,
\;\;\;P:=P_{\e}\circ A\,,
\end{equation}
where $\overline{\omega}(\th)$ is defined in \eqref{LinearFreqWW}
and ${\bf \Omega}$ in \eqref{HamLinearReal}.
We look for an embedded invariant torus (recall \eqref{decomposition})
\[
i \,:\, \mathbb{T}^{\nu} \,\to \,\mathbb{T}^{\nu}\times\mathbb{R}^{\nu}\times H_{S}^{\perp}\,,
\qquad \vphi\mapsto i(\vphi)=(\theta(\vphi),I(\vphi),z(\vphi)) 
\]
of the Hamiltonian vector field
$X_{H_{\e}}=(\pa_{I}H_\e, -\pa_{\theta}H_\e, \mathbb{P}_{S}^{\perp}J\nabla_{z}H_{\e})$
filled by traveling quasi-periodic solutions with Diophantine frequency
$\omega\in \R^{\nu}$.
We remark that the embedding $i(\vphi)$ in the original coordinates
(see \eqref{actionanglevar})
reads
\[
(\eta(t,x),\psi(t,x))=U(\omega t,x)\,,\qquad U(\vphi,x)=\e A(\theta(\vphi),I(\vphi),z(\vphi))\,,
\qquad\vphi\in \T^{\nu}\,. 
\]
Instead of looking directly for quasi-periodic solutions of  $ X_{H_\varepsilon} $
we look for quasi-periodic solutions 
of the family of modified Hamiltonians, 
where $\alpha\in \R^\nu$ are additional parameters,
\begin{equation}\label{Halpha}
H_\alpha := \cN_\alpha + \varepsilon P \,, 
\quad \cN_\alpha:= \alpha \cdot I + \tfrac12 \left( {\bf \Omega}z,z \right)_{L^2} \, .
\end{equation}
Actually, we construct the invariant tori of $H_{\alpha, \varrho}$ which is defined as
\begin{equation}\label{HHalphaMumu}
H_{\alpha, \varrho}:=H_{\alpha}  + \varrho \cdot \theta\,. 
\end{equation}
We consider   the nonlinear operator
\begin{align}\label{F_op}
\cF(i,\alpha, \varrho) & := \cF(\omega,\th,\varepsilon;i,\alpha, \varrho) := \omega\cdot\pa_\vphi i(\vphi) - X_{H_{\alpha, \varrho}}(i(\vphi)) \notag \\
& = \begin{pmatrix}
\omega\cdot \pa_\vphi \theta(\vphi)  - \alpha -\varepsilon \pa_I P(i(\vphi)) 
\\
\omega\cdot \pa_\vphi I(\vphi) + \varrho+ \varepsilon \pa_\theta P(i(\vphi)) 
\\
\omega\cdot \pa_\vphi z(\vphi)   - \, \mathbb{P}_{S}^{\perp}J ( {\bf \Omega} z(\vphi) +\varepsilon \nabla_{z} P(i(\vphi)) )  
\end{pmatrix} \, . 
\end{align}
If $\cF(i,\alpha, \varrho)=0$, then 
the embedding $\vphi\mapsto i(\vphi)$ is an invariant torus for the Hamiltonian vector field $X_{H_{\alpha, \varrho}}$, filled with quasi-periodic solutions with frequency $\omega$.

Each Hamiltonian $H_{\alpha, \varrho}$ in \eqref{Halpha} 
is invariant under the translations $\tau_\vs$, $\vs \in \R^{d}$, 
defined in  \eqref{def:vec.tau}:
\begin{equation}
\label{Halpha.symm}
H_{\alpha, \varrho}\circ \tau_\vs = H_{\alpha, \varrho} \, ,  \quad \forall\, \vs \in \R^d \, . 
\end{equation}
We look for a traveling torus embedding 
$  i (\vphi) = $ $ ( \theta(\vphi), I(\vphi), z(\vphi)) $, 
according to Definition \ref{def:travelembedd},
namely 
satisfying  (see \eqref{cond:travelemb}-\eqref{eq:travelemb})
\begin{equation}\label{RTTT}
\theta(\vphi)-\tV\vs=\theta(\vphi-\mathtt{V}\vs)\,,\quad I(\vphi)=I(\vphi-\mathtt{V}\vs)\,,
\quad
z(\vphi,x+\vs)=z(\vphi-\tV \vs,x) \,, \quad \forall \,\vs \in \R^d \, . 
\end{equation}
Note that, by \eqref{F_op} and \eqref{Halpha.symm}, 
the operator $ \cF$ is momentum preserving.
The norm of the periodic components of the embedded torus
\begin{equation}\label{ICal}
\cI (\vphi):= i(\vphi)-(\vphi,0,0) := \left( \Theta(\vphi), I(\vphi), w(\vphi) \right)\,, 
\quad \Theta(\vphi):= \theta(\vphi)-\vphi\,,
\end{equation}
is 
\begin{equation}\label{norma.ical}
\norm{ \cI }_s^{k_0,\gamma} := 
\norm{\Theta}_{H_\vphi^s}^{k_0,\gamma} 
+ \norm{I}_{H_\vphi^s}^{k_0,\gamma} + \norm{w}_s^{k_0,\gamma} \,,
\end{equation}
where 
$\norm{w}_s^{k_0,\gamma} =\norm{\eta}_s^{k_0,\gamma} 
+\norm{\psi}_s^{k_0,\gamma} $.
We define
\begin{equation}\label{k0_def}
k_0:= m_{0} + 2
\end{equation}
and $m_0 \in \N $ is the index of non-degeneracy provided 
by Proposition \ref{lem:transversality}, which only depends 
on the linear unperturbed frequencies. 
We will often omit to write the dependence 
of the various constants with respect to $k_0$, which
is considered as an absolute constant.
We look for quasi-periodic solutions of frequency $\omega$
belonging to a $\delta$-neighborhood (independent of $\varepsilon$)
\begin{equation}\label{definition mathtt Omega}
\mathtt{\Omega} := \big\{ \omega \in \R^\nu \ : \
\dist \big( \omega, \overline{\oo}[\th_1,\th_2] \big) < \delta \big\} \,, 
\quad \delta >0 \,,
\end{equation}
of the curve
$\overline{\oo}[\th_1,\th_2]$ defined by \eqref{LinearFreqWW}.
The next theorem
provides, for $ \varepsilon $ small enough, 
 a solution $ (i_\infty, \alpha_\infty, \varrho_\infty)(\omega,\th; 	\varepsilon) $ of the nonlinear operator
 $ {\cF}(\omega,\th,\varepsilon; i,\alpha, \varrho) = 0 $ for all the values of
 $ (\omega, \th) $ in the Cantor like set   $ \mathcal{G}_\infty(\gamma)$  suitably constructed in the theorem  below.
Actually, it will be shown that $\varrho_\infty=0$, 
implying that we find an invariant torus $i_\infty$ of the Hamiltonian $H_{\alpha_\infty}$.

\begin{thm}  {\bf (Theorem of hypothetical conjugation)}\label{NMT}
Consider the tangential sites $S$ in \eqref{TangentialSitesWW} (recall also Hyp. \ref{hyp:lattice}) 
and let $\tau\geq1$.
There exist positive constants ${\rm a_0},\bar{\sigma}, C,\kappa_{1}$ such that, 
 for any 
$\mathfrak{s} \geq \bar\sigma$, there exists $\e_0=\e_0(\mathfrak{s})\in (0,1)$ such that 
for all $\e\in(0,\e_0)$, any $\gamma=\e^{\rm a}$, ${\rm a}\in(0,{\rm a}_0)$
the following holds. 
There exist
\begin{enumerate}
\item 
a $k_0$-times differentiable function of the form 
$ \alpha_\infty :  \, \mathtt{\Omega} \times [\th_1,\th_2] \mapsto \R^\nu $, 
\begin{align}\label{alpha_infty}
& \alpha_\infty(\omega,\th) := \omega + r_\varepsilon(\omega,\th)  
\quad \text{ with } \quad  |r_\varepsilon|^{k_0,\gamma} \leq C \varepsilon \gamma^{-1} \, ;
\end{align}
moreover the map 
$\alpha_\infty(\,\cdot\,,\th)$ from $\mathtt{\Omega}$ into its image $\alpha_\infty(\mathtt{\Omega},\th)$ 
is invertible with inverse of the form
\begin{equation}\label{inv_alpha100}
\begin{aligned}
& \beta = \alpha_\infty(\omega,\th)  \quad \Leftrightarrow \quad
   \omega = \alpha_\infty^{-1}(\beta,\th) = \beta+\breve{r}_\varepsilon(\beta,\th) \, , \\&
 \abs{ \breve{r}_\varepsilon }^{k_0,\gamma} \leq C\varepsilon\gamma^{-1} \,;
\end{aligned}
\end{equation}

\item	
a family of embedded traveling tori $i_\infty (\vphi) $ (cfr. \eqref{RTTT}), 
defined for all $(\omega,\th)\in\mathtt{\Omega} \times[\th_1,\th_2] $, satisfying
\begin{equation}\label{i.infty.est}
\| i_\infty (\vphi) -(\vphi,0,0) \|_{\mathfrak{s}}^{k_0,\gamma} \leq C \varepsilon\gamma^{-1} \, ;
\end{equation}

\item a Cantor-like set $\mathcal{G}_{\infty}(\gamma)\subseteq \mathtt{\Omega}\times[\th_1,\th_2]$
satisfying 
\begin{equation}\label{megaMisura}
\big| \big\{ \th \in [\th_1, \th_2] 
\ : \  \big( \alpha_\infty^{-1}( \overline{\omega}(\th),\th ),\th \big) 
\in {\mathcal G}_\infty(\e^{\rm a}) \big\}\big|\to \th_2-\th_1\,,\quad {\rm as}\;\;\; \e\to0\,,
\end{equation}
\end{enumerate}
such that for any $(\omega,\th)\in\mathcal{G}_{\infty}(\gamma)$ 
the function $i_\infty(\vphi):= i_\infty(\omega,\th,\varepsilon;\vphi)$ is a 
	solution of 
	\[
	\cF(\omega,\th,\varepsilon; (i_\infty,\alpha_\infty, 0)(\omega,\th))=0\,.
	\] 
As a consequence, the embedded torus $\vphi\mapsto i_\infty(\vphi)$ is invariant for 
the Hamiltonian vector field $X_{H_{\alpha_\infty(\omega,\th)}}$ as it is filled by  
quasi-periodic traveling wave solutions  with frequency $\omega$.	
\end{thm}

%The next goal is to deduce Theorem \ref{thm:main0} from Theorem \ref{NMT}. 
The next sections are devote to the proof of the result above. In sections 
\ref{subsec:measest}-\ref{sec:linearstability} we will show how to deduce 
Theorem \ref{thm:main0} from Theorem \ref{NMT}.

\section{Approximate inverse}\label{sec:approx_inv}
In order to implement a convergent Nash-Moser scheme 
that leads to a solution of $\cF(i,\alpha, \varrho)=0$, where 
$ \cF (i, \alpha, \varrho) $ is the nonlinear operator  defined in \eqref{F_op},  
we construct an \emph{approximate right inverse} 
of the linearized operator
\begin{equation*}
\di_{i,\alpha, \varrho}\cF(i_0,\alpha_0, \varrho_0)[\whi,\wh\alpha, \wh\varrho] = 
\omega\cdot \pa_\vphi \whi 
- \di_i X_{H_\alpha}\left( i_0(\vphi) \right)[\whi] - \left(\wh\alpha,0,0\right) 
+ (0, \wh\varrho, 0) \,.
\end{equation*}
Note that $\di_{i,\alpha, \varrho}\cF(i_0,\alpha_0, \varrho_0)=\di_{i,\alpha,\varrho}\cF(i_0)$ 
is independent of $\alpha_0$.
We assume that the torus 
$ i_0 (\vphi) = ( \theta_0 (\vphi), I_0 (\vphi), z_0 (\vphi)) $ 
is traveling,  according to  \eqref{RTTT} and along this section, we assume the following hypothesis, 
which is verified by the approximate solutions obtained 
at each step of the Nash-Moser Theorem in section \ref{sec:NaM}. From now one we fix arbitrary parameters $s_0, \bar s$ satisfying
\begin{equation}\label{Sone}
\bar s \gg  s_0 \geq \nu + d + 3 
\end{equation}
and $k_0$ as in \eqref{k0_def}.
\begin{itemize}
\item {\sc ANSATZ.} The map 
$(\omega,\mathtt h)\mapsto \cI_0(\omega,\th) = i_0(\vphi;\omega,\th)- (\vphi,0,0)$ 
is $k_0$-times differentiable with respect to the parameters 
$(\omega,\th)\in \mathtt \Omega\times [\th_1,\th_2]$ where $\mathtt \Omega$ is defined in \eqref{definition mathtt Omega} and, 
for some $\s:=\s(\tau,\nu, k_0)>0$, $\gamma\in (0,1)$,
\begin{equation}\label{ansatz}
\norm{\cI_0}_{s_0+\s}^{k_0,\gamma} 
+ \abs{ \alpha_0-\omega }^{k_0,\gamma} \leq C \varepsilon \gamma^{-1} \,.
\end{equation} 
\end{itemize}

In the sequel we shall assume  the smallness condition,  
for some $\tk := \tk (\tau,\nu, k_0)>0$, 
\begin{equation}\label{epsilonpiccoloKK}
\varepsilon\gamma^{-\tk} \ll 1 \,.
\end{equation}
\noindent
We closely follow the  strategy presented in \cite{BB} 
and implemented for the water waves 
equations in \cite{BM20,BBHM, BFM}. 
As shown in \cite{BFM}, this construction preserves
the momentum preserving properties needed 
for the search of traveling waves and the estimates are very similar.

\smallskip
First of all, we state  tame estimates for the composition 
operator induced by the Hamiltonian vector field 
$X_{P}= ( \pa_I P , - \pa_\theta P, \mathbb{P}_{S}^\perp J \nabla_{z} P )$ 
in \eqref{F_op}.

We first estimate the composition operator induced by $ v(\teta, y) $  defined in \eqref{definizione.v}.
Since the functions $ I_j  \mapsto \sqrt{\xi_j + I_j}  $, $\theta \mapsto {\rm cos}(\theta)$, 
$\theta \mapsto {\rm sin}(\theta)$ 
are analytic  for $|I| \leq r $ small, 
the composition Lemma \ref{Moser norme pesate} implies that, for all  $ \Theta, y \in H^s(\T^\nu, \R^\nu )$,  
$   \| \Theta \|_{s_0},  \| y \|_{s_0} \leq r $, setting $\theta(\varphi) := \varphi + \Theta (\varphi)$, 
\begin{equation}\label{stima k0 gamma v per il campo}
\| \partial_\theta^\alpha \partial_I^\beta v(\theta(\cdot)\,, 
I(\cdot)) \|_s^{k_0, \gamma} \leq_s 1 + \| \cI\|_{s}^{k_0, \gamma}\,, 
\qquad \forall \alpha, \beta \in \N^\nu\,, \ \  |\alpha| + |\beta| \leq 3\,.
\end{equation}

\begin{lemma}{\bf (Estimates of the perturbation $P$)} \label{XP_est}
There exist 
$\wh\s=\wh\s(k_0)$ large enough such that if  $\cI(\vphi)$, defined in \eqref{ICal}, 
satisfies, 
$
\norm{ \cI }_{s_{0}+\wh\s}^{k_0,\gamma}\leq 1$ then, 
for any $ s \geq s_0 
$, 
\begin{equation}\label{stima.Xp}
\norm{ X_{P}(i) }_{s}^{k_0,\gamma} \lesssim_s 1 
+ \norm{ \cI }_{s+\wh\s}^{k_0,\gamma}\,,
\end{equation} 
and, for all $\whi:= (\wh\theta,\whI,\whz)$,
\begin{align}
\|\di_i X_{P}(i)[\whi] \|_{s}^{k_0,\gamma} &\lesssim_s 
\|\whi\|_{s + \wh\s}^{k_0,\gamma} 
+ \| \cI \|_{s+\wh\s}^{k_0,\gamma}
\| \whi\|_{s_0 +  \wh\s }^{k_0,\gamma} \,, \label{stima.dXp}
\\
\norm{ \di_i^2 X_{P}(i)[\whi,\whi] }_{s}^{k_0,\gamma} 
&\lesssim_s 
\|\whi\|_{s + \wh\s}^{k_0,\gamma}\|\whi\|_{s_0 +  \wh\s}^{k_0,\gamma} 
+ \| \cI \|_{s+\wh\s}^{k_0,\gamma} ( \|\whi\|_{s_0+  \wh\s}^{k_0,\gamma} )^2 \,.\label{stima.dXpquadro}
\end{align}
\end{lemma}

\begin{proof}
By the definition \eqref{hamsez4}, $P = P_\e \circ A$, where $A$ is defined in \eqref{actionanglevar} and $P_\e$ is defined in \eqref{P.epsilon}. Hence 
\begin{equation}\label{campo hamiltoniano perturbazione P}
\begin{aligned}
& X_P =   \begin{pmatrix} [\partial_I v(\theta, I)]^T \nabla P_\e(A(\theta, I, z)) \\
 - [\partial_\theta v(\theta, I)]^T \nabla P_\e(A(\theta, I, z)) \\
  \mathbb{P}_{{ S}}^\bot J \nabla P_\e(A(\theta, I, z))  
  \end{pmatrix}\,
\end{aligned}
\end{equation}
where $\mathbb{P}_{ S}^\bot$ is the $L^2$-projector on the space $H_{{ S}}^\bot$ defined in \eqref{decomposition}.
Now  $ \nabla P_\e = - J X_{P_\e} $ %(see \eqref{WW-riscalato})
 where 
$  X_{P_\e} $ is the explicit Hamiltonian vector field in  \eqref{campo hamiltoniano X P epsilon}. 
Since 
\[
\begin{aligned}
\| \eta\|_{s_0 + \wh\s}^{k_0, \gamma} \leq \e \| A(\theta(\cdot), I(\cdot), z(\cdot, \cdot)) \|_{s_0 + \wh\s}^{k_0, \gamma} 
& \leq C(s_0) \e (1 + \| \cI \|_{ s_0 + \wh\s}^{k_0, \gamma})\,,
\end{aligned}
\]
then the  smallness condition of Lemma \ref{stime tame dirichlet neumann} is fulfilled
for $ \e $ small (w.r.t. $s\geq s_0$).
Thus, by the tame estimate
 \eqref{stima tame dirichlet neumann} for the Dirichlet-Neumann operator, 
 the interpolation inequality \eqref{p1-pr},  and \eqref{stima k0 gamma v per il campo},  we get 
$$
\begin{aligned}
\| \nabla P_\e(A(\theta(\cdot), I(\cdot), z(\cdot, \cdot))) \|_s^{k_0, \gamma} 
& \leq_s  \| A(\theta(\cdot), I(\cdot), z(\cdot, \cdot))\|_{s + \wh\s}^{k_0, \gamma}  \leq_s 1 + \| \cI\|_{s + \wh\s}^{k_0, \gamma}\, .
\end{aligned}
$$ 
Hence  \eqref{stima.Xp} % for $X_P$ in \eqref{campo hamiltoniano perturbazione P} 
follows by \eqref{campo hamiltoniano perturbazione P},
 interpolation and  \eqref{stima k0 gamma v per il campo}. %  on $v(\theta, y)$.
  The estimates \eqref{stima.dXp}, \eqref{stima.dXpquadro} for 
  ${\rm d}_i X_P$ and ${\rm d}_i^2 X_P$ follow
 by differentiating the expression of $X_P$ in \eqref{campo hamiltoniano perturbazione P} 
 and applying the estimates \eqref{stima tame derivata dirichlet neumann}, 
 \eqref{stima tame derivata seconda dirichlet neumann} 
 on the Dirichlet-Neumann operator, the estimate \eqref{stima k0 gamma v per il campo} on 
 $v(\theta, y)$ and using the interpolation inequality \eqref{p1-pr}. 
\end{proof}

We first modify the approximate torus $i_0 (\vphi) $ to obtain a nearby 
isotropic torus $i_\delta (\vphi) $, namely such that the pull-back 1-form  $i_\delta^*\Lambda $  is closed, 
where $\Lambda$ is the Liouville 1-form defined in 
\eqref{liouville}.  
Consider the pull-back $ 1$-form 
\begin{align}\label{ak}
i_0^*\Lambda & = \sum_{k=1}^{\nu} a_k(\vphi) \di \vphi_k \, , \quad
a_k(\vphi) := -\big( [ \pa_\vphi \theta_0(\vphi) ]^T I_0(\vphi) \big)_k +\tfrac12 
\big( J^{-1} z_0(\vphi), \pa_{\vphi_k} z_0(\vphi) \big)_{L^2} \, , 
\end{align} 
and define $A_{kj}(\vphi)  := \pa_{\vphi_k} a_j(\vphi) - \pa_{\vphi_j}a_k(\vphi) $.
Let us define the ``error function''
\[
Z(\vphi)   
:= \cF(i_0,\alpha_0, \varrho_0)(\vphi) = \omega\cdot \pa_\varphi i_0(\vphi) 
- X_{H_{\alpha_0, \varrho_0}}(i_0(\vphi)) \,,
\]
which measure how the embedding $i_0$ is close to 
being invariant for $H_{\alpha,\varrho}$ in \eqref{HHalphaMumu}.
The next lemma, which is a consequence of Lemma $6.1$ in \cite{KdVAut}, 
proves that if $i_0$ is a solution of the equation 
$\mathcal{F}=0$ (see \eqref{F_op}), 
then the parameter $ \varrho_0$ has to be naught, 
hence the embedded torus $i_0$ 
supports a quasi-periodic solution of the 
``original'' system with Hamiltonian  $H_{\alpha_0}$. 

\begin{lemma}{(Lemma $6.1$ in \cite{KdVAut})}\label{Lemma6.1DP}
We have
$\lvert  \varrho_0 \rvert^{\gamma, \calO_0}\lesssim \lVert Z \rVert_{s_0}^{\gamma, \calO_0}$.
In particular, if $\mathcal{F}(i_0, \alpha_0,\varrho_0)=0$ then 
$\varrho_0=0$ and the torus 
$i_0(\varphi)$ is invariant for the vector field $X_{H_{\alpha_0}}$.
\end{lemma}

The next lemma follows as in Lemma 5.3 in  \cite{BBHM} and 
Lemma 6.2 in \cite{BFM}. 
\begin{lemma} {\bf (Isotropic torus).} \label{torus_iso}
The torus $i_\delta(\vphi):= ( \theta_0(\vphi),I_\delta(\vphi),z_0(\vphi) )$, 
defined by
\[
I_\delta(\vphi):= I_0(\vphi) + [ \pa_\vphi \theta_0(\vphi) ]^{-T}\rho(\vphi) \,, 
\quad \rho = (\rho_j)_{j=1, \ldots,\nu} \, , 
\quad 
\rho_j(\vphi):= \Delta_\vphi^{-1} \sum_{k=1}^{\nu}\pa_{\vphi_k}A_{kj}(\vphi)\,,
\]
is isotropic. Moreover, there is $\sigma:= \sigma(\nu,\tau)$ such that, for all 
$ s \geq s_0 $, 
\begin{align}
\norm{ I_\delta-I_0 }_s^{k_0,\gamma} &\lesssim_s 
\norm{\cI_0}_{s+1}^{k_0,\gamma}     \label{ebb1} \, , 
\quad 
\norm{ I_\delta-I_0 }_s^{k_0,\gamma}  
\lesssim_s \gamma^{-1}
\big( \norm{Z}_{s+\sigma}^{k_0,\gamma} 
+\norm{Z}_{s_0+\sigma}^{k_0,\gamma} 
\norm{ \cI_0 }_{s+\sigma}^{k_0,\gamma} \big) \,, 
\\
\norm{ \cF(i_\delta,\alpha_0, \varrho_0) }_s^{k_0,\gamma} &\lesssim_s  
\norm{Z}_{s+\sigma}^{k_0,\gamma} +\norm{Z}_{s_0+\sigma}^{k_0,\gamma} 
\norm{ \cI_0 }_{s+\sigma}^{k_0,\gamma} \, , 
\quad 
\|\di_i(i_\delta)[\whi] \|_{s_1} 
\lesssim_{s_1} \| \whi \|_{s_1+1}  \, ,   \label{ebb3} 
\end{align}
for  $  s_1 \leq s_0 + \s $ (cfr. \eqref{ansatz}).
Furthermore  $i_\delta(\vphi)$
is a traveling torus,  cfr.  \eqref{RTTT}. 
\end{lemma}

In view of Lemma \ref{torus_iso}, we
first find an  approximate inverse of the linearized operator 
$\di_{i,\alpha, \varrho}\cF(i_\delta)$. 
We introduce the  diffeomorphism 
$G_\delta:(\phi,y,\tw) \rightarrow (\theta,I,z)$ 
of the phase space $ \T^\nu\times \R^\nu \times H_S^\perp$, 
\begin{equation}\label{Gdelta}
\begin{pmatrix}
\theta \\ I \\ z
\end{pmatrix} := 
G_\delta 
\begin{pmatrix}
\phi \\ y \\ \tw
\end{pmatrix} := 
\begin{pmatrix}
\theta_0(\phi) \\ 
I_\delta(\phi) + \left[ \pa_\phi \theta_0(\phi) \right]^{-T}y 
+ \left[(\pa_\theta\wtz_0)(\theta_0(\phi))  \right]^T J^{-1} \tw \\
z_0(\phi) + \tw
\end{pmatrix}\,,
\end{equation}
where $\wtz_0(\theta):= z_0(\theta_0^{-1}(\theta))$.
It is proved in Lemma $2$ of \cite{BB} that $G_\delta$ is symplectic, 
because the torus $i_\delta$ is isotropic (Lemma \ref{torus_iso}). 
In the new coordinates, $i_\delta$ 
is the trivial embedded torus $(\phi,y,\tw)=(\phi,0,0)$.
Moreover, the diffeomorphism $G_\delta$ in \eqref{Gdelta} 
is momentum preserving, 
in the sense that  (Lemma 6.3 in \cite{BFM})
$\tau_\vs \circ G_\delta = G_\delta \circ \tau_\vs $, 
$ \forall \,\vs \in \R^d $, 
where $\tau_\vs $ is defined in \eqref{def:vec.tau}.

Under the symplectic diffeomorphism $G_\delta $, the Hamiltonian vector field $X_{H_{\alpha, \varrho}}$ changes into
\[
X_{K_{\alpha, \varrho}} = \left(DG_\delta  \right)^{-1} X_{H_{\alpha, \varrho}} \circ G_\delta 
\qquad 
{\rm where} \qquad K_{\alpha, \varrho} := H_{\alpha, \varrho} \circ G_\delta \,,
\]
is momentum preserving, in the sense that
$ K_{\alpha,\varrho}\circ \tau_\vs = K_{\alpha,\varrho}$,  $ \forall\, \vs \in \R^d  $.
We shall write $K$ instead of 
$K_{\alpha, \varrho}$ in order to simplify the notation.
The Taylor expansion of $K$ at the trivial torus $(\phi,0,0)$ is
\begin{equation}\label{taylor_Kalpha}
\begin{aligned}
K_{\alpha,\varrho}(\phi,y,\tw) =& \ K_{00}(\phi,\alpha) + K_{10}(\phi,\alpha) \cdot y 
+ ( K_{01}(\phi,\alpha),\tw )_{L^2} + \tfrac12 K_{20}(\phi) y\cdot y 
\\& 
+ ( K_{11}(\phi)y,\tw )_{L^2} + \tfrac12 ( K_{02}(\phi)\tw,\tw )_{L^2} 
+ K_{\geq 3}(\phi,y,\tw)+ \varrho\cdot \theta_0 (\phi)\,,
\end{aligned}
\end{equation}
where $K_{\geq 3}$ collects all terms at least cubic in the variables $(y,\tw)$.
Here $K_{00}\in \R$, $K_{10}\in\R^\nu$, 
$K_{01}\in H_S^\perp $, whereas $K_{20} $ is a $ \nu \times \nu $ symmetric matrix, 
$K_{11}\in\cL ( \R^\nu,H_S^\perp )$ and 
$ K_{02} $ is a self-adjoint operator  acting on $H_S^\perp $.

\noindent
The Hamilton equations associated to \eqref{taylor_Kalpha} are 
\begin{equation}\label{hameq_Kalpha}
\begin{cases}
\dot\phi =   K_{10}(\phi,\alpha) + K_{20}(\phi)y + [K_{11}(\phi)]^T \tw 
+ \pa_y K_{\geq 3}(\phi,y,\tw) 
\vspace{0.4em}
\\
\dot y = - \pa_\phi K_{00}(\phi,\alpha) - [\pa_\phi K_{10}(\phi,\alpha)]^T y 
- [\pa_\phi K_{01}(\phi,\alpha)]^T \tw  -[\pa_\phi(\theta_0(\phi))]^T\varrho 
\\
\qquad\quad- \pa_\phi\left( \tfrac12 K_{20}(\phi)y\cdot y + \left( K_{11}(\phi)y,\tw \right)_{L^2} 
+ \tfrac12 \left( K_{02}(\phi)\tw,\tw \right)_{L^2} 
+ K_{\geq 3}(\phi,y,\tw) \right) 
\vspace{0.4em}
\\
\dot\tw =  J \, \left( K_{01}(\phi,\alpha)+ K_{11}(\phi)y 
+ K_{02}(\phi)\tw + \nabla_{\tw} K_{\geq 3}(\phi,y,\tw) \right) \, , 
\end{cases} 
\end{equation}
where $\pa_\phi K_{10}^T $ is the $\nu\times\nu$ 
transposed matrix and 
$\pa_\phi K_{01}^T , K_{11}^T: H_S^\perp\rightarrow\R^\nu$ 
are defined by the duality relation 
$ (\pa_\phi K_{01}[\wh\phi],\tw  )_{L^2}=\wh\phi\cdot [\pa_\phi K_{01} ]^T \tw $ 
for any $\wh\phi\in\R^\nu$, $\tw\in H_S^\perp$. 

The terms $K_{00}, K_{01}$, $K_{10} -  \alpha$ 
in the Taylor expansion \eqref{taylor_Kalpha} vanish at an exact solution,
indeed,  
arguing as  in Lemma 5.4 in \cite{BBHM},  
there is $ \sigma := \sigma (\nu, \tau) > 0 $, such that, 
for all $ s \geq s_0 $, 
\begin{equation}\label{Kcoeff_est}
\begin{aligned}
\norm{ \pa_\phi K_{00}(\cdot , \alpha) }_s^{k_0,\gamma} 
+ \norm{ K_{10}(\cdot ,\alpha)-\alpha}_s^{k_0,\gamma} 
&+ \norm{ K_{01}( \cdot ,\alpha) }_s^{k_0,\gamma} 
\\&
\lesssim_s  
\norm{Z}_{s+\sigma}^{k_0,\gamma} 
+ \norm{Z}_{s_0+\sigma}^{k_0,\gamma} 
\norm{ \cI_0 }_{s+\sigma}^{k_0,\gamma} \, .  
\end{aligned}
\end{equation}
Under the linear change of variables
\begin{equation*}
DG_\delta(\vphi,0,0)\begin{pmatrix}
\wh\phi\, \\ \why \\ \wh\tw
\end{pmatrix}:= 
\begin{pmatrix}
\pa_\phi \theta_0(\vphi) & 0 & 0 
\\ 
\pa_\phi I_\delta(\vphi) & [\pa_\phi \theta_0(\vphi)]^{-T} 
&  [(\pa_\theta\widetilde{z}_0)(\theta_0(\vphi))]^T  J^{-1} 
\\
\pa_\phi z_0(\vphi) & 0 & {\rm Id}
\end{pmatrix}
\begin{pmatrix}
\wh\phi \, \\ \why \\ \wh\tw
\end{pmatrix} \,,
\end{equation*}
the linearized operator $\di_{i,\alpha, \varrho}\cF(i_\delta)$ 
is approximately transformed into the one obtained 
when one linearizes the Hamiltonian system \eqref{hameq_Kalpha} 
at $(\phi,y,\tw) = (\vphi,0,0)$, 
differentiating also in $\alpha$ at $\alpha_0$ 
and changing $\pa_t \rightsquigarrow \omega\cdot \pa_\vphi$, 
namely
\begin{equation}\label{lin_Kalpha}
\resizebox{.91\hsize}{!}{$
\begin{pmatrix}
\widehat \phi  \\
\widehat y    \\ 
\widehat \tw \\
\widehat \alpha\\
\widehat \varrho
\end{pmatrix} \mapsto
\begin{pmatrix}
\omega\cdot \pa_\vphi \wh\phi - \pa_\phi K_{10}(\vphi)[\wh\phi] 
- \pa_\alpha K_{10}(\vphi)[\wh\alpha] 
- K_{20}(\vphi)\why - [K_{11}(\vphi)]^T \wh\tw  
\vspace{0.2em}\\
\omega\cdot \pa_\vphi\why + \pa_{\phi\phi}K_{00}(\vphi)[\wh\phi]
+ \pa_\alpha\pa_\phi K_{00}(\vphi)[\wh\alpha] 
+ [\pa_\phi K_{10}(\vphi)]^T \why + [\pa_\phi K_{01}(\vphi)]^T  \wh\tw 
+ [\pa_\phi(\theta_0(\phi))]^T\widehat{\varrho}  
\vspace{0.2em}\\
\omega\cdot \pa_\vphi \wh\tw - J \,  \big( \pa_\phi K_{01}(\vphi)[\wh\phi] 
+ \pa_\alpha K_{01}(\vphi)[\wh\alpha] + K_{11}(\vphi) \why 
+ K_{02}(\vphi) \wh\tw \big)   
\end{pmatrix}\,.
$} 
\end{equation}
In order to construct an ``approximate''
 inverse of \eqref{lin_Kalpha}, 
we need that
\begin{equation}\label{Lomegatrue}
\cL_\omega := 
\mathbb{P}_S^{\perp} \left( \omega\cdot \pa_\vphi - J K_{02}(\vphi) \right)|_{H_S^\perp}
\end{equation}
is ``invertible'' on traveling waves.
Let us now define (recall \eqref{decomposition})
\begin{equation}\label{spazioprodperp}
H_{\perp}^s:=H_{\perp}^s(\T_{*}^{\nu+d};\R^2):= H^s(\T_{*}^{\nu+d};\R^2)\cap H_S^\perp\,.
\end{equation}
We make the following assumption:
\begin{itemize}
\item[(AI)]  {\bf Invertibility of $\cL_\omega$}: 
{\it There exist positive real numbers  $\bar{\mu},\s$ (depending on $k_0,\tau$)
%$ \s(\tb) $, $ \ta $, $ p  $, $ K_0 $, $\kappa_1$ 
and 
a subset 
$\mathtt{\Lambda}_o \subset \tD\tC(\gamma,\tau)\times [\th_1,\th_2]$ 
such that, for all $(\omega,\th) \in \mathtt{\Lambda}_o$, 
and any $\bar{s}\geq s_0+\bar{\mu}+\s$ the following holds.
%the operator $\cL_\omega$ may be decomposed as
%\begin{equation}\label{Lomega}
%\cL_\omega = \cL_\omega^< + \cR_\omega + \cR_\omega^\perp \,,
%\end{equation}
%where, 
For any traveling wave function 
$ g\in H_{\perp}^{s+\s}$
(i.e. satisfying \eqref{condtraembedd}), $s_0\leq s\leq \bar{s}-\bar{\mu}-\s$ and 
for any $  (\omega,\th) \in \mathtt{\Lambda}_o $, 
there is a traveling wave solution 
$ h \in H_{\perp}^{s} $ 
of $ \cL_\omega h = g$ satisfying, for all $s_0\leq s\leq \bar{s} - \bar{\mu}-\s$,  
\begin{equation}\label{almi4}
\norm{ \cL_\omega^{-1}g }_s^{k_0,\gamma} \lesssim_{\bar{s}}
\gamma^{-1}\big( \norm g_{s+\s}^{k_0,\gamma}
+ \norm g_{s_0+\s}^{k_0,\gamma}
\norm{ \cI_0}_{s+\bar{\mu}+\s}^{k_0,\gamma} \big) \,.
\end{equation}	
}
\end{itemize}
This assumption shall be verified, in section \ref{quasi invertibilita},
at each 
step of the Nash-Moser  iteration.

In order to find an almost approximate inverse 
of the linear operator in \eqref{lin_Kalpha} 
(and, as a  consequence, of the operator  $\di_{i,\alpha, \varrho}\cF(i_\delta)$), 
it is sufficient to  invert the operator
\begin{equation}\label{Dsys}
\D\big[ \wh\phi,\why,\wh\tw,\wh\alpha, \wh\varrho \big]:=
\begin{pmatrix}
\omega\cdot \pa_\vphi\wh\phi - \pa_\alpha K_{10}(\vphi)[\wh\alpha] 
- K_{20}(\vphi)\why- K_{11}^T(\vphi)\wh\tw 
\\
\omega\cdot \pa_\vphi \why +\pa_\alpha\pa_\phi K_{00}(\vphi)[\wh\alpha] 
+ [\pa_\phi \theta_0(\phi)]^T\wh\varrho 
\\
\cL_\omega^< \wh\tw -  J \left( \pa_\alpha K_{01}(\vphi)[\wh\alpha] 
+ K_{11}(\vphi)\why \right)
\end{pmatrix}\,,
\end{equation}
obtained neglecting in \eqref{lin_Kalpha} the terms 
$\pa_\phi K_{10}$, $\pa_{\phi\phi}K_{00}$, $\pa_\phi K_{00}$, $\pa_\phi K_{01}$ 
(they vanish at an exact solution by  \eqref{Kcoeff_est}).
%and the small remainders $\cR_\omega$, $\cR_\omega^\perp$ 
%appearing in \eqref{Lomega}. 

As in  section 6 of \cite{BFM} we have  the following result, 
where we denote
$ \normk{(\phi,y,\tw,\alpha, \varrho)}{s}:= 
\max \big\{ \normk{(\phi,y,\tw)}{s}, 
\abs{\alpha}^{k_0,\gamma}, \abs{\varrho}^{k_0,\gamma} \big\} $ 
(see- \cite[Proposition 6.5]{BFM}).
\begin{prop}\label{Dsystem}
Assume \eqref{ansatz}  and {\rm (AI)}. Then, 
for all $(\omega,\th)\in\mathtt{\Lambda}_o $, 
for any traveling wave embedding  $ g =(g_1,g_2,g_3)$ (i.e. satisfying \eqref{RTTT}), 
there exists a unique solution 
$\D^{-1}g:= ( \wh\phi,\why,\wh\tw,\wh\alpha, \wh\varrho)$
of  $\D ( \wh\phi,\why,\wh\tw,\wh\alpha, \wh\varrho) = g $
where $( \wh\phi,\why,\wh\tw)$ is a traveling wave embedding.  
Moreover, for any $s_0\leq s\leq \bar{s}-\bar{\mu}-\s$, 
\[
\normk{\D^{-1}g}{s} \lesssim_{\bar{s}} 
\gamma^{-1}\big( \normk{g}{s+\s}
+\normk{\cI_0}{s+\bar{\mu}+\s}\normk{g}{s_0+\s} \big)\,.
\]
\end{prop}
Finally, we conclude that the operator
\begin{equation}\label{bT0}
\bT_0 := \bT_0(i_0):=  
( D\wtG_\delta )(\vphi,0,0) \circ \D^{-1} \circ (D G_\delta ) (\vphi,0,0)^{-1}
\end{equation}
is an approximate right inverse for $\di_{i,\alpha, \varrho}\cF(i_0)$, 
where 
\[
\wtG_\delta(\phi,y,\tw,\alpha, \varrho) := \left( G_\delta(\phi,y,\tw),\alpha, \varrho \right)
\]
is the identity on the $\alpha$-component. 
Arguing exactly as in  Theorem 6.6 in \cite{BFM}, we deduce the following. 

\begin{thm} {\bf (Approximate inverse)} \label{alm.approx.inv}
Assume {\rm (AI)}.  Then there are 
$\bar\sigma :=\bar\sigma(\tau,\nu,k_0)>0$ and $\kappa_1=\kappa_1(\tau,
\nu,k_0)$ such that, if \eqref{ansatz} holds with 
$\s=\bar{\mu}+\bar\sigma$ and \eqref{epsilonpiccoloKK} holds for some $\tk\gg \kappa_1$, then, 
for all $(\omega,\th)\in\mathtt{\Lambda}_o$ 
and for any   $g:=(g_1,g_2,g_3)$ satisfying \eqref{RTTT},
	 the operator $\bT_0$ defined in \eqref{bT0} satisfies, for all $s_0 \leq s \leq \bar{s} 
	 - \bar{\mu}- \bar\sigma$,
	\begin{equation}\label{tame-es-AI}
	\normk{\bT_0 g}{s} \lesssim_{\bar{s}} \gamma^{-1} \big( \normk{g}{s+\bar\sigma} +\normk{\cI_0}{s+\bar{\mu}+\bar\sigma}\normk{g}{s_0+\bar\sigma}  \big)\,.
	\end{equation}
	Moreover, the first three components of $\bT_0 g $  form a 
	traveling wave embedding (i.e. \eqref{RTTT} holds). 
	Finally, $\bT_0$ is an  approximate right inverse of $\di_{i,\alpha,\varrho}\cF(i_0)$, namely
	\begin{equation}\label{splitting per approximate inverse}
	  \cP(i_0) :=\di_{i,\alpha, \varrho}\cF(i_0) \circ \bT_0 - {\rm Id} \,,
	%+ \cP_\omega(i_0)+\cP_\omega^\perp(i_0)\,,
	\end{equation}
	where for all $s_0 \leq s \leq \bar{s}-\bar{\mu}-\bar\sigma$,
	\begin{align}
	\normk{\cP g}{s} 
	& \lesssim_{\bar{s}} \gamma^{-1}
	\Big(  \normk{\cF(i_0,\alpha_0)}{s_0+\bar\sigma}\normk{g}{s+\bar\sigma}  \nonumber 
	\\
	& \qquad + \,   \big(  \normk{\cF(i_0,\alpha_0)}{s+\bar\sigma}+\normk{\cF(i_0,\alpha_0)}{s_0+\bar\sigma}\normk{\cI_0}{s+\bar{\mu}+\bar\sigma}  \big)\normk{g}{s_0+\bar\sigma}  \Big)\, .\label{pfi0} 
%	\\
%	\normk{\cP_\omega g}{s} & \lesssim_{\bar{s}} \varepsilon\gamma^{-\kappa_1-1} N_{\tn-1}^{-\ta} \big( \normk{g}{s+\bar\sigma}+ \normk{\cI_0}{s+\mu(\tb)+\bar\sigma}\normk{g}{s_0+\bar\sigma}  \big)\, , \label{pfi1}  \\
%	\normk{\cP_\omega^\perp g}{s_0} & \lesssim_{\bar{s},b} \gamma^{-1} K_\tn^{-b} \left( \normk{g}{s_0+\bar\sigma+b}+\normk{\cI_0}{s_0+\mu(\tb)+b+\bar\sigma}\normk{g}{s_0+\bar\sigma} \right)\,, \quad  \forall\,b>0\,,   \label{pfi2} \\
%	\normk{\cP_\omega^\perp g}{s} & \lesssim_{\bar{s}} \gamma^{-1}\big(  \normk{g}{s+\bar\sigma}+ \normk{\cI_0}{s+\mu(\tb)+\bar\sigma}\normk{g}{s_0+\bar\sigma} \big) \,. 
%	\label{pfi3} 
	\end{align}
\end{thm}

%\red{TORNANDO SU CORREGGERE E CAMBIARE}
\section{Symmetrization of the linearized operator at the highest order}\label{linearizzato siti normali}

In order to write an explicit  expression of the linear operator 
${\mathcal L}_\omega$ defined in \eqref{Lomegatrue}
we have to express the operator $ K_{02}(\phi) $ in \eqref{taylor_Kalpha} 
in terms of the original water waves Hamiltonian vector field.  

\begin{lemma} \label{thm:Lin+FBR}
The operator  $ K_{02}(\phi) $ is 
\begin{equation}\label{K 02}
K_{02}(\phi) =\mathbb{P}_S^{\perp} \partial_u \nabla_u H(T_\delta(\phi)) + \e R(\phi) \,,
\end{equation}
where $ H $ is the water waves Hamiltonian defined in \eqref{HamiltonianWW} 
%(with gravity constant $ g =1 $ and depth 
%$ h $ replaced by $ \mathtt h $), 
evaluated at the torus 
\begin{equation}\label{T delta}
T_\delta(\phi) := \e A(i_\delta(\phi)) = \e A(\theta_0(\phi), I_\delta(\phi), z_0(\phi) ) =
\e  v (\theta_0 (\phi), I_\delta(\phi)) +  \e z_0 (\phi) 
\end{equation}
with $ A (\theta, I, z ) $, $ v (\teta, I )$ defined in \eqref{actionanglevar}.
%The operator $ K_{02}(\phi) $ is even and reversible.  
The remainder $ R(\phi) $ has the ``finite dimensional" form 
\begin{equation}\label{forma buona resto}
R(\phi)[h] = {\mathop \sum}_{j \in S} \big(h\,,\,g_j \big)_{L^2_x} \chi_j\,, \quad \forall h \in H_{S}^\bot \, ,  
\end{equation}
for functions $ g_j, \chi_j \in H_{S}^\bot  $ which satisfy the tame estimates: 
for some $ \sigma:= \sigma(\tau, \nu) > 0 $, 
$ \forall s \geq s_0 $, 
\begin{equation}\label{stime gj chij}
\begin{aligned}
\| g_j\|_s^{k_0, \gamma} +\| \chi_j\|_s^{k_0, \gamma} &\lesssim_s 
1 + \| {\mathcal I}_\delta\|_{s + \s}^{k_0, \gamma}\,,
\\
 \| \pa_{i}g_j[\widehat \imath]\|_s +\| \partial_i \chi_j
[\widehat \imath]\|_s& \lesssim_s 
\| \widehat \imath \|_{s + \s}+ \| {\mathcal I}_\delta\|_{s + \s} \| \widehat \imath\|_{s_0 + \s}\, . 
\end{aligned}
\end{equation}
\end{lemma}

\begin{proof}
The lemma follows as in Lemma 6.1 in \cite{BM20}, 
(see also Lemma 7.1 in \cite{BFM}-\cite{BFM2}).
\end{proof}

By Lemma \ref{thm:Lin+FBR} the linear operator $ {\mathcal L}_\omega $ defined in \eqref{Lomegatrue} has the form 
\begin{equation}\label{representation Lom}
{\mathcal L}_\omega =  \mathbb{P}_S^{\perp} ( {\mathcal L} + \e R) \mathbb{P}_S^{\perp}\,,
\qquad 
{\rm where} \qquad {\mathcal L} := 
\omega \cdot \partial_\vphi   - J \partial_u  \nabla_u 
H (T_\delta(\vphi))
\end{equation}
is obtained linearizing the original water waves system at the torus
$ u = (\eta, \psi) = T_\delta(\vphi) $ 
defined in \eqref{T delta}, changing $ \pa_t \rightsquigarrow \omega \cdot \partial_\vphi $. 
The operator $\mathcal{L}$ above has an explicit structure in terms of 
the Dirichlet-Neumann operator $G(\eta)$  in \eqref{eq:112a},
whose properties have been discussed in Sections \ref{sec:DNsezione} and \ref{sec:pseudoDN}.

In particular, we recall that 
the map $(\eta,\psi)\to G(\eta)\psi$ is linear with respect to $\psi$
and nonlinear with respect to the profile $\eta$.
The derivative with respect to $\eta$ (which is called ``shape derivative'')
is given by the formula  in \eqref{shapeDer}.
Recall also that by $V(\eta,\psi), B(\eta,\psi)$ 
we denoted 
the horizontal and vertical components 
of the velocity field at the free interface, see formul\ae\, \eqref{def:V}-\eqref{form-of-B}.
Notice that,
by the shape derivative formula \eqref{shapeDer},  
the linearized operator of \eqref{eq:113}
at $(\eta,\psi)(\vphi,x)$ is given by 
\begin{equation}\label{linWW}
\mathcal{L}:=\mathcal{L}(\vphi):=\omega\cdot\pa_{\vphi}+
\left(
\begin{matrix}
V\cdot\nabla+{\rm div}(V)+G(\eta)B & -G(\eta) \\
1 +B{\rm div}(V)+BG(\eta)B & V\cdot\nabla-BG(\eta)
\end{matrix}
\right)\,.
\end{equation}
We recall the set of \emph{tangential sites} $S$ 
defined in \eqref{TangentialSitesWW}
and
that 
 $\mathbb{P}_S, \mathbb{P}_S^{\perp}$ are the orthogonal projections on the spaces 
 $H_S, H_S^\bot$ defined in \eqref{decomposition}.

\smallskip

We are going to make several transformations, whose aim is to conjugate the linearized 
operator to a constant coefficients operator, up to a remainder 
that is small in size and regularizing at a conveniently high order.

\smallskip

For the sequel we will always assume the following ansatz 
(satisfied by the approximate solutions obtained along the nonlinear  Nash-Moser iteration of Section
\ref{sec:NaM}): 
for some constants $\mu_0 :=\mu_0(\tau,\nu,k_0)>0$, $\gamma\in (0,1)$, (cfr. 
Lemma \ref{torus_iso})
\begin{equation}\label{ansatz_I0_s0}
\norm{ \cI_0 }_{s_0+\mu_0} ^{k_0,\gamma} \, ,  \ 
\norm{ \cI_\delta }_{s_0+\mu_0}^{k_0,\gamma}  
\leq 1 
\end{equation} 
where we recall that $s_0$ satisfies \eqref{Sone}. 
In order to estimate the variation of the eigenvalues with respect to the approximate invariant torus, 
we need also to estimate the variation with respect to the torus $i(\vphi)$ 
in another low norm $\|\cdot\|_{p}$ 
for all Sobolev indexes $p$ such that
\begin{equation}\label{ps0}
p+\sigma_0 \leq s_0 +\mu_0 \,, \quad \text{ for some } \ \sigma_0:=\sigma_0(\tau,\nu)>0 \,. 
\end{equation}
Thus, by \eqref{ansatz_I0_s0}, we have
$ \norm{ \cI_0 }_{p+\sigma_0}^{k_0,\gamma} $,
$  \norm{ \cI_\delta }_{p+\sigma_0}^{k_0,\gamma} \leq 1 $.
The constants $\mu_0$ and $\sigma_0$ represent the \emph{loss of derivatives} accumulated along the reduction procedure of the next sections. What is important is that they are independent of the Sobolev index $s$. 

 In general, the loss of derivatives $\s_0$ 
 needed to estimate the Lipschitz variation (which is independent of $k_0$) is smaller w.r.t. $\mu_0$. 
 Nevertheless,
 to lighten the notation, 
 along Sections \ref{linearizzato siti normali}-\ref{sym.low.order}-\ref{sec:redulower},
 we shall denote by $\sigma := \sigma(k_0, \tau, \nu) > 0$ a constant 
 (which possibly increases from lemma to lemma) representing the loss 
 of derivatives along the finitely many steps of the reduction procedure, without  making 
 distinction between $\s_0$ and $\mu_0$. We will measure the Lipschitz variation in the norm
 $p$ satisfying \eqref{ps0} by meaning that $\mu_0$ in the ansatz \eqref{ansatz_I0_s0}
is taken, w.l.o.g., larger than $\s$.

\medskip
As a consequence of Moser composition Lemma \ref{Moser norme pesate}, 
the Sobolev norm of the function $ u = T_{\delta} $ defined in \eqref{T delta} (see also \eqref{actionanglevar}) 
satisfies, $ \forall s \geq s_0 $,  
\begin{equation}\label{tame Tdelta}
\| u \|_s^{k_0,\gamma} = \| \eta \|_s^{k_0,\gamma}  + 
\| \psi \|_s^{k_0,\gamma} \leq \e C(s)  \big(  1 + \| {\mathcal I}_0 \|_{s}^{k_0, \gamma})  \,,
\end{equation}
and
\begin{equation}\label{derivata i T delta}
%\| \partial_i u [\hat \imath] \|_{p} \lesssim_{s_1} \e \| \hat \imath\|_{p}, \qquad
\| \Delta_{12} u \|_{p} \lesssim_{p} \e \| i_2 - i_1 \|_{p} \,,
\end{equation}
where we denote $ \Delta_{12} u:= u(i_2) - u(i_1)$ (see \eqref{deltaunodue}).
%; we will systematically use this notation.
In the next sections we shall also assume % (and it will not be reminded) 
that, for some $ \kappa := \kappa (\tau, \nu) > 0 $, we have  
\begin{equation}\label{smallcondepsilon}
\e \gamma^{- \kappa} \leq \delta (\bar{s})  \, ,
\end{equation}
where $ \delta (\bar{s}) > 0 $ is a constant small enough and  $ \bar{s} $ 
is a high regularity index $\bar{s}=s_0+\widetilde{\s}$  that will be fixed below in Section
\ref{sec:NaM}.
We recall that $ {\mathcal I}_0 := {\mathcal I}_0 (\omega, \th )  $ 
is defined for all  $ (\omega, \th) \in \R^\nu \times [\th_1, \th_2] $ by the extension procedure
that we perform along the Nash-Moser nonlinear iteration. Moreover, all the functions appearing in 
$ {\mathcal L } $ in \eqref{linWW} are $ {\mathcal C}^\infty $ in $ (\vphi, x) $ as 
the approximate torus 
$ u = (\eta, \psi) =  T_\delta (\vphi)  $. 
In particular, we have the following:
there is  $\sigma \gg 0$ large enough such that, for any $\bar{s}\geq s_0+\s$,
if the smallness conditions \eqref{ansatz_I0_s0}-\eqref{smallcondepsilon}
 hold (with $\mu_0\geq \s$) 
 the functions $V(\vphi, x)$ and $B(\vphi, x)$ satisfy  
 %for $s \geq s_0$ 
the tame estimates 
%(recall Lemma \ref{stime tame dirichlet neumann})
\begin{equation}\label{stime tame V B inizio lineariz}
\begin{aligned}
\| V \|_s^{k_0, \gamma}, \| B \|_s^{k_0, \gamma} 
&\lesssim_s \e (1 + \| {\mathcal I}_0 \|_{s + \sigma}^{k_0, \gamma})\,,\;\;\;s_0\leq s\leq \bar{s}\,, \\
\| \Delta_{12}V \|_p, \| \Delta_{12}B \|_p &\lesssim_p \e  \|i_1-i_2\|_{p+\sigma}\,,
\end{aligned}
\end{equation}
with $p$ as in \eqref{ps0}. This follows by Lemma \ref{stime tame dirichlet neumann}.
By \eqref{linWW} and Lemma \ref{thm:Lin+FBR}, the linear operator ${\mathcal L}_\omega : H_S^\bot \to H_S^\bot$ takes the form 
\begin{equation}\label{cal L omegaWW}
\begin{aligned}
\mathcal{L}_\omega&:=\mathcal{L}_{\oo}(\vphi)
\\&:=\mathbb{P}_S^{\perp}\omega\cdot\pa_{\vphi}+
\mathbb{P}_S^{\perp} \left(
\begin{matrix}
V\cdot\nabla+{\rm div}(V)+G(\eta)B & -G(\eta) \\
1 +B{\rm div}(V)+BG(\eta)B & V\cdot\nabla-BG(\eta)
\end{matrix}
\right) \mathbb{P}_S^{\perp} + \e R \,,
\end{aligned}
\end{equation}
where the finite rank operator $R$ is given in Lemma \ref{thm:Lin+FBR}.

The next lemma is fundamental for our scope.

\begin{lemma}\label{algStrucHamlinear}
Consider functions $(\eta,\psi)\in S_{\mathtt{V}}$ (see \eqref{SV})
and the 
linearized operator 
$\mathcal{L}$ in \eqref{linWW}.
Then $\mathcal{L}$ 
is \emph{Hamiltonian} and \emph{momentum preserving}
according to Def. \ref{operatoreHam}.
The same holds true for the operator $\mathcal{L}_{\omega}$ in \eqref{cal L omegaWW}.
\end{lemma}

\begin{proof}
It follows by  Lemma $7.9$ in \cite{FGtrave}. See also Lemma $7.3$ in \cite{BFM}.
\end{proof}

 The aim of the following sections is to 
conjugate the linearized operator $\mathcal{L}_{\omega}$ in \eqref{cal L omegaWW}
to a constant coefficients operator, up to a regularizing 
remainder.
This will be achieved by applying several transformations
which clearly depends nonlinearly on the point $(\eta,\psi)$ 
on which  we linearized.

We first start by symmetrizing the operator $\mathcal{L}_{\omega}$ at the highest order.
This is the content of the following subsections.
In particular, we remark that sections \ref{sec:goodunknown}, \ref{sec:coordinatecomplesse}
are quite standard and basically independent of the space dimension. 
On the contrary, starting from subsection \ref{sec:almoststraightening},  the reduction procedure is more delicate and the %straightening of the first order vector field is 
deeply affected by the dimension.

\subsection{Good unknown of Alinhac}\label{sec:goodunknown}
The aim of this section is to rewrite 
the operator $\mathcal{L}_{\omega}$ 
in \eqref{cal L omegaWW} in terms of the so called ``good unknown`` of Alinhac 
(more precisely its symplectic correction). This will be done in 
Proposition \ref{lemma:7.4}. As we will see, these 
coordinates are the correct ones
in order to diagonalize, at the highest order, the operator $\mathcal{L}_{\omega}$.
We shall first prove some preliminary results.
Following \cite{AB15}, \cite{BM20}
we conjugate the linearized operator $\mathcal{L}$ in \eqref{linWW}
by the operator (recall \eqref{decomposition})
\begin{equation}\label{flussoG}
\begin{aligned}
\mathcal{G}_{B}&=\mathbb{P}_S^{\perp}\mathcal{G}\mathbb{P}_S^{\perp}={\rm Id}+{\bf B}
\;:\; H_{S}^{\perp}\to H_{S}^{\perp}\,,
\\
{\bf B}&:=\mathbb{P}_S^{\perp}\left( 
\begin{matrix} 0 & 0\\ B & 0
\end{matrix}
\right)\mathbb{P}_S^{\perp}\,,
\qquad
\mathcal{G} :=\left( 
\begin{matrix} 1 & 0\\ B & 1
\end{matrix}
\right)\,,
\end{aligned}
\end{equation}
where $B$ is the real valued function in \eqref{form-of-B}.
Recalling also \eqref{def:V}
define the functions
\begin{equation}\label{func:a}
a(\vphi,x):=(\omega\cdot\pa_{\vphi}B)+V\cdot\nabla B\,, \quad b(\vphi, x) := {\rm div}(V)\,.
\end{equation}
In the following lemma we study some properties of the map $\mathcal{G}_{B}$.
\begin{lemma}\label{lem:mappGood}
There exists $\sigma \gg 0$ large enough such that, for  any $\bar{s} > s_0$,
if \eqref{smallcondepsilon} and \eqref{ansatz_I0_s0} hold (with $\mu_0 \geq \sigma$), then the following holds. 
 The following statements hold for any $s_0\leq s\leq \bar{s}$:  
\begin{align}
\| a \|_s^{k_0, \gamma}\,,\, \| b \|_s^{k_0, \gamma} &
 \lesssim_s \e (1 + \| {\mathcal I}_0 \|_{s + \sigma}^{k_0, \gamma})\,, \label{stime a b dopo alinac} 
 \\
 \| \Delta_{12}a \|_p, \| \Delta_{12}b\|_p &\lesssim_p \e \|i_1-i_2\|_{p+\sigma}\,,
 \label{stime a b dopo alinac2}
 \\
\fM_{(\mathcal{G}^{\pm 1}-\id)}(s) %+\fM_{(\mathcal{G}^{\pm 1}-\id)^*}(s) 
&\lesssim_{\bar{s}} \eps (1+ \|\cI_0 \|_{s+ \sigma}^{k_0,\gamma})\,,\label{cost tame good}
\\
\| \Delta_{12} (\mathcal{G}^{\pm 1}) h\|_p %+ \| \Delta_{12} (\mathcal{G}^{\pm 1})^* h\|_p 
&\lesssim_{p} \eps \|h\|_p \|i_1-i_2\|_{p+\sigma}\,,\label{delta12 good}
\end{align}
for $p$ as in \eqref{ps0}.
The functions $a$ and $b$ belongs to $S_{\mathtt{V}}$.
Moreover, the map $\mathcal{G}_{B}$ is symplectic on $H_S^\bot$
and momentum preserving. Finally the maps $\mathcal{G}_{B}^{\pm1}$ satisfy  the same estimates
of $\mathcal{G}$.
\end{lemma}
\begin{proof}
We shall apply estimates \eqref{stime tame V B inizio lineariz} on the functions $V,B$ with 
$\bar{s}\rightsquigarrow \bar{s}+1$. Then
estimates \eqref{stime a b dopo alinac}-\eqref{stime a b dopo alinac2} 
follow by \eqref{stime tame V B inizio lineariz}
and the tame estimates given by Lemma \ref{lemma:LS norms}.
Recalling that we linearized on a tavelling wave $\eta,\psi\in S_{\mathtt{V}}$ (see \ref{SV})
we have that $a,b\in S_{\mathtt{V}}$ thanks to the algebraic properties on $B,V$ given by Theorem 
\ref{lemma totale dirichlet neumann}.
Estimates \eqref{cost tame good}-\eqref{delta12 good}
follow by Lemma \ref{lemma: action Sobolev}.
A direct calculation shows that ${\mathcal G}_B$ is symplectic on $H_S^\bot$, since it is the restriction of the good unknown of Alhinac (which is symplectic on $L^2$) to the symplectic subspace $H_S^\bot$. A direct calculations shows that 
(recall \eqref{flussoG})
${\bf B}^n = 0$ for any $n \geq 2$, implying that ${\mathcal G}_B$ is invertible and 
\begin{equation}\label{inverso cal GB}
{\mathcal G}_B^{- 1} = \mathbb{P}_S^{\perp} - {\bf B} = \mathbb{P}_S^{\perp} \mathcal{G}^{-1}
%\begin{pmatrix}
%1 & 0 \\
%- B & 1\end{pmatrix} 
\mathbb{P}_S^{\perp}\,,\qquad \mathcal{G}^{-1}=\sm{1}{0}{-B}{1}\,.
\end{equation}
Therefore, the thesis follows. 
\end{proof}

We now study the conjugate of the operator $\mathcal{L}_{\omega}$.
\begin{prop}{\bf (Symplectic good unknown).}\label{lemma:7.4}
The conjugate the operator $\mathcal{L}_{\omega}$ in \eqref{cal L omegaWW}
under the map $\mathcal{G}_{B}$ in \eqref{flussoG} has the form
\begin{equation}\label{primamappa}
\mathcal{L}_1=\mathcal{G}_{B}^{-1}\mathcal{L}_{\omega}
\mathcal{G}_{B}=
\mathbb{P}_S^{\perp}\omega \cdot \partial_\vphi + \mathbb{P}_S^{\perp} \left(\begin{matrix}
V\cdot\nabla+ b(\vphi, x) & -G(\eta) \\
1 + a(\vphi, x) & V\cdot\nabla
\end{matrix}\right)\mathbb{P}_S^{\perp} +\e {\mathcal R}_1\,,
\end{equation}
where $a,b$ are in \eqref{func:a} and the remainder $\mathcal{R}_1$ is 
of the form \eqref{forma buona resto}
for some functions $\chi_j^{(1)}, g_j^{(1)}$
satisfying the estimates
\begin{equation}\label{stimeRR11}
\begin{aligned}
\| \chi_j^{(1)} \|_s^{k_0, \gamma}\,,\, \| g_j^{(1)} \|_s^{k_0, \gamma}& \lesssim_s 
 1 + \| {\mathcal I}_0 \|_{s + \sigma}^{k_0, \gamma}\,, \quad \forall s_0 \leq s \leq \bar{s}\,. 
\\
 \| \Delta_{12}g_j^{(1)}\|_p +\| \Delta_{12} \chi_j\|_p& \lesssim_p 
\| i_1-i_2 \|_{p + \s}\, ,
\end{aligned}
\end{equation}
$p$ as in \eqref{ps0}.
 Finally, the operator $\mathcal{L}_1$ is 
Hamiltonian and momentum preserving.
\end{prop}

\begin{proof}
Recalling that $\mathbb{P}_S^{\perp}={\rm Id}-\mathbb{P}_S$ and using 
\eqref{cal L omegaWW} (see also \eqref{linWW})
and  \eqref{flussoG}, \eqref{inverso cal GB}, we get
\begin{equation}\label{tarditardi}
\begin{aligned}
{\mathcal L}_1  := {\mathcal G}_B^{- 1} \circ {\mathcal L}_\omega \circ {\mathcal G}_B 
%\\&
&=\mathbb{P}_S^{\perp}
\mathcal{G}^{-1}
\mathbb{P}_S^{\perp}
\mathcal{L}\mathbb{P}_S^{\perp}
\mathcal{G}
\mathbb{P}_S^{\perp}
+\e {\mathcal G}_B^{- 1} R {\mathcal G}_B 
\\&=
\mathbb{P}_S^{\perp}
\mathcal{G}^{-1}
\mathcal{L}
\mathcal{G}
\mathbb{P}_S^{\perp}+\mathcal{R}_1
\end{aligned}
\end{equation}
%\begin{equation}\label{tarditardi}
%\begin{aligned}
%{\mathcal L}_1 & := {\mathcal G}_B^{- 1} \circ {\mathcal L}_\omega \circ {\mathcal G}_B 
%\\&=\mathbb{P}_S^{\perp}
%\begin{pmatrix}
%1 & 0 \\
%- B & 1
%\end{pmatrix}\mathbb{P}_S^{\perp}
%\mathcal{L}\mathbb{P}_S^{\perp}
%\begin{pmatrix}
%1 & 0 \\
%B & 1
%\end{pmatrix}\mathbb{P}_S^{\perp}
%+\e {\mathcal G}_B^{- 1} R {\mathcal G}_B 
%\\&=
%\mathbb{P}_S^{\perp}\begin{pmatrix}
%1 & 0 \\
%- B & 1
%\end{pmatrix}
%\mathcal{L}
%\begin{pmatrix}
%1 & 0 \\
%B & 1
%\end{pmatrix}\mathbb{P}_S^{\perp}+\mathcal{R}_1
%\end{aligned}
%\end{equation}
where we defined 
\begin{equation}\label{def a b cal R1}
\begin{aligned}
\mathcal{R}_1&:=\e {\mathcal G}_B^{- 1} R {\mathcal G}_B -
\mathbb{P}_S^{\perp}
\mathcal{G}^{-1}
\mathbb{P}_S
\mathcal{L}\mathbb{P}_S^{\perp}
\mathcal{G}
\mathbb{P}_S^{\perp}
%\\&
-
\mathbb{P}_S^{\perp}
\mathcal{G}^{-1}
\mathcal{L}\mathbb{P}_S
\mathcal{G}
\mathbb{P}_S^{\perp}\,.
\end{aligned}
\end{equation}
Now, using  \eqref{linWW}, \eqref{flussoG}, \eqref{inverso cal GB} 
and an explicit computation, we get
\[
\begin{aligned}
\mathcal{G}^{-1}
\mathcal{L}
\mathcal{G}
&=
\begin{pmatrix}
0 & 0 \\
\omega \cdot \partial_\vphi B & 0
\end{pmatrix} + 
\left(\begin{matrix}
V\cdot\nabla+ {\rm div}(V) & -G(\eta) \\
1 + V\cdot\nabla B & V\cdot\nabla
\end{matrix}\right)\,.
\end{aligned}
\]
%\[
%\begin{aligned}
%\mathcal{G}^{-1}
%\mathcal{L}
%\mathcal{G}
%&=
%\begin{pmatrix}
%1 & 0 \\
%- B & 1
%\end{pmatrix}
%\begin{pmatrix}
%\omega\cdot\pa_{\vphi} & 0 \\ 0 & \omega\cdot\pa_{\vphi}
%\end{pmatrix}
%\begin{pmatrix}
%1 & 0 \\
%B & 1
%\end{pmatrix}
%\\&+
%\begin{pmatrix}
%1 & 0 \\
%- B & 1
%\end{pmatrix}
%\left(
%\begin{matrix}
%V\cdot\nabla+{\rm div}(V)+G(\eta)B & -G(\eta) \\
%1 +B{\rm div}(V)+BG(\eta)B & V\cdot\nabla-BG(\eta)
%\end{matrix}
%\right)
%\begin{pmatrix}
%1 & 0 \\
%B & 1
%\end{pmatrix}
%\\&=
%\begin{pmatrix}
%0 & 0 \\
%\omega \cdot \partial_\vphi B & 0
%\end{pmatrix} + 
%\left(\begin{matrix}
%V\cdot\nabla+ {\rm div}(V) & -G(\eta) \\
%1 + V\cdot\nabla B & V\cdot\nabla
%\end{matrix}\right)\,.
%\end{aligned}
%\]
The latter expression, together with \eqref{tarditardi} and \eqref{func:a},
implies formula \eqref{primamappa} with the remainder $\mathcal{R}_1$ 
defined in \eqref{def a b cal R1}.
%one then has that 
%\[
%\begin{aligned}
%{\mathcal L}_1 & := {\mathcal G}_B^{- 1} \circ {\mathcal L}_\omega \circ {\mathcal G}_B 
%\\&
% = \omega \cdot \partial_\vphi  + \mathbb{P}_S^{\perp} \Big(\begin{pmatrix}
%0 & 0 \\
%\omega \cdot \partial_\vphi B & 0
%\end{pmatrix} + \begin{pmatrix}
%1 & 0 \\
%- B & 1
%\end{pmatrix} \mathcal{L}_\omega \begin{pmatrix}
%1 & 0 \\
% B & 1
%\end{pmatrix}  \Big)\mathbb{P}_S^{\perp} + {\mathcal R}_1   
%\\&
%= \omega \cdot \partial_\vphi + \mathbb{P}_S^{\perp} \left(\begin{matrix}
%V\cdot\nabla+ b(\vphi, x) & -G(\eta) \\
%1 + a(\vphi, x) & V\cdot\nabla
%\end{matrix}\right) \mathbb{P}_S^{\perp} + {\mathcal R}_1
%\end{aligned} 
%\]
%where we used that  $\mathbb{P}_S \mathbb{P}_S^{\perp} = 0$, \eqref{func:a}, and where
%\begin{equation}\label{def a b cal R1}
%\begin{aligned}
%& {\mathcal R}_1 := \e {\mathcal G}_B^{- 1} R {\mathcal G}_B 
%-  \mathbb{P}_S^{\perp} \begin{pmatrix}
%0 & 0 \\
%- B & 0
%\end{pmatrix} \mathbb{P}_S \cL_\omega  \mathbb{P}_S^{\perp}\begin{pmatrix} 1 & 0 \\
%B & 1
%\end{pmatrix} \mathbb{P}_S^{\perp}  \\
%& - \mathbb{P}_S^{\perp} \begin{pmatrix}
%1 & 0 \\
%- B & 1
%\end{pmatrix} \cL_\omega \mathbb{P}_S \begin{pmatrix} 0 & 0 \\
%B & 0
%\end{pmatrix} \mathbb{P}_S^{\perp} \,.
%\end{aligned}
%\end{equation}
Using the estimates \eqref{stime tame V B inizio lineariz}, 
Lemma \ref{lemma totale dirichlet neumann} and the fact that $\mathbb{P}_S$ is a finite rank operator, 
one has that the remainder $\cR_{1}$ has the finite rank form \eqref{forma buona resto}
and that it satisfies the estimates \eqref{stimeRR11}.
The operator $\mathcal{L}_1$ is Hamiltonian and momentum preserving 
since, by Lemma \ref{lem:mappGood}, the map $\mathcal{G}_{B}$
is symplectic and momentum preserving.
\end{proof}

\subsection{Complex formulation of Water waves}\label{sec:coordinatecomplesse}
We want to rewrite the operator $\cL_1$ in \eqref{primamappa} 
in the complex coordinates \eqref{CVWW}.

\begin{prop}\label{prop operatore cal L2}
Let  $\beta_{0} \in \N$, $M \in \N$, $\mathtt c \geq k_0 + \beta_0 + 3$ with $M \geq \mathtt c$. Then there is $\s=\s(\beta_{0}, M,k_0)\gg1$ 
large enough such that, for any $\bar{s} > s_0$,
if \eqref{smallcondepsilon}, \eqref{ansatz_I0_s0}-\eqref{ps0} hold for $\mu_0\gg\s$ then
(recall \eqref{wild}, \eqref{CVWW}, \eqref{decomp siti tangenziali coordinate complesse}, 
\eqref{primamappa}) one has that
\begin{equation}\label{elle2}
{\mathcal L}_2 := {\mathcal C} \circ {\mathcal L}_1 \circ {\mathcal C}^{- 1} =
 {\bf \Pi}_S^\bot\omega\cdot\partial_{\vphi}+ 
 {\bf \Pi}_S^\bot\Big(  \opw(\ii V\cdot \xi) \uno + \ii E \Omega(D)  
+ {\mathcal A}_2 + {\mathcal B}_2 + {\mathcal R}_2\Big) {\bf \Pi}_S^\bot
\end{equation}
where $E$ is in \eqref{ham1}, $\uno$ in \eqref{egogenerator}
and ${\mathcal A}_2, {\mathcal B}_2, {\mathcal R}_2$ satisfy the following properties. The operator ${\mathcal A}_2 = {\rm Op}^W \big( A_2(\vphi, x, \xi) \big) \in OPS^{\frac12}$ has the form 
\begin{equation}\label{cal A2 Hamiltoniano 0}
 \begin{aligned}
& A_2(\vphi, x, \xi) = \ii \begin{pmatrix}
a_2(\vphi, x, \xi) & b_2(\vphi, x, \xi) \\
- b_2(\vphi, x, \xi) & - a_2(\vphi, x,  \xi)
\end{pmatrix}\,, \\
& a_2(\vphi, x, \xi) =  \overline{a_2(\vphi, x,  \xi)}, \quad b_2(\vphi, x, \xi) = \overline{b_2(\vphi, x,  \xi)}\,, \\
& a_2(\vphi, x, \xi) = a_2(\vphi, x, - \xi), \quad b_2(\vphi, x, \xi) = b_2(\vphi, x, - \xi)\,,
\end{aligned}
\end{equation}
and satisfy the estimates 
\begin{equation}\label{stima cal A 2 new}
\begin{aligned}
\| {\mathcal A}_2 \|_{\frac12, s, \alpha}^{k_0, \gamma} &\lesssim_{s,M,\alpha} 
\e (1 + \| {\mathcal I}_0 \|_{s + \sigma }^{k_0, \gamma})\,, 
\quad \forall s_0 \leq s \leq \bar{s}\,, \quad \forall \alpha \in \N\,,
\\
\| \Delta_{12}{\mathcal A}_2 \|_{\frac12, p, \alpha}&\lesssim_{p,M,\alpha} \e \|i_1-i_2\|_{p+\sigma}\,,
\end{aligned}
\end{equation}
%for $p$ as in \eqref{ps0},
while the operator ${\mathcal B}_2 \in OPS^0$ satisfies the estimates
\begin{equation}\label{stima cal B 2 new}
\begin{aligned}
\| {\mathcal B}_2 \|_{0, s, \alpha}^{k_0, \gamma} &\lesssim_{s,M,\alpha} 
\e (1 + \| {\mathcal I}_0 \|_{s + \sigma }^{k_0, \gamma})\,, 
\quad \forall s_0 \leq s \leq \bar{s}\,, \quad \forall \alpha \in \N\,,
\\
\| \Delta_{12}{\mathcal B}_2 \|_{0, p, \alpha}&\lesssim_{p,M,\alpha} \e \|i_1-i_2\|_{p+\sigma }\,.
\end{aligned}
\end{equation}
for $p$ as in \eqref{ps0}.\\
The remainder ${\mathcal R}_2$ satisfies the following property. For any $\beta\in\N^{\nu}$, $|\beta|\le \beta_{0}$,
one has that  the operator 
$ \langle D \rangle^{M}\partial_\vphi^\beta {\mathcal R}_2 \langle D \rangle^{- \mathtt c}$ 
is a ${\mathcal D}^{k_0}$-tame operator  (see Def. \ref{Dksigmatame}) satisfying 
\begin{equation}\label{stima tame cal R2 new}
\begin{aligned}
{\mathfrak M}_{ \langle D \rangle^{M}
 \partial_\vphi^\beta{\mathcal R}_2 \langle D \rangle^{- \mathtt c}}(s) 
&\lesssim_{\bar{s},M} \e (1 + \| {\mathcal I}_0 \|_{s + \sigma})\,, 
\quad \forall s_0 \leq s \leq \bar{s}\,,
\\
\|  \langle D \rangle^{M} 
\partial_\vphi^\beta \Delta_{12}{\mathcal R}_2\langle D \rangle^{- \mathtt c} 
\|_{\cL(H^{p} , H^{p})}& \lesssim_{\bar{s},M}\e\|i_1-i_2\|_{p+\s}\,,
\end{aligned}
\end{equation}
$p$ as in \eqref{ps0}.
Moreover, the operators
${\mathcal L}_2$,  
$ {\mathcal A}_2$,  ${\mathcal B}_2$ and $\cR_{2}$, are real-to-real, 
Hamiltonian and momentum preserving.
\end{prop}

\begin{proof}
We reason as in Proposition $8.3$ in \cite{FGtrave}. 
First of all, note that the operator $\mathcal{L}_1$ in \eqref{primamappa}
depends on the functions $a,b$ estimated in Lemma \ref{lem:mappGood}, and on the Dirichlet-Neumann operator 
$G(\eta)$.
In order to use the pseudo-differential expansion \eqref{pseudoespansione}, we 
apply Theorem \ref{lemma totale dirichlet neumann} with 
$\beta_0\in \N$, $M\geq \mathtt{c}$ large enough
and $\bar{s}\rightsquigarrow \bar{s}+\s$ for some $\s\geq \s_1$. The smallness conditions 
\eqref{smallcondepsilon}, \eqref{ansatz_I0_s0}, together with \eqref{tame Tdelta}, 
imply \eqref{smalleta}.

We  note that the map $\mathcal{C}$ 
commutes with the projector  ${\bf \Pi}_S^\bot$, see \eqref{decomp siti tangenziali coordinate complesse2}.
Moreover, since ${\mathcal C}$ is  $\vphi$-independent 
${\mathcal C} \circ \omega \cdot \partial_\vphi  \circ {\mathcal C}^{- 1} = \omega \cdot \partial_\vphi$.
Recalling \eqref{primamappa}, \eqref{func:a}, and using Theorem \ref{lemma totale dirichlet neumann} we define 
the operators
\begin{equation}\label{ariafredda1}
\begin{aligned}
A & := V \cdot \nabla + b= \opw(\ii V\cdot \xi) + \frac{b}{2}  \,, 
\\
B & :=- G(\eta) \stackrel{\eqref{pseudoespansione}}{=}  
- |D | \tanh(\mathtt h |D|) - {\mathcal A}_G - {\mathcal B}_G - {\mathcal R}_G\,, \\
C & := 1 + a  \,, \qquad
D  := V \cdot \nabla = \opw(\ii V\cdot \xi) - \frac{b}{2}  \,. 
\end{aligned}
\end{equation}
By an explicit computation one can check that 
\[
{\mathcal L}_2   := {\mathcal C} \circ {\mathcal L}_1 \circ {\mathcal C}^{- 1} 
= {\bf \Pi}_{S}^{\perp}\omega \cdot \partial_\vphi + {\bf \Pi}_{S}^{\perp}{\mathcal T}{\bf \Pi}_{S}^{\perp}\,, 
\quad {\mathcal T} := \begin{pmatrix}
{\mathcal T}_1 & {\mathcal T}_2 \\
\overline{\mathcal T}_2 & \overline{\mathcal T}_1
\end{pmatrix}\,,
\]
where (recall  \eqref{coniugio complesse simmetrizzate}) 
the operator $ {\mathcal T}_1$, ${\mathcal T}_2$ have the form
\[
\begin{aligned}
{\mathcal T}_{1} = & \frac12 \Big(  \Omega^{- \frac12} A \Omega^{\frac12} 
+ \Omega^{\frac12} D \Omega^{- \frac12} + \ii \Omega^{\frac12} C \Omega^{\frac12} 
- \ii \Omega^{- \frac12} B \Omega^{- \frac12} \Big)\,, 
\\
{\mathcal T}_{2} = & \frac12 \Big( \Omega^{- \frac12} A \Omega^{\frac12} 
- \Omega^{\frac12} D \Omega^{- 1} + \ii \Omega^{- \frac12} B \Omega^{- \frac12} 
+ \ii \Omega^{\frac12} C \Omega^{\frac12}   \Big)\,,
\end{aligned}
\]
 with $A,B,C,D$ in \eqref{ariafredda1} and $\Omega$ in \eqref{wild}.
 We now study the structure of the operator $\mathcal{T}_1$ above,
 by studying each
summand  separately.

\medskip
\noindent
Let us consider the term $-\ii \Omega^{- \frac12} B \Omega^{- \frac12}$.
Fix $N:=M-\tc+2 $. % large w.r.t. $M\geq \mathtt{c}$ (to be chosen later).
Recalling \eqref{formula composizione pseudo} and \eqref{espanzionecomposizioneesplicita}, we have
\[
\begin{aligned}
-\ii \Omega^{- \frac12} B \Omega^{- \frac12}&\stackrel{\eqref{ariafredda1}}{=}
\ii \Omega^{- \frac12} |D | \tanh(\mathtt h |D|) \Omega^{- \frac12}
+\ii \Omega^{- \frac12} \Big(\mathcal{A}_{G}+\mathcal{B}_{G}\Big) \Omega^{- \frac12}+
+\ii \Omega^{- \frac12}\mathcal{R}_{G} \Omega^{- \frac12}
\\&\stackrel{\eqref{wild}, \eqref{pseudoespansione}}{=}
\ii \Omega(D)+
\ii \opw\big(\Omega^{- \frac12}(\x)\#^{W} \big(a_{G}+b_{G}\big)\#^{W}\Omega^{- \frac12}(\x)  \big)
+\ii \Omega^{- \frac12}\mathcal{R}_{G} \Omega^{- \frac12}
\\&
\stackrel{\eqref{espanzionecomposizioneesplicita}}{=}
\ii \Omega(D)+
\ii \opw\big(\widetilde{a}_{1/2} \big)
+\widetilde{\mathcal{R}}_1+\widetilde{\mathcal{R}}_2\,,
\end{aligned}
\]
where we defined (recall formula \eqref{restoerreenne})
\[
\begin{aligned}
\widetilde{a}_{1/2}=&\ \Omega^{- \frac12}(\x)\#^{W}_{N} \big(a_{G}+b_{G}\big)\#^{W}_{N}\Omega^{- \frac12}(\x) \,,
\\
\widetilde{\mathcal{R}}_1:=&\ \ii \opw(\mathtt{q}_1+\mathtt{q}_2)\,,
\qquad \;\; 
&&\widetilde{\mathcal{R}}_2:=\ii \Omega^{- \frac12}\mathcal{R}_{G} \Omega^{- \frac12}\,,
\\
\mathtt{q}_1:=&\ \Omega^{- \frac12}(\x)\#   r_{N}\Big(a_{G}+b_{G}\,,\, \Omega^{- \frac12}(\x) \Big) \,,  \qquad
&&\mathtt{q}_2:= r_{N}\Big(\Omega^{- \frac12}(\x)\, ,\, (a_{G}+b_{G})\#^{W}_{N}\Omega^{- \frac12}  \Big)\,.
\end{aligned}
\]
 Using estimate \eqref{stima:sharptroncato} in  Lemma \ref{lemmacomposizione}
 (using also that $\Omega(\x)$ is a Fourier multiplier), one has that there is some $\s=\s(N)\gg1$ such that
 \[
 \|\widetilde{a}_{1/2}\|_{1/2,s,\alpha}^{k_0, \gamma}\lesssim_{s,N,\al} \|a_{G}\|_{1/2,s+N,\alpha+N}^{k_0, \gamma}
 +
  \|b_{G}\|_{0,s+N,\alpha+N}^{k_0, \gamma}
  \stackrel{\eqref{stimefinalisimboliDN}}{\lesssim_{s,N,\al}} 
  \e (1 + \| {\mathcal I}_0 \|_{s + \sigma }^{k_0, \gamma})\,, 
  \quad \forall\, \alpha\in\N\,.
 \]
 This implies the first estimate in \eqref{stima cal A 2 new}. The estimates on the Lipschitz variation
 follow similarly  recalling Remark \ref{rmk:deltaunodueop}.
 
 \medskip
 \noindent
 We now claim that the remainders $\widetilde{\mathcal{R}}_i$, $i=1,2$, satisfy the bounds 
 \eqref{stima tame cal R2 new}.
 First of all, by composition and using estimates \eqref{stimefinaliRR} one gets the claim for 
 $\widetilde{\mathcal{R}}_2$.
 By applying  estimates \eqref{stimasharp} and 
 \eqref{stima:resto composizione} in  Lemma \ref{lemmacomposizione}
 with $\alpha=0$ and recalling \eqref{stimefinalisimboliDN}, 
 one gets that there is some $\s=\s(N)\gg1$ such that 
 \[
 \|\mathtt{q}_i\|_{-N+\frac{1}{2},s,0}^{k_0, \gamma}\lesssim_{s,N} 
 \e (1 + \| {\mathcal I}_0 \|_{s + \sigma }^{k_0, \gamma})\,,\;\;\;i=1,2\,.
 \]
 Therefore, by the estimate above and 
 using Lemma \ref{constantitamesimbolo} with $m_1=-M$, $m_2=\mathtt{c}$ and recalling
 $N= M-\mathtt{c}+2$ , we have that $\widetilde{\mathcal{R}}_1$ 
 satisfies \eqref{stima tame cal R2 new}.
 
 The estimates for the other summands in $\mathcal{T}_1$ and $\mathcal{T}_2$
 follows similarly, so that one gets
% 
%  
%  For $\mathcal{T}_2$ one can reason similarly.
% \eqref{espanzionecomposizioneesplicita}
%
% 
%Moreover, by using Lemmata \ref{lemma: action Sobolev}, \ref{constantitamesimbolo}, \ref{lem:pseudocomposizioSHARP}\,,
%we obtain that
\begin{equation*}
\begin{aligned}
{\mathcal T}_1 & :=  \ii \Omega(D) + \opw(\ii V\cdot \xi)
 + \ii {\rm Op}^W(a_2) + {\rm Op}^W(c_2) + {\mathcal R}_2^{(1)} \,, \\
{\mathcal T}_2 & := \ii {\rm Op}^{W} (b_2) + {\rm Op}^W(d_2) + {\mathcal R}_2^{(2)} \,,
\end{aligned}
\end{equation*}
where $a_2,b_2\in S^{\frac{1}{2}}$,
$c_2,d_2\in S^0$
satisfy the estimate 
\eqref{stima cal A 2 new}-\eqref{stima cal B 2 new} and the remainders ${\mathcal R}_2^{(1)}, {\mathcal R}_2^{(2)}$ 
satisfy the estimate \eqref{stima tame cal R2 new}. 
The properties \eqref{cal A2 Hamiltoniano 0} is verified. 
The last assertion follows since the map $\mathcal{C}$ is symplectic and momentum preserving, 
by Remark \ref{rmk:algsimboli}
and by the properties of the operator $\mathcal{L}_1$ and $G(\eta)$ (see condition \eqref{contagion1}).
\end{proof}

\begin{rmk}\label{caffe5}
Recalling the product spaces in \eqref{spazioprodperp}, with a slight abuse of notation,
we define the spaces (of complex valued functions)
\begin{equation}\label{spaziograssetto}
{\bf H}_{\perp}^{s}:=H^s(\T_{*}^{\nu+2};\C^2)\cap {\bf H}_S^\perp\,,
\end{equation}
%$\mathcal{U}$ is the subspace \eqref{calUUU}.
with ${\bf H}_S^\perp$ in \eqref{decomp siti tangenziali coordinate complesse}.
We also remark that the operator $\mathcal{L}_2$ in \eqref{elle2}
maps ${\bf H}_{\perp}^{s}$ into ${\bf H}_{\perp}^{s-1}$
and it is real-to-real according to Definition \ref{operatorerealtoreal}, since it has 
the form \eqref{divano}-\eqref{coniugio complesse simmetrizzate}.
 \end{rmk}

\subsection{Straightening of the transport operator}\label{sec:almoststraightening}
We now consider the operator $ \cL_2 $ in \eqref{elle2}.
	The aim of this section is to conjugate the variable  coefficients quasi-periodic 
transport operator
\begin{equation}\label{LTR}
	\cL_{\rm TR} := \omega\cdot\pa_\vphi \uno+  \opw\big(\ii V\cdot\x\big)\uno\,,
\end{equation}	
to a constant coefficients transport operator 
$\omega\cdot\pa_\vphi\uno + \opw(\ii \mathtt{m}_{1}\,\cdot \xi) \uno$. 
%
%where $\tn\in\N_0$ and the 
%scale $(N_{\tn})_{\tn\in\N_0}$ is defined, for $N_0>1$, by
%\begin{equation}\label{scala.strai}
%N_{\tn}:=N_0^{\chi^{\tn}}\,, \quad \chi=3/2\,, \quad N_{-1}:=1\,.
%\end{equation} 
%Such small remainder  is left 
%because we assume  only finitely many 
%non-resonance conditions, see \eqref{tDtCn}. 	

%Let us define %\red{DA VERIFICARE}
%\begin{equation}\label{tbta}
%%\tb:= [\ta] + 2 \in \N \,, \quad 
%\ta \ge 3(\tau_1+1)  \geq 1 \,, \quad \tau_1:= k_0 +(k_0+1)\tau\,.
%\end{equation}
%\red{spostare \eqref{scala.strai}, \eqref{tbta} vicino al teorema C.2}
The main result of this section is the following.

\begin{prop}{\bf (Reduction of the transport term).}\label{riduzione trasporto}
Let  $\beta_{0} \in \N$, $M \in \N$, $\mathtt c \geq 2(k_0 + \beta_0 ) + 3$ with $M \geq \mathtt c$.
%There exist $\tau_2(\tau,\nu)> \tau_1(\tau,\nu)+1 $ and 
There exists $\s:=\s(\beta_0,M ,k_0)\gg1$ 
	such that, for all $\bar{s} \geq s_0  $, 
	there is
	%are $N_0:=N_0(\bar{s})\in\N$ and 
	$\delta := \delta(\bar{s},\tau,M,\beta_0) \in (0,1) $ such that, if
	\eqref{ansatz_I0_s0}
holds with $\mu_0\gg \s$ and if (see \eqref{smallcondepsilon})
\begin{equation}\label{smallness.raddrizzo}
\varepsilon \gamma^{-1} < \delta  \,,
% N_0^{\tau_2} \varepsilon \gamma^{-1} < \delta    \,,
\end{equation}
	the following holds true. 
There exist a constant vector 	$\mathtt{m}_{1} := \mathtt{m}_{1}(\omega,\th) \in \R^d $,
operators $\mathcal{A}_3\in OPS^{1/2}$, $\mathcal{B}_3 \in OPS^0$, a 
$\mathcal{D}^{k_0}$-tame operator ${\mathcal R}_3$
and a
real-to-real, invertible, symplectic and momentum preserving 
(see \eqref{spaziograssetto})
map ${\mathcal{E}}_{\perp} : 
{\bf H}_{\perp}^{s}\to {\bf H}_{\perp}^s$
defined for all $(\omega,\mathtt{h})\in \mathtt \Omega \times[\th_1,\th_2]$ such that the following holds.
For any $(\omega,\th) $ belonging  the set
		\begin{equation}\label{tDtCn}
			\begin{aligned}
				\tT\tC_{\infty}(\gamma,\tau) 
				& :=\tT\tC_{\infty}(\gamma,\tau; i_0)\\
				&:= \Big\{ (\omega,\th)\in \mathtt \Omega \times[\th_1,\th_2] \,:\, 
				|(\omega-\tV\mathtt{m}_{1})\cdot\ell| 
				\geq  4\gamma\jap{\ell}^{-\tau} \,\ \forall\,\ell\neq 0 \Big\}
			\end{aligned}
		\end{equation}
one has the conjugation
(recall \eqref{elle2})
\begin{equation}\label{op cal L3 new}
\begin{aligned}
{\mathcal L}_3 &:= 
{\mathcal E}_{\perp}^{- 1} \circ {\mathcal L}_2 \circ {\mathcal E}_{\perp}
\\&= {\bf \Pi}_{S}^{\perp}\omega \cdot \partial_\vphi \uno+
 {\bf \Pi}_{S}^{\perp}\Big(
 \opw (\ii
\mathtt{m}_{1}\cdot\xi) \uno
+ \ii E\Omega(D)
+\cA_{3}+ \mathcal{B}_{3}+ {\mathcal R}_3
\Big) {\bf \Pi}_{S}^{\perp}\,.
\end{aligned}
\end{equation}
Moreover for any $(\omega,\mathtt{h})\in \mathtt\Omega\times[\th_1,\th_2]$
one has the following estimates:
	
\noindent
$(i)$ the constant vector $\mathtt{m}_1$ satisfies 
\begin{equation}\label{beta.FGMP.est}
\begin{aligned}
& |\mathtt{m}_{1}|^{k_0,\gamma}\lesssim\varepsilon  \,,
\quad  
| \Delta_{12} \mathtt{m}_{1} |  \lesssim 
\varepsilon \norm{i_1-i_2}_{s_0 +\sigma}  \,.
\end{aligned}
\end{equation}

\noindent
$(ii)$ The operator ${\mathcal A}_3 = {\rm Op}^W \big( A_3(\vphi, x, \xi) \big) \in OPS^{\frac12}$ has the form 
\begin{equation}\label{cal A3 Hamiltoniano}
 \begin{aligned}
& A_3(\vphi, x, \xi) = \ii \begin{pmatrix}
a_3(\vphi, x, \xi) & b_3(\vphi, x, \xi) \\
- b_3(\vphi, x, \xi) & - a_3(\vphi, x,  \xi)
\end{pmatrix}\,, \\
& a_3(\vphi, x, \xi) =  \overline{a_3(\vphi, x,  \xi)}, \qquad b_3(\vphi, x, \xi) = \overline{b_3(\vphi, x,  \xi)}\,, \\
& a_3(\vphi, x, \xi) = a_3(\vphi, x, - \xi), \quad \ b_3(\vphi, x, \xi) = b_3(\vphi, x, - \xi)\,,
\end{aligned}
\end{equation}
and satisfies the estimates 
\begin{equation}\label{stima cal A 3 new}
\begin{aligned}
\| \cA_{3} \|_{\frac12, s, \alpha}^{k_0, \gamma}  &\lesssim_{s,M,\al} 
\e\g^{-1}(1 + \| {\mathcal I}_0 \|_{s + \sigma}^{k_0, \gamma})\,, \quad \forall s_0 \leq s \leq \bar{s}\,, \quad \forall \alpha \in \N\,,
\\
\| \Delta_{12} \cA_{3} \|_{\frac12, p, \alpha}&\lesssim_{p,M,\al} \e \g^{-1}  \|i_1-i_2\|_{p+\sigma}\,,
\end{aligned}
\end{equation}
while 
the operator ${\mathcal B}_3 \in OPS^0$ satisfies the estimates
\begin{equation}\label{stima cal B 3 new}
\begin{aligned}
\| {\mathcal B}_3 \|_{0, s, \alpha}^{k_0, \gamma} &\lesssim_{s,M,\al} 
\e \g^{-1} (1 + \| {\mathcal I}_0 \|_{s + \sigma}^{k_0, \gamma})\,, 
\quad \forall s_0 \leq s \leq \bar{s}\,, \quad \forall \alpha \in \N\,,
\\
\| \Delta_{12}{\mathcal B}_3 \|_{0, p, \alpha}&\lesssim_{p,M,\al} \e \g^{-1} \|i_1-i_2\|_{p+\sigma}\,,
\end{aligned}
\end{equation}
for any $p$ as in \eqref{ps0}.

\noindent
$(iii)$
The remainder ${\mathcal R}_3$ satisfies the following property. 
For any $|\beta|\leq \beta_0$ one has that the operator 
%For any $M \in \N$, $M \le N - \red{2}(\tb + k_{0})-3$, 
%and $|\beta|\leq\tb$,
%one has that 
$ \langle D \rangle^{M} \partial_\vphi^\beta{\mathcal R}_3 \langle D \rangle^{-\mathtt{c}}$ 
is a ${\mathcal D}^{k_0}$-tame operator  (see Def. \ref{Dksigmatame}) 
satisfying 
\begin{equation}\label{stima tame cal R3 new}
\begin{aligned}
{\mathfrak M}_{ \langle D \rangle^{M}
\partial_\vphi^\beta {\mathcal R}_3 \langle D \rangle^{-\mathtt{c}}}(s) 
&\lesssim_{\bar{s},M} \e \gamma^{-1}  (1 + \| {\mathcal I}_0 \|_{s + \sigma})\,, 
\quad \forall s_0 \leq s \leq \bar{s}\,,
\\
\|  \langle D \rangle^{M} 
\partial_\vphi^\beta \Delta_{12}{\mathcal R}_3\langle D \rangle^{-\mathtt{c}} 
\|_{\cL(H^{p} , H^{p})}& \lesssim_{\bar{s},M}\e \gamma^{-1}  \|i_1-i_2\|_{p+\s}\,.
\end{aligned}
\end{equation}

\noindent
$(iv)$
The operators
${\mathcal E}_{\perp}^{\pm1} $ are ${\mathcal D}^{k_0}$-$(k_0+1)$-tame, the operators
$\mathcal{E}_{\perp}^{\pm1}-\id$ %, (\mathcal{E}_{\perp}^{\pm1}-\id)^{*}$
are ${\mathcal D}^{k_0}$-$(k_0+2)$-tame and 
satisfy
\begin{equation}\label{stimamappaEperp}
\begin{aligned}
\mathfrak{M}_{\mathcal{E}_{\perp}^{\pm1}}(s)
&\lesssim_{\bar{s}}
1+\|\mathcal{I}_{0}\|_{s+\s}^{k_0,\gamma}\,,
\qquad
\mathfrak{M}_{(\mathcal{E}_{\perp}^{\pm1}-{\rm Id})}(s)
%\,,\;
%\mathfrak{M}_{(\mathcal{E}_{\perp}^{\pm1}-{\rm Id})^{*}}(s)
&\lesssim_{\bar{s}}
\e \g^{-1}(1+\|\mathcal{I}_{0}\|_{s+\s}^{k_0,\gamma})\,,
\\
\|\Delta_{12}\mathcal{E}_{\perp}^{\pm1}h\|_{p}&\lesssim_{\bar{s}}
\e\gamma^{-1}\|h\|_{p+\s}\|i_1-i_{2}\|_{p+\s}\,.
%\quad \mathcal{A}\in \{\mathcal{E}_{\perp}^{\pm1}, (\mathcal{E}_{\perp}^{\pm1})^{*}\}\,.
\end{aligned}
\end{equation}
Finally, the operators $\mathcal{L}_3$,
$\mathcal{A}_{3}$, $\cB_{3}$ and $\cR_{3}$ are real-to-real, 
Hamiltonian and momentum preserving.
\end{prop}

The rest of the section is devoted to the proof of the result above.

\noindent
We start by  conjugating the  $ \cL_{\rm TR} $ in \eqref{LTR} 
through a \emph{symplectic} 
(Definition \ref{def simplettica complesse}) transformation
induced by the diffeomorphism of the torus
$\mathbb{T}_{\Gamma}^{d}\ni x\mapsto x+\breve\beta(\vphi,x)$
for some smooth, small quasi-periodic traveling wave $\beta:\T^\nu\times \T^d_\Gamma\to \R $
to be determined.
%
%Given some smooth function $\beta(\vphi,x)$ (to be determined) we define
%\begin{equation}\label{defcE}
%\cE:= \begin{pmatrix}
%	\cA & 0 \\ 0 & \cA
%\end{pmatrix} \, , \quad 
%\cE^{-1}:=\begin{pmatrix}
%	\cA^{-1} & 0 \\ 0 & \cA^{-1}
%\end{pmatrix}
%\end{equation}
%where
%\[
%\cA h(\varphi, x)= \cM \circ T_{\beta} h(\varphi, x) \,, \quad \varphi\in\T^\nu, x\in\T_{\Gamma}^d\,,
%\]
%where $\cM$ is the multiplication operator
%defined as \begin{equation}\label{def:CCalso}
%\cM h(\varphi, x) = \sqrt{\det(\uno+\nabla_x \beta(\varphi, x))} h(\varphi, x)\\
%\end{equation}
%and $T_{\beta}$ is induced by the $\vphi-$dependent diffeomorphism 
%$y=x+\beta(\vphi, x)$ of the torus $\T^d_\Gamma$, 
%for some small quasi-periodic traveling wave $\beta:\T^\nu\times \T^d_\Gamma\to \R $.

The next lemma provide a suitable choice of the function $\beta$.

\begin{lemma}\label{conju.tr}
{\bf (Straightening of the transport operator)}
Under the assumptions of Proposition \ref{riduzione trasporto}
there exist
a  constant vector $\mathtt{m}_{1} := \mathtt{m}_{1}(\omega,\th) \in \R^d $, satisfying
\eqref{beta.FGMP.est} and
a quasi-periodic traveling wave
%$\breve{\beta}(\vphi,x)$,
 $\beta(\vphi,x):=\beta_{\infty}(\vphi,x)$, 
defined for all $(\omega,\th)\in \mathtt\Omega\times[\th_1,\th_2]$,
and 
for some $\sigma=\sigma(\nu,k_0)>0$, satisfying
\begin{equation}\label{beta.FGMP.estBETA}
\begin{aligned}
\normk{\beta}{s} &\lesssim_{s} 
\varepsilon\gamma^{-1} (1+\normk{\cI_0}{s+\sigma}) \,, 
\quad  \forall \, s_0 \leq s \leq \bar{s} \, , 
\\
    \| \Delta_{12}\beta\|_{p} &\lesssim_{p} \e \gamma^{-1}\|i_1-i_2\|_{p+\sigma}\,,
\end{aligned}			
\end{equation}
for any $p$ as in \eqref{ps0},  such that the following holds.

\noindent
(i) The diffeomorphism $x \mapsto x + \beta(\vphi, x)$ is invertible  with inverse 
$y \mapsto y + \breve \beta(\vphi, y)$ where $\breve{\beta}$ 
satisfies bound like \eqref{beta.FGMP.estBETA}.

\noindent
(ii) Let us define the map, for $\varphi\in\T^\nu, x\in\T_{\Gamma}^d$,
\begin{equation}\label{def:acheck}
\breve{\cA} h(\varphi, x)= \breve{\cM} \circ T_{\breve\beta} h(\varphi, x) \,, 
\qquad {\rm with}\qquad 
\begin{aligned}
T_{\breve\beta} h(\varphi, x) &= h(\varphi, x+ \breve{\beta}(\varphi, x))\,,
\\
\breve{\cM} h(\varphi, x) &= \sqrt{\det(\uno+\nabla_x \breve{\beta}(\varphi, x))} h(\varphi, x)\,.
\end{aligned}
\end{equation}
For any $(\omega,\th)\in \tT\tC_{\infty}(2\gamma,\tau)$ in \eqref{tDtCn}
the operator $\cL_{\rm TR}$ in \eqref{LTR} is conjugated 
%under the map $\mathcal{E}$ in \eqref{defcE}, to
\begin{equation}\label{coniugaLTR}
\breve{\cE}\circ \cL_{\rm TR} \circ\breve{\cE}^{-1} = 
\omega\cdot\pa_\vphi + \opw (\ii \mathtt{m}_{1}\cdot \xi)  \uno\,,
\qquad \breve{\mathcal{E}}:=\begin{pmatrix}
\breve{\cA} & 0 \\ 0 & \breve{\cA}
\end{pmatrix}\,.
\end{equation}
Finally, the maps $\breve{\mathcal{E}}^{\pm1}$ satisfy similar 
estimates as in \eqref{stimamappaEperp}.
\end{lemma}
\begin{proof}
The proof of the result is postponed  to Appendix \ref{dim:lemma9.8}.
\end{proof}

\begin{rmk}\label{rmk:EedEbreve}
We remark the following facts.

\noindent
$(i)$ The map $\breve{\mathcal{A}}$ in Lemma \ref{conju.tr}
can be seen as $\breve{\mathcal{A}}\equiv(\breve{\mathcal{A}}^{\tau})_{|\tau=1}$ where
$\breve{\mathcal{A}}^{\tau}$ is the flow at time $\tau\in [0,1]$ of the problem 
\eqref{gianduiotto} with generator given by $b\rightsquigarrow \breve{b}$ defined as
\begin{equation}\label{breveBB}
\breve{b}(\tau; \varphi, x):=(\uno+\tau \nabla_x\breve{\beta}(\varphi, x))^{-T}\breve\beta(\varphi, x)\,.
\end{equation}

\noindent
$(ii)$
Consider the function $\beta$ given by Lemma \ref{conju.tr}.
By defining 
\begin{equation}\label{defcE}
\cE:= \begin{pmatrix}
	\cA & 0 \\ 0 & \cA
\end{pmatrix} \, , 
%\quad 
%\cE^{-1}:=\begin{pmatrix}
%	\cA^{-1} & 0 \\ 0 & \cA^{-1}
%\end{pmatrix}
\qquad
{\cA} h(\varphi, x)= {\cM} \circ T_{\beta} h(\varphi, y) \,, 
\end{equation}
where
\begin{equation}\label{def:acheck2}
\begin{aligned}
T_{\beta} h(\varphi, y) &= h(\varphi, y+ {\beta}(\varphi, y))\,,
\qquad 
{\cM} h(\varphi, y) = \sqrt{\det(\uno+\nabla_y {\beta}(\varphi, y))} h(\varphi, y)\,,
\end{aligned}
\end{equation}
one can note that $\mathcal{\mathcal{E}}^{-1}=\breve{\mathcal{E}}$.
\end{rmk}

\begin{lemma}\label{mappaEsymplectic}
There exists a real-to-real, momentum preserving, invertible 
map $\breve\cE_{\perp}:  {\bf H}_{{\perp}}^{s}
\to {\bf H}_{{\perp}}^{s} $
which is  symplectic on ${\bf H}_{{\perp}}^{s}$
%(i.e. w.r.t. the symplectic form \eqref{estesaPois})
such that
\[
\breve\cE_{\perp}=\Pi_{S}^{\perp} \breve\cE \Pi_{S}^{\perp}\circ(\uno +\mathcal{R})\,,
\]
where $\mathcal{R}$ is a finite rank operator of the form \eqref{forma buona resto}
satisfying \eqref{giallo2}-\eqref{giallo3}. In particular, one has
\begin{equation}\label{peperoni}
\begin{aligned}
\mathfrak{M}_{(\breve{\cE}_{\perp}^{\pm1}-{\rm Id})}(s)
%+\mathfrak{M}_{(\breve{\cE}_{\perp}^{\pm}-{\rm Id})^{*}}(s)
&\lesssim_{\bar{s}}
\e \g^{-1}(1+\|\mathcal{I}_{0}\|_{s+\s}^{k_0,\gamma}) 
\quad  \forall \, s_0 \leq s \leq \bar{s} \, , 
\\
\|\Delta_{12}(\breve{ \cE}_{\perp}^{\pm1}-{\rm Id})h\|_{p}
%+\|\Delta_{12}(\breve{\cE}_{\perp}^{\pm}-{\rm Id})^{*}h\|_{p}
&\lesssim_{\bar{s}}
\e \g^{-1} \|h\|_{p}\|i_1-i_{2}\|_{p+\s}\,,
\end{aligned}
\end{equation}
for any $p$ as in \eqref{ps0}.
\end{lemma}
\begin{proof} 
Recalling Remark \ref{rmk:EedEbreve}-(i),
the map $\breve\cE$ in Lemma \ref{conju.tr}
is constructed as the flow generated by $\opw(\ii \breve{b} \cdot \xi \uno)$ (see \eqref{gianduiotto}) 
with symbol $\breve{b}$ defined in 
\eqref{breveBB}. Therefore, the thesis follows just taking $\breve{\cE}_{\perp}$ as the flow generated
by the truncation of $\opw(\ii \breve{b} \cdot \xi \uno)$ (see \eqref{natale})
and 
by applying 
Lemma \ref{differenzaFlussi} (recalling also Remark \ref{differenzaFlussiRMK}).
The bounds \eqref{peperoni} follow recalling by the fact that $\breve{\cE}$ satisfies estimates like 
\eqref{stimamappaEperp}.
The map is momentum preserving since the generator preserves the subspace $S_{\tV}$
(see Lemma \ref{lem:momentoflusso} and Remark \ref{rmk:EedEbreve}-(i)).
\end{proof}

\begin{proof}[Proof of Proposition \ref{riduzione trasporto}]
Consider the map $\breve{\cE}_{\perp}$ in Lemma \ref{mappaEsymplectic}
and recall Remark \ref{rmk:EedEbreve}-(i).
First of all, by applying Lemma \ref{georgiaLem} and recalling 
\eqref{elle2}, we notice that
\begin{equation}\label{trota1}
{\mathcal L}_3 := 
\breve{\cE}_{\perp} \circ {\mathcal L}_2 \circ \breve{\cE}_{\perp}^{- 1}=
\Pi_{S}^{\perp}\mathcal{P}\Pi_{S}^{\perp}+Q_{3}\,,
\end{equation}
where $Q_{3}$ is a   finite rank operators satisfying the estimates \eqref{stima tame cal R3 new} and
where 
\[
\mathcal{P}:=
\breve\cE \circ\Big(\omega\cdot\pa_{\vphi}+
%+V\cdot\nabla \uno
 \opw(\ii V\cdot \xi) \uno 
 + \ii E \Omega(D)  + {\mathcal A}_2 
 + {\mathcal B}_2 
 + {\mathcal R}_2\Big) \circ \breve{\cE}^{- 1}\,,
\]
where $\breve{\cE}$ is given by Lemma \ref{conju.tr}.
We also remark that, thanks to the properties of $\breve{\cE}$, the operator $\mathcal{L}_{3}$
is real-to-real, Hamiltonian and momentum preserving.
The same holds true for the operator $\mathcal{P}$ acting on the whole space $H^{s}\times H^{s}$.
Therefore, we only have to study the structure of $\mathcal{P}$ which has the form
\begin{align}
{\mathcal P} & =
\breve{\cE} \circ \cL_{TR} \circ \breve{\cE}^{- 1}
\label{pepe1}
\\&
+\breve{\cE} \circ  \Big(  \ii E \Omega(D)  + {\mathcal A}_2 
 + {\mathcal B}_2 \Big)  \circ \breve{\cE}^{- 1}
\label{pepe2}
 \\&
+ \breve{\cE} \circ {\mathcal R}_2 \circ \breve{\cE}^{- 1}\,. \label{pepe3}
\end{align}
We analyze separately all the terms in the above expression. 

\smallskip
\noindent
{\bf Analysis of \eqref{pepe1}.}
By Lemma \ref{conju.tr} we have that 
\[
\eqref{pepe1}= \omega\cdot\pa_\vphi +\opw( \ii\mathtt{m}_{1}\cdot\xi)
%+\mathcal{F}_1
\,,
\]
where $\mathtt{m}_{1}$ satisfies \eqref{beta.FGMP.est}.

\smallskip
\noindent
{\bf Analysis of \eqref{pepe2}.}  
%This is the most crucial term.
We start by considering the conjugate of $\ii E\Omega(D)$.
Recalling \eqref{def:acheck}-\eqref{coniugaLTR} we note that
\[
\ii \breve{\mathcal{E}}E \Omega(D)\breve{\mathcal{E}}^{-1}=\ii \left(\begin{matrix}
\breve{\mathcal{A}}\opw(\Omega(\x))\breve{\mathcal{A}}^{-1} & 0 
\\ 0 & -\breve{\mathcal{A}}\opw(\Omega(\x))\breve{\mathcal{A}}^{-1}
\end{matrix}\right)\,.
\]
In view of Remark \ref{rmk:EedEbreve}-(i) one has that $\breve{\mathcal{A}}$ can be seen as the time one flow map 
of \eqref{gianduiotto} with $\breve{b}$ in \eqref{breveBB}.
Therefore, we shall apply Theorem \ref{quantitativeegorov}
with $\beta\rightsquigarrow \breve{\beta}$, $w(\vphi,x,\x)\rightsquigarrow \Omega(\x)$.
We obtain that
\begin{equation}\label{coniugoOmega}
\breve{\mathcal{A}}\opw(\Omega(\x))\breve{\mathcal{A}}^{-1}=
\opw(\Omega(\x)+q_0+q_1) + R\,,
\end{equation}
where $R$ is some smoothing remainder, $q_1\in S^{0}$ and 
\[
q_0:=\Omega\Big( (\uno+\nabla_y{\beta}(1, \varphi, y))^{T}_{|y=x+\breve{\beta}(\varphi, x)}\xi \Big)-\Omega(\x)\,.
\]
Using then fact that $\|\Omega(\x)\|_{1/2,s,\alpha}\lesssim_{s,\alpha}1$ and estimates
\eqref{parlare}, \eqref{beta.FGMP.estBETA}, we deduce
\[
\|q_1\|_{1/2,s,\alpha}^{k_0,\gamma}+
\|q_0\|_{-1/2,s,\alpha}^{k_0,\gamma}
\lesssim_{s,M,\alpha}\varepsilon\gamma^{-1} (1+\normk{\cI_0}{s+\sigma}) \,,
\;\;\;\forall\, \alpha\in\N\,.
\]
This implies estimates \eqref{stima cal A 3 new}-\eqref{stima cal B 3 new}.
The bounds on the Lipschitz variation follow similarly by using \eqref{troppo}.
In the same way, using \eqref{francia1}-\eqref{francia2} one gets that the remainder 
$R$ in \eqref{coniugoOmega} satisfies the bounds \eqref{stima tame cal R3 new}.
Let us now consider the conjugate  $\breve{\cE}\mathcal{A}_2\breve{\cE}^{- 1}$.
Recalling \eqref{coniugaLTR} and  \eqref{cal A2 Hamiltoniano 0}, we have
\begin{equation}\label{bandalarga1}
\breve{\cE}\mathcal{A}_2\breve{\cE}^{- 1}=
\left(
\begin{matrix}
\breve{\mathcal{A}}\opw(a_2)\breve{\mathcal{A}}^{-1} &
 \breve{\mathcal{A}}\opw(b_2)\breve{\mathcal{A}}^{-1}
\\[0.3em]
-\breve{\mathcal{A}}\opw(b_2)\breve{\mathcal{A}}^{-1} & -\breve{\mathcal{A}}\opw(a_2)\breve{\mathcal{A}}^{-1}
\end{matrix}
\right)\,.
\end{equation}
Let us consider $\breve{\mathcal{A}}\opw(a_2)\breve{\mathcal{A}}^{-1}$. 
We apply Theorem \ref{quantitativeegorov} with $w\rightsquigarrow a_2$ and we get
\[
\breve{\mathcal{A}}\opw(a_2)\breve{\mathcal{A}}^{-1}=\opw(\widetilde{a}_2)+\widetilde{R}\,,
\]
where $\widetilde{a}_2\in S^{1/2}$ is some symbol satisfying 
\begin{equation}\label{cappotto1}
\begin{aligned}
\|\widetilde{a}_2\|_{1/2,s,\alpha}^{k_0,\gamma}
&\stackrel{\eqref{parlare}}{\lesssim_{s,M,\alpha}}
\|a_2\|_{1/2,s,\alpha+\sigma_1}^{k_{0}, \g} + \sum_{s}^* \|a_2\|_{m,k_1,\alpha+k_2+\sigma_1}^{k_{0}, \g} 
\|\breve\beta\|^{k_{0}, \g}_{k_3+\sigma_1}
\\&
\stackrel{\eqref{stima cal A 2 new}, \eqref{beta.FGMP.estBETA}}{\lesssim_{s,M,\alpha}}
\e (1 + \| {\mathcal I}_0 \|_{s + \sigma }^{k_0, \gamma})
+\sum_{s}^*  \e (1 + \| {\mathcal I}_0 \|_{k_1+ \sigma }^{k_0, \gamma})
\e\gamma^{-1} (1 + \| {\mathcal I}_0 \|_{k_3+ \sigma }^{k_0, \gamma})\,,
\end{aligned}
\end{equation}
for some $\s\gg s_0$ large enough. 
Now, recalling the definition of $\sum_{s}^{*}$ above of Theorem \ref{quantitativeegorov}
we note that $k_1+k_3\leq s$ since $k_1+k_2+k_3=s$.
Moreover, using the smallness assumption \eqref{ansatz_I0_s0}, we deduce that,
in the case $k_1=0$ or $k_3=0$,  
one has
$\| {\mathcal I}_0 \|_{k_1+ \sigma }^{k_0, \gamma} \| {\mathcal I}_0 \|_{k_3+ \sigma }^{k_0, \gamma}\lesssim_{s}
 \| {\mathcal I}_0 \|_{s+ \sigma }^{k_0, \gamma} $.
 On the other hand, in the case $1\leq k_1,k_3$ and $k_1+k_3\leq s$, using 
 interpolation estimates
on Sobolev spaces (see Lemma \ref{lemma:interpolation}) one gets 
%
%
%Now, using that $k_1+k_2+k_3=s$ and  interpolation estimates
%on Sobolev spaces (see Lemma \ref{lemma:interpolation}) one gets 
\[
 \| {\mathcal I}_0 \|_{k_1+ \sigma }^{k_0, \gamma} \| {\mathcal I}_0 \|_{k_3+ \sigma }^{k_0, \gamma}
 \stackrel{\eqref{2202.2}}{\lesssim_{s}}
  \| {\mathcal I}_0 \|_{k_1+k_3+ \sigma }^{k_0, \gamma} \| {\mathcal I}_0 \|_{ \sigma }^{k_0, \gamma}
  \lesssim_{s} \| {\mathcal I}_0 \|_{s+ \sigma }^{k_0, \gamma}\,.
\]
The latter bound, together with \eqref{cappotto1}, implies 
that $\widetilde{a}_2$ satisfies \eqref{stima cal A 3 new}. The estimate on the Lipschitz variation follows using 
\eqref{troppo}. By reasoning similarly and using bounds \eqref{francia1}-\eqref{francia2}, one concludes that 
$\widetilde{R}$ satisfies \eqref{stima tame cal R3 new}.
The conjugate of $\opw(b_2)$ in \eqref{bandalarga1}
and $\breve{\cE}\mathcal{B}_2\breve{\cE}^{- 1}$ in \eqref{pepe2}
can be treated similarly (recalling  \eqref{stima cal B 2 new}).
%
%\noindent
%By applying Theorem \ref{quantitativeegorov} and estimates 
% \eqref{stima cal A 2 new}, \ref{stima cal B 2 new} 
 Hence, one gets that 
\[
\eqref{pepe2}=\ii E \Omega(D) + \cA_{3} + \cB_{3} +  {\mathcal F}_1
\]
where $\cA_{3}\in OPS^{1/2}$ satisfies \eqref{stima cal A 3 new}, $\cB_{3}\in OPS^{0}$ satisfies \eqref{stima cal B 3 new} 
and  $\mathcal{F}_1$ is some smoothing reminder satisfying  \eqref{stima tame cal R3 new}.
Moreover, by \eqref{cal A2 Hamiltoniano 0} we have that \eqref{cal A3 Hamiltoniano} holds.

\smallskip
\noindent
{\bf Analysis of \eqref{pepe3}.}  
By Lemma \ref{conj.NEWsmoothresto} ({applied with  $m=1$}), using estimates \eqref{stima tame cal R2 new} on $\mathcal{R}_2$, we have that 
\[
\eqref{pepe3} = \cF_{2}
\]
with $\cF_{2}$ a smoothing reminder satisfying  \eqref{stima tame cal R3 new}.

By the discussion above, formula \eqref{trota1} and setting 
$\mathcal{E}_{\perp}:=\breve{\mathcal{E}}_{\perp}^{-1}$ one gets the conjugation 
\eqref{op cal L3 new}.

 \smallskip
 \noindent
 In view of Remark \ref{rmk:algsimboli} we have that $\mathcal{A}_{3}$ and $\cB_{3}$ are momentum preserving thanks to Lemma \ref{conju.tr} 
 and Proposition \ref{prop operatore cal L2}. Then by difference also $\cR_{3}$ is so.
 The discussion above implies the thesis.
%
%The claimed statement then follows by defining ${\mathcal A}_3 := {\mathcal F}_1 + {\mathcal F}_2 + {\mathcal F}_3$ and ${\mathcal R}_3 := \Phi_{\mathcal M}^{- 1} \circ {\mathcal R}_2 \circ  \Phi_{\mathcal M}$. 
\end{proof}

\section{Symmetrization of the linearized operator at the lower order}\label{sym.low.order}

\subsection{Block-Diagonalization at order 1/2}\label{sym.low.orderZERO}
In this subsection we block diagonalize the order $1/2$ 
appearing in the operator ${\mathcal L}_3$ in \eqref{op cal L3 new}. Let $\chi_0 \in C^\infty(\R^d, \R)$ be a cut-off function satisfying 
\begin{equation}\label{cut-offCHIZERO}
\begin{aligned}
& 0 \leq \chi_0\leq 1, \quad \chi_0(\xi) = \chi_0(- \xi)\,, \\
&  \chi_0(\xi) = 0 \quad \text{if} \quad |\xi| \leq 1\,, 
%\\& 
\qquad 
\chi_0(\xi) = 1 \quad \text{if} \quad |\xi| \geq 2\,. 
\end{aligned}
\end{equation}
We split the operator  ${\mathcal A}_3 = {\rm Op}(A_3)$ in \eqref{cal A3 Hamiltoniano} as 
\begin{equation}
A_3(\vphi, x, \xi) = \chi_0(\xi) A_3(\vphi, x, \xi) + \big( 1 - \chi_0(\xi)\big) A_3(\vphi, x, \xi)\,. 
\end{equation}
Since $1 - \chi_0 \equiv 0$ for $|\xi| \geq 2$ and by using \eqref{stima cal A 3 new}, one obtains that 
\begin{equation}\label{A3 cut off resto}
\begin{aligned}
& (1 - \chi_0) A_3 \in S^{- \infty}\,, \quad  \\
& \| (1 - \chi_0) A_3 \|_{- N, s, \alpha}^{k_0, \gamma} \lesssim_{N, s, \alpha} \e \gamma^{- 1} (1 + \| {\mathcal I}_0 \|_{s + \sigma}^{k_0, \gamma})\,, \quad 
s_0\leq  s \leq \bar{s}\,, \quad \forall N, \alpha \in \N\,, \\
&  \| (1 - \chi_0) \Delta_{12} A_3 \|_{- N, p, \alpha} \lesssim_{N, p, \alpha} \e \gamma^{- 1} \| i_1 - i_2 \|_{p + \sigma}, \quad \forall N, \alpha \in \N \,,
\end{aligned}
\end{equation}
for any $p \geq s_0$ as in \eqref{ps0} , for some $\sigma \gg 0$ large enough. Then, we write
\begin{equation}\label{caffe10}
\begin{aligned}
\ii E \Omega(\xi) + A_3(\vphi, x, \xi) & = C(\vphi,x,\x) + (1 - \chi_0(\xi)) A_3(\vphi, x, \xi)\,, \\
C(\vphi, x, \xi) & := \ii E\Omega(\x) + \chi_0(\xi) A_3(\vphi, x, \xi)\,,
\end{aligned}
\end{equation}
(where $\Omega(\x)$ is in \eqref{Omega D cut off}) and we shall diagonalize the matrix $C(\vphi, x, \xi)$ (since $(1 - \chi_0) A_3 \in S^{- \infty}$). 
More precisely, the matrix $C(\vphi, x, \xi)$ has eigenvalues 
\begin{equation}\label{autovalori matrice C 1/2}
\begin{aligned}
 \lambda_\pm(\vphi, x, \xi) &  = \pm \ii  \lambda(\vphi, x, \xi), \quad \lambda(\vphi, x, \xi)  :=  { \sqrt{\Delta(\vphi, x, \xi)}} \\
 \Delta (\vphi, x, \xi) & :=  \Big(\Omega(\xi)   + \chi_0(\xi)a_3(\vphi, x, \xi) \Big)^2   -  \chi_0(\xi)^2b_3(\vphi, x, \xi)^2\,.
\end{aligned}
\end{equation}
The following Lemma holds. 
\begin{lemma}\label{lemma autovalori decoupling ordine alto}
One has that the symbol $\lambda(\omega, \mathtt h; \vphi, x, \xi) \in S^{\frac12}$ is a real valued, momentum preserving symbol defined for any $(\omega, \mathtt h) \in \mathtt \Omega \times [\mathtt h_1, \mathtt h_2]$ and 
$$
\begin{aligned}
 \lambda (\vphi, x, \xi)  &= \Omega(\xi) + r_\lambda(\vphi, x, \xi), \quad r_\lambda \in S^{\frac12}\,, \\
 r_\lambda(\vphi, x, \xi) &\equiv 0 , \quad \text{for } \quad |\xi| \leq 1\,, \\
 \lambda(\vphi, x, \xi) &= \lambda(\vphi, x, - \xi), \quad \forall (\vphi, x, \xi) \in \T^\nu \times \T^d_\Gamma \times \R^d\,. 
\end{aligned}
$$
%Moreover one has
%\begin{equation}\label{prop simmetria lambda +- e fpm}
%\begin{aligned}
%&\lambda_{\pm}(\vphi, x, \xi) = \lambda_\pm(\vphi, x, - \xi)\,, \quad
%   \lambda_-(\vphi, x, \xi) = - \lambda_+(\vphi, x,  \xi)\,. 
%\end{aligned}
%\end{equation}
Furthermore, for any $\bar{s} \geq s_0$, $\alpha \in \N$, 
$$
\begin{aligned}
 \| \lambda - \Omega \|_{\frac12, s, \alpha}^{k_0, \gamma} 
&= \|r_\lambda \|_{\frac12, s, \alpha}^{k_0, \gamma}   
\lesssim_{s,M, \alpha} 
\e \gamma^{- 1}(1 + \|{\mathcal I}_0 \|_{s + \sigma}^{k_0, \gamma})\,, 
\;\;\;s_0\leq s\leq \bar{s}\,,
\\
 \| \Delta_{12} \lambda   \|_{\frac12, p, \alpha} &= \| \Delta_{12} r_\lambda   \|_{\frac12, p, \alpha} 
\lesssim_{p, M,\alpha} 
\e \gamma^{- 1} \| i_1 - i_2 \|_{p + \sigma}\,,
\end{aligned}
$$
for any $p \geq s_0$ as in \eqref{ps0} , for some $\sigma \gg 0$ large enough, and where $M$ 
is as in Proposition \ref{riduzione trasporto}.
\end{lemma}
\begin{proof}
By the properties of the cut-off function $\chi_0$ in \eqref{cut-offCHIZERO}, one has that $\Omega (\xi) \gtrsim \langle \xi \rangle^{\frac12} \gtrsim 1$, for $|\xi| \geq 1$ which is the support of $\chi_0$. 
Hence, 
 the desired decomposition of $\lambda$ follows by \eqref{autovalori matrice C 1/2}. 
 Since $a_{3}$ and $b_{3}$ are momentum preserving then $\lambda$ is so. 
 The claimed estimates follow by applying \eqref{stima cal A 3 new}.
 \end{proof}
%\red{OPZIONE: considera il cut-off \eqref{cut-offCHIZERO} dire voglio cancellare la matrice di simboli
%$\Omega(\x)+\chi_0(\x)A_{3}$}
By a direct calculation, the normalized eigenvectors $u_{\pm}(\vphi, x, \xi)$ of the matrix $C(\vphi, x, \xi)$ are given by 
\begin{equation}\label{f pm lambda pm}
\begin{aligned}
u_+(\vphi, x, \xi) & = \begin{pmatrix}
f_{+}(\vphi, x, \xi) \\
f_{-}(\vphi, x, \xi)
\end{pmatrix} , \quad u_-(\vphi, x, \xi)  = \begin{pmatrix}
f_{-}(\vphi, x, \xi)  \\
f_{+}(\vphi, x, \xi)
\end{pmatrix}
\end{aligned}
\end{equation}
where 
\begin{equation}\label{def f g u + u -}
\begin{aligned}
f_{+}(\vphi, x, \xi) & := \dfrac{\Omega(\xi) + \chi_0(\xi)a_3(\vphi, x,  \xi) + \lambda(\vphi, x, \xi)}{\sqrt{\big( \Omega(\xi) + \chi_0(\xi)a_3(\vphi, x,  \xi) + \lambda(\vphi, x, \xi)\big)^2 - \chi_0(\xi)^2b_3(\vphi, x, \xi)^2}}  
\\
f_{-} (\vphi, x, \xi) & := 
- \dfrac{{\chi_0(\xi)b_3(\vphi, x, \xi)}}{\sqrt{\Big( \Omega(\xi) + \chi_0(\xi)a_3(\vphi, x,  \xi) 
+ \lambda(\vphi, x, \xi)  \Big)^2 - \chi_0(\xi)^2b_3(\vphi, x, \xi)^2}} 
\end{aligned}
\end{equation}
and we define the matrix 
\begin{equation}\label{defmatrixM}
F(\vphi, x, \xi) :=  \begin{pmatrix}
f_{+}(\vphi, x, \xi) & f_{-}(\vphi, x, \xi) \\
f_{-}(\vphi, x, \xi) & f_{+} (\vphi, x, \xi)
\end{pmatrix}\,. 
\end{equation}
In the next lemma we shall verify that the latter expressions are well defined and that $f_+, f_- \in S^0$ with $f_+ - 1 = O(\e \gamma^{- 1})$ and $f_- = O(\e \gamma^{- 1})$.

\begin{lemma}\label{lem:lambdino}
The symbols $ f_{\pm}$ in \eqref{autovalori matrice C 1/2} 
and \eqref{def f g u + u -}
are real valued, momentum preserving (according to Def. \ref{def:mompressimbo}) 
symbols respectively  and $S^{0}$ defined for any $(\omega,\th)\in \Omega \times [\mathtt h_1, \mathtt h_2]$. 
Moreover, under the assumptions of Prop. \ref{riduzione trasporto}, they satisfy the estimates 
\begin{equation}\label{stime lambda pm fpm}
\begin{aligned}
 \| f_{+} - 1 \|_{0, s, \alpha}^{k_0, \gamma}\,,  \| f_{-}  \|_{0, s, \alpha}
& \lesssim_{s,M, \alpha} \e \g^{-1} (1 + \| {\mathcal I}_0 \|_{s + \sigma}^{k_0, \gamma})\,, 
\quad s_0\leq s \leq \bar{s} \,,
 \\
  \| \Delta_{12}f_{\pm} \|_{0, p, \alpha}
 &\lesssim_{p,M,\al} \e  \g^{-1} \|i_1-i_2\|_{p+\sigma}\,
 \end{aligned}
\end{equation}
for any $p \geq s_0$ as in \eqref{ps0} , for some $\sigma \gg 0$ large enough. 
 The symbols $f_{\pm}$ satisfy the symmetry condition
  \begin{equation}\label{prop simmetria lambda +- e fpm}
   f_{\pm}(\vphi, x, \xi) = f_{\pm}(\vphi, x, - \xi)\,.
   \end{equation}
The matrix $F(\vphi, x, \xi)$ in \eqref{defmatrixM} is invertible and it satisfies 
\begin{align}
&{\rm det}\Big( F(\vphi, x, \xi)\Big) = f_{+}(\vphi, x, \xi)^2 - f_{-}(\vphi, x, \xi)^2 = 1\,.\label{prop determinante M}
\end{align}
Moreover $F^{\pm 1}$ are close to the identity and they satisfy the estimates
\begin{equation}\label{stime M M inv ordine alto}
\begin{aligned}
\| F^{\pm1} - \uno \|_{0, s, \alpha}^{k_0, \gamma}&\lesssim_{s, M,\alpha} 
\e \g^{-1} (1 + \| {\mathcal I}_0 \|_{s + \sigma}^{k_0, \gamma})\,,\quad s_0\leq s \leq \bar{s}\,,
\\
 \| \Delta_{12}F^{\pm1} \|_{0, p, \alpha}
 &\lesssim_{p,M,\al} \e  \g^{-1} \|i_1-i_2\|_{p+\sigma}\,,
\end{aligned}
\end{equation}
for any $p \geq s_0$ as in \eqref{ps0} , for some $\sigma \gg 0$ large enough. 
 Finally ,
\begin{align}
&F(\vphi, x, \xi)^{- 1} \Big[\ii E \Omega(\xi) + A_3(\vphi,x, \xi) \Big]F(\vphi, x, \xi)  \\
&  \stackrel{\eqref{autovalori matrice C 1/2}}{=}
%\stackrel{\eqref{prop simmetria lambda +-}}{=} 
\begin{pmatrix}
\ii \lambda(\vphi, x, \xi) & 0 \\
0 & - \ii \lambda(\vphi, x, \xi)
\end{pmatrix} + R_{\infty, \lambda}(\vphi,x, \xi)  \,,\label{C M inv M}
\end{align}
where $R_{\infty, \lambda} \in S^{- \infty}$ and 
\begin{equation}\label{stima R infty lambda dec ord max}
\begin{aligned}
\| R_{\infty, \lambda} \|_{- N, s, \alpha}^{k_0, \gamma} & \lesssim_{N, s,M, \alpha} 
\e \gamma^{- 1} (1 + \| {\mathcal I}_0 \|_{s + \sigma}^{k_0, \gamma})\,, 
\quad s_0\leq s \leq \bar{s}\,, %\quad \forall N, \alpha \in \N\,, 
\\
\| \Delta_{12} R_{\infty, \lambda} \|_{- N, p, \alpha} & \lesssim_{N, p, M,\alpha}  \e \gamma^{- 1} \| i_1 - i_2 \|_{p + \sigma}, 
\end{aligned}
\end{equation}
for any $p$ as in \eqref{ps0} and for any $\alpha, N \in \N$, for some some $\sigma\gg0$ large enough. 
\end{lemma}
\begin{proof}
Fix any $\bar{s}\geq s_0$ and assume the smallness condition of Prop. \ref{riduzione trasporto}.
By using the property of the cut off function $\chi_0$ in \eqref{cut-offCHIZERO}, by Lemma \ref{lemma autovalori decoupling ordine alto}, by the estimates \eqref{stima cal A 3 new} and the first estimate in \eqref{stimasharpSTANDARD}, one deduces that 
\[
\begin{aligned}
& \big( \Omega + \chi_0 a_3 + \lambda\big)^2 - \chi_0^2b_3^2 = 4 \Omega^2 + d\,, 
%\\& 
\qquad 
\Omega(\xi) + \chi_0a_3(\vphi, x,  \xi) + \lambda = 2 \Omega + e \,,
\end{aligned}
\]
where $d \in S^1, e \in S^{\frac12}$ and they satisfy 
\[
d(\vphi,x, \xi), e(\vphi,x, \xi) \equiv 0 \quad \text{for} \quad |\xi| \leq 1\,,
\]
and, for some $\s\gg1$ large, 
\begin{align}
\| d \|_{1, s, \alpha}^{k_0, \gamma}\,,\, \| e \|_{\frac12, s, \alpha}^{k_0, \gamma} 
&\lesssim_{s, M,\alpha} 
\e \gamma^{- 1} 
\big(1 + \| {\mathcal I}_0 \|_{s + \sigma}^{k_0, \gamma} \big)\,, 
\quad s_0\leq  s \leq \bar{s}\,, 
\label{stimeed1}
\\
\| \Delta_{12} d \|_{1, p, \alpha}\,,\;
\| \Delta_{12} e \|_{\frac12, p, \alpha}&\lesssim_{p,M,\al} \e \g^{-1}  \|i_1-i_2\|_{p+\sigma}\,,
\label{stimeed2}
\end{align} 
for any $\alpha \in \N$ and for $p$ as in \eqref{ps0}.
Since $\Omega (\xi) \gtrsim \langle \xi \rangle^{\frac12} \gtrsim 1$ for $|\xi| \geq 1$, 
one deduces that
$d (4 \Omega)^{- 2} \in S^{0}$, $e (2\Omega)^{- 1} \in S^{0}$ with estimates, using \eqref{stimeed1},
\begin{equation}\label{sist din merda putrefatta 0}
\begin{aligned}
& \| d (4 \Omega)^{- 2}\|_{0, s, \alpha}^{k_0, \gamma}\,,\, 
\|  e (2 \Omega)^{- 1} \|_{0, s, \alpha}^{k_0, \gamma} 
\lesssim_{s,M, \alpha} 
\e \gamma^{- 1}(1 + \| {\mathcal I}_0 \|_{s + \sigma}^{k_0, \gamma})\,, 
\quad s_0\leq s \leq \bar{s}\,, % \quad \forall \alpha \in \N\,. 
\end{aligned}
\end{equation}
for all $\alpha \in \N$.
Therefore, we shall write 
\[
\sqrt{\big( \Omega + \chi_0 a_3 + \lambda\big)^2 - \chi_0^2b_3^2} = 2 \Omega + f 
\]
where $f$ is a symbol satisfying the following: one has $f\in S^{\frac{1}{2}}$, $f (2 \Omega)^{- 1} \in S^0$
and 
\begin{equation}\label{sist din merda putrefatta 1}
\begin{gathered}
 f (\vphi, x, \xi) \equiv 0, \quad \forall |\xi| \leq 1\,, 
\\ 
\| f \|_{\frac12, s, \alpha}^{k_0, \gamma} \,,\;\;  \| f (2 \Omega)^{- 1}\|_{0, s, \alpha}^{k_0, \gamma} 
\lesssim_{s,M, \alpha} \e \gamma^{- 1}(1 + \| {\mathcal I}_0 \|_{s + \sigma}^{k_0, \gamma})\,, 
\end{gathered}
\end{equation}
for any $s_0\leq s\leq \bar{s}$ and 
any $\alpha\in \N$.
The latter properties clearly implies for $\e \gamma^{- 1} \ll 1$ small enough that if $|\xi| \geq 1$ that 
\[
g:=
\sqrt{\big( \Omega + \chi_0 a_3 + \lambda \big)^2 - \chi_0^2b_3^2} \gtrsim \langle \xi \rangle^{\frac12} \gtrsim 1\,.
\]
Hence one can check that $\chi_0 g^{-1}\in S^{-\frac{1}{2}}$ with estimate
\[
\|\chi_0 g^{-1}\|_{-\frac{1}{2},s,\alpha}\lesssim_{s,M,\alpha}
1 + \| {\mathcal I}_0 \|_{s + \sigma}^{k_0, \gamma}\,, \;\;\; s_0\leq  s \leq \bar{s}\,, \quad \forall \alpha \in \N\,. 
\]
%$$
%\begin{aligned}
%%& \frac{\chi_0}{\sqrt{\big( \Omega + \chi_0 a_3 + \lambda \big)^2 - \chi_0^2b_3^2}} \in S^{- \frac12}\,, 
%%\\& 
%\Big\| \frac{\chi_0}{\sqrt{\big( \Omega + \chi_0 a_3 + \lambda \big)^2 - \chi_0^2b_3^2}}  \Big\|_{- \frac12, s, \alpha}^{k_0, \gamma} \lesssim_{s , \alpha} 1 + \| {\mathcal I}_0 \|_{s + \sigma}^{k_0, \gamma}, \quad \forall s \geq s_0, \quad \forall \alpha \in \N\,. 
%\end{aligned}
%$$
The latter estimate, together with the estimate \eqref{stima cal A 3 new} and the first estimate in \eqref{stimasharpSTANDARD} imply the claimed bound \eqref{stime lambda pm fpm} on $f_{-}$. For $f_+$, we write
\[
\begin{aligned}
f_+ & = \frac{\Omega(\xi) + \chi_0a_3(\vphi, x,  \xi) + \lambda}{ \sqrt{\big( \Omega + \chi_0 a_3 + \lambda \big)^2 - \chi_0^2b_3^2}}  = \frac{2 \Omega + e}{2 \Omega + f } = \frac{1 + e(2 \Omega)^{- 1}}{1 + f(2 \Omega)^{- 1}}
\end{aligned}
\]
and hence by the estimates \eqref{sist din merda putrefatta 0} and \eqref{sist din merda putrefatta 1}, 
by the first estimate in  \eqref{stimasharpSTANDARD} and by the mean value theorem, 
one gets the claimed bound on $f_+ - 1$. 
The estimates on $\Delta_{12} f_{\pm}$ follow by similar arguments using 
also \eqref{stimeed2} and reasoning as above. 
The claimed bounds \eqref{stime M M inv ordine alto} on $F$ are a 
direct application of the bounds \eqref{stime lambda pm fpm} on $f_+, f_-$.  
The properties \eqref{C M inv M}-\eqref{stima R infty lambda dec ord max}
 follow by a direct calculation by using \eqref{A3 cut off resto}, 
 \eqref{f pm lambda pm}-\eqref{defmatrixM}, the estimates 
 \eqref{stime M M inv ordine alto} and the first estimate in \eqref{stimasharpSTANDARD}. 
%The estimates \eqref{stime lambda pm fpm}, \eqref{stime M M inv ordine alto}
%follows 
%by using Lemma \ref{lemmacomposizione} and
%%an explicit computation 
%by using 
%\eqref{stima cal A 3 new} and the explicit expressions \eqref{autovalori matrice C 1/2} 
%\eqref{def f g u + u -} and \eqref{defmatrixM}.
The algebraic properties \eqref{prop simmetria lambda +- e fpm} 
follows by using \eqref{cal A3 Hamiltoniano}.
The symbols $f_{\pm}$ satisfy \eqref{invariantsymbol} 
(i.e. are momentum preserving) since  the operator $\mathcal{A}_3$ is 
momentum preserving by Proposition 
\ref{riduzione trasporto}.
\end{proof}

Hence the operator ${\mathcal F} = \opw(F)$ 
is a good candidate the block decouple the order $1/2$.
On the other hand, ${\mathcal F}$ is not symplectic, 
hence we have to construct a symplectic correction that makes the same job. 
So in order to prove the main result of this section (see Proposition \ref{prop operatore cal L3})
we need some preliminary results.

\begin{lemma}\label{lemma M  exp P}
There exists a
map $\Phi_{\mathcal F}:  H^{s}\times H^{s}\to H^{s}\times H^{s} $
which is real-to-real, 
 invertible, symplectic 
 (according to Def. \ref{def simplettica complesse})
 and momentum preserving 
 %which is a correction of the map ${\mathcal M}$ 
 %(up to a one-smoothing remainder) in \eqref{def cal M} 
 satisfying the following properties. 
 Recalling \eqref{defmatrixM}, one has
\begin{equation}\label{espaphiMM}
\begin{aligned}
\Phi_{\mathcal F} & = {\rm Op}^W\Big(F(\vphi, x, \xi) \Big) + {\mathcal R}_{\mathcal F}\,, \quad \Phi_{\mathcal F}^{- 1} & = {\rm Op}^W\Big(F(\vphi, x, \xi)^{- 1} \Big) + {\mathcal Q}_{\mathcal F}
\end{aligned}
\end{equation}
for some momentum preserving operators 
${\mathcal R}_{\mathcal F},  {\mathcal Q}_{\mathcal F} \in OPS^{- 1}$ 
satisfying for some $\s=\s(\alpha)$ the following estimates
\begin{equation}\label{caffe3}
\begin{aligned}
\| {\mathcal R}_{\mathcal F} \|_{- 1, s, \alpha}^{k_0, \gamma} \,,
\,\| {\mathcal Q}_{\mathcal F} \|_{- 1, s, \alpha}^{k_0, \gamma} &\lesssim_{s,M, \alpha} 
\e \g^{-1} (1 + \| {\mathcal I}_0 \|_{s + \sigma}^{k_0, \gamma})\,, \quad \forall s_0 \leq s \leq \bar{s}\,,
\\
 \| \Delta_{12}{\mathcal R}_{\mathcal F} \|_{0, p, \alpha}\,, 
  \| \Delta_{12}{\mathcal Q}_{\mathcal F} \|_{0, p, \alpha}
 &\lesssim_{p,M,\al} \e \g^{-1} \|i_1-i_2\|_{p+\sigma}\,,
\end{aligned}
\end{equation}
for any $p\geq s_{0}$ as in \eqref{ps0}.%, for some $\sigma\gg0$ large enough.
\end{lemma}

\begin{proof}
We look for $\Phi_{\mathcal F}$ of the form 
$\Phi_{\mathcal F} = {\rm exp}(\opw(G(\vphi, x, \xi)))$, 
where the matrix of symbols $G(\vphi, x, \xi) \in S^0$ is of the form 
\begin{equation}\label{caffe6}
G(\vphi, x, \xi) = \begin{pmatrix}
0 & g(\vphi, x, \xi) \\
g(\vphi, x, \xi) & 0
\end{pmatrix}
\end{equation}
for some symbol $g\in S^{0}$ which satisfies 
\begin{equation}\label{prop p enunciato}
g(\vphi, x, \xi) = \overline{g(\vphi, x, \xi)}\,, \qquad  g(\vphi, x, \xi) = g(\vphi, x, - \xi)\,.
\end{equation}
A direct calculation shows that 
\begin{equation}\label{caffe4}
\Phi_{\mathcal F} = \begin{pmatrix}
{\mathcal G}_1 & {\mathcal G}_2 \\
\overline{\mathcal G}_2 & \overline{\mathcal G}_1
\end{pmatrix}\,,
\qquad {\mathcal G}_1 := \sum_{n \geq 0} \dfrac{\opw(g)^{2n}}{(2 n)!}\,, 
\qquad {\mathcal G}_2 := \sum_{n \geq 0} \dfrac{\opw(g)^{2 n + 1} }{(2 n + 1)!}\,.
\end{equation}
Then one has that 
\begin{equation}\label{caffe2}
{\mathcal G}_1 = \opw\Big( g_1(\vphi, x, \xi) \Big) 
+ {\mathcal R}_1\,, \quad {\mathcal G}_2 = \opw\Big( g_2(\vphi, x, \xi)\Big) + {\mathcal R}_2\,,
\end{equation}
with ${\mathcal R}_1, {\mathcal R}_2 \in OPS^{- 1}$ and 
\[
\begin{aligned}
g_1(\vphi, x, \xi) &:= \sum_{n \geq 0} \dfrac{g(\vphi, x, \xi)^{2 n}}{(2 n)!}=\cosh \big(g(\vphi, x, \xi) \big)\,, 
\\
 g_2(\vphi, x, \xi) &:= \sum_{n \geq 0} \dfrac{g(\vphi, x, \xi)^{2 n + 1} }{(2 n + 1)!}
 = \sinh \big(g(\vphi, x, \xi) \big)\,.
\end{aligned}
\]
%This implies that 
%\[
%g_1(\vphi, x, \xi) = \cosh \big(g(\vphi, x, \xi) \big)\,, 
%\quad g_2(\vphi, x, \xi) = \sinh \big(g(\vphi, x, \xi) \big)\,. 
%\]
%Then, we want to choose the symbol $p$ in such a way that 
%\[
%g_1(\vphi, x, \xi) = f_{+}(\vphi, x, \xi) \quad \text{and} \quad g_2(\vphi, x, \xi) = f_{-}(\vphi, x, \xi)
%\]
%namely
Recalling \eqref{defmatrixM}, we want to choose the symbol $g$ in such a way that 
\begin{equation}\label{equazioni per p}
\cosh\big(g(\vphi, x, \xi) \big) = f_{+}(\vphi, x, \xi)\,, 
\quad \sinh \big(g(\vphi, x, \xi) \big) = f_{-}(\vphi, x, \xi)\,. 
\end{equation}
Since $f_{+}^2 - f_{-}^2 = 1$ identically (see \eqref{prop determinante M}), 
the two equations above are compatible and then we can choose 
\[
g(\vphi, x, \xi) = {\rm arcsinh}\big( f_{-}(\vphi, x, \xi)\big)\,.
\]
Clearly, one has that $g(\vphi, x, \xi)$ is a real symbol and symmetric w.r.t. $\xi$ 
(using \eqref{prop simmetria lambda +- e fpm}), 
namely the property \eqref{prop p enunciato} holds. 
It is also momentum preserving by the properties of $f_-$ given 
by Lemma \ref{lem:lambdino}.
So it is the flow $\Phi_{\mathcal{F}}$. Finally, by \eqref{prop p enunciato} and \eqref{prop determinante M}, we have that $\Phi_{\mathcal{F}}$ is symplectic.
Moreover, by explicit computations and recalling \eqref{stime lambda pm fpm},
one gets  the estimates 
\[
\| g \|_{0, s, \alpha}^{k_0, \gamma} \lesssim_{s,M, \alpha} 
\e \g^{-1} (1 + \|{\mathcal I}_0 \|_{s + \sigma}^{k_0, \gamma})\,,\qquad
 \| \Delta_{12} g \|_{0, p, \alpha}
 \lesssim_{p,M,\al} \e \g^{-1}\|i_1-i_2\|_{p+\sigma}\,.
\]
The bounds above, together with Lemma \ref{lemma.potenze.sharp},
%reasoning as  in Lemma $5.2$ in \cite{FIloc} (see also \cite{FI2}) 
and the smallness of $\e$,
imply that  
the remainders ${\mathcal R}_1, {\mathcal R}_2 \in OPS^{- 1}$ 
in \eqref{caffe2} satisfy estimates like \eqref{caffe3}.
In conclusion, formul\ae\, \eqref{caffe4}, \eqref{caffe2}, \eqref{defmatrixM} 
and the equations \eqref{equazioni per p},
imply that
\[
\Phi_{{\mathcal F}}=\opw\big(F(\vphi,x,\x)\big)+\mathcal{R}_{\mathcal{F}}\,,
\qquad \mathcal{R}_{\mathcal{F}}:=
\begin{pmatrix}
{\mathcal R}_1 & {\mathcal R}_2 \\
\overline{\mathcal R}_2 & \overline{\mathcal R}_1 
\end{pmatrix}\,.
\]
By similar arguments one can prove that $\Phi_{{\mathcal F}}^{- 1}$ has the expansion
\[
\Phi_{{\mathcal F}}^{- 1} = {\rm Op}^W\big( F(\vphi, x, \xi)^{- 1} \big) 
+ {\mathcal Q}_{{\mathcal F}}
\]
where ${\mathcal Q}_{{\mathcal F}} \in OPS^{- 1}$ and it satisfies the estimates
\eqref{caffe3}.
The claimed statement has then been proved. 
\end{proof}

We now construct a correction of the map $\Phi_{\mathcal{F}}$ 
%(which is symplectic on $L^{2}$)
which is symplectic on the symplectic  subspace ${\bf H}^{s}_{\perp}$
(see Remark \ref{caffe5}).

\begin{lemma}\label{mappaPhiMsymplectic}
There exists a real-to-real, momentum preserving, invertible 
map $\Phi_{\mathbb F}:  {\bf H}^{s}_{\perp}
\to {\bf H}^{s}_{\perp} $
which is  symplectic on ${\bf H}^{s}_{\perp}$
%(i.e. w.r.t. the symplectic form \eqref{estesaPois})
such that
\[
\Phi_{\mathbb{F}}={\bf \Pi}_{S}^{\perp} \Phi_{\mathcal{F}}{\bf \Pi}_{S}^{\perp}\circ(\uno +\mathcal{R})\,,
\]
where $\mathcal{R}$ is a finite rank operator of the form \eqref{forma buona resto}
satisfying \eqref{giallo2}-\eqref{giallo3}. In particular, one has for some $\s=\s(\alpha)$ the following estimates
\begin{equation}\label{AVB3mappa2PRO}
\begin{aligned}
\mathfrak{M}_{({\Phi}_{\mathbb{F}}^{\pm1}-{\rm Id})}(s)
%+\mathfrak{M}_{({ \Phi}_{\mathbb{F}}^{\pm}-{\rm Id})^{*}}(s)
&\lesssim_{\bar{s},M}
\e \g^{-1} (1+\|\mathcal{I}_{0}\|_{s+\s}^{k_0,\gamma})\,,\;\;\;s_0\leq s\leq \bar{s}\,,
\\
\|\Delta_{12}({ \Phi}_{\mathbb{F}}^{\pm1}-{\rm Id})h\|_{p}
%+\|\Delta_{12}({ \Phi}_{\mathbb{F}}^{\pm}-{\rm Id})^{*}h\|_{p}
&\lesssim_{\bar{s},M}
\e \g^{-1}  \|h\|_{p}\|i_1-i_{2}\|_{p+\s}\,,
\end{aligned}
\end{equation}
for any $p$ as in \eqref{ps0}.
\end{lemma}
\begin{proof}
The map $\Phi_{\mathcal{F}}$ in Lemma \ref{lemma M  exp P}
is constructed as the flow generated by $\opw(G(\vphi,x,\x))$ with symbol in 
\eqref{caffe6}. Therefore, the thesis follows just taking $\Phi_{\mathbb{F}}$ as the flow generated
by the truncation of $\opw(G(\vphi,x,\x))$ (see \eqref{natale})
and 
by applying 
Lemma \ref{differenzaFlussi} (recalling also Remark \ref{differenzaFlussiRMK}).
The bounds \eqref{AVB3mappa2PRO} follow recalling 
\eqref{stime M M inv ordine alto}, \eqref{caffe3}.
The map is momentum preserving since the generator preserves the subspace $S_{\tV}$
by Lemma \ref{lemma M  exp P} (see Lemma \ref{lem:momentoflusso} and Remark \ref{rmk:algflussoEgo}).
\end{proof}

We are now in position to state the main result of the section.
\begin{prop}\label{prop operatore cal L3}
Let  $\beta_{0} \in \N$, $M \in \N$, $\mathtt c \geq 2(k_0 + \beta_0 ) + 3$ with $M \geq \mathtt c$.
%Let $\bar{s} > s_0$, $\tb, N \in \N$. 
For any $\alpha\in \N$ there exists $\s=\s(\beta_0,M, \alpha)\gg 1$ such that for any $\bar{s}>s_0$, if 
condition \eqref{ansatz_I0_s0}
holds with $\mu_0\gg\s$ and if $\e$ is small enough (see \eqref{smallcondepsilon})
then the following holds.
The symplectic map $\Phi_{\mathbb F}$ in Lemma \ref{mappaPhiMsymplectic} 
transforms the Hamiltonian operator ${\mathcal L}_3$ 
in Proposition \ref{riduzione trasporto} into the Hamiltonian operator ${\mathcal L}_4$ given by 
\begin{equation}\label{elle4new}
\begin{aligned}
{\mathcal L}_4 & := \Phi_{\mathbb F}^{- 1} \circ {\mathcal L}_3 \circ \Phi_{\mathbb F} 
\\&=
{\bf \Pi}_{S}^{\perp}\omega\cdot\partial_{\vphi} 
+ {\bf \Pi}_{S}^{\perp}\Big( \opw(\ii \mathtt{m}_{1}\cdot \xi) \uno + \ii E \opw(\lambda) 
+ {\mathcal B}_4 + {\mathcal R}_4\Big){\bf \Pi}_{S}^{\perp}
\end{aligned}
\end{equation}
where $\lambda \in S^{\frac{1}{2}}$  is the symbol in \eqref{autovalori matrice C 1/2} satisfying
the properties given by Lemma \ref{lem:lambdino},
and the operators ${\mathcal B}_4, {\mathcal R}_4$ are defined for all
$(\omega,\th)\in \mathtt \Omega \times[\th_1,\th_2]$ and satisfy the following properties. 
The operator ${\mathcal B}_4  \in OPS^{0}$ satisfies the estimates 
\begin{equation}\label{stima cal B4 new}
\begin{aligned}
\| {\mathcal B}_4 \|_{0, s, \alpha}^{k_0, \gamma} 
&\lesssim_{s,M,\alpha} \e \g^{-1} (1 + \| {\mathcal I}_0 \|_{s + \sigma}^{k_0, \gamma})\,, 
\quad \forall s_0 \leq s \leq \bar{s}\,,
\\
  \| \Delta_{12}{\mathcal B}_{4} \|_{0, p, \alpha}\, 
 &\lesssim_{p,M,\al} \e \g^{-1} \|i_1-i_2\|_{p+\sigma}\,,
\end{aligned}
\end{equation}
for any $p$ as in \eqref{ps0}.
The smoothing remainder ${\mathcal R}_4$ satisfies the following property. 
%For any $M \in \Z$, $M \le N - 2( \tb + k_{0})-3$,  and for 
For any $\beta \in \N^\nu$, $|\beta|\le \beta_0$,
one has that $ \langle D \rangle^{M}\partial_\vphi^\beta {\mathcal R}_4 \langle D \rangle^{- \mathtt{c}}$ 
is a ${\mathcal D}^{k_0}$-tame operator satisfying 
\begin{equation}\label{stima tame cal R4 new}
\begin{aligned}
{\mathfrak M}_{\langle D \rangle^{M}\partial_\vphi^\beta  {\mathcal R}_4 \langle D \rangle^{- \mathtt{c}}}(s) 
&\lesssim_{\bar s,M} 
\e  \g^{-1} (1 + \| {\mathcal I}_0 \|_{s + \sigma})\,, \quad \forall s_0 \leq s \leq \bar{s}\,,
\\
\|  \langle D \rangle^{M} 
\partial_\vphi^\beta \Delta_{12}{\mathcal R}_4\langle D \rangle^{- \mathtt{c}} 
\|_{\cL(H^{p} , H^{p})}& \lesssim_{\bar{s},M}\e  \g^{-1}  \|i_1-i_2\|_{p+\s}\,,
\end{aligned} 
\end{equation}
for any $p$ as in \eqref{ps0}.
Finally the operators $\mathcal{L}_4$, $\mathcal{B}_{4}$ and $\cR_{4}$ are real-to-real, 
Hamiltonian and momentum preserving.
\end{prop}

\begin{proof}
Consider the map $\Phi_{\mathbb{F}}$ in Lemma \ref{mappaPhiMsymplectic}.
First of all, by applying Lemma \ref{georgiaLem} and recalling 
\eqref{op cal L3 new}, we notice that
\[
{\mathcal L}_4 := 
\Phi_{{\mathbb F}}^{- 1} \circ {\mathcal L}_3 \circ \Phi_{\mathbb F}=
\Pi_{S}^{\perp}\mathcal{P}\Pi_{S}^{\perp}+Q_4\,,
\] 
where $Q_4$ is a
finite rank operators satisfying the estimates \eqref{stima tame cal R4 new} and 
where 
\[
\mathcal{P}:=
\Phi_{{\mathcal F}}^{- 1} \circ\Big(\omega\cdot\pa_{\vphi}
+  \opw(\ii \mathtt{m}_{1} \cdot \xi) \uno
+ \ii E \Omega(D)  
+ {\mathcal A}_3 + {\mathcal B}_3 + {\mathcal R}_3\Big) \circ \Phi_{\mathcal F} \,,
\]
where $\Phi_{\mathcal{F}}$ is given by Lemma \ref{lemma M  exp P}.
We also remark that, thanks to the properties of $\Phi_{\mathbb{F}}$, the operator $\mathcal{L}_{4}$
is real-to-real, Hamiltonian and momentum preserving.
The same holds true for the operator $\mathcal{P}$ acting on the whole space $H^{s}\times H^{s}$.
Therefore we only have to study the structure of $\mathcal{P}$ which has the form
\begin{align}
{\mathcal P} & =
\omega\cdot\pa_{\vphi}\nonumber
\\&
+\Phi_{{\mathcal F}}^{- 1} \circ (\omega \cdot \partial_\vphi  \Phi_{\mathcal F} )
\label{cacio1}
\\&+ \Phi_{\mathcal F}^{- 1} \circ \Big( \ii E \Omega(D) +  {\mathcal A}_3 \Big) \circ \Phi_{\mathcal F} 
+
\Phi_{\mathcal F}^{- 1} \circ \Big(   \opw(\ii \mathtt{m}_{1} \cdot \xi) \uno
 + {\mathcal B}_{3} \Big) \circ \Phi_{\mathcal F} 
\label{cacio2}
 \\&
 + \Phi_{\mathcal F}^{- 1} \circ {\mathcal R}_{3} \circ \Phi_{\mathcal F}\,.\label{cacio3}
\end{align}
We analyze separately all the terms in the above expression. 

\smallskip
\noindent
{\bf Analysis of \eqref{cacio1}.}
By Lemma \ref{lemma M  exp P} and by the composition estimates
in Lemma \ref{lemmacomposizione}
one can easily check that $\eqref{cacio1}=:\mathcal{Q}_1$
for some $\mathcal{Q}_1\in OPS^{0}$
satisfying \eqref{stima cal B4 new}. 

%(using that $\omega \in \mathtt \Omega$ which is open and bounded, implying that $\| \omega \cdot \partial_\vphi  \Phi_{\mathcal F} \|_{0, s, \alpha}^{k_0, \gamma} \lesssim \| \Phi_{\mathcal F} \|_{0, s+ 1, \alpha}^{k_0, \gamma}$ for any $s$ and $\alpha$).

\smallskip
\noindent
{\bf Analysis of \eqref{cacio2}.} This is the most crucial term.

\noindent
By applying Lemma \ref{lemma M  exp P}, using the expansion \eqref{espaphiMM},
the composition Lemma \ref{lemmacomposizione}, estimates 
\eqref{stime tame V B inizio lineariz}, \eqref{stima cal A 3 new}, \eqref{stime M M inv ordine alto}, 
\eqref{caffe3} (with possibly a larger $\bar{s}$),
we get that 
\[
\eqref{cacio2}=
\opw(F^{- 1}) \circ (  \opw(\ii \mathtt{m}_{1} \cdot \xi) \uno
+ \ii E \Omega(\xi) 
+ \cA_3
 + {\mathcal B}_{3} ) \circ \opw(F) + {\mathcal Q}_2
\]
where $\mathcal{Q}_2$ is some pseudo-differential operator in $OPS^{0}$
satisfying \eqref{stima cal B4 new}.
Similarly, using symbolic calculus and applying again Lemma \ref{lemmacomposizione} (and recalling \eqref{prop determinante M})
we deduce that 
$$
\begin{aligned}
\opw(F^{- 1}) \circ (  & \opw(\ii \mathtt{m}_{1} \cdot \xi) \uno
+ \ii E \Omega(\xi) 
+ \cA_{3}
 + {\mathcal B}_{3} ) \circ \opw(F)  
\\
& =  \opw(\ii \mathtt{m}_{1} \cdot \xi) \uno
+  {\rm Op}^W \Big( F^{- 1} \big( \ii E \Omega(\xi) + A_3\big) F \Big) + {\mathcal Q}_3
\\&
\stackrel{\eqref{caffe10},\eqref{C M inv M}}{=}
 \opw(\ii \mathtt{m}_{1} \cdot \xi) \uno
+\opw\big(  \ii E \lambda(\vphi, x, \xi)\big) +{\mathcal Q}_3 \,,
\end{aligned}
$$
where ${\mathcal Q}_3 \in OPS^0$ and satisfies the estimates
\eqref{stima cal B4 new}.

\smallskip
\noindent
{\bf Analysis of \eqref{cacio3}.}  
Let us consider the operator
$\Phi_{\mathcal F}^{- 1} 
\circ {\mathcal R}_{3} \circ  \Phi_{\mathcal F}$ where $\Phi_{\mathcal F}^{\pm1}$ are given in Lemma
\ref{lemma M  exp P}.
Moreover we have
\[
\begin{aligned}
\langle D \rangle^{M}\partial_\vphi^\beta  
&\Phi_{\mathcal F}^{- 1} 
\circ {\mathcal R}_{3} \circ  \Phi_{\mathcal F}
 \langle D \rangle^{- \mathtt{c}}
 \\&=
 \sum_{\beta_1+\beta_2+\beta_3=\beta}C_{\beta_1\beta_2\beta_3}
 \langle D \rangle^{M}\partial_\vphi^{\beta_1}  
\Phi_{\mathcal F}^{- 1}  \langle D \rangle^{-M}
\circ  \langle D \rangle^{M}\partial_\vphi^{\beta_2}  {\mathcal R}_{3}
 \langle D \rangle^{-\mathtt{c}}
 \circ \langle D \rangle^{\mathtt{c}} \partial_\vphi^{\beta_3}   \Phi_{\mathcal F}
 \langle D \rangle^{- \mathtt{c}}\,.
 \end{aligned}
\]
Therefore the first bound in \eqref{stima tame cal R4 new}
follows using the composition Lemma \ref{composizione operatori tame AB}
together with Lemma \ref{constantitamesimbolo}
applied on the expansions \eqref{espaphiMM} of $\Phi_{\mathcal{F}}^{\pm1}$.
The estimate on the Lipschitz variation follows similarly.
%By applying 
%Lemma \ref{conj.NEWsmoothresto} (with  $m=0$) and recalling \eqref{stima tame cal R3 new}, we deduce that 
%the operator $\Phi_{\mathcal F}^{- 1} 
%\circ {\mathcal R}_{3} \circ  \Phi_{\mathcal F}$ satisfies the bounds \eqref{stima tame cal R4 new}.

 \smallskip
 \noindent
 In view of Remark \ref{rmk:algsimboli} we have that the pseudo-differential operator 
 $\mathcal{B}_4:= {\mathcal Q}_1 + {\mathcal Q}_2 + {\mathcal Q}_3$
 is momentum preserving thanks to Lemma \ref{lemma M  exp P} 
 and Proposition \ref{riduzione trasporto}. Moreover, since by Lemma \ref{lemma M  exp P}, $\Phi_{\cF}$ is real-to-real and symplectic, then $\cB_{4}$ is real-to-real and Hamiltonian and therefore, by difference, $\cR_{4}$ is so. The discussion above implies the thesis.
%
%The claimed statement then follows by defining ${\mathcal A}_3 := {\mathcal F}_1 + {\mathcal F}_2 + {\mathcal F}_3$ and ${\mathcal R}_3 := \Phi_{\mathcal M}^{- 1} \circ {\mathcal R}_2 \circ  \Phi_{\mathcal M}$. 
\end{proof}

\subsection{Block-Diagonalization at lower orders}\label{subsec:loweroffdiag}
In this subsection we iteratively 
 block-decouple the operator ${\mathcal L}_4$ in \eqref{elle4new}  
 (which is already block-decoupled at the orders $ 1$ and $ 1/ 2 $)  
 up to a smoothing remainder. 
 The main result of the section is the following.

\begin{prop}{\bf (Block-decoupling up to a smoothing remainder).}\label{blockTotale}
%Let $N, \tb \in \N$.
Let  $\beta_{0} \in \N$, $M \in \N$, $\mathtt c \geq 2(k_0 + \beta_0 ) + 3$ with $M \geq \mathtt c$.
Then for any $\alpha\in \N$ there exists $\sigma = \sigma(\beta_0, M, \alpha) \gg 1$ such that for any $\bar{s}\geq s_0$, if 
condition \eqref{ansatz_I0_s0}
holds with $\mu_0\gg\s$ and if $\e$ is small enough (see \eqref{smallcondepsilon})
then the following holds.
There exists a
real-to-real, invertible, symplectic and momentum preserving 
map $\mathcal{S} : {\bf H}^{s}_{{\perp}}\to {\bf H}^{s}_{{\perp}}$, defined for all
$(\omega,\th)\in \mathtt \Omega\times[\th_1,\th_2]$,
such that  
(recall \eqref{elle4new})
\begin{equation}\label{op cal L5}
\begin{aligned}
{\mathcal L}_5 &:= 
{\mathcal S}^{- 1} \circ {\mathcal L}_4 \circ {\mathcal S}
\\&={\bf \Pi}_{S}^{\perp}\omega \cdot \partial_\vphi +
{\bf \Pi}_{S}^{\perp}\Big(
    \opw(\ii \mathtt{m}_{1} \cdot \xi) \uno  +\ii E\opw(\lambda)+ {\mathcal B}_5 + {\mathcal R}_5
\Big){\bf \Pi}_{S}^{\perp}
\end{aligned}
\end{equation}
where
\[
\mathcal{B}_5=\opw\big(B_{5}(\vphi,x,\x)\big)\,,\qquad B_{5}(\vphi,x,\x)=\ii \left(
\begin{matrix}
b_{5}(\vphi,x,\x) & 0 \\ 0 & - b_{5}(\vphi,x,-\x)
\end{matrix}
\right)
\]
for some symbol $b_{5}(\vphi, x, \xi)$ and smoothing operator ${\mathcal R}_5$ 
defined for all $(\omega,\th)\in \mathtt \Omega\times[\th_1,\th_2]$
 satisfing
the following properties. 

\noindent
{\bf (Symbols)}
The symbol $b_{5} \in S^{0}$  satisfies the estimates
\begin{equation}\label{stima d dopo decoupling}
\begin{aligned}
\| b_{5}  \|_{0, s, \alpha}^{k_0, \gamma} &\lesssim_{s, M,\alpha} 
\e \g^{-1} (1 + \| {\mathcal I}_0 \|_{s + \sigma}^{k_0, \gamma})\,, 
\quad \forall s_0 \leq s \leq \bar{s}\,,
\\
  \| \Delta_{12}b_{5} \|_{0, p, \alpha}\, 
 &\lesssim_{p,M,\al} \e \g^{-1} \|i_1-i_2\|_{p+\sigma}\,,
\end{aligned}
\end{equation}
for any $p$ as in \eqref{ps0}.

\noindent
{\bf (Remainder)}
%For any $M \in \N$, $M \le N - 2(\tb +k_{0})-3$,  and for 
For any $\beta \in \N^\nu$, $|\beta|\le \beta_0$,
one has that $\partial_\vphi^\beta \langle D \rangle^{M} {\mathcal R}_5 \langle D \rangle^{- \mathtt{c}}$ 
is a ${\mathcal D}^{k_0}$-tame operator satisfying 
\begin{equation}\label{stima tame cal R5 new}
\begin{aligned}
{\mathfrak M}_{ \langle D \rangle^{M} \partial_\vphi^\beta{\mathcal R}_5 \langle D \rangle^{- \mathtt{c}}}(s) 
&\lesssim_{\bar{s},M} \e \g^{-1}(1 + \| {\mathcal I}_0 \|_{s + \sigma}^{k_0, \gamma})\,, 
\quad \forall s_0 \leq s \leq \bar{s}\,,
\\
\|  \langle D \rangle^{M} 
\partial_\vphi^\beta \Delta_{12}{\mathcal R}_5\langle D \rangle^{- \mathtt{c}} 
\|_{\cL(H^{p} , H^{p})}& \lesssim_{\bar{s},M}\e \g^{-1}\|i_1-i_2\|_{p+\s}\,,
\end{aligned}
\end{equation}
for any $p$ as in \eqref{ps0}.
Moreover, the map 
${\mathcal S} $ satisfies 
\begin{equation}\label{stimamappaS}
\begin{aligned}
\mathfrak{M}_{(\mathcal{S}^{\pm1}-{\rm Id})}(s)&\lesssim_{\bar{s},M}
\e(1+\|\mathcal{I}_{0}\|_{s+\s}^{k_0,\gamma})\,,
\\
\|\Delta_{12}(\mathcal{S}^{\pm1}-{\rm Id})h\|_{p}&\lesssim_{\bar{s},M}
\e\|h\|_{p}\|i_1-i_{2}\|_{p+\s}\,,
\end{aligned}
\end{equation}
for any $p$ as in \eqref{ps0}.
Finally the operators $\mathcal{L}_5$, $\mathcal{B}_{5}$ and $\cR_{5}$ are real-to-real, 
Hamiltonian and momentum preserving.
\end{prop}

\begin{proof}
The proof of this proposition is made by performing an inductive normal form procedure. 
First, we shall describe the iterative step of such iteration. 

\smallskip

\noindent
{\sc INDUCTION STEP.}

\noindent
At the $k$-th step, $k\geq0$, of such a procedure we deal with a Hamiltonian operator of the form 
\begin{equation}\label{cal L 4 n block decoupling}
\begin{aligned}
{\mathcal L}_4^{(k)} & = 
{\bf \Pi}_{S}^{\perp}\omega \cdot \partial_\vphi 
\\&+
{\bf \Pi}_{S}^{\perp}\Big(
 \opw(\ii \mathtt{m}_{1} \cdot \xi) \uno 
+\ii E\opw(\lambda)+ {\mathcal B}^{(k)}_4+\mathcal{C}_{4}^{(k)}+ {\mathcal R}_4^{(k)}
\Big){\bf \Pi}_{S}^{\perp}\,,
\end{aligned}
\end{equation}
where $\cB_{4}^{(k)}, {\mathcal C}_4^{(k)}, {\mathcal R}_4^{(k)}$ satisfy the following properties. 
Given $\beta_0, M ,\mathtt{c}\in \N$, $M\geq\mathtt{c}\geq 2(\beta_0+k_0)+3$, 
for any $\alpha\in \N$ there exists a non-decreasing sequence 
$\sigma_k \equiv \sigma_k(\beta_0, M, \alpha)$ 
such that 
such that if 
condition \eqref{ansatz_I0_s0}
holds with $\mu_0=\s_{k}$ and if $\e$ is small enough (see \eqref{smallcondepsilon}) then
the following hold. 
\begin{itemize}
\item{\bf $({\mathcal P}1)_k$} 
The operator 
${\mathcal B}_4^{(k)} = \opw(B_4^{(k)}(\vphi, x, \xi)) \in OPS^{0}$ is a real-to-real, 
Hamiltonian operator of the form 
\[
\mathcal{B}^{(k)}_4=\opw\big(B^{(k)}_{4}(\vphi,x,\x)\big)\,,\qquad B^{(k)}_{4}(\vphi,x,\x)=\ii \left(
\begin{matrix}
b_{4}^{(k)}(\vphi,x,\x) & 0 \\ 0 & -b_{4}^{(k)}(\vphi,x,-\x)
\end{matrix}
\right)
\]
where
$b_{4}^{(k)} \in S^{0}$ 
satisfies
\[
\begin{aligned}
\|  b_{4}^{(k)} \|_{0, s, \alpha}^{k_0, \gamma} 
&\lesssim_{s, \alpha, k,M} 
\e \g^{-1} (1 + \| {\mathcal I}_0 \|_{s + \sigma_k}^{k_0, \gamma}), \quad \forall s_0 \leq s \leq \bar{s}\,,
\\
  \| \Delta_{12}b_{4}^{(k)}\|_{0, p, \alpha}\,, 
 &\lesssim_{p,\al, k,M} \e\g^{-1} \|i_1-i_2\|_{p+\sigma_k}\,,
\end{aligned}
\]
for any $p$ as in \eqref{ps0}.
\item{\bf $({\mathcal P}2)_k$} The operator 
${\mathcal C}_4^{(k)} = \opw(C_4^{(k)}(\vphi, x, \xi)) \in OPS^{- \frac{k}{2}}$ is a real-to-real, 
Hamiltonian operator of the form 
\begin{equation}\label{proCCk}
\begin{gathered}
 {C}_4^{(k)}(\vphi, x, \xi) = \ii \begin{pmatrix}
0& c_4^{(k)}(\vphi, x, \xi) \\[0.3em]
- \overline{c_4^{(k)}(\vphi, x, \xi)} & 0
\end{pmatrix}\,, 
%\\& 
\quad c_4^{(k)}(\vphi, x, \xi) = c_4^{(k)}(\vphi, x, - \xi)\,, 
\\
 \| { C}_4^{(k)} \|_{- \frac{k}{2}, s, \alpha}^{k_0, \gamma} 
\lesssim_{s, \alpha, k,M} \e \g^{-1} (1 + \| {\mathcal I}_0 \|_{s + \sigma_k}^{k_0, \gamma})\,, 
\quad \forall s_0 \leq s \leq \bar{s}\,,
\\
  \| \Delta_{12}{ C}_4^{(k)}\|_{0, p, \alpha}\,
 \lesssim_{p,\al,k,M} \e \g^{-1} \|i_1-i_2\|_{p+\sigma_k}\,,
\end{gathered}
\end{equation}
for any $p$ as in \eqref{ps0}.
\item{$({\mathcal P}3)_k$.} % For any $M \in \N$, $M \le N - 2(\tb + k_{0})-3$,  
For any $\beta \in \N^\nu$, $|\beta|\le \beta_0$,
 the operator  $  \langle D \rangle^{M}\partial_\vphi^\beta{\mathcal R}_4^{(k)} \langle D \rangle^{-\mathtt{c}}$ 
 is ${\mathcal D}^{k_0}$- tame operator and it satisfies the estimates
\begin{equation}\label{punkina1}
\begin{aligned}
{\mathfrak M}_{ \langle D \rangle^{M}\partial_\vphi^\beta{\mathcal R}_4^{(k)} 
\langle D \rangle^{- \mathtt{c}}}(s) 
&\lesssim_{\bar{s}, k,M} \e \g^{-1} (1 + \| {\mathcal I}_0 \|_{s + \sigma_k}^{k_0, \gamma})\,, 
\quad \forall s_0 \leq s \leq \bar{s}\,,
\\
\|  \langle D \rangle^{M} 
\partial_\vphi^\beta \Delta_{12}{\mathcal R}_4^{(k)}\langle D \rangle^{- \mathtt{c}} 
\|_{\cL(H^{p} , H^{p})}& \lesssim_{\bar{s},k,M}\e \g^{-1}\|i_1-i_2\|_{p+\s}\,,
\end{aligned}
\end{equation}
for any $p$ as in \eqref{ps0}.
\item{$({\mathcal P}4)_k$.} For $k\geq1$
there exists a momentum preserving symbol $p_{k-1}\in S^{-\frac{k}{2}}$
such that
\begin{equation}\label{coniugiostepk}
\mathcal{L}_{4}^{(k)}
={\bf \Psi}_{\mathcal{T}_{k-1}}^{-1}\circ \mathcal{L}_{4}^{(k-1)}\circ{\bf \Psi}^{1}_{\mathcal{T}_{k-1}}
\end{equation}
where ${\bf \Psi}_{\mathcal{T}_{k-1}}^{\theta}$, $\theta\in[0,1]$ 
is the flow at time 
$\theta\in [0,1]$ of 
\begin{equation}\label{flow-jth}
\partial_{\theta} {\bf \Psi}^{\theta}_{\mathcal{T}_{k-1}} (\vphi) = 
 {\bf \Pi}_{S}^{\perp}\mathcal{P}_{k-1}
{\bf \Pi}_{S}^{\perp}{\bf \Psi}_{\mathcal{T}_{k-1}}^{\theta}(\vphi) \, ,
 \quad {\bf \Psi}_{\mathcal{T}_{k-1}}^{0}(\vphi) = {\rm Id} \, ,  
\end{equation}
with
\begin{equation}\label{generatore-jth}
 {\mathcal P}_{k-1} := \opw\begin{pmatrix}
0 &  p_{k-1}(\vphi, x, \xi) \\[0.3em]
- \overline{ p_{k-1}(\vphi, x, \xi) } & 0
\end{pmatrix}\,,\quad \quad p_{k-1}(\vphi, x, \xi) = p_{k-1}(\vphi, x, - \xi)\,.
\end{equation}

\item{$({\mathcal P}5)_k$.} The operator $\mathcal{L}_{4}^{(k)}$ is real-to-real, Hamiltonian and momentum 
preserving as well as the operators
$\mathcal{B}^{(k)}_{4}$, ${\mathcal C}_4^{(k)}$ and ${\mathcal R}_{4}^{(k)}$.
\end{itemize}

\smallskip
At $k=0$ we have that the operator $\mathcal{L}_{4}$ in \eqref{elle4new}
has the form $\mathcal{L}_{4}^{(0)}$ in \eqref{cal L 4 n block decoupling}
with ${\mathcal B}_4^{(0)}:={\mathcal B}_4$,
${\mathcal C}_4^{(0)}\equiv0$ and ${\mathcal R}_4^{(k)}:={\mathcal R}_4$.
The properties $({\mathcal P}1)_0$, $({\mathcal P}2)_0$, $({\mathcal P}3)_0$, $({\mathcal P}5)_0$
follow by Proposition \ref{prop operatore cal L3} while $({\mathcal P}4)_0$ in empty.

Let us now assume  $({\mathcal P}j)_k$ with $k\geq1$ and $j=1,\ldots,5$.
In order to prove inductively the assertions above
we want to choose the symbol $p_k$ in order to cancel out 
the off-diagonal term $c_{4}^{(k)}$ 
up to a remainder of order $- (k + 1)/2$. 
We want to choose $p_{k}$ in such a way that 
\[
- 2 \lambda(\vphi, x, \xi) p_k(\vphi, x, \xi) + \ii c_{4}^{(k)}(\vphi, x, \xi) = S^{- \infty} \quad 
\text{(smoothing symbol)}\,. 
\]
Note that by Lemma \ref{lemma autovalori decoupling ordine alto}, one deduces that 
$$
\begin{gathered}
 \lambda(\vphi, x, \xi) = \Omega(\xi) + r_\lambda(\vphi, x, \xi), \quad r_\lambda \in S^{\frac12}\,,  \\
 |r_\lambda(\vphi, x, \xi)| \lesssim \e \langle \xi \rangle^{\frac12}, \quad \forall (\vphi, x, \xi) \in \T^\nu \times \T^d_\Gamma \times \R^d\,. 
\end{gathered}
$$
Moreover, by \eqref{cutofffunct}, \eqref{Omega D cut off}, if $|\xi | \geq 1$, one has that 
$$
|\lambda(\vphi, x, \xi)| \gtrsim \langle \xi \rangle^{\frac12}  - C \e \langle \xi \rangle^{\frac12} \gtrsim \langle \xi \rangle^{\frac12}
$$
for $\e \ll 1$ small enough. 
Consider the even  cut-off function $\chi_0 \in C^{\infty}(\R^{d};\R)$ defined in \eqref{cut-offCHIZERO}.
% such that 
%\begin{equation}
%\begin{aligned}
%& 0 \leq \chi_0\leq 1, \quad \chi_0(\xi) = \chi_0(- \xi)\,, \\
%&  \chi_0(\xi) = 0 \quad \text{if} \quad |\xi| \leq 1\,, \qquad 
% \chi_0(\xi) = 1 \quad \text{if} \quad |\xi| \geq 2\,. 
%\end{aligned}
%\end{equation}
% and on the support of the cut-off function one has $\chi_0$, $|\lambda(\vphi, x, \xi) | \gtrsim \langle \xi \rangle^{\frac12} \gtrsim 1$. 
We define 
\begin{equation}\label{def pk decoupling}
p_k(\vphi, x, \xi) :=  \dfrac{ \ii \chi_0(\xi)c_{4}^{(k)}(\vphi, x, \xi)}{2 \lambda(\vphi, x, \xi)}\,.
\end{equation}
By using the properties of $\lambda$ (recall the cut-off functions $\chi$ in \eqref{cutofffunct} \eqref{Omega D cut off} and $\chi_{0}$ in \eqref{cut-offCHIZERO}) one can check that $\frac{\chi_0}{\lambda} \in S^{- \frac12}$. Moreover, by the estimates of Lemma \eqref{lemma autovalori decoupling ordine alto} on 
$\lambda$ and $({\mathcal P}2)_k$,  one has  that $p_k \in S^{- \frac{k + 1}{2}}$ and that it satisfies the estimates
\begin{equation}\label{stima pk decoupling}
\begin{aligned}
%\| {\mathcal P}_k \|_{- \frac{k + 1}{2}, s, \alpha}^{k_0, \gamma} 
%\simeq 
\| p_k \|_{- \frac{k + 1}{2}, s, \alpha}^{k_0, \gamma} 
&\lesssim_{s,k, \alpha,M} 
\e \g^{-1}(1 + \| {\mathcal I}_0 \|_{s + \sigma_k}^{k_0, \gamma})\,, 
\quad \forall s_0 \leq s \leq \bar{s}\,,
\\
  \| \Delta_{12} p_{k}\|_{0, p, \alpha}\, 
 &\lesssim_{p,k,\al,M} 
 \e \g^{-1} \|i_1-i_2\|_{p+\sigma_k}\,.
\end{aligned}
\end{equation}
Moreover, in view of \eqref{def pk decoupling}, one has that $p_{k}$ is the solution of
\begin{equation}\label{pizza decoupling 0}
- 2 \lambda(\vphi, x, \xi) p_k(\vphi, x, \xi) + \ii \chi_0(\xi)c_{4}^{(k)}(\vphi, x, \xi) = 
\ii (1 - \chi_0(\xi)) c_{4}^{(k)}(\vphi, x, \xi) \in S^{- \frac{ k + 1}{2}}\,.
\end{equation}
By the estimates \eqref{proCCk} in $({\mathcal P}2)_k$ and
observing that $1 - \chi_0(D) \in OPS^{- \infty} $, then 
one gets that for any $s_0 \leq s \leq \bar{s}$, 
\begin{equation}\label{pizza decoupling 1}
\begin{aligned}
\| (1 - \chi_{0}) c_{4}^{(k)} \|_{- \frac{k + 1}{2}, s, \alpha}^{k_0, \gamma} 
&\lesssim_{s, k,\alpha} 
\| c_{4}^{(k)} \|_{- \frac{k}{2}, s, \alpha}^{k_0, \gamma} 
\lesssim_{s,k, \alpha,M} 
\e \g^{-1} (1 + \| {\mathcal I}_0 \|_{s + \sigma_k}^{k_0, \gamma})\,,
\\
  \| \Delta_{12}(1 - \chi_{0}) c_{4}^{(k)}\|_{-\frac{k+1}{2}, p, \alpha}\, 
 &\lesssim_{p,k,\al,M} \e \g^{-1}\|i_1-i_2\|_{p+\sigma_k}\,.
 \end{aligned}
\end{equation}
By $({\mathcal P}5)_k$ (recall Def. \ref{def:mompressimbo})
we deduce that also the symbol $p_k$ is momentum preserving.
By Lemma \ref{lemma:buonaposFlussi}
we have that the flow map ${\bf \Psi}_{\mathcal{T}_{k}}$ in \eqref{flow-jth} is well-posed
and satisfies estimates like \eqref{stimamappaS}
by using \eqref{stima pk decoupling}.
 In particular, it is 
real-to-real 
Hamiltonian and momentum preserving (see Lemma \ref{lem:momentoflusso} and Remark \ref{rmk:algflussoEgo}).
Setting $\mathcal{T}_{k}:=\exp(\ii \mathcal{P}_{k})$,
which is the solution of \eqref{flussosimboloGenerico} with generator $\ii\mathcal{P}_{k}$,
 and, applying 
Lemmata \ref{differenzaFlussi} and \ref{georgiaLem} (recall Remark. \ref{differenzaFlussiRMK}),
we deduce that 
 the conjugate operator $\mathcal{L}_{3}^{(k+1)}$ in 
\eqref{coniugiostepk} has the form
\begin{equation}\label{nokia1}
\mathcal{L}_{4}^{(k+1)}={\bf \Pi}_{S}^{\perp}{\mathcal T}_k^{- 1} \circ {\mathcal L}_4^{(k)} 
\circ {\mathcal T}_k{\bf \Pi}_{S}^{\perp}+\mathcal{Q}\,,
\end{equation}
for some finite rank operator $\mathcal{Q}$ satisfying bounds like \eqref{punkina1}. 

We then compute 
${\mathcal T}_k^{- 1} \circ {\mathcal L}_4^{(k)} \circ {\mathcal T}_k$ 
and we analyze separately all the terms arising in the conjugation. 
By applying Lemma \ref{lemma:buonaposFlussi} 
and the estimates \eqref{stime tame V B inizio lineariz}, 
\eqref{stima pk decoupling}, \eqref{stime lambda pm fpm},
$({\mathcal P}1)_k$-$({\mathcal P}3)_k$, and recalling
\eqref{cal L 4 n block decoupling}
one obtains the following expansions: for some constant $\mu \gg 0$
\begin{align}
{\mathcal T}_k^{- 1} \circ \omega \cdot \partial_\vphi \circ {\mathcal T}_k 
& 
= \omega \cdot \partial_\vphi + {\mathcal F}_1\,, 
\label{Tk om d vphi decoupling}
\\
{\mathcal T}_k^{- 1} \circ    \opw(\ii \mathtt{m}_{1} \cdot \xi) 
\uno  \circ {\mathcal T}_k 
& =     \opw(\ii \mathtt{m}_{1} \cdot \xi) \uno 
\nonumber
\\&
+  \ii [   \opw(\ii \mathtt{m}_{1} \cdot \xi) \uno \,,\, {\mathcal P}_k] 
+ {\mathcal F}_2\,,
\label{Tk V dot nabla decoupling}
\\
{\mathcal T}_k^{- 1} \circ\big( \ii E\opw(\lambda)
+ {\mathcal B}^{(k)}_4\big)\circ {\mathcal T}_k 
&=
\ii E\opw(\lambda)+ {\mathcal B}^{(k)}_4\nonumber
\\&+\big[\mathcal{P}_k,\ii E\opw(\lambda)+ {\mathcal B}^{(k)}_4\big]
+\mathcal{F}_3\,,\label{Tk lambda dk decoupling}
\\
{\mathcal T}_k^{- 1} \circ {\mathcal C}_4^{(k)} \circ {\mathcal T}_k 
& = {\mathcal C}_{4}^{(k)} + {\mathcal F}_4\,, 
\qquad\quad \qquad {\mathcal F}_4 \in OPS^{- k - \frac32}\,,
\end{align}
for some pseudo-differential operators
\[
{\mathcal F}_1 \in OPS^{- \frac{k + 1}{2}}\,,
\qquad 
{\mathcal F}_2, {\mathcal F}_3 \in OPS^{- k - 2}\,,
\qquad
{\mathcal F}_4 \in OPS^{- k - \frac{3}{2}}\,,
\]
satisfying the estimates
\[
\| {\mathcal F}_1 \|_{- \frac{k + 1}{2}, s, \alpha}^{k_0, \gamma} \,,
\| {\mathcal F}_j \|_{-k-2 , s, \alpha}^{k_0, \gamma} \,,
\| {\mathcal F}_4 \|_{- k-\frac{3}{2}, s, \alpha}^{k_0, \gamma} 
 \lesssim_{s,k, \alpha,M} 
 \e \g^{-1}(1 + \|{\mathcal I}_0 \|_{s + \sigma_k + \mu})\,, 
\]
for $s_0 \leq s \leq \bar{s}$, $j=2,3$, and
 \[
  \| \Delta_{12}{\mathcal F}_1\|_{- \frac{k + 1}{2}, p, \alpha}\,, 
    \| \Delta_{12}{\mathcal F}_j\|_{- k-2, p, \alpha}\,, 
      \| \Delta_{12}{\mathcal F}_1\|_{- k-\frac{3}{2}, p, \alpha}\, 
 \lesssim_{p,k,\al,M} 
 \e \g^{-1} \|i_1-i_2\|_{p+\sigma_k+\mu}\,.
\]
for $j=2,3$. Moreover, by defining 
$\widetilde{\mathcal R}_4^{(k + 1)} := {\mathcal T}_k^{- 1} \circ {\mathcal R}_4^{(k)} \circ {\mathcal T}_k$, 
one has that,
using Lemma \ref{conj.NEWsmoothresto} (recall also item $(ii)$ in Remark \ref{rmk:restimatrici}),
 for any $\beta \in \N^\nu$, $|\beta| \leq \beta_0$, $M\in\N, M\geq \mathtt{c}\geq   2(\tb+ k_{0})+3$, the operator 
 $ \langle D \rangle^{M} \partial_\vphi^\beta{\mathcal R}_4^{(k+ 1)} \langle D \rangle^{- \mathtt{c}}$ 
 is ${\mathcal D}^{k_0}$-tame with tame constant satisfying, for $s_0 \leq s \leq \bar{s}$,
\begin{equation}\label{stima cal R4 k+ 1 decoupling}
\begin{aligned}
{\mathfrak M}_{ \langle D \rangle^{M} \partial_\vphi^\beta\widetilde{\mathcal R}_4^{(k + 1)} 
\langle D \rangle^{- \mathtt{c}}}(s)  
&\lesssim_{\bar{s}, k, M} 
\e \g^{-1} (1 + \| {\mathcal I}_0 \|_{s + \sigma_k + \mu}^{k_0, \gamma})\,, 
\\
\|  \langle D \rangle^{M} 
\partial_\vphi^\beta \Delta_{12}\widetilde{\mathcal R}_{4}^{(k + 1)} \langle D \rangle^{- \mathtt{c}} 
\|_{\cL(H^{p} , H^{p})}
& 
\lesssim_{p,k,M}
\e \g^{-1}\|i_1-i_2\|_{p+\s_{k}+\mu}\,.
\end{aligned}
\end{equation}
It remains only to analyze the two commutators in \eqref{Tk V dot nabla decoupling}, \eqref{Tk lambda dk decoupling}. By applying Lemma \ref{commutatorsimbol} and the estimates \eqref{stime tame V B inizio lineariz}, \eqref{stima pk decoupling}, recalling \eqref{generatore-jth} one gets that 
the operator ${\mathcal F}_5  := \ii [   \opw(\ii \mathtt{m}_{1} \cdot \xi) \uno  \,,\,{\mathcal P}_k ]$
satisfies 
\begin{equation}\label{commutatore V nabla Pk decoupling}
\begin{aligned}
\| {\mathcal F}_5 \|_{- \frac{k + 1}{2}, s, \alpha}^{k_0, \gamma} 
& \lesssim_{ s,k, \alpha,M} 
\e \g^{-1} (1 + \| {\mathcal I}_0 \|_{s + \sigma_k + \mu}^{k_0, \gamma})\,, 
\quad \forall s_0 \leq s \leq \bar{s}\,,
\\
    \| \Delta_{12}{\mathcal F}_5\|_{-\frac{k+1}{2}, p, \alpha}\,
 &\lesssim_{p,k,\al,M} 
 \e \g^{-1}\|i_1-i_2\|_{p+\sigma_k+\mu}\,.
\end{aligned}
\end{equation}
Moreover
\[
\begin{aligned}
\big[\mathcal{P}_k,\ii E\opw(\lambda)&+ {\mathcal B}^{(k)}_4\big]
 = \begin{pmatrix}
0 & {\mathcal B}_k \\
\overline{\mathcal B}_k & 0
\end{pmatrix}\,, \\
& {\mathcal B}_k := - \Big(  \opw(p_k) \circ \overline{{\rm Op}^W(\lambda + b_{4}^{(k)})} 
+ {\rm Op}^W(\lambda + b_{4}^{(k)}) \circ \opw(p_k) \Big)\,. 
\end{aligned}
\]
Then by applying Lemma  \ref{commutatorsimbol} and the estimates  
\eqref{stima pk decoupling}, \eqref{stime lambda pm fpm} 
(using the property that $\lambda$ is real and symmetric w.r.t. $\xi$) 
one gets that 
\begin{equation}\label{espansione cal Ck decoupling}
\begin{aligned}
\big[\mathcal{P}_k,\ii E\opw(\lambda)+ {\mathcal B}^{(k)}_4\big]
& =\opw \begin{pmatrix}
0 & - 2 \lambda p_k \\
- \overline{2 \lambda p_k} & 0
\end{pmatrix} + {\mathcal F}_6, \quad {\mathcal F}_6 \in OPS^{- \frac{k + 1}{2}}\,,
\end{aligned}
\end{equation}
for some operator $\mathcal{F}_6$ satisfying 
\[
\begin{aligned}
 \| {\mathcal F}_6 \|_{- \frac{k + 1}{2}, s, \alpha} 
 &\lesssim_{ s,k, \alpha,M} 
 \e \g^{-1} (1 + \| {\mathcal I}_0 \|_{s + \sigma_k + \mu}^{k_0, \gamma})\,, 
 \quad \forall s_0 \leq s \leq \bar{s}\,,
 \\
     \| \Delta_{12}{\mathcal F}_6\|_{-\frac{k+1}{2}, p, \alpha}\, 
 &\lesssim_{p,k,\al,M} 
 \e \g^{-1} \|i_1-i_2\|_{p+\sigma_k+\mu}\,.
 \end{aligned}
\]
Finally,  by \eqref{pizza decoupling 0}, \eqref{pizza decoupling 1},
we deduce
\begin{equation}\label{cal F7 decoupling}
\begin{aligned}
{\mathcal F}_7 & :=
\opw \begin{pmatrix}
0 & - 2 \lambda p_k \\
- \overline{2 \lambda p_k} & 0
\end{pmatrix}
+ {\mathcal C}_3^{(k)} \in OPS^{- \frac{k + 1}{2}}\,, 
\end{aligned}
\end{equation}
where $F_{7}\in OPS^{-\frac{k+1}{2}}$ satisfies 
\[
\begin{aligned}
 \| {\mathcal F}_7 \|_{- \frac{k + 1}{2}, s, \alpha}^{k_0, \gamma} 
 &\lesssim_{ s,k, \alpha,M} 
 \e \g^{-1} (1 + \| {\mathcal I}_{0} 
 \|_{s + \sigma_k + \mu}^{k_0, \gamma})\,, 
 \quad \forall s_0 \leq s \leq \bar{s}\,,
 \\
     \| \Delta_{12}{\mathcal F}_7\|_{-\frac{k+1}{2}, p, \alpha}\,
 &\lesssim_{p,k,\al,M} 
 \e \g^{-1} \|i_1-i_2\|_{p+\sigma_k+\mu}\,.
 \end{aligned}
\]
We now define the $2\times 2$ matrix of  pseudo-differential 
operators ${\mathcal F} := \sum_{i = 1}^7 {\mathcal F}_i  \in OPS^{- \frac{k + 1}{2}}$ 
and by summarizing all the previous expansions and estimates, 
one obtains that 
\[
\begin{aligned}
 \| {\mathcal F} \|_{- \frac{k + 1}{2}, s, \alpha}^{k_0, \gamma} 
 &\lesssim_{ s,k, \alpha,M} 
 \e \g^{-1} (1 + \| {\mathcal I}_{0} 
 \|_{s + \sigma_k + \mu}^{k_0, \gamma})\,, 
 \quad \forall s_0 \leq s \leq \bar{s}\,,
 \\
     \| \Delta_{12}{\mathcal F}\|_{-\frac{k+1}{2}, p, \alpha}\, 
 &\lesssim_{p,k,\al,M} \e \g^{-1} \|i_1-i_2\|_{p+\sigma_k+\mu}\,,
 \end{aligned}
\]
and 
\[
{\mathcal T}_k^{- 1} \circ {\mathcal L}_4^{(k)} \circ {\mathcal T}_k  
=\omega\cdot\pa_{\vphi}+
  \opw(\ii \mathtt{m}_{1} \cdot \xi) \uno 
  +\ii E\opw(\lambda)+ {\mathcal B}^{(k)}_4+\mathcal{F}+ \widetilde{{\mathcal R}}_4^{(k+1)}\,.
\]
The claimed statements 
$({\mathcal P}1)_{k + 1}$-$({\mathcal P}4)_{k + 1}$
 then follow recalling \eqref{nokia1} and 
by taking $\sigma_{k + 1} \geq \sigma_k + \mu$. 
The last item $({\mathcal P}5)_{k + 1}$
follows by the inductive assumption and since the map 
${\bf \Psi}_{\mathcal{T}_{k}}$
generated by choosing $p_k$ as in \eqref{def pk decoupling}
is real-to-real, symplectic (on the restricted subspace ${\bf H}_{S}^{\perp}$)
and momentum preserving.
\end{proof}

\begin{proof}[{\bf Proof of Proposition \ref{blockTotale} concluded}.]
We assume \eqref{smallcondepsilon} with $\kappa\gg 2(M-\mathtt{c})$.
%let $k \in \{ 0, \ldots, 2(M-\mathtt{c}) \}$.  
We define
\begin{equation*}
{\mathcal S} := {\bf \Psi}^{1}_{\mathcal{T}_{0}}\circ\ldots\circ{\bf \Psi}^{1}_{\mathcal{T}_{2(M-\mathtt{c})}}\,.
\end{equation*}
By composition (using \eqref{stima pk decoupling}) and the estimates on the map 
${\bf \Psi}^{1}_{\mathcal{T}_{i}}$,
one can prove that the map $\mathcal{S}$ satisfies 
the bounds \eqref{stimamappaS}
for $\sigma \gg \sigma_k$, for any $k \in \{0, \ldots, 2 (M-\mathtt{c})  \}$ large enough. 
By iterating the induction step described above, we set
${\mathcal L}_5 \equiv {\mathcal L}_4^{(2 (M-\mathtt{c}))} 
= {\mathcal S}^{- 1} \circ {\mathcal L}_4 \circ {\mathcal S}$.
In view of \eqref{cal L 4 n block decoupling} we define $\mathcal{B}_{5}\equiv\mathcal{B}_{4}^{(2(M-\mathtt{c}))} $.
On the other hand, by \eqref{proCCk}
and Lemma \ref{constantitamesimbolo} we have that $\cC_{4}^{(2(M-\mathtt{c}))}$ can be absorbed into the remainder
satisfying \eqref{stima tame cal R5 new}, as well as $\mathcal{R}_{4}^{(2(M-\mathtt{c}))}$, 
see \eqref{punkina1}.
The discussion above implies that ${\mathcal L}_5$ 
satisfies the properties \eqref{op cal L5}-\eqref{stima tame cal R5 new}. 
The proof of the proposition is then concluded. 
\end{proof}

\section{Reduction  to constant coefficients up to smoothing remainders}\label{sec:redulower}

%\subsection{Reduction of lower orders}\label{subsec.redulower}
In this section we iteratively reduce to constant coefficients the linear operator 
${\mathcal L}_5$ in \eqref{op cal L5} up to a smoothing remainder. 
We split this procedure in two parts. The first part concerns 
the reduction to constant coefficient of the order $1/2$ whereas 
the second part is about the reduction of the lower order terms $0, - 1/2, \ldots$. 

\subsection{Reduction at order 1/2}\label{sec:riduzione diagonal ordine 1 2}
Here we normalize the symbol $\lambda$, 
by removing the dependence from $(\vphi, x)$ up to an arbitrarily regularizing remainder. 
The main result of this subsection is the following.

\begin{prop}{\bf (Reduction to constant coefficients of the order $1/2$.)}\label{riduzione diagonal ordine 1 2}
%Let $N, \beta_{0} \in \N$. Then for any $\alpha\in\N$ there exists $\s=\s(N, \tb, \alpha)$  such that if 
Let  $\beta_{0} \in \N$, $M \in \N$, $\mathtt c \geq 3(k_0 + \beta_0 ) + 3$ with $M \geq \mathtt c$.
Then for any $\alpha\in \N$ there exists $\sigma = \sigma(\beta_0, M, \alpha) \gg 1$ such that, for any $\bar{s}\geq s_0$, if 
condition \eqref{ansatz_I0_s0}
holds with $\mu_0\geq \s$ and if $\e$ is small enough (see \eqref{smallcondepsilon})
then the following holds.
There exist a symbol $\mathtt{m}_{6}\in S^{\frac{1}{2}}$ independent of $(\vphi,x)$, 
an operator  $\mathcal{B}_6 \in OPS^0$, a 
$\mathcal{D}^{k_0}$-tame operator ${\mathcal R}_6$
and 
a real-to-real, invertible symplectic and momentum preserving map 
${\bf \Psi}^{\pm 1} : {\bf H}^s_{\perp}\to {\bf H}^s_{\perp}$, $s_0 \leq s \leq \bar{s}$ ,
defined for all $(\omega,\mathtt{h})\in \mathtt \Omega\times[\th_1,\th_2]$ such that the following holds.
For any $(\omega,\mathtt{h})\in 	\tT\tC_{\infty}(\gamma,\tau) $ (see \eqref{tDtCn})
one has the conjugation 
(recall \eqref{op cal L5})
\begin{equation}\label{op cal L6}
\begin{aligned}
{\mathcal L}_6 & = {\bf \Psi}^{- 1} \circ {\mathcal L}_5 \circ {\bf \Psi}\\
&= {\bf \Pi}_{S}^{\perp}\omega \cdot \partial_\vphi + 
{\bf \Pi}_{S}^{\perp}\Big( 
\ii\opw(\mathtt{m}_{1}  \cdot \xi) \uno  
+ \ii E\opw(\mathtt{m}_{6})+\cB_{6} + \cR_{6} \Big) {\bf \Pi}_{S}^{\perp}\,,
\end{aligned}
\end{equation}

where $\mathtt{m}_{1}$ is given in Proposition \ref{riduzione trasporto} .
Moreover, for any $(\omega,\mathtt{h})\in \mathtt \Omega\times[\th_1,\th_2]$
one has the following:

\noindent
$(i)$ the Fourier multiplier  $\mathtt{m}_{6}\in S^{\frac{1}{2}}$ is a real valued, momentum preserving symbol and satisfies 
\[
\mathtt{m}_{6}( \xi)=\mathtt{m}_{6}(-\xi);
\]
and  (recall \eqref{wild})
\begin{equation}\label{stima.m.6}
\begin{aligned}
\| \mathtt{m}_{6} - \Omega  \|_{\frac12, s, \alpha}^{k_0, \gamma}
&\lesssim_{s, M, \alpha} 
\e {\gamma^{-1} }(1 + \| {\mathcal I }_0 \|_{s + \sigma}^{k_0, \gamma})\,, 
\quad \forall s_0 \leq s \leq \bar{s}\,,
\\
  \| \Delta_{12}\mathtt{m}_{6}\|_{\frac{1}{2}, p, \alpha}
 &\lesssim_{p,M,\al} 
 \e {\gamma^{-1} } \|i_1-i_2\|_{p+\sigma}\,,
\end{aligned}
\end{equation}
for any $p$ as in \eqref{ps0}; 
The operator 
$\mathcal{B}_{6}$ has the form
\begin{equation}\label{bressan}
\mathcal{B}_6=\opw\big(B_{6}(\vphi,x,\x)\big)\,,\qquad B_{6}(\vphi,x,\x)=\ii \left(
\begin{matrix}
b_{6}(\vphi,x,\x) & 0 \\ 0 & -b_{6}(\vphi,x,-\x)
\end{matrix}
\right)\,,
\end{equation}
for some momentum preserving symbol $b_6\in S^0$. Moreover, it
satisfies
 the estimates
\begin{equation}\label{stima.b.6}
\begin{aligned}
\| b_6   \|_{0, s, \alpha}^{k_0, \gamma} 
&\lesssim_{s, M, \alpha} 
\e {\gamma^{-2} }(1 + \| {\mathcal I }_0 \|_{s + \sigma}^{k_0, \gamma})\,, 
\quad \forall s_0 \leq s \leq \bar{s}\,,
\\
    \| \Delta_{12}b_{6}\|_{0, p, \alpha}
 &\lesssim_{p,M,\al} 
 \e {\gamma^{-2} } \|i_1-i_2\|_{p+\sigma}\,,
\end{aligned}
\end{equation}
for any $p$ as in \eqref{ps0}; 

\noindent
$(ii)$
%For any $M \in \N$, $M \le N - 3(\tb+ k_{0})-3$,  and for any 
For any $\beta \in \N^\nu$, $|\beta|\le \beta_{0}$,
one has that the operator 
$\langle D \rangle^{M} \partial_\vphi^\beta  {\mathcal R}_6 \langle D \rangle^{- \tc}$ 
is a ${\mathcal D}^{k_0}$-tame and satisfies
\begin{equation}\label{stima tame cal R6 new}
\begin{aligned}
{\mathfrak M}_{ \langle D \rangle^{M}  \partial_\vphi^\beta {\mathcal R}_6 \langle D \rangle^{- \tc}}(s) 
&\lesssim_{\bar{s},M} 
\e  {\gamma^{-2}} (1 + \| {\mathcal I}_0 \|_{s + \sigma}^{k_0, \gamma})\,, 
\quad \forall s_0 \leq s \leq \bar{s}\,,
\\
\|  \langle D \rangle^{M} 
\partial_\vphi^\beta \Delta_{12}{\mathcal R}_6\langle D \rangle^{- \tc} 
\|_{\cL(H^{p} , H^{p})}& \lesssim_{\bar{s},M}
\e  {\gamma^{-2} }\|i_1-i_2\|_{p+\s}\,,
\end{aligned}
\end{equation}
for any $p$ as in \eqref{ps0}.

\noindent
$(iii)$
The maps ${\bf \Psi}^{\pm 1}$ are
 ${\mathcal D}^{k_0}$-$(k_0/2)$-tame,
 the operators
${\bf \Psi}^{\pm 1}-\id$ 
are ${\mathcal D}^{k_0}$-$((k_0+1)/2)$-tame and 
satisfy, for all $s_0 \leq s \leq \bar{s}$,
%\begin{equation}\label{stimamappaPsi}
%\begin{aligned}
%&\mathfrak{M}_{({\bf \Psi}^{\pm}-{\rm Id})}(s)
%%+\mathfrak{M}_{({\bf \Psi}^{\pm}-{\rm Id})^{*}}(s)
%\lesssim_{\bar{s},M}
%\e\g^{-2}(1+\|\mathcal{I}_{0}\|_{s+\s}^{k_0,\gamma})\,, 
%\quad  \forall \, s_0 \leq s \leq \bar{s} \, , 
%\\
%&\|\Delta_{12}({\bf \Psi}^{\pm}-{\rm Id})h\|_{p}
%%+\|\Delta_{12}({\bf \Psi}^{\pm}-{\rm Id})^{*}h\|_{p}
%\lesssim_{p,M}
%\e \g^{-2}\|h\|_{p}\|i_1-i_{2}\|_{p+\s}\,.
%\end{aligned}
%\end{equation}

\begin{equation}\label{stimamappaPsi}
\begin{aligned}
\mathfrak{M}_{{\bf \Psi}^{\pm}}(s)
&\lesssim_{\bar{s},M}
1+\|\mathcal{I}_{0}\|_{s+\s}^{k_0,\gamma}\,,
\qquad
\mathfrak{M}_{({\bf \Psi}^{\pm}-{\rm Id})}(s)
%\,,\;
&\lesssim_{\bar{s},M}
\e \g^{-2}(1+\|\mathcal{I}_{0}\|_{s+\s}^{k_0,\gamma})\,,
\\
\|\Delta_{12}{\bf \Psi}^{\pm}h\|_{p}&\lesssim_{\bar{s},M}
\e\gamma^{-2}\|h\|_{p+\s}\|i_1-i_{2}\|_{p+\s}\,,
\end{aligned}
\end{equation}
for any $p$ as in \eqref{ps0}.

\noindent
Finally, the operators $\mathcal{L}_6$, $\mathcal{B}_{6}$ and $\cR_{6}$ are real-to-real, 
Hamiltonian and momentum preserving.
\end{prop}

In order to prove Proposition \ref{riduzione diagonal ordine 1 2}, we first need the following result, 
which will be used along the whole section.

\begin{lemma}{\bf (Homological Equation).}\label{proposizione equazione omologica generale riduzione in ordini}
Let $m \in \R$, $a(\lambda; \vphi, x, \xi) \in S^m$ be 
a momentum preserving symbol depending in a Lipschitz way on the embedding $i(\vphi)$. 
Then the following hold. 

\noindent
$(i)$ The average w.r.t. $\vphi$ of $a$ is independent of $x$, namely 
\[
\langle a \rangle_\vphi := \frac{1}{(2 \pi)^\nu} \int_{\T^\nu} a(\vphi, x, \xi)\, d \vphi\,,
\]
is a Fourier multiplier $\langle a \rangle_\vphi(\xi)$. Moreover,
\begin{equation}\label{phon}
\| \langle a \rangle_\vphi \|_{m, s, \alpha}^{k_0, \gamma} \lesssim
 \| a \|_{m, s_0, \alpha}^{k_0, \gamma}, \quad \forall s,\alpha 
\geq 0\,.  
\end{equation}
Moreover if $a$ is real then $\langle a \rangle_\vphi$ is real. 

\noindent
$(ii)$ 
There exists a symbol $g= g(\omega, \th; \vphi, x,\xi) \in S^{m}$ with zero-average and satisfying the conservation of momentum, defined for all $(\omega, \th)\in \mathtt\Omega \times [\mathtt h_1, \mathtt h_2]$, satisfying the estimates 
\begin{equation}\label{stima.soluzione.omologica}
\begin{aligned}
\| g \|_{m, s, \alpha}^{k_0, \gamma} &\lesssim \gamma^{- 1} \| a \|_{m, s + \sigma, \alpha}^{k_0, \gamma}, \quad \forall s \geq 0\,, \\
\| \Delta_{12}g \|_{m,p,\alpha}  &\lesssim 
 \eps \gamma^{-2} \| i_{1} - i_{2}\|_{p+\sigma}\max_{j=1,2} \|a(i_{j})\|_{m,p+\s, \alpha}+
\gamma^{- 1} \| \Delta_{12}a \|_{m,p+\s,\alpha}\,,\;\;p\geq0\,,
\end{aligned}
\end{equation}
for some $\sigma \equiv \sigma (k_0, \tau) \gg 0$ large enough, such that
for $\lambda= (\omega, \mathtt h) \in \tT\tC_{\infty}(\gamma,\tau)$ (recall \eqref{tDtCn})
it is the solution of  
\begin{equation}\label{omologica riduzione in ordini generale}
\big( \omega \cdot \partial_\vphi + \mathtt{m}_{1}\cdot \nabla \big) g(\vphi, x, \xi) + a(\vphi, x, \xi) = \langle a \rangle_\vphi(\xi)\,.
\end{equation}
%Let $\lambda= (\omega, \mathtt h) \in \tT\tC_{\infty}(\gamma,\tau)$ (recall \eqref{tDtCn}). Then there exists a unique solution $g(\lambda; \vphi, x, \xi) \in S^m$, with zero-average and satisfying the conservation of momentum, of the equation \footnote{Note that, since $\mathtt{m}_{1}$ is constant in $(\vphi, x)$, we have that $\mathtt{m}_{1}\cdot \nabla = \opw(\ii \mathtt{m}_{1} \cdot \xi)$.}
%\begin{equation}\label{omologica riduzione in ordini generale}
%\big( \omega \cdot \partial_\vphi + \mathtt{m}_{1}\cdot \nabla \big) g(\vphi, x, \xi) + a(\vphi, x, \xi) = \langle a \rangle_\vphi(\xi)\,.
%\end{equation}
%The symbol $g$ satisfies the estimates 
%\begin{equation}\label{stima.soluzione.omologica}
%\begin{aligned}
%\| g \|_{m, s, \alpha}^{k_0, \gamma} &\lesssim \gamma^{- 1} \| a \|_{m, s + \sigma, \alpha}^{k_0, \gamma}, \quad \forall s \geq 0\,, \\
%\| \Delta_{12}g \|_{m,p,\alpha}  &\lesssim \gamma^{- 1} \| \Delta_{12}a \|_{m,p+\s,\alpha}
%\end{aligned}
%\end{equation}
%for some $\sigma \equiv \sigma (k_0, \tau) \gg 0$ large enough. 
Finally, if $a$ is real, then $g$ is real. 
\end{lemma}

\begin{proof}
$(i)$ Since $a$ is momentum preserving, one has that (recall Lemma \ref{lem:mompressimbolo})
\[
a (\vphi, x, \xi) = \widetilde{a} (\vphi - \mathtt V x\,,\, \xi) = \sum_{k \in \Z^\nu} \widehat{\widetilde{a}} (k, \xi) 
e^{\ii k \cdot \big( \vphi - \mathtt V x \big)}\,. 
\]
By using this and computing the time average of $a$, namely 
\[
\langle a (\cdot, x, \xi)\rangle_\vphi := \frac{1}{(2 \pi)^\nu} \int_{\T^\nu} a (\vphi, x, \xi)\, d \vphi\,,
\]
then we obtain that 
\[
\begin{aligned}
\frac{1}{(2 \pi)^\nu} \int_{\T^\nu} a (\vphi, x, \xi)\, d \vphi & = \frac{1}{(2 \pi)^\nu} 
\sum_{k \in \Z^\nu} \widehat{\widetilde{a}}  (k, \xi)  \int_{\T^\nu} e^{\ii k \cdot \big( \vphi - \mathtt V x \big)}\, d \vphi 
\\
&=  \frac{1}{(2 \pi)^\nu}\sum_{k \in \Z^\nu} \widehat{\widetilde{a}}  (k, \xi)\, 
e^{- \ii \mathtt V x \cdot k}  \int_{\T^\nu} e^{\ii k \cdot  \vphi }\, d \vphi = \widehat{\widetilde{a}}  (0, \xi)\,. 
\end{aligned}
\]
Therefore, we have the bound 
$\|\langle a \rangle_\vphi\|_{m, s, \alpha}^{k_0, \gamma} 
= \|\widehat{\widetilde{a}}(0, \cdot)\|_{m, s, \alpha}^{k_0, \gamma} 
\lesssim \| \widetilde{a} \|_{m, s_0, \alpha}^{k_0, \gamma}$ , for any $s\geq 0$. 

\noindent
$(ii)$ We shall now construct solutions of the equation \eqref{omologica riduzione in ordini generale}. We look for a momentum preserving symbol $g$, namely $g(\vphi, x, \xi) \equiv \widetilde{g}(\vphi - \mathtt V x, \xi)$. Hence we have 
\[
\begin{aligned}
\big( \omega \cdot \partial_\vphi +  \mathtt{m}_{1} \cdot \nabla \big) g(\vphi, x, \xi) 
& = \sum_{k \in \Z^\nu} \widehat{\widetilde{g}}(k, \xi) \big( \omega \cdot \partial_\vphi 
+ \mathtt{m}_{1} \cdot \nabla \big) \big[ e^{\ii  \vphi \cdot k} e^{- \ii \mathtt V x \cdot k} \big]    
\\& 
= \sum_{k \in \Z^\nu} \ii \Big( \omega \cdot k - \mathtt{m}_{1} \cdot {\mathtt V}^T k \Big) 
\widehat{\widetilde{g}}(k, \xi) e^{\ii \vphi \cdot k } e^{- \ii  x \cdot {\mathtt V}^T k} 
\\& 
= \sum_{k \in \Z^\nu} \ii \big( \omega -  {\mathtt V} \mathtt{m}_{1}\big) \cdot k \,
\widehat{\widetilde{g}}(k, \xi) e^{\ii k \cdot \big( \vphi - \mathtt V x \big)}
\,. 
\end{aligned}
\]
Therefore, the \eqref{omologica riduzione in ordini generale} in terms of $\widetilde{a}$ and  $\widetilde{g}$ reads as
\begin{equation}\label{solsolsol}
\big( \omega -  {\mathtt V} \mathtt{m}_{1}\big)  \cdot \pa_{\theta} \widetilde{g}(\theta, \xi) 
+  \widetilde{a} (\theta, \xi) =  \langle \widetilde{a} \rangle_\theta(\xi)\,.
\end{equation}

Recalling the operator $( (\omega - \tV \mathtt{m}_{1}  )\cdot \pa_{\theta})_{ext}^{- 1}$ 
defined in \eqref{def ompaph-1 ext2}
%Lemma \ref{lemma:WD2} 
we define the traveling function $\widetilde g$ as follow
\begin{equation}\label{piripicchio}
 \widetilde g :=  \Big( (\omega - \tV \mathtt{m}_{1}  )\cdot \pa_{\theta}\Big)_{ext}^{- 1}\Big[ \langle \widetilde{a} \rangle_\theta(\xi) - \widetilde{a}(\theta, \xi) \Big] \,.
\end{equation}
Lemma \ref{lemma:WD2} applies. Hence 
the estimates \eqref{stima.soluzione.omologica} follows by \eqref{piripicchio},
\eqref{2802.2nuova}
using  that $g(\vphi, x, \xi) \equiv \widetilde{g}(\vphi - \mathtt V x, \xi)$.
%The estimates on the Lipschitz variation can be checked using the equation
%\eqref{omologica riduzione in ordini generale} and reasoning as above.
We finally estimate the Lipschitz variation of the solution $g$.
\[
\begin{aligned}
\| \Delta_{12} \widetilde g \|_{m, p, \alpha} 
&=  \| \Delta_{12}\big[ ( (\omega - \tV \mathtt{m}_{1}  )\cdot \pa_{\theta})^{- 1}_{ext}
 [ \langle \widetilde{a} \rangle_\theta(\xi) - \widetilde{a}(\theta, \xi) ]\big]\|_{m, p, \alpha}
 \\
 &\le \| \Delta_{12} \big( (\omega - \tV \mathtt{m}_{1}(i_{1})  )\cdot \pa_{\theta}\big)^{- 1}_{ext}
 \big[ \langle \widetilde{a} \rangle_\theta(\xi) - \widetilde{a}(\theta, \xi)  \big] \|_{m, p, \alpha}
 \\
 &+  \| \big( (\omega - \tV \mathtt{m}_{1}  )\cdot \pa_{\theta}\big)^{- 1}_{ext}
 \Delta_{12} \big[ \langle \widetilde{a} \rangle_\theta(i_{2}; \xi) - \widetilde{a}(i_{2}; \theta, \xi) \big] \|_{m, p, \alpha}
 \\
 & \stackrel{\eqref{2802.2nuova},\eqref{2802.2nuova2}}{\lesssim} 
 \eps \gamma^{-2} \| i_{1} - i_{2}\|_{p+\sigma} 
 \max_{j=1,2}\|a(i_j)\|_{m, p+2 \tau +1, \alpha} +  \gamma^{-1} \|\Delta_{12} a\|_{m, p+\sigma, \alpha} \,.
 \end{aligned}
\]

%Then
%for any value of the parameter $\lambda= (\omega , \mathtt h) \in\tT\tC_{\infty}(2\gamma,\tau)$ 
%the solution of the above equation is the periodic function 
%\begin{equation}\label{soluzione omologica g ordine 12}
%\begin{aligned}
% \widetilde g(\theta, \xi) &  := \Big( (\omega - \tV \mathtt{m}_{1}  )\cdot \pa_{\theta}\Big)^{- 1}
% \Big[ \langle \widetilde{a} \rangle_\theta(\xi) - \widetilde{a}(\theta, \xi) \Big]  \,,
%\end{aligned}
%\end{equation}
%which we extend to the whole parameter space $\mathtt{\Lambda}_{0}$ (recall \eqref{def:setlambda0}) by setting 
%\[
% \widetilde g_{ext} :=  \Big( (\omega - \tV \mathtt{m}_{1}  )\cdot \pa_{\theta}\Big)_{ext}^{- 1}\Big[ \langle \widetilde{a} \rangle_\theta(\xi) - \widetilde{a}(\theta, \xi) \Big] 
%\]
%via the operator $( (\omega - \tV \mathtt{m}_{1}  )\cdot \pa_{\theta})_{ext}^{- 1}$ 
%defined in Lemma \ref{lemma:WD2}. %with $\oo \rightsquigarrow  \oo - \mathtt V \mathtt{m}_{1} $. 
%For simplicity we still denote by $\widetilde g$ this extension.
Recalling the definition of $( (\omega - \tV \mathtt{m}_{1}  )\cdot \pa_{\theta})_{ext}^{- 1}$  in \eqref{def ompaph-1 ext2} and the cut-off function in \eqref{cut off simboli 1}
one has that for any value of the parameter $\lambda= (\omega , \mathtt h) \in\tT\tC_{\infty}(2\gamma,\tau)$ the function 
$\widetilde g$ is solution of \eqref{solsolsol}.
\end{proof}

\begin{proof}[{\bf Proof of Proposition \ref{riduzione diagonal ordine 1 2}}]
 First of all we consider a real symbol $g(\vphi, x, \xi) \in S^{\frac12}$ (to be determined) and we consider the flow $\Psi^\tau(\vphi)$ which is defined by 
\begin{equation}\label{def Phi tau ordine 1 2}
\begin{cases}
\partial_\tau \Psi^\tau_g(\vphi)  = \ii \Pi_{S}^{\perp}{\rm Op}^W \big( g(\vphi, x, \xi) \Big) \circ [\Pi_{S}^{\perp}\Psi_g^\tau(\vphi)] \\
\Psi^0_g(\vphi) = {\rm Id}\,. 
\end{cases}
\end{equation}

By Lemma \ref{differenzaFlussi} we have that 
the flow $\Psi^\tau_{g}(\vphi)$ is the correction of the flow $\Phi^\tau_{g}$of the equation
\begin{equation}\label{tagliata}
\partial_\tau \Phi^\tau_g(\vphi)  = \ii \opw \Big( g(\vphi, x, \xi) \Big) \circ \Phi_g^\tau(\vphi)\,, \quad 
\Phi^0_g(\vphi) = {\rm Id}\,,
\end{equation}
up to a finite rank operator. More precisely,
\[
\Psi_{g}^{\tau} =\Pi_{S}^{\perp} \Phi_{g}^{\tau} \Pi_{S}^{\perp}\circ(\uno +\mathcal{R})\,,
\]
where $\Phi_{g}$ is the $1$-time flow of the equation \eqref{tagliata} and
$\mathcal{R}$ is a finite rank operator of the form \eqref{forma buona resto}
satisfying \eqref{giallo2}-\eqref{giallo3}.

We choose the symbol $g$ in such a way that 
\begin{equation}\label{equazione omologica g ordine 1 2}
 \big( \omega \cdot \partial_\vphi + \mathtt{m}_{1} \cdot \nabla \big) g(\vphi, x, \xi) 
+ \lambda(\vphi, x, \xi)-\Omega(\x)  = \langle \lambda  -\Omega\rangle_\vphi(\xi)\,. 
\end{equation}
We shall use the conservation of momentum of $\lambda$ and we apply Lemma
\ref{proposizione equazione omologica generale riduzione in ordini}
with $a\rightsquigarrow   \lambda(\vphi,x,\x)-\Omega(\x)$.
We deduce the existence of a symbol 
$g \in S^{\frac12}$, defined for all $(\omega,\mathtt{h})\in \R^{\nu}\times [\th_1,\th_2]$, with the following properties:

\noindent
$(i)$ using Lemma \ref{lemma autovalori decoupling ordine alto} to estimate the symbol $\lambda-\Omega(\x)$,
%using \eqref{stime lambda pm fpm}, 
the bound \eqref{stima.soluzione.omologica}
and taking $\e\gamma^{-1}\ll1$\,, we get, for all $(\omega,\mathtt{h})\in \mathtt\Omega \times [\mathtt h_1, \mathtt h_2]$
\begin{equation*}%\label{stime Omega 6 generatore g}
\begin{aligned}
%& \| \langle \lambda_5 - \Omega \rangle_\vphi \|_{\frac12, s, \alpha}^{k_0, \gamma}  \lesssim_{s, \alpha} \e\gamma^{-1} (1 + \| {\mathcal I }_0 \|_{s + \sigma}^{k_0, \gamma})\,, \\
& \| g \|_{\frac12, s , \alpha}^{k_0, \gamma} 
\lesssim_{s,M, \alpha} 
\e \gamma^{- 2} \big( 1 + \| {\mathcal I}_{0} \|_{s + \sigma}^{k_0, \gamma}\big)\,, \quad s_0 \leq s \leq \bar{s}\,,\\
& \| \Delta_{12}g \|_{\frac12, p, \alpha} \lesssim_{p,M,\al} 
\e \gamma^{- 2} \|i_{2}-i_{1}\|_{p+\s}\,,
\end{aligned}
\end{equation*}
for some $\s\gg1$ large enough and for any $p$ as in \eqref{ps0}.

\noindent
$(ii)$ Since by 
Lemma \ref{lemma autovalori decoupling ordine alto} the symbol $\lambda$ is even in $\x$
(and also $\Omega(\x)$), also $g$ is so.

\noindent
$(iii)$ For $(\omega, \mathtt h) \in \tT\tC_{\infty}(\gamma,\tau)$ (recall \eqref{tDtCn})
such symbol $g(\vphi,x,\x)$ satisfies the equation \eqref{equazione omologica g ordine 1 2}.

%
%.  by the estimates \eqref{stime lambda pm fpm}, 
%one obtains that satisfies \eqref{phon} and that there is 
%$g \in S^{\frac12}$, defined for all $(\omega,\mathtt{h})\in \R^{\nu}\times [\th_1,\th_2]$, such that 
%\begin{equation*}%\label{stime Omega 6 generatore g}
%\begin{aligned}
%%& \| \langle \lambda_5 - \Omega \rangle_\vphi \|_{\frac12, s, \alpha}^{k_0, \gamma}  \lesssim_{s, \alpha} \e\gamma^{-1} (1 + \| {\mathcal I }_0 \|_{s + \sigma}^{k_0, \gamma})\,, \\
%& \| g \|_{\frac12, s , \alpha}^{k_0, \gamma} \lesssim_{s, \alpha} \e \gamma^{- 2} \big( 1 + \| {\mathcal I}_{0} \|_{s + \sigma}^{k_0, \gamma}\big)\,, \quad s_0 \leq s \leq \bar{s}\,,\\
%& \| \Delta_{12}g \|_{p} \lesssim_{p} \e \gamma^{- 2} \|i_{2}-i_{1}\|_{p}\,. 
%\end{aligned}
%\end{equation*}

We now set 
\begin{equation}\label{def cal B}
{\bf\Psi }:= \begin{pmatrix}
 \Psi_g & 0 \\
 0 & \overline{\Psi_g}
 \end{pmatrix}, \quad \Psi_g \equiv \Psi_g^1\,.
% \qquad 
% {\bf\Phi }:= \begin{pmatrix}
% \Phi_g & 0 \\
% 0 & \overline{\Phi_g}
% \end{pmatrix}, \quad \Phi_g \equiv \Phi_g^1\,.
\end{equation} 

In particular, the map ${\bf \Psi}$ is a real-to-real,  invertible operator symplectic on
$ {\bf H}^s_{\perp}$, $s_0 \leq s \leq \bar{s}$
an it is momentum preserving since the generator is a momentum preserving symbol.
Moreover, by Lemmata \ref{lemma:buonaposFlussi}, \ref{differenzaFlussi} 
(recalling also Remark \ref{differenzaFlussiRMK}) 
it also satisfies the estimates in \eqref{stimamappaPsi}.
We now compute the conjugated operator 
%${\mathcal L}_6 :=
${\bf \Psi}^{- 1} \circ {\mathcal L}_5 \circ {\bf \Psi}$.
First of all, recalling \eqref{op cal L5} (see also \eqref{autovalori matrice C 1/2}),  
we define 
\begin{equation}\label{scalarstep56}
{\bf P}_1:=\left(\begin{matrix}\mathcal{P}_1 & 0 \\0 & \overline{\mathcal{P}_1}\end{matrix}\right)\,,
\qquad
\mathcal{P}_1:=\omega\cdot\pa_{\vphi}+\mathtt{m}_1\cdot\nabla
+\ii \opw(\lambda+b_{5}(\vphi,x,\x))\,.
\end{equation}
Hence we write 
\[
\begin{aligned}
{\bf \Psi}^{- 1} \circ {\mathcal L}_5 \circ {\bf \Psi} &=
{\bf \Psi}^{- 1} \circ{\bf \Pi}_{S}^{\perp} {\bf P}_1 {\bf \Pi}_{S}^{\perp}  \circ {\bf \Psi}
+
{\bf \Psi}^{- 1} \circ{\bf \Pi}_{S}^{\perp}\cR_{5} {\bf \Pi}_{S}^{\perp}\circ {\bf \Psi}
\\&
=\left(\begin{matrix}\mathcal{P}_2 & 0 \\0 & \overline{\mathcal{P}_2}\end{matrix}\right)
+{\bf \Psi}^{- 1} \circ{\bf \Pi}_{S}^{\perp}\cR_{5}{\bf \Pi}_{S}^{\perp} \circ {\bf \Psi}
\end{aligned}
\] 
where we defined (recall \eqref{projComplex})
\[
\mathcal{P}_2:=\Psi_{g}^{-1}\circ {\Pi}_{S}^{\perp}\mathcal{P}_{1}{ \Pi}_{S}^{\perp}\circ \Psi_{g}\,.
\]
By applying Lemma \ref{georgiaLem} we notice that
\[
\mathcal{P}_2=\Pi_{S}^{\perp}
\Phi_{g}^{-1}\mathcal{P}_{1}\Phi_{g}
\Pi_{S}^{\perp}+\mathcal{Q}_{6}\,,
\]
where $\cQ_{6}$ is a finite rank operator satisfying the estimates \eqref{stima tame cal R6 new}.
Moreover, by \eqref{scalarstep56}, the 
construction of $g$ and by using Lemma \ref{flussi coniugi eccetera}  (recall also estimates  on $\lambda$
in Lemma 
\ref{lemma autovalori decoupling ordine alto} and \eqref{stima d dopo decoupling} on $b_{5}$),
we have that
\[
\Phi_{g}^{-1}\mathcal{P}_{1}\Phi_{g}=
\oo\cdot \partial_\vphi + \mathtt{m}_{1}\cdot\nabla 
+\ii  \opw(  \mathtt{q} + b_6 )+\widetilde{R}_1
\]
where 
\begin{equation}\label{def:simboloqq}
\mathtt{q}:=\mathtt{q}(\vphi,x,\x):=
\big( \omega \cdot \partial_\vphi + \mathtt{m}_{1} \cdot \nabla \big) g(\vphi, x, \xi) + \lambda(\vphi, x, \xi) 
%\stackrel{\eqref{equazione omologica g ordine 1 2}}{\equiv}
%\Omega(\x)  + \langle \lambda  -\Omega\rangle_\vphi(\xi)\,,
\end{equation}
and where the symbol 
$b_6 \in S^{0}$ satisfyies estimates like in \eqref{stima.b.6} 
and the remainder $\widetilde{R}_{1}$ 
%is a smoothing in space operator 
satisfies \eqref{stima tame cal R6 new}.
We remark that the symbol $\mathtt{q}$ is even in $\x$ and real valued 
thanks to the properties of $\lambda $ and $g$
constructed above.
We also define the matrix of pseudo-differential operators $\mathcal{B}_{6}$ as in \eqref{bressan}.

Let us now consider
the term  ${\bf \Psi}^{- 1} \circ{\bf \Pi}_{S}^{\perp}\cR_{5}{\bf \Pi}_{S}^{\perp} \circ {\bf \Psi}=  \widetilde{\cR}_{3}$
which is the only one containing ``off-diagonal'' terms (which are smoothing operators).
This term satisfies 
\eqref{stima tame cal R6 new} in view of \eqref{stima tame cal R5 new} and using Lemma 
\ref{conj.NEWsmoothresto} and Remark \ref{rmk:restimatrici}-(i).

The discussion above implies that 
\[
\begin{aligned}
{\bf \Psi}^{- 1} &\circ {\mathcal L}_5 \circ {\bf \Psi}= {\bf \Pi}_{S}^{\perp}\omega \cdot \partial_\vphi + 
{\bf \Pi}_{S}^{\perp}\Big( 
\ii\opw(\mathtt{m}_{1}  \cdot \xi) \uno  
+ \ii E\opw(\mathtt{q})+\cB_{6} + \cR_{6} \Big) {\bf \Pi}_{S}^{\perp}\,,
\end{aligned}
\]
 where $\mathtt{q}$ is in \eqref{def:simboloqq} and where we set
 \[
 \cR_{6}=\widetilde{\cR}_{3}+\sm{\mathcal{Q}_{6} + \widetilde{R}_{1}}{0}{0}{\overline{\mathcal{Q}_{6}}+\overline{\widetilde{R}_{1}}}\,.
 %+\left(\begin{matrix} \mathcal{Q}_{6} + \widetilde{R}_{1}& 0 \\ 0 
 %& \overline{\mathcal{Q}_{6}}+\overline{\widetilde{R}_{1}}\end{matrix}\right)\,.
 \]
  Since the map ${\bf \Psi}$ is generated by the real, momentum preserving symbol $g$, 
  then it is real-to-real, Hamiltonian and momentum 
  preserving and therefore ${\bf \Psi}^{- 1} \circ {\mathcal L}_5 \circ {\bf \Psi}$ 
  is so.
  In view of Remark \ref{rmk:algsimboli} we have also that the pseudo-differential term
 $\mathcal{B}_6$ is
 real-to-real, Hamiltonian and momentum preserving.
Finally, by difference $\cR_{6}$ satisfies the same algebraic properties. 
In order to conclude the proof, we define
 $\mathtt{m}_{6}:= \jap{\lambda}_{\vphi}(\xi)$ and the operator 
 (defined for all values $(\omega,\th)\in \mathtt{\Omega}\times[\th_1,\th_2]$)
 \[
 {\mathcal L}_6 = {\bf \Pi}_{S}^{\perp}\omega \cdot \partial_\vphi + 
{\bf \Pi}_{S}^{\perp}\Big( 
\ii\opw(\mathtt{m}_{1}  \cdot \xi) \uno  
+ \ii E\opw(\mathtt{m}_{6})+\cB_{6} + \cR_{6} \Big) {\bf \Pi}_{S}^{\perp}\,.
 \]
 First of all, note that by 
 \eqref{phon} and Lemma \ref{lemma autovalori decoupling ordine alto} for the symbol $\lambda$,
 one has that the symbol $\mathtt{m}_6$ satisfies the estimates \eqref{stima.m.6}. 
 Finally, by \eqref{equazione omologica g ordine 1 2}, \eqref{def:simboloqq} and
 taking $(\omega, \th)\in \tT\tC_{\infty}(\gamma,\tau)$, than one has that 
 %(recall \eqref{def:simboloqq})
  $\mathtt{q}(\vphi,x,\x)\equiv \mathtt{m}_{6}(\x)$,
 implying the conjugation result in \eqref{op cal L6}.
  The discussion above implies the thesis.
\end{proof}

\subsection{Reduction up to smoothing remainders}\label{sec:reductionlower}
In the next proposition we reduce to constant coefficients the linear operator ${\mathcal L}_6$ up to an arbitrarily smoothing remainder. 
\begin{prop}{\bf (Reduction to constant coefficients of the lower order terms.)}\label{riduzione diagonal ordini bassi} 
Let  
%and fix $\tb, N \in \N$ such that $N\gg 3(\tb+k_0)+3$. 
 $\beta_{0} \in \N$, $M \in \N$, $\mathtt c \geq 3(k_0 + \beta_0 ) + 3$ with $M \geq \mathtt c$.
Then for any $\alpha\in \N$ there exists $\sigma = \sigma(\beta_0, M, \alpha) \gg 1$ such that, for any 
$\bar{s}\ge s_{0}$, if 
%Then for any $\alpha\in\N$ there exist $\sigma \equiv \sigma(\tb, N, \alpha) \gg 0$, 
%$\mathtt c \equiv \mathtt c (N) \gg 0$ large enough such that if 
condition \eqref{ansatz_I0_s0}
holds with $\mu_0\gg\s$ and if \eqref{smallcondepsilon} holds with $\kappa\ge 2(M + \tc )+2$ then the following holds. 
There exist
a symbol $\mathtt{m}_{7}\in S^{\frac{1}{2}}$ independent of $(\vphi,x)$, 
 a 
$\mathcal{D}^{k_0}$-tame operator ${\mathcal R}_7$ and
a real-to-real, invertible, symplectic and momentum preserving map 
$\cW : {\bf H}^s_{\perp} \to {\bf H}^s_{\perp}$, $s_0 \leq s \leq \bar{s}$ 
defined for all $(\omega,\mathtt{h})\in \mathtt \Omega\times[\th_1,\th_2]$ such that the following holds.
For any $(\omega,\mathtt{h})\in 	\tT\tC_{\infty}(\gamma,\tau) $ (see \eqref{tDtCn})
one has the conjugation 
(recalling \eqref{op cal L6})
\begin{equation}\label{op cal L7}
\begin{aligned}
{\mathcal L}_7 & = \cW^{- 1} \circ {\mathcal L}_6 \circ \cW
\\
&= {\bf \Pi}_{S}^{\perp}\omega \cdot \partial_\vphi + 
{\bf \Pi}_{S}^{\perp}\Big( \mathtt{m}_{1}  \cdot \nabla \uno  
 + \ii \opw\left(
 \begin{matrix}
 \mathtt{m}_{7}(\x) & 0 \\0 & -\mathtt{m}_{7}(-\x)\end{matrix}
 \right)+ \cR_{7} \Big) {\bf \Pi}_{S}^{\perp}
\end{aligned}
\end{equation}
where $\mathtt{m}_{1}$  is given in Proposition \ref{riduzione trasporto} and 
satisfies estimates in \eqref{beta.FGMP.est}. 
Moreover, for any $(\omega,\mathtt{h})\in \mathtt \Omega\times[\th_1,\th_2]$
one has the following estimates:

 \noindent 
 $(i)$  The Fourier multiplier $\mathtt{m}_{7}\in S^{\frac12}$ is a real valued, momentum preserving symbol.
%  such that
%\[
%\mathtt{m}_{7}( \xi)=\mathtt{m}_{7}(-\xi);
%\]
% the operator $\cC_{7}$ has the form
% \begin{equation}\label{trapani}
%\mathcal{C}_7=\opw\big(C_{7}(\vphi,x,\x)\big)\,,\quad C_{7}(\vphi,x,\x)=\left(
%\begin{matrix}
% c_{7}(\vphi,x,\x) & 0 \\ 0 &  - c_{7}(\vphi,x,-\x)
%\end{matrix}
%\right)\,,
%\end{equation}
%for some momentum preserving symbols  $c_{7}\in S^{-\rho}$. 
Moreover, it 
satisfies
 the estimates (recall \eqref{wild})
\begin{equation}\label{stima lambda 7}
\begin{aligned}
\| \mathtt{m}_{7} - \Omega   \|_{\frac12, s, \alpha}^{k_0, \gamma} 
%\,,\,\| c_7   \|_{-\rho, s, \alpha}^{k_0, \gamma} 
&\lesssim_{s, M, \alpha} 
\e {\gamma^{- 2(M - \tc )-1}} (1 + \| {\mathcal I }_0 \|_{s + \sigma}^{k_0, \gamma})\,, 
\quad \forall s_0 \leq s \leq \bar{s}\,,
\\
  \| \Delta_{12}\mathtt{m}_{7}\|_{\frac12, p, \alpha}
  %\,,\,  \| \Delta_{12}c_{7}\|_{-\rho, p, \alpha}
 &\lesssim_{p,M,\al} 
 \e {\gamma^{- 2(M - \tc )-1}} \|i_1-i_2\|_{p+\sigma}\,,
\end{aligned}
\end{equation}
for any $p$ as in \eqref{ps0}.

\noindent
$(ii)$
%For any $M \in \N$, $M \le N - 3(\tb+ k_{0})-3$,  and for any $\beta \in \N^\nu$, $|\beta|\le \beta_{0}$,
For any $\beta \in \N^\nu$, $|\beta|\le \beta_{0}$,
one has that the operator 
$\partial_\vphi^\beta \langle D \rangle^{M} {\mathcal R}_7 \langle D \rangle^{- \tc}$ 
is a ${\mathcal D}^{k_0}$-tame and satisfies
\begin{equation}\label{stima tame cal R7 new}
\begin{aligned}
{\mathfrak M}_{ \langle D \rangle^{M} \partial_\vphi^\beta{\mathcal R}_7 \langle D \rangle^{- \tc}}(s) 
&\lesssim_{\bar{s},M} \e {\gamma^{- 2(M - \tc )-2}}(1 + \| {\mathcal I}_0 \|_{s + \sigma}^{k_0, \gamma})\,, 
\quad \forall s_0 \leq s \leq \bar{s}\,,
\\
\|  \langle D \rangle^{M} 
\partial_\vphi^\beta \Delta_{12}{\mathcal R}_7\langle D \rangle^{- \tc} 
\|_{\cL(H^{p} , H^{p})}& \lesssim_{\bar{s},M}\e  {\gamma^{- 2(M - \tc )-2}} \|i_1-i_2\|_{p+\s}\,,
\end{aligned}
\end{equation}
for any $p$ as in \eqref{ps0}.

\noindent
$(iv)$
The maps $\cW^{\pm 1}$ are
 ${\mathcal D}^{k_0}$-tame and satisfy
\begin{equation}\label{stimamappaPsi2}
\begin{aligned}
&\mathfrak{M}_{(\cW^{\pm 1}-{\rm Id})}(s)
%+\mathfrak{M}_{(\cW^{\pm 1}-{\rm Id})^{*}}(s)
\lesssim_{\bar{s},M}
\e {\gamma^{- 2(M - \tc )-2}}(1+\|\mathcal{I}_{0}\|_{s+\s}^{k_0,\gamma})\,, 
\quad  \forall \, s_0 \leq s \leq \bar{s} \, , 
\\
&\|\Delta_{12}(\cW^{\pm 1}-{\rm Id})h\|_{p}
%+\|\Delta_{12}(\cW^{\pm 1}-{\rm Id})^{*}h\|_{p}
\lesssim_{\bar{s},M}
\e  {\gamma^{- 2(M - \tc )-2}} \|h\|_{p}\|i_1-i_{2}\|_{p+\s}\,.
\end{aligned}
\end{equation}

\noindent
Finally the operators $\mathcal{L}_7$, $\mathcal{P}_{7}$ and $\mathcal{R}_{7}$ are real-to-real, 
Hamiltonian and momentum preserving.
\end{prop}

\begin{proof}
The proof is based on an iterative procedure in which we remove the dependence on 
$(\vphi, x)$ (in decreasing order) from the symbol $b_{6}(\vphi, x, \xi)$ 
given in Proposition \ref{riduzione diagonal ordine 1 2}. 
In order to achieve this purpose, it is enough to implement a 
normal form reduction scheme on the operator $\mathcal{L}_6$ defined 
in \eqref{op cal L6}.

We shall describe the induction step of such an iterative procedure. 
At the $k$-th step, $k\geq0$, we deal with the operator $\mathcal{L}_6^{(k)}$ which has the form 
\begin{equation}\label{forma L6 (k)}
\begin{aligned}
\mathcal{L}_6^{(k)} &= {\bf \Pi}_{S}^{\perp}\omega \cdot \partial_\vphi 
+{\bf \Pi}_{S}^{\perp}\Big(\mathtt{m}_{1}  \cdot \nabla 
+ \ii {\rm Op}^W \left(\begin{matrix}\mathtt{m}_{6}^{(k)}(\xi) & 0 \\0 & -\mathtt{m}_{6}^{(k)}(-\xi)\end{matrix}\right) 
+ \mathcal{B}_{6}^{(k)}
+\mathcal{R}_{6}^{(k)}
\Big){\bf \Pi}_{S}^{\perp}
\,,
\end{aligned}
\end{equation}
defined for all $(\omega,\mathtt{h})\in \mathtt\Omega \times [\mathtt h_1, \mathtt h_2]$ and satisfying 
 the following properties. 
Given $M\in \N$, there exists a non-decreasing sequence $\sigma_k=\sigma_k(M)$ 
such that if $\| \cI_{0} \|_{s_{0}+\sigma_k}^{k_{0}, \g}\le 1$ , $\eps \g^{- k-2}\lesssim_{s}1$, the following hold. 
\begin{itemize}
\item{$({\mathcal P}1)_k$} The Fourier multiplier $\mathtt{m}_{6}^{(k)}\equiv\mathtt{m}_{6}^{(k)}(\x) \in S^{\frac12}$ 
is real and it satisfies the estimates
\begin{equation}
\begin{aligned}
\| \mathtt{m}_{6}^{(k)} - \Omega \|_{\frac12, s, \alpha}^{k_0, \gamma} 
&\lesssim_{s,k, \alpha,M} 
\e {\g^{-2-(k-1)}}(1 + \| {\mathcal I}_0 \|_{s + \sigma_k}^{k_0, \gamma}), \quad s_0 \leq s \leq \bar{s}
\\
\|\Delta_{12} \mathtt{m}_{6}^{(k)}  \|_{\frac12, p, \alpha}
&\lesssim_{p, k,\al,M} \e{\g^{-2-(k-1)}}\|i_{1} - i_{2} \|_{s + \sigma_k}\,,
\end{aligned}
\end{equation}
for any $p$ as in \eqref{ps0}.
\item{$({\mathcal P}2)_k$} 
$\mathcal{B}_{6}^{(k)}=\opw\big(B_{6}(\vphi,x,\x)\big)$
with
\begin{equation}\label{Zippoblu}
B_{6}(\vphi,x,\x):=\ii\left(\begin{matrix} b_{6}^{(k)}(\vphi, x, \xi) & 0 
\\ & -b_{6}^{(k)}(\vphi, x, -\xi) \end{matrix}\right)\,,
\end{equation}
where 
the symbol $b_{6}^{(k)} \in S^{- \frac{k}{2}}$ is real, momentum preserving and it satisfies the estimates 
\begin{equation}\label{stimeB6kappino}
\begin{aligned}
\| b_{6}^{(k)} \|_{- \frac{k}{2}, s, \alpha}^{k_0, \gamma}
& \lesssim_{s, k,\alpha,M} 
\e\g^{-2-k} (1 + \| {\mathcal I}_0 \|_{s + \sigma_k}^{k_0, \gamma}), \quad s_0 \leq s \leq \bar{s} \\
\|\Delta_{12} b_{6}^{(k)}  \|_{-\frac{k}{2}, p, \alpha} 
&\lesssim_{p, k,\al,M} 
\e \g^{-2-k}\|i_{1} - i_{2} \|_{p+ \sigma_k}\,,
\end{aligned}
\end{equation}
for any $p$ as in \eqref{ps0}. 

\item{$({\mathcal P}3)_k$} 
%For any $M \in \N$, $M \le N - 3(\tb + k_{0})-3$,  
%and for any $\beta \in \N^\nu$, $|\beta|\le \tb$,
For any $\beta \in \N^\nu$, $|\beta|\le \beta_{0}$,
 the operator  $  \langle D \rangle^{M}
 \partial_\vphi^\beta{\mathcal R}_6^{(k)} \langle D \rangle^{-\tc}$ 
 is ${\mathcal D}^{k_0}$- tame operator and it satisfies the estimates
\begin{equation}\label{punkina1KKK}
\begin{aligned}
{\mathfrak M}_{ \langle D \rangle^{M}
\partial_\vphi^\beta{\mathcal R}_6^{(k)} 
\langle D \rangle^{- \tc}}(s) 
&\lesssim_{\bar{s}, k,M} 
\e \g^{-2-k} (1 + \| {\mathcal I}_0 \|_{s + \sigma_k}^{k_0, \gamma})\,, 
\quad \forall s_0 \leq s \leq \bar{s}\,,
\\
\|  \langle D \rangle^{M} 
\partial_\vphi^\beta \Delta_{12}{\mathcal R}_6^{(k)}&\langle D \rangle^{- \tc} 
\|_{\cL(H^{p} , H^{p})}
\lesssim_{\bar{s},k,M}
\e \g^{-2-k}\|i_1-i_2\|_{p+\s}\,,
\end{aligned}
\end{equation}
for any $p$ as in \eqref{ps0}.

\item{$({\mathcal P}4)_k$} For $k\geq1$, 
there exists a real valued momentum preserving symbol $f_{k}\in S^{-\frac{k}{2}}$,
defined for all $(\omega, \mathtt h) \in \mathtt\Omega \times [\mathtt h_1, \mathtt h_2]$,
for any $(\omega, \mathtt h) \in \tT\tC_{\infty}(\gamma,\tau)$ (recall \eqref{tDtCn})
such that
\begin{equation}\label{coniugiostepkDIAGO}
\mathcal{L}_{6}^{(k)}
=\mathcal{W}_{k-1}^{-1}\circ \mathcal{L}_{6}^{(k-1)}\circ\mathcal{W}_{k-1}^{1}
\end{equation}
where $\mathcal{W}_{k-1}^{\theta}$, $\theta\in[0,1]$ 
is the flow at time 
$\theta\in [0,1]$ of 
\begin{equation}\label{flow-jthDIAGO}
\begin{aligned}
\partial_{\theta} \mathcal{W}_{k-1}^{\theta} (\vphi) &= 
 {\bf \Pi}_{S}^{\perp}
 \opw\begin{pmatrix}
\ii f_{k-1}(\vphi,x,\x) &  0 \\
0 & -\ii f_{k-1}(\vphi,x,-\x)
\end{pmatrix}
{\bf \Pi}_{S}^{\perp}\mathcal{W}_{k-1}^{\theta}(\vphi) \, ,
\\
\mathcal{W}_{k-1}^{0}(\vphi) &= {\rm Id} \, . 
 \end{aligned}
\end{equation}

\item{$({\mathcal P}5)_k$} The operator $\mathcal{L}_{6}^{(k)}$ is real-to-real, Hamiltonian and momentum 
preserving as well as the operators
$\mathcal{B}^{(k)}_{6}$,   and ${\mathcal R}_{6}^{(k)}$.
\end{itemize}

\smallskip
At $k=0$ we have that the operator $\mathcal{L}_{6}$ in \eqref{op cal L6}
has the form $\mathcal{L}_{6}^{(0)}$ in \eqref{forma L6 (k)}
with $\mathtt{m}_6^{(0)}:=\mathtt{m}_{6}$,
${\mathcal B}_6^{(0)}:=\mathcal{B}_{6}$ (see \eqref{bressan}) 
and ${\mathcal R}_6^{(k)}:={\mathcal R}_6$.
The properties $({\mathcal P}1)_0$, $({\mathcal P}2)_0$, $({\mathcal P}3)_0$, $({\mathcal P}5)_0$
follow by Proposition \ref{riduzione diagonal ordine 1 2} while $({\mathcal P}4)_0$ in empty.

\smallskip
Let us now assume  $({\mathcal P}j)_k$ with $k\geq1$ and $j=1,\ldots,5$.
In order to prove inductively the assertions above,
we want to choose the symbol $f_k$  in order to eliminate the $(\vphi,x)$-dependence of 
the diagonal term $b_{6}^{(k)}$ 
up to a remainder of order $- (k + 1)/2$. 

 First of all, we consider a real symbol $f_{k}(\vphi, x, \xi) \in S^{-\frac{k}{2}}$ (to be determined) 
 and we consider the flow $\Psi_{k}^\tau(\vphi)$ which is defined by 
\begin{equation}\label{def Phi2 tau ordine  - k 2}
\begin{cases}
\partial_\tau \Psi_{k}^\tau(\vphi)  = 
\ii \Pi_{S}^{\perp}\opw\Big( f_{k}(\vphi, x, \xi) \Big) \circ [\Pi_{S}^{\perp}\Psi_k^\tau(\vphi)] \\
\Psi^0_k(\vphi) = {\rm Id}\,. 
\end{cases}
\end{equation}
In this way, see \eqref{flow-jthDIAGO}, we shall write
\begin{equation}\label{graal4}
\mathcal{W}_{k}^{\tau}:=\left(\begin{matrix}\Psi_{k}^\tau & 0 \\ 0 & \overline{\Psi_{k}^\tau} \end{matrix}\right)\,,\quad \tau\in[0,1]\,.
\end{equation}
Moreover,
by Lemma \ref{differenzaFlussi} we have that 
the flow $\Psi^\tau_{k}(\vphi)$ is the correction of the flow $\Phi^\tau_{k}$ of the equation
\begin{equation}\label{tagliata2}
\partial_\tau \Phi^\tau_k(\vphi)  = \ii \opw \Big( f_k(\vphi, x, \xi) \Big) \circ \Phi_k^\tau(\vphi)\,, \quad 
\Phi^0_k(\vphi) = {\rm Id}
\end{equation}
up to a finite rank operator. More precisely,
\[
\Psi_{k}^{\tau} =\Pi_{S}^{\perp} \Phi_{k}^{\tau} \Pi_{S}^{\perp}\circ(\uno +\mathcal{R})\,,
\]
where 
%$\Phi_{k}$ is the $1$-time flow of the equation \eqref{tagliata2} and
$\mathcal{R}$ is a finite rank operator of the form \eqref{forma buona resto}
satisfying \eqref{giallo2}-\eqref{giallo3}.

We chose the symbol $f_{k}$ in such a way that 
\begin{equation}\label{omologica fk ordini bassi}
%\begin{aligned}
 \big( \omega \cdot \partial_\vphi + \mathtt{m}_{1} \cdot \nabla \big) f_k(\vphi, x, \xi) 
 + b_{6}^{(k)}(\vphi, x, \xi) = \langle b_{6}^{(k)} \rangle_\vphi (\xi)\,.
%\end{aligned}
\end{equation}
We shall use the conservation of momentum of $b_{6}^{(k)}$ and hence by applying Lemma \ref{proposizione equazione omologica generale riduzione in ordini},
one obtains that $\langle b_{6}^{(k)} \rangle_\vphi (\xi)$ is real and it satisfies
\begin{equation}\label{stima.media.b6k}
\|\langle b_{6}^{(k)} \rangle_\vphi \|_{- \frac{k}{2}, s, \alpha}^{k_0, \gamma} 
\lesssim_{s, k,\al,M} 
\e\g^{-2-k} (1 + \| {\mathcal I}_0 \|_{s_0 + \sigma_k}^{k_0, \gamma})\,, 
\quad \forall s, \alpha \geq 0
\end{equation}
and
 $f_k \in S^{- \frac{k}{2}} $ is a 
real, momentum preserving symbol satisfying for some $\s\gg 0$ large enough,
and any $(\omega, \mathtt h) \in \mathtt\Omega \times [\mathtt h_1, \mathtt h_2]$, the estimates
\begin{equation}\label{stima fk ordini bassi}
\begin{aligned}
 \| f_k \|_{- \frac{k}{2}, s, \alpha}^{k_0, \gamma} 
 & \lesssim_{s, k, \alpha,M} 
 \e \gamma^{- 2 - (k +1)} \big( 1 + \| {\mathcal I}_0 \|_{s + \sigma_k + \sigma}^{k_0, \gamma} \big)\,, 
 \quad s_0 \leq s \leq \bar{s}\,,
 \\
  \| \Delta_{12} f_k \|_{- \frac{k}{2}, p, \alpha}^{k_0, \gamma} 
 & \lesssim_{p, k, \alpha,M} 
 \e \gamma^{- 2 - (k +1)} \| i_{1} - i_{2}\|_{p + \sigma_k + \sigma}\,,
\end{aligned}
\end{equation}
for $p$ as in \eqref{ps0}. Moreover, by the smallness conditions \eqref{ansatz_I0_s0}-\eqref{smallcondepsilon} and the \eqref{stima fk ordini bassi},
then Lemmata \ref{lemma:buonaposFlussi}, \ref{differenzaFlussi} apply and one gets
\begin{equation}\label{stima mappa phi k}
\begin{aligned}
\mathfrak{M}_{( \Psi_k^{\pm 1}-{\rm Id})}(s)
%+\mathfrak{M}_{( \Psi_k^{\pm 1}-{\rm Id})^{*}}(s)
&\lesssim_{\bar{s},k,M}
\e {\gamma^{- 2 - (k +1)}} (1+\|\mathcal{I}_{0}\|_{s+\s_{k}+\s}^{k_0,\gamma})\,, 
\quad  \forall \, s_0 \leq s \leq \bar{s} \, , 
\\
\|\Delta_{12}( \Psi_k^{\pm 1}-{\rm Id})h\|_{p}
%+\|\Delta_{12}( \Psi_k^{\pm 1}-{\rm Id})^{*}h\|_{p}
&\lesssim_{\bar{s},k,M}
\e  {\gamma^{- 2 - (k +1)}} \|h\|_{p}\|i_1-i_{2}\|_{p+\s}\,,
\end{aligned}
\end{equation}
for $p$ as in \eqref{ps0}.
In view of the properties of $f_{k}$ the map $\mathcal{W}_{k}^{\tau}$ is Hamiltonian and momentum preserving. 
We then define 
\begin{equation}\label{def Omega 6 k + 1}
 \mathtt{m}_{6}^{(k + 1)} (\xi) := \mathtt{m}_{6}^{(k)}(\xi) + \langle b_{6}^{(k)} \rangle_\vphi (\xi)\,.\end{equation}
The symbol $ \mathtt{m}_{6}^{(k + 1)} (\xi)$  is real, since $\mathtt{m}_{6}^{(k)}(\xi)$ is real by hypothesis 
 and by \eqref{stima.media.b6k} it satisfies 
\begin{equation}\label{stima Omega 6 k + 1}
\begin{aligned}
 \| \mathtt{m}_{6}^{(k + 1)} - \Omega \|_{\frac12, s, \alpha} 
 &\lesssim_{s, \alpha, k,M} \e \gamma^{- 2 - k} \big( 1 + \| {\mathcal I}_0 \|_{s + \sigma_k + \sigma}^{k_0, \gamma} \big) , \quad \forall s, \alpha \geq 0\,
 \\
 \|\Delta_{12} \mathtt{m}_{6}^{(k+1)}  \|_{\frac12, p, \alpha}^{k_0, \gamma} 
 &\lesssim_{p, k,\al,M} \e \g^{-2-k}\|i_{1} - i_{2} \|_{s + \sigma_k}\,,
 \end{aligned}
\end{equation}
for some $\sigma \gg 0$ large enough.
We now compute the conjugation 
$\mathcal{W}_k^{- 1} \circ {\mathcal L}_6^{(k)} \circ \mathcal{W}_k$ 
where the operator ${\mathcal L}_6^{(k)}$ is defined in \eqref{forma L6 (k)}. By
Lemma \ref{georgiaLem} (recall Remark \ref{differenzaFlussiRMK})
we deduce that 
 the conjugated operator  has the form
 \begin{equation}\label{graal2}
\begin{aligned}
&{\bf \Pi}_{S}^{\perp}{\mathcal T}_k^{- 1} 
\Big(\omega\cdot\vphi+\mathtt{m}_{1}\cdot\nabla+\mathcal{B}_{6}^{(k)}
+  \ii {\rm Op}^W \left(\begin{matrix}\mathtt{m}_{6}^{(k)}(\xi) & 0
\\
0 & -\mathtt{m}_{6}^{(k)}(-\xi)
\end{matrix}\right)
\Big)
{\mathcal T}_k{\bf \Pi}_{S}^{\perp}
\\&\qquad\qquad
+{\bf \Pi}_{S}^{\perp}{\mathcal T}_k^{- 1} \mathcal{R}_{6}^{(k)}{\mathcal T}_k{\bf \Pi}_{S}^{\perp}
+\mathcal{Q}\,,
\end{aligned}
\end{equation}
for some finite rank operator $\mathcal{Q}$ satisfying bounds like \eqref{punkina1KKK}
and where we defined
\[
{\mathcal T}_k:=\left(\begin{matrix}\Phi_{k} & 0 \\[0.1em] 0 & \overline{\Phi_k{}} \end{matrix}\right)\,.
\]
By applying Lemma \ref{flussi coniugi eccetera}, using $({\mathcal P}1)_k, ({\mathcal P}2)_k,  ({\mathcal P}3)_k$, 
one obtains that 
\begin{equation}\label{coniugi L6 (k)}
\begin{aligned}
 \Phi_k^{- 1} \circ \Big( \omega \cdot \partial_\vphi +  \mathtt{m}_{1}  \cdot \nabla \Big) \circ \Phi_k  
&=  \omega \cdot \partial_\vphi +  \mathtt{m}_{1} \cdot \nabla 
\\&+ \ii {\rm Op}^W \Big( \big( \omega \cdot \partial_\vphi  +  \mathtt{m}_{1}  \cdot \nabla \big) f_k
 + r_1\Big) +\mathcal{Q}_{1}^{(k+1)}\,, 
 \\
 \Phi_k^{- 1} \circ \ii {\rm Op}^W\Big( \mathtt{m}_{6}^{(k)}(\xi) + b_{6}^{(k)}  \Big) \circ \Phi_k 
 &= \ii {\rm Op}^W \Big( \mathtt{m}_{6}^{(k)}(\xi) + b_{6}^{(k)}+ r_2 \Big)+\mathcal{Q}_{2}^{(k+1)}
\\
 {\mathcal T}_k^{-1}\circ \mathcal{R}_{6}^{(k)}\circ {\mathcal T}_k  & =\mathcal{Q}_{3}^{(k+1)}\,,
\end{aligned}
\end{equation}
where  $r_1\in S^{-k-1}, r_2 \in S^{- \frac{k + 1}{2}}$ are 
real, momentum preserving\footnote{this follows from the fact that $\Phi_{k}$ is symplectic, 
$ \mathtt{m}_{6}^{(k)}(\xi), b_{6}^{(k)}, \mathtt{m}_{1}$ are real and by Remark \ref{rmk:algsimboli}}, 
 and satisfy the estimates
\eqref{stimeB6kappino}
with $k\rightsquigarrow k+1$. The remainders $\mathcal{Q}_{2}^{(k+1)},\mathcal{Q}_{3}^{(k+1)}$
satisfy \eqref{punkina1KKK} thanks to \eqref{restoconiugazione4} and \eqref{stima fk ordini bassi}.   
The term $\mathcal{T}_{k}^{-1}\circ \mathcal{R}_{6}^{(k)}
\circ
{\mathcal T}_k$ satisfies \eqref{punkina1KKK} with $k\rightsquigarrow k+1$ 
by Lemma \ref{conj.NEWsmoothresto} (recalling also item $(i)$ in Remark \ref{rmk:restimatrici}) 
and the inductive assumption.

Then, by using \eqref{graal2}, \eqref{coniugi L6 (k)}, together with 
%\small{\eqref{omologica fk ordini bassi}, 
\eqref{stima fk ordini bassi}, 
one obtains that, for any $(\omega,\mathtt{h})\in \mathtt{\Omega}\times[\th_1,\th_2]$,
\begin{equation}\label{mathtt L 6 (k + 1) finale}
\begin{aligned}
\mathcal{W}_k^{- 1}& \circ {\mathcal L}_6^{(k)} \circ \mathcal{W}_k= {\bf \Pi}_{S}^{\perp}\omega \cdot \partial_\vphi 
 \\&+
 {\bf \Pi}_{S}^{\perp}\Big(
  \mathtt{m}_{1} \cdot \nabla
  + \ii {\rm Op}^W \left(\begin{matrix} \mathtt{q}^{(k + 1)}(\xi) & 0
  \\ 0 & -\mathtt{q}^{(k + 1)}(-\xi)
  \end{matrix}\right) 
  +\mathcal{B}_{6}^{(k+1)}+
  \mathcal{R}_{6}^{(k+1)}
 \Big)\Pi_{S}^{\perp}\,,  
\end{aligned}
\end{equation}
where 
\begin{equation}\label{quetto}
\tq^{(k + 1)} := \mathtt{m}_{6}^{(k)}(\xi)+  \big( \omega \cdot \partial_\vphi  +  \mathtt{m}_{1}  \cdot \nabla \big) f_k +  b_{6}^{(k)}\,,
\end{equation}
and where $  \mathcal{R}_{6}^{(k+1)}$ is some remainder satisfying \eqref{punkina1KKK} at $k+1$.
In order to conclude the proof, we define the operator 
 (defined for all values $(\omega,\th)\in \mathtt{\Omega}\times[\th_1,\th_2]$)
 \begin{equation}\label{mathtt L 6 (k + 1) finale2}
 {\mathcal L}_6^{(k+1)} = {\bf \Pi}_{S}^{\perp}\omega \cdot \partial_\vphi + 
{\bf \Pi}_{S}^{\perp}\Big( 
\ii\opw(\mathtt{m}_{1}  \cdot \xi) \uno  
+ \ii E\opw(\mathtt{m}_{6}^{(k+1)})+\cB_{6}^{(k+1)} + \cR_{6}^{(k+1)} \Big) {\bf \Pi}_{S}^{\perp}\,,
 \end{equation}
where $\cB_{6}^{(k+1)}$ has the form \eqref{Zippoblu} with $k\rightsquigarrow k+1$ and 
$b_{6}^{(k + 1)} := r_1 + r_2$.

In conclusion, by \eqref{mathtt L 6 (k + 1) finale} and \eqref{mathtt L 6 (k + 1) finale2}, using
 the homological equation \eqref{omologica fk ordini bassi} and \eqref{quetto}, recalling \eqref{def Omega 6 k + 1}, 
we deduce that for any 
$(\omega, \th)\in \tT\tC_{\infty}(\gamma,\tau) $,
 than one has that
$\tq^{(k + 1)}(\vphi, x, \xi) = \mathtt{m}_{6}^{(k+1)}(\xi)$.
This implies the 
 the conjugation result in \eqref{forma L6 (k)} with $k\rightsquigarrow k+1$.
%where $\mathtt{m}_{6}^{(k + 1)}(\xi)$ is in \eqref{def Omega 6 k + 1} and satisfies \eqref{stima Omega 6 k + 1},
%the estimates \eqref{stima mappa phi k} on the map
%and Lemmata \ref{composizione operatori tame AB}-\ref{lemma operatore e funzioni dipendenti da parametro};
%The term
% $  \mathcal{R}_{6}^{(k+1)}$ is some remainder satisfying 
%\eqref{punkina1KKK} at $k+1$ and where we set
%$b_{6}^{(k + 1)} := r_1 + r_2$.
The discussion above implies $({\mathcal P}j)_k$, $j=1,\ldots,5$.

\bigskip

\noindent
{\sc Proof of Proposition  \ref{riduzione diagonal ordini bassi} concluded.} 

\noindent
We perform $N= 2(M - \tc)$ steps of the previous precedure
%Let $k \in \{ 1, \ldots, 2N  \}$ and let us 
and we define (recall \eqref{def Phi2 tau ordine  - k 2}-\eqref{graal4})
$\cW:=\cW_0\circ\cdots\circ\cW_{N}$.
By composition and using \eqref{stima mappa phi k}, one 
gets \eqref{stimamappaPsi2}
for $\sigma \gg \sigma_{N}$ large enough. 
Moreover, the map $\cW$ is 
real-to-real,
Hamiltonian and, since every $\Psi_{k}$ is a flow generated by a momentum preserving symbol, it is momentum preserving.

By iterating the induction step described above, starting from $\mathcal{L}_{6}^{(0)}\equiv\mathcal{L}_{6}$, 
one then obtain that 
${\mathcal L}_7 \equiv {\mathcal L}_6^{(N)} = \cW^{- 1} \circ {\mathcal L}_6 \circ \cW$ is of the form \eqref{op cal L7} with $\mathtt{m}_7 := \mathtt{m}_{6}^{(N)}$ and 
$$
{\mathcal R}_7 := \cW^{- 1} \circ {\mathcal R}_6 \circ \cW+ \ii{\rm Op}^W \begin{pmatrix}
b_{6}^{(N)}(\vphi, x, \xi) & 0 \\
0 & -  b_{6}^{(N)}(\vphi, x, - \xi)
\end{pmatrix}
$$
The estimates \eqref{stima tame cal R7 new} on ${\mathcal R}_7$ follows by \eqref{punkina1KKK} with $k=N$,
and by \eqref{stimeB6kappino} with $k=N$, together with Lemma \ref{constantitamesimbolo}. 
Finally, it is Hamiltonian and momentum preserving by difference. 
This concludes the proof.
\end{proof}

\medskip
In view of 
Proposition \ref{lemma:7.4}, \ref{prop operatore cal L2}, \ref{riduzione trasporto}, \ref{prop operatore cal L3},
\ref{blockTotale}, \ref{riduzione diagonal ordine 1 2} and \ref{riduzione diagonal ordini bassi}
we have that, for  any $\lambda=(\omega , \th)\in \tT\tC_{\infty}(2\gamma,\tau)$
(see \eqref{tDtCn})
the real, Hamiltonian and  momentum preserving operator $\cL_{\oo}$ in \eqref{cal L omegaWW} is
conjugated, 
under the  real, symplectic and momentum preserving map
\begin{equation}\label{def.calkappa}
\cK:= \cG_{B}{\mathcal C}^{- 1} \cE_{\perp} \Phi_{\mathbb{F}} \cS {\bf \Psi} \cW % \cW {\bf \Psi} \cS \Phi_{\mathbb{F}} \cE_{\perp} \cC^{-1} \cG_{B}
\end{equation}
to the real, Hamiltonian and momentum preserving operator $\cL_{7}$ defined in \eqref{op cal L7}.
The next theorem summarizes the main result of sections \ref{linearizzato siti normali}, \ref{sym.low.order}
and of Propositions \ref{riduzione diagonal ordine 1 2}-\ref{riduzione diagonal ordini bassi}.

\begin{thm}\label{red.Lomega.smooth.rem}
{\bf (Reduction of $\cL_\oo$ up to smoothing remainders)} 
Fix 
$\alpha_{*}\in \N$,
 $\beta_{0} \in \N$, $M \in \N$, $\mathtt c \geq 3(k_0 + \beta_0 ) + 3$ with $M \geq \mathtt c$. 
There exists $\s=\s(M, \beta_{0}, k_0, \alpha_{*})$
such that for any $\bar{s}\geq s_0$ there is $\delta(\bar{s},M,\beta_0,\alpha_*)$ such that
if \eqref{smallcondepsilon} holds with $\kappa=2(M-\mathtt{c})+3$ and \eqref{ansatz_I0_s0} holds with 
$\mu_0\geq\s$
then the following holds.
For all $ \lambda = (\oo, \th) \in \tT\tC_{\infty}(\gamma,\tau)  $ (see \eqref{tDtCn}),
the operator $ \cL_\oo $ in \eqref{representation Lom} is conjugated 
by the map $\cK $ defined in \eqref{def.calkappa} 
 to the real, Hamiltonian and momentum preserving operator
$\cL_7 $ of the form
\begin{equation}\label{operatoreintero}
\cL_7 = \mathcal{K}^{-1}\mathcal{L}_{\omega}\mathcal{K}
=  \omega \cdot \partial_\vphi {\mathbb I}_{\perp}+ \ii \cD_\bot + \cR_\bot 
\end{equation}
where ${\mathbb I}_{\perp}$ denotes the identity map of ${\bf H}_{S}^{\perp}$ 
(see \eqref{decomp siti tangenziali coordinate complesse}) and:

\noindent
$(i)$ {\bf (Normal form).} The operator 
$ \cD_\bot $ is diagonal and has  the form 
\begin{equation}\label{proprieta cal D botAA}
\begin{aligned}
{\mathcal D}_{\perp}& := \left(\begin{matrix}
\mathtt{D}_{\perp} & 0
\\  0 &-\overline{\mathtt{D}}_{\perp}
\end{matrix}
\right)\,,
%E {\rm diag}_{j \in S_0^c} \mu_0 (j)\,, 
\qquad 
\begin{aligned}
\mathtt{D}_{\perp}&:={\rm diag}_{j \in S_0^c} \mu_0 (j)\,,
% \mathtt{m}_{1} \cdot j + \mathtt{m}_{7}(j)\,, 
\end{aligned}
\end{aligned}
\end{equation}
where the eigenvalues $\mu_0 (j)$, $j\in S_0^{c}$, are real valued, 
are defined for any $\lambda = (\omega, \mathtt h) \in \mathtt \Omega \times [\mathtt h_1, \mathtt h_2]$
and admit the expansion
\begin{equation}\label{proprieta cal D bot}
\mu_0 (j) :=  \Omega(j) + \mathtt m_{1} \cdot j + \mathfrak{r}_0(j)\,,
\qquad
\mathfrak{r}_0(j):=\mathtt{m}_{7}(j)-\Omega(j)\,,
\end{equation}
where $\Omega(j)\equiv\Omega(j;\mathtt{h})$ is in \eqref{wild}, $\mathtt{m}_{1} $  in \eqref{beta.FGMP.est}, 
the symbol 
$\mathtt{m}_{7}\in S^{\frac{1}{2}}$ satisfies \eqref{stima lambda 7} for any $0\leq \alpha\leq\alpha_{*}$.
Moreover, for any $\lambda = (\omega, \mathtt h) \in \mathtt \Omega \times [\mathtt h_1, \mathtt h_2]$, 
the following estimates hold:

\begin{itemize}
\item[a)]
for any $j \in S_0^c$ one has
\begin{align}
|\partial_\lambda^k (\mu_0 - \Omega)(j; \lambda) | &\lesssim_{M}  
\e \gamma^{-2(M-\mathtt{c})-1 - |k|} |j|, \quad \forall k \in \N^{\nu + 1}, \quad |k| \leq k_0\,,
%\quad \forall \lambda \in \R^\nu \times [\mathtt h_1, \mathtt h_2]\,.  
\label{stima partial mu 0 pre rid}
\\
|\Delta_{12} \mu_0(j) | &\lesssim_{M} 
\e \gamma^{-2(M-\mathtt{c})-1} |j| \| i_1 - i_2 \|_{s_0 + \s}\,.
\label{stima delta 12 mu 0 pre rid}
\end{align}

\item[b)]
for any $j, j' \in S_0^c$, $j \neq j'$, and for all $|k| \leq k_0$  one has 
\begin{equation}\label{stima tipo simbolo mu 0 j j'}
\begin{aligned}
& |\partial_{\mathtt h}^k (\Omega(j) - \Omega(j'))| \lesssim_{M} |j - j'|\,,  
\\
& |\pa_{\lambda}^{k}(\mu_0 - \Omega)(j) - \pa_{\lambda}^{k}(\mu_0 - \Omega)(j'))|\,, 
\,|\pa_{\lambda}^{k}(\mathfrak{r}_0(j) - \mathfrak{r}_0(j'))|
\lesssim_{M}  \e \gamma^{-2(M-\mathtt{c})-1 - |k|}  |j - j'|\,.
\end{aligned}
\end{equation}
As a consequence, for any $j, j' \in S_0^c$, $j \neq j'$, $|k| \leq k_0$, one has 
\begin{equation}\label{stima tipo simbolo mu 0 j j'a}
|\pa_{\lambda}^{k}(\mathtt m_7(j) - \mathtt m_7(j'))|\,,\, |\pa_{\lambda}^{k}(\mu_0(j) - \mu_0(j'))|\lesssim  |j - j'|\,. 
\end{equation}

\end{itemize}

\smallskip
\noindent
$(ii)$ {\bf (Remainder).} The smoothing remainder 
${\mathcal R}_{\perp}$ is a Hamiltonian, momentum preserving operator, it is defined for any $\lambda = (\omega , \mathtt h) \in \mathtt \Omega \times [\mathtt h_1, \mathtt h_2]$ and satisfies the following:
  %for any $\tm \in \N$ such that $\tm \leq N - 3(\tb + k_{0})-3$, 
  for any $\beta \in \N^\nu$, $|\beta| \leq \beta_0$, 
  the linear operator $\langle D \rangle^M \partial_\vphi^\beta {\mathcal R}_{\perp}\jap{D}^{-\mathtt{c}}$ 
  is ${\mathcal D}^{k_0}$-tame with tame constant satisfying
 \begin{align}
{\mathfrak M}_{ \langle D \rangle^{M} \partial_\vphi^\beta{\mathcal R}_\bot 
\langle D \rangle^{-\mathtt{c}}}(s) 
&\lesssim_{\bar{s},M} 
\e {\gamma^{-2(M-\mathtt{c})-2}}(1 + \| {\mathcal I}_0 \|_{s + \s}^{k_0, \gamma})\,, 
\quad \forall s_0 \leq s \leq \bar{s}\,,
\label{tame riducibilita cal R bot iniziale}
\\
\| \langle D \rangle^{M}\partial_\vphi^\beta \Delta_{12} \cR_\bot \langle D \rangle^{-\mathtt{c}} \|_{\cL(H^{p})}
&\lesssim_{\bar{s},M} \e \gamma^{- 2(M-\mathtt{c})-2}  \| i_1-i_2 \|_{p + \s}\,,
\label{stima tame derivata cal R bot}
\end{align}
for any $p$ as in \eqref{ps0}.

\smallskip
\noindent
$(iii)$ {\bf (Transformation).} The map $\mathcal{K}$ in \eqref{def.calkappa} is invertible and, for any $s_0\leq s\leq \bar{s}$, one has 
\begin{equation}\label{stimamappacomplessiva}
\begin{aligned}
\|\mathcal{K}^{\pm1}h\|_{s}^{k_0,\gamma}&\lesssim_{\bar{s},M} \|h\|_{s+\s}^{k_0,\gamma}
+ \| {\mathcal I}_0 \|_{s + \s}^{k_0, \gamma}\|h\|_{s_0+\s}^{k_0,\gamma}\,,
\\
\|(\Delta_{12}\mathcal{K}^{\pm1})h\|_{p}&\lesssim_{\bar{s},M}\e\gamma^{-2(M-\mathtt{c})-2}\|i_1-i_2\|_{p+\s}\|h\|_{p+\s}\,,
\end{aligned}
\end{equation}
for any $p$ as in \eqref{ps0}.
\end{thm}

\begin{proof}
Formula \eqref{operatoreintero} with the expansion 
\eqref{proprieta cal D botAA}-\eqref{proprieta cal D botAA} follow by \eqref{op cal L7}
 by setting
 \[
 \begin{aligned}
 \ii \cD_\bot&:=\mathtt{m}_{1}  \cdot \nabla \uno  +  \ii \opw\left(
 \begin{matrix}
 \mathtt{m}_{7}(\x) & 0 \\0 & -\mathtt{m}_{7}(-\x)\end{matrix}
 \right)\,,
 \qquad
  \mathcal{R}_{\perp}:={\bf \Pi}_{S}^{\perp}\mathcal{R}_{7}{\bf \Pi}_{S}^{\perp}\,.
 \end{aligned}
 \]
 The estimates \eqref{tame riducibilita cal R bot iniziale}-\eqref{stima tame derivata cal R bot} 
 on $\mathcal{R}_{\perp}$
 follow 
 %respectively by \eqref{stima lambda 7} and 
 \eqref{stima tame cal R7 new}, while 
 the estimates on the map $\mathcal{K}$ in \eqref{stimamappacomplessiva}
 follow by composition, recalling Lemma \ref{lemma operatore e funzioni dipendenti da parametro},
the bounds \eqref{cost tame good}-\eqref{delta12 good}, \eqref{stimamappaEperp}, \eqref{AVB3mappa2PRO},
 \eqref{stimamappaS}, \eqref{stimamappaPsi}, \eqref{stimamappaPsi2}
 and the smallness condition on $\e$ in \eqref{smallcondepsilon}.
 This proves items $(ii)$ and $(iii)$.
 It remains to prove  item $(i)$.
 
 \smallskip
 \noindent
 {\sc Proof of \eqref{stima partial mu 0 pre rid}.}
In view of \eqref{proprieta cal D bot}, one has to estimate both 
$\mathtt{m}_{1} (\lambda) \cdot j$ and $\tm_{7}(j)$, $j\in S_0^c$. 
By the estimate \eqref{beta.FGMP.est}, one has that  for any $k \in \N^{\nu + 1}$, $|k| \leq k_0$, 
\[
|\partial_\lambda^k \mathtt{m}_{1}(\lambda) \cdot j| 
\leq 
\gamma^{- |k|}|\mathtt{m}_{1}|^{k_0, \gamma} |j| \lesssim \e \gamma^{- |k|} |j|\,. 
% \lesssim |j|\,. 
\]
We now estimate $(\tm_{7} - \Omega)$ (recall \eqref{dispersionLaw}). 
%A direct calculation shows that for any $n \in \N$, $\mathtt h \in [\mathtt h_1, \mathtt h_2]$, 
%$|\partial_{\mathtt h}^n \Omega(j; \mathtt h)| \lesssim_n |j|^{\frac12}\,. $
Moreover, for any $k \in \N^{\nu + 1}$, $|k| \leq k_0$, 
$\lambda \in {\mathtt\Lambda}_{0}\subseteq \R^\nu \times [\mathtt h_1, \mathtt h_2]$, 
one has that 
\[
\begin{aligned}
|\partial_\lambda^k (\tm_7(j; \lambda) - \Omega(j; \lambda))| 
& = 
\gamma^{- |k|} |j|^{\frac12} \Big( \gamma^{|k|} |\partial_\lambda^k (\tm_7(j; \lambda) - \Omega(j; \lambda))| 
|j|^{- \frac12} \Big) 
\\& 
\leq 
\gamma^{- |k|} |j|^{\frac12} \sup_{|k| \leq k_0} \sup_{\lambda \in \mathtt \Omega \times [\mathtt h_1, \mathtt h_2]} 
\sup_{ \xi \in \R^d}  \gamma^{|k|} |\partial_\lambda^k (\tm_7(\xi; \lambda) - \Omega(\xi; \lambda))| 
\langle \xi \rangle^{- \frac12} 
\\& 
\leq 
\gamma^{- |k|} |j|^{\frac12} \| \tm_7 - \Omega \|^{k_{0}, \gamma}_{\frac12, s_0, 0} 
\stackrel{\eqref{stima lambda 7}, \eqref{ansatz_I0_s0} }{\lesssim_{M}} 
\e \gamma^{- |k|-2(M-\mathtt{c})-1} |j|^{\frac12} \,.
%\lesssim |j|^{\frac12}\,,
\end{aligned}
\]
%provided that $\e \gamma^{- k_0- 2(\mathtt{m}-\mathtt{c}_1)-1} \ll 1$. 
%The claimed estimate has then been proved. 
The discussion above implies the bound \eqref{stima partial mu 0 pre rid}.

\smallskip
\noindent
{\sc Proof of \eqref{stima tipo simbolo mu 0 j j'}, \eqref{stima tipo simbolo mu 0 j j'a}.} The first estimate in \eqref{stima tipo simbolo mu 0 j j'} follows by a direct calculation by using the explicit expression of 
$\Omega(j)$  in \eqref{wild}. 
Let us prove the other estimates. One has that for any $j, j' \in S_0^c$
\[
\mathfrak{r}_0(j) - \mathfrak{r}_0(j') = (\tm_7 - \Omega)(j) - (\tm_7 - \Omega)(j')
\]
Since, by \eqref{stima lambda 7}, $\partial_\lambda^k \mathfrak{r}_0$, $|k| \leq k_0$, is a Fourier multiplier in $S^{1/2} $, one has that 
$\partial_\lambda^k\pa_\xi\mathfrak{r}_0 \in S^{- \frac12} \subset S^0$. Moreover, using   the bound 
\eqref{stima lambda 7}, one gets that for any $|k| \leq k_0$,
we get
\begin{equation}\label{sistemi dinamici merdaccia}
\sup_{\xi \in \R^d}|\partial_\lambda^k\pa_\xi \mathfrak{r}_0 (\xi)| 
\leq C\e \gamma^{- |k|-2(M-\mathtt{c})-1}
\end{equation}
for some constant $C > 0$. By using the 
mean value theorem, one then has that for any $j, j' \in S_0^c$, $j \neq j'$, $|k| \leq k_0$ 
\[
\begin{aligned}
|\partial_\lambda^k(\mu_0 - \Omega)(j) - \partial_\lambda^k(\mu_0 - \Omega)(j'))| 
& \leq |\partial_\lambda^k\mathtt{m}_{1}| | |j - j'| 
+ | \partial_\lambda^k \mathfrak{r}_0(j) - \partial_\lambda^k \mathfrak{r}_0(j')|  
\\
&  \stackrel{\eqref{beta.FGMP.est}, \eqref{sistemi dinamici merdaccia}}{\lesssim_{M}} 
\big(\e \gamma^{- |k|} + \e \gamma^{- |k|-2(\mathtt{m}-\mathtt{c}_1)-1}\big) |j - j'| 
\\
& \lesssim_{M} \e \gamma^{- |k|-2(\mathtt{m}-\mathtt{c}_1)-1} |j - j'|\,. 
\end{aligned}
\]
This latter chain of inequalities proves the second estimate in \eqref{stima tipo simbolo mu 0 j j'}. 
The estimate \eqref{stima tipo simbolo mu 0 j j'} (using that 
$\e \gamma^{- k_0-2(M-\mathtt{c})-1} \ll 1 $) trivially implies \eqref{stima tipo simbolo mu 0 j j'a}.

\smallskip
\noindent
{\sc Proof of \eqref{stima delta 12 mu 0 pre rid}} It follow as done above and by using \eqref{beta.FGMP.est} 
and \eqref{stima lambda 7}.
\end{proof}

\section{Perturbative reduction scheme}\label{sec:KAMreducibility}

%\subsection{KAM Reducibility}\label{sec:KAMreducibility}

In this section we diagonalize the operator 
\begin{equation}\label{operatorenonintero}
\bL_{0}:=\cL_7=\oo \cdot \pa_{\vphi} + \ii \cD_{\perp}+ {\mathcal R}_{\perp}\,,
\end{equation}
obtained by Theorem \ref{red.Lomega.smooth.rem}.
 We diagonalize the operator $\bL_{0}$ through a KAM reducibility scheme. 
We want to implement a reducibility scheme in which we deal with 
${\mathcal D}^{k_0}$-modulo tame operators (see Def. \ref{Dk0modulotame}). 
In order to do this we need the following auxiliary  lemmata.  

Given $\tau > 0$\footnote{the parameter $\tau$ will be fixed in the following, see \eqref{param_small_meas},
depending only on the dimensions $\nu,d$ and on the non-degeneracy index $m_0$ given by Proposition \ref{lem:transversality}.}
we fix the following constants:
\begin{equation}\label{costanti.riducibilita}
\begin{aligned}
&\tau_3 := 2(k_0 + \tau(k_0 + 1))\,, &&\ta_{1} := 3 \tau_{3}^2 + \chi (2 (\mathtt m - \mathtt c_1) + 2(\tau + 1)(\tau + 2) + 4)\,,  \\
&  \tb_{1} := 3\big(\lfloor \ta_{1} \rfloor + 2\big),  &&\beta_0:=\mathtt{b}_1+ \nu + d + 2\,,
\\
&  \mathtt{c}_1:=3(\beta_0+k_0)+4\,,
 &&\mathtt{m}:= \mathtt{c}_1+\tau_3+\tau\,.
%\\&
% \mathtt m \gg \mathtt b_{1} + \tk_0 + 2 + \tau_{3} + \tau\,,
%\qquad \mathtt k_0 := 2s_0 +  3k_0 + 4\,, 
\end{aligned}
\end{equation}
%and we set
%\begin{equation}\label{Scelta.N}
%N \ge \tb_1+ \tk_{0} + \tm \in \N \,  .
%\end{equation}
We define
 \begin{equation}\label{definizione.mu.b}
 \mu (\mathtt{b}_{1}) := % \s(\mathtt{m}, 2s_0 + \mathtt b_{1},\nu,d)\,, 
 \s(\mathtt{m}, \beta_0, \alpha_* , \ta_{1})\,,
\end{equation}
where the constant $\s=\s(\mathtt{m}, \beta_0, \alpha_{*}, \ta_{1})$ 
is given by Theorem \ref{red.Lomega.smooth.rem}
applied with $M\rightsquigarrow\mathtt{m}$, 
$\mathtt{c}\rightsquigarrow \mathtt{c}_1$, $\ta \rightsquigarrow \mathtt{a}_1$,
$\beta_0$ in \eqref{costanti.riducibilita}
and $\alpha_{*}=1$.
The condition \eqref{costanti.riducibilita} 
means that one has to perform a 
sufficiently large number of regularizing 
steps in Sections \ref{sym.low.order} and \ref{sec:redulower}.
As a consequence of the estimates of 
Theorem \ref{red.Lomega.smooth.rem}
in the next Lemma 
we show that 
$ \cR_\bot$ satisfies 
suitable \emph{modulo}-tame estimates.

\begin{lemma}\label{inizializzazione modulo tame red}
Let $\mathtt{b}_{1}$, $\mathtt m$, $\mathtt{c}_1$ as in \eqref{costanti.riducibilita}.
Then the operators 
$\langle D \rangle^{\mathtt m} \langle \pa_{\vphi}\rangle^{\mathfrak{d}}{\mathcal R}_\bot 
 \langle D \rangle^{- \mathtt{c}_{1} } $, with $\mathfrak{d}=0,\tb_1$ 
 are ${\mathcal D}^{k_0}$-modulo tame operators whose tame constants 
 satisfy the following estimates: 
\begin{equation}\label{stime mod tame primo resto rid}
\begin{aligned}
{\mathfrak M}^\sharp_{\langle D \rangle^{\mathtt m} \langle \partial_\vphi \rangle^{\mathfrak{d}} 
{\mathcal R}_\bot \langle D \rangle^{- \mathtt{c}_{1}}}(s)
 &\lesssim_{\bar{s}} \e {\gamma^{-2-2(\mathtt{m}-\mathtt{c}_1) }} 
 \big(1 + \| {\mathcal I}_{0} \|_{s +  \mu (\mathtt{b}_{1})}^{k_0, \gamma} \big)\,, 
\end{aligned}
\end{equation}
for $\mathfrak{d}=0,\tb_{1}$ and $s_0 \leq s \leq \bar{s}$.
\end{lemma}
\begin{proof}
In view of the choices \eqref{costanti.riducibilita} we shall apply 
Lemma \ref{lem: Initialization astratto}. Then
the estimate \eqref{stime mod tame primo resto rid} follows recalling 
\eqref{tame riducibilita cal R bot iniziale}.
\end{proof}
The main result of this section is the following.

 \begin{thm}\label{Teorema di riducibilita} {\bf (KAM reducibility)}
Assume that $\| {\mathcal I}_{0} \|_{s_0 + \mu_0}^{k_0, \gamma} \leq 1$ with 
$ \mu_0 \gg s_0 +   \mu (\mathtt{b}_{1})$. 
Then for any $\bar{s} > s_0$ there exist $K_0 := K_0(\bar{s}) > 0$, 
$ 0 < \delta_0 := \delta_0(\bar{s})  \ll 1 $,  so that if 
\begin{equation}\label{KAM smallness condition1}
K_0^{\tau_3} \e \gamma^{-k_0- 2(\mathtt{m}-\mathtt{c}_1) - 3} \leq \delta_0\,, 
%\qquad N_0 \gamma \geq \mathtt C_0(\bar s)
\end{equation}
for $\tau_3$ as in \eqref{costanti.riducibilita}, 
%\eqref{KAM smallness condition1} is fullfilled, 
then the following holds. 

\smallskip
\noindent
$(i)$ {\bf (Eigenvalues).} For any 
$ (\omega, \mathtt h) $ in $\mathtt \Omega \times [\mathtt h_1, \mathtt h_2]$
there exists a sequence\footnote{Whenever it is not strictly necessary, we shall drop the 
dependence on $\lambda$ and $i_0$.}
\begin{equation}\label{autovalMuinfinity}
\begin{aligned}
\mu_\infty(j)&= \mu_\infty(j;\,\cdot\,): \mathtt \Omega \times [\mathtt h_1, \mathtt h_2]\to \R \,,
\\
 \mu_\infty(j)  &\equiv \mu_\infty(j; \lambda, i_0(\lambda)) 
= \mu_0(j; \lambda , i_0(\lambda)) + r_\infty(j; \lambda, i_0(\lambda))\,, \quad j\in S_0^{c}\,,
\end{aligned}
\end{equation}
where $\mu_0(j)$ are in \eqref{proprieta cal D bot} and 
\begin{align}
 | r_\infty(j)|^{k_0, \gamma} 
& \lesssim_{\bar{s}}   \e {\gamma^{-2-2(\mathtt{m}-\mathtt{c}_1)}}  |j|^{- {\mathtt m} + \mathtt{c}_{1} }\,,
\label{stimaAutovaloriFinalissimi1}
\\
 | \Delta_{12} r_\infty(j)| & \lesssim_{\bar{s}}   \e {\gamma^{-2-2(\mathtt{m}-\mathtt{c}_1)}}  
 |j|^{- {\mathtt m} + \mathtt{c}_{1} } \| i_1 - i_2 \|_{s_0 + \mu(\mathtt b_1)}\,. 
 \label{stimaAutovaloriFinalissimi2}
\end{align}

\smallskip
\noindent
$(ii)$ {\bf (Conjugacy).}
There exists
 a real-to-real, symplectic, momentum preserving, invertible transformation 
 $\mathcal{U}_{\infty}: {\bf H}^{s}_{{\perp}}\to {\bf H}^{s}_{{\perp}}$
defined for all $ (\omega, \mathtt h) $ in $\mathtt \Omega \times [\mathtt h_1, \mathtt h_2]$, 
such that for all
$ \lambda = (\omega, \mathtt h) $ belonging to the set (recall \eqref{def:DCgt}-\eqref{tDtCn})
 \begin{equation}\label{Cantor set}
 	\begin{aligned} 
 \Lambda_\infty^{\gamma}\equiv \Lambda_\infty^{\gamma}(i_0) 
 %\equiv \Omega_\infty^\gamma(v) 
 	:= \Big\{ \lambda = &(\omega, \mathtt h) \in \tD\tC (2 \gamma,\tau) 
	\cap {\mathtt T}{\mathtt C}_\infty( \gamma, \tau)(i_0) \, : \,%{\mathcal O}_\infty^\gamma 
 		\\&  |\omega \cdot \ell + \mu_\infty(j; \lambda, i_0(\lambda)) \pm  \mu_\infty(j'; \lambda, i_0(\lambda))| \geq \frac{2  \,\gamma}{\langle \ell \rangle^\tau | j' |^\tau } \,, \\
 		&  \forall \ell \in \Z^\nu \setminus \{ 0 \}\,, \, 	j,j' \in S_0^c\,, \ \  {\mathtt V}^T\ell  + j \pm j' =0\Big\}\,. 
 	\end{aligned}
 \end{equation}

the following conjugation formula holds (recall \eqref{operatorenonintero}):
\begin{equation}\label{cal L infinito}
\begin{aligned}
{\mathbf L}_{\infty} &:= {\mathcal U}_{\infty}^{- 1} {\mathbf L}_0 {\mathcal U}_{\infty}= 
\omega \cdot \partial_\vphi  + \ii {\mathcal D}_{\infty} \,,
\\ {\mathcal D}_{\infty} & := \begin{pmatrix}
\mathtt D_{\infty} & 0 \\
0 & - \overline{\mathtt D}_{\infty}
\end{pmatrix}\,,\quad \mathtt D_{\infty}  := {\rm diag}_{j \in S_0^c} \mu_\infty (j)\,.
\end{aligned}
\end{equation}
%where $ \mu_\infty (j)$ are given by Lemma \ref{lemma blocchi finali}.
Moreover, the operators
$ {\mathcal U}_{\infty}^{\pm 1} -  {\rm Id}_\bot  $ are $ {\mathcal D}^{k_0}$-modulo-tame 
with a modulo-tame constant satisfying
\begin{equation}\label{stima Phi infinito}
{\mathfrak M}^\sharp_{{\mathcal U}_{\infty}^{\pm 1} -  {\rm Id}} (s) 
\lesssim_{\bar{s}} 
K_0^{\tau_3} \e  {\gamma^{-3-2(\mathtt{m}-\mathtt{c}_1)}} (1  + 
\| {\mathcal I}_{0} \|_{s  + \mu (\mathtt{b}_{1}) }^{k_0, \gamma})\,, \quad \forall s_0 \leq s \leq \bar{s} \, .
\end{equation} 
\end{thm}

The rest of the section is devoted to the proof of the Theorem above.

\subsection{Iterative scheme}

In view of Lemma \ref{inizializzazione modulo tame red} we now 
 diagonalize the operator $\bL_0$ in \eqref{operatorenonintero}
by means of a perturbative scheme adapted to $\mathcal{D}^{k_0}$-modulo-tame operators.
We perform the reducibility of $\bL_{0}$ along the scale
\begin{equation}\label{EnneEnne}
 K_{- 1} := 1, \quad K_n := K_0^{\chi^n}\,, 
 \quad K_0>0\,,
 %\quad \chi \in (1, 2)\,, 
 \quad \chi:=\frac32\,,
\end{equation}
for some $K_0>0$ to be chosen later (large enough and independent of the diophantine constant $\gamma$), 
requiring inductively at each step 
the non-resonance conditions in \eqref{Omega nu + 1 gamma}.
We shall prove the following.
\begin{prop}\label{iterazione riducibilita}
{\bf (Iterative lemma).} 
Let $\bar s \ge s_{0} + \mu(\mathtt{b}_{1})$ (recall \eqref{costanti.riducibilita}, \eqref{definizione.mu.b}).
%and fix $N$ as in \eqref{Scelta.N}.
There exist 
$ K_0 := K_0 (\bar{s})\,,\, \mathtt C_0(\bar s) \gg 0$, $\delta_0 := \delta_0(\bar{s}) \ll 1$ 
such that, if \eqref{ansatz_I0_s0} holds with $\mu_{0}= \mu(\mathtt{b}_{1})$ 
and \eqref{KAM smallness condition1} is fulfilled
%\begin{equation}\label{KAM smallness condition1}
%K_0^{\tau_3} \e \gamma^{-k_0- 2(\mathtt{m}-\mathtt{c}_1) - 3} \leq \delta_0\,, 
%%\qquad N_0 \gamma \geq \mathtt C_0(\bar s)
%\end{equation}
%for $\tau_3$ as in \eqref{costanti.riducibilita} 
then, for all $ n \in \N $, %$n = 0, 1 , \ldots, \overline{\tn}$: 
\begin{itemize}
\item[${\bf(S1)_{n}}$] 
There exists a real operator
\begin{equation}\label{cal L nu}
\bL_{n} := \omega \cdot \partial_\vphi + \ii {\mathcal D}_n  + {\mathcal R}_n \,, \quad 
{\mathcal D}_n  := \begin{pmatrix}
\mathtt D_n & 0 \\
0 & - \overline{\mathtt D}_n 
\end{pmatrix}\,,\quad \mathtt D_n  := {\rm diag}_{j \in S_0^c} \mu_n (j)\,,
\end{equation}
defined for all $\lambda=  (\omega, \mathtt h) \in \mathtt \Omega \times [\mathtt h_1, \mathtt h_2]$ 
where $ \mu_n(j ) $ 
are $ k_0 $ times differentiable functions of the form 
\begin{equation}\label{mu j nu}
\mu_n (j; \lambda) := \mu_0 (j; \lambda ) + r_n (j ; \lambda )
\in \R
\end{equation}
where $\mu_0(j)$ are defined in \eqref{proprieta cal D bot}, 
and $r_n(j)$ satisfies the estimate
\begin{equation}\label{stima rj nu}
 \quad |r_n(j)|^{k_0, \gamma} \lesssim_{\bar{s}} 
 \e \gamma^{-2-2(\mathtt{m}-\mathtt{c}_1)}  |j|^{- \mathtt m + \mathtt{c}_{1} }\,,\ \  \forall j \in S_0^c
\end{equation}
and, for $n \geq 1$,
\begin{equation}\label{vicinanza autovalori estesi}
|\mu_n(j) - \mu_{n - 1}(j) |^{k_0, \gamma} = |r_n(j) - r_{n - 1}(j) |^{k_0, \gamma} 
\lesssim_{\bar{s}}
K_{n - 2}^{- {\mathtt a}_{1}}
\e {\gamma^{-2-2(\mathtt{m}-\mathtt{c}_1)}} 
 |j|^{- {\mathtt m} + \mathtt{c}_{1} } \,. 
\end{equation}
The remainder
\begin{equation}\label{forma cal R nu}
{\mathcal R}_n := \begin{pmatrix}
\mathtt R_{n, 1} & \mathtt R_{n, 2} \\
\overline{\mathtt R}_{n, 2} & \overline{\mathtt R}_{n, 1}
\end{pmatrix}
\end{equation}
is a Hamiltonian, momentum preserving and satisfies the following properties.

\noindent
The operators $\langle D \rangle^\tm {\mathcal R}_n \langle D \rangle^{- \mathtt{c}_1} $ 
and $ \langle D \rangle^{\mathtt m}  \langle \partial_{\vphi } \rangle^{\mathtt{b}_{1}}  
{\mathcal R}_n \langle D \rangle^{- \mathtt{c}_1}$ are 
$ {\mathcal D}^{k_0} $-modulo-tame 
and there exists a constant $C_* = C_*(\bar{s})  \gg  0$ large enough such that, for any $s \in [s_0, \bar{s}] $, 
\begin{equation}\label{stima cal R nu}
\begin{aligned}
 {\mathfrak M}_{ \langle D \rangle^\tm {\mathcal R}_n \langle D \rangle^{- \mathtt{c}_1}}^\sharp (s)
 &\leq C_* K_{n - 1}^{- \ta_{1}} \e 
{\gamma^{-2-2(\mathtt{m}-\mathtt{c}_1)}} 
(1 + \| {\mathcal I_{0} } \|_{s + \mu(\mathtt{b}_{1})}^{k_0, \gamma}) \,, 
 \\ 
  {\mathfrak M}_{\langle D \rangle^{\mathtt m}  \langle \partial_{\vphi } \rangle^{\mathtt{b}_{1}}  
{\mathcal R}_n 
\langle D \rangle^{- \mathtt{c}_1}}^\sharp ( s) 
&\leq  C_*  K_{n - 1} \e{\gamma^{-2-2(\mathtt{m}-\mathtt{c}_1)}} 
(1 + \| {\mathcal I_{0} } \|_{s + \mu(\mathtt{b}_{1})}^{k_0, \gamma})\,.
 \end{aligned}
\end{equation}
Define the sets $ \Lambda_n^\gamma $ by
 (recall \eqref{tDtCn} and \eqref{def:DCgt}) 
 \begin{equation}\label{Omega nu + 1 gamma}
\begin{aligned}
 \Lambda_0^\gamma & := {\mathtt T}{\mathtt C}_{\infty}( \gamma, \tau; i_0)
 \cap \mathtt{DC}( 2\gamma, \tau) \quad \text{and for} \quad n \geq 1\,, 
 \\
\Lambda_n^{\gamma}  &:=  \Lambda_n^{\gamma} (i_{0}) 
\\&:=  
\Big\{ \lambda = (\omega, \mathtt h) \in \Lambda_{n - 1}^\gamma  \ : \ 
%\\&  
|\omega \cdot \ell  + \mu_{n - 1}(j) \pm \mu_{n - 1} (j')| \geq  
\frac{\gamma}{ \langle j' \rangle^{ \tau}\langle \ell \rangle^{\tau}}\,,
\\
  & \qquad\qquad  \forall \ell \in \Z^\nu\,, \, 0 <  |\ell  | \leq K_{n - 1}, \ \  j, j' \in S_0^c\,, 
  \mathtt V^T \ell + j \pm j' = 0\,, 
\Big\} \, . 
%\Lambda_n^{\gamma, +}  &:=  \Lambda_n^{\gamma, +} (i_{0}) 
%\\&:=  
%\Big\{ \lambda = (\omega, \mathtt h) \in \Lambda_{n - 1}^\gamma  \ : \ 
%  %\\& 
% |\omega \cdot \ell  + \mu_{n - 1}(j) + \mu_{n - 1}(j')| \geq  
%\frac{\gamma (|j|^{\frac12} + |j'|^{\frac12})}{\langle \ell \rangle^{\tau}}\,,
%\\
%& \qquad\qquad  \forall \ell \in \Z^\nu\,, \,  |\ell|  \leq N_{n - 1}\,,  \ j, j' \in S_0^c\,, \mathtt V^T \ell + j + j' = 0  \Big\} \, . 
\end{aligned}
\end{equation}
For $ n \geq 1 $, there exists a symplectic, momentum preserving transformation
defined for all $ (\omega, \mathtt h) $ in $\mathtt \Omega \times [\mathtt h_1, \mathtt h_2]$, of the form 
\begin{equation}\label{Psi nu - 1}
 {\Phi}_{n - 1} := {\rm exp}(\Psi_{n - 1})\,, \quad 
{\Psi}_{n - 1} := \begin{pmatrix}
\mathtt\Psi_{n - 1, 1} & \mathtt\Psi_{n - 1, 2} \\
\overline{\mathtt \Psi}_{n - 1, 2} & \overline{\mathtt \Psi}_{n - 1, 1}
\end{pmatrix}\,
\end{equation}
such that for all $\lambda = (\omega, {\mathtt h}) \in \Lambda_n^\gamma $ 
the following conjugation formula holds:
\begin{equation}\label{coniugionu+1}
{\mathbf L}_n = { \Phi}_{n - 1}^{-1} {\mathbf L}_{n - 1} {\Phi}_{n - 1}\, .
\end{equation}
The operators $ A_{n - 1} \in \{ \langle D \rangle^{{\mathtt m}} \Psi_{n - 1} \langle D \rangle^{- {\mathtt m}},\, \langle D \rangle^{ \mathtt{c}_{1}} \Psi_{n - 1} \langle D \rangle^{- { \mathtt{c}_{1}} } \}$
and $ \langle \partial_{\vphi} \rangle^{\mathtt{b}_{1}} A_{n - 1}$, 
are $ {\mathcal D}^{k_0} $-modulo-tame on $\mathtt\Omega \times [\mathtt h_1, \mathtt h_2]$ 
with modulo-tame constants satisfying, for all $s \in [s_0, \bar{s}] $,
($ \tau_3$ and $\mathtt a_{1} $ are defined in \eqref{costanti.riducibilita} )
\begin{align}\label{tame Psi nu - 1}
  {\mathfrak M}_{A_{n - 1}}^\sharp (s) 
&\lesssim_{\bar{s}}
K_{n - 1}^{\tau_3} K_{n - 2}^{- \mathtt a_{1}} 
\e  {\gamma^{-3-2(\mathtt{m}-\mathtt{c}_1)}} 
(1 + \| {\mathcal I }_{0} \|_{s + \mu(\mathtt{b}_{1})}^{k_0, \gamma}) \, , 
\\
 {\mathfrak M}_{\langle \partial_\vphi \rangle^{\mathtt{b}_{1}}A_{n - 1} }^\sharp  (s) 
&\lesssim_{\bar{s}}
K_{n - 1}^{\tau_3} K_{n - 2} 
\e  {\gamma^{-3-2(\mathtt{m}-\mathtt{c}_1)}} 
(1 + \| {\mathcal I }_{0} \|_{s + \mu(\mathtt{b}_{1})}^{k_0, \gamma}) \, . 
\label{tame Psi nu - 1 vphi x b}
\end{align}

\item[${\bf(S2)_{n}}$] Let $ i_1(\omega, \mathtt h) $, $i_2(\omega, \mathtt h) $ satisfy \eqref{ansatz_I0_s0}
Then % for all $\nu = 0, \ldots n$, 
for all $(\omega, \mathtt h) \in \mathtt \Omega \times [\mathtt h_1, \mathtt h_2]$
%\begin{color}{red} (OPPURE su tutto $\R^{|S^+|} \times [\th_1,\th_2]$ ??) \end{color}
%with $\gamma_1, \gamma_2 \in [\gamma/2, 2 \gamma]$, 
the following estimates hold
\begin{align}\label{stima R nu i1 i2}
 \! \! \| | \langle D \rangle^{\mathtt m} 
\Delta_{12} {\mathcal R}_n \langle D \rangle^{- \mathtt{c}_{1} } | \|_{{\mathcal L}(H^{s_0})}  
&\lesssim_{\bar{s}}  
K_{n - 1}^{- \mathtt a_{1}}
\e  {\gamma^{-2-2(\mathtt{m}-\mathtt{c}_1)}} 
 \|i_1 - i_2\|_{s_0 +  \mu(\mathtt{b}_{1})}, 
\\ 
 \label{stima R nu i1 i2 norma alta}
\! \! \|  |\langle D \rangle^{\mathtt m} \langle \partial_{\vphi }  \rangle^{\mathtt{b}_{1}}  
\Delta_{12} {\mathcal R}_n \langle D \rangle^{- \mathtt{c}_{1}}  | \|_{{\mathcal L}(H^{s_0})}
&\lesssim_{\bar{s}} 
K_{n - 1}
\e  {\gamma^{-2-2(\mathtt{m}-\mathtt{c}_1)}}  \|i_1 - i_2\|_{ s_0 +  \mu(\mathtt{b}_{1})}\,.
\end{align}
Moreover, for $n \geq 1$, for all $j \in S^c_0$, 
\begin{align}\label{r nu - 1 r nu i1 i2}
 \big|\Delta_{12}(r_n(j) - r_{n - 1}(j)) \big|
 &\lesssim_{\bar{s}} 
  K_{n - 2}^{- \mathtt a_{1}} 
 \e  {\gamma^{-2-2(\mathtt{m}-\mathtt{c}_1)}}  
 |j|^{- {\mathtt m} + \mathtt{c}_{1} }   
 \| i_1-i_2 \|_{ s_0  + \mu(\mathtt{b}_{1})}\,, 
 \\
 \ |\Delta_{12} r_{n}(j)| &\lesssim_{\bar{s} }  
\e   {\gamma^{-2-2(\mathtt{m}-\mathtt{c}_1)}}|j|^{- {\mathtt m} + \mathtt{c}_{1} }  
\| i_1-i_2 \|_{ s_0  + \mu(\mathtt{b}_{1})}\,. \label{r nu i1 - r nu i2}
\end{align}

\end{itemize}
\end{prop}

%\subsubsection{Proof of Proposition \ref{iterazione riducibilita}}

Proposition \ref{iterazione riducibilita} is proved by induction by constructing 
$\Psi_n$, $\Phi_n$, ${\mathbf L}_{n + 1}$
etc.,  iteratively.
In subsection \ref{sec:homologicalKAM} we show how to construct the suitable change of coordinates 
$\Phi_{n}$.
Then in subsection \ref{sec:inductionKAM} we prove some properties of the 
conjugate operator ${\mathbf L}_{n + 1}$.
At the end of the section we will prove by induction that 
 the properties
 ${\bf (S1)_{n }}, {\bf (S2)_{n }}$ in Proposition \ref{iterazione riducibilita} holds for any $n$.
 %at the step $n+1$.

\subsubsection{The homological equation}\label{sec:homologicalKAM}

Consider the operator ${\bf L}_{n}$ in \eqref{cal L nu} at the step $n$.
We want to determine $\Psi_n$ in such a way that it solves the homological equation 
\begin{equation}\label{equazione omologica}
\begin{aligned}
& \omega \cdot \partial_\vphi \Psi_n +  \ii [\mathcal{D}_n, \Psi_n] + \Pi_{K_n} {\mathcal  R}_n 
= [{\mathcal  R}_n] \,, 
\qquad    
[{\mathcal R}_n] := \begin{pmatrix}
\widehat{\mathtt R}_{n, 1}(0) & 0 \\
0 &\overline{\widehat{\mathtt R}_{n, 1}(0)}
\end{pmatrix}\, ,
\end{aligned}
\end{equation}
where the projector $ \Pi_{K_n} $ is defined in \eqref{proiettore-oper} %and 
%$ \Pi_{K_n}^\bot = 
%{\rm Id } - \Pi_{K_n}$ 
and where we  assumed that the normal form $\mathcal{D}_{n}$
and the remainder $\mathcal{R}_{n}$ fulfil the inductive assumptions ${\bf (S1)_n}, {\bf (S2)_n}$.
The solution of \eqref{equazione omologica} is given by 
\begin{equation}\label{shomo1}
\begin{aligned}
& \Psi_n = \begin{pmatrix}
\mathtt\Psi_{n, 1} & \mathtt\Psi_{n, 2} \\
\overline{\mathtt\Psi}_{n, 2} & \overline{\mathtt \Psi}_{n, 1} 
\end{pmatrix}\,, \\
& \widehat {\mathtt\Psi}_{n, 1}(\ell  )_j^{j'} := \begin{cases}
 \dfrac{\widehat{\mathtt R}_{n, 1}(\ell)_j^{j'} }{\ii (\omega \cdot \ell  + \mu_n (j) - \mu_n (j'))} \qquad \forall 
0 < | \ell  | \leq K_{n} \, , \, j, j'\in S_0^c 
\\
 \qquad \qquad \qquad \qquad \qquad \qquad  \mathtt V^T \ell + j - j' = 0\,, \\
0 \qquad \text{otherwise}\, , 
\end{cases}   
\\
& \widehat{\mathtt \Psi}_{n, 2}(\ell  )_j^{j'} := \begin{cases}
 \dfrac{\widehat{\mathtt R}_{n, 2}(\ell)_j^{j'} }{\ii (\omega \cdot \ell  + \mu_n (j) +  \mu_n (j'))} \qquad \forall 
| \ell  | \leq K_{n} \, , \, j, j'\in S_0^c \\
 \qquad \qquad \qquad \qquad \qquad \qquad  \mathtt V^T \ell + j + j' = 0\,, \\
0 \qquad \text{otherwise}\, .
\end{cases}  
\end{aligned}
\end{equation}
Note that the denominators in \eqref{shomo1} do not vanish
for $ \lambda = (\omega, \mathtt h) \in \Lambda_{n +1}^\gamma $ 
(cf. \eqref{Omega nu + 1 gamma}). In the next lemma we provide some quantitative 
estimates on the transformation $\Psi_n$ defined above. 

\begin{lemma}[\bf Estimates on the homological equation]\label{stime eq omologica}
Let $\tau_3$ as in \eqref{costanti.riducibilita} and assume \eqref{KAM smallness condition1}.
 The operator $\Psi_n$ defined for 
$\lambda = (\omega, \mathtt h) \in  \Lambda_{n + 1}^\gamma$ 
admits an extension to the whole parameter space 
$\mathtt \Omega \times [\mathtt h_1, \mathtt h_2]$, 
which we still denote by $\Psi_n$, which satisfies the following property. 

\noindent
The operators $ A_n \in \{ \langle D \rangle^{ \mathtt m} \Psi_n 
\langle D \rangle^{-  \mathtt m}\,,\, 
\langle D \rangle^{ \mathtt c_1} \Psi_n \langle D \rangle^{- \mathtt c_1}\,,\, \Psi_n \}$, 
$\langle \partial_{\vphi} \rangle^{\mathtt{b}_{1}} A_n$ 
are ${\mathcal D}^{k_0}$-modulo tame and, for any $s_0 \leq s \leq \bar{s}$, 
satisfy the estimates 
\begin{equation}\label{stima Pn s lemma homo}
\begin{aligned}
 {\mathfrak M}^\sharp_{A_n}(s) 
&
\lesssim
K_n^{\tau_3} \gamma^{-1} 
 {\mathfrak M}_{ \langle D \rangle^\tm {\mathcal R}_n \langle D \rangle^{- \mathtt{c}_1}}^\sharp (s)\,,
\\
 {\mathfrak M}^\sharp_{\langle \partial_\vphi \rangle^{\mathtt{b}_{1}}A_n }(s)  
&
 \lesssim 
 K_n^{\tau_3}\gamma^{-1}
  {\mathfrak M}_{\langle D \rangle^{\mathtt m}  \langle \partial_{\vphi } \rangle^{\mathtt{b}_{1}}  
{\mathcal R}_n 
\langle D \rangle^{- \mathtt{c}_1}}^\sharp ( s)\,.
\end{aligned}
\end{equation}
%In particular, for $s = s_0$, using the ansatz \eqref{ansatz_I0_s0} 
%(with $\mu_0\rightsquigarrow \mu(\tb_1)$),
%%$\| {\mathcal I}_{0} \|_{s_0 + \mu(\mathtt{b}_{1})}^{k_0, \gamma} \leq 1$, 
%one has 
%\begin{equation}\label{stima Pn s0 lemma homo}
%\begin{aligned}
% {\mathfrak M}^\sharp_{A_n}(s_0)\,,\, {\mathfrak M}^\sharp_{ \Psi_n }(s_0) 
% &
% \lesssim_{\bar{s}}
% K_n^{\tau_3}  K_{n - 1}^{- \ta_1} \e  {\gamma^{-3-2(\mathtt{m}-\mathtt{c}_1)}}\,,
%  \\
% {\mathfrak M}^\sharp_{\langle \partial_\vphi \rangle^{\mathtt{b}_{1}}A_n}(s_0)
% & \lesssim_{\bar{s}}
% K_n^{\tau_3} K_{n - 1} \e  {\gamma^{-3-2(\mathtt{m}-\mathtt{c}_1)}} \,. 
%\end{aligned}
%\end{equation}
Assume that ${i}_1(\lambda), {i}_2(\lambda)$, $\lambda \in \R^\nu \times [\mathtt h_1, \mathtt h_2]$ 
are two tori satisfying \eqref{ansatz_I0_s0}. Then, for any 
$\lambda = (\omega, \mathtt h) \in \mathtt \Omega \times [\mathtt h_1, \mathtt h_2]$, 
%$\gamma_1 , \gamma_2 \in [\gamma/2\,,\, 2 \gamma]$,  
the following estimates hold: 
\begin{equation}\label{stime delta 12 Psin}
\begin{aligned}
% \| |\Delta_{12}A_n|\|_{{\mathcal L}(H^{s_0})} 
%&\lesssim_{\bar{s}}
%K_n^{ \tau_3} \gamma^{- 1} \| |\langle \partial_\vphi \rangle^{\mathtt b_1} \mathtt R_1 
%\langle D \rangle^{- \mathtt c_1}| \|_{{\mathcal L}(H^{s_0})} \| i_1 - i_2 \|_{s_0 + \mu(\mathtt b_1)} 
%\\& 
%\qquad 
%+ K_n^{ \tau_3} \gamma^{- 1} \| |\langle \partial_\vphi \rangle^{\mathtt b_1} \Delta_{12} \mathtt R_1 \langle D \rangle^{- \mathtt c_1}| \|_{{\mathcal L}(H^{s_0})} \,,
% \\
 \| | \langle \partial_\vphi \rangle^{p} \Delta_{12} A_n |\|_{{\mathcal L}(H^{s_0})} 
&\lesssim_{\bar{s}}
K_n^{ \tau_3} \gamma^{- 1} 
\mathfrak{R}_{n}(s_0,p)
%\| |\langle \partial_\vphi \rangle^{p} 
%\langle D \rangle^{ \mathtt{m}}\mathcal{R}_n 
%\langle D \rangle^{- \mathtt c_1}| \|_{{\mathcal L}(H^{s_0})} 
\| i_1 - i_2 \|_{s_0 + \mu(\mathtt b_1)} 
\\& 
\qquad 
+ K_n^{ \tau_3} \gamma^{- 1} \| |\langle D \rangle^{ \mathtt{m}} \langle \partial_\vphi \rangle^{p} 
\Delta_{12}\mathcal{R}_n 
\langle D \rangle^{- \mathtt c_1}| \|_{{\mathcal L}(H^{s_0})} \,,
\end{aligned}
\end{equation}
for $p=0,\mathtt{b}_1$ and 
where
\begin{equation}\label{bassaGotR0}
\mathfrak{R}_{n}(s_0,p):=\sup_{j=1,2}\| |
\langle D \rangle^{ \mathtt{m}}
\langle\pa_{\vphi}\rangle^{p}\mathcal{R}_n(i_{j})
\langle D \rangle^{- \mathtt c_1}| \|_{{\mathcal L}(H^{s_0})}\,,\qquad p=0,\mathtt{b}_1\,.
\end{equation}
Finally, $\Psi_{n}$ is Hamiltonian and momentum preserving.
\end{lemma}

\begin{proof}
To simplify notations in this proof we drop the index $n$, 
namely we write $\Psi, {\mathcal R}, \mu(j)$, etc. 
instead of ${\Psi}_n, {\mathcal R}_n, \mu_n(j)$, etc. 

\noindent
We prove the claimed statement for the operator $\mathtt \Psi_1$. 
The statement for $\mathtt \Psi_2$ can be proved in a similar fashion. 

\smallskip
\noindent
\emph{Small divisors.}
First let us analyze the $k_0$-times differentiable function 
\[
f_{\ell j j'}(\lambda) \equiv f_{\ell j j'}(\lambda , i_0(\lambda)) := \omega \cdot \ell + \mu (j; \lambda, i_0(\lambda)) - \mu(j'; \lambda, i_0(\lambda))\,, 
\qquad \lambda = (\omega, \mathtt h) \in \mathtt \Omega \times [\mathtt h_1, \mathtt h_2]\,,
\]
where $\ell \in \Z^\nu {\setminus\{0\}}$, $j, j' \in S_0^c$ and $\mathtt V^T \ell + j - j' = 0$. 

\noindent 
We first estimate $\mu(j; \lambda)$, $j \in S_0^c$. 
By Theorem \ref{red.Lomega.smooth.rem}-$(i)$
%Lemma \ref{lemma autovalore 0 pre red} 
and by the induction hypotheses \eqref{mu j nu}, \eqref{stima rj nu}, 
one has that for any $k \in \N^{\nu + 1}$, $|k| \leq k_0$,
\begin{equation}\label{stima.mu.j}
\begin{aligned}
|\partial_\lambda^k \mu(j ;  \lambda)|  
& \leq 
|\partial_\lambda^k \mu_0(j; \lambda)| 
+ |\partial_\lambda^k r(j; \lambda)| 
\lesssim |j|  + \gamma^{- |k|} |r(j)|^{k_0, \gamma} 
\\& 
\lesssim
|j| + C_{\bar{s}}\e \gamma^{- |k| -2-2(\mathtt{m}-\mathtt{c}_1)} |j|^{- \mathtt m + \mathtt c_1} 
\lesssim |j|\,,
\end{aligned}
\end{equation}
for some $C_{\bar{s}}>0$, 
using that $\e \gamma^{- k_0 -2-2(\mathtt{m}-\mathtt{c}_1)}  \ll 1$ (since $\gamma \in (0, 1)$, 
$\mathtt m > \mathtt c_1$). 
Similarly if $i_1(\lambda), i_2(\lambda)$, $\lambda \in \mathtt\Omega \times [\mathtt h_1, \mathtt h_2]$ 
are two tori, arguing as above, by Theorem \ref{red.Lomega.smooth.rem}-$(i)$
%Lemma \ref{lemma autovalore 0 pre red} 
and by \eqref{r nu i1 - r nu i2}, one can show that 
\begin{equation}\label{delta i2 mu j eq omologica}
|\Delta_{12} \mu(j)| = |\mu (j; \lambda, i_1(\lambda)) - \mu (j; \lambda, i_2(\lambda))  | 
\lesssim_{\bar{s}} \e   {\gamma^{-2-2(\mathtt{m}-\mathtt{c}_1)}}
|j|  \| i_1-i_2 \|_{ s_0  + \mu(\mathtt{b}_{1})}
\end{equation}
The latter estimates \eqref{stima.mu.j}, \eqref{delta i2 mu j eq omologica} imply that 
\begin{equation}\label{stima.f.j}
\begin{aligned}
 |\partial_\lambda^k f_{\ell j j'}(\lambda)| &\lesssim |\ell| + |j| + |j'| \lesssim |\ell| |j||j'|\,, 
\qquad \forall k \in \N^\nu\,, \quad 0 < |k| \leq k_0\,, 
\\
 |\Delta_{12} f_{\ell j j'}|& \lesssim\bar{s} 
\e   {\gamma^{-2-2(\mathtt{m}-\mathtt{c}_1)}} |j| |j'|  \| i_1-i_2 \|_{ s_0  + \mu(\mathtt{b}_{1})}\,.
\end{aligned}
\end{equation}
Then, by defining 
\[
\Lambda (\ell, j, j') := \Big\{ \lambda = (\omega , \mathtt  h) \in \mathtt\Omega \times [\mathtt h_1, \mathtt h_2]
: |f_{\ell j j'}(\lambda)| \geq \rho  \Big\}\,, \qquad \rho := \frac{\gamma}{2 \langle \ell \rangle^\tau |j'|^\tau }\,,
\]
one has that 
\[
\frac{1}{f_{\ell j j'}} : \Lambda (\ell, j, j')  \to \R, 
\qquad 
\lambda \mapsto \frac{1}{f_{\ell j j'}(\lambda, i_0(\lambda))} \equiv \frac{1}{f_{\ell j j'}(\lambda)} \,,
\] 
is $k_0$-times differentiable 
and for any $k \in \N^{\nu + 1}$, $|k| \leq k_0$, for any $\lambda \in \Lambda (\ell, j, j')$, 
it satisfies the estimate 
\[
\Big| \partial_\lambda^k \Big( \dfrac{1}{f_{\ell j j'}(\lambda)}\Big) \Big|
\lesssim 
\dfrac{|\ell|^{ |k|} |j|^{|k|}|j'|^{|k|}}{\rho^{|k| + 1} } 
\lesssim 
\langle \ell \rangle^{k_0 + \tau(k_0 + 1)}|j|^{k_0} |j'|^{k_0 + \tau (k_0 + 1)} \gamma^{- (|k| + 1)}\,.
\]
By recalling \eqref{def-gl}, we define 
\[
g_{\ell j j'}(\lambda) \equiv g_{\ell j j'}(\lambda, i_0(\lambda)) 
:= h_\rho\big( f_{\ell j j'}(\lambda, i_0(\lambda)) \big)\,, 
\qquad 
\lambda \in \mathtt \Omega \times [\mathtt h_1, \mathtt h_2]\,,
\]
and by applying Lemmata \ref{lemma:cut-off sd}, \ref{lemma lip estensione lip} 
(recall the estimates \eqref{stima.f.j}), one obtains the following properties:  
\begin{equation}\label{prop extension omologica}
\begin{aligned}
 g_{\ell j j'}(\lambda, i_0(\lambda)) &= \dfrac{1}{f_{\ell j j'}(\lambda, i_0(\lambda))}\,, \qquad \forall \lambda \in \Lambda (\ell, j, j')\,, 
\\
 |\partial_\lambda^k g_{\ell j j'}(\lambda)| 
 &\lesssim
 \langle \ell \rangle^{\tau_1} |j|^{\tau_1} |j'|^{\tau_1} \gamma^{- (|k| + 1)} \,, 
 \qquad  \forall \lambda \in \mathtt \Omega \times [\mathtt h_1, \mathtt h_2]\,, \quad \forall |k| \leq k_0 \\
 & \qquad \text{where} \qquad \tau_1 := k_0 + \tau(k_0 + 1)\,. %\\ 
\end{aligned}
\end{equation}
From the latter estimate, we deduce
\begin{equation}\label{prop extension omologicaBIS}
\begin{aligned}
 |g_{\ell j j'}|^{k_0, \gamma} 
 &\lesssim 
 \langle \ell \rangle^{\tau_1}|j|^{\tau_1} |j'|^{\tau_1} \gamma^{- 1}\,,
 \\
 |\Delta_{12} g_{\ell j j'}| 
 & = 
 | g_{\ell j j'}(\lambda, i_1(\lambda)) -  g_{\ell j j'}(\lambda, i_2(\lambda))|  
 \\&  
 \lesssim_{\bar{s}} 
 \e   {\gamma^{-4-2(\mathtt{m}-\mathtt{c}_1)}} \langle \ell \rangle^{2 \tau} |j| |j'|^{2 \tau + 1}  
 \| i_1-i_2 \|_{ s_0  + \mu(\mathtt{b}_{1})}\,. 
 \end{aligned}
\end{equation}
%implying that 
%\begin{equation}\label{prop extension omologica 2}
%  |g_{\ell j j'}|^{k_0, \gamma} \lesssim \langle \ell \rangle^{\tau_1}|j|^{\tau_1} |j'|^{\tau_1} \gamma^{- 1}\,. 
%\end{equation}

\smallskip
\noindent
\emph{Estimates on the transformation.}
Then, we extend $\mathtt \Psi_1 (\lambda ) \equiv \mathtt \Psi_1(\lambda, i_0(\lambda))$ for any $\lambda \in \mathtt \Omega \times [\mathtt h_1, \mathtt h_2]$, by defining 
\[
\widehat{\mathtt \Psi}_1(\ell)_j^{j'} (\lambda, i_0(\lambda)) :=
\begin{cases}
g_{\ell j j'} (\lambda, i_0(\lambda))  \widehat{\mathtt R}_1(\ell)_j^{j'}(\lambda, i_0(\lambda))  
\quad \text{if} \quad 0 < |\ell| < K\,,\quad \mathtt V^T \ell + j - j' = 0\,, \\
0 \qquad \text{otherwise\,.}
\end{cases}
\]
We now provide estimates on the operators
\begin{equation}\label{sigaro1}
{\mathcal F} := \langle D \rangle^{ \mathtt m} 
\langle \partial_\vphi \rangle^{\mathtt{b}_{1}} 
\mathtt \Psi_1 \langle D \rangle^{- \mathtt m}\,,
\qquad
{\mathcal G} := \langle D \rangle^{\mathtt c_1} 
\langle \partial_\vphi \rangle^{\mathtt{b}_{1}} 
\mathtt \Psi_1 \langle D \rangle^{ - \mathtt c_1}\,,
\end{equation}
being the estimates for the case $\mathtt{b}_1=0$ similar.
\medskip

\noindent
{\sc Estimate of $\mathcal{F}$ in \eqref{sigaro1}}.
One clearly has 
\[
\widehat{\mathcal F}(\ell)_j^{j'} 
= \langle \ell \rangle^{\mathtt{b}_{1}}\langle j \rangle^{ \mathtt m} g_{\ell j j'} 
\widehat{\mathtt R}_1(\ell)_j^{j'} \langle j' \rangle^{- \mathtt m}\,,
\]
and for any $k \in \N^{\nu + 1}$, $|k| \leq k_0$, 
\begin{equation}\label{mortal kombat omologico 1}
\begin{aligned}
 |\partial_\lambda^k \widehat{\mathcal F}(\ell)_j^{j'} |  
&\lesssim 
\sum_{k_1 + k_2 = k} \langle j \rangle^{ \mathtt m} |\partial_\lambda^{k_1}g_{\ell j j'}|  
\langle \ell \rangle^{\mathtt{b}_{1}}
|\partial_\lambda^{k_2}\widehat{\mathtt R}_1(\ell)_j^{j'}| \langle j' \rangle^{- \mathtt m} 
\\& 
\stackrel{\eqref{prop extension omologica}}{\lesssim}  
\langle \ell \rangle^{\tau_1} \langle j \rangle^{\tau_1}  \langle j' \rangle^{\tau_1 + \mathtt c_1 - \mathtt m}  
\sum_{k_1 + k_2 = k}   \gamma^{- (|k_1| + 1)} \langle j \rangle^{\mathtt m} \langle \ell \rangle^{\mathtt{b}_{1}}  |\partial_\lambda^{k_2}\widehat{\mathtt R}_1(\ell)_j^{j'}| \langle j' \rangle^{- \mathtt c_1}  \,. 
\end{aligned}
\end{equation}
The conservation of momentum $\mathtt V^T \ell + j - j' = 0$ implies that $|j - j'| \lesssim |\ell|$, 
therefore, using that 
$\langle j \rangle \lesssim \langle j - j' \rangle + \langle j' \rangle
 \lesssim \langle j - j' \rangle \langle j' \rangle 
\lesssim \langle \ell \rangle \langle j'\rangle$, 
one gets that 
\begin{equation}\label{mortal kombat omologico 2}
\begin{aligned}
& \langle \ell \rangle^{\tau_1} \langle j \rangle^{\tau_1}  
\langle j' \rangle^{\tau_1 + \mathtt c_1 - \mathtt m}  
\lesssim \langle \ell \rangle^{2 \tau_1} \langle j' \rangle^{2 \tau_1 + \mathtt c_1 - \mathtt m}  
\lesssim  \langle \ell \rangle^{2 \tau_1}\,,
 \end{aligned}
\end{equation}
since by \eqref{costanti.riducibilita} (recall also the definition of $\tau_1$ in \eqref{prop extension omologica}), one has that $\mathtt m \geq 2 \tau_1 + \mathtt c_1$. 
Hence, by \eqref{mortal kombat omologico 1}, \eqref{mortal kombat omologico 2} 
and using that $|\ell| \leq K$, one gets for any $k \in \N^{\nu + 1}$, $|k| \leq k_0$, 
the estimate
\begin{equation}\label{mortal kombat omologico 3}
\begin{aligned}
&|\partial_\lambda^k \widehat{\mathcal F}(\ell)_j^{j'} |  
\lesssim K^{2 \tau_1} 
\gamma^{- (|k| + 1)}\sum_{k_1 + k_2 = k}   
\gamma^{|k_2|} \langle j \rangle^{\mathtt m} \langle \ell \rangle^{\mathtt{b}_{1}} 
|\partial_\lambda^{k_2}\widehat{\mathtt R}_1(\ell)_j^{j'}| \langle j' \rangle^{- \mathtt c_1} \,. 
\end{aligned}
\end{equation}
Therefore, for any $k \in \N^{\nu + 1}$, $|k| \leq k_0$, $s_0 \leq s \leq \bar{s}$ 
one has  
\begin{equation*}
\begin{aligned}
& \| |\partial_\lambda^k {\mathcal F}| h \|_s^2  
\leq \sum_{\ell, j} \langle \ell, j \rangle^{2 s} \Big( \sum_{\ell', j'} 
|\partial_\lambda^k \widehat{\mathcal F}(\ell - \ell')_j^{j'}| |\widehat h(\ell', j')|\Big)^2  
\\&  \hspace{-0.2cm}
\stackrel{\eqref{mortal kombat omologico 3}}{\lesssim}  K^{4 \tau_1} 
\gamma^{- 2 (|k| + 1)} 
\sum_{|k| \leq k_2} \gamma^{2 |k_2|}\sum_{\ell, j}   \langle \ell, j \rangle^{2 s} 
\Big(  \sum_{\ell', j'} \langle j \rangle^{\mathtt m} \langle \ell - \ell'\rangle^{\mathtt{b}_{1}}  
|\partial_\lambda^{k_2}\widehat{\mathtt R}_1(\ell)_j^{j'}| \langle j' \rangle^{- \mathtt c_1}   
|\widehat h(\ell', j')|\Big)^2 
\\&
\lesssim K^{4 \tau_1} 
\gamma^{- 2 (|k| + 1)} \sum_{|k| \leq k_2} \gamma^{2 |k_2|} 
\| | \langle \partial_\vphi \rangle^{\mathtt{b}_{1}} \langle D \rangle^{\mathtt m} 
\partial_\lambda^{k_2} {\mathtt R}_1  \langle D \rangle^{- \mathtt c_1}| |h|\|^2_s 
\\& 
\lesssim
K^{4 \tau_1}\gamma^{- 2 (|k| + 1)} 
\Big( {\mathfrak M}^\sharp_{\langle \partial_\vphi \rangle^{\mathtt{b}_{1}} 
\langle D \rangle^{\mathtt m} {\mathtt R}_1 \langle D \rangle^{- \mathtt c_1} }(s) 
\| |h| \|_{s_0} + {\mathfrak M}^\sharp_{\langle \partial_\vphi \rangle^{\mathtt{b}_{1}} 
\langle D \rangle^{\mathtt m}  {\mathtt R}_1 
\langle D \rangle^{- \mathtt c_1} }(s_0) \| |h| \|_{s} \Big)^2\,.
\end{aligned}
\end{equation*}
Using that $\| |h| \|_s = \| h \|_s$, % and recalling \eqref{costanti.riducibilita}, 
we have that the latter bounds implies 
\eqref{stima Pn s lemma homo} for the operator $\mathcal{F}$ in \eqref{sigaro1}.

\noindent
To estimate the Lipschitz variation we reason similarly.
We consider the coefficients
\[
\Delta_{12} \widehat{\mathcal F}(\ell)_j^{j'} 
= \widehat{\mathcal F}(\ell)_j^{j'}(\lambda, i_1(\lambda)) 
- \widehat{\mathcal F}(\ell)_j^{j'}(\lambda, i_2(\lambda))\,.
\] 
Using that $\e   \gamma^{-3-2(\mathtt{m}-\mathtt{c}_1)} \leq 1$ 
(recall the anstatz \eqref{smallcondepsilon})
one has that 
\begin{equation}\label{stima delta 12 mathcal F matrix}
\begin{aligned}
&|\Delta_{12} \widehat{\mathcal F}(\ell)_j^{j'}|  
  \leq 
|\Delta_{12} g_{\ell j j'}|  \langle \ell \rangle^{\mathtt{b}_{1}}\langle j \rangle^{ \mathtt m} 
|\widehat{\mathtt R}_1(\ell)_j^{j'}(\lambda, i_1(\lambda))| \langle j' \rangle^{- \mathtt m}  
\\& 
\qquad + |g_{\ell j j'}(\lambda, i_2(\lambda))| \langle \ell \rangle^{\mathtt{b}_{1}}\langle j \rangle^{ \mathtt m} 
|\Delta_{12} \widehat{\mathtt R}_1(\ell)_j^{j'}| \langle j' \rangle^{- \mathtt m}  
\\& \hspace{-0.35cm}
\stackrel{\eqref{prop extension omologicaBIS}}{\lesssim_{\bar{s}}}  
\frac{\e\langle \ell \rangle^{2 \tau} \langle j \rangle }{\gamma^{4+2(\mathtt{m}-\mathtt{c}_1)}}
%\e   {\gamma^{-4-2(\mathtt{m}-\mathtt{c}_1)}} 
 \langle j' \rangle^{2 \tau + 1 + \mathtt c_1 - \mathtt m}  
\langle \ell \rangle^{\mathtt{b}_{1}}\langle j \rangle^{ \mathtt m} 
|\widehat{\mathtt R}_1(\ell)_j^{j'}(\lambda, i_1(\lambda))| \langle j' \rangle^{- \mathtt c_1}  
\| i_1-i_2 \|_{ s_0  + \mu(\mathtt{b}_{1})}   
\\& 
\qquad 
+  \langle \ell \rangle^{\tau_1} |j|^{\tau_1} |j'|^{\tau_1 + \mathtt c_1 - \mathtt m} \gamma^{-  1}  
\langle \ell \rangle^{\mathtt{b}_{1}}\langle j \rangle^{ \mathtt m} 
|\Delta_{12} \widehat{\mathtt R}_1(\ell)_j^{j'}| \langle j' \rangle^{- \mathtt c_1} 
\\& 
{\lesssim}  
\langle \ell \rangle^{\tau_1} |j|^{\tau_1} |j'|^{\tau_1 + \mathtt c_1 - \mathtt m} \gamma^{-  1}  
\langle \ell \rangle^{\mathtt{b}_{1}}\langle j \rangle^{ \mathtt m} 
|\widehat{\mathtt R}_1(\ell)_j^{j'}(\lambda, i_1(\lambda))| \langle j' \rangle^{- \mathtt c_1}  
\| i_1-i_2 \|_{ s_0  + \mu(\mathtt{b}_{1})}  
\\& 
\qquad + \langle \ell \rangle^{\tau_1} |j|^{\tau_1} |j'|^{\tau_1 + \mathtt c_1 - \mathtt m} \gamma^{-  1}  
\langle \ell \rangle^{\mathtt{b}_{1}}\langle j \rangle^{ \mathtt m} 
|\Delta_{12} \widehat{\mathtt R}_1(\ell)_j^{j'}| \langle j' \rangle^{- \mathtt c_1} 
\\& \hspace{-0.75cm}
\stackrel{\eqref{mortal kombat omologico 2}, |\ell| \leq K}{\lesssim} 
K^{2 \tau_1} \gamma^{-  1}  \langle \ell \rangle^{\mathtt{b}_{1}}\langle j \rangle^{ \mathtt m} 
|\widehat{\mathtt R}_1(\ell)_j^{j'}(\lambda, i_1(\lambda))| \langle j' \rangle^{- \mathtt c_1}  
\| i_1-i_2 \|_{ s_0  + \mu(\mathtt{b}_{1})} 
 \\& 
 \qquad + K^{2 \tau_1} \gamma^{-  1}  \langle \ell \rangle^{\mathtt{b}_{1}}\langle j \rangle^{ \mathtt m} 
 |\Delta_{12} \widehat{\mathtt R}_1(\ell)_j^{j'}| \langle j' \rangle^{- \mathtt c_1}\,.
\end{aligned}
\end{equation}
By \eqref{stima delta 12 mathcal F matrix}
one can check that $\| |\Delta_{12} {\mathcal F}| \|_{{\mathcal L}(H^{s_0})} $ satisfies the estimate
\eqref{stime delta 12 Psin}.
%
%
%Using that $\| |h| \|_s = \| h \|_s$, the latter inequality implies that 
%\[
%\begin{aligned}
%{\mathfrak M}^\sharp_{\mathcal F}(s) 
%& \lesssim 
%K_n^{2 \tau_1} \gamma^{- 1}  {\mathfrak M}^\sharp_{\langle \partial_\vphi \rangle^{\mathtt{b}_{1}} 
%\langle D \rangle^{\mathtt m} {\mathcal R}_n \langle D \rangle^{- \mathtt c_1} }(s)  
%\\& 
%\stackrel{\eqref{stima cal R nu}}{\lesssim_{\bar{s}}} 
%K_n^{2 \tau_1 } K_{n - 1} \e{\gamma^{-3-2(\mathtt{m}-\mathtt{c}_1)}}  
%(1 + \| {\mathcal I}_{0} \|_{s + \mu(\mathtt{b}_{1})}^{k_0, \gamma})\,.
%\end{aligned}
%\] 
%The corresponding estimate for $s = s_0$ follows by the latter one, 
%using also the ansatz \eqref{ansatz_I0_s0}. Similarly, by using \eqref{stima delta 12 mathcal F matrix}, and the induction estimates \eqref{stima cal R nu}, \eqref{stima R nu i1 i2 norma alta}, one obtains that 
%\[
%\begin{aligned}
%\| |\Delta_{12} {\mathcal F}| \|_{{\mathcal L}(H^{s_0})} 
%& \lesssim 
%K_n^{2 \tau_1} \gamma^{- 1} \| |\langle \partial_\vphi \rangle^{\mathtt b_1} \mathtt R_1 
%\langle D \rangle^{- \mathtt c_1}| \|_{{\mathcal L}(H^{s_0})} \| i_1 - i_2 \|_{s_0 + \mu(\mathtt b_1)} 
%\\& 
%\qquad 
%+ K_n^{2 \tau_1} \gamma^{- 1} \| |\langle \partial_\vphi \rangle^{\mathtt b_1} \Delta_{12} \mathtt R_1 \langle D \rangle^{- \mathtt c_1}| \|_{{\mathcal L}(H^{s_0})} 
%\\& 
%\lesssim 
%K_n^{2 \tau_1 } N_{n - 1} \e{\gamma^{-3-2(\mathtt{m}-\mathtt{c}_1)}}  \| i_1 - i_2 \|_{s_0 + \mu(\mathtt b_1)}\,. 
%\end{aligned} 
%\]

\medskip

\noindent
{\sc Estimate of $\mathcal{G}$ in \eqref{sigaro1}}.
One clearly has 
\[
\widehat{\mathcal G}(\ell)_j^{j'} = \langle \ell \rangle^{\mathtt{b}_{1}}
\langle j \rangle^{ \mathtt c_1} g_{\ell j j'} 
\widehat{\mathtt R}_1(\ell)_j^{j'} \langle j' \rangle^{- \mathtt c_1}\,,
\]
and for any $k \in \N^{\nu + 1}$, $|k| \leq k_0$, 
\begin{equation}\label{mortal kombat omologico 1a}
\begin{aligned}
 |\partial_\lambda^k \widehat{\mathcal G}(\ell)_j^{j'} |  
&\lesssim 
\sum_{k_1 + k_2 = k} \langle j \rangle^{ \mathtt c_1} |\partial_\lambda^{k_1}g_{\ell j j'}|  
\langle \ell \rangle^{\mathtt{b}_{1}}|\partial_\lambda^{k_2}\widehat{\mathtt R}_1(\ell)_j^{j'}| 
\langle j' \rangle^{- \mathtt c_1} 
\\& \hspace{-0.2cm}
\stackrel{\eqref{prop extension omologica}}{\lesssim}  
\langle \ell \rangle^{\tau_1} \langle j \rangle^{\tau_1 + \mathtt c_1 - \mathtt m}  
\langle j' \rangle^{\tau_1 }  \sum_{k_1 + k_2 = k}   \gamma^{- (|k_1| + 1)} 
\langle j \rangle^{\mathtt m} \langle \ell \rangle^{\mathtt{b}_{1}}  
|\partial_\lambda^{k_2}\widehat{\mathtt R}_1(\ell)_j^{j'}| \langle j' \rangle^{- \mathtt c_1}  \,. 
\end{aligned}
\end{equation}
The conservation of momentum $\mathtt V^T \ell + j - j' = 0$ implies that $|j - j'| \lesssim |\ell|$, 
therefore, using that 
$\langle j' \rangle \lesssim \langle j - j' \rangle + \langle j \rangle 
\lesssim \langle j - j' \rangle \langle j \rangle \lesssim \langle \ell \rangle \langle j\rangle$, 
one gets that 
\begin{equation}\label{mortal kombat omologico 2a}
\begin{aligned}
& \langle \ell \rangle^{\tau_1} \langle j \rangle^{\tau_1 + \mathtt c_1 - \mathtt m}  \langle j' \rangle^{\tau_1 } 
 \lesssim 
 \langle \ell \rangle^{2 \tau_1} \langle j \rangle^{2 \tau_1 + \mathtt c_1 - \mathtt m} 
  \lesssim  
  \langle \ell \rangle^{2 \tau_1}
 \end{aligned}
\end{equation}
since by \eqref{costanti.riducibilita} (recalling also the definition of $\tau_1$ 
in \eqref{prop extension omologica}), 
one has that $\mathtt m \geq 2 \tau_1 + \mathtt c_1$. 
Hence, by \eqref{mortal kombat omologico 1a}, \eqref{mortal kombat omologico 2a} 
and using that $|\ell| \leq K$, 
one gets for any $k \in \N^{\nu + 1}$, $|k| \leq k_0$, 
the estimate
\begin{equation}\label{mortal kombat omologico 3a}
\begin{aligned}
&|\partial_\lambda^k \widehat{\mathcal G}(\ell)_j^{j'} |  
\lesssim 
K^{2 \tau_1} \gamma^{- (|k| + 1)}
\sum_{k_1 + k_2 = k}   \gamma^{|k_2|} \langle j \rangle^{\mathtt m} 
\langle \ell \rangle^{\mathtt{b}_{1}} |\partial_\lambda^{k_2}\widehat{\mathtt R}_1(\ell)_j^{j'}| 
\langle j' \rangle^{- \mathtt c_1} \,. 
\end{aligned}
\end{equation}
Therefore, for any $k \in \N^{\nu + 1}$, $|k| \leq k_0$, $s_0 \leq s \leq \bar{s}$ 
one has 
\begin{equation*}
\begin{aligned}
& \| |\partial_\lambda^k {\mathcal G}| h \|_s^2  
\leq 
\sum_{\ell, j} \langle \ell, j \rangle^{2 s} \Big( \sum_{\ell', j'} |\partial_\lambda^k 
\widehat{\mathcal G}(\ell - \ell')_j^{j'}| |\widehat h(\ell', j')|\Big)^2  
\\& \hspace{-0.2cm}
\stackrel{\eqref{mortal kombat omologico 3a}}{\lesssim} 
K^{4 \tau_1} 
\gamma^{- 2 (|k| + 1)} \sum_{|k| \leq k_2} \gamma^{2 |k_2|}\sum_{\ell, j}   
\langle \ell, j \rangle^{2 s} \Big(  \sum_{\ell', j'} \langle j \rangle^{\mathtt m} 
\langle \ell - \ell'\rangle^{\mathtt{b}_{1}}  
|\partial_\lambda^{k_2}\widehat{\mathtt R}_1(\ell)_j^{j'}| 
\langle j' \rangle^{- \mathtt c_1}   |\widehat h(\ell', j')|\Big)^2 
\\& 
\lesssim 
K^{4 \tau_1} 
\gamma^{- 2 (|k| + 1)} \sum_{|k| \leq k_2} \gamma^{2 |k_2|} 
\| | \langle \partial_\vphi \rangle^{\mathtt{b}_{1}} \langle D \rangle^{\mathtt m} 
\partial_\lambda^{k_2} {\mathtt R}_1  \langle D \rangle^{- \mathtt c_1}| |h|\|^2_s 
\\& 
\lesssim
K^{4 \tau_1}\gamma^{- 2 (|k| + 1)} 
\Big( {\mathfrak M}^\sharp_{\langle \partial_\vphi \rangle^{\mathtt{b}_{1}} 
\langle D \rangle^{\mathtt m} {\mathtt R}_1 
\langle D \rangle^{- \mathtt c_1} }(s) \| |h| \|_{s_0} 
+ {\mathfrak M}^\sharp_{\langle \partial_\vphi \rangle^{\mathtt{b}_{1}} 
\langle D \rangle^{\mathtt m}  {\mathtt R}_1 \langle D \rangle^{- \mathtt c_1} }(s_0) 
\| |h| \|_{s} \Big)^2\,,
\end{aligned}
\end{equation*}
from which we deduce the bound 
\eqref{stima Pn s lemma homo} for the operator $\mathcal{G}$ in \eqref{sigaro1}.

\noindent
To estimate the Lipschitz variation
we consider the coefficients 
\[
\Delta_{12} \widehat{\mathcal G}(\ell)_j^{j'} 
= \widehat{\mathcal G}(\ell)_j^{j'}(\lambda, i_1(\lambda)) 
- \widehat{\mathcal G}(\ell)_j^{j'}(\lambda, i_2(\lambda))\,.
\] 
Using that $\e   \gamma^{-3-2(\mathtt{m}-\mathtt{c}_1)} \leq 1$ 
(recall the anstatz \eqref{smallcondepsilon})
one has that 
\begin{equation}\label{stima delta 12 mathcal G matrix}
\begin{aligned}
&|\Delta_{12} \widehat{\mathcal G}(\ell)_j^{j'}|    \leq 
|\Delta_{12} g_{\ell j j'}|  \langle \ell \rangle^{\mathtt{b}_{1}}\langle j \rangle^{ \mathtt c_1} 
|\widehat{\mathtt R}_1(\ell)_j^{j'}(\lambda, i_1(\lambda))| \langle j' \rangle^{- \mathtt c_1}  
\\& 
\quad 
+ |g_{\ell j j'}(\lambda, i_2(\lambda))| \langle \ell \rangle^{\mathtt{b}_{1}}\langle j \rangle^{ \mathtt c_1} 
|\Delta_{12} \widehat{\mathtt R}_1(\ell)_j^{j'}| \langle j' \rangle^{- \mathtt c_1}  
\\ 
& \hspace{-0.35cm} \stackrel{\eqref{prop extension omologicaBIS}}{\lesssim_{\bar{s}}}  
\frac{\e\langle \ell \rangle^{2 \tau} 
\langle j \rangle^{\mathtt c_1 - \mathtt m} }{\gamma^{4+2(\mathtt{m}-\mathtt{c}_1)}}
%\e   {\gamma^{-4-2(\mathtt{m}-\mathtt{c}_1)}} \langle \ell \rangle^{2 \tau} 
%\langle j \rangle^{\mathtt c_1 - \mathtt m}  
\langle j' \rangle^{2 \tau + 1 }  \langle \ell \rangle^{\mathtt{b}_{1}}\langle j \rangle^{ \mathtt m} 
|\widehat{\mathtt R}_1(\ell)_j^{j'}(\lambda, i_1(\lambda))| \langle j' \rangle^{- \mathtt c_1}  
\| i_1-i_2 \|_{ s_0  + \mu(\mathtt{b}_{1})}   
\\& 
\qquad 
+  \langle \ell \rangle^{\tau_1} |j|^{\tau_1 + \mathtt c_1 - \mathtt m} |j'|^{\tau_1 } \gamma^{-  1}  
\langle \ell \rangle^{\mathtt{b}_{1}}\langle j \rangle^{ \mathtt m} 
|\Delta_{12} \widehat{\mathtt R}_1(\ell)_j^{j'}| \langle j' \rangle^{- \mathtt c_1} 
\\& 
{\lesssim}  
\langle \ell \rangle^{\tau_1} |j|^{\tau_1 + \mathtt c_1 - \mathtt m} |j'|^{\tau_1 } \gamma^{-  1}  
\langle \ell \rangle^{\mathtt{b}_{1}}\langle j \rangle^{ \mathtt m} 
|\widehat{\mathtt R}_1(\ell)_j^{j'}(\lambda, i_1(\lambda))| \langle j' \rangle^{- \mathtt c_1}  
\| i_1-i_2 \|_{ s_0  + \mu(\mathtt{b}_{1})}  
\\& 
\qquad 
+  \langle \ell \rangle^{\tau_1} |j|^{\tau_1 + \mathtt c_1 - \mathtt m} |j'|^{\tau_1 }\gamma^{-  1}  
\langle \ell \rangle^{\mathtt{b}_{1}}\langle j \rangle^{ \mathtt m} 
|\Delta_{12} \widehat{\mathtt R}_1(\ell)_j^{j'}| \langle j' \rangle^{- \mathtt c_1} 
\\&  \hspace{-0.75cm}
\stackrel{\eqref{mortal kombat omologico 2a}, |\ell| \leq K}{\lesssim} 
K^{2 \tau_1} \gamma^{-  1}  \langle \ell \rangle^{\mathtt{b}_{1}}\langle j \rangle^{ \mathtt m} 
|\widehat{\mathtt R}_1(\ell)_j^{j'}(\lambda, i_1(\lambda))| \langle j' \rangle^{- \mathtt c_1}  
\| i_1-i_2 \|_{ s_0  + \mu(\mathtt{b}_{1})}  
\\& 
\qquad 
+ K^{2 \tau_1} \gamma^{-  1}  \langle \ell \rangle^{\mathtt{b}_{1}}\langle j \rangle^{ \mathtt m} 
|\Delta_{12} \widehat{\mathtt R}_1(\ell)_j^{j'}| \langle j' \rangle^{- \mathtt c_1}\,. 
\end{aligned}
\end{equation}
By \eqref{stima delta 12 mathcal G matrix}
one can check that $\| |\Delta_{12} {\mathcal G}| \|_{{\mathcal L}(H^{s_0})} $ satisfies the estimate
\eqref{stime delta 12 Psin}.
%
%
%Using that $\| |h| \|_s = \| h \|_s$, the latter inequality implies that 
%\[
%\begin{aligned}
%{\mathfrak M}^\sharp_{\mathcal G}(s) 
%& \lesssim 
%K_n^{2 \tau_1} 
%\gamma^{- 1}  {\mathfrak M}^\sharp_{\langle \partial_\vphi \rangle^{\mathtt{b}_{1}} 
%\langle D \rangle^{\mathtt m} {\mathcal R}_n \langle D \rangle^{- \mathtt c_1} }(s) 
% \\& 
% \stackrel{\eqref{stima cal R nu}}{\lesssim_{\bar{s}}} 
% K_n^{2 \tau_1} K_{n - 1} \e{\gamma^{-3-2(\mathtt{m}-\mathtt{c}_1)}}
% (1 + \| {\mathcal I}_{0} \|_{s + \mu(\mathtt{b}_{1})}^{k_0, \gamma})\,.
%\end{aligned}
%\]
%The corresponding estimate for $s = s_0$ follows by the latter one, 
%using also the ansatz \eqref{ansatz_I0_s0}. 
%Similarly, by using \eqref{stima delta 12 mathcal G matrix},  \eqref{stima cal R nu}, 
%\eqref{stima R nu i1 i2 norma alta} one obtains the estimate 
%\[
%\| |\Delta_{12} {\mathcal G}| \|_{{\mathcal L}(H^{s_0})} 
%\lesssim 
%K_n^{2 \tau_1 } K_{n - 1} \e{\gamma^{-3-2(\mathtt{m}-\mathtt{c}_1)}}  
%\| i_1 - i_2 \|_{s_0 + \mu(\mathtt b_1)}\, .
%\]
%
%\medskip
%
%\noindent
%The estimates of ${\mathfrak M}^\sharp_{\langle D \rangle^{\mathtt m} \Psi_ 
%\langle D \rangle^{- \mathtt m}}(s)\,,\, {\mathfrak M}^\sharp_{\langle D 
%\rangle^{\mathtt c_1} \Psi_ \langle D \rangle^{-\mathtt c_1}}(s)\,, \,  {\mathfrak M}^\sharp_{\Psi}(s)\,,\,  {\mathfrak M}^\sharp_{\langle \partial_\vphi \rangle^{\mathtt{b}_{1}} \Psi }(s)$ 
%can be done by similar arguments. 
%The proof of the claimed estimate is then concluded. 

Recalling Definition \ref{operatoreHam} and Lemma \ref{lemma:hamfou} and \eqref{shomo1}, it is easy to show that since ${\mathcal R}$
is a Hamiltonian operator, then also $\Psi$ is a Hamiltonian operator. 
Moreover, by the definition of $\Psi$, one can immediately recognize 
that it is a momentum preserving operator, since ${\mathcal R}$ is momentum preserving. 
This concludes the proof.
%The proof of the estimates \eqref{stime delta 12 Psin} can be done by very 
%similar arguments and therefore it is omitted. 
\end{proof}

Let us define the time one flow map 
\begin{equation}\label{defPsi}
\Phi_n:= \exp( \Psi_n ) =\mathcal{F}_1\,,\qquad 
 {\mathcal F}_\tau := {\rm exp}(\tau \Psi_n), \quad \tau \in [- 1, 1]\,.
\end{equation}
In the following we need the following smallness assumption on the remainder $\mathcal{R}_n$:
one has
\begin{equation}\label{Rnpiccolopiccolo}
C(\bar{s})K_{n}^{\tau_{3}}K_{n-1}\gamma^{-1}  {\mathfrak M}_{\langle D \rangle^{\mathtt m}   
{\mathcal R}_n 
\langle D \rangle^{- \mathtt{c}_1}}^\sharp ( s_0)\leq 1\,,
\end{equation}
for some $C(\bar{s})>0$ large enough.
The condition above will be verified inductively along the proof of Proposition 
\ref{iterazione riducibilita}.
We have the following.
\begin{lemma}{\bf (Estimates on the flow map).}\label{lem:on the flow map}
Under the assumptions \eqref{KAM smallness condition1}, \eqref{Rnpiccolopiccolo} one has the following.
The map $\mathcal{F}_{\tau}$, $\tau \in [- 1, 1]$ in \eqref{defPsi} satisfies,
 for any $s_0 \leq s \leq \bar{s}$, $\lambda = (\omega, \mathtt h) \in \mathtt \Omega \times [\mathtt h_1, \mathtt h_2]$ 
 and $\mathtt e = \mathtt m, \mathtt c_1$ the estimates
\begin{equation}\label{stima exp tau Psi n}
\begin{aligned}
{\mathfrak M}^\sharp_{\langle D \rangle^{\mathtt e} {\mathcal F}_\tau \langle D \rangle^{- \mathtt e}} (s) 
 &\lesssim
 1+
  K_n^{\tau_3}\gamma^{-1}
  {\mathfrak M}_{\langle D \rangle^{\mathtt m}   
{\mathcal R}_n 
\langle D \rangle^{- \mathtt{c}_1}}^\sharp ( s)\,,
\\
{\mathfrak M}^\sharp_{\langle \pa_{\vphi} \rangle^{\mathtt{b}_1}\langle D \rangle^{\mathtt e} 
{\mathcal F}_\tau \langle D \rangle^{- \mathtt e}} (s) 
 &\lesssim
 1+
  K_n^{\tau_3}\gamma^{-1}
  {\mathfrak M}_{\langle D \rangle^{\mathtt m}  \langle \partial_{\vphi } \rangle^{\mathtt{b}_{1}}  
{\mathcal R}_n 
\langle D \rangle^{- \mathtt{c}_1}}^\sharp ( s)
\\&\qquad
+  K^{2\tau_3}\gamma^{-2}
{\mathfrak M}_{\langle D \rangle^{\mathtt m}  \langle \partial_{\vphi } \rangle^{\mathtt{b}_{1}}  
{\mathcal R}_n 
\langle D \rangle^{- \mathtt{c}_1}}^\sharp ( s_0)
{\mathfrak M}_{\langle D \rangle^{\mathtt m}  
{\mathcal R}_n 
\langle D \rangle^{- \mathtt{c}_1}}^\sharp ( s)\,.
\end{aligned}
\end{equation}
Moreover, given 
$i_1(\lambda), i_2(\lambda)$ two tori satisfying \eqref{ansatz_I0_s0} one has
\begin{equation}\label{stima delta 12 exp Psi n}
\begin{aligned}
 \| | \langle \partial_\vphi \rangle^{p} \langle D \rangle^{ \mathtt e}
 \Delta_{12} &{\mathcal F}_\tau 
 \langle D \rangle^{- \mathtt e} |\|_{{\mathcal L}(H^{s_0})}  
\\& \lesssim
K_n^{ \tau_3} \gamma^{- 1}
\mathfrak{R}_{n}(s_0,p)
%\| |\langle \partial_\vphi \rangle^{p} 
%\langle D \rangle^{ \mathtt{m}}\mathcal{R}_n 
%\langle D \rangle^{- \mathtt c_1}| \|_{{\mathcal L}(H^{s_0})} 
\| i_1 - i_2 \|_{s_0 + \mu(\mathtt b_1)} 
\\& 
\qquad 
+ K_n^{ \tau_3} \gamma^{- 1} \| |\langle \partial_\vphi \rangle^{p} 
\Delta_{12}\langle D \rangle^{ \mathtt{m}} \mathcal{R}_n 
\langle D \rangle^{- \mathtt c_1}| \|_{{\mathcal L}(H^{s_0})} \,,
\end{aligned}
\end{equation}
with $p=0,\mathtt{b}_1$ and recalling \eqref{bassaGotR0}. The estimates above are uniform in $\tau\in[-1,1]$.
Finally, the map $\cF_{\tau}$ is symplectic and momentum preserving.
\end{lemma}

\begin{proof}
First of all, by the estimate \eqref{stima Pn s lemma homo} in 
Lemma \ref{stime eq omologica}, recalling the choice of the parameters in 
\eqref{costanti.riducibilita}, the ansatz \eqref{ansatz_I0_s0},
and the assumption 
\eqref{Rnpiccolopiccolo} one obtains that 
\begin{equation}\label{so Psi n Psi D m bla}
\begin{aligned}
& {\mathfrak M}^\sharp_{\langle D \rangle^{ \mathtt e} \Psi_n \langle D \rangle^{-  \mathtt e}}(s_0) 
\lesssim 
K_n^{\tau_3}  \gamma^{-1}
{\mathfrak M}_{\langle D \rangle^{\mathtt m}   
{\mathcal R}_n 
\langle D \rangle^{- \mathtt{c}_1}}^\sharp ( s_0)
\stackrel{\eqref{Rnpiccolopiccolo}}{\ll}1\,,
%K_{n - 1}^{- \ta_1} \e  {\gamma^{-3-2(\mathtt{m}-\mathtt{c}_1)}} 
%(1 + \| {\mathcal I }_{0} \|_{s_0 + \mu(\mathtt{b}_{1})}^{k_0, \gamma})
% \leq 1\,,
 \quad  \mathtt e = \mathtt m, \mathtt c_1\,.
\end{aligned}
\end{equation} 
%by the smallness condition \eqref{KAM smallness condition1}, 
%and by the inductive assumption \eqref{stima cal R nu} on $\mathcal{R}_n$,
%one obtains that 
%\begin{equation}\label{so Psi n Psi D m bla}
%\begin{aligned}
%& {\mathfrak M}^\sharp_{\langle D \rangle^{ \mathtt e} \Psi_n \langle D \rangle^{-  \mathtt e}}(s_0) 
%\leq C 
%K_n^{\tau_3}  K_{n - 1}^{- \ta_1} \e  {\gamma^{-3-2(\mathtt{m}-\mathtt{c}_1)}} 
%(1 + \| {\mathcal I }_{0} \|_{s_0 + \mu(\mathtt{b}_{1})}^{k_0, \gamma})
% \leq 1\,,\quad  \mathtt e = \mathtt m, \mathtt c_1\,.
%\end{aligned}
%\end{equation} 
Hence, the smallness condition \eqref{piccolezza neumann tamea} is verified. Therefore,
Lemma \ref{serie di neumann per maggiorantia} implies, 
for any $s_0 \leq s \leq \bar{s}$ and $e = \mathtt m, \mathtt c_1$, that
\begin{align}
{\mathfrak M}^\sharp_{\langle D \rangle^{\mathtt e} {\mathcal F}_\tau \langle D \rangle^{- \mathtt e}} (s) 
 &
\lesssim 1 + {\mathfrak M}^\sharp_{\langle D \rangle^{ \mathtt e} \Psi_n 
\langle D \rangle^{- \mathtt e}}(s)\,,
\label{stima exp tau Psi n1}
\\
{\mathfrak M}^\sharp_{\langle \partial_\vphi \rangle^{\mathtt{b}_{1}} (\langle D \rangle^{ \mathtt e} 
{\mathcal F}_\tau \langle D \rangle^{- \mathtt e})} (s) 
&
\lesssim 
1 + {\mathfrak M }_{ \langle \partial_{\vphi} \rangle^{\mathtt{b}_{1}} 
(\langle D \rangle^{ \mathtt e} \Psi_n \langle D \rangle^{ - \mathtt e} )}^\sharp(s) 
\nonumber
\\
&   
\quad+ {\mathfrak M}_{\langle D \rangle^{ \mathtt e}\Psi_n \langle D \rangle^{- \mathtt e}}^\sharp(s) \, 
{\mathfrak M }_{ \langle \partial_{\vphi} \rangle^{\mathtt{b}_{1}}  
(\langle D \rangle^{ \mathtt e}\Psi_n \langle D \rangle^{ - \mathtt e} )}^\sharp(s_0) \,.
\label{stima exp tau Psi n2}
\end{align}
The bounds \eqref{stima exp tau Psi n1}, \eqref{stima exp tau Psi n2}, together with 
 \eqref{stima Pn s lemma homo}, 
imply  \eqref{stima exp tau Psi n}.
Reasoning similarly, by applying Lemmata \ref{serie di neumann per maggiorantia}, 
\ref{stime eq omologica} (see \eqref{stime delta 12 Psin}), with the ansatz \eqref{ansatz_I0_s0},
by the smallness conditions \eqref{KAM smallness condition1}, \eqref{Rnpiccolopiccolo},  
one gets \eqref{stima delta 12 exp Psi n}.
Since $\Psi_{n}$ is Hamiltonian and momentum preserving, then $\cF_{\tau}$ is symplectic and momentum preserving.
 \end{proof}

\subsubsection{The new operator}\label{sec:inductionKAM}

We conjugate $ {\mathbf L}_n $ by the symplectic time one flow map in \eqref{defPsi}
generated by a Hamiltonian operator $ \Psi_n $ acting in $ {\bf H}^s_\bot $ (see \eqref{spaziograssetto})
given by Lemma \ref{stime eq omologica}.
By a Lie expansion we get 
\begin{equation}\label{Ltra1}
\begin{aligned}
 {\mathbf L}_{n + 1} 
 &:= 
 \Phi_n^{- 1}  {\mathbf L }_n \Phi_n = \Phi_n^{- 1} ( \omega \cdot \partial_\vphi 
 + \ii {\mathcal D}_n ) \Phi_n + \Phi_n^{- 1} {\mathcal R}_n \Phi_n  
 \\& 
 = \omega \cdot \partial_\vphi + \ii {\mathcal D}_n 
+  \omega \cdot \partial_\vphi \Psi_n 
 + \ii [{\mathcal D}_n, \Psi_n ] + {\mathcal R}_n 
 \\&
  + \int_0^1 {\rm exp}( - \tau \Psi_n)[{\mathcal R}_n, \Psi_n ] {\rm exp}( \tau \Psi_n)\, d \tau 
 \\& 
+ \int_0^1 (1 - \tau) {\rm exp}( - \tau \Psi_n) \Big[\omega \cdot \partial_\vphi \Psi_n 
 + \ii [{\mathcal D}_n, \Psi_n ]  , \Psi_n \Big] {\rm exp}( \tau \Psi_n)\, d \tau \,.
\end{aligned}
\end{equation}
In view of the homological equation \eqref{equazione omologica} (and recalling 
$ \Pi_{K_n} $ in \eqref{proiettore-oper}) we have that \eqref{Ltra1} becomes
\begin{equation}\label{forma cal P+ stima quadratica iterazione}
\begin{aligned}
 \bL_{n + 1} &= \Phi_n^{- 1} \bL_n \Phi_n 
 =  \omega \cdot \partial_\vphi + \ii {\mathcal D}_{n + 1}  + {\mathcal R}_{n + 1} \,, 
 \\
 {\mathcal D}_{n + 1} &:= {\mathcal D}_n -\ii [{\mathcal R}_n]\,, 
 \\
 {\mathcal R}_{n + 1} &:=
 %\stackrel{\eqref{equazione omologica}}{=}
 \Pi_{N_n}^\bot {\mathcal R}_n +  \int_0^1 ({\mathcal Q}_\tau + {\mathcal V}_\tau)\, d \tau \,, 
\end{aligned}
\end{equation}
where we defined, for $\tau\in[0,1]$,
\begin{equation}\label{cal R n + 1 iterazione 2}
\begin{aligned}
{\mathcal Q}_\tau & :=  {\rm exp}( - \tau \Psi_n)[{\mathcal R}_n, \Psi_n] {\rm exp}( \tau \Psi_n)\,, 
\\
{\mathcal V}_\tau & :=   (1 - \tau) {\rm exp}( - \tau \Psi_n) 
\Big[ [{\mathcal R}_n] - \Pi_{K_n} {\mathcal R}_n, \Psi \Big] {\rm exp}( \tau \Psi_n)\,. 
\end{aligned}
\end{equation}

In the next lemma we give some properties of the new normal form part ${\mathcal D}_{n + 1}$. 

\begin{lemma}{\bf (The new normal form).}\label{lemma nuova parte diagonale}
Under the assumptions \eqref{KAM smallness condition1}-\eqref{Rnpiccolopiccolo} one has the following.
The operator ${\mathcal D}_{n + 1}$ in \eqref{forma cal P+ stima quadratica iterazione} has the form 
\begin{equation}\label{cal D+ 0}
{\mathcal D}_{n + 1} = \begin{pmatrix}
\mathtt D_{n + 1} & 0 \\
0 & - \overline{\mathtt D_{n + 1}}
\end{pmatrix} \quad \text{where} 
\quad \mathtt D_{n + 1} := \mathtt D_n -\ii \widehat{\mathtt R}_{n, 1}(0)\,. 
\end{equation}
The operator $\mathtt D_{n+1}$ is diagonal, with real eigenvalues
\begin{equation}\label{cal D+ 1}
\begin{aligned}
& \mathtt D_{n + 1} := {\rm diag}_{j\in S_0^c} \mu_{n + 1}(j), 
\quad \mu_{n + 1}(j) := \mu_n(j) -\ii \widehat{\mathtt R}_{n, 1}(0)_j^j = \mu_0(j) + r_{n + 1}(j)\,,  
\\& 
r_{n + 1}(j) = r_n(j) -\ii \widehat{\mathtt R}_{n, 1}(0)_j^j \,, \quad  j\in S_0^c\,.  
\end{aligned}
\end{equation}
Moreover one has
\begin{equation}\label{cal D+ 1b}
\begin{aligned}
|\mu_{n + 1}(j) - \mu_n(j)|^{k_0, \gamma}  
& = |r_{n + 1}(j) - r_n(j)|^{k_0, \gamma}
 %\lesssim |j|^{- \mathtt m + \mathtt{b}_{1} + \mathtt k_0} {\mathfrak M}^\sharp_{\langle D \rangle^{\mathtt m} 
 %{\mathcal R}_n \langle D \rangle^{- \mathtt{b}_{1} - \mathtt k_0} }(s_0) 
 %\\&
 \lesssim |j|^{- \mathtt m + \mathtt c_1} 
  {\mathfrak M}^\sharp_{\langle D \rangle^{\mathtt m} 
{\mathcal R}_n \langle D \rangle^{- \mathtt c_1} }(s_0) \,.
% K_{n - 1}^{- \ta_1} 
%\e{\gamma^{-2-2(\mathtt{m}-\mathtt{c}_1)}} \,, 
% \\
%|r_{n + 1}(j)|^{k_0, \gamma} & \lesssim_{\bar{s}} 
%|j|^{- \mathtt m + \mathtt c_1} \e{\gamma^{-2-2(\mathtt{m}-\mathtt{c}_1)}}  \,. 
\end{aligned}
\end{equation}
Let $i_1(\lambda), i_2 (\lambda)$ be two tori satisfying \eqref{ansatz_I0_s0}. 
Then for all $\lambda = (\omega, \mathtt h) 
\in \mathtt\Omega \times [\mathtt h_1, \mathtt h_2]$
%with $\gamma_1, \gamma_2 \in [\gamma/2, 2 \gamma]$, 
the following estimates hold: 
\begin{equation}\label{Delta 12 r step n + 1}
\begin{aligned}
 \big|\Delta_{12}(r_{n + 1}(j) - r_n(j)) \big|  
&\lesssim
|j|^{- {\mathtt m} + \mathtt c_1} 
   \|| \langle D \rangle^{\mathtt m} \Delta_{12} {\mathcal R}_n \langle D \rangle^{-\mathtt c_1} | 
   \|_{{\mathcal L}(H^{s_0})}  \,.
%
%\e{\gamma^{-2-2(\mathtt{m}-\mathtt{c}_1)}}
%|j|^{- {\mathtt m} + \mathtt c_1}   K_{n - 1}^{- \ta_1}  
%\| i_1-i_2 \|_{ s_0  + \mu(\mathtt{b}_{1})}\,, 
%\\ 
%\ |\Delta_{12} r_{n + 1}(j)| & \lesssim_{\bar{s}}  
%\e{\gamma^{-2-2(\mathtt{m}-\mathtt{c}_1)}}  |j|^{- {\mathtt m} + \mathtt c_1}  
%\| i_1-i_2 \|_{ s_0  + \mu(\mathtt{b}_{1})}\,. 
\end{aligned}
\end{equation}
\end{lemma}

\begin{proof} 
In view of the assumptions \eqref{KAM smallness condition1}-\eqref{Rnpiccolopiccolo}
we can apply Lemma \ref{lem:on the flow map}. 
Hence the new operator ${\bf L}_{n+1}$ in \eqref{forma cal P+ stima quadratica iterazione}
is well-defined.
We start by prooving \eqref{cal D+ 1}, \eqref{cal D+ 1b}.
Since ${\mathcal R}_n$ in \eqref{forma cal R nu}
%\[
%{\mathcal R}_n = \begin{pmatrix}
%{\mathtt R}_{n, 1} & {\mathtt R}_{n, 2} \\
%\overline{\mathtt R_{n, 2}} & \overline{\mathtt R_{n, 1}}
%\end{pmatrix}
%\]
 is a Hamiltonian operator (recall Def. \ref{operatoreHam}), one gets that $-\ii  \mathtt R_{n, 1}$ is self-adjoint, 
 implying that $-\ii \widehat{\mathtt R}_{n, 1}(0)_j^j \in \R$ for any $j\in S_0^c$. 
 Moreover since ${\mathcal R}_n$ and then ${\mathtt R}_{n, 1}$ 
 is a momentum preserving operator (recall Remark \ref{verofood}), one has that 
 \[
 \widehat{\mathtt R}_{n, 1}(\ell)_j^{j'} \neq 0\quad \Rightarrow\;\; \mathtt V^T \ell + j - j' = 0
 \;\;\;{\rm then}\;\;\;\widehat{\mathtt R}_{n, 1}(0)_j^{j'} \neq 0
 \quad \Rightarrow \;\;\;j=j'\,,
 \]
% $\widehat{\mathtt R}_{n, 1}(\ell)_j^{j'} \neq 0$
%  implies that $\mathtt V^T \ell + j - j' = 0$ and hence $\widehat{\mathtt R}_{n, 1}(0)_j^{j'} \neq 0$ 
%  implies that $j = j'$, 
  meaning that $\widehat{\mathtt R}_{n, 1}(0)$ is the diagonal operator 
  $\widehat{\mathtt R}_{n, 1}(0) = {\rm diag}_{j\in S_0^c} \widehat{\mathtt R}_{n, 1}(0)_j^j$. 
  Hence by the above considerations and recalling \eqref{forma cal P+ stima quadratica iterazione}, 
  the properties \eqref{cal D+ 0}, \eqref{cal D+ 1} easily follow. 
  Moreover, we note that
\begin{equation*}
\begin{aligned}
& |\mu_{n + 1}(j) - \mu_n(j)|^{k_0, \gamma} = |r_{n + 1}(j) - r_n(j)|^{k_0, \gamma} 
= |\widehat{\mathtt R}_{n, 1}(0)_j^j|^{k_0, \gamma}  
\\& 
\qquad \leq 
|j|^{- \mathtt m + \mathtt c_1} \sup_{k \in S_0^c } 
|\widehat{\mathtt R}_{n, 1}(0)_k^k|^{k_0, \gamma} |k|^{\mathtt m - \mathtt c_1} 
\lesssim 
|j|^{- \mathtt m + \mathtt c_1} {\mathfrak M}^\sharp_{\langle D \rangle^{\mathtt m} 
{\mathcal R}_n \langle D \rangle^{- \mathtt c_1} }(s_0) \,,
%\\& 
%\qquad \stackrel{\eqref{stima cal R nu}, \eqref{ansatz_I0_s0}}{\lesssim_{\bar{s}}} 
% |j|^{- \mathtt m + \mathtt c_1} 
% K_{n - 1}^{- \ta_1} \e{\gamma^{-2-2(\mathtt{m}-\mathtt{c}_1)}} \,. 
\end{aligned}
\end{equation*}
and
\[
\begin{aligned}
  | \Delta_{12}(r_{n + 1}(j) - r_n(j))| & = |\Delta_{12}\widehat{\mathtt R}_{n, 1}(0)_j^j|    
  \leq 
  |j|^{- \mathtt m + \mathtt c_1} \sup_{k \in S_0^c } 
  |\Delta_{12}\widehat{\mathtt R}_{n, 1}(0)_k^k| |k|^{\mathtt m -  \mathtt c_1} 
   \\& 
   \lesssim 
   |j|^{- \mathtt m +  \mathtt c_1} 
   \|| \langle D \rangle^{\mathtt m} \Delta_{12} {\mathcal R}_n \langle D \rangle^{-\mathtt c_1} | 
   \|_{{\mathcal L}(H^{s_0})}  \,,
   \end{aligned}
\]
which imply \eqref{cal D+ 1b} and \eqref{Delta 12 r step n + 1}.
\end{proof}

\smallskip

We now estimate the new remainder ${\mathcal R}_{n + 1}$ given in 
\eqref{forma cal P+ stima quadratica iterazione}-\eqref{cal R n + 1 iterazione 2}. 

\begin{lemma}{\bf (The new remainder).}\label{lemma nuova parte resto}
Under the assumptions \eqref{KAM smallness condition1}-\eqref{Rnpiccolopiccolo} one has the following.
The remainder ${\mathcal R}_{n + 1}$ in 
\eqref{forma cal P+ stima quadratica iterazione}-\eqref{cal R n + 1 iterazione 2}
satisfies the following estimates: 
%\eqref{stima cal R nu}, \eqref{stima R nu i1 i2} and \eqref{stima R nu i1 i2 norma alta}
%with $n\rightsquigarrow n+1$.
\begin{equation}\label{stimaastraatta1}
\begin{aligned}
{\mathfrak M}^\sharp_{\langle D \rangle^{\mathtt m} {\mathcal R}_{n + 1} 
\langle D \rangle^{- \mathtt c_1} }(s) 
& \lesssim 
K_n^{- \mathtt{b}_{1}} 
{\mathfrak M}^\sharp_{
\langle D \rangle^{\mathtt m} 
\langle \partial_\vphi \rangle^{\mathtt{b}_{1}} 
{\mathcal R}_n\langle D \rangle^{- \mathtt c_1} }(s) 
\\& 
+ K_n^{\tau_3} \gamma^{- 1} {\mathfrak M}^\sharp_{\langle D \rangle^{\mathtt m} 
{\mathcal R}_n \langle D \rangle^{- \mathtt c_1} }(s) 
{\mathfrak M}^\sharp_{\langle D \rangle^{\mathtt m} {\mathcal R}_n \langle D \rangle^{- \mathtt c_1} }(s_0)\,, 
\\
{\mathfrak M}^\sharp_{
\langle D \rangle^{\mathtt m} \langle \partial_\vphi \rangle^{\mathtt{b}_{1}} 
 {\mathcal R}_{n + 1} \langle D \rangle^{- \mathtt c_1} } (s) 
& \lesssim 
{\mathfrak M}^\sharp_{
\langle D \rangle^{\mathtt m} 
\langle \partial_\vphi \rangle^{\mathtt{b}_{1}} 
{\mathcal R}_n \langle D \rangle^{- \mathtt c_1} }(s)  
\\& 
+ K_n^{\tau_3} \gamma^{- 1} {\mathfrak M}^\sharp_{\langle D \rangle^{\mathtt m} 
{\mathcal R}_n \langle D \rangle^{- \mathtt c_1} }(s)
{\mathfrak M}^\sharp_{
\langle D \rangle^{\mathtt m} 
\langle \partial_\vphi \rangle^{\mathtt{b}_{1}} 
{\mathcal R}_n \langle D \rangle^{- \mathtt c_1} }(s_0)\,. 
\end{aligned}
\end{equation}
and given 
$i_1(\lambda), i_2(\lambda)$ two tori satisfying \eqref{ansatz_I0_s0} one has
\begin{equation}\label{stimaastratta2}
\begin{aligned}
  \| | \langle D \rangle^{\mathtt m}& \Delta_{12} {\mathcal R}_{n+1}
\langle D \rangle^{- \mathtt c_1} | \|_{{\mathcal L}(H^{s_0})}
\lesssim
K_n^{- \mathtt{b}_{1}} \| |
\langle D \rangle^{\mathtt m}
\langle \partial_\vphi \rangle^{\mathtt{b}_{1}} 
\Delta_{12} 
 {\mathcal R}_n 
\langle D \rangle^{- \mathtt c_1} | \|_{{\mathcal L}(H^{s_0})} 
\\&+
K_n^{ \tau_3} \gamma^{- 1} 
(\mathfrak{R}_{n}(s_0,0))^2
 \| i_1 - i_2 \|_{s_0 + \mu(\mathtt b_1)} 
\\& 
+ K_n^{ \tau_3} \gamma^{- 1} 
\mathfrak{R}_{n}(s_0,0)
\| |
\langle D \rangle^{ \mathtt{m}}
\Delta_{12} \mathcal{R}_n 
\langle D \rangle^{- \mathtt c_1}| \|_{{\mathcal L}(H^{s_0})}\,,
\end{aligned}
\end{equation}

\begin{equation}\label{stimaastratta3}
\begin{aligned}
&  \| |
\langle D \rangle^{\mathtt m}
 \langle\pa_{\vphi}\rangle^{\mathtt{b}_1} \Delta_{12} 
{\mathcal R}_{n+1} 
\langle D \rangle^{- \mathtt c_1} | \|_{{\mathcal L}(H^{s_0})}
\lesssim
\|  |
\langle D \rangle^{\mathtt m}
\langle \partial_\vphi \rangle^{\mathtt{b}_{1}} \Delta_{12} 
  {\mathcal R}_n 
\langle D \rangle^{- \mathtt c_1} | \|_{{\mathcal L}(H^{s_0})}
\\&+
K_n^{ \tau_3} \gamma^{- 1} 
\mathfrak{R}_{n}(s_0,0)\mathfrak{R}_{n}(s_0,\mathtt{b_1})
\| i_1 - i_2 \|_{s_0 + \mu(\mathtt b_1)} 
\\& 
\qquad 
+ K_n^{ \tau_3} \gamma^{- 1} 
\mathfrak{R}_{n}(s_0,\mathtt{b}_1)
\| |
\langle D \rangle^{ \mathtt{m}}
\langle\pa_{\vphi}\rangle^{\mathtt{b}_1}
\Delta_{12} \mathcal{R}_n 
\langle D \rangle^{- \mathtt c_1}| \|_{{\mathcal L}(H^{s_0})}
\\& 
\qquad 
+ K_n^{ \tau_3} \gamma^{- 1} 
\mathfrak{R}_{n}(s_0,0)
\| |
\langle D \rangle^{ \mathtt{m}} 
\langle\pa_{\vphi}\rangle^{\mathtt{b}_1}\Delta_{12}
\mathcal{R}_n 
\langle D \rangle^{- \mathtt c_1}| \|_{{\mathcal L}(H^{s_0})}\,,
\end{aligned}
\end{equation}
where $\mathfrak{R}_{n}(s_0,q)$, $q=0,\mathtt{b}_1$, is defined in \eqref{bassaGotR0}.
%\begin{equation}\label{bassaGotR}
%\mathfrak{R}_{n}(s_0,q):=\sup_{j=1,2}\| |\langle\pa_{\vphi}\rangle^{q}
%\langle D \rangle^{ \mathtt{m}}\mathcal{R}_n(i_{j})
%\langle D \rangle^{- \mathtt c_1}| \|_{{\mathcal L}(H^{s_0})}\,,\qquad q=0,\mathtt{b}_1\,.
%\end{equation}
\end{lemma}

\begin{proof}
%\noindent
%{\sc Proof of \eqref{stima cal R nu} at the step $n + 1$.}
We estimate separately the terms $\Pi_{K_n}^\bot {\mathcal R}_n$ and the terms 
${\mathcal Q}_\tau, {\mathcal V}_\tau$ (uniformly w. r. to $\tau \in [0, 1]$) in 
\eqref{forma cal P+ stima quadratica iterazione}. 

\smallskip
\noindent
{\sc Estimate of $\Pi_{K_n}^\bot {\mathcal R}_n$.} 
A direct application of Lemma \ref{lemma:smoothing-tame} 
allows to deduce that for any $s_0 \leq s \leq \bar{s}$, 
\begin{equation}\label{stima Pi N bot cal Rn0}
\begin{aligned}
\mathfrak{M}^\sharp_{\Pi_{K_n}^\bot \langle D \rangle^{\mathtt m} 
{\mathcal R}_n \langle D \rangle^{- \mathtt c_1} }(s) 
& \leq 
K_n^{- \mathtt{b}_{1}} {\mathfrak M}^\sharp_{\langle D \rangle^{\mathtt m} 
\langle \partial_\vphi \rangle^{\mathtt{b}_{1}} {\mathcal R}_n \langle D \rangle^{- \mathtt c_1} }(s) \,, 
\\
\mathfrak{M}^\sharp_{\Pi_{K_n}^\bot \langle D \rangle^{\mathtt m} 
\langle \partial_\vphi \rangle^{\mathtt{b}_{1}}{\mathcal R}_n \langle D \rangle^{- \mathtt c_1} }(s) 
&\leq {\mathfrak M}^\sharp_{\langle D \rangle^{\mathtt m} 
\langle \partial_\vphi \rangle^{\mathtt{b}_{1}} {\mathcal R}_n \langle D \rangle^{- \mathtt c_1} }(s)\,.
\end{aligned}
\end{equation}

%
%\begin{equation}\label{stima Pi N bot cal Rn0}
%\begin{aligned}
%\mathfrak{M}^\sharp_{\Pi_{K_n}^\bot \langle D \rangle^{\mathtt m} 
%{\mathcal R}_n \langle D \rangle^{- \mathtt c_1} }(s) 
%& \leq 
%K_n^{- \mathtt{b}_{1}} {\mathfrak M}^\sharp_{\langle D \rangle^{\mathtt m} 
%\langle \partial_\vphi \rangle^{\mathtt{b}_{1}} {\mathcal R}_n \langle D \rangle^{- \mathtt c_1} }(s) \,, 
%\\&
%\stackrel{\eqref{stima cal R nu}_{n}}{\lesssim}
%K_n^{- \mathtt{b}_{1}}K_{n - 1}^{- \ta_{1}} \e 
%{\gamma^{-2-2(\mathtt{m}-\mathtt{c}_1)}} 
%(1 + \| {\mathcal I_{0} } \|_{s + \mu(\mathtt{b}_{1})}^{k_0, \gamma}) \,.
%\end{aligned}
%\end{equation}
%Similarly, using the second estimate in the inductive assumption \eqref{stima cal R nu}$_{n}$
%one gets
%\begin{equation}\label{stima Pi N bot cal Rn}
%\begin{aligned}
%\mathfrak{M}^\sharp_{\Pi_{K_n}^\bot \langle D \rangle^{\mathtt m} 
%\langle \partial_\vphi \rangle^{\mathtt{b}_{1}}{\mathcal R}_n \langle D \rangle^{- \mathtt c_1} }(s) 
%& 
%\lesssim
%K_{n - 1} \e{\gamma^{-2-2(\mathtt{m}-\mathtt{c}_1)}} 
%(1 + \| {\mathcal I_{0} } \|_{s + \mu(\mathtt{b}_{1})}^{k_0, \gamma})\,.
%\end{aligned}
%\end{equation}

\smallskip
\noindent
{\sc Estimate of ${\mathcal Q}_\tau, {\mathcal V}_\tau$.} 
Actually the terms ${\mathcal Q}_\tau$ and ${\mathcal V}_\tau$ can be estimated similarly, 
hence we only estimate ${\mathcal Q}_\tau$.  
We write (recall the notation in \eqref{defPsi}) 
\[
\begin{aligned}
{\mathcal Q}_\tau & =  {\mathcal F}_{- \tau}[{\mathcal R}_n, \Psi_n] {\mathcal F}_\tau  =  {\mathcal F}_{- \tau}{\mathcal R}_n \Psi_n  {\mathcal F}_\tau -  {\mathcal F}_{- \tau}\Psi_n {\mathcal R}_n {\mathcal F}_\tau 
\end{aligned}
\]
and hence 
\begin{equation}\label{moltiplica e dividi rangle langle}
\begin{aligned}
& \langle D \rangle^{\mathtt m}  {\mathcal Q}_\tau \langle D \rangle^{- \mathtt c_1} 
\\
&  
=  \big( \langle D \rangle^{\mathtt m}{\mathcal F}_{- \tau} \langle D \rangle^{- \mathtt m} \big)
\big(\langle D \rangle^{\mathtt m} {\mathcal R}_n \langle D \rangle^{- \mathtt c_1} \big)
 \big( \langle D \rangle^{\mathtt c_1} \Psi_n \langle D \rangle^{- \mathtt c_1}  \big)
  \big(  \langle D \rangle^{\mathtt c_1}{\mathcal F}_\tau  
  \langle D \rangle^{- \mathtt c_1} \big) 
   \\
&  
-  \big( \langle D \rangle^{\mathtt m}{\mathcal F}_{- \tau} \langle D \rangle^{- \mathtt m} \big) 
\big( \langle D \rangle^{\mathtt m} \Psi_n \langle D \rangle^{- \mathtt m} \big) 
\big( \langle D \rangle^{\mathtt m} {\mathcal R}_n \langle D \rangle^{- \mathtt c_1}  \big) 
\big( \langle D \rangle^{\mathtt c_1}  {\mathcal F}_\tau 
\langle D \rangle^{- \mathtt c_1} \big) \,.
\end{aligned}
\end{equation}
%Recall the smallness \eqref{Rnpiccolopiccolo}.
%First of all note that, by \eqref{costanti.riducibilita}, \eqref{KAM smallness condition1}, the anstatz
%\eqref{ansatz_I0_s0} and the inductive assumption \eqref{stima cal R nu},
%one has 
%\begin{equation}\label{tetraedro bla}
%K_n^{\tau_3} K_{n - 1} \gamma^{-1}{\mathfrak M}^\sharp_{\langle D \rangle^{\mathtt m} 
%{\mathcal R}_n \langle D \rangle^{- \mathtt{c}_{1} } }(s_0) 
%\stackrel{\eqref{ansatz_I0_s0}, \eqref{stima cal R nu}}{\lesssim} 
%K_n^{\tau_3} K_{n - 1}^{1 - \ta_1} \e{\gamma^{-3-2(\mathtt{m}-\mathtt{c}_1)}}    \ll 1\,,
%\end{equation} 
Moreover, using the composition  
Lemma \ref{interpolazione moduli parametri}
one gets
\[
\begin{aligned}
&
\mathfrak{M}^{\sharp}_{
\big( \langle D \rangle^{\mathtt{m}} \Psi_n \langle D \rangle^{- \mathtt{m}}  \big)
\big(\langle D \rangle^{\mathtt m} {\mathcal R}_n \langle D \rangle^{- \mathtt c_1} \big)}+
\mathfrak{M}^{\sharp}_{\big(\langle D \rangle^{\mathtt m} {\mathcal R}_n \langle D \rangle^{- \mathtt c_1} \big)
 \big( \langle D \rangle^{\mathtt c_1} \Psi_n \langle D \rangle^{- \mathtt c_1}  \big)}
 \\&
 \qquad \lesssim\sup_{\mathtt{e}=\mathtt{m},\mathtt{c}_1}\Big(
 {\mathfrak M}^\sharp_{\langle D \rangle^{\mathtt m} 
{\mathcal R}_n \langle D \rangle^{- \mathtt{c}_{1} } }(s_0) 
{\mathfrak M}^\sharp_{\langle D \rangle^{\mathtt{e}} 
\Psi_n \langle D \rangle^{-\mathtt{e} } }(s) +
 {\mathfrak M}^\sharp_{\langle D \rangle^{\mathtt m} 
{\mathcal R}_n \langle D \rangle^{- \mathtt{c}_{1} } }(s) 
{\mathfrak M}^\sharp_{\langle D \rangle^{\mathtt{e}} 
\Psi_n \langle D \rangle^{- \mathtt{e} } }(s_0) \Big)
\\&\qquad
\stackrel{\eqref{stima Pn s lemma homo}}{\lesssim}
K_{n}^{\tau_3}\gamma^{-1} {\mathfrak M}^\sharp_{\langle D \rangle^{\mathtt m} 
{\mathcal R}_n \langle D \rangle^{- \mathtt{c}_{1} } }(s_0) 
 {\mathfrak M}^\sharp_{\langle D \rangle^{\mathtt m} 
{\mathcal R}_n \langle D \rangle^{- \mathtt{c}_{1} } }(s) \,.
\end{aligned}
\]
By combining the estimate above with \eqref{stima exp tau Psi n} on the flow map 
and with \eqref{Rnpiccolopiccolo},
%\eqref{tetraedro bla}
one gets
\begin{equation}\label{stima cal Q tau}
{\mathfrak M}^\sharp_{\langle D \rangle^{\mathtt m} {\mathcal Q}_\tau 
\langle D \rangle^{- \mathtt c_1} }(s) 
 \lesssim 
K_n^{\tau_3} \gamma^{- 1} {\mathfrak M}^\sharp_{\langle D \rangle^{\mathtt m} 
{\mathcal R}_n \langle D \rangle^{- \mathtt c_1} }(s) 
{\mathfrak M}^\sharp_{\langle D \rangle^{\mathtt m} {\mathcal R}_n \langle D \rangle^{- \mathtt c_1} }(s_0)
\,.
\end{equation}
Reasoning similarly one obtains that
\begin{equation}\label{stima cal Q tauALTA}
\begin{aligned}
{\mathfrak M}^\sharp_{ 
\langle D \rangle^{\mathtt m}
\langle \partial_\vphi \rangle^{\mathtt{b}_{1}} 
 {\mathcal Q}_\tau \langle D \rangle^{- \mathtt c_1} }(s) 
& \lesssim 
{\mathfrak M}^\sharp_{
\langle D \rangle^{\mathtt m}
\langle \partial_\vphi \rangle^{\mathtt{b}_{1}} 
 {\mathcal R}_n \langle D \rangle^{- \mathtt c_1} }(s)   
\\
 + K_n^{\tau_3} \gamma^{- 1} &{\mathfrak M}^\sharp_{\langle D \rangle^{\mathtt m} 
{\mathcal R}_n \langle D \rangle^{- \mathtt c_1} }(s)
{\mathfrak M}^\sharp_{
\langle D \rangle^{\mathtt m}
\langle \partial_\vphi \rangle^{\mathtt{b}_{1}} 
 {\mathcal R}_n \langle D \rangle^{- \mathtt c_1} }(s_0) \,,
\end{aligned}
\end{equation}
and that the operator $\mathcal{V}_{\tau}$
satisfies similar estimates as \eqref{stima cal Q tau}-\eqref{stima cal Q tauALTA}.
By recalling the expression of ${\mathcal R}_{n + 1}$ in \eqref{cal R n + 1 iterazione 2},
and by collecting the 
estimates \eqref{stima Pi N bot cal Rn0}, 
\eqref{stima cal Q tau}, \eqref{stima cal Q tauALTA}
one deduces the bounds \eqref{stimaastraatta1} on ${\mathcal R}_{n + 1}$.

\medskip
\noindent
We now shortly describe how to estimate $\Delta_{12} {\mathcal R}_{n + 1}$ 
(recall the expression given in \eqref{forma cal P+ stima quadratica iterazione}-\eqref{cal R n + 1 iterazione 2}).
First of all by  Lemma \ref{lemma:smoothing-tame}, one has that 
\begin{equation}\label{Delta 12 rRn bot}
\begin{aligned}
\| |  \Pi_{K_n}^\bot \langle D \rangle^{\mathtt m}
\Delta_{12} {\mathcal R}_n 
&\langle D \rangle^{- \mathtt c_1} | \|_{{\mathcal L}(H^{s_0})} 
%\\& 
\lesssim 
K_n^{- \mathtt{b}_{1}} \| |
\langle D \rangle^{\mathtt m}
\langle \partial_\vphi \rangle^{\mathtt{b}_{1}} \Delta_{12} 
 {\mathcal R}_n \langle D \rangle^{- \mathtt c_1} 
| \|_{{\mathcal L}(H^{s_0})} 
\\
\| |
 \Pi_{K_n}^\bot \langle D \rangle^{\mathtt m} 
 \langle \partial_\vphi \rangle^{\mathtt{b}_{1}} 
\Delta_{12}
 {\mathcal R}_n 
&\langle D \rangle^{- \mathtt c_1} | \|_{{\mathcal L}(H^{s_0})} 
%\\&
 \lesssim 
\|  |
\langle D \rangle^{\mathtt m} 
\langle \partial_\vphi \rangle^{\mathtt{b}_{1}} \Delta_{12} 
 {\mathcal R}_n 
\langle D \rangle^{- \mathtt c_1} | \|_{{\mathcal L}(H^{s_0})} \,.
\end{aligned}
\end{equation}
Now, as done before, we provide an estimate for 
$\Delta_{12}\mathcal{Q}_{\tau}$, being the one for $\Delta_{12}\mathcal{V}_{\tau}$ similar.
By reasoning as in \eqref{moltiplica e dividi rangle langle},
using   repeatedly  the 
triangular inequality, estimates
\eqref{stime delta 12 Psin}, \eqref{stima delta 12 exp Psi n}, and 
\eqref{Rnpiccolopiccolo},
%\eqref{tetraedro bla},
one gets (recalling the notation in \eqref{bassaGotR0})
\begin{equation}\label{stime delta 12 cal Q tau cal V tau}
\begin{aligned}
  \| | \langle D \rangle^{\mathtt m} \Delta_{12} {\mathcal Q}_\tau 
&\langle D \rangle^{- \mathtt c_1} | \|_{{\mathcal L}(H^{s_0})}
%\\&
\lesssim
K_n^{ \tau_3} \gamma^{- 1} 
(\mathfrak{R}_{n}(s_0,0))^2
\| i_1 - i_2 \|_{s_0 + \mu(\mathtt b_1)} 
\\& 
\qquad 
+ K_n^{ \tau_3} \gamma^{- 1} 
\mathfrak{R}_{n}(s_0,0)
\| |
\langle D \rangle^{ \mathtt{m}} 
\Delta_{12}
\mathcal{R}_n 
\langle D \rangle^{- \mathtt c_1}| \|_{{\mathcal L}(H^{s_0})}\,,
\end{aligned}
\end{equation}
and
\begin{equation}\label{stime delta 12 cal Q tau cal V tauALTA}
\begin{aligned}
  \| |\langle D \rangle^{\mathtt m} 
  &\langle\pa_{\vphi}\rangle^{\mathtt{b}_1} \Delta_{12} {\mathcal Q}_\tau 
\langle D \rangle^{- \mathtt c_1} | \|_{{\mathcal L}(H^{s_0})}
\\&\lesssim
K_n^{ \tau_3} \gamma^{- 1} 
\mathfrak{R}_{n}(s_0,0)\mathfrak{R}_{n}(s_0,\mathtt{b}_1)
%\||
%\langle D \rangle^{ \mathtt{m}}\mathcal{R}_n 
%\langle D \rangle^{- \mathtt c_1}| \|_{{\mathcal L}(H^{s_0})}
%\| |
%\langle\pa_{\vphi}\rangle^{\mathtt{b}_1}
%\langle D \rangle^{ \mathtt{m}}\mathcal{R}_n 
%\langle D \rangle^{- \mathtt c_1}| \|_{{\mathcal L}(H^{s_0})} 
\| i_1 - i_2 \|_{s_0 + \mu(\mathtt b_1)} 
\\& 
\qquad 
+ K_n^{ \tau_3} \gamma^{- 1} 
\mathfrak{R}_{n}(s_0,
\mathtt{b}_1)
%\| |\langle\pa_{\vphi}\rangle^{\mathtt{b}_1}
%\langle D \rangle^{ \mathtt c_1}\mathcal{R}_n 
%\langle D \rangle^{- \mathtt c_1}| \|_{{\mathcal L}(H^{s_0})}
\| |
\langle D \rangle^{ \mathtt{m}} \Delta_{12}\mathcal{R}_n 
\langle D \rangle^{- \mathtt c_1}| \|_{{\mathcal L}(H^{s_0})}
\\& 
\qquad 
+ K_n^{ \tau_3} \gamma^{- 1} 
\mathfrak{R}_{n}(s_0,0)
%\| |
%\langle D \rangle^{ \mathtt{m}}\mathcal{R}_n 
%\langle D \rangle^{- \mathtt c_1}| \|_{{\mathcal L}(H^{s_0})}
\| |
\langle D \rangle^{ \mathtt{m}} 
\langle\pa_{\vphi}\rangle^{\mathtt{b}_1}
\Delta_{12}\mathcal{R}_n 
\langle D \rangle^{- \mathtt c_1}| \|_{{\mathcal L}(H^{s_0})}\,.
\end{aligned}
\end{equation}
By collecting \eqref{Delta 12 rRn bot}, \eqref{stime delta 12 cal Q tau cal V tau}, 
\eqref{stime delta 12 cal Q tau cal V tauALTA}
we deduce that $\mathcal{R}_{n+1}$ satisfies \eqref{stimaastratta2}-\eqref{stimaastratta3}.
\end{proof}

In view of Lemmata \ref{stime eq omologica}, \ref{lem:on the flow map}, \ref{lemma nuova parte diagonale} 
and \ref{lemma nuova parte resto}
we are in position to give the proof of the iterative lemma.

 \begin{proof}[{\bf Proof of Proposition \ref{iterazione riducibilita}}]
 We argue by induction.
 
 \smallskip
\noindent {\bf Initialization.} The items ${\bf(S1)}_{0}$, ${\bf(S2)}_{0}$ on the operator ${\bf L}_0$
in \eqref{operatorenonintero}
are proved separately.

\noindent
{\sc Proof of ${\bf(S1)}_{0}$}. Properties \eqref{cal L nu}-\eqref{stima rj nu} for $ n = 0 $
 follow by \eqref{proprieta cal D bot} with  $r_0(j) = 0$. 
Furthermore, \eqref{stima cal R nu} holds for $ n = 0 $ in view of Lemma \ref{inizializzazione modulo tame red},
%\begin{lemma}\label{lem: Initialization} 
%$ {\mathfrak M}_0^\sharp (s) $,  $ {\mathfrak M}_0^\sharp ( s, \mathtt{b}_{1}) \lesssim_{ \mathtt{b}_{1}}    
%{\mathfrak M}_0 (s, {\mathtt{b}_{1}})  $ where ${\mathfrak M}_0 (s, {\mathtt{b}_{1}})$ is defined in \eqref{stima mathfrak M s0 b}. 
%\end{lemma}
\noindent

\noindent
{\sc Proof of ${\bf(S2)}_0$}. Estimates \eqref{stima R nu i1 i2} and \eqref{stima R nu i1 i2 norma alta} follow similarly 
by \eqref{stima tame cal R7 new}. 
%\shu{se vuoi fare con dettagli vedi Berti Montalto standing waves}

\medskip

\noindent
{\bf Induction step.}
Assume that ${\bf (S1)_n}, {\bf (S2)_n}$ holds true. We shall prove ${\bf (S1)_{n + 1}}, {\bf (S2)_{n + 1}}$. 

\noindent
First of all, in view of the inductive assumption on $\mathcal{R}_n$ we have
\[
K_{n}^{\tau_{3}}K_{n-1}\gamma^{-1}  {\mathfrak M}_{\langle D \rangle^{\mathtt m}   
{\mathcal R}_n 
\langle D \rangle^{- \mathtt{c}_1}}^\sharp ( s_0)
\stackrel{\eqref{stima cal R nu}_{n}}{\leq} C_{*}(\bar{s})
K_{n}^{\tau_{3}}K_{n - 1}^{1- \ta_{1}} \e 
{\gamma^{-3-2(\mathtt{m}-\mathtt{c}_1)}} 
(1 + \| {\mathcal I_{0} } \|_{s_0 + \mu(\mathtt{b}_{1})}^{k_0, \gamma})\ll_{\bar{s}}1
\]
by using  \eqref{costanti.riducibilita}, the smallness \eqref{KAM smallness condition1} and  
the anstatz
\eqref{ansatz_I0_s0}. Then  the smallness \eqref{Rnpiccolopiccolo} is verified and  Lemmata
\ref{stime eq omologica}, \ref{lem:on the flow map}, \ref{lemma nuova parte diagonale}, 
\ref{lemma nuova parte resto}
apply.

\noindent
{\sc Proof of \eqref{tame Psi nu - 1}$_{n+1}$-\eqref{tame Psi nu - 1 vphi x b}$_{n+1}$.}
Consider the map $\Psi_{n}$ (which solves \eqref{equazione omologica}) 
given by  Lemma \ref{stime eq omologica}.
Therefore, the bounds 
 \eqref{tame Psi nu - 1}-\eqref{tame Psi nu - 1 vphi x b} at the step 
 $n + 1$ follow by estimates 
 \eqref{stima Pn s lemma homo}
 and the inductive assumption
 \eqref{stima cal R nu} on $\mathcal{R}_{n}$.
 
 We consider now the new operator ${\bf L}_{n+1}=\Phi_{n}^{-1}{\bf L}_{n}\Phi_{n}$
 in \eqref{forma cal P+ stima quadratica iterazione} where $\Phi_{n}$
 is the flow map given by Lemma \ref{lem:on the flow map}.
 
 \smallskip
 \noindent
 {\sc Proof of \eqref{mu j nu}$_{n+1}$-\eqref{vicinanza autovalori estesi}$_{n+1}$, 
 \eqref{r nu - 1 r nu i1 i2}$_{n+1}$-\eqref{r nu i1 - r nu i2}$_{n+1}$ }.
Consider  the new eigenvalues $\mu_{n+1}(j)$  given 
by 
 Lemma \ref{lemma nuova parte diagonale}. 
 The expansion \eqref{mu j nu}$_{n+1}$ follows by \eqref{cal D+ 1}.
 By estimates \eqref{cal D+ 1b}-\eqref{Delta 12 r step n + 1}
 and the inductive assumption we get
 \[
 \begin{aligned}
 |\mu_{n + 1}(j) - \mu_n(j)|^{k_0, \gamma}  
& = |r_{n + 1}(j) - r_n(j)|^{k_0, \gamma}
 \lesssim |j|^{- \mathtt m + \mathtt c_1} 
  {\mathfrak M}^\sharp_{\langle D \rangle^{\mathtt m} 
{\mathcal R}_n \langle D \rangle^{- \mathtt c_1} }(s_0)
\\&
\stackrel{\eqref{stima cal R nu}}{\lesssim_{\bar{s}}}
 |j|^{- \mathtt m + \mathtt c_1} 
 K_{n - 1}^{- \ta_{1}} \e 
{\gamma^{-2-2(\mathtt{m}-\mathtt{c}_1)}} 
(1 + \| {\mathcal I_{0} } \|_{s_{0} + \mu(\mathtt{b}_{1})}^{k_0, \gamma})
\\&\stackrel{\eqref{ansatz_I0_s0}}{\lesssim_{\bar{s}}}
 |j|^{- \mathtt m + \mathtt c_1} 
 K_{n - 1}^{- \ta_{1}} \e 
{\gamma^{-2-2(\mathtt{m}-\mathtt{c}_1)}} \,,
 \end{aligned}
 \]
 which is the \eqref{vicinanza autovalori estesi} at the step $n+1$.
 Similarly by \eqref{Delta 12 r step n + 1} and the inductive assumption ${\bf(S2)_{n}}$
 we get
 \begin{equation*}
\begin{aligned}
  | \Delta_{12}(r_{n + 1}(j) - r_n(j))| & = |\Delta_{12}\widehat{\mathtt R}_{n, 1}(0)_j^j|    
   \lesssim 
   |j|^{- \mathtt m +  \mathtt c_1} 
   \|| \langle D \rangle^{\mathtt m} \Delta_{12} {\mathcal R}_n \langle D \rangle^{-\mathtt c_1} | 
   \|_{{\mathcal L}(H^{s_0})}  
   \\&  
   \stackrel{\eqref{stima R nu i1 i2}}{\lesssim_{\bar{s}}}  
   |j|^{- \mathtt m + \mathtt c_1} K_{n - 1}^{- \ta_1} \e{\gamma^{-2-2(\mathtt{m}-\mathtt{c}_1)}}
   \| i_1-i_2\|_{s_0 + \mu(\mathtt{b}_{1})} \,,
\end{aligned}
\end{equation*}
which is the \eqref{r nu - 1 r nu i1 i2} at the step $n+1$.
Let us prove \eqref{stima rj nu}$_{n+1}$ and \eqref{r nu i1 - r nu i2}$_{n+1}$.
We use a telescopic argument.
Indeed $r_{n + 1}(j) = \sum_{i = 0}^n r_{i + 1}(j) - r_i(j)$  (recall that by definition $r_0(j) := 0$) and 
\[
\begin{aligned}
|r_{n + 1}(j)|^{k_0, \gamma} & \leq \sum_{i = 0}^n | r_{i + 1}(j) - r_i(j) |^{k_0, \gamma} 
\lesssim_{\bar{s}} \e{\gamma^{-2-2(\mathtt{m}-\mathtt{c}_1)}}  |j|^{- \mathtt m +\mathtt c_1}
\sum_{i = 0}^{+ \infty} K_{i - 1}^{- \ta_1} 
\\&
\lesssim_{\bar{s}} \e{\gamma^{-2-2(\mathtt{m}-\mathtt{c}_1)}} |j|^{- \mathtt m + \mathtt c_1}\,,
\end{aligned}
\]
since $\sum_{i \geq 0} K_{i - 1}^{- \ta_1}$ is convergent, see \eqref{costanti.riducibilita}. 
Similarly, using \eqref{r nu - 1 r nu i1 i2}$_{n+1}$ we have
\[
\begin{aligned}
|\Delta_{12} r_{n + 1}(j)| & \leq \sum_{i = 0}^n |\Delta_{12}(r_{i + 1}(j) - r_i(j))| 
\lesssim_{\bar{s}} 
\e{\gamma^{-2-2(\mathtt{m}-\mathtt{c}_1)}}  |j|^{- \mathtt m + \mathtt c_1} 
\sum_{i = 0}^{n}  K_{i - 1}^{- \ta_1} \| i_1-i_2\|_{s_0 + \mu(\mathtt{b}_{1})} 
\\& 
\lesssim_{\bar{s}} 
\e{\gamma^{-2-2(\mathtt{m}-\mathtt{c}_1)}}  |j|^{- \mathtt m + \mathtt c_1} 
\sum_{i \geq 0} K_{i - 1}^{- \ta_1} \| i_1-i_2\|_{s_0 + \mu(\mathtt{b}_{1})} 
\\& 
\lesssim_{\bar{s}} 
\e{\gamma^{-2-2(\mathtt{m}-\mathtt{c}_1)}}  |j|^{- \mathtt m +\mathtt c_1} 
\| i_1-i_2\|_{s_0 + \mu(\mathtt{b}_{1})}
\end{aligned}
\]
since, as we already observed, the series $\sum_{i \geq 0} K_{i - 1}^{- \ta_1}$ is convergent. 
 
 \smallskip
 \noindent
 {\sc Proof of \eqref{stima cal R nu}$_{n+1}$.} Let us estimate the remainder $\mathcal{R}_{n+1}$
given in Lemma \ref{lemma nuova parte resto}.
By estimates
\eqref{stimaastraatta1} and the inductive assumption \eqref{stima cal R nu} on $\mathcal{R}_{n}$
we get
\[
\begin{aligned}
&{\mathfrak M}^\sharp_{\langle D \rangle^{\mathtt m} {\mathcal R}_{n + 1} 
\langle D \rangle^{- \mathtt c_1} }(s) 
\\
& \qquad\qquad\leq 
C(\bar{s}) \e{\gamma^{-2-2(\mathtt{m}-\mathtt{c}_1)}} \Big( K_n^{- \mathtt{b}_{1}} K_{n - 1}   
+ K_n^{\tau_3} K_{n - 1}^{- 2 \ta_1} \e{\gamma^{-2-2(\mathtt{m}-\mathtt{c}_1)}}   \gamma^{- 1} \Big) 
(1 + \| {\mathcal I}_{0} \|_{s + \mu(\mathtt{b}_{1})}^{k_0, \gamma})  
\\
&\qquad\qquad \leq 
C_* K_n^{- \ta_1} \e{\gamma^{-2-2(\mathtt{m}-\mathtt{c}_1)}}
(1 + \| {\mathcal I}_{0} \|_{s + \mu(\mathtt{b}_{1})}^{k_0, \gamma})\,,
\end{aligned}
\]
provided that 
\begin{equation}\label{parametri conv iterazione step red}
C(\bar{s}) K_n^{- \mathtt{b}_{1}} K_{n - 1}  \leq \frac{C_*}{2} K_n^{- \ta_1} \,, 
\quad 
C(\bar{s}) K_n^{\tau_3} K_{n - 1}^{- 2 \ta_1} \e{\gamma^{-3-2(\mathtt{m}-\mathtt{c}_1)}}  
\leq 
\frac{C_*}{2} K_n^{- \ta_1}\,. 
\end{equation}
These two conditions are fullfilled by \eqref{costanti.riducibilita}, 
by the smallness condition \eqref{KAM smallness condition1} 
and by taking $C_* \equiv C_* (\bar{s}) \gg 0$ large enough. 
Similarly, using the second estimate in \eqref{stimaastraatta1}, 
one has that
\[
\begin{aligned}
&{\mathfrak M}^\sharp_{
\langle D \rangle^{\mathtt m}
\langle \partial_\vphi \rangle^{\mathtt{b}_{1}} 
 {\mathcal R}_{n + 1} \langle D \rangle^{- \mathtt c_1} } (s) 
\\ 
%\lesssim {\mathfrak M}^\sharp_{\langle \partial_\vphi \rangle^{\mathtt{b}_{1}} 
%\langle D \rangle^{\mathtt m} {\mathcal R}_n \langle D \rangle^{- \mathtt c_1} }(s)  
%\\
%&+ K_n^{\tau_3} \gamma^{- 1} {\mathfrak M}^\sharp_{\langle D \rangle^{\mathtt m} 
%{\mathcal R}_n \langle D \rangle^{- \mathtt c_1}  }(s)
%{\mathfrak M}^\sharp_{\langle \partial_\vphi \rangle^{\mathtt{b}_{1}} 
%\langle D \rangle^{\mathtt m}  {\mathcal R}_n \langle D \rangle^{- \mathtt c_1}  }(s_0) 
%\\
&\qquad  \lesssim_{\bar{s}}  
 \e{\gamma^{-2-2(\mathtt{m}-\mathtt{c}_1)}}  \big( K_{n - 1} + K_n^{\tau_3} K_{n - 1}^{1 - \ta_1}  \e{\gamma^{-2-2(\mathtt{m}-\mathtt{c}_1)}} \gamma^{- 1} 
  \big) (1 + \| {\mathcal I}_{0} \|_{s + \mu(\mathtt{b}_{1})}^{k_0, \gamma})  
\\
&\qquad  \leq
C(\bar{s})  K_{n - 1}   \e{\gamma^{-2-2(\mathtt{m}-\mathtt{c}_1)}}
(1 + \| {\mathcal I}_{0} \|_{s + \mu(\mathtt{b}_{1})}^{k_0, \gamma}) 
\\
&\qquad \leq C_* K_n  \e{\gamma^{-2-2(\mathtt{m}-\mathtt{c}_1)}} 
(1 + \| {\mathcal I}_{0} \|_{s + \mu(\mathtt{b}_{1})}^{k_0, \gamma}) \,,
\end{aligned}
\]
by taking $C_* \equiv C_* (\bar{s}) \gg 0$ large enough.
This concludes the proof of the estimate \eqref{stima cal R nu} at the step $n + 1$. 

\smallskip
\noindent
{\sc Proof of \eqref{stima R nu i1 i2}$_{n+1}$, \eqref{stima R nu i1 i2 norma alta}$_{n+1}$. } 
By estimate \eqref{stimaastratta2}, the inductive assumption \eqref{stima cal R nu}, 
\eqref{stima R nu i1 i2}-\eqref{stima R nu i1 i2 norma alta}
at the step $n$, and \eqref{KAM smallness condition1},
we get 
\[
\begin{aligned}
  \| |  \Delta_{12}&\langle D \rangle^{\mathtt m} {\mathcal R}_{n+1}
\langle D \rangle^{- \mathtt c_1} | \|_{{\mathcal L}(H^{s_0})}
\\&\lesssim_{\bar{s}}
\big(
K_n^{- \mathtt{b}_{1}}K_{n-1}+
\e  {\gamma^{-3-2(\mathtt{m}-\mathtt{c}_1)}}
K_{n}^{\tau_3}K_{n-1}^{-2\mathtt{a}_1}  \big)
\e  {\gamma^{-2-2(\mathtt{m}-\mathtt{c}_1)}} K_{n - 1} \|i_1 - i_2\|_{ s_0 +  \mu(\mathtt{b}_{1})}\,,
\end{aligned}
\]
which implies \eqref{stima R nu i1 i2}$_{n+1}$ recalling also \eqref{costanti.riducibilita}.
Similarly using \eqref{stimaastratta3}, the inductive assumptions,
 one deduces
\[
\begin{aligned}
  \| | \langle D \rangle^{\mathtt m} 
  &\langle\pa_{\vphi}\rangle^{\mathtt{b}_1}\Delta_{12} 
{\mathcal R}_{n+1} 
\langle D \rangle^{- \mathtt c_1} | \|_{{\mathcal L}(H^{s_0})}
\\&\lesssim_{\bar{s}}
\e  {\gamma^{-2-2(\mathtt{m}-\mathtt{c}_1)}} \|i_1 - i_2\|_{ s_0 +  \mu(\mathtt{b}_{1})}
\big(
K_{n-1}+
\e  {\gamma^{-3-2(\mathtt{m}-\mathtt{c}_1)}}
K_{n}^{\tau_3}K_{n-1}^{1-\mathtt{a}_1}  \big)
\end{aligned}
\]
which implies \eqref{stima R nu i1 i2 norma alta} at the step $n+1$ recalling 
again \eqref{costanti.riducibilita} and \eqref{KAM smallness condition1}
(by taking $C_* \equiv C_*(\bar{s}) \gg 0$ eventually larger).
This proves
 ${\bf (S1)_{n + 1}, (S2)_{n + 1}}$ and concludes the proof of 
 Proposition \ref{iterazione riducibilita}.
 \end{proof}

\subsection{Convergence}
In this section we prove the convergence of the scheme given by Proposition 
\ref{iterazione riducibilita} and we provide the proof of the main result of 
Theorem \ref{Teorema di riducibilita}.
 
% In this section we prove that the KAM reducibility scheme for the operator ${\mathcal L}_{e}^{(2)}$, whose iterative step is described in Proposition \ref{prop riducibilita}, is convergent under the smallness condition \eqref{KAM smallness condition} with the final operator being diagonal with purely imaginary eigenvalues.
First of all we show that 
 the eigenvalues of the normal form operator $\mathcal{D}_{n}$ in \eqref{cal L nu} 
 are Cauchy sequences and so they converge to some limit.
 \begin{lemma}\label{lemma blocchi finali}
 For any $j \in S_0^c$, the sequence 
 $\{  \mu_n(j) \equiv \mu_n(j; \lambda, i_0(\lambda))\}_{n\in\N}$, 
 $\lambda \in \mathtt \Omega \times [\mathtt h_1, \mathtt h_2]$ 
 defined in \eqref{mu j nu} (recall also \eqref{proprieta cal D bot}),
 converges to some limit
\begin{equation}\label{def cal N infty nel lemma}
\begin{aligned}
\mu_\infty(j)&= \mu_\infty(j;\,\cdot\,):  \mathtt\Omega \times [\mathtt h_1, \mathtt h_2] \to \R \,,
\\
 \mu_\infty(j)  &\equiv \mu_\infty(j; \lambda, i_0(\lambda)) 
= \mu_0(j; \lambda , i_0(\lambda)) + r_\infty(j; \lambda, i_0(\lambda))\,, 
\end{aligned}
\end{equation}
 satisfying the following estimates % the following estimates hold: 
\begin{equation}\label{stime forma normale limite}
\begin{aligned}
  |  \mu_\infty(j) -  \mu_n(j) |^{k_0, \gamma} &= |  r_\infty(j) - r_n(j) |^{k_0, \gamma} 
\\& \lesssim_{\bar{s}}
 K_{n - 1}^{- \mathtt a_1} \e {\gamma^{-2-2(\mathtt{m}-\mathtt{c}_1)}}  |j|^{- {\mathtt m} + \mathtt{c}_{1} }  \,, 
 \\ 
 | r_\infty(j)|^{k_0, \gamma} 
& \lesssim_{\bar{s}}   \e {\gamma^{-2-2(\mathtt{m}-\mathtt{c}_1)}}  |j|^{- {\mathtt m} + \mathtt{c}_{1} }\,.
\end{aligned}
\end{equation}
Moreover, let $i_1(\lambda), i_2(\lambda)$, $\lambda \in \mathtt\Omega \times [\mathtt h_1, \mathtt h_2]$ 
be two tori satisfying \eqref{ansatz_I0_s0}. Then 
\begin{equation}\label{stime forma normale limitea}
 \begin{aligned}
 |\Delta_{12} (\mu_\infty(j) -  \mu_n(j))| &= |\Delta_{12}(r_\infty(j) - r_n(j))|  
 \\&
 \lesssim_{\bar{s}} 
 K_{n - 1}^{- \mathtt a_1} \e {\gamma^{-2-2(\mathtt{m}-\mathtt{c}_1)}}  |j|^{- {\mathtt m} + \mathtt{c}_{1} } 
 \| i_1 - i_2 \|_{s_0 + \mu(\mathtt b_1)}\,, 
 \\
 | \Delta_{12} r_\infty(j)| & \lesssim_{\bar{s}}   \e {\gamma^{-2-2(\mathtt{m}-\mathtt{c}_1)}}  
 |j|^{- {\mathtt m} + \mathtt{c}_{1} } \| i_1 - i_2 \|_{s_0 + \mu(\mathtt b_1)}\,. 
 \end{aligned}
\end{equation}
 \end{lemma}
 
% \begin{proof}
% 	By Proposition \ref{prop riducibilita}, in particular by \eqref{reversibility reality auto} and \eqref{lambdaestesi}, we have that the sequence $\{\widetilde \mu_n(j;\omega)\}_{n\in\N}\subset \im\,\R$ is Cauchy on the closed set $\tD\tC(2\gamma,\tau)$,  therefore it is convergent for any $ \omega \in \tD\tC(2\gamma,\tau)$. The estimates in \eqref{stime forma normale limite} follow then by a telescoping argument with the estimate \eqref{lambdaestesi}. 
% \end{proof}
\begin{proof}
For any torus $i_0(\lambda)$ satisfying \eqref{ansatz_I0_s0},
by using \eqref{vicinanza autovalori estesi} in Proposition \ref{iterazione riducibilita}
%By Proposition \ref{iterazione riducibilita} for any torus $i_0(\lambda)$ satisfying \eqref{ansatz_I0_s0} 
%and by
%%	\eqref{reversibility reality auto} and
%\eqref{vicinanza autovalori estesi}, 
we have that the sequence 
$\{ \mu_n(j; \lambda, i_0(\lambda))\}_{n\in\N}\subset \R$ is Cauchy on the  set 
$\mathtt \Omega \times [\mathtt h_1, \mathtt h_2]$,  therefore it is convergent. 
The estimates in \eqref{stime forma normale limite}, \eqref{stime forma normale limitea}
%	 , resp. the estimate in \eqref{auto.infty.delta12}, 
follow then by a telescoping argument, using the estimates 
\eqref{vicinanza autovalori estesi}, \eqref{r nu - 1 r nu i1 i2}.	 
\end{proof}

 We now define  the set $\Lambda_\infty^\gamma$ of the non-resonance conditions for the final eigenvalues  
as in \eqref{Cantor set} with $\mu_{\infty}(j)$ given by Lemma \ref{lemma blocchi finali}.
We show that $\Lambda_\infty^\gamma$ is included in the sets \eqref{Omega nu + 1 gamma}
at each step of the iteration.

 \begin{lemma}\label{prima inclusione cantor}
 We have $\Lambda_\infty^\gamma \subseteq \cap_{n \geq 0} \, \Lambda_n^\gamma$ 
(see \eqref{Cantor set} and \eqref{Omega nu + 1 gamma}).
 \end{lemma}
 
 \begin{proof}
 	%\red{DA RINCONTROLLARE TUTTO, USO DEL MOMENTO SOPRATTUTTO!}
 We prove by induction that $\Lambda_\infty^\gamma \subseteq \Lambda_n^\gamma$ 
 for any integer $n \geq 0$. 
 The statement is trivial for $n=0$,
 since 
 $\Lambda_0^\gamma := \tD\tC (2 \gamma,\tau) \cap {\mathtt T}{\mathtt C}_\infty( \gamma, \tau) $ 
 (see Proposition \ref{iterazione riducibilita}).
We now assume by induction that 
$\Lambda_\infty^\gamma \subseteq \Lambda_n^\gamma$ for some $n \geq 0$ 
and we  show that $\Lambda_\infty^\gamma \subseteq \Lambda_{n + 1}^\gamma$. 
Let $\Lambda \in \Lambda_\infty^\gamma$, $\ell \in \Z^\nu \setminus \{ 0 \}$, 
$j, j' \in\Gamma^* \setminus \{ 0 \}$, with $ \mathtt V^T \ell +j-j'=0$ and $|\ell|  \leq K_n$. 
By \eqref{stime forma normale limite}, \eqref{Cantor set}, we compute
 \[
 \begin{aligned}
 | \,\omega \cdot \ell + \mu_n(j) - \mu_n(j')|   &\geq 
 | \omega \cdot \ell + \mu_\infty(j) - \mu_\infty(j')| 
 \\
 &\qquad 
 - |\mu_\infty(j) - \mu_n(j)|  - |\mu_\infty(j') - \mu_n(j')| 
 \\& 
 \geq \frac{2 \gamma}{\langle \ell \rangle^\tau | j' |^\tau } 
 - C K_{n - 1}^{- \mathtt a_1} \e {\gamma^{-2-2(\mathtt{m}-\mathtt{c}_1)}}
 \big(  |j|^{- {\mathtt m} + \mathtt{c}_{1}} + |j'|^{- {\mathtt m} + \mathtt{c}_{1}} \big)  
\\& 
\geq \frac{\gamma}{\langle \ell \rangle^\tau | j' |^\tau } \,,
\end{aligned}
\]
 for some positive constant $C>0$, provided
 \begin{equation}\label{frittata di maccheroni 10}
 C  K_{n - 1}^{- \ta_1}\e {\gamma^{-3-2(\mathtt{m}-\mathtt{c}_1)}} 
 \langle \ell \rangle^\tau |j'|^\tau 
 \big(  |j|^{- {\mathtt m} + \mathtt{c}_{1}} + |j'|^{- {\mathtt m} + \mathtt{c}_{1}} \big) \leq 1\,.
 \end{equation}
 Recalling that $\mathtt{m} >  \tau+\mathtt{c}_1$ by \eqref{costanti.riducibilita}, 
we have that $|j'|^{\tau +\mathtt{c}_1- \mathtt{m}} \leq 1$. 
Moreover, using that $ |\ell| \leq K_n$ and the momentum condition 
$\mathtt V^T  + j-j'=0$, we have that
 	\begin{equation*}
 		\begin{aligned}
 			& |j'| \leq |j| + |j - j'| \leq |j| + K_n \lesssim K_n |j|\,, 
 		\end{aligned}
 	\end{equation*}
	and therefore, using again that $\mathtt m > \tau+\mathtt{c}_1$, we have 
	\[
	| j' |^\tau |j|^{- \mathtt{m}+\mathtt{c}_1}  
	\lesssim K_n^\tau |j|^{\tau+\mathtt{c}_1 - \mathtt{m}} 
	\lesssim K_n^\tau\,.
	\]
We then deduce that 
\[
C  K_{n - 1}^{- \ta_1}\e {\gamma^{-3-2(\mathtt{m}-\mathtt{c}_1)}} 
\langle \ell \rangle^\tau |j'|^\tau 
\big(  |j|^{- {\mathtt m} + \mathtt{c}_{1}} + |j'|^{- {\mathtt m} + \mathtt{c}_{1}} \big)
\leq 
C_0 K_n^{2 \tau} K_{n - 1}^{- \ta_1} \varepsilon \gamma^{-3-2(\mathtt{m}-\mathtt{c}_1)}
\]
for some large constant $C_0 \gg 0$. 
Hence, the condition \eqref{frittata di maccheroni 10} is verified since 
\[
C_0 K_n^{2 \tau} K_{n - 1}^{- \ta_1} \varepsilon \gamma^{-3-2(\mathtt{m}-\mathtt{c}_1)} \leq 1\,,
\]
recalling \eqref{costanti.riducibilita}, 
and the smallness condition \eqref{KAM smallness condition1}. 
The proof for the Second Melnikov conditions with $+$ sign can be done similarly.
 Thus, we conclude that $\omega \in \Lambda_{n + 1}^\gamma$ and the claim is proved. 
 \end{proof}

We are in position to prove our main result.

\begin{proof}[{\bf Proof of Theorem \ref{Teorema di riducibilita}}]
By applying Proposition \ref{iterazione riducibilita}
we define the sequence of invertible maps 
 \begin{equation}\label{defUn}
 	{\mathcal U}_n := \Phi_0 \circ \Phi_1 \circ \ldots \circ \Phi_n \,, \quad n \in \N\,,
 \end{equation}
 defined for $(\omega, \mathtt h) \in \R^\nu \times [\mathtt h_1, \mathtt h_2] $.
 By estimates \eqref{tame Psi nu - 1}-\eqref{tame Psi nu - 1 vphi x b}
 and reasoning as in Theorem $7.5$ in \cite{BM20}
 one can check that, for any $n\geq1$
 \[
 \begin{aligned}
    {\mathfrak M}_{{\mathcal U}_n}^\sharp (s) &\lesssim_{\bar{s}}
    1+K_0^{\tau_3}\e  {\gamma^{-3-2(\mathtt{m}-\mathtt{c}_1)}} 
    (1 + \| {\mathcal I }_{0} \|_{s + \mu(\mathtt{b}_{1})}^{k_0, \gamma}) \, , 
 \\
   {\mathfrak M}_{{\mathcal U}_n-{\mathcal U}_{n-1}}^\sharp (s) 
&\lesssim_{\bar{s}} 
K_{n - 1}^{\tau_3} K_{n - 2}^{- \mathtt a_{1}} 
\e  {\gamma^{-3-2(\mathtt{m}-\mathtt{c}_1)}} 
(1 + \| {\mathcal I }_{0} \|_{s + \mu(\mathtt{b}_{1})}^{k_0, \gamma}) \, , 
\end{aligned}
 \]
from which one can deduce that there exists 
$\mathcal{U}_{\infty}=\lim_{n\to\infty}\mathcal{U}_{n}$
in the topology induced by the operator norm. In particular, 
\eqref{stima Phi infinito} holds for $\mathcal{U}_{\infty}$.
The estimate on the inverse $\mathcal{U}_{\infty}^{-1}$ by a Neumann series argument.

\noindent
Moreover, by the inclusion given by Lemma \ref{prima inclusione cantor},
and using \eqref{coniugionu+1}, one deduces, for any $n\geq1$,
the conjugation $\mathbf{L}_{n}=\mathcal{U}_{n}^{-1}\mathbf{L}_0\mathcal{U}_n$
where $\mathbf{L}_{n}$ given  in \eqref{cal L nu}.
Passing to the limit one gets the conjugated operator in \eqref{cal L infinito}.
The bounds \eqref{stimaAutovaloriFinalissimi1}-\eqref{stimaAutovaloriFinalissimi2} follow
by Lemma \ref{lemma blocchi finali}.
\end{proof}

\subsection{Inversion of the linearized operator
}\label{quasi invertibilita}
%Let us fix any $\bar{\tn}\in\N$. Take parameters as in \eqref{tbta} and \eqref{costanti.riducibilita}.
%Assume \eqref{ansatz_I0_s0} and \eqref{ps0} with $\mu_{0} = \mu(\tb_{1})$ given in \eqref{definizione.mu.b}.  
%Let $N_{0}$ % \ge \max\{ N_{0}^{(1)}\,, N_{0}^{(2)} \}$ where $N_{0}^{(1)}$ 
%bigger that the
%one appearing 
%in the smallness condition of Proposition \ref{riduzione trasporto} %and $N_{0}^{(2)}$ 
%and the one appearing in the smallness condition of Proposition \ref{iterazione riducibilita}. 
%%
%Assume 
%\[
%N_{0}^{\overline{\tau}} \eps \gamma^{-2N-3} \mathtt{C}(\bar{s}, \tb) \le 1
%\]
%for $\mathtt{C}(\bar{s}, \tb)$ large enough and $\overline{\tau} \ge \max \{ \tau_{3}, \tau_{2} \}$ such that the smallness conditions \eqref{smallness.raddrizzo} and \eqref{KAM smallness condition1} hold.
%
In this section we conclude the proof of the invertibility assumption \textbf{(AI)} (see \eqref{almi4}).
Recall the linearized operator $\mathcal{L}_{\omega}$ in \eqref{representation Lom}.
By applying  Theorem \ref{red.Lomega.smooth.rem}, one has that 
and 
 \[
 \cL_{\oo}= \cK \cL_{7} \cK^{-1} \stackrel{\eqref{operatorenonintero}}{=}
  \cK {\bf L}_0 \cK^{-1}\,, %+ \cK \cP_{\bot, \bar{\tn}} \cK^{-1}\,, 
  \qquad \forall  \lambda = (\oo, \th) \in \tT\tC_{\infty}(\gamma,\tau) 
 \]
 where  $\cL_{7}$ is defined in \eqref{op cal L7} and 
 $\cK$ is defined in \eqref{def.calkappa}. 
Then by  applying   Theorem \ref{Teorema di riducibilita}
%\red{We will prove that \textbf{(AI)} 
%(see \eqref{Lomega}-\eqref{almi3}) hold with $\ta =\ta_{1} $ (recall  \eqref{costanti.riducibilita}).}
(see also Lemma \ref{prima inclusione cantor}) for any $\lambda\in  \Lambda_{\infty}^\gamma \subset \tT\tC_{\infty}(\gamma,\tau)  $, we have that there exists an invertible, 
real, symplectic and momentum preserving map $\cU_{\infty}$ such that
 $\bL_{0}=\cU_{\infty} \bL_{\infty} \cU_{\infty}^{-1}$
where $\bL_{0}$ is defined in \eqref{operatorenonintero}
and $ \bL_{\infty}$ in \eqref{cal L infinito}.
In conclusion we have obtained
\begin{equation}\label{def.Psi.barn}
\begin{aligned}
\cL_{\oo} &  = {\bf \Phi}_{\infty} \bL_{\infty}  {\bf \Phi}_{\infty}^{-1} \,, 
\qquad {\bf \Phi}_{\infty}:= \cK \cU_{\infty}\,. 
\end{aligned}
\end{equation}
%\\
%Thus, since $\bL_{0}$ is obtained by $\cL_{7}$ neglecting $\cP_{\perp, \bar{\tn}}$ (see \eqref{operatoreintero}), we have that for any $\lambda = (\omega, \mathtt h) \in  \Lambda_{\bar{\tn}}^\gamma$, 
%\begin{equation}\label{def.Psi.barn}
%\begin{aligned}
%\cL_{\oo} & = \cK {\bf L}_0 \cK^{-1} + \cK \cP_{\bot, \bar{\tn}} \cK^{-1} \\
%& = \cK \cU_{\bar{\tn}} \bL_{\bar{\tn}} \cU_{\bar{\tn}}^{-1} \cK^{-1} + \cK \cP_{\bot, \bar{\tn}} \cK^{-1} \\
%& = {\bf \Phi}_{\bar{\tn}} (\bL_{\bar{\tn}} + \cU_{\bar{\tn}}^{-1} \cP_{\perp, \bar{\tn}}\cU_{\bar{\tn}} ) {\bf \Phi}_{\bar{\tn}}^{-1} \,, \qquad {\bf \Phi}_{\bar{\tn}}:= \cK \cU_{\bar{\tn}} 
%\end{aligned}
%\end{equation}
%By Lemmata \ref{lemma operatore e funzioni dipendenti da parametro} and \ref{A versus |A|},
% by estimate \eqref{stima Phi infinito} with the smallness condition \eqref{KAM smallness condition1} 
% %and $\tau_{3}> \tau_{0}$ (see Proposition \ref{iterazione riducibilita})
%the operators $ \cU_{\bar\tn}^{\pm 1}$ satisfy, for all $s_0 \leq s \leq \bar{s}$, 
%$\|  \cU_{\bar \tn}^{\pm 1} h\|_s^{k_0, \gamma} \lesssim_{\bar s} \| h \|_s^{k_0, \gamma} 
%+ \| \mathcal{I}_{0}\|_{s + \mu(\mathtt{b}_1)}^{k_0, \gamma} \| h \|_{s_0 }^{k_0, \gamma}$. 
Therefore, by Lemmata \ref{lemma operatore e funzioni dipendenti da parametro} and \ref{A versus |A|},
 recalling \eqref{stima Phi infinito}, \eqref{stimamappacomplessiva} and
\eqref{KAM smallness condition1} (recall that $K_0 $ is independent of $\gamma$), one gets that the
operators ${\bf \Phi}_{\infty}^{\pm 1}$ satisfy, for all $s_0 \leq s \leq \bar{s}$,
\begin{equation}\label{stime.Psi.barn}
\| {\bf \Phi}_{\infty}^{\pm 1} h \|_s^{k_0, \gamma} \lesssim_{\bar{s}} \| h \|_{s + \sigma}^{k_0, \gamma}  + 
\|  \mathcal{I}_{0} \|_{s + \mu(\mathtt{b}_1)}^{k_0, \gamma} \| h \|_{s_0 + \sigma}^{k_0, \gamma} \,,
\end{equation}
for some $\sigma = \sigma(k_0, \tau, \nu) > 0$. 

\smallskip
\noindent
If $\lambda=(\oo, \th)$ satisfies some first order Melnikov non-resonance conditions, namely if $\lambda$ belongs to the set
\begin{equation}\label{prime.di.Melnikov}
\begin{aligned}
\Lambda^{\gamma, I}_{\infty}:= \Lambda^{\gamma, I}_{\infty}(i_{0})&:= 
\big\{\lambda \in \mathtt \Omega \times [\th_{1}, \th_{2} ] : | \oo\cdot \ell + \mu_{\infty}(j)| 
\ge 2\gamma \jap{\ell}^{-\tau} |j|^{- \tau}\,, 
\\&\qquad\qquad\qquad\qquad  \ell\in\Z^{\nu}\setminus\{0\}, j \in S_{0}^{c}\,, \tV^{T}\ell + j= 0\big\}\,,
\end{aligned}
\end{equation}
then the operator ${\mathbf L}_{\infty}$ is invertible with suitable estimates. 
This is the content of the following lemma.

\begin{lemma}%[First order Melnikov non-resonance conditions]
For all $\lambda\in \Lambda^{\gamma, I}_{\infty}$ 
the operator ${\mathbf L}_{\infty}$ in \eqref{cal L infinito} is invertible 
on the subspace of the traveling waves $\tau_\vs  h(\vphi,x)=u(\vphi,x+\vs)$, $\vs\in\R$, such that $h\in {\bf H}_{S}^{\perp}$
and there is an extension of the 
inverse operator (that we denote in the same way) to the whole $\mathtt \Omega  \times [\mathtt h_1, \mathtt h_2] $ 
satisfying the estimate
\begin{equation}\label{stima.inv.bL}
\| {\mathbf L}_{\infty}^{-1} h \|^{k_{0}, \gamma}_{s} \lesssim \gamma^{- 1} \| h \|^{k_{0}, \gamma}_{s+ \s}
\end{equation}
where $\s= k_{0} + \tau(k_{0} + 1)$ is the constant in \eqref{Diophantine-1} with $k_{0} = k+1$. 
Moreover, ${\mathbf L}_{\infty}^{-1} $ is a momentum preserving operator.
\end{lemma}

\begin{proof}
By \eqref{stima tipo simbolo mu 0 j j'}, \eqref{stime forma normale limite}, using the smallness condition \eqref{KAM smallness condition1}, arguing as in \eqref{stima.f.j}, one has 
$| \partial_\lambda^k( \oo \cdot \ell + \mu_{\infty}(j) ) | 
\lesssim  \langle \ell \rangle |j|$ for all $1 \leq |k| \leq k_0$.
Hence, the thesis follows by applying Lemma 
%$B.4$ in \cite{BBHM}
\ref{lemma:cut-off sd} 
 to the function $f_{\ell j}(\lambda) = \oo \cdot \ell + \mu_{\infty}(j ; \lambda) $
with $M = C \langle \ell \rangle |j|$ and 
$\rho =  \frac{\gamma}{2} \langle \ell  \rangle^{- \tau}$, $\mathtt \Lambda_0 = \mathtt \Omega \times [\mathtt h_1, \mathtt h_2]$ and we obtain \eqref{stima.inv.bL}.  
\end{proof}
%
%Standard smoothing properties imply that 
%the operator $\cR_{\bar{\tn}}^\bot $ defined in \eqref{suddivisione} satisfies for any traveling wave $h\in H_{S}^{\perp}$, for all $ b  > 0$,   
%\begin{equation}\label{stima.RnPerp}
%\| \cR_{\bar{\tn}}^\bot h \|_{s_0}^{k_0, \gamma} \lesssim K_n^{- b} \| h \|_{s_0 + b + 1}^{k_0, \gamma}\,,\quad \| \cR_{\bar{\tn}}^\bot h\|_s^{k_0, \gamma} \lesssim \| h \|_{s + 1}^{k_0, \gamma} \, .
%\end{equation}
By \eqref{def.Psi.barn}, 
Theorem \ref{Teorema di riducibilita},  \eqref{operatoreintero}, 
and estimates \eqref{stima.inv.bL}, %\eqref{stima.RnPerp}, 
\eqref{stime.Psi.barn}, 
we deduce the following theorem.

\begin{thm}\label{inversione parziale cal L omega}
{\bf (Invertibility of $ \cL_\oo $)}
Assume \eqref{ansatz}. 
Recall the parameters in \eqref{tbta}, \eqref{costanti.riducibilita}. 
Let $\bar{s} > s_0$, and assume the smallness condition \eqref{KAM smallness condition1}.
Then for all 
\begin{equation}\label{Melnikov-invert}
(\omega, \th) \in  {\bf \Lambda}_{\infty}^{\g}  := {\bf \Lambda}_{\infty}^{\g} (i) 
:= { \Lambda}_{\infty}^\gamma  \cap  { \Lambda}_{\infty}^{\gamma, I}
\end{equation}
(see \eqref{Cantor set}, \eqref{prime.di.Melnikov}) 
the operator $ \cL_\omega$ defined in \eqref{Lomegatrue} (see also \eqref{representation Lom})
is invertible and satisfies the 
estimate \eqref{almi4}
for some 
$\sigma := \sigma(k_0, \tau, \nu) > 0$ and with 
with $ \overline \mu  \equiv \mu (\mathtt b_1) $  defined in \eqref{definizione.mu.b}.
Notice that these latter estimates hold on the whole $\mathtt\Omega \times [\mathtt h_1, \mathtt h_2]$. 
 \end{thm}
 
This result implies the assumption ${\bf (AI)}$ and, by setting  (see \eqref{costanti.riducibilita})
\begin{equation}\label{kappa1}
\kappa_1\equiv k_0+2(\mathtt{m}-\mathtt{c}_1)+3\,,
\end{equation}
allows to deduce Theorem \ref{alm.approx.inv}, which is
the key step for a Nash-Moser iterative scheme.

\section{The Nash-Moser nonlinear iteration}\label{sec:NaM} 
In this section we prove Theorem \ref{NMT}. 
It will be a consequence of Theorem \ref{iterazione-non-lineare} below
where we construct iteratively a sequence of better and better approximate solutions 
of the equation $ \cF( i, \alpha, \mu) = 0$, 
with $\cF(i,\alpha, \mu)$ defined in \eqref{F_op}. 
We consider the finite-dimensional subspaces 
$$
E_n := \Big\{ \cI (\vphi ) = (\Theta , I , z) (\vphi) , \ \  
\Theta = \Pi_n \Theta, \ I = \Pi_n I, \ z = \Pi_n z \Big\}
$$
where $ \Pi_n  $ is the projector 
\begin{equation}\label{truncation NM}
\Pi_n := \Pi_{N_n} : z(\vphi,x) = \sum_{\ell  \in \Z^\nu,  j \in {S}_0^c} z_{\ell, j} e^{\ii (\ell  \cdot \vphi + j\cdot x)} 
\ \mapsto \Pi_n z(\vphi,x) 
:= \sum_{|(\ell ,j)| \leq N_n} 
z_{\ell,  j} e^{\ii (\ell  \cdot \vphi + j\cdot x)}  
\end{equation}
with $ N_n = N_0^{\chi^n} $, $N_0>0$ (chosen below large with respect to the diophantine constant $\gamma$) 
and 
we denote with the same symbol 
$ \Pi_n p(\vphi) :=  \sum_{|\ell | \leq N_n}  p_\ell e^{\ii \ell  \cdot \vphi} $. 
We define $ \Pi_n^\bot := {\rm Id} - \Pi_n $.  
The projectors $ \Pi_n $, $ \Pi_n^\bot$ satisfy the smoothing properties \eqref{p2-proi}, \eqref{p3-proi} for the weighted Whitney-Sobolev norm $\| \cdot \|_s^{k_0, \gamma}$ 
defined in \eqref{def norm Lip Stein uniform}.

\bigskip 

In view of the Nash-Moser Theorem \ref{iterazione-non-lineare} 
we introduce the following constants: 
\begin{align}
&  \mathtt e := 
\mu({\mathtt b}_1) +  4 \sigma_1  + 2\,, \qquad 
{\mathtt a}_2 := \chi^2 (\tau + 2(\tau + 2)(\tau + 1) + \kappa_1 + 4) + 3 \mathtt e + 1 \label{costanti nash moser} 
\\
% {\rm max}\{3 \mathtt e + 1, \chi p \red{(2\tau^2+7\tau+4)} + \chi(\mu(\mathtt b_{1}) + 2 \sigma_1) +1 \}\,,  
& \mathtt a_3 := \chi^{- 1} \mathtt a_2  - \mu(\mathtt b_{1}) - 2 \sigma_1 \,, 
\quad 
\mu_1 := 3 \mathtt e + 1, 
\qquad {\mathtt b}_2 := \mathtt e + \chi^{- 1} \mu_1 + \mathtt a_2 + 1 ,
%\qquad \chi = 3/ 2,
\label{costanti nash moser 1} 
\\
& \chi = 3/ 2\,,\quad \sigma_1 := \max \{ \bar \sigma\,, \wh\s \},
\label{costanti nash moser 2}
\end{align}
where $\bar \sigma := \bar \sigma(\tau, \nu, k_0) > 0$ is defined in Theorem \ref{alm.approx.inv}, 
$ \wh\s $ is the largest loss of regularity in the estimates of the Hamiltonian vector field $X_P$ in Lemma \ref{XP_est}, 
$\mu(\mathtt b_{1})$ is defined in \eqref{definizione.mu.b}. 

In this section, given  
$  W = ( \cI, \al, \varrho ) $ where
$   \cI = \cI (\lambda) $ %  \in H^{s}_\vphi  \times H^{s}_\vphi \times H^{s} $ 
is the periodic component of a torus as in \eqref{ICal}, and $ \alpha = \alpha (\lambda) \in \R^\nu $
we denote $ \|  W \|_{s}^{k_0, \gamma} := \max\{ \|  \cI \|_{s}^{k_0, \gamma} ,  |  \alpha |^{k_0, \gamma} 
{, |\varrho|^{k_{0}. \gamma}} \} $, 
where  $ \|  \cI \|_{s}^{k_0, \gamma}  $ is defined in \eqref{norma.ical}.
Since all the parameters in \eqref{costanti nash moser}-\eqref{costanti nash moser 2} are fixed,
in the following we shall only write $\lesssim$ in the estimates omitting 
the  dependence on such parameters.

\begin{thm}{\bf (Nash-Moser).} \label{iterazione-non-lineare} 
Fix any $s_0$ as in \eqref{Sone}.
There exist $ \delta_0$, $ C_* > 0 $, such that, if
\begin{equation}\label{nash moser smallness condition}  
\begin{aligned}
&N_0^{\tau_4} \e 
%\g^{-k_0-2(\tm - \tc_{1}) - 3} 
< \delta_0, 
\quad  \tau_4 := \max \Big\{ \tau_3 + \frac{k_0 + \kappa_1 + 1}{2},  \mathtt e +  {\mathtt a}_2  \Big\}\,,  
% \quad \red{N:= \tm - \tc_{1}}
\\
&N_0 := \gamma^{- 2}, 
\quad \gamma:= \e^a, 
\quad 0 < a < \frac{1}{2 \tau_4}\,,
\end{aligned}
\end{equation}
%where the constant $N$ is defined in \eqref{Scelta.N} and
where $ \tau_3 := \tau_3(\tau, \nu)$ is  defined in \eqref{costanti.riducibilita},
%Theorem \ref{iterazione riducibilita}, 
and $\kappa_1$ is given in Theorem \ref{alm.approx.inv},
then, for all $ n \geq 0 $: 
\begin{itemize}
\item[$(\cP1)_{n}$] 
there exists a $k_0$ times differentiable function $\tilde W_n : \mathtt \Omega  \times [\th_1, \th_2] 
\to  \R^\nu  {\times \R^\nu} \times E_{n -1}$, $ \lambda = (\omega, \th) \mapsto \tilde W_n (\lambda) 
:=  (\tilde \cI_n, \tilde \alpha_n - \omega {,\tilde \varrho_{n}}) $, for  $ n \geq 1$,  
and $\tilde W_0 :=(0,0, 0) $ (${E_{-1}:=\{0\}}$), satisfying 
\begin{equation}\label{ansatz induttivi nell'iterazione}
\| \tilde W_n \|_{s_0 + \mu({\mathtt b}_{1}) + \sigma_1}^{k_0, \gamma} \leq C_*   \e  \gamma^{-1}\,. 
\end{equation}
Let $\tilde U_n := U_0 + \tilde W_n$ where $ U_0 := (\vphi,0,0, \omega {, 0})$. %\red{io direi $U_0 := (\vphi,0,0, \omega,0)$}.
The difference $\tilde H_n := \tilde U_{n} - \tilde U_{n-1}$, $ n \geq 1 $,  satisfies
\begin{equation} \label{Hn}
\|\tilde H_1 \|_{s_0 + \mu({\mathtt b}_{1}) + \sigma_1}^{k_0, \gamma} \leq	 
C_* \e \gamma^{- 1} \,,\quad \| \tilde H_{n} \|_{ s_0 + \mu({\mathtt b}_{1}) + \sigma_1}^{k_0, \gamma} \leq C_* \e \gamma^{-1} N_{n - 1}^{- \mathtt a_3} \,,\quad \forall n \geq 2. 
\end{equation}
\item[$(\cP2)_{n}$]   
Setting $ {\tilde \imath}_n := (\vphi, 0, 0) + \tilde \cI_n $, we define 
\begin{equation}\label{def:cal-Gn}
\cG_{0} := \mathtt{\Omega} \times [\th_1, \th_2]\,, \quad 
\cG_{n+1} :=  \cG_n \cap {\bf \Lambda}_{\infty}^{  \gamma_n}({\tilde \imath}_n)
\,, \quad n \geq 0 \, , 
\end{equation}
where $  {\bf \Lambda}_{\infty}^{  \gamma_n}({\tilde \imath}_n) $ is defined in \eqref{Melnikov-invert} and $\gamma_n := \gamma(1 + 2^{- n})$, $n \geq 0$. 
Then, for all $\lambda \in \cG_n$,  setting $ N_{-1} := 1 $, we have 
\begin{equation}\label{P2n}
\| \cF(\tilde U_n) \|_{ s_{0}}^{k_0, \gamma}  \leq C_* \e N_{n - 1}^{- {\mathtt a}_2} \, .
\end{equation}
\item[$(\cP3)_{n}$] \emph{(High norms).} For all $\lambda \in \cG_n$ we have that
\begin{equation}\label{P3n}
\| \tilde W_n \|_{ s_{0}+ {\mathtt b}_2}^{k_0, \gamma} 
\leq C_* \e   N_{n - 1}^{\mu_1}\,.
\end{equation}
\end{itemize}
\end{thm}

\begin{proof}
To simplify notation, in this proof we denote $\| \ \|^{k_0, \gamma}$ by $\| \ \|$. 

\smallskip

{\sc Step 1:} \emph{Proof of} $(\cP1, \cP 2, \cP 3)_0$. The estimates \eqref{ansatz induttivi nell'iterazione} and \eqref{P3n} for $n=0$ are trivial. We now prove \eqref{P2n}.
Recalling \eqref{F_op},  by \eqref{stima.Xp} we have that for any $s\geq s_{0}$ 
\begin{equation}\label{stima cal F U0 nel lemma}
\|\cF(U_0) \|_s = \| \cF(\vphi, 0, 0, \oo, 0)\|_{s} =  \e \| X_{P}(\vphi, 0, 0)\|_{s} \lesssim_{s}  \e 
\end{equation}
then $(\cP1, \cP 2, \cP 3)_0$ hold taking $C_* $ large enough.
%One has $\|\cF(U_0) \|_s \lesssim \e$ 
%by \eqref{F_op}, \eqref{stima.dXp}, 
%then take $C_* $ large enough.

\smallskip

{\sc Step 2:} \emph{Assume that $(\cP1, \cP 2, \cP 3)_n$ hold for some $n \geq 0$, 
and prove $(\cP1, \cP 2, \cP3)_{n+1}$.}
By \eqref{nash moser smallness condition}
\begin{equation}\label{telefono1}
N_{0}^{\tau_{4}}\eps = \gamma^{- 2 \tau_4} \eps  = \eps^{-2 a \tau_{4}  +1} \le \delta_{0}\,.
\end{equation}
%for $\eps$ small enough.\\
We are going to define the next approximation $ \tilde U_{n+1} $ by a modified Nash-Moser scheme.
To this aim, we prove the approximate invertibility of the linearized operator 
\[
L_n := L_n(\lambda) := d_{i,\alpha,\varrho} \cF({\tilde \imath}_n(\lambda)) 
\]
by applying Theorem \ref{alm.approx.inv} to ${ L}_n(\lambda) $.
To prove that the inversion assumptions \eqref{almi4} 
hold, we apply Theorem \ref{inversione parziale cal L omega} with $ i = \tilde \imath_n $.
Let $K_{0}$ be the constant given by Theorem \eqref{Teorema di riducibilita}, we must verify he smallness condition \eqref{KAM smallness condition1}. 
By taking $\delta_{0}$ in \eqref{nash moser smallness condition} small enough, by \eqref{telefono1} one has
$\varepsilon$ small enough in such a way that $N_{0}\geq K_{0}$. %This is possible since $N_0 = \gamma^{- 2}= \varepsilon^{-2a}$ (see \eqref{nash moser smallness condition}). 
Then by the relations  $\kappa_1=2(\tm - \tc_{1})+2$ (see  \eqref{kappa1}) and by using
 \eqref{nash moser smallness condition} , one has that
%(and recalling the relations  $\kappa_1=2(\tm - \tc_{1})+2$ in \eqref{kappa1} and $N_0 = \gamma^{- 2}$ in \eqref{nash moser smallness condition}),
%the smallness condition \eqref{KAM smallness condition1} 
%holds for $ \e $ small enough, indeed
\[
K_0^{\tau_3} \e \gamma^{-k_0- \kappa_1 - 1}
{\leq} N_0^{\tau_3 + \frac{k_0 + \kappa_1 + 1}{2}} \e 
\leq N_0^{\tau_4} \e 
\stackrel{\eqref{nash moser smallness condition}}{\leq}\delta_0\,.
\]
Therefore Theorem \ref{inversione parziale cal L omega} applies, 
and we deduce that \eqref{almi4} hold 
for all $ \lambda \in {\bf \Lambda}_{\infty}^{ \gamma_n}(\tilde \imath_n)$, see \eqref{Melnikov-invert}.

Now we apply Theorem \ref{alm.approx.inv} to the linearized operator 
$ L_n(\lambda) $  with $ \mathtt{\Lambda}_o = {\bf \Lambda}_{\infty}^{ \gamma_n}(\tilde \imath_n)$ and 
\begin{equation}\label{valore finalissimo S}
\bar s = s_0 + \mathtt b_2
\end{equation}
where $\mathtt b_2$ 
is defined in \eqref{costanti nash moser 1}.
%(such a value of the Sobolev index $ S $ is the biggest one which is taken for the convergence of all the scheme). 
It implies the existence of
an approximate inverse ${\bf T}_n  := { \bf T} (\lambda, {\tilde \imath}_n(\lambda))$ 
satisfying
\begin{align}
\| {\bf T}_n g  \|_s 
 &\lesssim \gamma^{-1} \big( \| g \|_{s + \sigma_1} 
+ \| \tilde \cI_n \|_{s + \mu({\mathtt b}_{1}) + \sigma_1 } \| g \|_{s_0+ \sigma_1}\big) 
\quad 
\forall s_0 < s \leq \bar s - \mu(\tb_{1})-\bar\s \,, \label{stima Tn} 
\\
\| {\bf T}_n  g \|_{s_0} 
&\lesssim
 \gamma^{-1} \| g \|_{s_0 + \sigma_1} \label{stima Tn norma bassa}
\end{align}
because $\s_1 \geq \bar\s$ by \eqref{costanti nash moser 2}, 
where $\bar \s$ is the loss in \eqref{tame-es-AI} (use also that $\gamma \leq \gamma_n = \gamma(1 + 2^{- n}) \leq 2 \gamma$, for any $n \geq 0$).
For all $\lambda \in \cG_{n + 1} = \cG_n \cap {\bf \Lambda}_{\infty}^{\gamma_n}(\tilde \imath_n)$ (see \eqref{def:cal-Gn}), we define 
\begin{equation}\label{soluzioni approssimate}
\begin{aligned}
U_{n + 1} &:= \tilde U_n + H_{n + 1}\,, \qquad W_{n+1}:=\tilde{W}_{n}+H_{n+1}\,,
\\
H_{n + 1} &:=
( \widehat \cI_{n+1}, \widehat \alpha_{n+1}{, \widehat \mu_{n+1}}) :=  - {\bf \Pi}_{n } {\bf T}_n \Pi_{n } \cF(\tilde U_n) 
\in E_n \times \R^\nu  {\times \R^{\nu}}
\end{aligned}
\end{equation}
where  $ {\bf \Pi}_n $ is defined by (see \eqref{truncation NM})
\begin{equation*}%\label{proiettore modificato}
 {\bf \Pi}_n ({\cI}, \alpha, \varrho) := (\Pi_n \cI, \alpha, \varrho)\,,
 \quad {\bf \Pi}_n^\bot (\cI, \alpha, \varrho) := (\Pi_n^\bot \cI, 0{, 0})\,,\quad \forall (\cI, \alpha, \varrho)\,.
\end{equation*}
We show that the iterative scheme in \eqref{soluzioni approssimate} is rapidly converging. 
We write  
\begin{equation}\label{def:Qn} 
\cF(U_{n + 1}) =  \cF(\tilde U_n) + L_n H_{n + 1} + Q_n 
\end{equation} 
where $ L_n := d_{i,\alpha, \varrho} \cF(\tilde U_n)$ 
and $Q_n$ is defined by difference. 
Then, by the definition of $ H_{n+1} $ in \eqref{soluzioni approssimate}, 
we have %(recall also \eqref{proiettore modificato})
\begin{align}
\cF(U_{n + 1}) & = 
 \cF(\tilde U_n) - L_n {\bf \Pi}_{n } {\bf T}_n \Pi_{n } \cF(\tilde U_n) + Q_n \nonumber \\
 & = 
 \cF(\tilde U_n) - L_n  {\bf T}_n \Pi_{n } \cF(\tilde U_n) + L_n  {\bf \Pi}_n^\bot  {\bf T}_n \Pi_{n } \cF(\tilde U_n)
 + Q_n \nonumber\\ 
  & = 
 \cF(\tilde U_n) - \Pi_{n } \cF(\tilde U_n) - ( L_n  {\bf T}_n - {\rm Id})  \Pi_{n } \cF(\tilde U_n) + L_n  {\bf \Pi}_n^\bot  {\bf T}_n \Pi_{n } \cF(\tilde U_n)
 + Q_n \nonumber\\ 
%& =  \cF(\tilde U_n)  - \Pi_{n } L_n {\bf T}_n \Pi_{n }\cF(\tilde U_n) 
%+ ( L_n  {\bf \Pi}_n^\bot -  \Pi_n^\bot L_n ) {\bf T}_n \Pi_{n }\cF(\tilde U_n) + Q_n \nonumber\\
 & = \Pi_{n }^\bot \cF(\tilde U_n) + R_n + Q_n + P_n  
\label{relazione algebrica induttiva}
\end{align}
where 
\begin{equation}\label{Rn Q tilde n}
%R_n := (L_n  {\bf \Pi}_n^\bot -  \Pi_n^\bot L_n) {\bf T}_n \Pi_{n }\cF( \tilde U_n) \,,
 R_n := L_n  {\bf \Pi}_n^\bot  {\bf T}_n \Pi_{n } \cF(\tilde U_n)\,,
\qquad 
%P_n := - \Pi_{n } ( L_n {\bf T}_n - {\rm Id}) \Pi_{n } \cF( \tilde U_n)\,.
P_n := -( L_n  {\bf T}_n - {\rm Id})  \Pi_{n } \cF(\tilde U_n)\,.
\end{equation}
We first note that, for all $ \lambda \in \mathtt{\Omega} \times [\th_1, \th_2] $,  $s \geq s_0 $, 
by triangular inequality 
and by \eqref{F_op}, \eqref{stima.dXp}, \eqref{stima cal F U0 nel lemma},
\eqref{costanti nash moser 2}, \eqref{ansatz induttivi nell'iterazione} 
we have
\begin{equation}\label{F tilde Un W tilde n}
\| \cF(\tilde U_n)\|_{s} 
\lesssim_s \|\cF(U_0)\|_s + \| \cF(\tilde U_n) - \cF(U_0)\|_s 
\lesssim_s \e + \| \tilde W_n\|_{s + \sigma_1}  
\end{equation}
and, by \eqref{ansatz induttivi nell'iterazione}, \eqref{nash moser smallness condition}, 
\begin{equation}\label{gamma - 1 F tilde Un}
 \| \cF(\tilde U_n)\|_{s_0} \leq 1\, . 
\end{equation}

We now prove the following inductive estimates of Nash-Moser type.
\begin{lemma} \label{lemma:2017.0504.1}
For all $\lambda \in \cG_{n + 1}$ we have
\begin{align}
 \| \cF(U_{n + 1})\|_{s_0}  
&\lesssim
{N_n^{\mathtt e - {\mathtt b}_2}} (\e + \| \tilde W_n \|_{s_0 +\tb_2}) 
+ N_n^{\mathtt e} \| \cF(\tilde U_n)\|_{s_0}^2
\label{F(U n+1) norma bassa} 
\\ 
 \label{U n+1 alta}
\| W_{n + 1}\|_{s_0 + {\mathtt b}_2} &\lesssim
N_n^{\mathtt e}  (\e  +  \| \tilde W_n\|_{s_0 + {\mathtt b}_2}  )\, , \ n \geq 0 \, . 
\end{align}
where $\mathtt e$ and $\mathtt{b}_{2}$ are defined in \eqref{costanti nash moser}-\eqref{costanti nash moser 1}
\end{lemma}

\begin{proof} 
We first estimate $ H_{n +1} $ defined in  \eqref{soluzioni approssimate}.

\medskip

\noindent
{\bf Estimates of $ H_{n+1} $.}
 By \eqref{soluzioni approssimate} and \eqref{p2-proi}, %\eqref{p3-proi}, 
\eqref{stima Tn}, \eqref{stima Tn norma bassa}, \eqref{ansatz induttivi nell'iterazione}, $\gamma^{- 1} \leq \gamma^{- 2} = N_0 \leq N_n$,  we get 
\begin{align}
\|  H_{n + 1} \|_{s_0 + \tb_{2}} 
& \lesssim  \gamma^{- 1} 
\big( N_n^{\sigma_1} \|\cF(\tilde U_n) \|_{s_0 + {\mathtt b}_2} + 
N_n^{\mu({\mathtt b}_1) + 2 \sigma_1 } \|\tilde \cI_n \|_{s_0 + {\mathtt b}_2}\| \cF(\tilde U_n)\|_{s_0 }  \big)
\nonumber \\
& \stackrel{\eqref{F tilde Un W tilde n}, \eqref{gamma - 1 F tilde Un}}{\lesssim} 
N_n^{\mu(\mathtt b_1) + 2 \sigma_1 + 1 }  \big( \e  +  \| \tilde W_n \|_{s_0 +\tb_2} \big)\, ,  \label{H n+1 alta} 
\\
\label{H n+1 bassa}
\|  H_{n + 1}\|_{s_0} 
& \lesssim\gamma^{-1}N_{n}^{\sigma_1} \| \cF(\tilde U_n)\|_{s_0} \, .
\end{align}
Moreover, by \eqref{p2-proi}, \eqref{stima Tn}, \eqref{stima cal F U0 nel lemma}, one gets that 
\begin{equation}\label{stima H n + 1 per ansatz}
\begin{aligned}
\| H_1 \|_{s_0 + \mu(\mathtt b_1) + \sigma_1} & \lesssim \gamma^{- 1} \| {\mathcal F}(U_0) \|_{s_0 + \mu(\mathtt b_1) + 2 \sigma_1} \lesssim \e \gamma^{- 1}\,, \\
\| H_{n + 1} \|_{s_0 + \mu(\mathtt b_1) + \sigma_1} & \lesssim N_n^{\mu(\mathtt b_1) + \sigma_1} \| H_{n + 1} \|_{s_0} \stackrel{\eqref{H n+1 bassa}}{\lesssim} N_n^{\mu(\mathtt b_1) + 2 \sigma_1} \gamma^{- 1} \| \cF(\tilde U_n)\|_{s_0}\,, \quad n \geq 1\,. 
\end{aligned}
\end{equation}
Now we  estimate the terms $ Q_n $ in \eqref{def:Qn} and $ P_n , R_n $ in \eqref{Rn Q tilde n} in $ \| \ \|_{s_0} $ norm. 

\medskip

\noindent
{\bf Estimate of $ Q_n $.}
By 
\eqref{ansatz induttivi nell'iterazione},
\eqref{costanti nash moser 2},
\eqref{p2-proi},
\eqref{H n+1 bassa},
\eqref{P2n},
and since $3\s_1 - \mathtt{a}_2 \leq 0$ and $a < 1 / (1 + 3\s_1)$ 
(see \eqref{costanti nash moser}, \eqref{nash moser smallness condition}),
one has $\| \tilde W_n + t H_{n+1} \|_{s_0 + \wh\s} \lesssim 1$ for all $t \in [0,1]$.
Hence, by Taylor's formula, using
\eqref{def:Qn}, 
\eqref{F_op}, 
\eqref{stima.dXpquadro},
\eqref{H n+1 bassa},
\eqref{p2-proi},
and $\e \g^{-2} \leq 1$, we get
\begin{equation} \label{Qn norma bassa}
\| Q_n \|_{s_0} 
\lesssim \e \| H_{n+1} \|_{s_0+ \wh\s }^2 
\lesssim N_n^{4 \sigma_1} \| \cF(\tilde U_n) \|_{s_0}^2\, . 
\end{equation}
{\bf Estimate of $ P_n $.} 
By \eqref{splitting per approximate inverse}, we estimate
$ \cP({\tilde \imath}_n ) : = L_n {\bf T}_n - {\rm Id}$ and then $ P_n $ defined in \eqref{Rn Q tilde n}. 
By \eqref{p2-proi}-\eqref{p3-proi}, 
\begin{equation} \label{2017.0404.1}
\begin{aligned}
\| \cF(\tilde U_n) \|_{s_0 + \sigma_1} 
&\leq \|\Pi_n \cF(\tilde U_n) \|_{s_0  + \sigma_1} 
+ \|\Pi_n^\bot \cF(\tilde U_n) \|_{s_0  + \sigma_1} 
\\
&\leq N_n^{\sigma_1}  \| \cF(\tilde U_n)\|_{s_0} 
+ N_n^{ \sigma_{1} - \mathtt b_2} \|\cF(\tilde U_n) \|_{s_0 + \mathtt b_2}  \,.
%\leq K_n^{\sigma_1} ( \| \cF(\tilde U_n)\|_{s_0} 
%+ K_n^{- \mathtt b_2} \|\cF(\tilde U_n) \|_{s_0 + \mathtt b_2} ).
\end{aligned}
\end{equation}
By 
\eqref{pfi0},
\eqref{ansatz induttivi nell'iterazione},
\eqref{2017.0404.1}, 
and then \eqref{F tilde Un W tilde n},
\eqref{p2-proi}, $\gamma^{- 1} \leq \gamma^{- 2} = N_0 \leq N_n$,
we obtain
\begin{align} 
\| P_n  \|_{s_0}  
& \lesssim
\g^{-1} N_n^{\sigma_1 } \| \cF(\tilde U_n)\|_{s_0+\s_{1}} \| \cF(\tilde U_n)\|_{s_0} 
\nonumber
\\
&\lesssim
 N_n^{2 \sigma_1 + 1} \| \cF(\tilde U_n)\|_{s_0} 
( \| \cF(\tilde U_n)\|_{s_0} 
+ N_n^{- {\mathtt b}_2} \| \cF(\tilde U_n) \|_{s_0 + {\mathtt b}_2} ) 
\nonumber
\\ & 
\lesssim
 N_n^{2 \sigma_1 + 1} \| \cF(\tilde U_n)\|_{s_0} 
( \| \cF(\tilde U_n)\|_{s_0} + N_n^{\s_1 - {\mathtt b}_2}
(\e + \| \tilde W_n\|_{s_0 + \mathtt b_2} ) )
\nonumber
\\
& 
\stackrel{\eqref{gamma - 1 F tilde Un}}{\lesssim} N_n^{2 \sigma_1 + 1} 
\| {\mathcal F}(\tilde U_n) \|_{s_0}^2 
+ N_n^{3 \sigma_1 + 1- \mathtt b_2} (\e + \| \tilde W_n\|_{s_0 + \mathtt b_2} ) \,. 
\label{Q n 1 bassa}
\end{align}
%By 
%\eqref{pfi1}, 
%\eqref{ansatz induttivi nell'iterazione}, 
%\eqref{p2-proi}, $\gamma^{- 1} = K_0 \leq K_n$, 
%we have
%\begin{equation} \label{Q n omega bassa}
%\| P_{n, \omega} \|_{s_0}   \lesssim_{s_0 + \mathtt b_2} 
%\e \gamma^{- \kappa_1-1} N_{n - 1}^{- {\mathtt a}_1} K_n^{\sigma_1} \| \cF(\tilde U_n) \|_{s_0} \lesssim_{s_0 + \mathtt b_2} 
%\e  N_{n - 1}^{- {\mathtt a}_1} K_n^{\sigma_1 + \kappa_1 + 1} \| \cF(\tilde U_n) \|_{s_0}\,,
%\end{equation}
%where $\mathtt{a}_1$ is defined in \eqref{costanti.riducibilita}.
%By 
%\eqref{pfi2}, 
%\eqref{p2-proi},
%\eqref{costanti nash moser 1},
%\eqref{P2n},
%and then 
%\eqref{F tilde Un W tilde n}, 
%\eqref{p2-proi},
%we get
%\begin{align}
%\| P_{n, \omega}^\bot\|_{s_0} 
%& \lesssim_{s_0 + \mathtt b_2} \g^{-1} K_{n}^{-\tb_{2}} (\| \cF(\tilde U_n) \|_{s_0 + {\mathtt b}_2} 
%+  K_{n}^{\mu(\tb_1)+2\s_{1}} \| \cI \|_{s_{0}+\tb_{2}} \| \cF(\tilde U_n) \|_{s_0})
%\\
%& \lesssim_{s_0 + \mathtt b_2}  K_{n}^{\mu({\mathtt b}_1) + 2 \s_1 - {\mathtt b}_2} 
%\g^{-1} ( \| \cF(\tilde U_n) \|_{s_0 + {\mathtt b}_2} +  \| \tilde W_n \|_{s_0 + {\mathtt b}_2}) \nonumber \\ 
%& \lesssim_{s_0 + \mathtt b_2}  K_{n}^{\mu({\mathtt b}_1) + 3 \sigma_1 - {\mathtt b}_2}
%\g^{-1} (\e + \| \tilde W_n \|_{s_0 + \mathtt b_2 }). 
%\label{Q n omega bot bassa} 
%\end{align}

\noindent
{\bf Estimate of $ R_n $.} 
By 
\eqref{stima.dXp}, \eqref{costanti nash moser 2}, 
\eqref{ansatz induttivi nell'iterazione}, we get
\begin{equation}\label{stima commutatore modi alti norma bassa}
\| L_n H \|_{s_0} 
\lesssim \| H \|_{s_0 + \sigma_1}\,. 
\end{equation}
Hence, by 
\eqref{p3-proi}, \eqref{Rn Q tilde n},
\eqref{stima commutatore modi alti norma bassa},  
\eqref{stima Tn},
\eqref{ansatz induttivi nell'iterazione},
\eqref{p2-proi},
and then 
\eqref{F tilde Un W tilde n},
\eqref{p2-proi},
\eqref{gamma - 1 F tilde Un}, $\gamma^{- 1} \leq \gamma^{- 2} = N_0 \leq N_n$,
we get 
\begin{align} 
\| R_n\|_{s_0}  & \lesssim  \| {\bf \Pi}_n^\bot  {\bf T}_n \Pi_{n } \cF(\tilde U_n)  \|_{s_0 + \sigma_1} 
\lesssim 
N_n^{- \mathtt b_2 + \sigma_1} \| {\bf T}_n \Pi_{n } \cF(\tilde U_n) \|_{s_0 + \mathtt b_2} 
\nonumber
\\
& \lesssim 
N_n^{- \mathtt b_2 + \sigma_1}  \g^{-1} 
( \| \Pi_n \cF(\tilde U_n) \|_{s_0 + {\mathtt b}_2 + \sigma_1} 
+ \| \tilde W_n \|_{s_0 + {\mathtt b}_2 + \mu(\mathtt b_1) + \sigma_1} \| \Pi_n \cF(\tilde U_n) \|_{s_0 + \s_1} ) 
\nonumber 
\\ 
& \lesssim
 N_n^{- \mathtt b_2 + 3 \sigma_1 + \mu(\mathtt b_1) + 1} 
( \|  \cF(\tilde U_n) \|_{s_0 + {\mathtt b}_2 } 
+ \| \tilde W_n \|_{s_0 + {\mathtt b}_2 } \|   \cF(\tilde U_n) \|_{s_0 } ) 
\nonumber 
\\ 
& \lesssim
N_n^{ \mu({\mathtt b}_1) + 4 \s_1 + 1 - {\mathtt b}_2} 
(\e + \| \tilde W_n\|_{s_0 + \mathtt b_2} ). 
\label{stima Rn norma bassa}
\end{align}
We can finally estimate $ \cF(U_{n + 1}) $. % in $ \| \ \|_{s_0} $. 
By 
\eqref{relazione algebrica induttiva},
\eqref{p3-proi},
\eqref{F tilde Un W tilde n},
\eqref{stima Rn norma bassa},
\eqref{Qn norma bassa}, 
\eqref{Q n 1 bassa}, 
%  \eqref{nash moser smallness condition},  % (using $\e \gamma^{- 1} \leq 1$), 
\eqref{ansatz induttivi nell'iterazione}, 
we get \eqref{F(U n+1) norma bassa}. 
%By \eqref{soluzioni approssimate} and \eqref{stima Tn} we have bound \eqref{U n+1 alta} 
%for $W_1 := H_1$, namely 
%$$ 
%\| W_1 \|_{s_0+ {\mathtt b}_2} = \| H_1 \|_{s_0+ {\mathtt b}_2}  \lesssim_{s_0+ {\mathtt b}_2} \g^{-1} \| \cF(U_0)\|_{s_0+ {\mathtt b}_2 + \sigma_1} 
%\lesssim_{s_0+ {\mathtt b}_2}  \e \gamma^{- 1} \, .
%$$ 
Estimate \eqref{U n+1 alta} for $ W_{n+1} := \tilde W_n + H_{n+1} $, $ n \geq 1 $,  
follows by \eqref{H n+1 alta}. 
\end{proof}
Now that Lemma \ref{lemma:2017.0504.1} has been proved, 
we continue the proof of Theorem \ref{iterazione-non-lineare}. 
As a corollary of Lemma \ref{lemma:2017.0504.1} we get the following lemma, where for clarity 
we use the extended notation $ \| \  \|^{k_0, \gamma} $ (instead of $\| 	 \|$ used above). 

\begin{lemma}\label{lemma:quadra}
For all  $ \lambda \in \cG_{n + 1}$, we have
\begin{align}
& \qquad \quad \| \cF(U_{n + 1}) \|_{s_0}^{k_0, \gamma}  
\leq C_* \e N_n^{- \mathtt a_2} \label{stima F u n + 1 induttiva} \, , \qquad 
\| W_{n + 1} \|_{s_0 + \mathtt b_2}^{k_0, \gamma} 
\leq C_* \e  N_n^{\mu_1} \, ,  \\ % \label{stima W n + 1 tilde W n} 
& \label{stima H n+1 lemma}
\| H_1 \|_{s_0 + \mu(\mathtt b_1) + \sigma_1}^{k_0, \gamma} \leq C \e \g^{-1} \, , \quad
\| H_{n + 1}\|_{s_0 + \mu(\mathtt b_1) + \sigma_1}^{k_0, \gamma} 
\lesssim  \e \gamma^{- 1} N_n^{\mu(\mathtt b_1) + 2 \sigma_1 } N_{n - 1}^{- \mathtt a_2} \, ,
\  	\ n \geq 1 \,.
\end{align}
\end{lemma}

\begin{proof}
First note that, by  \eqref{def:cal-Gn},  if  
$ \lambda \in  \cG_{n + 1}$, then 
$ \lambda \in  \cG_n$ and so  \eqref{P2n} and the inequality in $(\cP3)_{n}$ hold. 
Then the first inequality in \eqref{stima F u n + 1 induttiva} follows by 
\eqref{F(U n+1) norma bassa},
$ (\cP2)_{n}$, $(\cP3)_{n}$,  $ \e  \ll 1 $ small, 
and by \eqref{costanti nash moser}, \eqref{costanti nash moser 1}.
%(see also Remark \ref{choice-a-b}). 
For $ n = 0 $ we use also  \eqref{nash moser smallness condition}. 

\noindent
The second inequality in \eqref{stima F u n + 1 induttiva}  follows directly from 
 \eqref{U n+1 alta}, $(\cP3)_{n}$, by the choice of $\mu_1$ 
 in \eqref{costanti nash moser 1} and by taking $N_0 \gg 1$ large enough.  
The estimates \eqref{stima H n+1 lemma} follow by 
\eqref{stima H n + 1 per ansatz} and by \eqref{P2n}. 
\end{proof}

By Theorem 
%$B.2$ in \cite{BBHM},
\ref{thm:WET}, 
we define a $ k_0 $ times differentiable extension $ \tilde H_{n + 1}$ 
of $ (H_{n + 1})_{|\cG_{n+1}} $ to the whole $ \R^\nu \times [\th_1, \th_2] $. For convenience we consider the restriction of $\tilde H_{n + 1}$ to $\mathtt \Omega \times [\mathtt h_1, \mathtt h_2]$.
Clearly this extension $\tilde H_{n + 1}$ satisfies the same bound for $ H_{n+1} $ in \eqref{stima H n+1 lemma} and therefore, by  the definition 
of $ {\mathtt a}_2 $ in \eqref{costanti nash moser}, the estimate \eqref{Hn} at $ n + 1 $ holds.   

Finally we define the functions
$$
\tilde W_{n+1} := \tilde W_{n} + \tilde H_{n + 1} \, , \quad 
\tilde U_{n + 1} := \tilde U_n + \tilde H_{n + 1} = U_0 + \tilde W_n + \tilde H_{n + 1} = U_0 + \tilde W_{n + 1}\,,
$$
which are defined for all $\lambda \in \mathtt \Omega \times [\th_1, \th_2] $ and satisfy 
$$
\tilde W_{n + 1} = W_{n + 1} \, ,  \quad 
\tilde U_{n + 1} = U_{n + 1} \, , \ \ \forall \lambda \in \cG_{n + 1} \, . 
$$ 
Therefore  $(\cP2)_{n + 1}$, $(\cP3)_{n + 1}$ are proved by Lemma \ref{lemma:quadra}.
Moreover by \eqref{Hn}, which has been proved up to the step $n + 1 $, 
%by using \eqref{stima H n+1 estensione Whitney}, 
we have 
$$
\| \tilde W_{n + 1} \|_{s_0 + \mu({\mathtt b}_1) + \sigma_1}^{k_0, \gamma} 
\leq {\mathop \sum}_{k = 1}^{n + 1} \| \tilde H_k \|_{s_0 + \mu({\mathtt b}_1) + \sigma_1}^{k_0, \gamma} 
\leq C_*  \e  \gamma^{-1}
$$ 
and thus \eqref{ansatz induttivi nell'iterazione} holds also 
at the step $n + 1$. This completes the proof of Theorem \ref{iterazione-non-lineare}.
\end{proof}

\subsection{Proof of Theorem \ref{NMT}}\label{proof.NMT}
 
Let $ \gamma = \e^{\rm a}  $ with $ {\rm a} \in (0, {\rm a_0}) $ and $ {\rm a_0} := 1 / (2\tau_4) $ where $\tau_4$ is defined in \eqref{nash moser smallness condition}.
Let $\kappa_1$ as in \eqref{kappa1} and define (recall \eqref{Sone}, \eqref{costanti nash moser 2} 
and \eqref{definizione.mu.b})
\[
\bar{\s}:=[(\nu + d)/2] +2+\mu(\mathtt{b}_1)+\s_1\,.
\]
Given any $\mathfrak{s}\geq \bar{\s}$ we apply Theorem \ref{iterazione-non-lineare}
with 
\begin{equation}\label{defSgotico}
s_0:=[(\nu + d)/2] +2+\mathfrak{s}-\bar{\s}\,.
\end{equation}
Note that thanks to this choice one has
$s_0+\mu(\mathtt{b}_1)+\s_1=\mathfrak{s}$ and that \eqref{Sone} is fulfilled.
Moreover, using $N_0=\gamma^{-2}=\e^{-2a}$ and the definition of $\tau_4$, one has
that 
the smallness condition given by the first inequality in 
\eqref{nash moser smallness condition} holds for $ 0 < \e < \e_0 $ 
small enough and Theorem \ref{iterazione-non-lineare} applies.   
By \eqref{Hn} the  sequence of functions 
\[
{\tilde W}_n = {\tilde U}_n - (\vphi,0,0, \omega, 0) 
:= (\tilde \cI_n, \tilde \alpha_n - \omega, \tilde\varrho_{n}) =
(\tilde \imath _n - (\vphi,0,0), \tilde \alpha_n - \omega, \tilde\varrho_{n}) 
\]
is a Cauchy sequence in $ \| \cdot \ \|_{s_0+\mu(\mathtt{b}_1)+\s_1}^{k_0, \gamma} $ 
and then it converges to a function
\[
W_\infty := \lim_{n \to + \infty} {\tilde W}_n  \, , \quad 
%\quad {\rm with} \quad
\]
with
\[
W_\infty : \mathtt \Omega \times [\th_1, \th_2] \to H^{s_0+\mu(\mathtt{b}_1)+\s_1}_\vphi  
\times H^{s_0+\mu(\mathtt{b}_1)+\s_1}_\vphi \times H^{s_0+\mu(\mathtt{b}_1)+\s_1}_{\vphi, x}
\times \R^\nu {\times \R^{\nu}}  \, . 
\]
We define 
\begin{equation}\label{QP:soluz}
U_\infty := (i_\infty, \alpha_\infty, \varrho_{\infty}) = (\vphi,0,0, \omega, 0) +  W_\infty \, .
\end{equation}
By \eqref{ansatz induttivi nell'iterazione} and \eqref{Hn} we also deduce that
\begin{equation}\label{U infty - U n}
\|  U_\infty -  U_0 \|_{s_0 + \mu(\mathtt b_1) + \sigma_1}^{k_0, \gamma} \leq C_* \e \gamma^{- 1} \,, 
\quad 
\| U_\infty - {\tilde U}_n \|_{s_0  + \mu({\mathtt b}_1) + \sigma_1}^{k_0, \gamma} 
\leq C_* \e \gamma^{-1} N_{n }^{- \mathtt a_3}\,, \ \  n \geq 1 \, .
\end{equation}
Moreover by Theorem \ref{iterazione-non-lineare}-$(\cP2)_n$, we deduce that 
(recall \eqref{F_op})
\[
\cF(\lambda;U_\infty(\lambda))  
= \cF(\lambda;i_{\infty},\alpha_{\infty}, \varrho_\infty) := 
\omega\cdot\pa_\vphi i_{\infty}(\vphi) - X_{H_{\alpha_{\infty}, \varrho_{\infty}}}(i_{\infty}(\vphi))=0
\]
%$ \cF(\lambda, U_\infty(\lambda)) = 0 $   
for all $ \lambda  $ belonging to 
\begin{equation}\label{defGinfty}
%{\mathcal G}_\infty := 
\bigcap_{n \geq 0} \cG_n = 
\cG_0 \cap \bigcap_{n \geq 0} 
 {\bf \Lambda}_{\infty}^{\gamma_n}(\tilde \imath_{n}) 
\stackrel{\eqref{Melnikov-invert}}{=} % , \eqref{Cantor set}, \eqref{prime di melnikov}}{=}
\cG_0 \cap \Big[ \bigcap_{n \geq 0}  {\Lambda}_{\infty}^{ \gamma_n}(\tilde \imath_{n}) \Big] \cap 
 \Big[ \bigcap_{n \geq 0}    \Lambda_{\infty}^{\gamma_n, I}(\tilde \imath_{n}) \Big]\, \,,
\end{equation}
where $\cG_0 =  \mathtt \Omega \times [\th_1, \th_2] $ is defined in \eqref{def:cal-Gn}.  
By the first inequality in \eqref{U infty - U n} and \eqref{defSgotico} 
we deduce
estimates \eqref{alpha_infty} and \eqref{i.infty.est} with $C\equiv C_{*}$.
By Lemma \ref{Lemma6.1DP}, we also deduce that actually $\varrho_{\infty}\equiv0$
and that  $i_{\infty}(\vphi)$ is invariant for the vector field 
$X_{H_{\alpha_\infty}} \equiv X_{H_{\alpha_\infty}, 0}$ . Then, we define the set of good
parameters $\mathcal{G}_{\infty}$ as  (recall \eqref{defGinfty} and \eqref{Melnikov-invert})
\begin{equation}\label{defGinftyBis}
\mathcal{G}_{\infty}:=\mathcal{G}_{\infty}(\gamma):= \bigcap_{n \geq 0} \cG_n \,.
\end{equation}
The discussion above 
implies the thesis of Theorem  \ref{NMT} up to the condition 
\eqref{megaMisura} on the measure of the set.
 This argument is discussed in the next subsection.

\subsection{Measure estimates: proof of Theorem \ref{NMT} concluded}\label{subsec:measest}
In this subsection we provide measure estimates on the set of good parameters $\mathcal{G}_{\infty}$
in \eqref{defGinftyBis} and prove the estimate \eqref{megaMisura}. 
More precisely we prove the following.
\begin{thm} {\bf (Measure estimates).}\label{MEASEST} 
Let
\begin{equation}\label{param_small_meas}
\tau > m_0 ( \nu +1) + 2 m_0 d \,, \qquad \gamma = \varepsilon^{a} 
\end{equation}
where $a$ is given in \eqref{nash moser smallness condition}, $m_0$ is the index of non-degeneracy 
given in Proposition \ref{lem:transversality}, 
%and $\kappa_1$ given by Theorem \ref{NMT}.
Then, for $ \varepsilon \in (0, \varepsilon_0) $ small enough,  the measure of the set
\begin{equation}\label{Gvare2}
\cC_\varepsilon := \big\{ \th \in [\th_1, \th_2] 
\ : \  \big( \alpha_\infty^{-1}( \overline{\omega}(\th),\th ),\th \big) 
\in {\mathcal G}_\infty \big\}
\end{equation}
satisfies $ | {\mathcal C}_\varepsilon | \rightarrow\th_2-\th_1$ as 
$\varepsilon\rightarrow 0$.
\end{thm}
\noindent
The result above is quite technical and involves several arguments which will be discussed in the rest of the section.
First of all, by the expansions \eqref{alpha_infty}-\eqref{inv_alpha100} in Thm. \ref{NMT}-$(1)$
we have

\begin{equation}\label{Om-per}
\overline{ \omega}_\varepsilon (\th):= \alpha_\infty^{-1}(\overline{\omega}(\th),\th) =
\overline{\omega}(\th) +\vec{r}_{\varepsilon}(\th) \,,
\end{equation}
where 
$\vec{r}_\varepsilon(\th) := \breve{r}_\varepsilon(\overline{\omega}(\th),\th) $
satisfies 
\begin{equation}\label{eq:tang_res_est}
|\pa_\th^k {\vec r}_\varepsilon (\th)| \leq C \varepsilon\gamma^{-(1+k)} \,, 
\quad \forall\,\abs k \leq k_0 \,, \ \text{uniformly on } [\th_1,\th_2] \,,
\end{equation}
where $k_0:= m_0+2$.
Therefore, the definition of ${\mathcal C}_\e$ in \eqref{Gvare2} becomes 
\begin{equation}\label{Gvare}
\cC_\varepsilon := \big\{ \th \in [\th_1, \th_2] 
\ : \  \big( \overline{ \omega}_\varepsilon (\th),\th \big) 
\in {\mathcal G}_\infty \big\}\,.
\end{equation}
We shall use in the sequel the following notation. For any function $f (\omega, \mathtt h)$, $(\omega , \mathtt h) \in \R^\nu \times [\mathtt h_1, \mathtt h_2]$, we write 
$$
f(\mathtt h) \equiv f(\overline{ \omega}_\varepsilon (\th)\,,\, \mathtt h), \quad \mathtt h \in [\mathtt h_1, \mathtt h_2]\,. 
$$
For any $n \geq 0$, we denote by (recall definition \eqref{def:cal-Gn})
\begin{equation}\label{mathcal Cn}
{\mathcal C}_n := \Big\{ \mathtt h \in [\mathtt h_1, \mathtt h_2] : (\overline{ \omega}_\varepsilon (\th)\,,\, \mathtt h) \in {\mathcal G}_n \Big\}\,,
\end{equation}
 and by recalling also \eqref{tDtCn}, we define for any $n \geq 0$
\begin{equation}\label{def mathcal En 0 a}
\begin{aligned}
{\mathcal E}_n^{(0)} & := \Big\{ \mathtt h \in {\mathcal C}_n: (\overline \omega_\e (\mathtt h), \mathtt h) \in  \tT\tC_{\infty}(\gamma_n,\tau; i_n)\Big\}\,. 
\end{aligned}
\end{equation}
Note that ${\mathcal C}_n \setminus {\mathcal E}_n^{(0)}$ is given by  
 \begin{equation}\label{def mathcal En 0}
\begin{aligned}
& {\mathcal C}_n \setminus {\mathcal E}_n^{(0)} = \bigcup_{\ell \in \Z^\nu \setminus \{ 0 \}} {\mathcal R}_n^{(0)}(\ell)\,, \\
& {\mathcal R}_n^{(0)}(\ell) := \Big\{ \mathtt h \in {\mathcal C}_n : |(\overline{ \omega}_\varepsilon (\th) - \mathtt V \mathtt m_1(i_n)) \cdot \ell| < \frac{4 \gamma_n }{\langle \ell \rangle^\tau} \Big\}\,.
\end{aligned}
\end{equation}
 We now analyze the complement set of ${\mathcal C}_\e$ in \eqref{Gvare}, which is given by 
\begin{equation}\label{inclusione Cn C n + 1 C epsilon}
[\mathtt h_1, \mathtt h_2] \setminus {\mathcal C}_\e \subseteq \bigcup_{n \geq 0} {\mathcal C}_n \setminus {\mathcal C}_{n + 1}\,.
\end{equation}

By \eqref{def:cal-Gn}, \eqref{tDtCn}\,,\, \eqref{Melnikov-invert}, \eqref{Cantor set}, \eqref{prime.di.Melnikov}, one gets that for any $n \geq 0$, 

\begin{equation}\label{splitting cal Gn meno cal G n+1}
\begin{aligned}
{\mathcal C}_n \setminus {\mathcal C}_{n + 1} & \subseteq ({\mathcal C}_n \setminus {\mathcal E}_n^{(0)}) \cup {\mathcal E}_n^{(1)} \cup {\mathcal E}_n^{(2, +)} \cup {\mathcal E}_n^{(2, -)} \,,
\end{aligned}
\end{equation}
\[
{\mathcal E}_n^{(1)} := \bigcup_{\begin{subarray}{c}
(\ell, j) \in \Z^\nu \times S_0^c \\
\mathtt V^T \ell + j = 0
\end{subarray}} {\mathcal R}_n^{(1)}(\ell, j)\,,
\qquad
{\mathcal E}_n^{(2, \s)} := \bigcup_{\begin{subarray}{c}
(\ell, j, j') \in (\Z^\nu \setminus \{ 0 \}) \times S_0^c \times S_0^c \\
\mathtt V^T \ell + j +\s j' = 0
\end{subarray}} {\mathcal R}_n^{(2, \s)}(\ell, j, j')\,,\;\;\;\s\in\{\pm\}\,,
\]
where the ``resonant'' sets ${\mathcal R}_n^{(1)}(\ell, j)$, ${\mathcal R}_n^{(2, \s)}(\ell, j, j')$, $\s\in\{\pm\}$
are defined as follows:
\begin{equation}\label{def mathcal En 1}
\begin{aligned}
%& {\mathcal E}_n^{(1)} := \bigcup_{\begin{subarray}{c}
%(\ell, j) \in \Z^\nu \times S_0^c \\
%\mathtt V^T \ell + j = 0
%\end{subarray}} {\mathcal R}_n^{(1)}(\ell, j)\,, \\
%& 
{\mathcal R}_n^{(1)}(\ell, j) := \Big\{ \mathtt h \in {\mathcal E}_n^{(0)} : |\overline{ \omega}_\varepsilon (\th) \cdot \ell + \mu_\infty(j; \mathtt h, i_n(\mathtt h))|  < \frac{2 \gamma_n }{\langle \ell \rangle^\tau |j|^\tau} \Big\}\,,
\end{aligned}
\end{equation}
\begin{equation}\label{def mathcal En 2 -}
\begin{aligned}
%& {\mathcal E}_n^{(2, -)} := \bigcup_{\begin{subarray}{c}
%(\ell, j, j') \in (\Z^\nu \setminus \{ 0 \}) \times S_0^c \times S_0^c \\
%\mathtt V^T \ell + j - j' = 0
%\end{subarray}} {\mathcal R}_n^{(2, -)}(\ell, j, j')\,, 
%\\
%& 
{\mathcal R}_n^{(2, -)}&(\ell, j, j') 
\\&
:= \Big\{ \mathtt h \in {\mathcal E}_n^{(0)} : |\overline{ \omega}_\varepsilon (\th) \cdot \ell 
+ \mu_\infty(j; \mathtt h, i_n(\mathtt h)) - \mu_\infty(j'; \mathtt h, i_n(\mathtt h)) |  
< \frac{2 \gamma_n }{\langle \ell \rangle^\tau |j'|^\tau} \Big\}\,,
\end{aligned}
\end{equation}
\begin{equation}\label{def mathcal En 2 +}
\begin{aligned}
%& {\mathcal E}_n^{(2, +)} := \bigcup_{\begin{subarray}{c}
%(\ell, j, j') \in (\Z^\nu \setminus \{ 0 \}) \times S_0^c \times S_0^c \\
%\mathtt V^T \ell + j + j' = 0
%\end{subarray}} {\mathcal R}_n^{(2, +)}(\ell, j, j')\,, \\
%& 
{\mathcal R}_n^{(2, +)}&(\ell, j, j') 
\\&:= \Big\{ \mathtt h \in {\mathcal E}_n^{(0)} : |\overline{ \omega}_\varepsilon (\th) \cdot \ell + \mu_\infty(j; \mathtt h, i_n(\mathtt h)) + \mu_\infty(j'; \mathtt h, i_n(\mathtt h)) |  < \frac{2 \gamma_n }{\langle \ell \rangle^\tau |j'|^\tau} \Big\}\,.
\end{aligned}
\end{equation}
Our first purpose is to compute the Lebesgue measure of ${\mathcal C}_n \setminus {\mathcal C}_{n + 1}$. We shall first prove the following lemma
\begin{lemma}\label{lemma inclusione cantor}
For any $n \geq 0$, the following inclusions hold:
\begin{equation}\label{inclusioni cantor per misura}
\begin{aligned}
&  {\mathcal C}_n \setminus {\mathcal E}_n^{(0)}  \subseteq  \bigcup_{|\ell| \geq N_{n }} {\mathcal R}_n^{(0)}(\ell)\,, 
%\\& 
\qquad\qquad  {\mathcal E}_n^{(1)} \subseteq \bigcup_{\begin{subarray}{c}
(\ell, j) \in \Z^\nu \times S_0^c \\
|\ell| \geq N_{n },\;\;
\mathtt V^T \ell + j = 0
\end{subarray}} {\mathcal R}_n^{(1)}(\ell, j) \,, \\
& {\mathcal E}_n^{(2, \s)} \subseteq \bigcup_{\begin{subarray}{c}
(\ell, j, j') \in (\Z^\nu \setminus \{ 0 \}) \times S_0^c \times S_0^c \\
|\ell| \geq N_{n }\,,\;\;
\mathtt V^T \ell + j +\s j' = 0
\end{subarray}} {\mathcal R}_n^{(2, \s)}(\ell, j, j')\,, \;\;\;\s\in\{\pm\}\,.
%\\
%& {\mathcal E}_n^{(2, +)} \subseteq \bigcup_{\begin{subarray}{c}
%(\ell, j, j') \in (\Z^\nu \setminus \{ 0 \}) \times S_0^c \times S_0^c \\
%|\ell| \geq N_{n } \\
%\mathtt V^T \ell + j - j' = 0
%\end{subarray}} {\mathcal R}_n^{(2, +)}(\ell, j, j')\
\end{aligned}
\end{equation}
\end{lemma}
\begin{proof}
We shall prove the inclusion involving ${\mathcal E}_n^{(2, -)}$. All the other inclusions can be proved by similar arguments (they are actually even easier). The proof is organized in several claims.

\medskip

\noindent
{\bf Claim 1.} Let us consider $\mathtt h \in {\mathcal E}_n^{(0)}$ with the constraints on the indexes
\begin{equation}\label{constraintClaim1}
0 < |\ell| \leq N_n\,,\qquad j, j' \in S_0^c\,,\;\;\;|j'| \geq N_n^q\,,\;\;\;q = 2 (\tau + 2)\,,\;\;\;
\mathtt V^T \ell + j - j' = 0\,.
\end{equation}
%$0 < |\ell| \leq N_n$, $j, j' \in S_0^c$, $|j'| \geq N_n^q$, $q = 2 (\tau + 2)$, $\mathtt V^T \ell + j - j' = '0$. 
Then  one has 
\begin{equation}\label{claim 3 inclusione cantor riducibilita a}
|\overline \omega_\e(\mathtt h) \cdot \ell + \mu_\infty(j; \mathtt h, i_n(\mathtt h)) - \mu_\infty(j'; \mathtt h,  i_n(\mathtt h))| \geq \frac{2 \gamma_n}{\langle \ell \rangle^\tau |j'|^\tau}\,. 
\end{equation}
This implies that ${\mathcal R}_n^{(2, -)}(\ell, j, j') = \emptyset$.

\medskip

\noindent
{\bf Proof of the claim 1.}
To simplify notations, we write 
$
\mu_\infty(j) \equiv \mu_\infty(j; \mathtt h) \equiv \mu_\infty(j: \mathtt h, i_n(\mathtt h))\,,
$
and 
we analyze 
\[
\delta_\infty(j) := \mu_\infty(j) - \mathtt{m}_{1} \cdot j \stackrel{\eqref{def cal N infty nel lemma}}{=}  \mu_0(j) - \mathtt{m}_{1} \cdot j + r_\infty(j) 
\stackrel{\eqref{proprieta cal D bot}}{=}
 \mathtt m_7(j) + r_\infty(j)\,.
\]
%for $j, j' \in S_0^c$ and $|j'| \geq N_n^q$, $q := 2 (\tau + 2)$, $|\ell| \leq N_n$, $\mathtt V^T \ell + j - j' = 0$, 
By conditions \eqref{constraintClaim1}
one has that $|j - j'| \leq C |\ell| \leq C N_n$ and by formulae \eqref{proprieta cal D bot}, \eqref{def cal N infty nel lemma}, 
one gets 
\begin{equation*}%\label{formule mu n j - m_1 dot ja}
\delta_\infty(j) - \delta_\infty(j') = \tm_7(j) - \tm_7(j') + r_\infty(j) - r_\infty(j')\,.
\end{equation*}
%By the estimate \eqref{uova fritte 2}, proved in Proposition \ref{iterazione riducibilita}-${\bf(S3)_{n}}$ one gets that
By the mean value theorem, there exists $\theta \in [0, 1]$ 
such that 
\begin{equation}\label{uova fritte 0}
|\tm_7(j) - \tm_7(j')| \leq |\nabla_\xi \tm_7(j' + \theta (j - j'))| |j - j'|
\end{equation}
and since $\mathtt m_7$ is a Fourier multiplier in $S^{\frac12}$, by the estimates \eqref{stima lambda 7}, one has that 
\begin{equation}\label{uova fritte 1}
|\nabla_\xi \tm_7(j' + \theta (j - j'))| \lesssim \langle j' + \theta (j - j') \rangle^{- \frac12}\,. 
\end{equation} 
Now since $|j'| \geq N_n^q$, $|j -j'| \leq C N_n$ and $p > 1$, one has that for any $\theta \in [0, 1]$
\begin{equation}\label{uova fritte 3}
\begin{aligned}
& |j' + \theta (j - j')| \geq |j'| - |j - j'| \geq N_n^q - C N_n = N_n^q (1 - C N_n^{1 - q}) \geq \frac12 N_n^q 
\end{aligned}
\end{equation}
by taking $N_0 \gg 0$ large enough. 
Then by \eqref{uova fritte 0}-\eqref{uova fritte 3} together with $|j - j' | \lesssim N_n$, one gets 
\begin{equation*}%\label{uova fritte 2a}
|\tm_7(j) - \tm_7(j')| \lesssim \frac{1}{N_n^{\frac{q}{2} - 1} }\,.
\end{equation*}
Moreover using that $|j'| \geq N_n^q$, $|j - j'| \leq C N_n$ and
$$
|j| = |j' +  (j - j')| \geq |j'| - |j - j'| \geq N_n^q - C N_n = N_n^q (1 - C N_n^{1 - q}) \geq \frac12 N_n^q\,,
$$
and the estimate \eqref{stime forma normale limite}  implies that 
\begin{equation*}%\label{uova fritte 4a}
|r_\infty(j) - r_\infty(j')| \leq |r_\infty(j)| + |r_\infty(j')| \lesssim  (|j|^{ - \mathtt m + \mathtt c_1} 
+ |j'|^{- \mathtt m + \mathtt c_1}) \lesssim  N_n^{ - q( \mathtt m - \mathtt c_1)}\,. 
\end{equation*}
Hence using that $q (\mathtt m- \mathtt c_1) > \frac{q}{2} - 1$ (see \eqref{costanti.riducibilita}), 
one gets the bound 
\begin{equation}\label{claim 2 S3 n+1a}
|\delta_\infty(j) - \delta_\infty(j')| \lesssim N_n^{- \frac{q}{2} + 1}\,.
\end{equation}

Now, since $\mathtt h \in {\mathcal E}_n^{(0)}$ (recall \eqref{def mathcal En 0 a}), one has that if $\ell \neq 0$, $\mathtt V^T \ell + j - j' = 0$, one has that 
\begin{equation}\label{prima melnikov per S3 n + 1a}
|\overline \omega_\e(\mathtt h) \cdot \ell +\mathtt{m}_{1}  \cdot (j - j')| \geq \frac{4 \gamma_n }{\langle \ell \rangle^\tau}\,. 
\end{equation}
Then, one obtains that 
$$
\begin{aligned}
|\overline \omega_\e(\mathtt h) \cdot \ell + \mu_\infty(j) - \mu_\infty(j')| & \geq |\overline \omega_\e(\mathtt h) \cdot \ell 
+ \mathtt{m}_{1} \cdot (j - j')| - |\delta_\infty(j) - \delta_\infty(j')|  \\
& \stackrel{\eqref{claim 2 S3 n+1a}, \eqref{prima melnikov per S3 n + 1a}}{\geq} \frac{4 \gamma_n}{\langle \ell \rangle^\tau} - \frac{C}{N_n^{\frac{q}{2} - 1}}  \geq \frac{2 \gamma_n}{\langle \ell \rangle^\tau}
\end{aligned}
$$
provided that 
\begin{equation}\label{tecnocasa merdaa}
N_n^{\frac{q}{2} - \tau - 1} \gamma  \gg 1, \quad \text{for any} \quad n \geq 0.
\end{equation} 
Since, $q = 2(\tau + 2)$, one has that $\frac{q}{2} - \tau - 1 = 1$, therefore the condition \eqref{tecnocasa merdaa} is fullfilled since  $N_0 = \gamma^{- 2}$ and $0 < \gamma \ll 1$ (see \eqref{nash moser smallness condition}). The proof of the {\bf claim 1} is then concluded. 

\medskip

{\bf Claim 2.} Let us consider $\mathtt h \in {\mathcal E}_n^{(0)}$ with the constraints on the indexes
\begin{equation}\label{constraintClaim2}
0 < |\ell| \leq N_n\,,\qquad j, j' \in S_0^c\,,\;\;\;|j'| \leq N_n^q\,,\;\;\;q = 2 (\tau + 2)\,,\;\;\;
\mathtt V^T \ell + j - j' = '0\,.
\end{equation} 
%$0 < |\ell| \leq N_n$, $j, j' \in S_0^c$, $|j'| \leq N_n^q$, $q = 2(\tau + 2)$, $\mathtt V^T \ell + j - j' = 0$. 
Then 
\begin{equation}\label{claim 3 inclusione cantor riducibilita b}
|\overline \omega_\e(\mathtt h) \cdot \ell + \mu_\infty(j; \mathtt h, i_n(\mathtt h)) - \mu_\infty(j'; \mathtt h,  i_n(\mathtt h))| \geq \frac{2 \gamma_n}{\langle \ell \rangle^\tau |j'|^\tau}\,. 
\end{equation}
This implies that ${\mathcal R}_n^{(2, -)}(\ell, j, j') = \emptyset$. 

\medskip

\noindent
{\bf Proof of the claim 2.} 
By the estimates \eqref{stima delta 12 mu 0 pre rid}, \eqref{stime forma normale limitea} (applied with $i_1 \equiv i_n$, $i_2 \equiv i_{n - 1}$),  for any value of the parameter $\mathtt h \in [\mathtt h_1, \mathtt h_2]$ and for any $j \in S_0^c$ (recalling that $\kappa_1 = 2(\mathtt m - \mathtt c_1) + 2$), using that $\gamma^{- 1} \leq \gamma^{- 2} = N_0 \leq N_n$, one has that 
\begin{equation*}%\label{L infty vn - v n - 1}
\begin{aligned}
|\mu_\infty(j; \mathtt h, i_n(\mathtt h)) - \mu_\infty(j; \mathtt h, i_{n - 1}(\mathtt h))| 
& \lesssim \e \gamma^{- \kappa_1}  |j| \| i_n - i_{n - 1} \|_{s_0 + \mu(\mathtt b_1)} \\
& \stackrel{\eqref{Hn}}{\lesssim} \e \gamma^{- \kappa_1 - 1}  |j|   N_{n - 1}^{- \mathtt a_3}    \lesssim N_n^{\kappa_1 + 1} N_{n - 1}^{- \mathtt a_3} \e |j| \,.
\end{aligned}
\end{equation*}
%If $|j'| \leq N_n^q$, $|\ell| \leq N_n$, $\mathtt V^T \ell + j - j' = 0$, one has that 
In view of \eqref{constraintClaim2} we deduce
\[
|j| \leq |j'| + |j - j'| \lesssim |j'| + |\ell| \lesssim N_n^q + N_n \lesssim N_n^q\,,
\]
and hence 
\begin{equation*}%\label{L infty vn - v n - 1}
\begin{aligned}
 |\mu_\infty(j; \mathtt h, i_n(\mathtt h)) - \mu_\infty(j; \mathtt h, i_{n - 1}(\mathtt h))|  
&\lesssim N_n^{q + \kappa_1 + 1} N_{n - 1}^{- \mathtt a_3} \e  
\\
  |\mu_\infty(j'; \mathtt h, i_n(\mathtt h)) - \mu_\infty(j'; \mathtt h, i_{n - 1}(\mathtt h))| 
  &\lesssim N_n^{q + \kappa_1 + 1} N_{n - 1}^{- \mathtt a_3} \e\,. 
\end{aligned}
\end{equation*}

Thus if $ \mathtt h \in {\mathcal C}_n$ % \cap {\mathcal O}_\infty^{\gamma_n}(v_n)$ 
 with $0 < |\ell| \leq N_{n}$ and $|j'| \leq N_n^q$, $\mathtt V^T \ell + j - j' = 0$, we get, for some constant $C > 0$ large enough, 
%\begin{equation}\label{bound bla so ri 0}
%\begin{aligned}
%|\lambda \im \omega \cdot \ell + \mu_\infty(j; w_n) - \mu_\infty(j'; w_n)| & \geq |\lambda \im \omega \cdot \ell + \mu_\infty(j; w_{n - 1}) - \mu_\infty(j'; w_{n - 1})| \\
%& \quad - 2 \sup_{j \in \Z^2 \setminus \{ 0 \}} |\mu_\infty(j; w_n) - \mu_\infty(j; w_{n - 1})| \\
%& \geq \frac{2 \lambda \gamma_{n - 1}}{\langle \ell \rangle^\tau |j'|^\tau} - C \Big( \tN_{n-1}^{2 \overline \sigma + 1} \tN_{n-2}^{-\mathtt a_1}
%+  \tN_{n - 1}^{- \mathtt a} \Big) \\
%& \geq \frac{2 \lambda \gamma_{n }}{\langle \ell \rangle^\tau |j'|^\tau}
%\end{aligned}
%\end{equation}
\begin{equation*} %\label{bound bla so ri 0}
		\begin{aligned}
		 & | \overline \omega_\e(\mathtt h) \cdot \ell  + \mu_\infty(j; \mathtt h, i_n(\mathtt h)) - \mu_\infty(j'; \mathtt h, i_n(\mathtt h))|  \\
		& \geq | \overline \omega_\e(\mathtt h) \cdot \ell  + \mu_\infty(j; \mathtt h, i_{n - 1}(\mathtt h)) - \mu_\infty(j'; \mathtt h, i_{n - 1}(\mathtt h))| \\
		& \quad - 
		|\mu_\infty(j; \mathtt h, i_n(\mathtt h)) - \mu_\infty(j; \mathtt h, i_{n - 1}(\mathtt h))|  - 
		|\mu_\infty(j'; \mathtt h, i_n(\mathtt h)) - \mu_\infty(j'; \mathtt h, i_{n - 1}(\mathtt h))|  \\
		& \geq \frac{2  \gamma_{n - 1}}{\langle \ell \rangle^\tau |j'|^\tau} - C N_n^{\kappa_1 + q + 1}  N_{n - 1}^{- \mathtt a_3}   \e  
		%\\&  
		\geq \frac{2  \gamma_{n }}{\langle \ell \rangle^\tau |j'|^\tau} \,, 
	\end{aligned}
\end{equation*}
provided  that 
\[
C  N_n^{\kappa_1 + q + 1}  N_{n - 1}^{- \mathtt a_3}   \e \leq \frac{2  (\gamma_{n - 1} - \gamma_n)}{\langle \ell \rangle^\tau |j'|^\tau}\,.
\]
Using that $|\ell| \leq  N_{n}$, $|j'| \leq N_n^q$, $\gamma_{n - 1} - \gamma_n = \gamma 2^{- n}$, $\gamma^{- 1} \leq \gamma^{- 2} = N_{0} \leq N_{n}$, the latter condition is implied by 
\begin{equation}\label{giobbe1}
	C' 2^{n } N_{n}^{\tau + \tau q + 1} N_n^{\kappa_1 + q + 1}  N_{n - 1}^{- \mathtt a_3}  \e    \leq 1\,, \quad n \geq 0\,, 
\end{equation}
for some constant $C' \geq C$. The condition \eqref{giobbe1} is verified by \eqref{costanti.riducibilita}, \eqref{costanti nash moser}, \eqref{costanti nash moser 1} (which imply that 
$\mathtt a_3 > \chi(\tau + (\tau + 1) q + \kappa_1  + 2 )$) and 
the smallness condition \eqref{nash moser smallness condition}.
%, since $\bar \tau > \chi(\tau + \tau^2 + 1) + 2 \bar \sigma + 1$, $\mathtt a > \chi(\tau + \tau^2 + 1)$ and $\mathtt a_1 > \chi^2(\tau + \tau^2 + 1) + \chi(2 \bar \sigma + 1)$.   
The proof of the {\bf claim 2} is then concluded. 

%\medskip

In conclusion, {\bf claim 1} and {\bf claim 2} imply the inclusion on the set ${\mathcal E}_n^{(2, -)}$ in \eqref{inclusioni cantor per misura}.
%$$
%{\mathcal E}_n^{(2, -)} \subseteq \bigcup_{\begin{subarray}{c}
%(\ell, j, j') \in (\Z^\nu \setminus \{ 0 \}) \times S_0^c \times S_0^c \\
%|\ell| \geq N_{n } \\
%\mathtt V^T \ell + j - j' = 0
%\end{subarray}} {\mathcal R}_n^{(2, -)}(\ell, j, j')\,. 
%$$
\end{proof}

The next step is to estimate the measure of the ``resonant sets" appearing in 
\eqref{def mathcal En 1}, \eqref{def mathcal En 2 -}, \eqref{def mathcal En 2 +}. 
We also denote, with a small abuse of notation (recall \eqref{proprieta cal D bot}, 
\eqref{def cal N infty nel lemma}, \eqref{wild}) %Lemma \ref{lemma autovalore 0 pre red}) 
for all $j\in S_0^c$, 
\begin{equation}\label{eq:final_eig_kappa}
\mu_\infty(j; \th) \equiv \mu_\infty(j ;  \mathtt h, i_n(\mathtt h)) = 
  \Omega(j; \mathtt h)
+ \fr_\infty(j; \mathtt h) 
\end{equation}
where $\fr_{\infty} (j; \mathtt h) := \mu_0(j; \mathtt h) - \Omega(j; \mathtt h) + r_\infty(j; \mathtt h)$. By the estimates \eqref{beta.FGMP.est}, \eqref{stima partial mu 0 pre rid}, \eqref{stima tipo simbolo mu 0 j j'}, \eqref{stime forma normale limite} and \eqref{stime forma normale limitea} (recall that $\kappa_1 = 2 (\mathtt m - \mathtt c_1) + 2$), we get the estimates for any $j \in S_0^c$ and for any $k \leq k_0$
\begin{align}
& |\partial_{\mathtt h}^k \mathtt m_1(\mathtt h)| \leq C \e \gamma^{- k}, \quad \forall \mathtt h \in [\mathtt h_1, \mathtt h_2]\,, \label{small_coeff_k0} \\
& |\pa_\th^k \fr_\infty(j; \mathtt h) |
\leq 
C \varepsilon\gamma^{-k - \kappa_1} |j|\, , \quad \forall j \in S_0^c, \quad \mathtt h \in [\mathtt h_1, \mathtt h_2]\,, \label{small_coeff_k} 
\\
 & |\pa_\th^k (\fr_\infty( j; \th) - \fr_\infty(j', \th)  )|
%\sup_{j\in S_0^c }|j|^{-\frac12}\abs{ \pa_\th^k \fr_j^\infty(\th) } 
\leq 
C \varepsilon\gamma^{-\kappa_1-k} \langle j - j' \rangle\,, 
\qquad  j, j' \in S^c_0\,, \quad \mathtt h \in [\mathtt h_1, \mathtt h_2] \,.  \label{small_rem_kb} 
\end{align}
In order to estimate the measure of the sets appearing in \eqref{def mathcal En 1}, \eqref{def mathcal En 2 -}, \eqref{def mathcal En 2 +},  
the key point is to prove that the perturbed frequencies satisfy 
transversality properties similar to the ones 
\eqref{imp1}-\eqref{imp4} satisfied by the unperturbed frequencies. 
By Proposition \ref{lem:transversality}, 
 \eqref{Om-per}, and the estimates \eqref{eq:tang_res_est}, 
\eqref{small_coeff_k0}, \eqref{small_coeff_k}, \eqref{small_rem_kb}, we deduce 
the following  lemma
(cfr. Lemma 5.5 in \cite{BFM}).  

\begin{lemma} {\bf (Perturbed transversality)} \label{lem:pert_trans}
	For $\varepsilon\in(0,\varepsilon_0)$ small enough and  $\forall\,\th\in[\th_1,\th_2]$, 
	\begin{align}
	& \max_{0\leq n \leq m_0}
	| \partial_\th^n \overline{\omega}_\varepsilon (\th)\cdot \ell | \geq \frac{\rho_0}{2} \jap{\ell} \,, 
	\quad \forall\,\ell\in\Z^\nu\setminus\{0\} \,; \label{eq:00_meln}
	\\
	&
	\max_{0\leq n \leq m_0} | \partial_\th^n 
	(\overline{\omega}_\varepsilon(\th)
	 - \tV \mathtt{m}_1 (\th)) \cdot \ell  | \geq \frac{\rho_0}{2}\jap{\ell} \, , 
	  \quad \forall \ell\in\Z^\nu \setminus \{0\} \label{eq:0_meln_pert}
	  \\	
	& 
	\begin{cases}
	\max_{0\leq n \leq m_0} \Big| \partial_\th^n \Big( \overline{\omega}_\varepsilon(\th)\cdot \ell + \mu_\infty(j; \th) \Big) \Big|
	\geq \frac{\rho_0}{2}	\jap{\ell} \, , 
	\\
	\tV^{T}\ell + j = 0 \,, \quad \ell\in\Z^\nu\,, \ j\in S_0^c \,;
	\end{cases} 
 \label{eq:1_meln_pert}
 \\
	&\begin{cases}
	\max_{0\leq n \leq m_0} \Big| \partial_\th^n \Big( \overline{\omega}_\varepsilon(\th)\cdot \ell + \mu_\infty(j; \th)  - \mu_\infty(j'; \th)   \Big) \Big| \geq \frac{\rho_0}{2}	\jap{\ell} 
	\\
	\tV^{T} \ell + j -j'= 0 \,, \quad \ell\in\Z^\nu \setminus \{ 0 \}\,, \ j,j'\in S_0^c \, ; 
	\end{cases} \label{eq:2_meln-_pert}
	\\
	&\begin{cases}
	\max_{0\leq n \leq m_0} \Big|  \partial_\th^n \Big( \overline{\omega}_\varepsilon(\th)\cdot \ell + \mu_\infty(j; \th)  + \mu_\infty(j'; \th)  \Big)  \Big| \geq \frac{\rho_0}{2} 	\jap{\ell} 
	\\
	\tV^{T}\ell + j + j'= 0 \,, \quad \ell\in\Z^\nu\,, \ j,   j'\in S_0^c  \, . 
	\end{cases} \label{eq:2_meln+_pert}
	\end{align}
	We recall that $\rho_0$ is the amount of non-degeneracy that has been defined in Proposition \ref{lem:transversality}.
\end{lemma}

\begin{proof}
We prove \eqref{eq:2_meln-_pert}. The proofs of \eqref{eq:00_meln}, \eqref{eq:0_meln_pert},  \eqref{eq:1_meln_pert} 
and \eqref{eq:2_meln+_pert} are similar.
By \eqref{Om-per} and \eqref{eq:final_eig_kappa} 
we have
\begin{equation*}
\begin{aligned}
\overline{\omega}_\varepsilon(\th)\cdot \ell +  \mu_\infty(j; \th)  - \mu_\infty(j'; \th)
&=\overline{\omega}(\th)\cdot\ell 
+\Omega(j,\th)-\Omega(j';\th)
\\&
+ \vec{r}_{\varepsilon}(\th) \cdot \ell 
+ \fr_\infty(j; \mathtt h)-\fr_\infty(j'; \mathtt h)\,.
\end{aligned}
\end{equation*}
Hence for any $0 \leq k \leq m_0$
\[
\begin{aligned}
\Big| \partial_{\mathtt h}^k \Big( \overline{\omega}_\varepsilon(\th)\cdot \ell 
&+ \mu_\infty(j; \th)  - \mu_\infty(j'; \th) \Big) \Big| 
\geq 
\Big| \partial_{\mathtt h}^k \Big( \overline{\omega}(\th)\cdot \ell 
+ \Omega(j;\th)-\Omega(j';\th) \Big) \Big| 
\\& 
-  |\partial_{\mathtt h}^k \vec{r}_{\varepsilon}(\th) \cdot \ell| 
- |\partial_{\mathtt h}^k (\fr_\infty(j ; \mathtt h)-\fr_\infty(j'; \mathtt h))| 
\\& 
\stackrel{\eqref{eq:tang_res_est}\,,\, \eqref{small_coeff_k}-\eqref{small_rem_kb}}{\geq}  
\Big|\partial_{\mathtt h}^k \Big(\overline{\omega}(\th)\cdot \ell + \Omega(j;\th)
-\Omega(j';\th) \Big) \Big| 
- C \e \gamma^{- 1- k} |\ell| 
 \\& 
- C \e \gamma^{- k  - \kappa_1} \langle j - j' \rangle
\end{aligned}
\]
and using the conservation of momentum $\mathtt V^T \ell + j - j' = 0$ 
(which implies that $|j - j'| \lesssim |\ell|$), one obtains that for any $0 \leq k \leq m_0$, 
$$
\begin{aligned}
\Big| \partial_{\mathtt h}^k \Big( \overline{\omega}_\varepsilon(\th)\cdot \ell + \mu_\infty(j; \th)  - \mu_\infty(j'; \th) \Big) \Big|& \geq \Big| \partial_{\mathtt h}^k \Big( \overline{\omega}(\th)\cdot \ell 
+ \Omega(j;\th)-\Omega(j';\th) \Big) \Big|  -  C_1 \e \gamma^{- \kappa_1 - m_0} \langle \ell \rangle
\end{aligned}
$$
for some large constant $C_1 > C> 0$. 
%A direct calculation using \eqref{dispersionLaw} shows that, for $\th\in [\th_1,\th_2]$,
%\begin{equation}\label{posacenere2}
%|\pa_{\th}^{k}\big(\Omega_j(\th)-\Omega_{j'}(\th)\big)|\leq C_{k}||j|^{\frac{1}{2}}-|j'|^{\frac{1}{2}} |\,,\quad \forall\, k\geq0\,.
%\end{equation}
%Furthermore, using that $\fr^\infty(\th;j)$ is a symbol satisfying \eqref{coeff_fin_small} and the mean value theorem, 
%we get
%\begin{equation}\label{posacenere3}
%|\pa_{\th}^{k} (\fr^\infty(\th;j)-\fr^\infty(\th;j'))|\leq C_{k}\e \gamma^{-\kappa_1-k}|j-j'|\,,\quad 0\leq k\leq m_0\,.
%\end{equation}
%By  \eqref{posacenere1}, \eqref{posacenere2}, \eqref{posacenere3}, using also 
%\eqref{eq:tang_res_est}, \eqref{small_coeff_k} and the momentum condition,
%we deduce that, for any $0\leq k\leq m_0$,
%\[
%\begin{aligned}
%| \partial_\th^k( \overline{\omega}_\varepsilon(\th)\cdot \ell + d_j^\infty(\th)-d_{j'}^\infty(\th) ) | &\geq 
%|\pa_{\th}^{k}(\overline{\omega}(\th)\cdot\ell 
%+\Omega_j(\th)-\Omega_{j'}(\th))|-C\e\gamma^{-k_1-m_0}\langle \ell\rangle
%\\&
%\stackrel{\eqref{imp3}}{\geq}\rho_0\langle\ell\rangle-C\e\gamma^{-k_1-m_0}\langle \ell\rangle
%\geq \rho\langle\ell\rangle/2\,,
%\end{aligned}
%\]
The latter estimate, together with \eqref{imp3} and by taking $\e\gamma^{-\kappa_1-m_0}\leq \rho_0/(2C_1)$ (see \eqref{param_small_meas}) implies that 
$$
{\rm max}_{0 \leq k \leq m_0} \Big\{ \Big| \partial_{\mathtt h}^k \Big( \overline{\omega}_\varepsilon(\th)\cdot \ell + \mu_\infty(j; \th)  - \mu_\infty(j'; \th) \Big) \Big| \Big\} \geq \frac{\rho_0 \langle \ell \rangle}{2}
$$
which is exactly the estimate \eqref{eq:2_meln-_pert}. 
\end{proof}

As a consequence of the transversality estimates 
\eqref{eq:00_meln}-\eqref{eq:2_meln+_pert} 
we get the following bounds for the resonant sets defined in \eqref{def mathcal En 0}, \eqref{def mathcal En 1}-\eqref{def mathcal En 2 +}. 

\begin{lemma} {\bf (Estimates of the resonant sets).} \label{lem:meas_res}
The measure of the sets \eqref{def mathcal En 0}, \eqref{def mathcal En 1}-\eqref{def mathcal En 2 +} satisfy
\begin{align*}
& |  {\mathcal R}_n^{(0)}(\ell) | \lesssim 
( \gamma\jap{\ell}^{-(\tau+1)} )^{\frac{1}{m_0}}  \,, 
 \qquad 
| {\mathcal R}_n^{(1)}(\ell, j) |\lesssim \big( \gamma |j|^{- \tau}\jap{\ell}^{-(\tau+1)} \big)^{\frac{1}{m_0}} \,,
\\
& | {\mathcal R}_n^{(2, -)}(\ell, j, j') |\,,\, | {\mathcal R}_n^{(2, +)}(\ell, j, j') |\lesssim 
\big( \gamma \jap{\ell}^{-(\tau+1) }  \jap{j'}^{-\tau} \big)^{\frac{1}{m_0}} \,.
\end{align*}
\end{lemma}

\begin{proof}
It follows by using  R\"ussmann Theorem $17.1$ in \cite{Russ},
which can be applied in view of Lemma \ref{lem:pert_trans}.
For more details
we refer to Lemma $5.6$ in \cite{BFM}.
\end{proof}

\begin{proof}[Proof of Theorem \ref{MEASEST} completed.]
We now estimate the measure of the set ${\mathcal C}_\e$ in  \eqref{Gvare2}. 
By the inclusion \eqref{inclusione Cn C n + 1 C epsilon}, one has to estimate the measure of the set ${\mathcal C}_n \setminus {\mathcal C}_{n + 1}$ for any $n \geq 0$. By \eqref{splitting cal Gn meno cal G n+1} 
\begin{equation*}%\label{conto stima misura 1}
|{\mathcal C}_n \setminus {\mathcal C}_{n + 1}| \leq |{\mathcal C}_n \setminus {\mathcal E}_n^{(0)}| + |{\mathcal E}_n^{(1)}| + |{\mathcal E}_n^{(2, +)}| +  |{\mathcal E}_n^{(2, -)}|\,. 
\end{equation*}
By applying Lemmata \ref{inclusioni cantor per misura}, \ref{lem:meas_res}, one has that  
$$
 |{\mathcal E}_n^{(2, -)}| \leq  \sum_{\begin{subarray}{c}
(\ell, j, j') \in (\Z^\nu \setminus \{ 0 \}) \times S_0^c \times S_0^c \\
|\ell| \geq N_{n } \\
\mathtt V^T \ell + j - j' = 0
\end{subarray}} |{\mathcal R}_n^{(2, -)}(\ell, j, j')|  \lesssim \sum_{\begin{subarray}{c}
(\ell, j, j') \in (\Z^\nu \setminus \{ 0 \}) \times S_0^c \times S_0^c \\
|\ell| \geq N_{n } \\
\mathtt V^T \ell + j - j' = 0
\end{subarray}} \big( \gamma \jap{\ell}^{-(\tau+1) }  \jap{j'}^{-\tau} \big)^{\frac{1}{m_0}}\,.
$$
The constraint $\mathtt V^T \ell + j - j' = 0$ implies that $|j| \lesssim |j'| + |\ell| \lesssim |j'| \langle \ell \rangle$, hence 
$$
\begin{aligned}
 |{\mathcal E}_n^{(2, -)}|  & \lesssim \sum_{\begin{subarray}{c}
(\ell, j, j') \in (\Z^\nu \setminus \{ 0 \}) \times S_0^c \times S_0^c \\
|\ell| \geq N_{n } \\
|j| \lesssim |j'| \langle \ell \rangle
\end{subarray}} \big( \gamma \jap{\ell}^{-(\tau+1) }  \jap{j'}^{-\tau} \big)^{\frac{1}{m_0}}  \\
& \lesssim \sum_{\begin{subarray}{c}
(\ell, j') \in (\Z^\nu \setminus \{ 0 \}) \times S_0^c  \\
|\ell| \geq N_{n } 
\end{subarray}} \langle \ell \rangle^d |j'|^d\big( \gamma \jap{\ell}^{-(\tau+1) }  
\jap{j'}^{-\tau} \big)^{\frac{1}{m_0}} \\
& \lesssim \gamma^{\frac{1}{m_0}} \sum_{|\ell| \geq N_n } \langle \ell \rangle^{- \frac{\tau + 1}{m_0} + d} \sum_{j' \in S_0^c} |j'|^{- \frac{\tau}{m_0} + d} \lesssim \gamma^{\frac{1}{m_0}} N_n^{- \frac{\tau + 1}{m_0} + d+\nu} 
\\
& \lesssim \gamma^{\frac{1}{m_0}} N_n^{- 1}
\end{aligned}
$$
by using that, by the choice of $\tau$ in \eqref{param_small_meas}, the series $
\sum_{\ell \in \Z^\nu } \langle \ell \rangle^{- \frac{\tau + 1}{m_0} + d}\,,\, \sum_{j' \in S_0^c}
 |j'|^{- \frac{\tau}{m_0} + d} < + \infty$ are convergent and $ \frac{\tau + 1}{m_0} - d-\nu > 1$. 
 By similar arguments one can show that 
$$
|{\mathcal E}_n^{(0)}|\,,\, |{\mathcal E}_n^{(1)}|\,,\, |{\mathcal E}_n^{(2, +)}| \lesssim \gamma^{\frac{1}{m_0}} N_n^{- 1}
$$
implying that 
$$
|{\mathcal C}_n \setminus {\mathcal C}_{n + 1}| \lesssim \gamma^{\frac{1}{m_0}} N_n^{- 1}\,.
$$
Therefore (recall again \eqref{inclusione Cn C n + 1 C epsilon}) 
$$
|[\mathtt h_1, \mathtt h_2] \setminus {\mathcal C}_\e| \leq \sum_{n \geq 0} |{\mathcal C}_n \setminus {\mathcal C}_{n + 1}| \lesssim \gamma^{\frac{1}{m_0}} \sum_{n \geq 0} N_n^{- 1} \lesssim \gamma^{\frac{1}{m_0}}
$$
since the series $\sum_{n \geq 0} N_n^{- 1}  < + \infty$ is convergent. The claimed statement then follows using that, by \eqref{param_small_meas}, $\gamma = \e^{ a}$. 
\end{proof}

\subsection{Linear stability and  proof of Theorem \ref{thm:main0}}\label{sec:linearstability}
We shall  prove the existence of quasi-periodic 
solutions of the original Hamiltonian system $ H_\varepsilon $ in \eqref{hamsez4}
by applying Theorem \ref{NMT}. 
%which is equivalent after a rescaling to  \eqref{eq:Ham_eq_zeta},
% and not of just of the Hamiltonian system generated 
% by the  modified Hamiltonian 
% $ H_{\alpha_\infty} $. 
 We proceed as follows. 
We recall that by  \eqref{alpha_infty}-\eqref{inv_alpha100}, the function $\alpha_\infty(\,\cdot\,,\th)$ from $\mathtt{\Omega}$ into its image $\alpha_\infty(\mathtt{\Omega},\th)$ is invertible and 
\begin{equation}\label{inv_alpha}
\begin{aligned}
& \beta = \alpha_\infty(\omega,\th) = \omega+r_\varepsilon(\omega,\th) \ \Leftrightarrow \\
&   \omega = \alpha_\infty^{-1}(\beta,\th) = \beta+\breve{r}_\varepsilon(\beta,\th) \, , \quad \abs{ {r}_\varepsilon }^{k_0,\gamma}\,,\, \abs{ \breve{r}_\varepsilon }^{k_0,\gamma} \leq C\varepsilon\gamma^{-1} \,.
\end{aligned}
\end{equation}
%We underline that the function $\alpha_\infty^{-1}(\,\cdot\,,\th)$ is the inverse of $\alpha_\infty(\,\cdot\,,\th)$ at any fixed value $\th$ in $[\th_1,\th_2]$. 
Then, for any $\beta\in\alpha_\infty({\mathcal G}_\infty)$ (see \eqref{defGinftyBis}), 
Theorem \ref{NMT} proves the existence of an embedded invariant torus filled by quasi-periodic solutions with Diophantine frequency $\omega=\alpha_\infty^{-1}(\beta,\th)$ for the Hamiltonian
\begin{equation*}
H_\beta = \beta \cdot I+  \tfrac12({\bf \Omega}z,z)_{L^2} + \varepsilon P \,.
\end{equation*}
Consider the curve of the unperturbed tangential frequency vector 
$ \overline{\oo}(\th) $ in \eqref{LinearFreqWW}. 
%Thanks to Theorem \ref{MEASEST} 
Thanks to the measure estimate  \eqref{megaMisura}
we have that for "most" values of $\th\in[\th_1,\th_2]$ the vector $(\alpha_\infty^{-1}(\overline{\omega}(\th),\th) ,\th )$ is in ${\mathcal G}_\infty$ 
(recall \eqref{defGinfty}-\eqref{defGinftyBis}),
obtaining
an embedded torus for the Hamiltonian $H_\varepsilon$ in  \eqref{hamsez4}, 
filled by quasi-periodic solutions with Diophantine frequency vector 
$\omega = \alpha_\infty^{-1}(\overline{\oo}(\th),\th) $.
% denoted  $ \wt \Omega $ in Theorem \ref{thm:main0}.
Thus $\varepsilon A(i_\infty(\omega t))$,  where $A$ 
is defined in \eqref{actionanglevar},
is a quasi-periodic traveling wave solution of the water waves equations 
\eqref{eq:113}-\eqref{HamiltonianWW}.
In order to conclude the proof of 
Theorem \ref{thm:main0} it remains to show the linear stability of the quasi-periodic solution.
This is done in the rest of this section (see Theorem \ref{thm:lin stab}).

In order to show the linear stability we first shall need to show that the second Melnikov conditions are satisfies at the torus $i_\infty$. Let us define 
\begin{equation*}
%\label{Melnikov toro finale}
{\mathcal C}(i_\infty) := \big\{ \mathtt h \in [\mathtt h_1, \mathtt h_2] : 
(\alpha_\infty^{- 1}(\bar \omega(\mathtt h), \mathtt h), \mathtt{h}) 
\in {\bf \Lambda}_\infty^\gamma(i_\infty) \big\}  
\end{equation*} 
(recall \eqref{tDtCn}, \eqref{Melnikov-invert}, \eqref{Cantor set}, \eqref{prime.di.Melnikov}).  We shall prove the following lemma. 
\begin{lemma}
One has the inclusion ${\mathcal C}_\e \subseteq {\mathcal C}(i_\infty)$ (recall \eqref{Gvare2}).
\end{lemma}
\begin{proof}
We need to show that if $\mathtt h \in {\mathcal C}_\e$, then 
\[
\alpha_\infty^{- 1}(\bar \omega(\mathtt h), \mathtt h) \in {\mathtt T}{\mathtt C}_\infty(\gamma, \tau)(i_\infty) \cap \Lambda_\infty^{\gamma, I}(i_\infty) \cap \Lambda_\infty^\gamma(i_\infty)\,. 
\]
First we shall prove that 
\[
\alpha_\infty^{- 1}(\bar \omega(\mathtt h), \mathtt h) \in  {\mathtt T}{\mathtt C}_\infty(\gamma, \tau)(i_\infty), 
\]
namely 
\begin{equation}\label{zero melnikov cantor finale}
\begin{aligned}
\Big|(\bar \omega_\e(\mathtt h) - \mathtt V \mathtt m_1(i_\infty)) \cdot \ell \Big| \geq 
4\gamma\langle\ell\rangle^{-\tau}\,,
%\frac{4 \gamma}{\langle \ell \rangle^\tau}\,, 
\qquad \forall \ell \neq 0\,. 
\end{aligned}
\end{equation}
Indeed, if $\mathtt h \in {\mathcal C}_\e$, one has that for any 
$n \geq 0$, $\ell \in \Z^\nu \setminus \{ 0 \}$, 
\[
\Big|(\bar \omega_\e(\mathtt h) 
- \mathtt V \mathtt m_1(i_n)) \cdot \ell \Big| \geq 
4\gamma_{n}\langle\ell\rangle^{-\tau}\,,
\qquad \gamma_n = \gamma(1 + 2^{- n})\,.
%\frac{4 \gamma_n}{\langle \ell \rangle^\tau}. 
\]
Let us fix some $\ell \in \Z^\nu \setminus \{ 0 \}$. Therefore, 
there exists $\bar n \in \N$ such that for any $n \geq \bar n$, $|\ell| \leq N_n$. As a consequence 
one has that 
\[
\begin{aligned}
\big|(\bar \omega_\e(\mathtt h) - \mathtt V \mathtt m_1(i_\infty)) \cdot \ell \big| 
& \geq 
\big|(\bar \omega_\e(\mathtt h) - \mathtt V \mathtt m_1(i_n)) \cdot \ell \big| 
- \big|\mathtt V(\mathtt m_1(i_\infty) - \mathtt m_1(i_n)) \cdot \ell \big|  
\\& 
\stackrel{\eqref{beta.FGMP.est}}{\geq} 
\frac{4 \gamma_n}{\langle \ell \rangle^\tau } 
- C_1\e N_n \| i_n - i_\infty \|_{s_0  + \mu({\mathtt b}_1) + \sigma_{1}}
 \\& 
  \stackrel{\eqref{U infty - U n}, \e \gamma^{- 1} \leq 1}{\geq} 
  \frac{4 \gamma_n}{\langle \ell \rangle^\tau } 
  - C_1\e N_n^{1 - \mathtt a_3}  
  \geq 
  \frac{4 \gamma}{\langle \ell \rangle^\tau} \,,
\end{aligned}
\]
since 
%(recall that $\gamma_n = \gamma(1 + 2^{- n})$) 
\[
2^n (4 \gamma)^{- 1}C_1 \e N_n^{1 - \mathtt a_3} \langle \ell \rangle^\tau \leq 1\,,
\]
and since $|\ell| \leq N_n$ for any $n \geq \bar n$. The latter condition is implied by
\[
2^n (4 \gamma)^{- 1}C_1 \e N_n^{\tau + 1 - \mathtt a_3} \leq 1\,. 
\]
This  condition above 
is satisfied for any $n \geq 0$, since by 
\eqref{costanti nash moser}-\eqref{costanti nash moser 2}, 
$\mathtt a_3 > \tau + 1$ and $\e \gamma^{- 1} \ll 1$. 
This proves \eqref{zero melnikov cantor finale}. 

\medskip

\noindent
We now show that $\alpha_\infty^{- 1}(\bar \omega(\mathtt h), \mathtt h) \in {\bf \Lambda}_\infty^\gamma(i_\infty)$, namely
\begin{equation*}%\label{seconde Melnikov cantor finale}
\begin{aligned}
& 
|\bar \omega_\e(\mathtt h) \cdot \ell + \mu_\infty(j; i_\infty) \pm  \mu_\infty(j'; i_\infty)| 
\geq \frac{2  \,\gamma}{\langle \ell \rangle^\tau | j' |^\tau } \,, 
\\
 &  \forall \ell \in \Z^\nu \setminus \{ 0 \}\,, \, 	j,j' \in S_0^c\,, \ \  {\mathtt V}^T\ell  + j \pm j' =0\,. 
\end{aligned}
\end{equation*}
We actually show the estimate with the $-$ sign. The one with $+$ is similar. 
Since $\mathtt h \in {\mathcal C}_\e$ then 
\begin{equation}\label{sec mel tutti i tori}
\begin{aligned}
& |\omega \cdot \ell + \mu_\infty(j; i_n) -  \mu_\infty(j'; i_n)| 
\geq 
\frac{2  \,\gamma}{\langle \ell \rangle^\tau | j' |^\tau } \,, \quad \forall n \geq 0\,,
\\
 &  \forall \ell \in \Z^\nu \setminus \{ 0 \}\,, \, 	j,j' \in S_0^c\,, \ \  {\mathtt V}^T\ell  + j - j' =0\,. 
\end{aligned}
\end{equation}
Let $\ell \neq 0$. Then there exists (as before) 
$\bar n \in \N$ such that for any $n \geq \bar n$, $|\ell| \leq N_n$. 
We distinguish two cases:

\smallskip
\noindent
{\bf Case 1 .} Let $j, j' \in S_0^c$, $|j'| \geq N_n^q$, $q = 2 (\tau + 2)$, $\mathtt V^T \ell + j - j' = '0$. 
Then  one has
\begin{equation*}%\label{claim 3 inclusione cantor riducibilita aa}
|\overline \omega_\e(\mathtt h) \cdot \ell + \mu_\infty(j;  i_\infty) - \mu_\infty(j'; i_\infty)| \geq \frac{2 \gamma}{\langle \ell \rangle^\tau |j'|^\tau}\,. 
\end{equation*}
%This implies that ${\mathcal R}_n^{(2, -)}(\ell, j, j') = \emptyset$.

\medskip

\noindent
{\bf Proof of the claim 1.}
We analyze 
\[
\delta_\infty(j) := \mu_\infty(j; i_\infty) - \mathtt{m}_{1}(i_\infty) \cdot j \,,
\]
for $j, j' \in S_0^c$ and $|j'| \geq N_n^q$, $q := 2 (\tau + 2)$, $|\ell| \leq N_n$, $\mathtt V^T \ell + j - j' = 0$, 
arguing as in the proof of \eqref{claim 2 S3 n+1a}.
One gets the bound 
\begin{equation}\label{claim 2 S3 n+1aa}
|\delta_\infty(j) - \delta_\infty(j')| \lesssim N_n^{- \frac{q}{2} + 1}\,.
\end{equation}
Moreover, by \eqref{zero melnikov cantor finale},  one has that 
\[
\begin{aligned}
|\overline \omega_\e(\mathtt h) \cdot \ell + \mu_\infty(j; i_\infty) - \mu_\infty(j'; i_\infty)| 
& \geq 
|\overline \omega_\e(\mathtt h) \cdot \ell 
+ \mathtt{m}_{1}(i_\infty) \cdot (j - j')| - |\delta_\infty(j) - \delta_\infty(j')|  
\\& 
\stackrel{\eqref{claim 2 S3 n+1aa}, \eqref{zero melnikov cantor finale}}{\geq} 
\frac{4 \gamma}{\langle \ell \rangle^\tau} 
- \frac{C}{N_n^{\frac{q}{2} - 1}}  
\geq 
\frac{2 \gamma}{\langle \ell \rangle^\tau}\,,
\end{aligned}
\]
provided that 
\begin{equation*}%\label{tecnocasa merdaa a}
N_n^{\frac{q}{2} - \tau - 1} \gamma  \gg 1, \quad \text{for any} \quad n \geq 0\,.
\end{equation*} 
Since, $q = 2(\tau + 2)$, one has that $\frac{q}{2} - \tau - 1 = 1$, 
therefore the condition \eqref{tecnocasa merdaa} is fulfilled since  
$N_0 = \gamma^{- 2}$ and $0 < \gamma \ll 1$ (see \eqref{nash moser smallness condition}). 
%This implies \eqref{claim 3 inclusione cantor riducibilita aa}.
The proof of the {\bf case 1} is then concluded. 

\medskip
\noindent
{\bf Case 2.}  For $j, j' \in S_0^c$, $|j'| \leq N_n^q$, $q = 2(\tau + 2)$, $\mathtt V^T \ell + j - j' = 0$. Then 
\begin{equation*}%\label{claim 3 inclusione cantor riducibilita ba}
|\overline \omega_\e(\mathtt h) \cdot \ell + \mu_\infty(j; i_\infty) - \mu_\infty(j'; i_\infty)| \geq \frac{2 \gamma}{\langle \ell \rangle^\tau |j'|^\tau}\,. 
\end{equation*}
%This implies that ${\mathcal R}_n^{(2, -)}(\ell, j, j') = \emptyset$. 

\medskip

\noindent
{\bf Proof of the case 2.} 
By the estimates \eqref{stima delta 12 mu 0 pre rid}, \eqref{stime forma normale limitea} (applied with 
$i_1 \equiv i_n$, $i_2 \equiv i_\infty$),  
for any value of the parameter 
$\mathtt h \in [\mathtt h_1, \mathtt h_2]$ and for any $j \in S_0^c$ (recalling that 
$\kappa_1 = 2(\mathtt m - \mathtt c_1) + 2$), using that 
$\gamma^{- 1} \leq \gamma^{- 2} = N_0 \leq N_n$, one has that 
\begin{equation*}%\label{L infty vn - v n - 1}
\begin{aligned}
|\mu_\infty(j; i_n) - \mu_\infty(j; i_\infty)| 
& \lesssim \e \gamma^{- \kappa_1}  |j| \| i_n - i_\infty\|_{s_0 + \mu(\mathtt b_1) + \sigma_1} 
\\
& \stackrel{\eqref{U infty - U n}}{\lesssim} \e \gamma^{- \kappa_1 - 1}  |j|   N_{n }^{- \mathtt a_3}    \lesssim N_n^{\kappa_1 + 1} N_{n }^{- \mathtt a_3} \e |j| \,.
\end{aligned}
\end{equation*}
If $|j'| \leq N_n^q$, $|\ell| \leq N_n$, $\mathtt V^T \ell + j - j' = 0$, one has that 
\[
|j| \leq |j'| + |j - j'| \lesssim |j'| + |\ell| \lesssim N_n^q + N_n \lesssim N_n^q
\]
hence 
\begin{equation*}%\label{L infty vn - v n - 1}
\begin{aligned}
 |\mu_\infty(j; i_n) - \mu_\infty(j; i_\infty)|  
&\lesssim N_n^{q + \kappa_1 + 1} N_{n}^{- \mathtt a_3} \e  
\\
  |\mu_\infty(j'; i_n) - \mu_\infty(j'; i_\infty)| 
  &\lesssim N_n^{q + \kappa_1 + 1} N_{n }^{- \mathtt a_3} \e\,. 
\end{aligned}
\end{equation*}
Thus, by using \eqref{sec mel tutti i tori}, we get, for some constant $C > 0$ large enough, 
\begin{equation*} %\label{bound bla so ri 0a}
		\begin{aligned}
		  | \overline \omega_\e(\mathtt h) &\cdot \ell  + \mu_\infty(j; i_\infty) - \mu_\infty(j'; i_\infty)|  \\
		& \geq | \overline \omega_\e(\mathtt h) \cdot \ell  + \mu_\infty(j; i_n) - \mu_\infty(j';i_n )| \\
		& \quad - 
		|\mu_\infty(j; i_n) - \mu_\infty(j; i_\infty)|  - 
		|\mu_\infty(j'; i_n) - \mu_\infty(j'; i_\infty)|  
		\\
		& \geq 
		\frac{2  \gamma_n}{\langle \ell \rangle^\tau |j'|^\tau} - C N_n^{\kappa_1 + q + 1 - \mathtt a_3}    \e  
		 \geq \frac{2  \gamma}{\langle \ell \rangle^\tau |j'|^\tau} \,, 
	\end{aligned}
\end{equation*}
provided (use that $\gamma_n - \gamma = \gamma 2^{- n}$)
\[
C 2^{n - 1}  N_n^{\kappa_1 + q + 1- \mathtt a_3} \langle \ell \rangle^\tau |j'|^\tau  \e \gamma^{- 1} \leq 1 \,.
\]
Using that $|\ell| \leq  N_{n}$, $|j'| \leq N_n^q$, $\gamma^{- 1} \leq \gamma^{- 2} = N_{0} \leq N_{n}$, 
the latter condition is implied by 
\begin{equation}\label{giobbe1a}
C' 2^{n - 1} N_{n}^{\tau + \tau q + 1} N_n^{\kappa_1 + q + 1}  N_{n }^{- \mathtt a_3}  \e    \leq 1\,, 
\quad n \geq 0\,, 
\end{equation}
for some constant $C' \geq C$. The condition \eqref{giobbe1a} is verified by \eqref{costanti.riducibilita}, \eqref{costanti nash moser}, \eqref{costanti nash moser 1} (which imply that 
$\mathtt a_3 > \chi(\tau + (\tau + 1) q + \kappa_1  + 2 ) > \tau + (\tau + 1) q + \kappa_1  + 2 \,$)
and the smallness condition \eqref{nash moser smallness condition}.  
The proof of the {\bf case 2} is then concluded. 

\medskip
\noindent
The fact that
\[
\mathtt h \in {\mathcal C}_\e \quad \Longrightarrow \quad 
\alpha_\infty^{- 1}(\bar \omega(\mathtt h), \mathtt h) \in \Lambda_\infty^{\gamma, I}(i_\infty)\,,
\]
can be done by similar arguments, hence it is omitted. The proof of the lemma is then concluded. 
\end{proof}
We prove that around each torus filled by the 
quasi-periodic solutions \eqref{QP:soluz}  
of the Hamiltonian system \eqref{eq:113} 
constructed in Theorem \ref{thm:main0}
there exist symplectic coordinates 
$ (\phi, y, \mathtt{w}) = (\phi, y, \eta, \psi) \in \T^\nu \times \R^\nu \times  H_{S}^\bot  $ 
(see \eqref{Gdelta}) 
in which the water waves Hamiltonian  reads
\begin{align}\label{weak-KAM-normal-form}
\widetilde\om \cdot y  + \frac12 K_{2 0}(\phi) y \cdot y +  \big( K_{11}(\phi) y , \mathtt{w}\big)_{L^2} 
+ \frac12 \big(K_{02}(\phi) \mathtt{w} , \mathtt{w} \big)_{L^2} + K_{\geq 3}(\phi, y, \mathtt{w}) 
\end{align}
where $ K_{\geq 3} $ collects the terms at least cubic in the variables $ (y, \mathtt{w} )$
(see \eqref{taylor_Kalpha} and note that, at a solution, 
one has $ \partial_\phi K_{00} = 0 $, $ K_{10} = \widetilde \omega $, 
$ K_{01} = 0 $ by Lemma %\ref{coefficienti nuovi}
$5.4$ in \cite{BBHM}). 
The $ (\phi, y) $ coordinates are modifications of the action-angle variables and 
$ \mathtt{w} $ is a translation of the cartesian variable $ z$ in the normal subspace, 
see \eqref{Gdelta}.
In these coordinates the quasi-periodic solution reads 
$ t \mapsto (\widetilde\om t , 0, 0 ) $ 
and  the corresponding  linearized water waves equations are
\begin{equation}\label{linear-torus-new coordinates}
\begin{cases}
\dot{ \phi} = K_{20}(\widetilde\omega t)[ y] + K_{11}^T(\widetilde\omega t)[ \mathtt{w}] \\
\dot{ y} = 0 \\
\dot{ \mathtt{w}} = J K_{02}(\widetilde\omega t)[ \mathtt{w}] + J K_{11}(\widetilde\omega t)[ y]\, .
\end{cases}
\end{equation}
The self-adjoint operator $  K_{02} (\widetilde\om t) $ is defined in 
%\eqref{KHG} 
\eqref{K 02}
and $ J K_{02} (\widetilde\om t) $ is the restriction to $ H_{S}^\bot $ of the linearized
water waves vector field  $ J \pa_u \nabla_u H (u (\widetilde\om t ))$ 
up to a finite dimensional remainder, see Lemma \ref{thm:Lin+FBR}.

\noindent
We have the following result of linear stability for the quasi-periodic solutions provided by Theorem \ref{thm:main0}. 
\begin{thm}{ \bf (Linear stability)} \label{thm:lin stab}
The quasi-periodic solutions \eqref{QP:soluz} of the pure gravity water waves system are linearly stable, meaning that for all $s$ belonging to a suitable interval $[s_1,s_2]$, for any initial datum 
\[
y(0) \in \R^\nu\,, \quad \mathtt{w}(0) \in H^{s - \frac14}_x \times H^{s + \frac14}_x\,, 
\]
the solutions $y(t)$, $\mathtt{w}(t)$ of system \eqref{linear-torus-new coordinates} satisfy 
\[
y(t) = y(0)\,, \qquad  
\| \mathtt{w} (t)\|_{H^{s - \frac14}_x \times H^{s + \frac14}_x} 
\leq 
C \big( \| \mathtt{w}(0)\|_{H^{s - \frac14}_x \times H^{s + \frac14}_x} + |y(0)| \big)\,, 
\quad \forall t\in\R\,. %\quad \forall s \in [s_1, s_2]\,. 
\]
\end{thm}

\begin{proof}
By \eqref{linear-torus-new coordinates},  
the actions $ y (t) = y(0) $ do not evolve in time and 
the third equation reduces to  the linear PDE
\begin{equation}\label{San pietroburgo modi normali}
\dot{\mathtt{w}} = J K_{02}(\widetilde\omega t)[ \mathtt{w}] + J K_{11}(\widetilde\omega t)[ y (0)] \, .
\end{equation}
Sections \ref{linearizzato siti normali}, \ref{sym.low.order}, \ref{sec:redulower} 
imply the existence of a transformation $(H^s_x \times H^s_x) \cap H_{S}^\bot  \to (H^{s - \frac14}_x \times H^{s + \frac14}_x) \cap H_{S}^\bot $ (see \eqref{def.Psi.barn}), 
bounded and invertible for all % on Sobolev spaces in a range 
$s \in [s_1, s_2]$,
% see \eqref{semiconiugio cal L8}, \eqref{final semi conjugation}, 
such that, in the new variables ${\mathtt w}_\infty$, 
the homogeneous equation $\dot{\mathtt{w}} = J K_{02}(\widetilde\omega t)[ \mathtt{w}]$ transforms into 
a system of infinitely many uncoupled % complex 
scalar and time independent ODEs of the form
\begin{equation}\label{san pietroburgo siti normali ridotta}
\partial_t {\mathtt w}_{\infty, j} = - \ii \mu_\infty(j) {\mathtt w}_{\infty, j} \, , 
\quad \forall j \in S_0^c \, , 
\end{equation}
where 
the eigenvalues $\mu_{\infty}(j)$ are (see \eqref{def cal N infty nel lemma}, \eqref{proprieta cal D bot}, \eqref{beta.FGMP.est})
\[
 \mu_{\infty}(j) := \mathtt{m}_{1} \cdot j + \mathtt{m}_{7}(j)+r_{\infty}(j)
\]
with $\mathtt{m}_1=O(\e)$, $\mathtt{m}_{7}$
 is a symbol in $S^{\frac{1}{2}}$ of size $\sim 1+O(\e \gamma^{-2(\mathtt{m}-\mathtt{c}_1)-1})$ 
 (see \eqref{stima lambda 7}),
 and (recall \eqref{stime forma normale limite})
$\sup_{j \in { S}_0^c} |j|^{{\mathtt m} - \mathtt{c}_{1}} |r_{\infty}(j)| = 
O(  \e {\gamma^{-2-2(\mathtt{m}-\mathtt{c}_1)}} ) $ 
see \eqref{costanti.riducibilita}. 
The above result is the {\it reducibility} of the linearized quasi-periodically time dependent equation 
$\dot{\mathtt{w}} = J K_{02}(\widetilde\omega t)[\mathtt{w}]$. 
The {\it Floquet exponents} of the quasi-periodic solution 
are the purely imaginary numbers 
$ \{ 0,  \ii \mu_{\infty}(j), j \in { S}_0^c  \}$
(the null Floquet exponent comes from the action component $\dot y = 0$).
Since $\mu_\infty(j)$ are real, the Sobolev norms of the solutions of 
\eqref{san pietroburgo siti normali ridotta} are constant.
This proves the linear stability result.
\end{proof}

\appendix

\section{Technical Lemmata}

\subsection{Symplectic flows and truncations}\label{sec:algebraflusso}

In this paper we shall consider the flow of the equation

\begin{equation}
\begin{cases}\label{natale}
\partial_\theta \Psi^\theta u= \im \Pi^\perp_S \big[\opw(a( \tau; \varphi, x, 	\xi))[\Pi^\perp_S \Psi^\theta u]\big]\\
\Psi^0 u=u
\end{cases}
\end{equation}
where $\Pi_{S}^{\perp}$ is defined in \eqref{projComplex} and $a$ is a symbol of the form 
\eqref{egogenerator} or \eqref{egogenerator2}.
In particular, we shall compare the map $\Psi^\tau, \tau\in[0,1]$ with the flow $\Phi^\tau$ of the equation
\eqref{flussosimboloGenerico}.
%Recalling \eqref{decomp siti tangenziali coordinate complesse}
%we consider the  symplectic form in the extended phase space 
%$(\vphi,Q,Z)\in \mathbb{T}^{\nu}\times\mathbb{R}^{\nu}\times {\bf H}_S^{\perp}$,
%(recall \eqref{poissonBraComp})
%\begin{align}
%\Omega_{e}(\varphi, Q, W)&:=d Q\wedge d\varphi -\ii d z\wedge d\bar{z}\,,\nonumber\\
%\{F_{\mathbb{C}}, G_{\mathbb{C}}\}_{e}&:=
%\partial_{Q} F_{\mathbb{C}} \partial_{\vphi} G_{\mathbb{C}}
%-\partial_{\vphi} F_{\mathbb{C}} \partial_{Q} G_{\mathbb{C}}+
%\{F_{\mathbb{C}},G_{\mathbb{C}}\}\label{estesaPois}
%\\
%&\;=
%\partial_{\varphi} F_{\mathbb{C}} \partial_{\eta} G_{\mathbb{C}}
%-\partial_{\eta} F_{\mathbb{C}} \partial_{\varphi} G_{\mathbb{C}}-\ii
%\int_{\mathbb{T}_{\Gamma}^{d}}
%(\nabla_{z}G_{\mathbb{C}}\nabla_{\bar{z}}F_{\mathbb{C}}
%-\nabla_{\bar{z}}G_{\mathbb{C}}\nabla_{{z}}F_{\mathbb{C}})dx\,.
%\nonumber
%\end{align}
We remark that
 the map $\Psi^\tau$ in \eqref{natale} (if well posed) is symplectic
on the restricted subspace.
Indeed it is the flow 
of the vector field 
associated %(w.r.t. the symplectic form in \eqref{estesaPois})
to the Hamiltonian
\[
S(\tau; \varphi, z)=\frac{1}{2}\int_{\T_{\Gamma}^d} 
\opw (a(\tau; \varphi, x, \xi))z\cdot \bar{z} \,dx, \quad Z=\vect{z}{\bar{z}}\in {\bf H}^\perp_S\,\,.
\]
In the following we study the properties of the flow in \eqref{natale}. Let us fix 
\[
s_0:=(\nu+d)/2+1\,.
\]

\begin{lemma}\label{differenzaFlussi}
Let $a$ satisfying the assumption of Lemma \ref{lemma:buonaposFlussi}.
Fix  $\rho,\beta_0\in \N$ with $\rho\geq \beta_0+k_0$.
%and $p\geq s_0$ then the following holds.
For any $m_1, m_2 \in \Z$  such that $m_1+m_2\geq \rho-(\beta_0+k_0)$, 
there exist $\s_1=\s_1(m_1,m_2,\rho,\beta_0,k_0)\geq \s_2(m_1,m_2,\rho,\beta_0,k_0)$ such that, 
for any $\bar{s}\geq s_0+\s_1$, %$s_0 > \frac{d + \nu}{2}$,
there exists $\delta=\delta(\s_1,\bar{s},\beta_0,k_0)$ such that 
 such that if
\begin{equation}\label{flow1}
\|a\|_{m,s_0+\s_1,\s_1}^{k_0,\gamma} \le \delta
\end{equation} 
then the following holds. 
Let $\Phi^{\tau}$ be the flow of the system \eqref{flussosimboloGenerico}
then the flow  $\Psi^{\tau}$ of \eqref{natale}
 is well defined 
for $|\tau|\le 1$ 
and one has
$\Psi^1= \Pi_S^\perp\Phi^{1}\Pi_S^\perp\circ(\uno +\mathcal{R})$
where $\mathcal{R}$ is an operator 
with the form 
\begin{equation}\label{FiniteDimFormDP10001}
\mathcal{R}(\varphi) z=\sum_{\lvert j \rvert\le C} \int_0^1 (w, g_j(\tau; \varphi))_{L^2}\,\chi_j(\tau; \varphi)\,d\tau\,,
\end{equation}
with $g_{j}^{(i)}, \chi_{j}^{(i)}$ functions satisfying 
$\lVert g^{(i)}_j \rVert_s^{k_0, \gamma}
+\lVert \chi_j^{(i)} \rVert_s^{k_0, \gamma}\lesssim_{s,m_1,m_2} 1$ and Lipschitz in the variable $i$.
Furthermore
the remainder $\cR$ satisfies the following:
%
%, let $\bar{s}>s_0, \tb\in\N$ and $\rho>(\tb+k_0)$ 
%then
% for any $s_0\le s\le \bar{s}$ the operator $\cR$ satisfies the following: \\

\noindent
$\bullet$ For any  $m_1+m_2\geq \rho-\beta_0-k_0$, 
and for any $b\in\N^{\nu}$ such that $|b|\leq \beta_0$, 
the operator $\jap{D}^{m_1} \pa_\vphi^{b} \cR \jap{D}^{m_2}$ 
is $\cD^{k_0}$-tame with tame constants satisfying (recall 
Def. \ref{Dksigmatame})
%\eqref{Mdritto})
\begin{equation}\label{giallo2}
%\mathbb{M}^{\gamma}_{\mathcal{R}}(s, \tb)
\mathfrak{M}_{\jap{D}^{m_1} \pa_\vphi^{b} \cR \jap{D}^{m_2}}(s)
\lesssim_{s,m_1,m_2}  \|a \|^{k_0, \gamma}_{m,s+\s_1,\s_1}\,.
\end{equation}

\noindent
$\bullet$ For any 
$m_{1}+m_{2}=\rho-\beta_0-1$ , the operator
$\jap{D}^{m_1} \pa_\vphi^b \Delta_{12} \cR \jap{D}^{m_2}$ satisfies 
\begin{equation}\label{giallo3}
\sup_{|b|\leq \beta_0}\sup_{m_1+m_2=\rho-\beta_0-1}
\|\jap{D}^{m_1} \pa_\vphi^b \Delta_{12} \cR \jap{D}^{m_2}\|_{\mathcal{L}(H^p;H^p)} \lesssim_{p,m_1,m_2}
\|\Delta_{12}a\|_{m,p+\s_2,\alpha}\,.
\end{equation}
\end{lemma}
\begin{proof}
One can follow almost word by word the proof of Lemma $C.1$ in \cite{FGP20}.
\end{proof}

As a consequence of Lemma \ref{differenzaFlussi} we have the following
\begin{lemma}\label{georgiaLem}
Fix $\beta_0\in \N$, $\mathtt{c}\geq m(\beta_0+k_0)+2$, and $M\geq \mathtt{c}$.
%Fix $N, \tb \in \N$ with
%$N- m (\tb+k_{0})\geq 2$ and $\alpha\in\N_0$.
Consider the operator
\[
\cL:=\Pi_{S}^{\perp}\mathcal{A}\Pi_{S}^{\perp}\,,\qquad \mathcal{A}:=\omega\cdot\pa_{\vphi}+\opw(d(\vphi,x,\x))+\mathcal{Q}
\]
where $\omega\in\mathtt{\Omega}\subseteq\R^{\nu}$,
 $d\in S^{n}$, 
$ n\in \R$ 
with  $d=d(\lambda, i(\lambda))$   is 
 $k_{0}$-times 
differentiable in $\lambda\in\mathtt{\Lambda}_{0}\subset\R^{\nu+1}$, 
Lipschitz in the variable $i$ and purely imaginary, and 
the operator $\mathcal{Q}$ has the  following properties:

\begin{itemize}
\item[(a)] for any $|\beta|\leq \beta_0$
% for any $M\leq N-\alpha(\tb+k_0)-2$ and any
%$b\in\N_0^{\nu}$ with $|b|\le \tb$,
the operator

\noindent
$\jap{D}^{M}\partial_{\vphi}^{\beta} \cQ \jap{D}^{-\mathtt{c}}$ is $\cD^{k_{0}}$-tame
(see Def. \ref{Dksigmatame}) with tame constants bounded for $s_0\leq s\leq \bar{s}$ for some $\bar{s}\gg1$
large;

\item[(b)] $\cQ=\cQ(i)$ depends  in a Lipschitz way on a  parameter $i$. 
For any $\beta\in\Z^{\nu}$ with $|\beta|\le \beta_{0}$, an some  $s_0\leq p\leq \bar{s}$
one has
$\jap{D}^{M}\partial_{\vphi}^{\beta} \Delta_{12}\cQ  \jap{D}^{-\mathtt{c}}\in\cL(H^{p}, H^{p})$.
%
% For  $M \le N-\alpha\tb-2$,
%any  $|b|\le \beta_{0}$, and   $s_0\leq p\leq \widehat{s}$
%one has
%$\partial_{\vphi}^{b}\jap{D}^{M} \Delta_{12}\cR  \jap{D}^{-\alpha\beta_0-2}\in\cL(H^{p}, H^{p})$.
\end{itemize} 
Let $\Psi^{\tau}, \Phi^{\tau}$ be respectively the flows of \eqref{natale} and \eqref{flussosimboloGenerico} 
with generator as in 
\eqref{egogenerator} or \eqref{egogenerator2}.
There exist constants
 $\s_{1}=\s_{1}(\beta_0, M,  k_{0},n,m)\geq \s_{2}=\s_{2}(\beta_0, M,n,m)$ 
such that if  $\bar{s}\geq s_0+\s_1$ the following holds. 
There exists $\delta=\delta(M,\beta_0, n, m,\bar{s})$ 
 such that if
\begin{equation*}%\label{flow1bertier}
\|a\|_{m,s_0+\s_1,\s_1}^{k_0,\gamma}+\|d\|_{n,s_0+\s_1,\s_1}^{k_0,\gamma}+
\mathbb{M}_{\mathcal{Q}}(s_0+\s_1, \beta_0) \le \delta
\end{equation*} 
where
\[
\mathbb{M}_{\mathcal{Q}}(s, \beta_0):=
\sup_{\substack{
0\leq |\beta|<\beta_0}}\mathfrak{M}_{\jap{D}^{M} \pa_\vphi^{\beta}  \cQ \jap{D}^{-\mathtt{c}}}(s)\,,
\]
then the following holds. 
One has
\begin{equation}\label{georgia100}
\Psi^{1}\cA(\Psi^{1})^{-1}=\Pi_{S}^{\perp}\Phi^{1}\cA(\Phi^{1})^{-1}\Pi_{S}^{\perp}
+\mathcal{R}\,,
\end{equation}
where $\mathcal{R}$ is a finite rank operator of the form \eqref{FiniteDimFormDP10001}
satisfying the following.
\begin{itemize}
\item
For 
 any 
  $|\beta|\le \beta_0$, one has that the operator \\
 $\jap{D}^{M} \partial_{\vphi}^{\beta}\cR \jap{D}^{-\mathtt{c}}$ is 
 $\cD^{k_{0}}$-tame and satisfies, for any $s_0 \leq s \leq \bar{s}-\s_1$, the estimate
\begin{equation}\label{giallo4}
\begin{aligned}
&{\mathfrak M}_{\jap{D}^{ M} \partial_{\vphi}^{\beta}\cR \jap{D}^{- \mathtt{c}}}(s) 
\\&\qquad\qquad \quad\lesssim_{m,n, s, M, \beta_0} 
 \|a \|^{k_0, \gamma}_{m,s+\s_1,\alpha}+
\|a \|^{k_0, \gamma}_{m,s_0+\s_1,\alpha}
(\|d\|^{k_0, \gamma}_{n,s+\s_1,\alpha}+\mathbb{M}^{\gamma}_{\mathcal{Q}}(s, \beta_0))\,;
\end{aligned}
\end{equation}

\item for  any 
$|\beta|\le \beta_0$ we have that, for any $p+\s_2\leq s_0+\s_1$,
\begin{equation}\label{giallo5}
\begin{aligned}
&\|\jap{D}^{M} \pa_\vphi^{\beta}  \Delta_{12} \cQ \jap{D}^{-\mathtt{c} }\|_{\mathcal{L}(H^p;H^p)}
\lesssim_{m,n,p,M, \beta_0} \|\Delta_{12}a\|_{m,p+\s_2,\s_2}(1+C_{p})\,,
\\&
C_{p}:=\|\Delta_{12}a\|_{m,p+\s_2,\s_2}+\sup_{|\beta|\leq \beta_0}
\|\jap{D}^{M} \pa_\vphi^\beta \Delta_{12} \cQ \jap{D}^{-\mathtt{c}}\|_{\mathcal{L}(H^p;H^p)}\,.
\end{aligned}
\end{equation}
\end{itemize}
If $m\leq 0$ the statements holds with $m=0$.
\end{lemma}

\begin{proof}
Formula \eqref{georgia100} follows by Lemma \ref{differenzaFlussi}
writing $\Psi^1= \Pi_S^\perp\Phi^{1}\Pi_S^\perp\circ(\rm{I} +\mathcal{R})$.
The estimates on the remainder $\mathcal{R}$ follow by reasoning as in 
the proof of Lemma $C.2$ in \cite{FGP20}
and by using Lemmata \ref{passaggiostweyl}, \ref{passaggioweylst}
%and Remark \ref{standard/Weyl}
to pass from the standard quantization
to the Weyl quantization.
\end{proof}

\begin{rmk}\label{differenzaFlussiRMK}
The same result of Lemmata \ref{differenzaFlussi}, \ref{georgiaLem} holds also in the case
$d(\vphi,x,\x)$ is a matrix of symbols in $S^{m'}$ (recall \eqref{matricidisimboli}),
$m'\leq 0$
and $\mathcal{Q}$ is a matrix of smoothing operators.
\end{rmk}

%\subsection{Green's function estimates: concluded}\label{app:green}
%We collect the proofs of some technical result of section \ref{sec:DNsezione}.

%\subsubsection{Proof of Corollary \ref{stima of di green forzante qualunque}}\label{app:greencorollaryQualunque}
%\begin{proof}[{\bf Proof of Corollary \ref{stima of di green forzante qualunque}}]
%The proof is postponed to Appendix \ref{app:proof4.13}.

%\end{proof}

%
%\subsubsection{Proof of Lemma \ref{lemma stima cal A partial x y m cal K vphi y}}\label{app:green2}
%%We give the proofs of some lemmata stated in section \ref{sec:Problemastriscianuova}.
%In this section we give the proof of Lemma \ref{lemma stima cal A partial x y m cal K vphi y}.
%%\subsection{Proof of Lemma \ref{lemma stima cal A partial x y m cal K vphi y}}\label{app:lemma3.8}
%The  result is quite technical. We first provide estimates on the 
%operator $\cT_m(\vphi, y)$ in \eqref{op fondamentale stime tame resto DN} 
%in the spaces $\mathcal{H}^{s}$ (see \eqref{def sobolev striscia}). 
%Then we conclude the estimates in $\mathcal{O}^{s}$. 
%We have the following.

\subsection{Tame estimates of the parabolic flow}\label{sec:backwardflow}
In this section we %give the proof of Lemma \ref{stima derivate lambda vphi flusso parabolico}. 
prove tame estimates for the flow $\cU(y)$ of the PDE in \eqref{calore striscia}
where $a(\vphi, x, \xi)=a(\lambda; \vphi, x, \xi)\in S^{1}$ is $k_{0}-$times 
differentiable with respect to the parameters $\lambda=(\omega, \th)$. 
The symbol $a:=a(i)$ may depend also on the “approximate” torus $i(\vphi)$ in a Lipschitz way.
Moreover, suppose that the symbol $a$ has the following structure (recall \eqref{opcutofffunct}):
\begin{equation}\label{propsimboloa.calore}
\begin{aligned}
 a(\vphi, x,  \xi ) &= \mathtt h|\xi| \chi(\x)+ r(\vphi, x, \xi)\,,  \quad r \in S^1 \,, 
\\
 \| r \|_{1, s, \alpha}^{k_0, \gamma}& \lesssim_{s, \alpha} \| \eta \|_{s + \sigma}^{k_0, \gamma}\,, 
\quad \forall s \geq s_0:=(\nu+d)/2+1\,,
\\
 \| \Delta_{12}r \|_{1, p, \alpha} &\lesssim_{p,\alpha}\| \eta_1-\eta_{2} \|_{p + \sigma}^{k_0, \gamma}\,, 
 \;\;\;p\geq s_0\,.
\end{aligned}
\end{equation}

We look for the solution of  \eqref{calore striscia} by a Galerkin approximation, as limit of the solutions of the truncated equations
\begin{equation}\label{calore.striscia.N}
\begin{cases}
\partial_\tau \cU(\vphi, \tau) =  - \Pi_{N}{\rm Op}\big(a(\vphi, x, \xi) \big) \Pi_{N}\cU(\vphi, \tau)\,, 
\quad  0 \leq \tau \leq 1\,, \\
\cU(\vphi, 0) = \Pi_{N} \,,
\end{cases}
\end{equation}
where, for any $N \in \N$,  we denote by 
$\Pi_N$  the $ L^2 $-orthogonal projector on the finite dimensional subspace
\[
E_N := \Big\{ u \in L^2(\T^d_\Gamma;\C) : 
u(x) = \sum_{|j| \leq N} u_j e^{\ii j\cdot x} \Big\} \, .
\]
%We  introduce the  ``paraproduct" decomposition 
%for the  product of two functions $ a, u : \T^{d}_{\Gamma} \to \C $,  
Given a function $ u : \T^{d}_{\Gamma} \to \C $,  we introduce the decomposition (recall \eqref{opcutofffunct})
\begin{equation}\label{paraproduct}
\op(a) u =\th|D| u+ T_r u + R_r u 
\end{equation}
where, recalling the notation in \eqref{def a chi fourier}-\eqref{bonyquantization},
\begin{equation}\label{Ta Ru}
T_r:=\opb(r)\,,\qquad R_r:=\opb(r-r_{\tilde\chi})\,.
\end{equation}
%\begin{equation}
%\begin{aligned} \label{Ta Ru}
%&  T_a u  := \sum_{k, \xi \in \Gamma^*\,,|k - \xi| \leq \frac12 |\xi|} \widehat a(k - \xi, \xi) \widehat u(\xi) e^{\ii k \cdot x}\,, \\
%&  R_a u  := \sum_{k, \xi \in \Gamma^*\,,|k - \xi| > \frac12  |\xi|} \widehat a(k - \xi, \xi) \widehat u(\xi) e^{\ii k \cdot x}= \sum_{k, \xi \in \Z\,,|k - \xi| <   2|\xi|} \widehat u(k - \xi) \widehat a(\xi, k-\xi) e^{\ii k\cdot x}\, .
%\end{aligned}
%\end{equation}
%Note that  
%\begin{equation}\label{simbolo.a0} 
%T_a = {\rm Op} ( a_0 (x, \xi) ) \quad {\rm with  } \quad
%a_0 (x, \xi) := {\mathop \sum}_{|\xi| \leq  2 |k-\xi|} \widehat a( \xi, k-\xi) e^{\ii k \cdot x } \, . 
%\end{equation}
%In order to do the energy estimates we need 
%to introduce a semi-norm on the class of symbols 
%$S^m$ that controls the regularity in $x$.
%
%\red{semplificare la def di $T_a$ vedi egorov. Le stime direttamente con $\|\eta\|$}
%\begin{defn}\label{norminaEasy}
%Let $a\in S^m$. For $\alpha\in\N, s\ge 0$ we define the norm
%\begin{equation}\label{norma simbolo solo in spazio}
%\|a\|_{m, H^s_x, \al}:= \max_{0\le |\beta| \le \alpha} \sup_{\xi\in\R^d} \|\pa_\xi^\beta a(\cdot, \xi)\|_{H^s_x} \jap{\xi}^{-m+\beta}\,.
%\end{equation}
%Analogously, we define the norm $\|a\|_{m, C^s_x, \alpha}:= \max_{0\le |\beta| \le \alpha} \sup_{\xi\in\R^d} \|\pa_\xi^\beta a(\cdot, \xi)\|_{C^s_x} \jap{\xi}^{-m+\beta}$. 
%\end{defn}

\noindent
%Let us fix $s_{1}\in \N$ with $s_{1}> \big\lfloor \frac{d}{2} \big\rfloor$.
By arguing as in Lemma \ref{lemma proprieta bony weil} 
(see also Lemma 2.21 in \cite{BM20}) and using \eqref{propsimboloa.calore}
one can prove the following estimates:
%for all $ s \geq 0   $,  following estimates  hold:
for all $s\in \R$ one has
\begin{equation}\label{stima Ta Ru1}
\| T_r u  \|_{H^s_x} \lesssim_{s} \| \eta \|_{s_0 + \sigma}^{k_0, \gamma} \| u \|_{H^{s+1}_x}\,, 
\end{equation}
and for any $s\geq 0$ one has
\begin{equation}\label{stima Ta Ru}
%\| T_a u  \|_{H^s_x} \leq C(s) \| a \|_{1, H^{s_{1}}_{x}, 0} \| u \|_{H^{s+1}_x}\,, \quad 
\| R_r u \|_{H_x^s} \lesssim_{s} \| \eta \|_{s + \sigma}^{k_0, \gamma} \| u \|_{L^{2}_x}  \,,
\end{equation}
(the operator $ u \mapsto R_r u $  is  smoothing)
for some $\s\gg1$.

\begin{rmk}\label{rmk:azionespaziOMO}
In view of the cut off $|k-\x|\leq 1/2|\x|$ in the definition of $T_r$ in \eqref{Ta Ru},
one can note that actually 
$T_r u=T_{r}\widetilde{u}$ where 
$\widetilde{u}:=u-u_0$, with $u_0=\frac{1}{(2\pi)^d}\int_{\T^{d}_{\Gamma}}u(x)dx$.
In other words the action of $T_r$ is invariant on the space of zero average functions.
%$T_{a}g\equiv0$ if $g$ is a constant function.

Therefore we note the following simple facts:

\noindent
$(i)$ %Writing $\widetilde{u}:=u-u_0$, with $u_0=\frac{1}{(2\pi)^d}\int_{\T^{d}_{\Gamma}}u(x)dx$, 
One has 
\[
(T_r u,u)_{L_x^{2}}=(T_r \widetilde{u}, \widetilde{u} )_{L^{2}_x}\,.
\]

\noindent
$(ii)$
For $s\in \R$, one has  
\begin{equation}\label{stimaspaziomo}
\| T_r u  \|_{\dot{H}^s_x} \leq C(s) \| \eta \|_{s_0 + \sigma}^{k_0, \gamma} \| u \|_{\dot{H}^{s+1}_x}\,, 
\qquad \forall\,u\in \dot{H}_{x}^{s+1}\,,
\end{equation}
where $\dot{H}_{x}^{s}$ is the standard  homogeneous Sobolev space. In particular, we recall that an equivalent norm on $\dot{H}_{x}^{s}$ is given by (recall \eqref{opcutofffunct})
\[
\|u\|_{\dot{H}_{x}^{s}}\sim_{s}\||D|^{s}u\|_{L^{2}_{x}}\,.
\]
\end{rmk}

%Indeed, 
%\begin{equation*}
%\begin{aligned}
%\| T_a u  \|^{2}_{H^s_x} 
%&= \sum_{k\in\Gamma^{*}} \left(  \sum_{|k-\xi|\le \frac12 |\xi|} \jap{k}^{s} \widehat{a}(k-\xi, \xi) \widehat{u}(\xi) \right)^{2}
%\\
%&\lesssim \sum_{k\in\Gamma^{*}} \left(  \sum_{|k-\xi|\le \frac12 |\xi|} \jap{\xi}^{-1} 
%|\widehat{a}(k-\xi, \xi)| \jap{\xi}^{s+1} |\widehat{u}(\xi)| \frac{ \jap{k-\xi}^{s_{1}}}{\jap{k-\xi}^{s_{1}}} 
%\right)^{2}
%\\
%&\lesssim_{s_{1}} \|a\|^{2}_{1, H_{x}^{s_{1}}, 0} \|u\|^{2}_{H_{x}^{s+1}}
%\end{aligned}
%\end{equation*}
%where in the last line we used Cauchy-Schwarz inequality. In the same way
%\begin{equation*}
%\begin{aligned}
%\| R_a u  \|^{2}_{H^s_x} 
%&= \sum_{k\in\Gamma^{*}} \left(  \sum_{|k-\xi|\le 2|\xi|} \jap{k}^{s} \widehat{u}(k-\xi) \widehat{a}(\xi, k-\xi) \right)^{2}
%\\
%&\lesssim \sum_{k\in\Gamma^{*}} \left(  \sum_{|k-\xi|\le 2|\xi|} \jap{\xi}^{s} \widehat{u}(k-\xi) \jap{k-\xi}^{-1} \jap{k-\xi}^{s_{1}+1} \frac{\widehat{a}(\xi, k-\xi)}{ \jap{k-\xi}^{s_{1}} } \right)^{2}
%\\
%&\lesssim_{s_{1}} \|a\|^{2}_{1, H_{x}^{s+s_{1}+1}, 0} \|u\|^{2}_{L_{x}^{2}}\,.
%\end{aligned}
%\end{equation*}

We need the following technical lemma.
\begin{lemma}\label{lemm:supergard}
Consider the symbol $a(\vphi,x,\x)$ in \eqref{propsimboloa.calore}.
There is $\s\gg1$ such that for any $s\geq 0$, there exist $C_s\gg1$, $\delta\ll1$ such that
the following holds.
If the smallness condition \eqref{smalleta} holds (for $\s_1\geq\s$) then
one has
\begin{equation}\label{eq:supergard}
\mathtt{h}\|\langle D\rangle^{s}h\|_{\dot{H}_x^{\frac{1}{2}}}^2
\leq 
2{\rm Re}(\langle D\rangle^{s}\op(a)h,\langle D\rangle^{s}h)_{L^2_x}+
C_s\Big(
\| h \|_{H^s_x}^2 
+ (\| \eta \|_{s + \sigma}^{k_0, \gamma})^2 \| h\|_{L^2_x  }^2\Big)\,,
\end{equation}
for any $h\in H^{s}_x$.
\end{lemma}
\begin{proof}
%First of all, 
%recalling Definition \ref{norminaEasy} and assumption
%\eqref{propsimboloa.calore} (see also \eqref{sobnormseparati} and \eqref{norma pseudodiff})
%one has
%\begin{equation}\label{piccolezzarpiccolo}
%\|r\|_{1,H_{x}^{s},\alpha}\lesssim_{s,\alpha}\|\eta\|_{s+\s}\,,
%\end{equation}
% for some $\s\gg 1$ large.

First of all we note that
\begin{align}
-2{\rm Re}(\langle D\rangle^{s}&\op(a)h,\langle D\rangle^{s}h)_{L^2_x} 
\nonumber
\\& =
-  \big( \langle D \rangle^s  \op( a) h , \langle D \rangle^s h \big)_{L^2_x}  
-  \big(\langle D \rangle^s  h ,\langle D \rangle^s  \op( a) h \big)_{L^2_x} \nonumber  \\ 
&\stackrel{\eqref{propsimboloa.calore}, \eqref{paraproduct}}{=} 
-  \th \big( \langle D \rangle^s  |D| h , \langle D \rangle^s h \big)_{L^2_x}  
-   \th  \big(\langle D \rangle^s  h ,\langle D \rangle^s |D| h \big)_{L^2_x} \label{termini0} \\ 
& -  \big( \langle D \rangle^s  T_{r} h , \langle D \rangle^s h \big)_{L^2_x}  
-  \big(\langle D \rangle^s  h ,\langle D \rangle^s  T_{r} h \big)_{L^2_x} \label{termini1}
\\
&-  \big(\langle D \rangle^s R_{r}  h ,\langle D \rangle^s h \big)_{L^2_x}
-  \big(\langle D \rangle^s  h ,\langle D \rangle^s R_{r} h \big)_{L^2_x} \label{termini2}\,. 
\end{align}    
Note that, recalling Remark \ref{rmk:azionespaziOMO}-$(ii)$, one has
\[
\eqref{termini0}=-2\mathtt{h}\|\langle D \rangle^s h\|_{\dot{H}_{x}^{\frac{1}{2}}}^2\,.
\]
We now provide estimates on the remaining terms in the inequality above.

\noindent
{\sc Estimate of \eqref{termini2}.}
Using Cauchy-Schwarz and Young inequalities we have that
\begin{equation}\label{terministima1}
\begin{aligned}
|\eqref{termini2}|\leq2\|R_{r}h\|_{H_{x}^{s}}\|h\|_{H_{x}^{s}}
&\stackrel{\eqref{stima Ta Ru}}{\lesssim_{s}}
\| \eta \|_{s + \sigma}^{k_0, \gamma} \|h\|_{L_{x}^2}\|h\|_{H_{x}^{s}}
\\&\lesssim_{s}(\| \eta \|_{s + \sigma}^{k_0, \gamma})^2\|h\|_{L_{x}^2}^2+\|h\|_{H_{x}^{s}}^2\,.
\end{aligned}
\end{equation}
\noindent
{\sc Estimate of \eqref{termini1}.} 
We start by proving that, given $v\in H^{s}$, the following holds for $s\geq0$:
\begin{equation}
 \big\| [\langle D \rangle^s, T_r ] v \big\|_{L^2_x} 
\lesssim_{s} \| \eta \|_{s_0 + \sigma}^{k_0, \gamma}  \|  v \|_{H^s_x} \label{lemma.commutatore.calore1}\,.
\end{equation}
We note that
\[
[\langle D\rangle^s, T_r ] v  = {\mathop \sum}_\xi 
{\mathop \sum}_{|j - \xi| \leq  \frac12 |\xi|} \psi(\xi, j) 
\widehat r(j - \xi, \xi) \jap{\xi}^{-1}\widehat v(\xi)  e^{\ii  j\cdot x } 
\]
where 
$ \psi(\xi, j) := (\langle j \rangle^s  - \langle\xi \rangle^s) \jap{\xi}   $. 
Since $|j - \xi| \leq \frac12 |\xi|$, 
then by triangular inequality $|j|\le \frac 32|\xi|$ and hence we have  
$ |\psi(\xi, j)| \lesssim_{s} \jap{j-\xi}  \langle \xi \rangle^{s}  $. 
Let $s_1:=d/2+1$. Then
using the Cauchy-Schwartz inequality we get
\begin{equation}\label{stimacommuTr}
\begin{aligned}
\| [\langle D \rangle^s, T_r ] v\|_{L^2_x}^2 & \leq
\sum_j \Big( \sum_{|j - \xi| \leq \frac12|\xi|} |\psi(\xi, j)| |\widehat r(j - \xi,\xi)| \jap{\xi}^{-1} |\widehat v(\xi)| \Big)^2 
\\
& \lesssim_{s}  \sum_j \Big( \sum_{|j - \xi| \leq \frac12|\xi|} \langle \xi \rangle^{s}  |\widehat r(j - \xi, \xi)| \jap{\xi}^{-1}
|\widehat v(\xi)| \frac{\langle j - \xi \rangle^{s_{1}+1}}{\langle j - \xi \rangle^{s_{1}}} \Big)^2  
\\
& \lesssim_{s} \sum_\xi \langle \xi \rangle^{2 s} |\widehat v(\xi)|^2 \sum_j \langle j - 
\xi \rangle^{2s_{1}+2} |\widehat r(j - \xi, \xi)|^2 \jap{\xi}^{-2} \\
& \lesssim_{s}
 %\| r \|_{1, H^{s_{1}+1}_x, 0}^2
( \sup_{\xi\in\R^d} \| r(\cdot, \xi)\|_{H^{s_1+1}_x} \jap{\xi}^{-1})^2
  \| v \|_{H^{s}_x}^2 	\, .
\end{aligned}
\end{equation}
Moreover, recalling \eqref{sobnormseparati} and \eqref{norma pseudodiff} and 
assumption
\eqref{propsimboloa.calore} one deduces 
\begin{equation*}%\label{piccolezzarpiccolo}
 \sup_{\xi\in\R^d} \| r(\cdot, \xi)\|_{H^s_x} \jap{\xi}^{-1}\lesssim_s\|r\|_{1,s,0}\lesssim_{s}\|\eta\|_{s+\s}\,,
\end{equation*}
 for some $\s\gg 1$ large. Therefore, the bound \eqref{stimacommuTr} implies 
 \eqref{lemma.commutatore.calore1}.
Now, recalling Remark \ref{rmk:azionespaziOMO} (setting $\widetilde{h}=h-h_0$), we note that
\begin{equation}\label{terministima2}
\begin{aligned}
|\eqref{termini1} |& \leq 
|\big( |D|^{-\frac{1}{2}}T_{r}  \langle D \rangle^s \widetilde{h} , 
|D|^{\frac{1}{2}}\langle D \rangle^s \widetilde{h} \big)_{L^2_x}  |
+
|\big([ \langle D \rangle^s , T_{r} ]h , \langle D \rangle^s h \big)_{L^2_x} |
\\
&+|\big(  |D|^{\frac{1}{2}}\langle D \rangle^s \widetilde{h} ,|D|^{-\frac{1}{2}}T_{r}  \langle D \rangle^s \widetilde{h} \big)_{L^2_x}  |
+
|\big( \langle D \rangle^s h, [ \langle D \rangle^s , T_{r} ]h  \big)_{L^2_x} |
\\&
\lesssim 
\| |D|^{-\frac{1}{2}}T_{r}  \langle D \rangle^s \widetilde{h}\|_{L^{2}_x}
\|\langle D\rangle^{s}h\|_{\dot{H}_x^{\frac{1}{2}}}
+\|h\|_{H_{x}^{s}}
\|[ \langle D \rangle^s , T_{r} ]h\|_{L^{2}_x}
\\&
\stackrel{\eqref{stimaspaziomo}, \eqref{lemma.commutatore.calore1}}{\lesssim_{s}}
%\|r\|_{1,H_{x}^{s_1},0}
\| \eta \|_{s_0 + \sigma}^{k_0, \gamma}  
\|\langle D\rangle^{s}h\|_{\dot{H}_x^{\frac{1}{2}}}^2
+\|h\|_{H_{x}^{s}}^2\| \eta \|_{s_0 + \sigma}^{k_0, \gamma}  \,.
%\|r\|_{1,H_{x}^{s_1+1},0}\,.
\end{aligned}
\end{equation}
By combining \eqref{termini0}-\eqref{termini2} with \eqref{terministima1}-\eqref{terministima2}
we have that there is some $C(s)\gg1 $ large such that 
\[
\begin{aligned}
-2{\rm Re}(\langle D\rangle^{s}\op(a)h,\langle D\rangle^{s}h)_{L^2_x}  & \leq
\Big(-2\mathtt{h}+C(s)\| \eta \|_{s_0 + \sigma}^{k_0, \gamma}  
%\|r\|_{1,H_{x}^{s_1},0}
\Big)
\|\langle D \rangle^s h\|_{\dot{H}_{x}^{\frac{1}{2}}}^2
\\&
+C(s)\big(1+\| \eta \|_{s_0 + \sigma}^{k_0, \gamma}  
%\| r\|_{1, H^{s_{1}+1}_x, 0}
\big)\|h\|_{H_{x}^{s}}^2
+C(s)(\| \eta \|_{s + \sigma}^{k_0, \gamma}  )^2
%\| r\|_{1, H^{s+s_{1}+1}_x, 0}^2
\|h\|_{L_{x}^2}^2\,.
%\\&
%\lesssim_{s}
%\| r\|_{1, H^{s+s_{1}+1}_x, 0}^2\|u_{N}\|_{L_{x}^2}^2+\|u_{N}\|_{H_{x}^{s}}^2\,.
\end{aligned}
\]
Therefore, %using also \eqref{piccolezzarpiccolo} and 
taking 
$\| \eta\|_{s_0+\s}^{k_0, \gamma}\ll1$ small enough w.r.t. $s$ and $\mathtt{h}$,
we get
\begin{equation*}%\label{uragano1}
\begin{aligned}
-2{\rm Re}(\langle D\rangle^{s}&\op(a)h,\langle D\rangle^{s}h)_{L^2_x}  
\\&\leq
-\mathtt{h}\|\langle D\rangle^{s}h\|_{\dot{H}_{x}^{\frac{1}{2}}}^2
+2C(s)\big(
\|h\|_{H_{x}^{s}}^2
+(\| \eta \|_{s + \sigma}^{k_0, \gamma})^2
%\| r\|_{1, H^{s+s_{1}+1}_x, 0}^2
\|h\|_{L_{x}^2}^2
\big)\,.
%\\&
%\lesssim_{s}
%\| r\|_{1, H^{s+s_{1}+1}_x, 0}^2\|u_{N}\|_{L_{x}^2}^2+\|u_{N}\|_{H_{x}^{s}}^2\,.
\end{aligned}
\end{equation*}
The latter estimate implies \eqref{eq:supergard}.
\end{proof}

We have the following.
\begin{prop}\label{Prop0-flow.calore}
%Assume  {$ \| r \|_{0, H_{x}^{s_1 + 1}, 0} \leq 1   $}. Then, 
For any $ \vphi \in \T^\nu $, for all $s \geq 0$   
the flow $ \cU_N(\vphi, \tau) $ of \eqref{calore.striscia.N} satisfies
\begin{align}
%& {\rm sup}_{\tau \in [0, 1] } \| \cU_N(\vphi, \tau) (u_0) \|_{H^s_x} \leq C \| u_0\|_{H^s_x}\,, 
%\qquad \qquad  \forall \ 0 \leq s \leq 1 \label{stima flusso PDE s 0 1 N.calore}
%\\   
& \label{stima tame Phi t N.calore}
{\rm sup}_{\tau \in [0, 1] } \| \cU_N(\vphi, \tau) (u_0) \|_{H^s_x} \leq  C(s)  \big( \| u_0 \|_{H^s_x}  
+ \| \eta \|_{s + \sigma}^{k_0, \gamma}
%\| r \|_{1,H^{s + s_{1}+1}_x,0}  
\| u_0\|_{L^{2}_x  }\big) \, , \ \  \forall  \ s \geq0 \, ,
\end{align}
uniformly for all $ N \in \N $.
The flow of  \eqref{calore striscia} is a linear bounded  operator  $\cU(\vphi, \tau) : H^s_x (\T^{d}_{\Gamma}) \to H^s_x (\T^{d}_{\Gamma})$ satisfying 
\begin{equation}\label{stima tame Phi t.calore}
{\rm sup}_{\tau \in [0, 1] } \|\cU(\vphi, \tau) (u_0) \|_{H^s_x} \leq  C(s)  \big( \| u_0 \|_{H^s_x}  
+ \| \eta \|_{s + \sigma}^{k_0, \gamma}
%\| r \|_{1,H^{s + s_{1}+1}_x,0}  
\| u_0\|_{L^{2}_x  }\big) \, , \quad \forall \ s\ge0\,.
\end{equation}
%\begin{align}\label{stima flusso PDE s 0 1.calore}
%& {\rm sup}_{\tau \in [0, 1] } \|\cU(\vphi, \tau) (u_0) \|_{H^s_x} \leq C \| u_0\|_{H^s_x}\,, 
%\qquad \qquad  \forall \ 0 \leq s \leq 1 \\   
%& \label{stima tame Phi t.calore}
%{\rm sup}_{\tau \in [0, 1] } \|\cU(\vphi, \tau) (u_0) \|_{H^s_x} \leq  C(s)  \big( \| u_0 \|_{H^s_x}  
%+ \| r \|_{1,H^{s + s_{1}+1}_x,0}  \| u_0\|_{H^{1}_x  }\big) \, , \quad \forall \ s\ge1\,.
%\end{align}
\end{prop}

\begin{proof}
%Before starting the proof we first note the following.  
 \noindent
{\sc Proof of %\eqref{stima flusso PDE s 0 1 N.calore}, 
\eqref{stima tame Phi t N.calore}}.
\\[1mm]
For any $N \in \N$, the equation \eqref{calore.striscia.N} is an ODE on the finite dimensional space $E_N$
which admits a unique solution $u_N(\tau) =  u_N(\lambda, \tau, \vphi, \cdot ) = \cU_N(u_0)(\tau) \in E_N $.
We now provides a priori estimates on the $H_{x}^{s}$-norm of the solution 
$ u_N(t) $ for any $s\geq0$.

\noindent
Recalling that $ \| u_N \|_{H^s_x}^2 = \| \langle D \rangle^s u_N \|_{L^2_x}^2 $ 
and using that $\Pi_{N}$ is self-adjoint, we deduce
\begin{align*}
\partial_\tau \| \langle D \rangle^s u_N \|_{L^2_x}^2 & \stackrel{\eqref{calore.striscia.N}}{=}
-  \big( \langle D \rangle^s \Pi_N \op( a) u_N , \langle D \rangle^s u_N \big)_{L^2_x}  
-  \big(\langle D \rangle^s  u_N ,\langle D \rangle^s \Pi_N \op( a) u_N \big)_{L^2_x} 
\\ 
& =-2{\rm Re}(\langle D\rangle^{s}\op(a)u_N,\langle D\rangle^{s}u_N)_{L^2_x}\,.
\end{align*}
By using estimate \eqref{eq:supergard} in Lemma \ref{lemm:supergard}
one gets
\begin{equation}\label{uragano1}
\begin{aligned}
\partial_\tau \| \langle D \rangle^s u_N \|_{L^2_x}^2
&\stackrel{\eqref{eq:supergard}}{\leq}
-\mathtt{h}\|\langle D\rangle^{s}u_{N}\|_{\dot{H}_{x}^{\frac{1}{2}}}^2
+2C(s)\big(
\|u_{N}\|_{H_{x}^{s}}^2
+(\| \eta \|_{s + \sigma}^{k_0, \gamma})^2
%\| r\|_{1, H^{s+s_{1}+1}_x, 0}^2
\|u_{N}\|_{L_{x}^2}^2
\big)\,.
\\&
\lesssim_{s}
(\| \eta \|_{s + \sigma}^{k_0, \gamma})^2
%\| r\|_{1, H^{s+s_{1}+1}_x, 0}^2
\|u_{N}\|_{L_{x}^2}^2+\|u_{N}\|_{H_{x}^{s}}^2\,.
\end{aligned}
\end{equation}

The bound above specialized for $s=0$  reduces to (using again 
%\eqref{piccolezzarpiccolo} and 
the smallness \eqref{smalleta}) 
$ \partial_\tau \| u_N \|_{L^2_x}^2  \leq C \|  u_N \|_{L^2_x}^2 $, which implies, by Gronwall 
inequality  ,
$ \| \cU_N(\vphi, \tau) (u_0)  \|_{L^2_x  } \leq C' \| u_0 \|_{L^2_x }$, $  \forall \tau \in [0, 1] $.  

For $ s > 0 $,  \eqref{uragano1} reduces to 
$ \partial_\tau \| u_N \|_{H^s_x}^2 \leq 
C(s) \big( (\| \eta \|_{s + \sigma}^{k_0, \gamma})^2
%\| r \|_{H^{s + s_{1}+1}_x}^2  
\| u_0 \|_{L^{2}_x }^2 + \|  u_N \|_{H^s_x}^2 \big) $ and 
the estimate \eqref{stima tame Phi t N.calore} follows by the  Gronwall inequality in differential form.

\noindent
 {\sc Proof of %\eqref{stima flusso PDE s 0 1.calore}, 
 \eqref{stima tame Phi t.calore}.} 
 Now we  pass to the limit $ N \to + \infty $. 
By \eqref{stima tame Phi t N.calore}
%\eqref{stima flusso PDE s 0 1 N.calore} 
the sequence of functions $u_N( \tau, \cdot )$ is bounded in $L^\infty_{ \tau} H^s_x$ and,  
 up to subsequences,  
 \begin{equation}\label{disuguaglianza liminf u uN.calore}
 u_N \stackrel{w^*}{\rightharpoonup} u \ \  \text{in } \  L^\infty_{ \tau} H^s_x \, , 
\qquad  \| u \|_{L^\infty_{ t} H^s_x} \leq \liminf_{N \to + \infty} \| u_N \|_{L^\infty_{ \tau} H^s_x}\,.
 \end{equation} 
  {\sc Claim}: {\it the sequence $ u_N \to u $ in ${\mathcal C}^0_{ \tau} H^{s}_x  \cap {\mathcal C}^1_{ \tau} H^{ s - 1}_x  $, 
 and $ u (\tau, x)$ solves the equation \eqref{calore striscia}}. 
  \\[1mm]
We first prove that $u_N$ is a Cauchy sequence in ${\mathcal C}^0_{ \tau}L^2_x$. Indeed, 
by \eqref{calore.striscia.N},  
    the difference $h_N := u_{N + 1} - u_N$ solves 
        \begin{align*}
    \partial_\tau h_N 
%    &= \pa_{\tau} u_{N + 1} - \pa_{\tau}u_N= 
%    - \Pi_{N+1}{\rm Op}\big(a(\vphi, x, \xi) \big)  u_{N + 1}
%     + \Pi_{N}{\rm Op}\big(a(\vphi, x, \xi) \big) u_N
%    \\
%    &
%=  - \Pi_{N+1}{\rm Op}\big(a(\vphi, x, \xi) \big)  u_{N + 1}
%     + \Pi_{N}{\rm Op}\big(a(\vphi, x, \xi) \big) u_N 
%     + \Pi_{N+1}{\rm Op}\big(a(\vphi, x, \xi) \big) u_N 
%     - \Pi_{N+1}{\rm Op}\big(a(\vphi, x, \xi) \big) u_N
%     \\
     &= - \Pi_{N+1}{\rm Op}\big(a \big) h_N  - (\Pi_{N + 1} - \Pi_N)\op\big(a \big) u_N 
     \\
      h_N(0) &= (\Pi_{N + 1} - \Pi_N) u_0 \, , 
   \end{align*} 
    and therefore, using also that $\Pi_{N + 1}$ is self-adjoint with respect to the $L^2$ scalar product,
    one gets
    \begin{align}
    \partial_t \| h_N(t)\|_{L^2_x}^2 & = (\partial_t h_N\,,\,h_N )_{L^2_x} + ( h_N\,,\,\partial_t h_N )_{L^2_x}  \nonumber
    \\
    & = (-{\rm Op}\big(a \big) h_N , h_N )_{L^2_x} +  ( h_N,- {\rm Op}\big(a \big) h_N )_{L^2_x} 
    \label{coccodrillo -1}
   \\ &
   \quad - ( (\Pi_{N + 1} - \Pi_N)\op\big(a \big) u_N , h_N )_{L^2_x} 
%	\nonumber\\
  %  & \quad 
  -  ( h_N, (\Pi_{N + 1} - \Pi_N)\op\big(a \big) u_N 
 )_{L^2_x} \,. \label{coccodrillo 0}
    \end{align}
    By using Lemma \ref{lemm:supergard} we deduce
    \begin{equation}\label{coccodrillo 22}
    \begin{aligned}
    \eqref{coccodrillo -1}&=-2{\rm Re}(\op(a)h_N,h_{N})_{L^2_x}
    \\&
    \stackrel{\eqref{eq:supergard}}{\leq}
    -\mathtt{h}\|h_{N}\|_{\dot{H}_x^{\frac{1}{2}}}^2+
    C_s\Big(
\| h_{N} \|^2_{L^2_x}  
+ (\| \eta \|_{s_0 + \sigma}^{k_0, \gamma})^2
%\| r \|_{1,H^{ s_{1}+1}_x,0}  
\| h_{N}\|^2_{L^2_x  }\Big)
\,\lesssim \|h_{N}\|_{L^2_x}^2
    \end{aligned}
    \end{equation}
    where, in the last inequality, we also used the smallness \eqref{smalleta}.
%Since  $\Pi_{N + 1}$ is self-adjoint with respect to the $L^2$ scalar product  
%  \begin{align}
%  &(- \Pi_{N+1}{\rm Op}\big(a(\vphi, x, \xi) \big) h_N , h_N )_{L^2_x} +  ( h_N,- \Pi_{N+1}{\rm Op}\big(a(\vphi, x, \xi) \big) h_N )_{L^2_x}   \nonumber\\ 
%&=(-{\rm Op}\big(a(\vphi, x, \xi) \big) h_N , h_N )_{L^2_x} +  ( h_N,-{\rm Op}\big(a(\vphi, x, \xi) \big) h_N )_{L^2_x} \le 0 \label{coccodrillo1}
%  \end{align}
  Moreover, by Cauchy-Schwarz inequality and Lemma \ref{lemma:smoothing}, we have
  \begin{equation}\label{coccodrillo 2}
  \begin{aligned}
|  \eqref{coccodrillo 0}|&\leq
2 \| (\Pi_{N + 1} - \Pi_N) \op(a) u_N  \|_{L^2_x} \| h_N\|_{L^2_x} 
\\&  \leq 
  \| h_N\|_{L^2_x}^2 + \| (\Pi_{N + 1} - \Pi_N) \op(a) u_N  \|_{L^2_x}^2  
  \\
  & \le  \| h_N\|_{L^2_x}^2 + \big( N^{-s_{1}} \| \op(a)  u_N  \|_{H^{s_{1}}_x} \big)^2 
  \\
&  \lesssim_{s_{1}}  \| h_N\|_{L^2_x}^2  
+  N^{-2s_{1}} \| u_0  \|_{H^{s_{1}+1}_x}^2  
  \end{aligned}
  \end{equation}
%  \begin{align}
%  &( (\Pi_{N + 1} - \Pi_N) \op(a) u_N, h_N )_{L^2_x} + ( h_N , (\Pi_{N + 1} - \Pi_N) \op(a) u_N )_{L^2_x} \nonumber\\
%  & \leq 2 \| (\Pi_{N + 1} - \Pi_N) \op(a) u_N  \|_{L^2_x} \| h_N\|_{L^2_x}  \nonumber \\
%&  \leq 
%  \| h_N\|_{L^2_x}^2 + \| (\Pi_{N + 1} - \Pi_N) \op(a) u_N  \|_{L^2_x}^2  \nonumber\\
%  & \le  \| h_N\|_{L^2_x}^2 + \big( N^{-s_{1}} \| \op(a)  u_N  \|_{H^{s_{1}}_x} \big)^2  \nonumber\\
%&  \lesssim_{s_{1}}  \| h_N\|_{L^2_x}^2  + \big( N^{-s_{1}} \| u_0  \|_{H^{s_{1}+1}_x} \big)^2  \label{coccodrillo 2}
%  \end{align}
where we used
Lemma \ref{lemma: action Sobolev} and  that $\| \eta \|_{s_0 + \sigma}^{k_0, \gamma}\leq 1$.
%$  \| r \|_{1, H^{s_{1}}_x, 0} \leq 1 $.  
Hence \eqref{coccodrillo 22}-\eqref{coccodrillo 2} 
imply that 
  \[
  \partial_\tau \| h_N(\tau)\|_{L^2_x}^2 
  \lesssim_{s_{1}} 
  \| h_N(\tau)\|_{L^2_x}^2 + N^{-2s_{1}} \| u_0  \|_{H^{s_{1}+1}_x}^2 \,,
  \]
 and, since $ \| h_N(0)\|_{L^2_x} \leq N^{- s_{1}} \| u_0\|_{H^{s_{1}}_x}$,  by Gronwall lemma we deduce that
$$
\| u_{N + 1} - u_N\|_{{\mathcal C}^0_{ \tau} L^2_x} = \sup_{\tau\in [0, 1]} \| u_{N + 1}( \tau, \cdot) - u_{N}(\tau,\cdot)\|_{L^2_x} 
\lesssim_{s_{1}} N^{-s_{1}} \| u_0\|_{H^{s_{1}+1}_x}\,.
 $$
 By fixing for instance $s_1=1$, by a telescopic argument, one gets that 
 %The above inequality implies that 
 $ u_N $ is a Cauchy sequence in 
 $ {C}^0_{ \tau} L^2_x $. 
Hence $ u_N \to \tilde u \in { C}^0_{ \tau} L^2_x $. By \eqref{disuguaglianza liminf u uN.calore} 
we have $ u = \tilde u  \in {C}^0_{ \tau} L^2_x \cap L^\infty_t H^s_x $. 
 Next, for any $ \bar s \in [0, s) $ we use the interpolation inequality 
  $$
  \| u_N - u \|_{L^\infty_{\tau} H^{\bar s}_x } \leq 
  \| u_N - u \|_{L^\infty_{ \tau} L^2_x }^{1 - \lambda} \| u_N - u \|_{L^\infty_{ \tau} H^{\bar s}_x }^\lambda \, ,
$$
and, since  $ u_N $ is bounded in $ L^\infty_\tau H^s_x $ (see %\eqref{stima flusso PDE s 0 1 N.calore}, 
\eqref{stima tame Phi t N.calore}),
$ u  \in L^\infty_\tau H^s_x $, 
and $ u_N \to u \in { C}^0_{ t} L^2_x $, we deduce that  $u_N \to u$ in each $ L^\infty_\tau H^{\bar s}_x$.
Since $u_N \in { C}^0_\tau H^{\bar s}_x$  are continuous  in $ \tau$, the limit function 
$u \in { C}^0_\tau H^{\bar s}_x $ is continuous as well. 
  Moreover  we also deduce that 
$$
  \partial_\tau u_N =  \Pi_N ( \op(a) u_N ) \to \op(a) u \quad  \text{in}\ \  { C}^0_{ t} H^{\bar s - 1}_x  
  \, , \quad \forall \bar s \in [0, s) \, . 
$$
 As a consequence $ u \in 
 {C}^1_\tau H^{\bar s - 1}_x $ 
 and $ \pa_\tau u =\op(a) u $  
  solves \eqref{calore striscia}. 
 
Finally,  
arguing as in \cite{Taylor}, Proposition 5.1.D, 
it follows that the function $\tau \to \| u(\tau)\|_{H^s_x}^2$ is Lipschitz.  Furthermore,  if $ \tau_n \to \tau $ then 
$ u(\tau_n ) \rightharpoonup u(\tau) $ weakly in $  H^s_x $, because 
 $u(\tau_n) \to u(\tau)$ in $H^{\bar s}_x$ for any $\bar s \in [0, s)$.
 As a consequence the sequence $u(\tau_n) \to u(\tau) $ strongly in $H^s_x$. 
 This proves that  $u \in {C}^0_\tau H^s_x$ and therefore 
 $\partial_\tau u = \op(a) u \in {\mathcal C}^0_\tau H^{s - 1}_x$. 
 \\[1mm]
  {\sc Uniqueness.} If $u_1, u_2 \in {C}^0_\tau H^{s}_x \cap { C}^1_t H^{s - 1}_x $, $s \geq 1 $,  
  are solutions of \eqref{calore striscia}, then $h: = u_1 - u_2$ solves 
  $$
  \partial_\tau h = \op(a) h\,, \qquad h(0) = 0\,.
  $$
  Arguing as in the proof of the estimate   \eqref{stima tame Phi t N.calore}
  %\eqref{diff-ine.calore} 
  we deduce the energy inequality 
  $  \partial_\tau \| h (\tau)\|_{L^2_x}^2 \leq C \| h(\tau)\|_{L^2_x}^2 $. Since $ h(0)= 0 $, 
Gronwall lemma implies that $\| h(\tau)\|_{L^2_x}^2 = 0$, for any $\tau \in [0, 1]$, i.e.  $h(\tau) = 0$. 
Therefore the problem \eqref{calore striscia} has a unique solution $u(t)$ that we denote by  $\cU(\vphi, \tau)(u_0)$. 
  The estimate %\eqref{stima flusso PDE s 0 1.calore}, 
  \eqref{stima tame Phi t.calore} then follows by
  \eqref{disuguaglianza liminf u uN.calore}, % \eqref{stima flusso PDE s 0 1 N.calore}, 
  \eqref{stima tame Phi t N.calore}, since $u_N(t) =\cU_{N}(\vphi, \tau)(u_0)$. 
\end{proof}

The following we collects some results about the flow $\mathcal{U}$ obtained   by following word by word 
 Appendix $A$ of \cite{BM20} 
(see in particular Lemma $A.3$, Proposition $A.5$, $A.7$).
%\red{scrivere che tutte le seguenti sono dimostrate in Berti Montalto Standing}

%In the next Proposition it is stated the smooth dependence of the flow with respect to parameters.

%\begin{lemma}\label{dipendenza liscia dai parametri.calore}
%Let $a(z, \cdot )$ % \in {\mathcal C}^\infty(\T)$ and 
%$ p_0 $-times differentiable, resp. ${ C}^{p_0} $,  with respect to $z \in {\mathcal B}_X $, 
%where ${\mathcal B}_X$ is 
%an open subset of a  Banach space $ X $. 
%Then, for any $ p \leq p_0 $, the flow map $ \cU(z; \tau)$, $ \tau \in [0, 1] $,  
%  is smooth in $ z $, more precisely,  the map 
%$$
%{\mathcal B}_X \ni z \mapsto \cU(z;\tau) 
%\in {\mathcal L}(H^{s}_x, H^{s - p - 1}_x) \, , \quad \forall s  \geq p +1  \, ,
%$$
%is $ p $-times differentiable, resp. $ { C}^p $.
% Moreover, for any $z \in {\mathcal B}_X $,  the derivative $\partial_z^p \Phi(z, t)$ is a 
% multilinear form from $X^p$ in ${\mathcal L}(H_x^s, H^{s - p}_x)$.
%\end{lemma} 
%
%\begin{lemma}\label{lemma finitezza norme stime flusso.calore}
%For any $|\beta| \leq \beta_0, |k| \leq k_0, \tau \in [0, 1], h \in { C}^\infty(\T_*^{\nu + d})$, 
%the function $ \partial_\lambda^k \partial_\vphi^\beta \cU (\vphi, \tau )  h $ is ${ C}^\infty(\T_*^{\nu + d})$. 
%\end{lemma}

\begin{prop} \label{Prop1-flow.calore}
%Let $s_0:=(\nu+d)/2+1$.
Assume that  
\begin{equation*}%\label{ansatz derivate vphi flusso.calore}
 \quad \| \eta \|_{ s_0 + \s_{1}} \leq 1 \, , 
 \quad  \| \eta \|_{s_0 + \s_{2} } \leq \delta (s)  
\end{equation*}
for some $\s_{1}=\s_{1}(s_{0})>0$, $\s_{2}=\s_{2}(s_{0})>0$ and for some $ \delta(s) > 0 $ small. Then the following tame estimates hold: 
\begin{align*}%\label{stima flusso norma bassa.calore}
%& {\rm sup}_{t \in [0, 1]}\| \cU(\tau) u_0 \|_s \leq C(s) \| u_0\|_s\,, 
%\qquad \quad  \forall s \in [0, s_0 + 1 ] \, , 
%\\
& 
%\label{stima flusso norme unite.calore}
{\rm sup}_{t \in [0, 1]} \|  \cU(\tau) u_0 \|_{s} 
\leq 
C(s) \big( \| u_0 \|_{s} +\| \eta \|_{s + {\s_{1}} } \| u_0\|_{s_0} \big)\,, 
\quad \forall s \geq s_0  \,.
\end{align*}
\end{prop}

\begin{prop} \label{Teorema totale partial vphi beta k D beta k Phi.calore}
Let $ \beta_0, k_0 \in \N$. 
For any $\beta, k \in \N^\nu$ with $|\beta| \leq \beta_0$, $|k| \leq k_0$,
for any $m_1, m_2 \in \R$ with $m_1 + m_2 = |\beta| + |k|$,
for any $s \geq s_0$, 
there exist constants $\s_{1}= \sigma_{1}(s_{0}, m_1) > 0$, 
$ \s_{2}=\sigma_{2}(|\beta|, |k|, m_1, m_2) > 0$, $\delta(s, m_1) > 0$ such that if
\begin{equation}\label{piccolezza a partial vphi beta k D beta k}
\| \eta \|_{ s_0 + \s_{1} } \leq \delta (s, m_1) \,, \quad  
\| \eta \|_{ s_0 + \sigma_{2}}^{k_0, \gamma} \leq 1\,, 
\end{equation} 
then the following estimates hold: 
\begin{equation*}%\label{copenaghen B omega.calore}
\sup_{t \in [0, 1]} \| \langle D \rangle^{- m_1} \partial_\lambda^k \partial_\vphi^\beta \cU(\tau) 
\langle D \rangle^{- m_2} h \|_s
%\\&
%\qquad\qquad\qquad
\lesssim_{s,\beta_0, k_0, m_1,m_2} \gamma^{- |k|} 
\Big( \| h \|_s + \| \eta \|_{s + \sigma_{2}}^{k_0, \gamma} \| h \|_{s_0} \Big)\,.
\end{equation*}
\end{prop}

\begin{prop}\label{flussoCk0.calore}
Assume \eqref{piccolezza a partial vphi beta k D beta k}. Then the flow $\cU( \lambda; \tau)$ of \eqref{calore striscia} is $\mathcal{D}^{k_{0}}$-$k_{0}$-tame 
(recall Definition \ref{Dksigmatame}), 
more precisely,
for all $ k \in \N^{\nu + 1} $, $ |k| \leq k_0 $,  $ s \geq s_0 $, it satisfies
%the flow $ \cU(\tau) := \cU(\lambda, \vphi, \tau)$ of \eqref{calore striscia} satisfies
%\shu{stime da rifare. O le fai con la perdita precisa di regolarità su a oppure esiste sigma.}
\begin{align*}
%\label{arido 1.calore}
& \sup_{t \in [0, 1]}\|  \partial_\lambda^{k} \cU(\tau) h \|_s    \lesssim_s  \gamma^{- |k|} 
\big( \| h \|_{s+ |k|} + \| \eta \|_{s + s_0   +  |k| + 1}^{k_0, \gamma} 
\|  h \|_{s_0 +  |k|}  \big)  \, , 
\\
& 
\sup_{t \in [0, 1]} \| \partial_\lambda^k (\cU(\tau)- {\rm I d}) h \|_s   \nonumber
\\&\qquad\qquad
\lesssim_s \gamma^{- |k|}\big( \| \eta \|_{s_0}^{k_0, \gamma} \|  h\|_{s +  |k|+1} + \| \eta \|_{s + s_0   +  |k|  + 1 + 1,0}^{k_0, \gamma} \| h \|_{s_0 + |k|+1}  \big) \, .  
%{arido 2.calore}
\end{align*}
\end{prop}

We consider also the dependence of the flow $\Phi$  with respect to the torus $ i := i(\vphi ) $.  

\begin{prop} 
%\label{lemma:tame derivate flusso ancora diverso derivate.calore}
Let $p > s_0$,  $ \beta_0 \in \N$. 
For any $\beta \in \N^\nu$, $|\beta| \leq \beta_0$, for any $m_1, m_2 \in \R$ satisfying $m_1 + m_2 = |\beta| + 1$ there exists 
a constant $\sigma(|\beta|) = \sigma(|\beta|, m_1, m_2) > 0$ such that if
$ \| \eta \|_{ p + \sigma(\beta_0)} \leq \delta (s_1) $
with $\delta(p) > 0$ small enough, then the following estimate holds:  
\begin{align*}
\sup_{t \in [0, 1]} \| \langle D \rangle^{- m_1}  \partial_\vphi^\beta \Delta_{12} \cU(\tau) \langle D \rangle^{- m_2} h \|_{p} 
\lesssim_{p} \|  \eta_1-\eta_2 \|_{ p+ \sigma(\beta_0)}  \| h \|_{p} \, ,
%\label{stima derivate flusso generalissima derivate i.calore}
\end{align*}
where $\Delta_{12}\cU:= \cU(\eta_2) - \cU(\eta_1)$.
% and $\Delta_{12}\eta:= \eta(i_2) - \eta(i_1)$.
\end{prop}
As a consequence of the discussion above we obtain Lemma \ref{stima derivate lambda vphi flusso parabolico}.

\section{Whitney differentiable functions} \label{sec:U}
In this Appendix we  recall  the notion of Whitney differentiable functions and 
the Whitney extension theorem, following the version of Stein \cite{Stein}. Then we 
 prove the lemmata stated in Section \ref{subsec:function spaces}.
 The following definition is the adaptation of the one in Section 2.3, Chapter VI of \cite{Stein}
to Banach-valued functions.

\begin{defn}{\bf (Whitney differentiable functions)}
\label{def:Lip F}
Let $F$ be a closed subset of $\R^n$, $n \geq 1$. 
Let $Y$ be a Banach space.
Let $k \geq 0$ be an integer, and $k < \rho \leq k+1$. 
We say that a function $f : F \to Y$ belongs to $\Lip(\rho,F,Y)$ if there exist 
functions 
\[
f^{(j)} : F \to Y, \quad j \in \N^n, \ 0 \leq |j| \leq k, 
\]
with $f^{(0)} = f$, and a constant $M > 0$ 
such that if $R_j(x,y)$ is defined by 
\begin{equation*} %\label{16 Stein}
%\resizebox{0.9\hsize}{!}{$
f^{(j)}(x) = \sum_{\ell \in \N^n : |j+\ell| \leq k} \frac{1}{\ell!} \, f^{(j+\ell)}(y) \, (x-y)^\ell
+ R_j(x,y), \quad x,y \in F, 
%$}
\end{equation*} 
then 
\begin{equation} \label{17 Stein}
\| f^{(j)}(x) \|_Y \leq M, \quad 
\| R_j(x,y) \|_Y \leq M |x-y|^{\rho - |j|} \, , \quad 
\forall x,y \in F, \ |j| \leq k \, .
\end{equation} 
An element of $\Lip(\rho,F,Y) $ 
is  in fact the collection $\{ f^{(j)} : |j| \leq k \}$. 
The norm of $ f \in \Lip(\rho,F,Y)$ is defined as the smallest $M$ 
for which the inequality \eqref{17 Stein} holds, namely
\begin{equation} \label{def norm Lip Stein}
\| f \|_{\Lip(\rho,F,Y)} := \inf \{ M > 0 : \text{\eqref{17 Stein} holds} \} \, .
\end{equation}
If $ F = \R^n $ by  $ \Lip(\rho, \R^n,Y) $  we shall mean the linear space of the functions $ f = f^{(0)} $ for which
there exist $ f^{(j)} = \pa_{x}^j f $, $ |j| \leq k $,  satisfying  \eqref{17 Stein}.
\end{defn}

Notice that, if $ F = \R^n $,  
the $ f^{(j)} $, $ | j | \geq 1 $,
 are uniquely determined by $ f^{(0)} $ (which is not the case  for a general $ F $ with for example isolated points).

In the case $F = \R^n$, $\rho = k+1$ and $ Y $ is  a Hilbert space, 
the space $\Lip(k+1,\R^n,Y)$ is isomorphic to the Sobolev space 
$W^{k+1,\infty}(\R^n,Y)$, with equivalent norms
\begin{equation} \label{equiv Sob Lip}
C_1 \| f \|_{W^{k+1,\infty}(\R^n,Y)}
\leq \| f \|_{\Lip(k+1,\R^n,Y)} 
\leq C_2 \| f \|_{W^{k+1,\infty}(\R^n,Y)}
\end{equation}
where $C_1, C_2$ depend only on $k,n$.
For $Y = \C$ this isomorphism is classical, see e.g. \cite{Stein}, 
and it is based on the Rademacher theorem concerning the a.e.\ differentiability of Lipschitz functions, and the fundamental theorem of calculus for the Lebesgue integral. 
Such a property may fail 
for a Banach valued function, 
but it holds for a Hilbert space, see Chapter 5 of \cite{BenyLind}
(more in general it holds if $ Y $ is reflexive or it satisfies the Radon-Nykodim property).

The following key result provides an extension of a Whitney differentiable function $f$ defined on a closed subset $F$ of $\R^n$ to the whole domain $\R^n$, with equivalent norm.

\begin{thm}{\bf (Whitney extension Theorem)}  \label{thm:WET}
Let $F$ be a closed subset of $\R^n$, $n \geq 1$,  
$Y$ a Banach space, 
$k \geq 0$ an integer, and $k < \rho \leq k+1$. 
There exists a linear continuous extension operator $\cE_k : \Lip(\rho,F,Y) \to \Lip(\rho, \R^n,Y)$ 
which gives an extension $\cE_k f \in \Lip(\rho, \R^n,Y)$ to any $f \in \Lip(\rho,F,Y)$. 
The norm of $\cE_k$ has a bound independent of $F$, 
\begin{equation} \label{bound Ek Stein}
\| \cE_k f \|_{\Lip(\rho,\R^n,Y)} \leq C \| f \|_{\Lip(\rho,F,Y)} \, , \quad \forall f \in \Lip(\rho,F,Y) \, ,
\end{equation}
where $C$ depends only on $n,k $ (and not on $ F , Y $). 
\end{thm}

\begin{proof}
This is Theorem 4 in Section 2.3, Chapter VI of \cite{Stein}. The proof in \cite{Stein}
is written for real-valued functions $f : F \to \R$, but it also holds for functions $f : F \to Y$ for any (real or complex) Banach space $Y$, \emph{with no change}. 
The extension operator $\cE_k$ is defined in formula (18) in Section 2.3, Chapter VI of \cite{Stein},
and it is linear by construction.
\end{proof}

Clearly, since $\cE_k f$ is an extension of $f$, one has 
\begin{equation} \label{basso Stein}
\| f \|_{\Lip(\rho,F,Y)} \leq \| \cE_k f \|_{\Lip(\rho,\R^n,Y)} 
\leq C \| f \|_{\Lip(\rho,F,Y)} \, .
\end{equation}
In order to extend a function defined on a closed set $F \subset \R^n$ 
with values in scales of Banach spaces (like $H^s(\T^{\nu+1})$),
we observe that the extension provided by Theorem \ref{thm:WET} 
does not depend on the index of the space (namely $s$).

\begin{lemma}  \label{lemma:2702.1}
Let $F$ be a closed subset of $\R^n$, $n \geq 1$, 
let $k \geq 0$ be an integer, and $k < \rho \leq k+1$.
Let $Y \subseteq Z$ be two Banach spaces. 
Then $\Lip(\rho,F,Y) \subseteq \Lip(\rho,F,Z)$. 
The two extension operators 
$\cE_k^{(Z)} : \Lip(\rho,F,Z) \to \Lip(\rho, \R^n,Z)$ 
and $\cE_k^{(Y)} : \Lip(\rho,F,Y) \to \Lip(\rho, \R^n,Y)$
provided by Theorem \ref{thm:WET} satisfy
\[
\cE_k^{(Z)} f = \cE_k^{(Y)} f \quad \forall f \in \Lip(\rho,F,Y) \, .
\]
As a consequence, we simply denote $\cE_k$ the extension operator.
\end{lemma}

\begin{proof}
The lemma follows directly by the construction of the extension operator $\cE_k$ in formula (18) in Section 2.3, Chapter VI of \cite{Stein}. 
The explicit construction relies on a nontrivial decomposition in cubes of the domain $\R^n$ only.
\end{proof}

\smallskip

Thanks to the equivalence \eqref{basso Stein}, Lemma \ref{lemma:2702.1}, and 
\eqref{equiv Sob Lip} which holds for functions valued in $ H^s $, 
classical interpolation and tame estimates for products, projections, and composition of Sobolev 
functions can be easily extended to Whitney differentiable functions.

The difference between the Whitney-Sobolev norm introduced in Definition \ref{def:Lip F uniform}
and the norm in Definition \ref{def:Lip F} (for $ \rho = k+1 $, $ n = \nu + 1 $, 
and target space $ Y = H^s ( \T^{\nu+d}_*,\C) $)
is the weight $ \g \in (0,1]$. 
Observe that the introduction of this weight simply amounts to the following rescaling $\cR_\g$: 
given $u = (u^{(j)})_{|j| \leq k}$, we define $\cR_\g u = U = (U^{(j)})_{|j| \leq k}$ as
\begin{equation} \label{resca}
\lambda = \g \mu, \qquad 
\g^{|j|} u^{(j)}(\lambda) 
= \g^{|j|} u^{(j)}(\g \mu) 
=: U^{(j)}(\mu) = U^{(j)}(\g^{-1} \lambda), \qquad 
U := \cR_\g u \, .
\end{equation}
Thus $u \in \Lip(k+1,F,s,\g)$ if and only if $U \in \Lip(k+1, \g^{-1} F ,s,1)$, with 
\begin{equation}\label{rescaling:gamma-1}
\| u \|_{s,F}^{k+1,\g} = \| U  \|_{s,\g^{-1} F}^{k+1,1} \,.
\end{equation}
Under the rescaling $\cR_\g$, \eqref{equiv Sob Lip} gives the equivalence of the two norms 
\begin{equation} \label{0203.1}
\| f \|_{W^{k+1,\infty, \gamma}(\R^{\nu+1}, H^s)} := 
\sum_{|\alpha| \leq k+1} \g^{|\alpha|} \| \pa_\lambda^\alpha f \|_{L^\infty(\R^{\nu+1}, H^s)} 
\sim_{\nu,k} \| f \|_{s, \R^{\nu+1}}^\kug \,.
\end{equation}
Moreover, given $u \in \Lip(k+1,F,s,\g)$, its extension 
\begin{equation} \label{Wg}
\tilde u := \cR_\g^{-1} \cE_k \cR_\g u \in \Lip(k+1, \R^{\nu+1}, s, \g) 
\quad \text{satisfies} \quad 
\| u \|^{k+1,\g}_{s,F} 
\sim_{\nu,k} \| \tilde u \|_{s, \R^{\nu+1}}^\kug \,.
\end{equation}

\begin{proof}[{\bf Proof of Lemma \ref{lemma:smoothing}}]
Inequalities \eqref{p2-proi}-\eqref{p3-proi} follow by 
\begin{align*}
 (\Pi_N u)^{(j)}(\lambda) = \Pi_N [u^{(j)}(\lambda)]\,, 
 \quad 
 R_j^{(\Pi_N u)}(\lambda, \lambda_0) = \Pi_N [ R_j^{(u)}(\lambda, \lambda_0) ]\,, 
\end{align*}
for all $0 \leq |j|  \leq k$, $\lambda, \lambda_0 \in F$, 
and the usual smoothing estimates
$ \| \Pi_N f \|_{s}  \leq N^\alpha \| f \|_{s-\alpha} $ and 
$ \| \Pi_N^\perp f \|_{s}  
\leq N^{-\alpha} \| f \|_{s + \alpha} $ for Sobolev functions.
\end{proof}

%\smallskip
\begin{proof}[{\bf Proof of Lemma \ref{lemma:interpolation}}]
Inequality \eqref{2202.3} follows from the classical interpolation inequality 
$\| u \|_s \leq \| u \|_{s_0}^\theta \| u \|_{s_1}^{1-\theta}$, 
$s = \theta s_0 + (1-\theta) s_1$ for Sobolev functions,
and from the Definition \ref{def:Lip F uniform} of 
Whitney-Sobolev norms, since 
\begin{align*}
\g^{|j|} \| u^{(j)}(\lambda) \|_s 
& \leq 
(\g^{|j|} \| u^{(j)}(\lambda) \|_{s_0})^\theta
(\g^{|j|} \| u^{(j)}(\lambda) \|_{s_1})^{1-\theta}
\leq 
(\| u \|_{s_0,F}^{k+1,\g})^\theta (\| u \|_{s_1,F}^{k+1,\g})^{1-\theta},
\\
\g^{k+1} \| R_j(\lambda,\lambda_0) \|_s 
& \leq 
(\g^{k+1} \| R_j(\lambda,\lambda_0) \|_{s_0})^\theta
(\g^{k+1} \| R_j(\lambda,\lambda_0) \|_{s_1})^{1-\theta}
\\&
\leq (\| u \|_{s_0,F}^{k+1,\g})^\theta (\| u \|_{s_1,F}^{k+1,\g})^{1-\theta}
|\lambda - \lambda_0|^{k+1-|j|}\,.
\end{align*}
Inequality \eqref{2202.2} follows from \eqref{2202.3} 
by using the asymmetric Young inequality
(like in the proof of Lemma 2.2 in \cite{BM20}).
\end{proof}

%\smallskip
\begin{proof}[{\bf Proof of Lemma \ref{lemma:LS norms}}]
By \eqref{0203.1}-\eqref{Wg}, the lemma follows from the corresponding inequalities 
for functions in $W^{k+1,\infty,\g}(\R^{\nu+1}, H^s)$, which are proved, for instance, 
in \cite{BM20}. % (formula (2.72), Lemma 2.30).
\end{proof}

\smallskip
For any $\rho > 0$, we define the  $ {\cC}^\infty $ function $ h_\rho : \R \to \R $,  
\begin{equation}\label{def-gl}
h_\rho (y) := \frac{\chi_\rho(y)}{y} = \frac{\chi(y \rho^{-1})}{y}  \, , 
\quad \forall y \in \R \setminus \{ 0 \}, 
\quad h_\rho(0) := 0  \, , 
\end{equation}
where $\chi$ is the cut-off function introduced in \eqref{cut off simboli 1}, 
and $\chi_\rho(y) := \chi (y / \rho)$.
Notice that the function $ h_\rho $ 
is of class $ \cC^\infty $ because  $ h_\rho (y) = 0 $ for $ |y| \leq \rho / 3 $.
Moreover 
by the properties of $ \chi $ in  \eqref{cut off simboli 1}
we have 
\begin{equation} \label{fatti d'acca rho}
h_\rho (y) = \frac{1}{ y} \, , \ \forall  |y| \geq \frac{2 \rho}{ 3} \, , \qquad 
| h_\rho(y) | \leq \frac{3}{\rho} \, , \ \forall y \in \R \, . 
\end{equation}
To prove Lemma \ref{lemma:WD}, we use the following preliminary lemma. 

\begin{lemma} 
\label{lemma:cut-off sd}
Let $\mathtt \Lambda_0 \subseteq \R^{\nu + 1}$, $f :\mathtt \Lambda_0 \to \R$ and  $\rho > 0$. 
Then the function 
\begin{equation} \label{2017.1805.1}
g(\lambda) := h_\rho (f(\lambda)), \quad \forall \lambda \in \Lambda_0 \, , 
\end{equation}
where $h_\rho$ is defined in \eqref{def-gl},
coincides with $ 1/ f (\lambda) $ 
on the set $ F := \{ \lambda \in\mathtt \Lambda_0  : |f(\lambda)| \geq \rho \}$. 

If the function $f $ is in $ W^{k+1,\infty}(\mathtt \Lambda_0 ,\R)$, with estimates
\begin{equation} \label{0103.1}
 | \pa_\lambda^\alpha f(\lambda) | \leq M \, ,  \quad  \forall \alpha \in \N^{\nu+1} \, , 
\ \ 
1 \leq |\alpha| \leq k+1 \, , 
\end{equation} 
for some $M \geq \rho$, 
then the function $g $ is in $ W^{k+1,\infty}( \mathtt \Lambda_0 ,\R)$ and 
\begin{equation} \label{0103.2}
 | \pa_\lambda^\alpha g(\lambda) | 
\leq C_\alpha \frac{M^{|\alpha|}}{\rho^{|\alpha|+1}}\, , 
\quad \forall \alpha \in \N^{\nu+1}, \ \ 
0 \leq |\alpha| \leq k+1\,.
\end{equation} 
\end{lemma}

\begin{proof}
By \eqref{fatti d'acca rho}, $g(\lambda) = 1 / f(\lambda)$ for all $\lambda \in F$.
The derivatives of $ h_\rho (y) $ are
\begin{align*}
h_\rho^{(m)} (y) = \pa_y^m \big( \frac{\chi_\rho(y)}{y} \big) 
& = \sum_{m_1 + m_2 = m } C_{m_1, m_2} 
\rho^{-m_1} (\pa_y^{m_1} \chi) (\rho^{-1} y) y^{-m_2 - 1}  \, , 
\quad m \geq 0,
\end{align*}
that we may bound  --- we have just to consider $ |y| \geq \rho / 3 $ 
(otherwise $ h_{\rho} (y) = 0 $) --- as 
\begin{equation}\label{pazh}
|h_\rho^{(m)} (y) | \lesssim_m \sum_{m_1 + m_2 = m} 
\rho^{-m_1} |y|^{-m_2 - 1} \lesssim_m \rho^{-m-1} \, , 
\quad \forall y \in \R \, . 
\end{equation}
Using the Fa\`a di Bruno formula, for $|\alpha| \geq 1$
we compute the derivatives of the composite function 
\[
\pa_\lambda^\alpha g(\lambda) 
= 
\pa_\lambda^\alpha  \, h_\rho ( f(\lambda)) 
=
\sum_{1\leq q \leq |\alpha|} 
\sum_{\substack{\s_1+ \ldots + \s_q = \alpha, \\ \s_j \neq 0, j=1, \ldots, q}} 
C_{q, \s_1, \ldots, \s_q} h_\rho^{(q)} ( f(\lambda)) 
\pa_\lambda^{\s_1} f(\lambda) \ldots \pa_\lambda^{\s_q} f(\lambda) \, , 
\]
and, using \eqref{pazh}, we get, for all $ |\alpha| \geq 1$,
\begin{align*}
|\pa_\lambda^\alpha g(\lambda)| & \leq C_{\alpha}
\sum_{1\leq q \leq |\alpha|} \sum_{\substack{\s_1+ \ldots + \s_q = \alpha, \\
 \s_j \neq 0, j=1, \ldots, q}} 
% C_{q, \gamma_1, \ldots, \gamma_q} 
\rho^{-q-1} 
|\pa_\lambda^{\s_1} f(\lambda)| \ldots |\pa_\lambda^{\s_q}  f(\lambda)| \\
& \stackrel{ \eqref{0103.1} } \leq 
C_\alpha \sum_{1\leq q \leq |\alpha|}  \rho^{-q-1} 
M^q \leq 
C_\alpha  \rho^{-|\alpha|-1} 
 M^{|\alpha|} \, .
\end{align*}
Formula \eqref{0103.2} for $\alpha = 0$ holds by \eqref{fatti d'acca rho}.
\end{proof}

To prove Lemma \ref{lemma:WD2}, we use the following preliminary lemma. 
\begin{lemma}\label{lemma lip estensione lip}
Let $X$ be a Banach space with norm $\| \cdot \|$ and let $f : B_R \to \C$ be a Lipschitz function where $B_R := \{ u \in X : \| u \| \leq R \}$. Let $\rho > 0$ and 
$$
\| f \|^{\rm lip} := \sup_{\begin{subarray}{c}
u_1, u_2 \in B_R \\
u_1 \neq u_2
\end{subarray}} \dfrac{|f(u_1) - f(u_2)|}{\| u_1 - u_2\|} < \infty\,. 
$$
Then the function 
\begin{equation}\label{g h rho per lip}
g(u) := h_\rho(f(u)), \quad u \in B_R
\end{equation}
(see \eqref{def-gl}) is a Lipschitz extension of the function 
$$
\frac{1}{f} : E_\rho := \{ u \in B_R : |f(u)| \geq \rho \} \to \C
$$
saisfying the Lipschitz estimate
$$
|g(u_1) - g(u_2)| \lesssim  \frac{\| f \|^{\rm lip}}{\rho^2} \| u_1 - u_2 \|, \quad \forall u_1, u_2 \in B_R\,. 
$$
\end{lemma}
\begin{proof}
Clearly $g = \frac{1}{f}$ on the set $E_\rho$, hence $g$ is an extension of $\frac{1}{f}$. Let $u_1, u_2 \in B_R$ with $\| u_1 - u_2 \| \geq \delta \rho$ (where $0 < \delta \ll 1$ is an absolute constant to be determines). One has that (recall \eqref{cut off simboli 1}, \eqref{g h rho per lip}, \eqref{def-gl})
$$
\begin{aligned}
|g(u_1) - g(u_2)| & \leq |g(u_1)| + |g(u_2)| \leq \Big| \dfrac{\chi\Big( f(u_1) \rho^{- 1}\Big)}{f(u_1)}\Big| + \Big| \dfrac{\chi\Big( f(u_2) \rho^{- 1}\Big)}{f(u_2)}\Big| \\
& \lesssim \rho^{- 1} \lesssim \rho^{- 2} \delta^{- 1} \rho \delta \lesssim \rho^{- 2} \delta^{- 1} \| u_1 - u_2 \|\,.
\end{aligned}
$$
Let us consider now $u_1, u_2 \in B_R$ such that $\| u_1 - u_2 \| \leq \delta \rho$. If $| f(u_1) | \leq \frac{\rho}{4}$, then 
$$
|f(u_2)| \leq |f(u_1)| + |f(u_2) - f(u_1)| \leq \frac{\rho}{4} + \| f \|^{\rm lip} \| u_1 - u_2 \| \leq \frac{\rho}{4} + \| f \|^{\rm lip} \delta \rho < \frac{\rho}{3}
$$
by taking $0 < \delta = \frac{1}{24 \| f \|^{\rm lip}} < \frac{1}{12 \| f \|^{\rm lip}}$. Hence 
$$
|f(u_1)|, |f(u_2)| < \frac{\rho}{3} \Longrightarrow \chi\Big( f(u_1) \rho^{- 1}\Big), \chi\Big( f(u_2) \rho^{- 1}\Big)  = 0 \Longrightarrow g(u_1 ), g(u_2 ) = 0
$$
and hence 
$$
|g(u_1) - g(u_2)| \lesssim \rho^{- 2} \delta^{- 1} \| u_1 - u_2 \| \,.
$$
If $|f(u_1)| \geq \frac{\rho}{3}$, then 
$$
\begin{aligned}
|f(u_2)| & \geq |f(u_1)| - |f(u_2) - f(u_1)| \geq \frac{\rho}{3} - \| f \|^{\rm lip} \| u_1 -u_2 \|  \\
& \geq \frac{\rho}{3} - \| f \|^{\rm lip} \delta \rho \geq \frac{\rho}{4}
\end{aligned}
$$
by taking $0 < \delta = \frac{1}{24 \| f \|^{\rm lip}} < \frac{1}{12 \| f \|^{\rm lip}}$. Thus 
$$
\begin{aligned}
|g(u_1) - g(u_2)| & \leq \dfrac{\Big| \chi\Big( f(u_1) \rho^{- 1}\Big) \Big| |f(u_1) - f(u_2)|}{|f(u_1)| |f(u_2)|} + \dfrac{\Big| \chi\Big( f(u_1) \rho^{- 1}\Big) - \chi\Big( f(u_2) \rho^{- 1}\Big) \Big| |f(u_1)| }{|f(u_1)| |f(u_2)|} \\
& \lesssim \rho^{- 2} |f(u_1) - f(u_2)| + \rho^{- 1} \dfrac{|f(u_1) - f(u_2)|}{|f(u_2)|} \\
& \lesssim \rho^{- 2} |f(u_1) - f(u_2)| \lesssim \rho^{- 2} \| f \|^{\rm lip} \| u_1 - u_2 \|\,. 
\end{aligned}
$$
The claimed statement has then been proved. 
\end{proof}

\begin{proof}[{\bf Proof of Lemma \ref{lemma:WD}}]
The function $(\ompaph)^{-1}_{ext} u$ defined in \eqref{def ompaph-1 ext} is
\[
\big( (\ompaph)^{-1}_{ext} u \big) (\lambda,\vphi,x) 
= - \ii \sum_{(\ell,j) \in \Z^{\nu+1}} 
g_\ell(\lambda) u_{\ell, j}(\lambda) \, e^{\ii (\ell \cdot \vphi + j x )} \, ,
\]
where $g_\ell(\lambda) = h_\rho( \om \cdot \ell )$ in \eqref{2017.1805.1}
with $\rho = \g \langle \ell \rangle^{-\t}$ and $f(\lambda) = \om \cdot \ell$.
The function $f(\lambda)$ satisfies \eqref{0103.1} 
with $M = |\ell|$. 
Hence $g_\ell(\lambda)$ satisfies \eqref{0103.2}, namely
\begin{equation} \label{0103.6}
 | \pa_\lambda^\alpha g_\ell(\lambda) | \leq C_k \g^{-1 - |\alpha|} \langle \ell \rangle^\mu 
\quad \forall \alpha \in \N^{\nu+1}, \ 0 \leq |\alpha| \leq k+1,
\end{equation}
where $\mu = k + 1 + (k+2) \t$ is defined in \eqref{Diophantine-1}.
One has 
\[
\pa_\lambda^\alpha ( g_\ell(\lambda) u_{\ell, j}(\lambda) ) 
= \sum_{\alpha_1 + \alpha_2 = \alpha} C_{\alpha_1, \alpha_2} (\pa_\lambda^{\alpha_1} g_\ell)(\lambda) 
(\pa_\lambda^{\alpha_2} u_{\ell, j})(\lambda), 
\]
whence, by \eqref{0103.6}, we deduce 
\[
\g^{|\alpha|} \| \pa_\lambda^{\alpha}((\ompaph)^{-1}_{ext} u) (\lambda) \|_s  
\leq C_k \g^{-1} \| u \|_{s+\mu}^{k+1,\g}
\]
and therefore \eqref{2802.2}.
The proof is concluded by observing that the restriction of 
$(\ompaph)^{-1}_{ext} u$ 
to $F$ gives $(\ompaph)^{-1} u$ as defined in \eqref{def:ompaph}, 
and \eqref{Diophantine-1} follows by \eqref{Wg}.
\end{proof}

\begin{proof}[{\bf Proof of Lemma \ref{lemma:WD2}}]
Let us check estimate \eqref{2802.2nuova}. 
The function $((\omega-\mathtt{V}\mathtt{m})\cdot\pa_{\vphi})^{-1}_{ext} u$ defined in \eqref{def ompaph-1 ext2} is
\[
\big( ((\omega-\mathtt{V}\mathtt{m})\cdot\pa_{\vphi})^{-1}_{ext}  u \big) (\lambda,\vphi) 
= - \ii \sum_{\ell\in \Z^{\nu}} 
g_\ell(\lambda) u_{\ell}(\lambda) \, e^{\ii \ell \cdot \vphi } \, ,
\]
where $g_\ell(\lambda) = h_\rho( (\omega-\mathtt{V}\mathtt{m}) \cdot \ell )$ 
in \eqref{2017.1805.1}
with $\rho = \g \langle \ell \rangle^{-\t}$ and $f(\lambda) = (\omega-\mathtt{V}\mathtt{m}) \cdot \ell$.
The function $f(\lambda)$ satisfies \eqref{0103.1} 
with $M = |\ell|$ recalling that, by assumption $|\mathtt{m}|^{k_0,\gamma}\lesssim\e$
and that $\e\gamma^{-k-1}\ll1$. 
Hence $g_\ell(\lambda)$ satisfies \eqref{0103.2}. Thereofere one concludes the proof follows 
Lemma \ref{lemma:WD}.

Let us now check \eqref{2802.2nuova2}.
Recalling definition \eqref{def ompaph-1 ext2} we have that
\[
 ((\omega-\mathtt{V}\mathtt{m}(i_1))\cdot \pa_\vphi )^{-1}_{ext} u
-
 ((\omega-\mathtt{V}\mathtt{m}(i_2))\cdot \pa_\vphi )^{-1}_{ext} u
 = - \ii \sum_{\ell\in \Z^{\nu}\setminus\{0\}} 
(g_\ell(i_1)-g_{\ell}(i_2) )u_{\ell}(\lambda) \, e^{\ii \ell \cdot \vphi }
\]
where
\[
g_\ell(i)=\frac{\chi( (\omega-\mathtt{V}\mathtt{m}(i))\cdot \ell \g^{-1} \langle \ell \rangle^{\t})}{(\omega-\mathtt{V}\mathtt{m}(i))\cdot\ell}=h_{\rho}(f(i))
\]
where $\rho=\gamma\langle \ell\rangle^{-\tau}$, $h_{\rho}$ in \eqref{def-gl}
and $f(i)=f_{\ell}(i)=(\omega-\mathtt{V}\mathtt{m})\cdot\ell$.
By the estimate \eqref{stimaemminodelta12} we deduce that 
\[
|f_{\ell}(i_1)-f_{\ell}(i_2)|\lesssim |\ell| \e\|i_1-i_2\|_{s_0+\s}\,.
\]
Then by  applying Lemma \ref{lemma lip estensione lip} we get
\[
|g_\ell(i_1)-g_\ell(i_2)|\lesssim\e \gamma^{-2} \langle\ell\rangle^{2\tau+1}\|i_1-i_2\|_{s_0+\s}\,,
\] 
and therefore \eqref{2802.2nuova2} follows.
\end{proof}

\begin{proof}[{\bf Proof of Lemma \ref{Moser norme pesate}}]
Given $u\in\Lip(k+1,F,s,\g)$, we consider its extension 
$ \tilde u$ belonging to the space $\Lip(k+1,\R^{\nu+1},s,\g)$ provided by \eqref{Wg}. 
Then we observe that the composition $\mathtt f(\tilde u)$ is an extension of $\mathtt f(u)$, and therefore one has the inequality 
\[
\| \mathtt f(u) \|_{s,F}^{k+1,\g} 
\leq \| \mathtt f(\tilde u) \|_{s,\R^{\nu+1}}^{k+1,\g} 
\sim \| \mathtt f(\tilde u) \|_{W^{k+1, \infty, \g}(\R^{\nu+1},H^s)}
\]
 by \eqref{0203.1}. 
Then \eqref{0811.10} follows by the Moser composition estimates 
for $\| \ \|_{s, \R^{\nu+1}}^{k+1,\g}$ 
(see for instance Lemma 2.31 in \cite{BM20}), 
together with the equivalence of the norms in \eqref{0203.1}-\eqref{Wg}. 
\end{proof}

\section{Straightening of a transport operator}\label{app:FGMP}

The main results of this appendix are Theorem \ref{thm:as} 
and Lemma \ref{conju.tr}.
The goal is to straighten the
linear quasi-periodic transport operator  
\begin{equation}\label{defX0}
\cT_0 :=\omega\cdot\pa_\vphi + V(\vphi, x) \cdot \nabla\,, 
\end{equation}
to a constant coefficient one 
$\omega\cdot\pa_\vphi + \mathtt{m}_{1}\cdot\nabla$
where $V$ is the function satisfying \eqref{stime tame V B inizio lineariz}.
%up to a small 
%term $p_{\tn}\cdot\nabla_x$, see \eqref{Xtn.gtn} and \eqref{stime.pn.w}. 
%We follow the scheme of Section 4 in \cite{BM21}.

We first recall the weighted-graded Sobolev norm, 
already defined in \eqref{def norm Lip Stein uniform}:
%\eqref{sobolevpesata}: 
for any $u=u(\lambda)\in H^{s}(\T^{\nu+d}_*)$,  $s \in \R$,  
$k_0$-times differentiable with respect to 
$\lambda=(\omega,\th)\in\R^\nu\times[\th_1,\th_2]$, 
we define  the weighted graded Sobolev norm 
\begin{equation*}
\absk{u}{s}:= \sum_{\substack{k\in\N^{\nu+1} \\ 0\leq |k|\leq k_0}} \gamma^{|k|}\sup_{\lambda\in\R^\nu\times[\th_1,\th_2]}\| \pa_\lambda^k u(\lambda) \|_{s-|k|}\,.
\end{equation*}
This norm satisfies usual tame and interpolation estimates. 
%The main reason to use this norm is the estimate \eqref{est:compo} for the composition 
% operator where there is no loss of $ k_0 $-derivatives on the highest norm
% $ | u |_{s}^{k_0,\upsilon} $, unlike 
%the corresponding estimate \eqref{est:compo-loss} for the $ \| \ \|_{s}^{k_0,\upsilon} $.
%This is used in a crucial way to prove   \eqref{ptn.s+b} and then 
%deduce the a-priori bound \eqref{stime.pn.w}
%for the divergence of the high norms of the functions $ p_{\tn} $. 
 In the following we consider
\begin{equation*}
\mathfrak{p}_0:= \frac12(\nu+d)+k_0 + 1\,.
\end{equation*}
We report the following estimates which can be proved by 
adapting the arguments of  \cite{BM20}.

\begin{lemma}\label{proprieta}
The following hold:
\\[1mm]
(i)  For any $s \in \R$, we have $\absk{u}{s}\leq \normk{u}{s}\leq \absk{u}{s+k_0}$. 
\\[1mm]
(ii)  For any  $s\geq \mathfrak{p}_0$, we have 
$
\absk{uv}{s}\leq C(s)\absk{u}{s}\absk{v}{\mathfrak{p}_0} 
+C(\mathfrak{s}_0)\absk{u}{\mathfrak{p}_0}\absk{v}{s}.
$ 
The tame constant $C(s) := C(s,k_0)$ is monotone in $s\geq \mathfrak{p}_0$.
\\[1mm]
(iii) For $ N \geq 1$ and $ \alpha \geq 0 $ 
we have 
$\absk{\Pi_N u}{s}\leq N^{\alpha}\absk{u}{s-\alpha}$  
and 
$\absk{\Pi_N^\perp u}{s}\leq N^{-\alpha}\absk{u}{s+\alpha}$,  $ \forall s \in \R $.
\\[1mm]
(iv) Let 
$\absk{\beta}{2\mathfrak{p}_0+1}\leq \delta $ small enough. 
Then the composition operator $\cB$ defined as 
\begin{equation}\label{diffeo.del.toro}
(\mathcal{B}u)(\vphi,x)=u(\vphi,x+\beta(\vphi,x))\,,
\end{equation}
satisfies the tame estimate, for any $s\geq \mathfrak{p}_0 +1 $,
\begin{equation}\label{est:compo}
\absk{\cB u}{s} \leq C(s) (\absk{u}{s} + \absk{\beta}{s} \absk{u}{\mathfrak{p}_0+1} ) \,. 
\end{equation}
The tame constant $C(s) := C(s,k_0) $ is monotone in $s\geq \mathfrak{p}_0$. 

Moreover the diffeomorphism $x\mapsto x + \beta(\vphi,x)$ 
is invertible and the inverse diffeomorphism $y \mapsto y + \breve{\beta}(\vphi,y)$  
satisfies, for any $s\geq \mathfrak{p}_0$, 
$\absk{\breve{\beta}}{s}\leq C(s) \absk{\beta}{s} $.
\\[1mm]
(v) For any $\epsilon>0$, $a_0, b_0\geq 0$ and $p,q>0$, 
there exists $C_\epsilon=C_\epsilon(p,q)>0$, with $C_1<1$, 
such that 
\[
\absk{u}{a_0+p}\absk{v}{b_0+q}\leq \epsilon\absk{u}{a_0+p+q}\absk{v}{b_0}
+ C_\epsilon\absk{u}{a_0}\absk{v}{b_0+p+q}\,.
\]
\end{lemma}

We now state the  straightening result of the quasi-periodic transport operator. 
%Remind that $ N_{\tn}:=N_0^{\chi^{\tn}} $, $ \chi=3/2 $, $ N_{-1}:=1 $, see \eqref{scala.strai}. 
We fix the constants
\begin{equation}\label{tbta}
\begin{aligned}
& \mathfrak s_0 := s_0 + 2 k_0 \geq 2 \mathfrak p_0 + 1 \quad \text{where} \quad s_0 \quad \text{is given in \eqref{Sone}} \\
& K_0 > 0, \quad \chi := \frac32, \quad K_{- 1} := 1, \quad K_n := K_0^{\chi^n}, \quad n \geq 0\,.  \\
%& \frak a := 3(\tau_0 + 1) + 2 \,, \quad \frak b := \frak a + 1, 
%\begin{equation}\label{tbta}
%\tb:= [\ta] + 2 \in \N \,, \quad 
&\ta \ge 3(\tau_1+1)  \geq 1 \,, \quad \tau_1:= k_0 +(k_0+1)\tau\,.
\end{aligned}
\end{equation}
 %recall that $\tau_0$ is defined in Lemma \ref{lemma:WD}. 
 The constant $K_0$ will be chosen in the following large enough. We remark that it will 
be independent of the diophantine constant  $\gamma$ appearing in \eqref{tDtCn} (see also \eqref{def:DCgt}).

\begin{thm}[{\bf Straightening of transport}]\label{thm:as}
Consider the quasi-periodic traveling operator $\cT_{0}$ in \eqref{defX0}. 
There exists $\sigma=\sigma(\tau, k_{0})>0$ large enough such that if \eqref{ansatz_I0_s0} holds then for any
 $\bar{s}>\mathfrak{s}_0$, there exist 
%$\tau_2> \tau_1 + 1+  \ta$, 
$\tau_{2}>\tau_{1}+1$ (with $\tau_1$ defined in \eqref{tbta}) and
$\delta := \delta(\bar{s},k_0, \tau) \in (0,1) $ and
$K_0 :=K_0(\bar{s},k_0) \in \N $ 
%(with $\tau_1$, $\ta$, $\tb$ defined in \eqref{tbta}) 
such that, for any $\gamma\in(0, 1)$, if 
\begin{equation}\label{small.V.as.AP}
K_0^{\tau_2} \, % \absk{p_0}{\mathfrak{s}_0+\sigma} \, 
\eps \gamma^{-1} \leq \delta \,,
\end{equation}
then, 
%there exists a real constant vector
%$\tm_{1}$ satisfying 
%\begin{equation}\label{tm.est}
%|\mathtt{m}_{1}|^{k_0,\gamma} \leq 2 %C(\mathfrak{s}_0, \tb) 
%\,  \absk{p_0}{\mathfrak{s}_0+\sigma} \,  
%\end{equation}
there exists a
	constant
 vector $\mathtt{m}_{1} := \mathtt{m}_{1}(\omega,\th) \in \R^d $
 and a  quasi-periodic traveling wave  
vector valued function  $\beta(\vphi,x)\in H^{s}$, $\mathfrak{s}_0\leq s\leq \bar{s}$, 
defined for all $ (\omega, \th) \in \mathtt\Omega \times [\mathtt h_1, \mathtt h_2]$ such that the following holds.
For any $(\omega,\th) $ in $\tT\tC_{\infty}(\gamma,\tau) $  (see \eqref{tDtCn})
one has the conjugation
\begin{equation}\label{coniugazione nel teo trasporto}
 \cB^{-1}\circ  \cT_{0}  \circ \cB 
 = \omega\cdot\pa_\vphi + \tm_{1} \cdot \nabla
\end{equation}
where $\cB$ is defined in \eqref{diffeo.del.toro}.
Moreover,  for any $(\omega,\th)\in\mathtt\Omega \times [\mathtt h_1, \mathtt h_2]$,  the following estimates hold:

\noindent
$(i)$ one has 
\begin{equation}\label{tm.est}
\begin{aligned}
& |\mathtt{m}_{1}|^{k_0,\gamma}\lesssim\varepsilon (1+\normk{\cI_0}{s_0+\sigma}) \,,
\\
&| \Delta_{12} \mathtt{m}_{1} |  \lesssim 
\varepsilon \norm{i_1-i_2}_{p+\sigma}  \,,
\end{aligned}
\end{equation}
%for any $p$ as in \eqref{ps0}, %uniformly in  $ \bar\tn$,

\noindent
$(ii)$
The  quasi-periodic traveling wave  
vector valued function  $\beta(\vphi,x)$, 
%defined for all $ (\omega, \th) \in \R^\nu \times [\th_1, \th_2] $,  
satisfies  for any $\mathfrak{s}_0\leq s \leq \bar{s}$,
\begin{equation}\label{gtn.est.better}
\begin{aligned}
\absk{\beta}{s}&\lesssim_{s} K_{0}^{\tau_1}\e\gamma^{-1} (1+\normk{\cI_0}{s+\sigma})\, ,
\\ 
\| \Delta_{12} \beta \|_{p}% \,,\, \| \Delta_{12} \breve \beta\|_{p} 
&\lesssim_{p} K_0^{\tau_1}\e \gamma^{- 1} \| i_1 - i_2 \|_{p + \sigma}\,,
\end{aligned}
%\quad   \|\Delta_{12} \beta\|_{p}\lesssim_{p} N_{0}^{\tau_1}\gamma^{-1}\eps \| i_{1} - i_{2} \|_{p+\sigma} \, , 
\end{equation}
where $p$ is as in   \eqref{ps0}.
%for some  constant $C(s) \geq 1 $ monotone in 
%$s\in[\mathfrak{s}_0,S]$,
% such that, defining the composition operators 
%\begin{equation*}
%(\cG u)(\vphi,x) := u(\vphi,x+ g(\vphi,x))\,, 
%\ \  (\cG^{-1}u)(\vphi,y) := u(\vphi,y+\breve{g}(\vphi,y))\,,
%\end{equation*}

% and $(\cB^{-1}u)(\vphi,y) := u(\vphi,y+\breve{\beta}(\vphi,y))$.

\noindent
$(iii)$
The maps $\cB, \cB^{- 1}$ are momentum preserving  and 
satisfy the tame estimates for every $\mathfrak{s}_{0}\le s \le \bar{s}$
\begin{equation}\label{stima tame cambio variabile rid trasporto}
\begin{aligned}
\| \cB^{\pm 1}  h\|_s^{k_0, \gamma} &\lesssim_{s} 
\| h \|_s^{k_0, \gamma} 
+ \| \cI_{0} \|_{s + \sigma}^{k_0, \gamma} \| h \|_{\mathfrak{s}_{0}}^{k_0, \gamma}\,,  
\\
 \| (\cB^{\pm 1} - {\rm Id})  h\|_s^{k_0, \gamma}%\,,\, \| (\cB^* - {\rm Id})  h\|_s^{k_0, \gamma} 
&\lesssim_{s} 
K_0^{\tau_1}\e \gamma^{- 1}\Big( \| \cI_{0}  \|_{\mathfrak{s}_{0} + \sigma}^{k_0, \gamma}
\| h \|_{s + 1}^{k_0, \gamma} + \| \cI_{0} \|_{s + \sigma}^{k_0, \gamma}
 \| h \|_{\mathfrak{s}_{0} + 1}^{k_0, \gamma} \Big)\,,
\\
 \| \Delta_{12} \cB^{\pm 1} h \|_{p}%\,,\, \| \Delta_{12} \cB^* h \|_{s_1}
& \lesssim_{p} K_0^{\tau_1} \e \gamma^{- 1} \| i_1 - i_2 \|_{p + \sigma} \| h \|_{p + 1}\,. 
\end{aligned}
\end{equation}
%Finally, let $p \geq s_0$ and assume that $i_1, i_2$ 
%satisfy \eqref{ansatz} with $\mu_0 =  \sigma$. 
%Then for any $\lambda = (\omega, \th) \in \tT\tC_{\infty}(2\gamma,\tau)$ one has % the following holds. 
%\begin{equation}\label{stime delta 12 prop trasporto}
%\begin{gathered}
% \| \Delta_{12} \beta \|_{p}% \,,\, \| \Delta_{12} \breve \beta\|_{p} 
%\lesssim_{p} N_0^{\tau_1}\e \gamma^{- 1} \| i_1 - i_2 \|_{p + \sigma}\,, \quad  
% \\
%%& |\Delta_{12} \mathtt m| \lesssim \e \| u_1 - u_2\|_{s_1 + \mu}\,, \\
% \| \Delta_{12} \cB^{\pm 1} h \|_{p}%\,,\, \| \Delta_{12} \cB^* h \|_{s_1}
% \lesssim_{p} N_0^{\tau_1} \e \gamma^{- 1} \| i_1 - i_2 \|_{p + \sigma} \| h \|_{p + 1}\,. 
%\end{gathered}
%\end{equation}
\end{thm}

In order to prove Proposition \ref{thm:as}, 
we closely follow \cite{FGMP19}. First we show the following iterative Lemma. 

\begin{lemma}\label{lemma.iterativo.straightening}
Let $\gamma \in (0, 1)$, $\bar{s} >\mathfrak{s}_0$ and $\ta , \tau_{1}$ as in \eqref{tbta} . Then there exist $\delta = \delta(\bar{s}, k_0, \tau, \nu) \in (0, 1)$, $K_0 = K_0(\bar{s}, k_0, \tau, \nu) > 0$, $\sigma = \sigma (k_0, \tau, \nu) > 0$ such that if \eqref{ansatz_I0_s0}, \eqref{small.V.as.AP} hold with $\mu_0 =\sigma$, then the following statements hold 
for all $n \geq 0$. 

\medskip

\noindent % For any $n \in \N$, there 
$ \bf (S1)_n$ There exists a linear operator 
\begin{equation}\label{def cal Tn transport}
\cT_n := \omega \cdot \partial_\vphi + \mathtt m_{1,n} \cdot \nabla + p_n(\vphi, x) \cdot \nabla
\end{equation}
defined for any $\lambda = (\omega, \th) \in  \mathtt\Omega \times [\mathtt h_1, \mathtt h_2]$ 
where $\mathtt m_{1,n}\in \R^{d}$ is a constant vector satisfying 
\begin{equation}\label{lalla100}
 |\mathtt m_{1,n} |^{k_0, \gamma} 
\lesssim \e\,,\qquad 
\text{and, if $n \geq 1$, } \ 
|\mathtt m_{1, n} - \mathtt m_{1, n - 1} |^{k_0, \gamma} 
\lesssim \e  K_{n - 1}^{- \ta}\,,
%| p_{n - 1}  |_{\mathfrak{s}_0}^{k_0, \gamma}
\end{equation}
and where $p_{n}(\vphi,x)$ is quasi-periodic traveling wave vector valued function satisfying 
for any $\mathfrak{s}_0 \leq s \leq \bar{s}$
\begin{equation}\label{stima mathtt mn an}
\begin{aligned}
& | p_n |_{s}^{k_0, \gamma} 
\leq C_* (s)  K_{n - 1}^{- \ta}\e (1+\normk{\cI_0}{s+\sigma}), 
\quad & & 
| p_n |_{s + \ta + 1}^{k_0, \gamma} 
\leq  C_*(s)K_{n - 1} \e  (1+\normk{\cI_0}{s+\sigma})\,,
%\quad \ \forall s_0 \leq s \leq \bar{s}\,, 
\end{aligned}
\end{equation}
for some constant $C_* (s) = C_*(s, k_0, \tau) > 0$. 
%\begin{equation} \label{correzione mn} % stava dentro \label{stima mathtt mn an}
%|\mathtt m_n - \zeta |^{k_0, \gamma} \lesssim \e, 
%\quad \ \text{and, if $n \geq 1$, } \ 
%|\mathtt m_n - \mathtt m_{n - 1} |^{k_0, \gamma} 
%\lesssim_{k_0} \| a_{n - 1}  \|_{s_0}^{k_0, \gamma}.
%\end{equation}

Define the sets 
$ \mathtt{\Lambda}_{0}^{\rm T}  :=\mathtt\Omega \times [\mathtt h_1, \mathtt h_2]$, 
and, for $n \geq 1$,  
\begin{equation}\label{set.nonres.tn}%\label{cal On gamma}
\begin{aligned}
\mathtt{\Lambda}_{n}^{\rm T} & :=	\mathtt{\Lambda}_{n}^{\gamma,\rm T}%(V_0)
\\&  
:=\big\{ (\omega,\th)\in\mathtt{\Lambda}_{n-1}^{\rm T} \,:\, 
|(\omega-\tV\mathtt{m}_{1,n-1})\cdot\ell|\geq \gamma \jap{\ell}^{-\tau} 
\ \forall\,\ell\in\Z^{\nu}\setminus\{0\} \,, \ |\ell|\leq K_{n-1}  \big\}\, . 
\end{aligned}
\end{equation}	

For $n \geq 1$, there exists an invertible diffeomorphism of the torus 
$\T^d_{\Gamma} \to \T^d_{\Gamma}$, $x \mapsto x + \beta_{n - 1}(\vphi, x)$ 
with inverse $\T^d_{\Gamma} \to \T^d_{\Gamma}$, $y \mapsto y + \check \beta_{n - 1}(\vphi, y)$,
defined for any $\lambda = (\omega, \th) \in \mathtt\Omega \times [\mathtt h_1, \mathtt h_2]$ 
and satisfying the estimates
\begin{equation}\label{stime alpha n tilde alpha n}
\begin{aligned}
& | \beta_{n - 1}|_s^{k_0, \gamma},
 | \breve \beta_{n - 1}|_s^{k_0, \gamma} \lesssim_{s} 
K_{n - 1}^{\tau_1} K_{n - 2}^{- \ta} \e \gamma^{- 1}
(1+\normk{\cI_0}{s+\sigma})
, \quad \forall \mathfrak{s}_0 \leq s \leq \bar{s}\,, \\
&  | \beta_{n - 1}|_{s + \ta + 1}^{k_0, \gamma},
 | \breve \beta_{n - 1}|_{s + \ta + 1}^{k_0, \gamma} 
 \lesssim_{s} 
 K_{n - 1}^{\tau_1} K_{n - 2} \e \gamma^{- 1}
(1+\normk{\cI_0}{s+\sigma})
,  \quad \forall \mathfrak{s}_0 \leq s \leq \bar{s} \,, 
\end{aligned}
\end{equation}
(with $\tau_1$ defined in \eqref{tbta}), such that the following holds.
By defining the operator 
$$
\cB_{n - 1} : h(\vphi, x) \mapsto h(\vphi, x + \beta_{n - 1}(\vphi, x))
$$
whose inverse inverse is 
$$
\cB_{n - 1}^{- 1} : h(\vphi, x) \mapsto h(\vphi, y + \breve \beta_{n - 1}(\vphi, y))
$$
then,  for any $\lambda = (\omega, \th) \in \mathtt{\Lambda}_{n}^{\rm T}$, 
one has
the conjugation
\begin{equation}\label{coniugazione cal A n - 1 cal T n - 1}
\cT_n = \cB_{n - 1}^{- 1}\cT_{n - 1} \cB_{n - 1}\,. 
\end{equation}
Furthermore, $p_n (\vphi, x), \beta_{n - 1}(\vphi, x)$ and $\breve \beta_{n - 1}(\vphi, x)$ are quasi-periodic traveling wave vector valued functions, implying that 
$\cT_n,% is a momentum preserving operator and
\cB_{n - 1}$ and $ \cB_{n - 1}^{- 1}$ are momentum preserving operators.

$ \bf (S2)_n$ Given two tori $i_1,i_2$ satisfying \eqref{ansatz_I0_s0} and setting 
% Given   $ p_{0}(i_1), p_{0}(i_2) $ let
 $ \Delta_{12} p_n :=  p_{n}(i_1)-p_{n}(i_2) $, we have that
%For any $p \in [s_0+1, \bar{s}] $, there exist $C(p)>0$ and 
% $\delta'(p)\in(0,1)$ such that
% if 
%%\begin{equation}
%%N_0^{\tau_2} \sup_{(\omega, \th) \in \R^\nu \times [\th_1, \th_2]}
%%\big( \| p_{0}(i_1) \|_{p + \ta+1}+\| p_{0}(i_2) \|_{p+\ta+1} \big)
%%\gamma^{-1} \leq \delta'(p)  \,,
%%\end{equation}  
%\begin{equation}\label{small12}
%N_0^{\tau_2} ( \|V(i_{1})\|_{p+\ta+1} +  \|V(i_{2})\|_{p+\ta+1})\gamma^{-1} \leq \delta'(p)  \,,
%\end{equation} 
%
%then, 
for all $(\omega, \th) \in \mathtt\Omega \times [\mathtt h_1, \mathtt h_2]$,
%$(\omega, \gamma) \in  \mathtt{\Lambda}_{\tn}^{\upsilon, \rm T}(p_{0,1}) \cap \mathtt{\Lambda}_{\tn}^{\upsilon,\rm T}(p_{0,2}) $, 
%\begin{align}
%\label{estp121}
%&\|\Delta_{12} p_n\|_{p-1}  \lesssim_{p}
%   \eps N_{n-1}^{-\ta} \| i_{1} - i_{2} \|_{p+\ta+1 +\sigma} \, , 
%\quad \|\Delta_{12} p_n\|_{p+\ta+1} \lesssim_{p}
%N_{n-1}\|\Delta_{12} p_0\|_{p+\ta+1} \\
%& |\Delta_{12}(\mathtt{m}_{1, n+1}-\mathtt{m}_{1, n})|
% \leq \|\Delta_{12}p_n\|_{s_0}\,, 
% \quad |\Delta_{12} \mathtt{m}_{1,n}| \lesssim_{p}
%  \|\Delta_{12} p_0\|_{s_0 }\,.\label{estm121}
%\end{align}	
\begin{align}
\label{estp121}
&\|\Delta_{12} p_n\|_{p}  \lesssim_{p}
  K_{n-1}^{-\ta} \e\| i_{1} - i_{2} \|_{s_0+\sigma} \, , 
%\quad \|\Delta_{12} p_n\|_{p+\ta+1} \lesssim_{p}
%\eps N_{n-1}\| i_{1} - i_{2} \|_{p+\ta+1+\sigma} 
\\
& |\Delta_{12}(\mathtt{m}_{1, n}-\mathtt{m}_{1, n-1})|
 \leq \|\Delta_{12}p_n\|_{s_0}\,, 
 \quad |\Delta_{12} \mathtt{m}_{1,n-1}| 
 \lesssim_{p}
 \eps \| i_{1} - i_{2} \|_{s_{0}+\sigma}\,. \label{estm121}
\end{align}
Moreover,  %for any $s \geq \mathfrak{s}_0$,
%\begin{align}
%& \|\Delta_{12} \beta_n\|_{p} \lesssim_{p} 
%\g^{-1}\big( \|\Pi_{N_{n}}\Delta_{12}p_n\|_{p+\tau} 
%+ \g^{-1} |\Delta_{12}\mathtt{m}_{1, n}
%|\|\Pi_{N_{n}} p_{n}(i_{2}) \|_{p+2\tau+1} \big)\,. \label{estg12}
%\end{align}
\begin{equation}\label{estg12}
\begin{aligned}
% \|\Delta_{12} \beta_n\|_{p} 
% &\lesssim_{p} 
%\g^{-1}\big( \eps N_{n-1}^{-\ta} \| i_{1} - i_{2} \|_{p + \sigma}%_{p + \tau +\ta+1 + \sigma}
%+ \g^{-1} 
% \eps \| i_{1} - i_{2} \|_{s_{0}+\sigma}
% \e N_{n - 1}^{- \ta} (1+\normk{\cI_0}{p + \sigma}%{p + 2\tau +1 +\sigma})
% \big)
% \\ 
%&
% \lesssim_{p} 
  \|\Delta_{12} \beta_{n-1}\|_{p}& \lesssim_{s}
  \e\gamma^{-1}K_{n-1}^{\tau_1}K_{n-2}^{-\mathtt{a}}   \e\gamma^{-1}\|i_1-i_2\|_{s_0+\s}\,.
% \e\gamma^{-2}\|i_1-i_2\|_{s_0+\s}\|\Pi_{N_n} p_n(i_2)\|_{p+2\tau+1}
% \\&+
% \gamma^{-1}\|\Pi_{N_n} \Delta_{12}p_n\|_{p+2\tau+1}
% %\eps \| i _{1} - i_{2}\|_{p + \sigma}%_{p+ \ta + 1 + \sigma}
 \end{aligned}
\end{equation}

where $p$ is as in   \eqref{ps0}.

\end{lemma}

\begin{rmk}
In Lemma \ref{lemma.iterativo.straightening} 
we prove that
the norms $ |p_n|_s^{k_0, \g} $ satisfy 
inequalities typical of
a Nash-Moser iterative scheme, which converges 
under the  smallness low norm condition \eqref{small.V.as.AP}.

\end{rmk}

\begin{proof}
{\sc Proof of the statement for $n = 0$.} 
The claimed statements for $n = 0$ follows directly by defining 
$p_{0} = V$.  Then \eqref{stima mathtt mn an} follows by \eqref{stime tame V B inizio lineariz}. 

\medskip

\noindent
{\sc Proof of the induction step.}
Now assume that the claimed properties hold for some $n \geq 0$ and let us prove them at the step $n + 1$. We look for a diffeomorphism of the torus $\T^d_{\Gamma} \to \T^d_{\Gamma}$, $x \mapsto x + \beta_n(\vphi, x)$, wih inverse given by $y \mapsto y + \breve \beta_n(\vphi, y)$ such that defining 
$$
\cB_n : h(\vphi, x) \mapsto h(\vphi, x + \beta_n(\vphi, x)), \quad \cB_n^{- 1} : h(\vphi, x) \mapsto h(\vphi, x + \breve \beta_n(\vphi, x))
$$
the operator 
%$\cT_{ n + 1} := \cA_n^{- 1} \cT_n \cA_n$ 
$\cA_n^{- 1} \cT_n \cA_n$ 
has the desired properties. One computes 
\begin{equation}\label{primo cal L n + 1}
\begin{aligned}
\cB_n^{- 1} \cT_n \cB_n
%\cT_{n + 1} 
& = \omega \cdot \partial_\vphi + \mathtt m_{1,n} \cdot \nabla 
+ \cB_n^{- 1} \big[ \omega \cdot \partial_\vphi \beta_{n} + \mathtt m_{1,n} \cdot \nabla \beta_{n} + p_n + p_n \cdot \nabla \beta_{n} \big] \cdot \nabla  \\
& = \omega \cdot \partial_\vphi + \mathtt m_{1,n} \cdot \nabla 
+ \cB_n^{- 1} \big[ \omega \cdot \partial_\vphi \beta_n + \mathtt m_{1,n} \cdot \nabla \beta_{n} + \Pi_{N_n}p_n \big] \cdot \nabla  + p_{n + 1} \cdot \nabla
\end{aligned}
\end{equation}
where 
\begin{equation}\label{def a n + 1}
p_{n + 1} := \cB_n^{- 1} f_n, \qquad 
f_n :=  \Pi_{K_n}^\bot p_n +  p_n \cdot \nabla \beta_{n},
\end{equation}
and the projectors $\Pi_{K_n}, \Pi_{K_n}^\bot$ are defined by \eqref{proiettore-oper}.
For any $(\omega, \th) \in \mathtt{\Lambda}_{n+1}^{\rm T}$, we solve the homological equation 
\begin{equation}\label{eq omologica trasporto}
\omega \cdot \partial_\vphi \beta_{n} + \mathtt m_{1, n} \cdot \nabla \beta_{n} + \Pi_{K_n}p_n = \langle p_n \rangle_{\vphi, x}
\end{equation}
where $ \jap{  p_{n} }_{\vphi,x}  $  is the average of $p_{n} $ defined as 
\begin{equation}\label{lalla4}
\jap{p_{n} }_{\vphi,x}:=\frac{1}{(2\pi)^{\nu}|\T_{\Gamma}^{d}|}
\int_{\T_{*}^{\nu+d}}p_{n}(\vphi,x)d\vphi dx\,,
\qquad \jap{p_{n}}(x):=\frac{1}{(2\pi)^{\nu}}\int_{\T^{\nu}}p_{n}(\vphi,x)d\vphi\,.
\end{equation}
Notice that, since $ p_{n} $ is a quasi-periodic traveling wave then 
$\jap{p_{n} }_{\vphi,x}=\jap{p_{n} }_{\vphi}$.
Moreover, since $p_{n}$ is a traveling wave and we look for traveling 
wave solution $\beta_{n}$ of \eqref{eq omologica trasporto}, we shall write, 
recalling  Def. \ref{def:quasitravelling} and \eqref{velocityvec1}-\eqref{velocityvec2} 
(see also subsection \ref{sec:travelfunction}),
\begin{equation}\label{gtn.def}
p_{n}(\vphi,x)=P_{n}(\vphi-\tV x)\,,\quad \beta_{n}(\vphi,x)=B_{n}(\vphi-\tV x)\,,
\end{equation}
for some (known) function $P_{n}$ and some (unknown) function $B_{n}$,
both belonging to $H^{s}(\T^{\nu};\C)$.
With this formalism, equation \eqref{eq omologica trasporto} reads
\[
 \big((\omega -\tV\mathtt{m}_{1, n})\cdot\pa_\vphi \big)	B_{n}(\vphi,x) 
 +\Pi_{K_{n}} P_{n} = \jap{p_{n}}_{\vphi} =\jap{P_{n}}_{\theta}
 :=\frac{1}{(2\pi)^{\nu}}\int_{\T^{\nu}}P_{n}(\theta)d\theta
\]
So we define
\begin{equation}\label{gtn.def2}
B_{n}(\theta) := -  \big((\omega -\tV\mathtt{m}_{1, n})\cdot\pa_\vphi \big)_{\rm ext}^{-1}
(\Pi_{K_{n}} P_{n}- \jap{P_{n}}_{\theta}) 
\end{equation}
where the operator $((\omega -\tV\mathtt{m}_{1, n})\cdot\pa_\vphi )_{\rm ext}^{-1}$ 
is introduced in \eqref{def ompaph-1 ext2}.
%\eqref{def ompaph-1 ext} with $\omega\rightsquigarrow (\omega -\tV\mathtt{m}_{1, n})$.
%The operator \eqref{primo cal L n + 1} takes the form 
We define 
\begin{equation}\label{secondo cal L n + 1}
\cT_{n + 1} 
:= \omega \cdot \partial_\vphi + \mathtt m_{1, n + 1} \cdot \nabla + p_{n + 1} \cdot \nabla 
\end{equation}
where 
\begin{equation}\label{def mathtt m n + 1}
\mathtt m_{1, n + 1} := \mathtt m_{1, n} + \langle p_n \rangle_{\vphi, x}\,.
\end{equation}
We observe that the function $\beta_{n}(\vphi,x)$
defined in \eqref{gtn.def} with $B_{n}$ given in \eqref{gtn.def2} is a quasi-periodic traveling wave 
defined for all $(\om,\th) \in \mathtt\Omega \times [\mathtt h_1, \mathtt h_2]$.
In particular,  for all $(\om,\th) \in \mathtt{\Lambda}_{n+1}^{\rm T}$, $\beta_{n}$ solves the equation
\eqref{eq omologica trasporto} and so, recalling \eqref{primo cal L n + 1} and \eqref{secondo cal L n + 1},
one has the conjugation $\cB_n^{-1} \cT_n \cB_n = \cT_{n+1}$.
%We observe that $\cT_{n+1}$ is defined for all $(\om,\th) \in \R^{\nu}\times [\th_1,\th_2]$, 
%and, for $(\om,\th) \in \mathtt{\Lambda}_{n+1}^{\rm T}$, 
%one has $\cB_n^{-1} \cT_n \cB_n = \cT_{n+1}$.
%Clearly $\beta_{n} (\vphi, x; \omega, \th)$ is $\cC^\infty$ in $(\vphi, x)$ 
%and $k_0$ times differentiable in $(\omega, \th) \in \R^{\nu + 1}$. 

Furthermore, by Lemma \ref{lemma:WD2}, and by the smoothing property \eqref{p2-proi}, for any $s \geq 0$, one has 
\begin{equation}\label{prima stima alpha n an}
\begin{aligned}
& | \beta_{n} |_s^{k_0, \gamma} \lesssim  \gamma^{- 1} | \Pi_{K_n} p_n |_{s + \tau_1}^{k_0, \gamma} 
\lesssim K_n^{\tau_1} \gamma^{- 1} | p_n|_s^{k_0, \gamma}\,,
%\\
%& | \nabla \beta_{n}|_s^{k_0, \gamma} \lesssim \gamma^{- 1}| \Pi_{N_n} p_n |_{s + \tau_1 + 1}^{k_0, \gamma}  \lesssim N_n^{\tau_1 + 1} \gamma^{- 1} \| p_n\|_s^{k_0, \gamma}\,. 
\end{aligned}
\end{equation}
which, together with 
the induction estimates on \eqref{stima mathtt mn an} on $p_n$, 
implies that, for any $s_0 \leq s \leq \bar{s}$,
\begin{equation}\label{stime alpha n + 1}
\begin{aligned}
 |  \beta_{n} |_s^{k_0, \gamma}
%\,,\, | \breve  \beta_{n} |_s^{k_0, \gamma} 
&\lesssim_{s} K_{n }^{\tau_1} K_{n - 1}^{-\ta} \e \gamma^{- 1}
(1+\normk{\cI_0}{s+\sigma}) \,, 
\\
 |  \beta_{n} |_{s + \ta +1}^{k_0, \gamma} 
%\,,\, | \breve  \beta_{n} |_{s + \ta +1}^{k_0, \gamma}
& \lesssim_{s} K_n^{\tau_1} K_{n - 1} \e \gamma^{- 1}
(1+\normk{\cI_0}{s+\sigma})\,.
\end{aligned}
\end{equation}
The bounds \eqref{stime alpha n + 1} are the estimates \eqref{stime alpha n tilde alpha n} at the step $n + 1$
for the function $\beta_{n}$.
Note that, using the definition of the constant $\mathtt{a}$ in \eqref{tbta} 
and the ansatz \eqref{ansatz}, 
from \eqref{stime alpha n + 1} with $s=  \mathfrak s_0 $ one deduces that
\begin{equation}\label{stima alpha n s0}
 | \beta_{n} |_{ \mathfrak{s}_0}^{k_0, \gamma}
 %\,,\, | \breve \beta_{n} |_{\mathfrak{s}_0}^{k_0, \gamma} 
 \lesssim K_0^{\tau_1} \e \gamma^{- 1}\,.
\end{equation}
By the smallness condition \eqref{small.V.as.AP} and 
\eqref{stima alpha n s0} we have that Lemma \ref{proprieta} applies
and so the diffeomorphism  $x \mapsto x + \beta_n(\vphi, x)$ is invertible  with inverse 
$y \mapsto y + \breve \beta_n(\vphi, y)$ satisfying 
\begin{equation}\label{lalla2}
  | \breve  \beta_{n} |_s^{k_0, \gamma} \lesssim_{s}  |  \beta_{n} |_s^{k_0, \gamma}\,.
\end{equation}
The latter bound and \eqref{stime alpha n + 1} imply the estimates 
\eqref{stime alpha n tilde alpha n} at the step $n + 1$
for  $\breve\beta_{n}$.
Moreover,  using  the smallness condition \eqref{small.V.as.AP} (choosing $\tau_2 > \tau_1$), together with Lemma \ref{lemma:LS norms} and the estimate \eqref{prima stima alpha n an} leads to the estimate
\begin{equation}\label{stima cal An pm 1}
| \cB_n^{\pm 1} h |_s^{k_0, \gamma} \lesssim_{s} | h |_s^{k_0, \gamma} 
+ K_n^{\tau_1} \gamma^{- 1} | p_n|_s^{k_0, \gamma}
 | h |_{\mathfrak{s}_0}^{k_0, \gamma}\,, \quad \forall \mathfrak{s}_0 \leq s \leq \bar{s} + \ta +1\,.  
\end{equation}

% The latter estimate, together with Lemma \ref{lemma:LS norms} and the induction estimates on \eqref{stima mathtt mn an} on $p_n$, imply that for any $s_0 \leq s \leq \bar{s}$
%\begin{equation}\label{stime alpha n + 1}
%\begin{aligned}
%& |  \beta_{n} |_s^{k_0, \gamma}\,,\, | \breve  \beta_{n} |_s^{k_0, \gamma} 
%\lesssim_{s} N_{n }^{\tau_1} N_{n - 1}^{-\ta} \e \gamma^{- 1}
%(1+\normk{\cI_0}{s+\sigma}) \,, \\
%& |  \beta_{n} |_{s + \ta +1}^{k_0, \gamma} \,,\,
% | \breve  \beta_{n} |_{s + \ta +1}^{k_0, \gamma}
% \lesssim_{s} N_n^{\tau_1} N_{n - 1} \e \gamma^{- 1}
%(1+\normk{\cI_0}{s+\sigma})
%\end{aligned}
%\end{equation}
%which are the estimates \eqref{stime alpha n tilde alpha n} at the step $n + 1$. 
On the solution $\beta_{n}(\vphi,x)=B_{n}(\vphi-\mathtt{V}x)$  (see \eqref{gtn.def2}) 
we also have estimates on the Lipschitz variation w.r.t. the embedding $i$.
More precisely, 
using \eqref{2802.2nuova2} in Lemma \ref{lemma:WD2}, and the inductive assumption \eqref{estm121} for 
$\Delta_{12}\mathtt{m}_{1,n}$ one gets
\begin{equation}\label{venezia1}
\begin{aligned}
 \|\Delta_{12} \beta_n\|_{p} &\lesssim_{s}
 \e\gamma^{-2}\|i_1-i_2\|_{s_0+\s}\|\Pi_{K_n} p_n(i_2)\|_{p+2\tau+1}+
 \gamma^{-1}\|\Pi_{K_n} \Delta_{12}p_n\|_{p+2\tau+1}
 \\&
 \stackrel{\eqref{p2-proi}}{\lesssim}
  \e\gamma^{-2} K_{n}^{2\tau+1}\|i_1-i_2\|_{s_0+\s}\| p_n(i_2)\|_{p}
  +\gamma^{-1}K_{n}^{2\tau+1}\| \Delta_{12}p_n\|_{p}
 \end{aligned}
\end{equation}
for any $p$ as in \eqref{ps0}.
%which is the \eqref{estg12}.

%
%
%
%Note that, using the definition of the constant $\frak a$ in \eqref{tbta} 
%and the ansatz \eqref{ansatz}, 
%from \eqref{stime alpha n + 1} with $s=s_0$ one deduces that
%\begin{equation}\label{stima alpha n s0}
% | \beta_{n} |_{\mathfrak{s}_0}^{k_0, \gamma}\,,\, | \breve \beta_{n} |_{\mathfrak{s}_0}^{k_0, \gamma} \lesssim N_0^{\tau_1} \e \gamma^{- 1}\,.
%\end{equation}

We now estimate the function $p_{n + 1}$ defined in \eqref{def a n + 1}. 
First, we estimate $f_{n }$. 
By \eqref{tbta}, \eqref{stima mathtt mn an} and using also the ansatz \eqref{ansatz} and the smallness condition \eqref{small.V.as.AP}, 
one has that 
\begin{equation}\label{a n + 1 Nn gamma}
K_n^{\tau_1 + 1} \gamma^{- 1} \| p_n \|_{s_0}^{k_0, \gamma} \leq 1\,.
\end{equation} 
By \eqref{p1-pr}, \eqref{p2-proi}, \eqref{p3-proi}, 
\eqref{prima stima alpha n an}, \eqref{a n + 1 Nn gamma}
one has 
\begin{equation}\label{stima f n + 1 trasporto}
\begin{aligned}
| f_n |_{\mathfrak{s}_0}^{k_0, \gamma} & \lesssim
 | p_n |_{\mathfrak{s}_0}^{k_0, \gamma}\,, \\
| f_{n } |_s^{k_0, \gamma} & \lesssim_{s} K_n^{- \ta -1}
 | p_n |_{s + \ta +1}^{k_0, \gamma} + K_n^{\tau_1 + 1}
  \gamma^{- 1} | p_n|_s^{k_0, \gamma}
   | p_n|_{\mathfrak{s}_0}^{k_0, \gamma}\,, 
\quad \forall s_0 < s \leq \bar s, \\
| f_{n }|_{s + \ta +1}^{k_0, \gamma} 
& \lesssim_{s} | p_n|_{s + \ta +1}^{k_0, \gamma} 
\big(1 + K_n^{\tau_1 + 1} \g^{-1} | p_n|_{\mathfrak{s}_0}^{k_0, \gamma} \big) 
\lesssim_{s} % \stackrel{\eqref{a n + 1 Nn gamma}}{\lesssim_s} 
| p_n |_{s + \ta +1}^{k_0, \gamma}, \quad \forall \mathfrak{s}_0 \leq s \leq \bar s. 
\end{aligned}
\end{equation}
Hence \eqref{stima cal An pm 1}-\eqref{stima f n + 1 trasporto} 
imply that, for any $\mathfrak{s}_0 \leq s \leq \bar s$, 
\begin{equation}\label{stime induttive a n + 1}
\begin{aligned}
| p_{n + 1} |_s^{k_0, \gamma} 
& \lesssim_{s} K_n^{- \ta -1} 
| p_n|_{s + \ta +1}^{k_0, \gamma} 
+ K_n^{\tau_1 + 1} \gamma^{- 1} \| p_n\|_s^{k_0, \gamma}
| p_n|_{\mathfrak{s}_0}^{k_0, \gamma}\,, \\
| p_{n + 1}|_{s + \ta +1}^{k_0, \gamma} 
& \lesssim_{s} | p_n |_{s + \ta +1}^{k_0, \gamma} 
\end{aligned}
\end{equation}
and using the definition of the constant $\ta$ in \eqref{tbta} and the induction estimates on $p_n$ one deduces the estimate \eqref{stima mathtt mn an} for $p_{n + 1}$. The estimates \eqref{stima mathtt mn an} for $\mathtt m_{n + 1}$ follows by its definition \eqref{def mathtt m n + 1}, by the induction estimate on $p_n$ and by using a telescoping argument. 
The proof of $\bf (S1)_{n+1}$ is complete.

\vspace{0.3em}
Let us now prove $\bf (S2)_{n+1}$  assuming $\bf (S1)_{n}$, $\bf (S2)_{n}$. 

\noindent
We start by considering $f_n$ in  the  \eqref{def a n + 1}\,.
Using the smoothing estimates \eqref{p3-proi} and the inductive assumption \eqref{estp121},
we get
\[
\begin{aligned}
\|\Delta_{12}\Pi_{K_{n}}^{\perp}p_n\|_{p}&\lesssim K_{n}^{-\mathtt{a}-1} 
\|\Delta_{12}p_{n}\|_{p+\mathtt{a}+1}
\lesssim \eps K_{n}^{-\mathtt{a}}\| i_{1} - i_{2} \|_{s_0+\sigma}\,,
%\\
%\|\Delta_{12} \Pi_{N_{n}}^{\perp}p_n\|_{p+\ta+1}& \lesssim_{p}
%\eps N_{n}\| i_{1} - i_{2} \|_{p+\ta+1+\sigma}\,.
\end{aligned}
\]
where in the last inequality we used that (recall \eqref{ps0}) $p\ll s_0+\s$ for some $\s>0$ large.
We also note that
\[
\Delta_{12}(p_n \cdot \nabla \beta_{n})=(\Delta_{12}p_{n})\cdot\nabla\beta_{n}(i_1)+
p_{n}(i_2)\cdot\nabla\Delta_{12}\beta_{n}\,.
\]
Therefore, using \eqref{p1-pr}, the inductive assumptions \ref{estp121} and \eqref{estg12},
estimates \eqref{stime alpha n + 1} (in low norm) and \eqref{stima mathtt mn an} on $p_{n}$,
one deduces
\[
\begin{aligned}
\|\Delta_{12}(p_n \cdot \nabla \beta_{n})\|_{p}&\lesssim
%\|(\Delta_{12}p_{n})\|_{p-1}\|\beta_{n}(i_1)\|_{s_0}
%+\|(\Delta_{12}p_{n})\|_{s_0}\|\beta_{n}(i_1)\|_{p}
%\\
%&+\|p_{n}(i_2)\|_{p-1}\|\Delta_{12}\beta_{n}\|_{s_0}
%+\|p_{n}(i_2)\|_{s_0}\|\Delta_{12}\beta_{n}\|_{p}
  \eps  K_{n-1}^{-2\ta} K_{n}^{\tau_1}\| i_{1} - i_{2} \|_{s_0 +\sigma}\e\gamma^{-1}
  \lesssim K_{n}^{-\mathtt{a}}\e\| i_{1} - i_{2} \|_{s_0 +\sigma}\,,
\end{aligned}
\]
where in the last inequality we used that $\e\gamma^{-1}K_{n-1}^{-2\ta} K_{n}^{\tau_1+\mathtt{a}}\ll1$
thanks to the choice of $\mathtt{a}$ in \eqref{tbta} 
and the smallness assumption \eqref{small.V.as.AP}.
By the estimates above one deduce that, for $p$ satisfying \eqref{def a n + 1}, 
$f_{n}$ in \eqref{def a n + 1} satisfies 
\[
\|\Delta_{12}f_{n}\|_{p}\lesssim K_{n}^{-\mathtt{a}}\e\| i_{1} - i_{2} \|_{s_0 +\sigma}\,.
\]
Now let us consider the function $p_{n+1}$ in \eqref{def a n + 1} which has the form
\begin{equation}\label{lalla3}
p_{n+1}(\vphi,x)=f_{n}(\vphi,x+\breve{\beta}_{n}(\vphi,x))\,.
\end{equation}
{\bf Claim.} One has
\begin{equation}\label{lalla1}
  \|\Delta_{12} \breve\beta_n\|_{p} \lesssim_{s}
K_{n}^{\tau_1}K_{n-1}^{-\mathtt{a}}   \e\gamma^{-1}\|i_1-i_2\|_{s_0+\s}\,.
\end{equation}
Indeed, recall that $\breve\beta_{n}$ is constructed as the inverse diffeomorphism of $x+\beta_{n}$, and so
by defining
$\beta_{n,j}:=\beta_{n}(i_j)$ and $\breve\beta_{n,j}:=\breve\beta_{n}(i_j)$ for $j=1,2$ one has 
\[
\Delta_{12}\breve\beta_{n}=\beta_{n,1}(\vphi,x+\breve\beta_{n,2})
-\beta_{n,1}(\vphi,x+\breve\beta_{n,1})+(\Delta\beta_{n})(\vphi,x+\breve\beta_{n,2})\,.
\]
Therefore we deduce \eqref{lalla1} using the inductive assumption \eqref{estg12} on $\beta_{n}$,
estimates \eqref{stime alpha n + 1}, \eqref{lalla2},  Lemma \ref{lemma:LS norms} 
and the smallness condition \eqref{small.V.as.AP}.
The estimates \eqref{estp121} on $p_{n+1}$ in \eqref{lalla3}
follows reasoning as above, using again Lemma \ref{lemma:LS norms} and estimate \eqref{lalla1}.
Let us check \eqref{estm121}.
Recalling 
\eqref{def mathtt m n + 1} we deduce the estimate on $\Delta_{12}\mathtt{m}_{1,n+1}$
by formula \eqref{lalla4} and estimates \eqref{estp121} on $p_{n}$. 
Similarly the first bound in \eqref{estm121}
with $n\rightsquigarrow n+1$ follows using \eqref{estp121} at the step $n+1$.

It remains to prove \eqref{estg12} at the step $n+1$.
We recall that the traveling wave $\beta_{n}(\vphi,x)$ is constructed in \eqref{gtn.def2}.
%The the traveling wave $\beta_{n+1}(\vphi,x)=B_{n+1}(\vphi-\mathtt{V}x)$ is constructed 
%(recall formula \eqref{gtn.def2}) 
%as
%\begin{equation}\label{gtn.def2n+1}
%B_{n+1}(\theta) := -  \big((\omega -\tV\mathtt{m}_{1, n+1})\cdot\pa_\vphi \big)_{\rm ext}^{-1}
%(\Pi_{N_{n}} P_{n+1}- \jap{P_{n+1}}_{\theta}) \,.
%\end{equation}
%with $p_{n+1}(\vphi,x)=P_{n+1}(\vphi-\mathtt{V}x)$ and where $\mathtt{m}_{1, n+1}$ 
%is the constant vector constructed at the step $n+1$ (see \eqref{def mathtt m n + 1} 
%with $n\rightsquigarrow n+1$). t
%In particular, we have already proved estimates 
%\eqref{estp121}-\eqref{estm121} on $p_{n+1}$ and $\mathtt{m}_{1, n+1}$.
Hence using estimate \eqref{venezia1} and the inductive assumptions
\eqref{estp121}-\eqref{estm121} on $p_{n}$  we get
\[
\begin{aligned}
 \|\Delta_{12} \beta_{n}\|_{p} &\lesssim_{s}
  \e\gamma^{-2} K_{n}^{2\tau+1}\| p_{n}(i_2)\|_{p}
  +\gamma^{-1}K_{n}^{2\tau+1}\| \Delta_{12}p_{n}\|_{p}
  \\&
 \stackrel{\eqref{stima mathtt mn an}_{n}, \eqref{estp121}_{n}}{ \lesssim}
   \e\gamma^{-2} K_{n}^{2\tau+1}\e K_{n-1 }^{- \ta}\|i_1-i_2\|_{s_0+\s}+
 \e \gamma^{-1}K_{n}^{2\tau+1}  K_{n-1}^{-\ta} \| i_{1} - i_{2} \|_{s_0+\sigma}
  \\&
\lesssim  \e\gamma^{-1} K_{n}^{2\tau+1} K_{n-1 }^{- \ta}\|i_1-i_2\|_{s_0+\s}
\end{aligned}
\]
where in the las inequality we used that $\e\gamma^{-1}\ll1$.  
The latter bound is the \eqref{estg12} at the step $n+1$.
The proof of Lemma \ref{lemma.iterativo.straightening} is concluded.
\end{proof}

We then define 
\begin{equation}\label{def widetilde cal An}
\widetilde\cB_n := \cB_0 \circ \cB_1 \circ \ldots \circ \cB_n, \quad \text{with inverse} \quad \widetilde\cB_n^{- 1} = \cB_n^{- 1} \circ \cB_{n - 1}^{- 1} \circ \ldots \circ \cB_0^{- 1}\,. 
\end{equation}

\begin{lemma}\label{lemma tilde cal An}
Let $\bar s > \mathfrak{s}_0$, $\gamma \in (0, 1)$. Then there exist $\delta = \delta(\bar s, k_0, \tau, \nu) \in (0, 1)$, $\sigma = \sigma(k_0, \tau, \nu) > 0$ such that if \eqref{ansatz} holds with $\mu_0 = \sigma$ and if \eqref{small.V.as.AP} holds, then the following properties hold.

\noindent
$(i)$
$$
\widetilde\cB_n h (\vphi, x) = h(\vphi, x + \alpha_n(\vphi, x))\,, \quad \widetilde\cB_n^{- 1} h(\vphi, y) = h(\vphi, x + \breve \alpha_n(\vphi, x))
$$ 
where, for all $(\omega, \th) \in \mathtt\Omega \times [\mathtt h_1, \mathtt h_2]$ and any $\mathfrak{s}_0 \leq s \leq \bar s$, 
\begin{equation}\label{stima beta n beta n - 1 trasporto}
\begin{aligned}
| \alpha_0 |_{s}^{k_0, \gamma}\,,\, | \breve \alpha_0|_{s}^{k_0, \gamma} 
& \lesssim_{s} K_{0}^{\tau_1} \e \gamma^{- 1} 
(1+\normk{\cI_0}{s+\sigma}) \,,
\\
| \alpha_n - \alpha_{n - 1} |_s^{k_0, \gamma}\,,\, 
|\breve \alpha_n - \breve \alpha_{n - 1} |_s^{k_0, \gamma} 
& \lesssim_{s} K_{n }^{\tau_1} N_{n - 1}^{- \ta}\e \gamma^{- 1} 
(1+\normk{\cI_0}{s+\sigma}), \quad n \geq 1. 
\end{aligned}
\end{equation}
As a consequence, 
\begin{equation}\label{bound solo beta n}
| \alpha_n |_s^{k_0, \gamma} \lesssim_{s} K_0^{\tau_1} \e \gamma^{- 1}
(1+\normk{\cI_0}{s+\sigma}), \quad \forall s_0 \leq s \leq \bar s\,. 
\end{equation}
Moreover, one has that
{\begin{equation}\label{bound delta12 beta n}
\begin{aligned}
| \Delta_{12} \alpha_n |_p&\lesssim_{p}  \e\gamma^{-1} K_0^{\tau_1}\|i_1-i_2\|_{s_0+\s} \,,
\\
\|\Delta_{12}(\alpha_{n}-\alpha_{n-1})\|_{p}&\lesssim_{p}
\e\gamma^{-1}K_{n}^{\tau_1}K_{n-1}^{-\mathtt{a}}\|i_1-i_2\|_{s_0+\s}\,,
\end{aligned}
\end{equation}}
for $p$ as in \eqref{ps0}.
Furthermore, $\alpha_n, \breve \alpha_n$ are quasi-periodic traveling vector valued functions. 

\noindent
$(ii)$ For any $ \leq s \leq \bar s$, 
the sequence $(\alpha_n)_{n \in \N}$ (resp. $ (\breve \alpha_n)_{n \in \N}$) 
is a Cauchy sequence with respect to the norm $| \cdot |_s^{k_0, \gamma}$ 
and it converges to some limit $\beta$ (resp. $\breve \beta$). 
Furthermore $\beta, \breve \beta $ are quasi-periodic traveling vector valued functions.  
and, for any $\mathfrak{s}_0 \leq s \leq \bar s$, $n \geq 0$, one has
\begin{equation}\label{stima beta n - beta infty}
\begin{aligned}
| \beta- \alpha_{n} |_s^{k_0, \gamma}\,,\, | \breve \beta - \breve \alpha_{n}|_s^{k_0, \gamma} 
& \lesssim_{s} 
K_{n + 1}^{\tau_1} K_{n }^{1 - \ta}\e \gamma^{- 1} (1+\normk{\cI_0}{s+\sigma})  \,, 
\\
| \beta |_s^{k_0, \gamma}, | \breve \beta |_s^{k_0, \gamma} 
& \lesssim_{s} 
K_0^{\tau_1}\e \gamma^{- 1} (1+\normk{\cI_0}{s+\sigma}) \,.
\end{aligned}
\end{equation}
Finally 
{\begin{equation}\label{bound delta12 beta infinito}
\begin{aligned}
| \Delta_{12} \beta |_p&\lesssim_{p}   K_0^{\tau_1} \e\gamma^{-1}\|i_1-i_2\|_{s_0+\s} \,,
\end{aligned}
\end{equation}}
for $p$ as in \eqref{ps0}.

\noindent
$(iii)$ Define 
\begin{equation}\label{def cal A infty pm 1}
\cB h (\vphi, x) := h(\vphi, x + \beta(\vphi, x)), \quad \text{with inverse} \quad \cB^{- 1} h(\vphi, y) = h(\vphi, y + \breve \beta(\vphi, y))\,. 
\end{equation}
Then, for any $\mathfrak{s}_0 \leq s \leq \bar s$, $\widetilde\cB_n^{\pm 1}$ converges pointwise in $H^s$ to $\cB^{\pm 1}$, namely $\lim_{n \to + \infty} \| \cB^{\pm 1 } h - \widetilde\cB_n^{\pm 1} h \|_s = 0$ for any $h \in H^s$. 
\end{lemma}

\begin{proof}
{\sc Proof of $(i)$.} We prove the Lemma arguing by induction. For $n = 0$, one has that $\widetilde\cB_0^{\pm 1} = \cB_0^{\pm 1}$ and we set $\alpha_0 := \beta_{0}$, $\breve \alpha_0 := \breve \beta_0$. Then the first estimate \eqref{stima beta n beta n - 1 trasporto} follows by \eqref{stime alpha n tilde alpha n} (applied with $n = 1$). The second statement in \eqref{stima beta n beta n - 1 trasporto} for $n = 0$ is empty. Now assume that the claimed statement holds for some $n \geq 0$ and let us prove it at the step $n + 1$. We prove the claimed statement for $\widetilde\cB_{n + 1}$ since the proof for the map $\widetilde\cB_{n + 1}^{- 1}$ is similar. Using that $\widetilde\cB_{n + 1} := \widetilde\cB_n \circ \cB_{n + 1}$, one computes that 
\begin{equation}\label{definizione induttiva beta n}
\begin{aligned}
\widetilde\cB_{n + 1} h(\vphi, x) = h(\vphi, x + \alpha_{n + 1}(\vphi, x))\,,\quad 
\alpha_{n + 1} := \alpha_{n } + \widetilde\cB_{n}[\beta_{n+1 }]\,. 
\end{aligned}
\end{equation}
We apply Lemma \ref{lemma.iterativo.straightening}.  
By the induction estimate \eqref{bound solo beta n} for $s = \mathfrak{s}_0$, 
using the ansatz \eqref{ansatz}, 
$| \alpha_n |_{\mathfrak{s}_0}^{k_0, \gamma} \lesssim K_0^{\tau_1} \e \gamma^{- 1}$. 
Then, by the smallness condition \eqref{small.V.as.AP}, 
we can apply Lemma \ref{lemma:LS norms} 
and \eqref{stime alpha n tilde alpha n}, 
\eqref{ansatz}, 
\eqref{small.V.as.AP}, 
\eqref{bound solo beta n}, 
obtaining that, for any $\mathfrak{s}_0 \leq s \leq \bar s$, 
\begin{equation}\label{stima beta n + 1 beta n}
\begin{aligned}
| \alpha_{n + 1} - \alpha_n |_s^{k_0, \gamma} 
& \lesssim_{s} | \beta_{n + 1} |_s^{k_0, \gamma} 
+ | \alpha_n |_s^{k_0, \gamma} | \beta_{n + 1} |_{\mathfrak{s}_0}^{k_0, \gamma}  
%\\ & \stackrel{\eqref{stime alpha n tilde alpha n}, \eqref{ansatz}, \eqref{condizione piccolezza rid trasporto}, \eqref{bound solo beta n}}{\lesssim_{s, k_0}} 
\lesssim_{s} K_{n + 1}^{\tau_1} K_n^{- \ta} \e \gamma^{- 1} (1+\normk{\cI_0}{s+\sigma}) , 
\end{aligned}
\end{equation}
which is \eqref{stima beta n beta n - 1 trasporto} at the step $n + 1$. 
The estimate \eqref{bound solo beta n} at the step $n + 1$
 follows by using a telescoping argument, since the series 
 $\sum_{n \geq 0} K_n^{\tau_1} K_{n - 1}^{- \ta} < \infty$.
Let us prove 
 estimate \eqref{bound delta12 beta n} at the step $n + 1$. 
 %follows by using \eqref{definizione induttiva beta n} and Lemma \eqref{lemma:LS norms}.
 First of all, by \eqref{definizione induttiva beta n}, we have that
 \[
 \alpha_{n+1}(i;\vphi,x)-\alpha_{n}(i;\vphi,x)=\beta_{n+1}(i;\vphi,x+\alpha_{n}(i;\vphi,x))\,,
 \]
 and so, taking two tori $i_1,i_2$, we have
 \[
 \begin{aligned}
 \Delta_{12}(\alpha_{n+1}-\alpha_n)&=\widetilde{B}_{n}(i_1)[\Delta_{12}\beta_{n+1}]
 \\&+
 \beta_{n+1}(i_2;\vphi,x+\alpha_{n}(i_1;\vphi,x))
 - \beta_{n+1}(i_2;\vphi,x+\alpha_{n}(i_2;\vphi,x))\,.
 \end{aligned}
 \]
 By the composition Lemma \ref{lemma:LS norms}, 
 estimates \eqref{stime alpha n tilde alpha n}$_{n+1}$, \eqref{estg12}$_{n+1}$
 to estimate $\beta_{n+1}$ and $\Delta_{12}\beta_{n+1}$, the inductive assumptions 
 \eqref{bound solo beta n}$_{n}$-\eqref{bound delta12 beta n}$_{n}$ 
 to bound $\alpha_n$ and $\Delta_{12}\alpha_n$,
 and using the smallness condition 
 \eqref{small.V.as.AP}, one gets
\[
\|\Delta_{12}(\alpha_{n+1}-\alpha_n)\|_{p}\lesssim_{p}
K_{n+1}^{\tau_1}K_{n}^{-\mathtt{a}} \e\gamma^{-1}\|i_1-i_2\|_{s_0+\s}\,,
\]
for all $p$ satisfying \eqref{ps0}.
By the estimate above, together with the  definition of the constant $\mathtt{a}$ in \eqref{tbta}, one deduces
\eqref{bound delta12 beta n} at the step $n+1$.

\medskip

\noindent
{\sc Proof of $(ii)$.} By item $(i)$ and in particular by estimate \eqref{stima beta n beta n - 1 trasporto} we have that the sequences $(\alpha_n)_{n \in \N}$ and $ (\breve \alpha_n)_{n \in \N}$ 
are Cauchy sequences with respect to the norm $| \cdot |_s^{k_0, \gamma}$ 
and therefore they converge. Let us define
\[
\beta:=\lim_{n\to\infty} \alpha_n\,, \quad \breve\beta:=\lim_{n\to\infty} \breve\alpha_n\,.
\]
Then \eqref{stima beta n - beta infty} follow by
the estimate \eqref{stima beta n beta n - 1 trasporto} and a telescoping argument. 
The bound \eqref{bound delta12 beta infinito} follows similarly using \eqref{bound delta12 beta n}.
The property of the quasi-periodic traveling wave functions follows from Lemma \ref{lemma.iterativo.straightening}, since this property is closed under limit.

\noindent
{\sc Proof of $(iii)$.} It follows by item $(ii)$
%, using the same arguments of the proof of Lemma 
%B6-$(i)$ in \cite{BHM}: 
reasoning as follows:
given $h \in H^s, \e_1 > 0$, 
there exists $N > 0$ (sufficiently large, depending on $\e_1,s,h$)   
such that $h_1 := \Pi_N^\bot h$ satisfies 
$\| \cB h_1 \|_s \leq \e_1 /4$, 
$\| \tilde \cB_n h_1 \|_s \leq \e_1 /4$
uniformly in $n$ (bound \eqref{bound solo beta n} is uniform in $n$).  
On the other hand, $h_0 := \Pi_N h$ satisfies 
\[ 
\begin{aligned}
\| (\cB - \tilde \cB_n) h_0 \|_s 
& \leq \int_0^1 \frac{d}{d\theta} \, 
h_0 \big( \varphi, x + \alpha_{n}(\varphi,x) + \theta [\beta(\varphi,x) - \alpha_n(\varphi,x)] \big) \, d\theta 
\\ 
& \lesssim_s 
\| \grad h_0 \|_s \| \beta - \alpha_n \|_{s_0} + \| \grad h_0 \|_{s_0} \| \beta - \alpha_n \|_s
\leq \e_1 / 2  
\end{aligned}
\]
for all $n \geq n_0$, for some $n_0$ depending on $\e_1, s, h$. 
The property of momentum preserving operators follows from Lemma \ref{lemma.iterativo.straightening}, since this property is closed under limit.
\end{proof}

\begin{lemma}\label{inclusione cantor diofantei}
$(i)$ The sequence $(\mathtt m_{1, n})_{n \in \N}$ given in Lemma 
\ref{lemma.iterativo.straightening} converges to 
$\tm_{1}\in\R^{d}$, uniformly for $(\omega, \th) \in \mathtt\Omega \times [\mathtt h_1, \mathtt h_2]$.
More precisely, it satisfies the bound 
\begin{equation}\label{lalla101}
\begin{aligned}
{\rm sup}_{(\omega, \th) \in \mathtt\Omega \times [\mathtt h_1, \mathtt h_2]} 
|\mathtt m_{1, n}(\omega, \th) - \tm_{1}(\omega, \th) |& \lesssim \e K_{n - 1}^{- \ta} \,,
\\
{\rm sup}_{(\omega, \th) \in \mathtt\Omega \times [\mathtt h_1, \mathtt h_2]} 
|\Delta_{12}(\mathtt m_{1, n}(\omega, \th) - \tm_{1}(\omega, \th)) |
&\lesssim \e K_{n - 1}^{- \ta}\,.
\end{aligned}
\end{equation}
$(ii)$ The following inclusion holds: $\tT\tC_{\infty}(\gamma,\tau)  \subseteq \cap_{n \geq 0} \mathtt{\Lambda}_{n}^{\rm T} $ (recall the definitions \eqref{tDtCn}, \eqref{set.nonres.tn}).
\end{lemma}

\begin{proof}
{\sc Proof of $(i)$.} 
By using \eqref{lalla100}  one has that the sequence $(\mathtt m_{1, n})_{n \in \N}$ 
is a Cauchy sequence satisfying 
\[
{\rm sup}_{(\omega, \th) \in \mathtt\Omega \times [\mathtt h_1, \mathtt h_2]} 
|\mathtt m_{1, n}(\omega, \th) - \tm_{1, n-1}(\omega, \th) | \lesssim \e K_{n - 2}^{- \ta} \,.
\]
Therefore, it converges to a vector $\tm_{1}\in \R^{d}$. Moreover, one has
\[
\begin{aligned}
|\tm_{1} - \tm_{1, n}| &= \sum_{k=n}^{\infty} |\tm_{1, k+1} - \tm_{1, k}| 
\leq \sum_{k=n}^{\infty} \eps K_{k-1}^{-\ta}
= \sum_{k=n-1}^{\infty} \eps K_{k}^{-\ta}
%= \sum_{k=n-1}^{\infty} \eps K_{0}^{-\ta\chi^{k}}
\lesssim\eps K_{n-1}^{-\ta}\,,
%\\
%&\leq \eps K_{n-1}^{-\ta}( 1+  \sum_{k=n}^{\infty} K_{0}^{-\ta(\chi^{k}-  \chi^{n-1} )})
%= \eps K_{n-1}^{-\ta}( 1+  \sum_{k=n}^{\infty} \frac{k^{2}}{k^{2}}K_{0}^{-\ta(\chi^{k}-  \chi^{n-1} )})
%\\
%&\leq  \eps K_{n-1}^{-\ta}( 1+ \sup_{k\geq n} \{k^{2} K_{0}^{-\ta(\chi^{k}-  \chi^{n-1} )} \} \sum_{k=n}^{\infty} \frac{1}{k^{2}}) \leq  C \eps K_{n-1}^{-\ta}\,,
\end{aligned}
\]
taking $K_0>0$ large.
This implies the first bound in \eqref{lalla101}. The second one follows reasoning similarly by using
\eqref{estm121}.

\medskip

\noindent
{\sc Proof of $(ii)$.} We prove the claimed inclusion by induction, i.e. we show that $\tT\tC_{\infty}(\gamma,\tau)  \subseteq \mathtt{\Lambda}_{n}^{\rm T} $ for any $n \geq 0$. 
For $n = 0$, the inclusion holds since $ \mathtt{\Lambda}_{0}^{\rm T}  := \mathtt\Omega \times [\mathtt h_1, \mathtt h_2]$ 
and $\tT\tC_{\infty}(\gamma,\tau)  \subseteq \mathtt\Omega \times [\mathtt h_1, \mathtt h_2]$ (see \eqref{tDtCn}).  
Now assume that 
$\tT\tC_{\infty}(\gamma,\tau)  \subseteq \mathtt{\Lambda}_{n}^{\rm T} $
 for some $n \geq 0$ 
and let us prove that 
$\tT\tC_{\infty}(\gamma,\tau)  \subseteq \mathtt{\Lambda}_{n+1}^{\rm T}$ . 
Let $(\omega, \th) \in \tT\tC_{\infty}(\gamma,\tau) $. 
By the induction hypothesis, $(\omega, \th)$ belongs to 
$ \mathtt{\Lambda}_{n}^{\rm T} $. 
Therefore, by item $(i)$, one has 
$|\mathtt m_{1, n}(\omega, \th) - \tm_{1}(\omega, \th) | \lesssim \e K_{n - 1}^{- \ta}$. 
Hence, for all $(\ell, j) \in \Z^{\nu}\times \Gamma^{*} \setminus \{(0,0)\}$, 
$|\ell| \leq K_n$, one has 
\[
\begin{aligned}
|(\omega - \tV \mathtt m_{1, n} )\cdot \ell | &
 \geq |(\omega - \tV \mathtt m_{1})\cdot \ell| - |\tV( \tm_{1} - \tm_{1, n})\cdot \ell | 
 \\
& \geq \frac{4 \gamma}{\langle \ell \rangle^\tau} - \e C K_n K_{n - 1}^{-\ta} 
\geq \frac{\gamma}{\langle \ell \rangle^\tau}
\end{aligned}
\]
provided $C K_n^{1+\tau} K_{n - 1}^{- \ta} \e \gamma^{- 1} \leq  1$.
%This holds for all $n \geq 0$ provided 
%\begin{equation} \label{provided trasporto}
%C N_0^{1+\tau} \e \g^{-1} \leq 1.
%\end{equation}
This confition is fulfilled by using
%Condition \eqref{provided trasporto} is fulfilled by taking $N_{0}$ sufficiently large, 
%and using 
\eqref{tbta} 
and the smallness condition \eqref{small.V.as.AP}. 
Thus,
one has that $(\omega, \th) \in \mathtt{\Lambda}_{n+1}^{\rm T}$, 
and the proof is concluded. 
\end{proof}

\begin{proof}[Proof of Proposition \ref{thm:as}]
 By recalling the definition \eqref{def widetilde cal An}, using \eqref{def cal Tn transport}, \eqref{coniugazione cal A n - 1 cal T n - 1}, one obtains  
\begin{equation}\label{A tilde n cal Tn}
\omega \cdot \partial_\vphi + \mathtt m_{1,n} \cdot \nabla + p_n \cdot \nabla = \cT_n = \widetilde\cB_{n - 1}^{- 1}\cT_0 \widetilde\cB_{n - 1}, \quad \forall \  (\omega, \th) \in \mathtt{\Lambda}_{n}^{\rm T}
\end{equation}
and by Lemma \ref{lemma tilde cal An}-$(i)$ one computes explicitely
\begin{equation}\label{A tilde n cal Tn 1}
\widetilde\cB_{n - 1}^{- 1}\cT_0 \widetilde\cB_{n - 1} 
= \omega \cdot \partial_\vphi 
+ \widetilde \cB_{n - 1}^{- 1} 
\big( \omega \cdot \partial_\vphi \alpha_{n - 1} + V 
+ V\cdot \nabla \alpha_{n - 1} \big) \cdot \nabla \,.
\end{equation}
Since $\widetilde\cB_{n - 1}$ is a change of variable, 
one has $\widetilde\cB_{n - 1}[c] = c$ for all constant $c \in \R$. 
Hence, by \eqref{A tilde n cal Tn}, \eqref{A tilde n cal Tn 1}, 
one obtains the identity
\begin{equation}\label{uguaglianza mn an betan}
\begin{aligned}
\mathtt m_{1, n} + \widetilde\cB_{n - 1} p_n = \omega \cdot \partial_\vphi \alpha_{n - 1} 
+ V +V \cdot \nabla \alpha_{n - 1}\,. 
\end{aligned}
\end{equation}

For any $(\om, \th) \in \tT\tC_{\infty}(\gamma,\tau)$, 
by Lemma \ref{inclusione cantor diofantei}, 
$\mathtt{m}_{1, n} \to \tm_{1}$ as $n \to \infty$. 
By \eqref{bound solo beta n}, 
\eqref{stima mathtt mn an}
and Lemma \ref{lemma:LS norms}, 
one has $| \widetilde \cB_{n-1} p_n |_{\mathfrak{s}_0} \lesssim | p_n |_{\mathfrak{s}_0} \to 0$
and 
\[
| \pa_\varphi \alpha_{n-1} - \pa_\varphi \beta |_{\mathfrak{s}_0} , \ 
| \grad \alpha_{n-1} - \grad \beta |_{\mathfrak{s}_0} 
\leq | \alpha_{n-1} - \beta |_{\mathfrak{s}_0+1} 
\to 0 \quad \ (n \to \infty).
\]
Hence, passing to the limit in norm $| \ |_{s_0}$ 
in the identity \eqref{uguaglianza mn an betan}, we obtain the identity
\begin{equation} \label{1303.1}
\ompaph \beta + V(\varphi,x)+V(\varphi,x)\cdot\grad \beta
%(\mathbb{I} + \grad \beta )   V(\varphi,x) 
= \tm_{1}
\end{equation}
in $H^{\mathfrak{s}_0}(\T^{\nu+d}_{*})$, 
and therefore pointwise for all $(\varphi,x) \in \T^{\nu+d}_{*}$,
for any $(\om,\th) \in  \tT\tC_{\infty}(\gamma,\tau)$.
As a consequence, 
\begin{align*}
\cB^{- 1} \cT_{0} \cB 
& = \omega \cdot \partial_\vphi  
+ \big\{ \cB^{-1} \big( \ompaph \beta + V+V\cdot\nabla\beta  %(\mathbb{I} + \grad \beta )  V 
\big) 
\big\} 
\cdot \grad 
%\\ & 
= \omega \cdot \partial_\vphi + \tm_{1} \cdot \nabla 
\end{align*}
for all $(\omega, \th) \in  \tT\tC_{\infty}(\gamma,\tau)$, 
which is \eqref{coniugazione nel teo trasporto}.
%The first estimate in \eqref{gtn.est.better}%, \eqref{stima tame cambio variabile rid trasporto} 
%follows from the estimates \eqref{stima beta n - beta infty} and by Lemma \ref{lemma:LS norms}. 
%
%It remains only to prove the second estimate in \eqref{gtn.est.better}. 
%
The estimates \eqref{gtn.est.better}, \eqref{stima tame cambio variabile rid trasporto} 
follow from the estimates \eqref{stima beta n - beta infty}-\eqref{bound delta12 beta infinito} 
and by Lemma \ref{lemma:LS norms}, while \eqref{tm.est} follow by Lemma \ref{inclusione cantor diofantei}.
Finally, the estimate for $\Delta_{12} \cB^{\pm 1}$ 
follows by using the estimate for $\Delta_{12} \beta$
using the mean value theorem and Lemma \ref{lemma:LS norms}. 
\end{proof}

%\begin{rmk}\label{rmk:tatamount}
%In view of Theorem 4.1 in \cite{FGMP19}, the result of Theorem \ref{thm:as} is tantamount to reduce the PDE
%\[
%\oo\cdot\pa_{\vphi} u+ V \cdot \nabla_{x} u=0
%\]
%to the equation with constant coefficients %(up to a exponentially small remind)
%\[
%\oo\cdot\pa_{\vphi} v+ \mathtt{m}_{1} \cdot \nabla v=0
%\]
%for $(\omega,\th)\in\tT\tC_{\infty}(\gamma,\tau)$ under the change of
%variable $\mathcal{B}$ defined in Lemma \ref{lemma tilde cal An} 
%which is composition operator induced bt the torus diffeomorphism $x\mapsto x + \beta(\vphi, x)$. 
%By Theorem  \ref{thm:as} and Corollary \ref{lemma tilde cal An} 
%we also deduce that actually $\beta$ is the solution of the equation
%\[
%\mathcal{B}^{-1}\big( \oo\cdot\partial_{\varphi} \beta + (\uno +\nabla_{x} \beta) V \big)= \mathtt{m}_{1} 
%\]
%\end{rmk}

%\begin{rmk}
%If the function $ V(\vphi,x)$ in \eqref{defX0} 
%is not  a quasi-periodic traveling wave, 
%the same kind of conjugation result holds requiring in
%\eqref{set.nonres.tn} 
%the non resonance conditions 
%\[
%| \omega \cdot \ell + \mathtt{m}_{1,\tn-1} j |\geq \gamma \jap{\ell}^{-\tau}, \  \forall\,
%(\ell, j) \in (\Z^{\nu} \times \Gamma^*) \setminus\{0\} \,, \ |(\ell,j)|\leq K_{n-1} \, . 
%\]
%\end{rmk}

\begin{proof}[{\bf Proof of Lemma \ref{conju.tr}}]\label{dim:lemma9.8}
First of all we recall that $\opw(\ii V\cdot \xi)= V\cdot \nabla +\frac{1}{2} {\rm div} V$.
We apply Theorem \ref{thm:as} and Lemma \ref{lemma tilde cal An} 
to  the transport operator  $  
%X_0 =
 \omega \cdot \partial_\varphi + V\cdot \nabla $, 
which has the form \eqref{defX0}.
% with $ p_0 =  V $.
 By % Lemma \ref{proprieta}-item 1, 
\eqref{stime tame V B inizio lineariz}, \eqref{ansatz_I0_s0},
and recalling that $K_0$ is chosen independently of  $\gamma$ in \eqref{tDtCn}, we have that 
the smallness conditions \eqref{small.V.as.AP}
hold by taking $\delta$ in \eqref{smallness.raddrizzo}
%$ K_0^{\tau_2} \varepsilon \gamma^{-1} $  
sufficiently small (recall also the ansatz \eqref{smallcondepsilon}).
%which is guaranteed by assumption \eqref{smallcondepsilon}.

Therefore there exist a constant $ \mathtt{m}_{1} \in \R^d $ and a 
quasi-periodic traveling wave $\beta(\vphi,x)$ (defined for all $(\omega,\th)\in \mathtt\Omega \times [\mathtt h_1, \mathtt h_2]$), 
 such that, for any $(\omega,\th) $ in $\tT\tC_{\infty}(\gamma,\tau)$
 we have, by Theorem \ref{thm:as} (see equation \eqref{coniugazione nel teo trasporto}) 
 and Lemma \ref{lemma tilde cal An}
 %Lemma \ref{lemma tilde cal An} and Remark \ref{rmk:tatamount},
\begin{equation}\label{pdeutile}
T_{\beta}^{-1} \big(\oo\cdot\partial_{\varphi} \beta + 
V+V\cdot\nabla\beta\big)
%(\uno +\nabla_{x} \beta)\big) V 
= \mathtt{m}_{1} \end{equation}
where we set $\mathcal{B}:=T_{\beta}$ as in \eqref{def:taubeta}.
The estimates 
\eqref{tm.est} and \eqref{stime tame V B inizio lineariz}  imply \eqref{beta.FGMP.est}  and \eqref{beta.FGMP.estBETA}.
Item $(i)$ of Lemma \ref{conju.tr} follows with $\breve{\beta}(\vphi,x)$ given by 
Lemma \ref{lemma tilde cal An}.

We consider the map $\breve{\mathcal{E}}$ in \eqref{coniugaLTR} and we recall that, in view of 
Remark \ref{rmk:EedEbreve}-(ii), one has 
$\breve{\mathcal{E}}^{-1}=\mathcal{E}$ with $\mathcal{E}$ 
given in \eqref{defcE}.
%with $\beta$ constructed above
%and we remark that it is the composition of 
%a torus diffeomorphism with a multiplication operator.
Therefore, recalling \eqref{coniugaLTR}-\eqref{defcE} and \eqref{def:acheck}-\eqref{def:acheck2},
by an explicit computation 
%recalling \eqref{defcE}-\eqref{def:CCalso} and \eqref{flussocorretto}-\eqref{belgio}, 
we have that
\[
\begin{aligned}
 \breve{\mathcal{A}} (\oo\cdot \partial_{\varphi}&+ \opw(\ii V\cdot \xi)) \breve{\mathcal{A}}^{-1}= 
 \mathcal{A}^{-1} (\oo\cdot \partial_{\varphi}+ \opw(\ii V\cdot \xi)) \mathcal{A}
 \\&
 = \mathcal{A}^{-1} (\oo\cdot \partial_{\varphi}+ V\cdot\nabla_{x}  + \frac{1}{2} {\rm div} V) \mathcal{A}
 \\&
% =T_{\beta}^{-1}\circ\frac{1}{\sqrt{\det(\uno+\nabla_x \beta)}}\circ 
% (\oo\cdot \partial_{\varphi}+ V\cdot\nabla_{x}  + \frac{1}{2} {\rm div} V) \circ \sqrt{\det(\uno+\nabla_x \beta)}
% \circ T_{\beta}
% \\&
% = T_{\beta}^{-1}\circ\Big(\oo\cdot \partial_{\varphi}+ V\cdot\nabla_{x}  + \frac{1}{2} {\rm div} V+
%\frac{\big[(\omega\cdot\pa_{\vphi}+V\cdot\nabla)\sqrt{\det(\uno+\nabla_x \beta)}\big]}{\sqrt{\det(\uno+\nabla_x \beta)}}
% \Big)T_{\beta}
% \\&
 =\omega\cdot\vphi+\mathtt{q_{1}}\cdot\nabla+\frac{1}{2}{\rm div}(\mathtt{q_{1}})\,,
 \end{aligned}
\]
where
\begin{equation}\label{belgio2}
\begin{aligned}
\mathtt{q_{1}}&:=T_{\beta}^{-1}\big(\omega\cdot\pa_{\vphi}\beta
+V+V\cdot\nabla\beta
\big)\,.
%+(\uno+\nabla\beta)V\big)\,,
%\\
%\mathtt{q_{2}}&:=\frac{1}{2}T_{\beta}^{-1}\Big({\rm div}V
%+{\rm Tr}\Big((\uno+\nabla\beta)^{-1}\big[(\omega\cdot\pa_{\vphi}+V\cdot\nabla)(\nabla\beta)\big]\Big)\Big)\,.
\end{aligned}
\end{equation}
%An explicit computation shows that
%${\rm div}(\mathtt{q_{1}})=2\mathtt{q_{2}}$. 
Therefore we conclude
\begin{equation*}
\begin{aligned}
 \mathcal{A}^{-1} (\oo\cdot \partial_{\varphi}+ \opw(\ii V\cdot \xi)) \mathcal{A}&=\omega\cdot\pa_{\vphi}+
 \mathtt{q_{1}}\cdot \nabla+\frac{1}{2}{\rm div}(\mathtt{q_{1}})\equiv
 \omega\cdot\pa_{\vphi}+\opw(\ii\mathtt{q_{1}}\cdot\x)
 \\&
 \stackrel{\eqref{pdeutile}, \eqref{belgio2}}{=}
 \omega\cdot\pa_{\vphi}+ \opw( \ii \mathtt{m}_{1}  \cdot \xi)\,,
 \end{aligned}
\end{equation*}
which implies formula \eqref{coniugaLTR}.
Since $\cL_{\rm TR}$ is Hamiltonian (Definition \ref{operatoreHam}), 
and the map $\breve\cE$ is symplectic, we have that $\breve\cE\cL_{\rm TR}\breve{\cE}^{-1}$ 
is Hamiltonian as well.  
The map  $\breve\cE$ is momentum preserving since the generator preserves the subspace 
$S_{\tV}$
by Corollary \ref{lemma tilde cal An}. % (recalling the form of $b$ in \eqref{gianduiotto}). 
Finally, by Lemma \ref{lemma:buonaposFlussi} one has that 
$\breve\cE$ satisfies \eqref{stimamappaEperp}.
\end{proof}

\end{document}